\documentclass[12pt,oneside,reqno]{amsbook}

\usepackage{amsmath,amsthm,amsfonts,amssymb,amsrefs,amscd,graphicx,nicefrac}
\allowdisplaybreaks  

\usepackage[titletoc]{appendix}
\usepackage[active]{srcltx}
\usepackage{dsfont,amscd}
\usepackage{mathrsfs}
\usepackage{amstext}
\usepackage[capposition=bottom]{floatrow}
\usepackage{setspace}
\usepackage{color,xcolor}
\usepackage{enumerate}
\usepackage{fancyhdr,url}
\usepackage[top=1in, bottom=1in, left=1in, right=1in]{geometry}

\usepackage{comment} 
\usepackage[english]{babel}
\usepackage{todonotes}
\usepackage{tikz,tikz-cd}
\usepackage{etoolbox}
\usepackage{xfrac,enumitem,soul}

\usepackage{imakeidx}
\makeindex[title=Index of Concepts,columns=3]

\usepackage[hidelinks]{hyperref}
\usepackage{setspace}\setdisplayskipstretch{}  
\definecolor{OliveGreen}{rgb}{0,0.6,0}

\newcommand{\Natural}{\mathbb{N}}
\newcommand{\Integer}{\mathbb{Z}}
\newcommand{\Rational}{\mathbb{Q}}
\newcommand{\Real}{\mathbb{R}}
\newcommand{\Complex}{\mathbb{C}}

\DeclareMathOperator{\im}{im}

\DeclareMathOperator{\diam}{diam}

\DeclareMathOperator{\Avg}{Avg}
\DeclareMathOperator{\Thin}{Thin}

\DeclareMathOperator{\dist}{dist}
\DeclareMathOperator{\Conv}{Conv}

\DeclareMathOperator{\Lip}{Lip}
\DeclareMathOperator{\CAT}{CAT}

\DeclareMathOperator{\Ord}{Ord}

\DeclareMathOperator{\Int}{Int}
\DeclareMathOperator{\Cl}{Cl}
\DeclareMathOperator{\Span}{Span}
\DeclareMathOperator{\FSR}{FSR}
\DeclareMathOperator{\Cone}{Cone}

\newcommand {\sr}{\stackrel}

\newcommand {\eqv}{\equiv}

\newcommand {\ol}{\overline}

\newcommand {\ra}{\rightarrow}
\newcommand {\Ra}{\Rightarrow}

\newcommand {\La}{\Leftarrow}

\newcommand {\Lra}{\Leftrightarrow}
\newcommand {\ral}{\longrightarrow}

\newcommand{\ha}{\hookrightarrow}
\newcommand{\hal}{\lhook\joinrel\longrightarrow}

\newcommand {\vep}{\varepsilon}
\newcommand {\vphi}{\varphi}
\newcommand {\del}{\partial}
\newcommand {\Ld}{\Lambda}
\newcommand {\ld}{\lambda}

\newcommand {\al}{\alpha}

\newcommand{\wt}{\widetilde}
\newcommand{\wh}{\widehat}

\newcommand {\bfr}{\begin{flushright}}
\newcommand {\efr}{\end{flushright}}
\newcommand {\bfl}{\begin{flushleft}}
\newcommand {\efl}{\end{flushleft}}

\newcommand {\nn} {\nonumber}

\newcommand {\txt}{\textrm}
\newcommand {\bd}{\begin{document}}
\newcommand {\ed}{\end{document}}

\newcommand {\bt}{\begin{tikzcd}}
\newcommand {\et}{\end{tikzcd}}
\newcommand {\be}{\begin{equation}}
\newcommand {\ee}{\end{equation}}
\newcommand {\bea}{\begin{eqnarray}}
\newcommand {\eea}{\end{eqnarray}}
\newcommand {\ba}{\begin{aligned}}
\newcommand {\ea}{\end{aligned}}

\newcommand{\bbib}{\begin{thebibliography}{99}}
\newcommand{\ebib}{\end{thebibliography}}
\newcommand {\bab}{\begin{abstract}}
\newcommand {\eab}{\end{abstract}}

\newcommand {\bc}{\begin{center}}
\newcommand {\ec}{\end{center}}
\newcommand {\bit}{\begin{itemize}}
\newcommand {\eit}{\end{itemize}}
\newcommand{\blue}{\textcolor{blue}}
\newcommand{\red}{\textcolor{red}}
\newcommand{\magenta}{\textcolor{magenta}}

\def\A{\mathcal{A}}\def\B{\mathcal{B}}\def\C{\mathcal{C}}\def\D{\mathcal{D}}\def\E{\mathcal{E}}\def\F{\mathcal{ F}}\def\G{\mathcal{G}}\def\H{\mathcal{H}}\def\I{\mathcal{I}}
\def\J{\mathcal{J}}\def\K{\mathcal{K}}\def\L{\mathcal{L}}\def\M{\mathcal{M}}\def\N{\mathcal{N}}\def\O{\mathcal{ O}}\def\P{\mathcal{P}}\def\Q{\mathcal{Q}}\def\R{\mathcal{R}}
\def\S{\mathcal{S}}\def\T{\mathcal{T}}\def\U{\mathcal{U}}\def\V{\mathcal{V}}\def\W{\mathcal{W}}\def\X{\mathcal{ X}}\def\Y{\mathcal{Y}}\def\Z{\mathcal{Z}}
\def\abar{\bar{a}}\def\bbar{\bar{b}}\def\cbar{\bar{c}}\def\dbar{\bar{d}}\def\ebar{\bar{e}}\def\fbar{\bar{f}}
\def\gbar{\bar{g}}\def\ibar{\bar{i}}\def\jbar{\bar{j}}\def\kbar{\bar{k}}\def\lbar{\bar{l}}
\def\mbar{\bar{m}}\def\nbar{\bar{n}}\def\obar{\bar{o}}\def\pbar{\bar{p}}\def\qbar{\bar{q}}\def\rbar{\bar{r}}
\def\sbar{\bar{s}}\def\tbar{\bar{t}}\def\ubar{\bar{u}}\def\vbar{\bar{v}}\def\wbar{\bar{w}}
\def\xbar{\bar{x}}\def\ybar{\bar{y}}\def\zbar{\bar{z}}
\def\Abar{\bar{A}}\def\Bbar{\bar{B}}\def\Cbar{\bar{C}}\def\Dbar{\bar{D}}\def\Ebar{\bar{E}}\def\Fbar{\bar{F}}
\def\Gbar{\bar{G}}\def\Hbar{\bar{H}}\def\Ibar{\bar{I}}\def\Jbar{\bar{J}}\def\Kbar{\bar{K}}\def\Lbar{\bar{L}}
\def\Mbar{\bar{M}}\def\Nbar{\bar{N}}\def\Obar{\bar{O}}\def\Pbar{\bar{P}}\def\Qbar{\bar{Q}}\def\Rbar{\bar{R}}
\def\Sbar{\bar{S}}\def\Tbar{\bar{T}}\def\Ubar{\bar{U}}\def\Vbar{\bar{V}}\def\Wbar{\bar{W}}
\def\Xbar{\bar{X}}\def\Ybar{\bar{Y}}\def\Zbar{\bar{Z}}
\def\0bar{\bar{0}}\def\1bar{\bar{1}}\def\2bar{\bar{2}}\def\3bar{\bar{3}}\def\4bar{\bar{4}}\def\5bar{\bar{5}}
\def\6bar{\bar{6}}\def\7bar{\bar{7}}\def\8bar{\bar{8}}\def\9bar{\bar{9}}
\def\cprod{\mathop{\prod\hspace{-0.457cm}\times}}
\def\tprod{\mathop{\prod\hspace{-0.457cm}\otimes}}
\def\dsum{\mathop{\sum\hspace{-0.458cm}\oplus}}
\def\Tr{\txt{Tr}}\def\tr{\txt{tr}}\def\str{\txt{Str}}
\def\bde{\begin{description}}\def\ede{\end{description}}

\newtheorem{thm}{Theorem}[section]
\newtheorem{lmm}[thm]{Lemma}
\newtheorem{dfn}[thm]{Definition}
\newtheorem{crl}[thm]{Corollary}
\newtheorem{prp}[thm]{Proposition}
\newtheorem{example}[thm]{Example}
\newtheorem{problem}[thm]{Problem}
\newtheorem{algorithm}[thm]{Algorithm}
\newtheorem{recall}[thm]{Recall}
\newtheorem{exercise}[thm]{Exercise}
\newtheorem{caution}[thm]{Caution}
\newtheorem{claim}[thm]{Claim}
\newtheorem{fact}[thm]{Fact}
\newtheorem{solution}[thm]{Solution}
\newtheorem{conjecture}[thm]{Conjecture}
\newtheorem{assumption}[thm]{Assumption}

\theoremstyle{definition}
\newtheorem{dfns}{Definitions}[section]
\newtheorem{examples}{Examples}[section]
\newtheorem{rmk}{Remark}[section]
\newtheorem{rmks}{Remarks}[section]
\newtheorem{problems}{Problems}[section]
\newtheorem{algorithms}{Algorithms}[section]
\newtheorem{question}{Question}[section]
\newtheorem{rquestion}{Research question}[section]
\newtheorem{questions}{Questions}[section]
\newtheorem{rquestions}{Research questions}[section]
\newtheorem{recalls}{Recalls}[section]
\newtheorem{exercises}{Exercises}[section]
\newtheorem{note}{Note}[section]
\newtheorem{notes}{Notes}[section]
\newtheorem{notation}{Notation}[section]
\newtheorem{convention}{Convention}[section]
\newtheorem{claims}{Claims}[section]
\newtheorem{facts}{Facts}[section]
\newtheorem{assumptions}{Assumptions}
\newtheorem{summary}{Summary}
\newtheorem{consequence}{Consequence}
\newtheorem{consequences}{Consequences}

\theoremstyle{definition}\newtheorem*{thm*}{Theorem}
\newtheorem*{lmm*}{Lemma}                     \newtheorem*{lmms*}{Lemmas}
\newtheorem*{dfn*}{Definition}                \newtheorem*{dfns*}{Definitions}
\newtheorem*{crl*}{Corollary}
\newtheorem*{prp*}{Proposition}
\newtheorem*{example*}{Example}               \newtheorem*{examples*}{Examples}
\newtheorem*{rmk*}{Remark}                    \newtheorem*{rmks*}{Remarks}
\newtheorem*{problem*}{Problem}               \newtheorem*{problems*}{Problems}
\newtheorem*{algorithm*}{Algorithm}           \newtheorem*{algorithms*}{Algorithms}
\newtheorem*{question*}{Question}             \newtheorem*{questions*}{Questions}
\newtheorem*{rquestion*}{Research question}
\newtheorem*{rquestions*}{Research questions}
\newtheorem*{recall*}{Recall}                 \newtheorem*{recalls*}{Recalls}
\newtheorem*{exercise*}{Exercise}             \newtheorem*{exercises*}{Exercises}
\newtheorem*{note*}{Note}                     \newtheorem*{notes*}{Notes}
\newtheorem*{notation*}{Notation}
\newtheorem*{convention*}{Convention}         \newtheorem*{conventions*}{Conventions}
\newtheorem*{caution*}{Caution}
\newtheorem*{claim*}{Claim}                   \newtheorem*{claims*}{Claims}
\newtheorem*{fact*}{Fact}                     \newtheorem*{facts*}{Facts}
\newtheorem*{assumption*}{Assumption}         \newtheorem*{assumptions*}{Assumptions}
\newtheorem*{summary*}{Summary}
\newtheorem*{consequence*}{Consequence}       \newtheorem*{consequences*}{Consequences}

\numberwithin{section}{chapter}
\numberwithin{equation}{chapter}

\makeatletter
\let\sv@thm\@thm
\def\@thm{\let\indent\relax\sv@thm}
\makeatother

\makeatletter
\let\@afterindenttrue\@afterindentfalse
\makeatother

\setul{1.5pt}{.4pt}

\AtBeginDocument{\renewcommand{\contentsname}{Table of contents}}        

\makeatletter
\patchcmd{\@tocline}
  {\hfil}
  {\leaders\hbox{\,.\,}\hfil}
  {}{}
\makeatother

\AtBeginDocument{\renewcommand{\bibname}{References}}

\makeatletter
\renewcommand{\l@chapter}{\@tocline{0}{8pt plus1pt}{0pt}{}{\bfseries}}
\makeatother

\DeclareMathOperator{\bigutimes}{\textstyle\bigcup\bigutimeskern\times}
\newcommand{\bigutimeskern}{%
  \mkern-17.5mu
  \mathchoice{}{}{\mkern0.2mu}{\mkern0.5mu}%
  }
\DeclareMathOperator{\utimes}{\cup\utimeskern\times}
\newcommand{\utimeskern}{%
  \mkern-13mu
  \mathchoice{}{}{\mkern0.2mu}{\mkern0.5mu}%
}

\newenvironment{myquote}[1]%
  {\list{}{\leftmargin=#1\rightmargin=#1}\item[]}%
  {\endlist}
\makeatletter
\@namedef{subjclassname@2020}{%
  \textup{2020} Mathematics Subject Classification}
\makeatother

\begin{document}
\pagestyle{plain}

\frontmatter
\makeatletter
\patchcmd{\@maketitle}{\newpage}{}{}{}  
\makeatother
\title{Metric Geometry of Finite Subset Spaces}
\author{\large ~\\[0.5em]
Earnest Akofor \\[3em]

Department of Mathematics\\[0.5em]
Syracuse University, Syracuse, NY 13244, USA \\[0.5em]
Email: eakofor@syr.edu\\[8em]

(PhD Dissertation)

%
%
%

}

\maketitle
\thispagestyle{empty}
{
\begin{center}
  \vspace*{3in}
  \copyright \hspace{5mm} Copyright 2020 Earnest Akofor\\[1.5em]
This thesis or part of it may be freely copied, transferred, edited, reproduced, shared, or published (for example, as part of an article or alone) in any form without permission.
\end{center}
}
\newpage
\vskip 2cm
\begin{center}
In memory of my father, King Barnabas Ngocho of Cameroon.
\end{center}

\newpage          
\vskip 2cm
{
\begin{center}
\textbf{Acknowledgments}\\[1.5em]
\end{center}
\begin{quote}
I would like to express profound thanks to Leonid Kovalev for his continuous insightful suggestions, technical checking/advice, and editing of the dissertation text throughout its preparation. I would also like to thank William Wylie for carefully reading through the thesis with several useful comments/suggestions to improve its readability.
\end{quote}
}
\newpage 
\newpage
\thispagestyle{empty}
  \begin{center}
    {\large Abstract}\\[1em]
  \end{center}
{
\begin{myquote}{0.2in}
If $X$ is a (topological) space, the $n$th finite subset space of $X$, denoted by $X(n)$, consists of $n$-point subsets of $X$ (i.e., nonempty subsets of cardinality at most $n$) with the quotient topology induced by the unordering map $q:X^n\ra X(n)$, $(x_1,\cdots,x_n)\mapsto\{x_1,\cdots,x_n\}$. That is, a set $A\subset X(n)$ is open if and only if its preimage $q^{-1}(A)$ is open in the product space $X^n$.

Given a space $X$, let $\H(X)$ denote all homeomorphisms of $X$. For any class of homeomorphisms $\C\subset\H(X)$, the $\C$-geometry of $X$ refers to the description of $X$ up to homeomorphisms in $\C$. Therefore, the topology of $X$ is the $\H(X)$-geometry of $X$. By a ($\C$-) geometric property of $X$ we will mean a property of $X$ that is preserved by homeomorphisms of $X$ (in $\C$). Metric geometry of a space $X$ refers to the study of geometry of $X$ in terms of notions of metrics (e.g., distance, or length of a path, between points) on $X$. In such a study, we call a space $X$ metrizable if $X$ is homeomorphic to a metric space.

Naturally, $X(n)$ always inherits some aspect of every geometric property of $X$ or $X^n$. Thus, the geometry of $X(n)$ is in general richer than that of $X$ or $X^n$. For example, it is known that if $X$ is an orientable manifold, then (unlike $X^n$) $X(n)$ for $n>1$ can be an orientable manifold, a non-orientable manifold, or a non-manifold. In studying geometry of $X(n)$, a central research question is ``If $X$ has geometric property $P$, does it follow that $X(n)$ also has property $P$?''. A related question is ``If $X$ and $Y$ have a geometric relation $R$, does it follow that $X(n)$ and $Y(n)$ also have the relation $R$?''.

Extensive work exists in the literature on the richness of the geometry of $X(n)$. Nevertheless, despite the fact that the spaces $X(n)$ considered in those investigations are metrizable (which is the case if and only if $X$ is itself metrizable) the important role of metrics has been mostly ignored. Consequently, the existing results mostly elucidate topological aspects of the geometry of $X(n)$.

The main goal of this thesis is to attempt to answer the above research question(s) for several geometric properties, with metrics playing a significant role (hence the title phrase ``\emph{Metric Geometry of ...}''). Some of the questions are relatively easy and will be answered completely. However, a question such as ``If a normed space $X$ is an absolute Lipschitz retract, does it follow that $X(n)$ is also an absolute Lipschitz retract?'' appears to require considerable effort and will be answered only partially. By the definition of an absolute Lipschitz retract, establishing the existence of Lipschitz retractions $X(n)$ $\ra$ $X(n-1)$ for all $n\geq 1$ would be a partial positive answer to this question.

Among other things, we will prove the following. If $X$ is a metrizable space, then so is $X(n)$. If a metric space $X$ is a snowflake, quasiconvex, or doubling then so is $X(n)$. If two spaces $X$ and $Y$ are (Lipschitz) homotopy equivalent, then so are $X(n)$ and $Y(n)$. If $X$ is a normed space (which is Lipschitz $k$-connected for all $k\geq 0$), then $X(n)$ is Lipschitz $k$-connected for all $k\geq 0$. If $X$ is a normed space, there exist (i) H\"older retractions $X(n)\ra X(n-1)$, (ii) Lipschitz retractions $X(n)\ra X(1),X(2)$, and (iii) Lipschitz retractions $X(n)\ra X(n-1)$ when $\dim X<\infty$ or $X$ is a Hilbert space.
\end{myquote}
}
{\let\thefootnote\relax}
\newpage
\doublespacing 
\tableofcontents
\listoffigures

\mainmatter
\bibliographystyle{plain}                     

\pagestyle{headings}
\setcounter{chapter}{-1}
\chapter{Introduction}\label{Intro}
\section{Basic setup}
If $X$ is a (topological) space, Definition \ref{PtSetTop}, the $n$th \emph{finite subset space} of $X$, denoted by $X(n)$, consists of $n$-point subsets of $X$ (i.e., nonempty subsets of cardinality at most $n$) with the \emph{quotient topology} induced by the \emph{unordering map} $q:X^n\ra X(n)$, $(x_1,\cdots,x_n)\mapsto\{x_1,\cdots,x_n\}$. That is, a set $A\subset X(n)$ is open if and only if its preimage $q^{-1}(A)$ is open in $X^n$ (with the product topology, Definition \ref{ProdTop}). Due to this description in terms of $X^n$, the finite subset space $X(n)$ is also called the $n$th \emph{symmetric product} of $X$ (see \cite{borsuk-ulam1} for example). A detailed discussion of finite subset spaces will begin in chapter \ref{FSS}.

Given a space $X$ (which is defined up to \emph{homeomorphism}, Definition \ref{HomeoImbed}), let $\H(X)$ denote all homeomorphisms of $X$. For any class of homeomorphisms $\C\subset\H(X)$, the $\C$-geometry of $X$ refers to the description of $X$ up to homeomorphisms in $\C$. In particular, the topology of $X$ is the same as the $\H(X)$-geometry of $X$.

For example, \emph{Lipschitz geometry} of a metric space $X$ involves the study of those properties of $X$ that are invariant under Lipschitz homeomorphisms of $X$, where we say a map of (between) metric spaces $f:X\ra Y$ is \ul{\emph{Lipschitz}} if there is a number $\ld\geq 0$, (the least of) which we denote by $\Lip(f)$, such that $d(f(x),f(x'))\leq \ld d(x,x')$ for all $x,x'\in X$. Similarly, \emph{biLipschitz geometry} of a metric space $X$ involves the study of those properties of $X$ that are invariant under biLipschitz homeomorphisms of $X$, where we say a map of metric spaces $f:X\ra Y$ is \ul{\emph{biLipschitz}} if there is a number $\ld>0$ such that $d(x,x')/\ld\leq d(f(x),f(x'))\leq \ld d(x,x')$ for all $x,x'\in X$.

By a ($\C$-) \emph{geometric property/invariant} of $X$ we will mean a property of $X$ that is preserved by homeomorphisms of $X$ (in $\C$). For example, it is a topological (or $\H(X)$-geometric) property whether or not a subspace $A\subset X$ is a \ul{\emph{retract}}, i.e., whether or not the identity map $id_A:A\ra A$, $a\mapsto a$ extends to a continuous map $r:X\ra A$ (which we call a \ul{\emph{retraction}} of $X$ onto $A$) in the sense that $r|_A=id_A$. As in \cite{BBI}, \emph{metric geometry} of a space $X$ refers to the study of geometry of $X$ in terms of notions of metrics (e.g., distance, or length of a path, between points) on $X$. In such a study, we call a space $X$ \emph{metrizable} if $X$ is homeomorphic to a \emph{metric space}, Definition \ref{MetricSp}.

\section{Motivation and related work}
{\flushleft From} a scientific point of view, when measurements or observations are made and a \emph{dataset} (i.e., data set) is obtained, one may choose to describe \emph{relationships} between variables in the dataset in a \emph{symmetric} way (i.e., without the need to express one of the variables as a function of the rest). Assume we have $p$ variables $z_1,...,z_p$ in the dataset, where the \emph{value set} of $z_i$ is $\Z_i$, and suppose the \emph{joint value set} $\Z_1\times\cdots\times\Z_p$ is contained in a metric space $X$. If the number of observations in the dataset is $n$, then the whole dataset is a point $x=(x_1,...,x_n)$, $x_i=(z_{1i},...,z_{pi})$, of $X^n$. Relationships between the $p$ variables can be described in a symmetric way through their appearance/effects in the \emph{distribution of clusters} in some \emph{rearrangement} of the dataset. We can view a retraction $r:X(n)\ra X(k)$ as a rule $r\circ q:X^n\sr{q}{\ral}X(n)\sr{r}{\ral} X(k)$ for making a $k$-cluster rearrangement (i.e., a rearrangement with at most $k$ clusters) from the dataset of $n$ observations. (footnote\footnote{For the special case where $k=n$ and $r=id$, the distribution of the trivial rearrangement $q(x)=\{x_1,...,x_n\}$ is given by the counts $n_i(x)=\big|\{j:x_j=x_i\}\big|$. In particular, if $|q(x)|=n$, i.e., $q(x)\in X(n)\backslash X(n-1)$, then $q(x)$ has a uniform distribution with $n_i(x)=1$ for all $i$.}).
We may require the retraction to be \emph{Lipschitz}, which ensures \emph{stability} of the clustering rule under metric-perturbations of the data, in the sense that the perturbation $d_H(r\circ q(x),r\circ q(x'))$ of a $k$-cluster rearrangement is no greater than a scalar multiple $\Lip(r\circ q)d(x,x')$ of the perturbation $d(x,x')$ of the dataset.

From a more natural point of view, $X(n)$ always inherits some aspect of every geometric property of $X$  or $X^n$. Thus, the geometry of $X(n)$ is in general richer than that of $X$ or $X^n$. For example, it is known, \cite[Remark on page 1123]{tuffley2002}, that if $X$ is an orientable manifold such as the circle $S^1$, then (unlike $X^n$) $X(n)$ for $n>1$ can be an orientable manifold (e.g., $S^1(3)$ is homeomorphic to $S^3$), a non-orientable manifold (e.g., $S^1(2)$ is homeomorphic to the M\"{o}bius band), or a non-manifold (e.g., $S^1(k)$, $k\geq 4$, is not a manifold). In studying geometry of $X(n)$, the following are central research questions.
\begin{rquestions*}
(i) If $X$ has geometric property $P$, does it follow that $X(n)$ also has property $P$? (ii) If $X$ and $Y$ have a geometric relation $R$, does it follow that $X(n)$ and $Y(n)$ also have the relation $R$?
\end{rquestions*}

Note that question (i) above can also be reversed in order to deduce a property of $X$ from a property of $X(n)$ for some $n>1$, for example, see \cite{IllanMart2017} and related work in \cite{CoronaEtal2017}.

Extensive work exists in the literature on the richness of the geometry of $X(n)$. It was shown in \cite{borsuk-ulam1} that $[0,1](n)$ imbeds in $\Real^n$ for $n=1,2,3$ but not for $n\geq 4$. For $n\geq 4$ however, $[0,1](n)$ can be imbedded in $\Real^N$ for a sufficiently large $N$ (see \cite{kovalev2015} and the references therein). Homotopy types of the spaces $S^1(n)$ have been considered in \cite{borsuk1949,bott1952,wu-wen47,tuffley2002,chinen-koyama2010}, where \cite{bott1952} is a correction of \cite{borsuk1949}. Related discussions can also be found in \cite{schori68,chinen2015}. Computation of the Euler characteristic of $X(n)$, for a manifold $X$, can be found in \cite{salazar2004}.

Nevertheless, despite the fact that the spaces $X(n)$ considered in these investigations are metrizable (which is the case if and only if $X$ is itself metrizable, Proposition \ref{FSSmzty}) the important role of metrics has been mostly ignored. Consequently, the existing results in the literature mostly elucidate topological aspects of the geometry of $X(n)$. Studies in which metrical properties have been considered include \cite{BorovIbrag2009,BorovEtal2010,chinen2018}.

\section{Objectives and contributions}
The main goal of this thesis is to attempt to answer the above research question (i) for several geometric properties with metrics playing a significant role in establishing general results (hence the phrase ``\emph{Metric Geometry}'' in the thesis title). Some of the questions, such as those on quasiconvexity and Lipschitz connectedness, are relatively easy and will be answered completely. However, a question such as ``\emph{If a normed space $X$ is an absolute Lipschitz retract, does it follow that $X(n)$ is also an absolute Lipschitz retract?}'' appears to require considerable effort and will be answered only partially. This question (Question \ref{MajQuest} in this thesis) appeared as Problem 1.4 in the collection of open problems ``\emph{AimPL: Mapping theory in metric spaces}'', published by the American Institute of Mathematics and available at \url{http://aimpl.org/mappingmetric} (see \cite{ProbList}). By the definition of an absolute Lipschitz retract (as a metric space $Z$ such that every containing metric space $Y\supset Z$ admits a Lipschitz retraction $r:Y\ra Z$), establishing the existence of Lipschitz retractions $X(n)$ $\ra$ $X(n-1)$ for all $n\geq 1$ would be a partial positive answer to the above question. Existing work in this direction includes \cite{kovalev2016,bacac-kovalev2016,akofor2019}.

Among other things, we will prove the following. In Proposition \ref{FSSmzty} we verify that if $X$ is a metrizable space, then so is $X(n)$. In Proposition \ref{SnwfSubSps}, Theorem \ref{QConvThm} (plus Corollary \ref{MetricQCQC}), and Lemma \ref{FinSubDoubl} we show (respectively) that if a metric space $X$ is a snowflake, quasiconvex, or doubling then so is $X(n)$. (footnote\footnote{Theorem \ref{QConvThm} (which proves that if $X$ is geodesic then $X(n)$ is $2$-quasiconvex) improves and generalizes \cite[Theorem 4.1]{BorovEtal2010} in which it was established that $\Real(n)$ is $4^n$-quasiconvex. Lemma \ref{FinSubDoubl} (which proves that if $X$ is a doubling metric space then so is $X(n)$) partly solves \cite[Problem 4.1]{BorovEtal2010} in which it was asked to show that $\Real(n)$ is doubling.}). In Lemma \ref{LipHEfssRecall}, we observe that if two spaces $X$ and $Y$ are (Lipschitz) homotopy equivalent, then so are $X(n)$ and $Y(n)$. In Theorem \ref{LipConSuff5}, we establish that if $X$ is a normed space (which is $24$-Lipschitz $k$-connected for all $k\geq 0$, Corollary \ref{LipConSuff4b}), then $X(n)$ is Lipschitz $k$-connected for all $k\geq 0$. Based on this, we deduce in Theorem \ref{FinDimLipRet} that if $X$ is a finite-dimensional normed space, then there exist Lipschitz retractions $X(n)\ra X(n-1)$.

It was asked in \cite[Question 3.4]{kovalev2016} whether Lipschitz retractions $r_n:X(n)\ra X(n-1)$ exist for all $n\geq2$ when $X$ is a Banach space. A related question in \cite[Remark 3.4]{bacac-kovalev2016} similarly asked whether strictly convex or uniformly convex Banach spaces admit such Lipschitz retractions. Also, in \cite[Question 3.2]{kovalev2016}, and later in \cite[Remark 3.4]{bacac-kovalev2016}, it was asked whether $\Lip(r_n)$ can be bounded above by a constant that is independent of $n$.

The results of this thesis provide partial answers and tools of investigation towards answering the above two questions. Based directly on \cite{akofor2019}, regarding the first question above (i.e., \cite[Question 3.4]{kovalev2016} and its counterpart in \cite[Remark 3.4]{bacac-kovalev2016}), we prove that if $X$ is a normed space, then we have (i) concrete Lipschitz retractions $X(3)\ra X(2)\ra X$ (Theorem \ref{LrThmV3}), (ii) Lipschitz retractions $X(n)\ra X,X(2)$ (Theorem \ref{LrThmVn2}), and (iii) retractions $X(n)\rightarrow X(n-1)$ that are H{\"o}lder-continuous on bounded sets (Theorem \ref{MainThm}). Similarly, Theorem \ref{LipConstGr} provides a partial negative answer to the second question above (i.e., \cite[Question 3.2]{kovalev2016} and its counterpart in \cite[Remark 3.4]{bacac-kovalev2016}) by showing that if $X$ is a normed space of dimension $\geq 2$, then $\Lip(r_n)$ must grow with $n$.

\section{Descriptive summary}\label{DescpSumm}
In the rest of the discussion, we begin with preliminaries in chapter/appendix \ref{Prelims}, where we introduce fundamental concepts and develop notation/terminology that will be used in the main chapters to follow. The results of section \ref{PrelimsTSC} on topological space concepts (to be used throughout), of section \ref{PrelimsAT} on algebraic topology (to be used in chapter \ref{LipConn} and in chapter \ref{FSRP2}, section \ref{GrwLipConst}), and of section \ref{PrelimsMET} on metrical extension theorems (to be used in chapter \ref{GQC}) are all essential to the ensuing discussion. On the other hand, the results of section \ref{PrelimsTET} on topological extension theorems and of section \ref{PrelimsMS} on metrizability of spaces are not strictly essential, and so a reader who is familiar with the results of those two sections can simply skim through that material.

The description of finite subset spaces begins in chapter \ref{FSS}, a thorough understanding of which is essential in the remaining chapters. Hausdorff distance and other tools such as diameter and minimum separation are discussed (in section \ref{FSShd} and in section \ref{FSSfs}, respectively), finite subset spaces of a metrizable space are shown to be metrizable (in section \ref{FSSmf}), a few inherited properties of finite subset spaces are noted (in section \ref{FSSsi}), and an interesting property of subsets of any metric space, namely, "gap reducing property" is established (in section \ref{FSSgi}).

In chapter \ref{GQC} we begin studying some aspects of metric geometry of finite subset spaces (for which section \ref{PrelimsMET} on metrical extension theorems is essential). We prove, in section \ref{GQCiq}, that finite subset spaces of a quasiconvex metric space are themselves quasiconvex. In sections \ref{GQCcg} and \ref{GQCcq}, we give full characterizations of geodesics and quasigeodesics in finite subset spaces, and discuss some consequences. Also discussed, in section \ref{GQCmb}, is the failure of metrical convexity, and of the binary intersection property, in finite subset spaces of normed spaces of dimension two and higher.

In chapter \ref{LipConn} we consider Lipschitz connectedness of finite subset spaces (for which section \ref{PrelimsAT} on algebraic topology is essential). We discuss Lipschitz homotopy and contractibility (in section \ref{LipConnLH1}), and prove that if metric spaces $X$ and $Y$ have the same (Lipschitz) homotopy type, then so do $X(n)$ and $Y(n)$. Using a notion of conical Lipschitz contractibility (from section \ref{LipConnLC1}), we review (in section \ref{LipConnLC2}) Lipschitz $k$-connectedness of a normed space and also establish Lipschitz $k$-connectedness of finite subset spaces of a normed space.

Chapter \ref{FSRP1} begins a study of the existence of Lipschitz retractions $X(n)\ra X(n-1)$, called the FSR (finite subset retraction) property, as an attempt at partially answering the question of whether finite subset spaces of an absolute Lipschitz retract must also be absolute Lipschitz retracts. We prove (in sections \ref{LipRetV3}, \ref{LipRetVn12}, \ref{Retractions}) that when $X$ is a normed space, there exist retractions $X(n)\rightarrow X(k)$, $n>k\geq 1$, that are (i) H{\"o}lder for $k=n-1$, (ii) Lipschitz for $k=1,2$, and (iii) Lipschitz for $k=n-1$ if $\dim X<\infty$ or $X$ is a Hilbert space. We also show (in section \ref{Selections}) that these retractions map each finite subset of $X$ into its convex hull, and discuss some consequences of this observation.

Chapter \ref{FSRP2} continues the study of the FSR property on the existence of Lipschitz retractions $X(n)\ra X(n-1)$. Using Lipschitz $k$-connectedness of finite subsets of normed spaces from chapter \ref{LipConn}, we prove (in section \ref{FSRP2lr}) that if $X$ is a finite dimensional normed space, then we have Lipschitz retractions $X(n)\ra X(n-1)$. Using existing results on the homotopy types of the spaces $S^1(n)$, we show (in section \ref{GrwLipConst}) that if $X$ is a normed space of dimension 2 or more, then any sequence of Lipschitz retractions $r_n:X(n)\ra X(n-1)$ satisfies the bound $4\pi\sqrt{2}\Lip(r_n)\geq n$. We also consider (in sections \ref{FSRP2fs}, \ref{FSRP2sn}) some non-existence results and further questions on the FSR property involving metric spaces more general than normed spaces. Section \ref{FSRP2cf} considers some counter examples and facts, not involving finite subset spaces, that can help in answering questions on the FSR property.

Apart from preliminaries on topology, the appendices contain additional information that is equally relevant to the discussion in the main chapters, but may require the introduction of terminology that differs significantly from that used in the main chapters.

\section{Prerequisites}
All relevant basic knowledge of topology is provided as preliminaries/appendices. However, throughout, we assume basic knowledge of set theory (see \cite{enderton1977,goldrei1996}), as well as basic knowledge of groups (mostly abelian groups), rings (mostly $\Integer$, $\Rational$, $\Real$, $\Complex$), and modules (mostly vector spaces over $\Real,\Complex$) found for example in \cite{artin-book}. We also assume basic knowledge of measure theory (especially calculus of absolutely continuous functions) that can be found in \cite{folland-book,wheeden-zygmund}, and of functional analysis (especially geometry of locally convex spaces) that can be found in \cite{conway}.

\section{Notation and conventions}
Our preferred set containment operations are $\subset$ for non-strict containment, and $\subsetneq$ for strict containment. However, when dealing with a partially ordered set (poset), we will for convenience also use $\subseteq$ (in place of $\subset$) for non-strict containment. That is $\subset,\subseteq$ will each denote ``non-strict containment'' while $\subsetneq$ will denote ``strict containment''.

Capital letters such as $A,B,\cdots$ will denote sets, meanwhile "curly" capital letters such as $\A,\B,\cdots$ will denote collections/families of sets. Lower case letters will denote elements of a set, i.e., an arbitrary element of a set $A$ is denoted by $a\in A$. If $f:X\ra Y$, $x\mapsto f(x)$ is a \emph{mapping} of sets (i.e., a rule that assigns to each element $x\in X$ a \emph{single} element $f(x)\in Y$), then the value of $f$ at $x\in X$ is denoted by $f(x)$ or $f_x$ or $f^x$, depending on convenience. The same notation applies to indexed collections. For example, we can write $\{C(\al)\}_{\al\in A}$ or $\C=\{C_\al\}_{\al\in A}$ or $\C=\{C^\al\}_{\al\in A}$. For countable collections, we will often write $\C=\{A_i\}_{i\in C}$, for a countable set $C$. That is, we mostly (but not always) use $\al,\beta,\gamma,\cdots$ for arbitrary indices and mostly (but not always) use $i,j,k,l,\cdots$ for countable indices. In cases where the countable set $C\subset\Integer$, we will often (but not always) use $m,n,r,s,t,\cdots$, e.g., $\C=\{A_n\}_{n=1}^\infty$.

Any specialized notation/convention that deviates from the above basic rules will be explicitly stated (along with any associated terminology) in definitions. Such specialized notation is common in all chapters/sections, and some chapters/sections adopt and use notation introduced in ``earlier'' sections (some important cases of which have already been indicated in section \ref{DescpSumm} above). If a chapter relies on preliminary notation/terminology other than that of section \ref{PrelimsTSC}, then this will be explicitly mentioned at the beginning of that chapter.

We will sometimes abbreviate certain frequent phrases as follows. \ul{~$:=$~} stands for \ul{``defined to be''}. \ul{~$\Ra$~} stands for \ul{``implies''}. \ul{\emph{iff} (or $\iff$)} stands for \ul{``if and only if''}. \ul{\emph{s.t.}} or \ul{~$:$~} stands for \ul{``such that''}. \ul{\emph{a.e.}} stands for \ul{``almost everywhere''}. \ul{\emph{a.e.f.}} stands for \ul{``all except finitely many''}. \ul{\emph{wlog}} stands for \ul{``without loss of generality''}. \ul{\emph{nbd}} stands for \ul{``neighborhood''}. \ul{\emph{resp.}} stands for \ul{``respectively''}. \ul{\emph{equiv.}} stands for \ul{``equivalently''}. \ul{\emph{const.}} stands for \ul{``constant''}. \ul{~$\ll$~} stands for \ul{``much less than''}. \ul{~$\gg$~} stands for \ul{``much greater than''}.

The symbol $\cong$ will denote both \ul{``isomorphism''} (of sets, groups, rings, modules, ...) and \ul{``homeomorphism''} (of spaces). Meanwhile, the symbol $\simeq$ will denote both \ul{``homotopy''} (of continuous maps and chain maps) and \ul{``homotopy equivalence''} (of spaces and chain complexes).

In a metric space $X$, if $A\subset X$, we will occasionally use both $B_R(A)$ and $N_R(A)$ to denote the set $\{x\in X:\dist(x,A)<R\}$, and similarly both $\ol{B}_R(A)$ and $\ol{N}_R(A)$ will denote the set $\{x\in X:\dist(x,A)\leq R\}$. The actual choice in any given occasion will depend on (notational) convenience.


\chapter{Finite Subset Spaces}\label{FSS}
We introduce subset spaces and discuss (in sections \ref{FSSmf}, \ref{FSSsi}) a few of their inherited properties such as metrizability in Proposition \ref{FSSmzty}, completeness in Proposition \ref{FSScmpl}, compactness in Corollary \ref{FSScmpc}, and snowflakiness in Proposition \ref{SnwfSubSps}. In sections \ref{FSShd}, \ref{FSSfs} we review basic properties of Hausdorff distance $d_H$, including Lipschitz properties (in Lemmas \ref{DiamCont}, \ref{DeltaCont}) of diameter and minimum distance (between elements of a finite set) that will be useful in later chapters. Also introduced in section \ref{FSSgi}, Proposition \ref{GRPprp}, is an interesting property of all metric spaces, namely, ``\emph{gap reducing property}''.

\begin{dfn}[\textcolor{blue}{\index{Finite subset space}{Finite subset spaces}, Unordering map, Total finite subset space}]
Let $X$ be a space and $n\geq 1$. The $n$th \ul{finite subset space} (or $n$th \ul{symmetric product}) of $X$ is the \ul{quotient space}
\bea
\textstyle X(n):={X^n\over\sim}=\Big\{x=\{x_1,...,x_n\}:(x_1,...,x_n)\in X^n\Big\}=\{x\subset X:|x|\leq n\},\nn
\eea
where the equivalence relation $\sim$ is defined as follows: For $(x_1,...,x_n),(y_1,...,y_n)\in X^n$,
\bea
(x_1,...,x_n)\sim (y_1,...,y_n)~~\txt{if}~~\{x_1,...,x_n\}=\{y_1,...,y_n\}.\nn
\eea
The \ul{unordering map} is the \ul{quotient map} $q:X^n\ra X(n)$, $(x_1,...,x_n)\mapsto\{x_1,...,x_n\}$.

The \ul{total finite subset space} (i.e., the space of all finite subsets) of $X$ is the union
\bea
\textstyle FS(X):=\bigcup_{n=1}^\infty X(n)=\{x\subset X:1\leq|x|<\infty\},\nn
\eea
where we declare a set $A\subset FS(X)$ open if $A\cap X(n)$ is open in $X(n)$ for each $n\geq 1$ (equivalently, $A\subset FS(X)$ is closed if $A\cap X(n)$ is closed in $X(n)$ for each $n\geq 1$).
\end{dfn}
Note that by the universal property of quotient maps (UPQM), every continuous map $f_n:X^n\ra X^{n-1}$ satisfying ``$u\sim v$ ~$\Ra$~ $q_{n-1}\circ f_n(u)=q_{n-1}\circ f_n(v)$ for all $u,v\in X^n$'' induces a unique continuous map $\wt{f}_n:X(n)\ra X(n-1)$ such that $\wt{f}_n\circ q=q\circ f_n$.
\bc\bt
\cdots\ar[r]& X^n\ar[d,"q_n"]\ar[rr,"f_n"]&& X^{n-1}\ar[d,"q_{n-1}"]\ar[r]& \cdots\\
\cdots\ar[r]& X(n)\ar[rr,dashed,"\exists!~\wt{f}_n"]&& X(n-1)\ar[r]& \cdots
\et\ec
Similarly, every continuous map $f:X\ra Y$, through the associated map $f_n:X^n\ra Y^n$, $(x_1,...,x_n)\mapsto(f(x_1),...,f(x_n))$, induces a unique continuous map $\wt{f}_n:X(n)\ra X(n)$ such that $\wt{f}_n\circ q^X=q^Y\circ f_n$ (and so, we say the unordering map $q:X^n\ra X(n)$ is \ul{natural}).
\bc\bt
X^n\ar[d,"q_n^X"]\ar[rr,"f_n"]&& Y^n\ar[d,"q_n^Y"]\\
X(n)\ar[rr,dashed,"\exists!~\wt{f}_n"]&& Y(n)
\et\ec

\begin{lmm}[\textcolor{blue}{Limit-open sets, Limit-closed sets, \index{Limit topology}{Limit topology}}, \textcolor{OliveGreen}{Metrizable union}]\label{MtzblUnion}
Let $(X_1,d_1)\subset (X_2,d_2)\subset\cdots$ be a sequence of metric spaces such that $d_{i+1}|_{X_i\times X_i}=d_i$ (i.e., $X_i$ is a \emph{subspace} of $X_{i+1}$ for each $i\geq 1$). Let $X:=\bigcup X_i$ be the space such that $A\subset X$ is open (resp. closed) iff $A\cap X_i$ is open (resp. closed) in $X_i$ for all $i$. Equivalently, open (resp. closed) sets $A\subset X$ have the form $A=\bigcup A_i$, with $A_i=A\cap X_i$ and open (resp. closed) in $X_i$. Let us call the open (resp. closed) sets \ul{limit-open} (resp. \ul{limit-closed}), and the topology, the \ul{limit-topology} on $X$.

Then the union space $X$ is metrizable, and the topology of $X$ is induced by the metric $d:=\bigcup d_i$ on $X$ given, for $x,x'\in X$, by
\bea
\textstyle d(x,x'):=d_i(x,x'),~~~~\txt{if}~~x,x'\in X_i.\nn
\eea
\end{lmm}
\begin{proof}
{\flushleft\ul{The $d$-topology contains the limit-topology}:} Let $A\subset X$ be limit-open, i.e., {\small $A_i:=A\cap X_i=\bigcup_{\al\in\Gamma_i}B^{d_i}_{r^i_\al}(c^i_\al)=\left(\bigcup_{\al\in\Gamma_i}B^d_{r^i_\al}(c^i_\al)\right)\cap X_i=\wt{A}_i\cap X_i$}, where {\small $\wt{A}_i:=\bigcup_{\al\in\Gamma_i}B^d_{r^i_\al}(c^i_\al)$} is a $d$-open set for each $i$. Since $A_i\subset A_{i+1}$ (as $X_i\subset X_{i+1}$), the $\wt{A}_i$ can be chosen such that $\wt{A}_i\subset\wt{A}_{i+1}$, say by replacing $\wt{A}_i$ with $\wt{A}'_i:=\wt{A}_i\cap\wt{A}_{i+1}$ (which is still $d$-open). Observe that {\small $A=\bigcup A_i=\bigcup\left(\wt{A}_i\cap X_i\right)\subset\left(\bigcup\wt{A}_i\right)\cap\left(\bigcup X_i\right)=\bigcup\wt{A}_i$}. On the other hand, if {\small $x\in\bigcup\wt{A}_i\subset X=\bigcup X_i$}, then for some $i,j$ such that $j\geq i$, we have
{\small $x\in\wt{A}_i\cap X_j\subset\wt{A}_j\cap X_j=A\cap X_j=A_j\subset A$}, and so {\small $\bigcup\wt{A}_i\subset A$}. Hence, {\small $A=\bigcup \wt{A}_i$}, which is $d$-open.

{\flushleft\ul{The limit-topology contains the $d$-topology}:} Given a $d$-ball $B_r(c)=B^d_r(c):=\{x:d(c,x)<r\}$ in $X$, let $A_i:=B_r(c)\cap X_i=\{x_i\in X_i:d(x_i,c)<r\}$. Let $c\in X_{i_c}$ (which implies $c\in X_j$ for all $j\geq i_c$). Then for all $j\geq i_c$,
\bea
&&\textstyle A_j=\{x_j\in X_j:d(x_j,c)<r\}=\{x_j\in X_j:d_j(x_j,c)<r\}=B^{d_j}_r(c),\nn\\
&&\textstyle~~\Ra~~B_r^d(c)=\bigcup_{i\geq 1}A_i~=~\bigcup_{1\leq i<i_c}A_i~\cup~\bigcup_{j\geq i_c}B^{d_j}_r(c),\nn
\eea
where for $1\leq i<i_c$, {\small $A_i=B^{d_{i_c}}_r(c)\cap X_i\subset X_i\subset X_{i_c}$} is open in $X_i$ since $B^{d_{i_c}}_r(c)$ is open in $X_{i_c}$.

\end{proof}

Alternative descriptions of the topology of $X(n)$ can be found in \cite{illanes-nadler}. We will show in Proposition \ref{FSSmzty} that $X(n)$, and hence $FS(X)$ (by Lemma \ref{MtzblUnion}), is metrizable whenever $X$ is metrizable.

\section{Hausdorff distance: Properties and estimates}\label{FSShd}
Throughout, we denote the collection of all subsets of a space $X$ by $\P(X)$, and nonempty subsets of $X$ by $\P^\ast(X)$. If $X$ is a metric space, $\B(X)$ and $\B^\ast(X)$ will denote bounded subsets and nonempty bounded subsets respectively. Similarly, $\B_c(X)$ and $\B_c^\ast(X)$ will denote bounded closed subsets and nonempty bounded closed subsets respectively.

If $A,B\subset X$ and $\vep>0$, we will denote the \emph{distance} between $A,B$ by $\dist(A,B):=\inf\{d(a,b):a\in A,b\in B\}$, the \emph{open $\vep$-neighborhood} of $A$ by $N_\vep(A):=\{x\in X~|~\dist(x,A)<\vep\}$, and the \emph{closed $\vep$-neighborhood} of $A$ by $\ol{N}_\vep(A):=\{x\in X~|~\dist(x,A)\leq\vep\}$.

The \emph{cardinality} of a set $A$ will be denoted by $|A|$.

\begin{dfn}[\textcolor{blue}{\index{Hausdorff distance $d_H$}{Hausdorff distance} $d_H$, Set distance $\dist_H$}] Let $X$ be a metric space. The Hausdorff distance on $\P^\ast(X)$ is the map $d_H:\P^\ast(X)\times\P^\ast(X)\ra[0,\infty]$ given by
{\small\begin{align}
\label{HausDistRep1}
d_H(A,B)&:=\max\left\{\sup_{a\in A}\inf_{b\in B}d(a,b)~,~\sup_{b\in B}\inf_{a\in A}d(a,b)\right\}\eqv\max\left\{\wt{d}(A,B),\wt{d}(B,A)\right\}\\
&=~\sup_{a,b}~\max\big\{\dist(a,B),\dist(A,b)\big\},\nn
\end{align}}
where ~$\wt{d}(A,B):=\sup_{a\in A}\dist(a,B)$ ~and ~$\wt{d}(B,A):=\sup_{b\in B}\dist(b,A)$.

If $\A,\B\subset \P^\ast(X)$, the $d_H$-metric space distance between $\A$ and $\B$ will be denoted by
\bea
\textstyle \dist_H(\A,\B):=\inf_{A\in\A,B\in\B}d_H(A,B).\nn
\eea
\end{dfn}
Alternative expressions for $d_H$ are given by part (b) of the following lemma.
\begin{lmm}\label{HausDistRmk}
Let $X$ be a metric space. Then we have the following.
\begin{enumerate}[leftmargin=0.8cm]
\item[(a)] $d_H$ induces a metric topology on the bounded closed subsets $\B_c^\ast(X)$ of $X$.
\item[(b)] Let $A,B\subset X$ be bounded sets, and $A_\vep:=N_\vep(A)$, $B_\vep:=N_\vep(B)$, $\ol{A}_\vep:=\ol{N}_\vep(A)$, $\ol{B}_\vep:=\ol{N}_\vep(B)$ the open/closed $\vep$-neighborhoods of $A$ and $B$ in $X$. Then
\begin{align}
\label{HausDistRep2}d_H(A,B)&=\inf\left\{\vep:A\subset\ol{B}_\vep,~B\subset\ol{A}_\vep\right\}=\inf\left\{\vep:A\cup B\subset\ol{A}_\vep\cap\ol{B}_\vep\right\}\\
\label{HausDistRep3}&=\inf\left\{\vep:A\subset B_\vep,~B\subset A_\vep\right\}=\inf\left\{\vep:A\cup B\subset A_\vep\cap B_\vep\right\}.
\end{align}
\item[(c)] In the representation (\ref{HausDistRep2}) of $d_H$ above, $d_H(A,B)$ is ``achieved'' in the sense that with $\rho:=d_H(A,B)$, we have $A\cup B\subset\ol{A}_\rho\cap\ol{B}_\rho$.
\end{enumerate}
\end{lmm}

\begin{proof}
\begin{enumerate}[leftmargin=0.8cm]
\item[(a)] Observe that $d_H(A,B)=d_H(B,A)$, and $d_H(A,B)=0$ $\iff$ $\ol{A}=\ol{B}$. Also,
\begin{align}
&\textstyle\wt{d}(A,B)\sr{(\ref{HausDistRep1})}{=}\sup_a\dist(a,B)\leq\sup_a\inf_c[d(a,c)+\dist(c,B)]\leq \sup_a\inf_c[d(a,c)+\wt{d}(C,B)]\nn\\
&\textstyle~~~~=\wt{d}(A,C)+\wt{d}(C,B),~~~~\txt{and similarly},~~\wt{d}(B,A)\leq\wt{d}(B,C)+\wt{d}(C,A).\nn
\end{align}
\item[(b)] Observe that $\wt{d}(A,B)=\inf\{\vep:\wt{d}(A,B)\leq\vep\}=\inf\{\vep:A\subset\ol{B}_\vep\}$, which implies
\begin{align}
d_H(A,B)&=\max\big\{\wt{d}(A,B),\wt{d}(B,A)\big\}=\inf\big\{\vep:~\wt{d}(A,B)\leq\vep,~\wt{d}(B,A)\leq\vep\big\}\nn\\
&=\inf\left\{\vep:A\subset\ol{B}_\vep,~B\subset\ol{A}_\vep\right\}=\inf\left\{\vep:A\cup B\subset\ol{A}_\vep\cap\ol{B}_\vep\right\},\nn
\end{align}
where we note that ~$\ol{(A\cup B)}_\vep=\ol{A_\vep\cup B_\vep}=\ol{A}_\vep\cup \ol{B}_\vep$ ~~and ~~$\ol{(A\cap B)}_\vep\subset\ol{A_\vep\cap B_\vep}\subset\ol{A}_\vep\cap \ol{B}_\vep$. Similarly, $\wt{d}(A,B)=\inf\{\vep:\wt{d}(A,B)<\vep\}=\inf\{\vep:A\subset B_\vep\}$, which implies
\begin{align}
d_H(A,B)&=\max\big\{\wt{d}(A,B),\wt{d}(B,A)\big\}=\inf\big\{\vep:~\wt{d}(A,B)<\vep,~\wt{d}(B,A)<\vep\big\}\nn\\
&=\inf\left\{\vep:A\subset B_\vep,~B\subset A_\vep\right\}=\inf\left\{\vep:A\cup B\subset A_\vep\cap B_\vep\right\}.\nn
\end{align}

\item[(c)] Let $\rho:=d_H(A,B)$. Then the closed sets $C_n:=\ol{A}_{\rho+1/n}\cap\ol{B}_{\rho+1/n}$ satisfy $C_n\supset C_{n+1}\supset A\cup B$ for all $n\geq 1$, and so $A\cup B\subset\bigcap\left(\ol{A}_{\rho+1/n}\cap\ol{B}_{\rho+1/n}\right)=\ol{A}_\rho\cap\ol{B}_\rho$.
\end{enumerate}
\end{proof}

A discussion of the above properties of $d_H$ can also be found in \cite[Sec. 7.3.1, p.252]{BBI}.

\begin{lmm}[\textcolor{blue}{\index{Diameter}{Diameter}}, \textcolor{OliveGreen}{2-Lipschitz property of diameter}]\label{DiamCont}
If $X$ is a metric space, the map $\diam:\B^\ast(X)\ra \Real$ given by $\diam(A):=\sup\{d(a,a'):a,a'\in A\}$ is $2$-Lipschitz with respect to the Hausdorff distance $d_H$.
\end{lmm}
\begin{proof}
Fix $\vep>0$ and bounded sets $A,B\in \B^\ast(X)$. Then there exist $a,a'\in A$ and $b,b'\in B$ such that $\diam(A)\leq d(a,a')+\vep$, $d(a,b)\leq \dist(a,B)+\vep\leq d_H(A,B)+\vep$, and $d(a',b')\leq \dist(a',B)+\vep\leq d_H(A,B)+\vep$. These three inequalities (together with the triangle inequality) in turn imply
\begin{equation*}
\diam(A)\leq d(a,b)+d(b,b')+d(a',b')+\vep\leq 2d_H(A,B)+\diam(B)+3\vep.
\end{equation*}
Hence, $|\diam(A)-\diam(B)|\leq 2 d_H(A,B)+3\vep$.
\end{proof}

\begin{lmm}[\textcolor{OliveGreen}{Factorization property}]\label{FactznProp}
$d_H(A\cup B,C\cup D)\leq \max\big(d_H(A,C),d_H(B,D)\big)$.
\end{lmm}
\begin{proof}
Let $\rho:=\max\big(d_H(A,C),d_H(B,D)\big)$, and pick any $\vep'>\rho$. Then because $d_H(A,C)\leq\rho<\vep'$ and $d_H(B,D)\leq\rho<\vep'$, we have the containments $A\subset \ol{C}_{\vep'}$, $C\subset \ol{A}_{\vep'}$, $B\subset \ol{D}_{\vep'}$, $D\subset\ol{B}_{\vep'}$, which imply $A\cup B\subset \ol{C}_{\vep'}\cup \ol{D}_{\vep'}=\ol{(C\cup D)}_{\vep'}$ and $C\cup D\subset \ol{A}_{\vep'}\cup \ol{B}_{\vep'}=\ol{(A\cup B)}_{\vep'}$. It follows that $d_H(A\cup B,C\cup D):=\inf\left\{\vep: A\cup B\subset\ol{(C\cup D)}_\vep,~C\cup D\subset\ol{(A\cup B)}_\vep\right\}\leq\vep'$ for all $\vep'>\rho$. Hence, ~$d_H(A\cup B,C\cup D)\leq\rho=\max\big(d_H(A,C),d_H(B,D)\big)$.
\end{proof}

\begin{lmm}[\textcolor{OliveGreen}{Hausdorff distance to a subset}]\label{SubsetDist} If $A\subset B$, then
\bea
\textstyle d_H(A,B)=\max_{b\in B\backslash A}\dist(b,A)\leq d_H(A,B\backslash A).\nn
\eea
\end{lmm}
\begin{proof}
The inequality is a corollary of the preceding result since $d(A,B)=d(A\cup A,A\cup B\backslash A)$. For the equality, $d_H(A,B)=\max\{\max_{a\in A}\dist(a,B),\max_{b\in B}\dist(b,A)\}=\max_{b\in B}\dist(b,A)=\max_{b\in B\backslash A}\dist(b,A)$.
\end{proof}

\begin{lmm}[\textcolor{OliveGreen}{Further estimates}]\label{DiamEstimates}
Let $X$ be a metric space and $A,B,C,D\in\B^\ast(X)$. Then
\begin{enumerate}[leftmargin=0.8cm]
\item[(a)] $d_H(A,B)\leq\max\big(\diam(A),\diam(B)\big)+\dist(A,B)$.
\item[(b)] $\diam(A\cup B)\leq \diam(A)+\diam(B)+\dist(A,B)$.
\item[(c)] $\diam(A\cup B)\leq \min\big(\diam(A),\diam(B)\big)+2d_H(A,B)$.
\item[(d)] $|\dist(A,B)-\dist(C,D)|\leq d_H(A,C)+d_H(B,D)$.
\end{enumerate}
\end{lmm}
\begin{proof}
(a) For any $\vep>0$, there exist $a_\vep,a'_\vep\in A$ and $b_\vep,b'_\vep\in B$ such that
\begin{align}
&\textstyle \sup_{a\in A}\dist(a,B)<\dist(a_\vep,B)+\vep,~~~~\sup_{b\in B}\dist(B,a)<\dist(A,b_\vep)+\vep,\nn\\
&\textstyle \dist(a'_\vep,B)<\dist(A,B)+\vep,~~~~\dist(A,b'_\vep)<\dist(A,B)+\vep.\nn
\end{align}
Therefore, with $\dist(a,B)=\min_bd(a,b)\leq \min_b[d(a,a')+d(a',b)]=d(a,a')+\dist(a',B)$,
{\small\begin{align}
&\textstyle d_H(A,B)=\max\left(\sup_{a\in A}\dist(a,B),\sup_{b\in B}\dist(A,b)\right)<\max\Big(\dist(a_\vep,B),\dist(A,b_\vep)\Big)+\vep\nn\\
&~~\leq\max\Big(|\dist(a_\vep,B)-\dist(a'_\vep,B)|+\dist(a'_\vep,B)~,~|\dist(A,b_\vep)-\dist(A,b'_\vep)|+\dist(A,b'_\vep)\Big)+\vep\nn\\
&~~\leq \max\Big(d(a_\vep,a'_\vep)~,~d(b_\vep,b'_\vep)\Big)+\dist(A,B)+2\vep.\nn
\end{align}}(b) Similarly, if $\diam(A\cup B)\neq\diam A$ and $\diam(A\cup B)\neq\diam B$, then for any $\vep>0$, there exist $a_\vep\in A$, $b_\vep\in B$ such that
\begin{align}
&\diam(A\cup B)\leq d(a_\vep,b_\vep)+\vep\leq d(a_\vep,a)+d(a,b)+d(b,b_\vep)+\vep,~~\txt{for all}~~a\in A,b\in B,\nn\\
&~~~~\leq \diam(A)+d(a,b)+\diam(B)+\vep,~~\txt{for all}~~a\in A,b\in B,\nn\\
&~~\Ra~~\diam(A\cup B)\leq\diam(A)+\dist(A,B)+\diam(B)+\vep.\nn
\end{align}
(c) Using the triangle inequality and the fact that diameter is 2-Lipschitz,
\begin{align}
&\diam(A\cup B)\leq \diam(B)+|\diam(A\cup B)-\diam(B)|\leq \diam(B)+2d_H(A\cup B,B)\nn\\
&~~\leq \diam(B)+2d_H(A,B),~~\txt{and by symmetry},~~\diam(A\cup B)\leq \diam(A)+2d_H(A,B),\nn\\
&~\Ra~\diam(A\cup B)\leq \min\big(\diam(A),\diam(B)\big)+2d_H(A,B).\nn
\end{align}
(d)  Finally, recall that by the triangle inequality, we have
\bea
|\dist(A,B)-\dist(C,D)|\leq |\dist(A,B)-\dist(C,B)|+|\dist(C,B)-\dist(C,D)|,\nn
\eea
and given any functions $\{f_a,g_c:Z\ra\Real~|~a\in A,c\in C\}$ (where $Z$ is any set), as done in the proof of Lemma \ref{MWetLmm},
\begin{align}
\textstyle\left|\inf_af_a(x)-\inf_cg_c(y)\right|\leq d_H\big(\{f_a(x)\},\{g_c(y)\}\big).\nn
\end{align}
Using this, we see that ~$|\dist(A,B)-\dist(C,B)|=\big|\inf_a\dist(a,B)-\inf_c\dist(c,B)\big|$ \\
$\leq d_H\big(\{\dist(a,B)\}_{a\in A},\{\dist(c,B)\}_{c\in C}\big)\leq d_H(A,C)$.
\end{proof}

Note that by writing $A\cup B=(A\backslash B)\cup(A\cap B)\cup(B\backslash A)$ in part (b) above, we get
\bea
\diam(A\cup B)\leq\diam(A\backslash B)+\max[\diam(A\cap B),\dist(A,B)]+\diam(B\backslash A).\nn
\eea

\section{Metrizability of finite subset spaces}\label{FSSmf}
Recall that a space is said to be metrizable iff it is homeomorphic to a metric space.

\begin{rmk}[\textcolor{OliveGreen}{Non-openness of the unordering quotient map}]
Let $X$ be a metric space and $d_H:X(n)\times X(n)\ra\Real$ the Hausdorff distance. Then the unordering map $q:X^n\ra \big(X(n),d_H\big)$, $(x_1,...,x_n)\mapsto\{x_1,..,x_n\}$ is in general not an open map.
\end{rmk}
\begin{proof}
Consider $q:\Real^3\ra\Real(3)$. Fix $u:=(1,1,2)$. Let $v:=q(u)=\{1,2\}$, and $E:=\{v_\vep:=\{1,2,2+\vep\}:0<\vep<\infty\}$. Then $d_H(v,v_\vep)=\vep$. Meanwhile, for all $\vep>0$,
\begin{align}
&\dist(u,q^{-1}(v_\vep))\geq 1,~~~\Ra~~~q^{-1}(E)\subset B_1(u)^c,~~~\Ra~~~B_1(u)\subset [q^{-1}(E)]^c=q^{-1}(E^c),\nn\\
&~~~\Ra~~~q\big(B_1(u)\big)\subset qq^{-1}(E^c)=E^c,\nn
\end{align}
which shows $q$ is not open since $q\big(B_1(u)\big)$ is not an open neighborhood of $q(u)$.
\end{proof}

Note that using Corollary \ref{OpnssCrit4} the above result is immediate, because with $u_\vep:=(1,2,2+\vep)$,
\bea
d_H\big(q(u),q(u_\vep)\big)=d_H(v,v_\vep)\ra 0~~~~\txt{but}~~~~\dist\Big(u,q^{-1}\big(q(u_\vep)\big)\Big)=\dist\Big(u,q^{-1}(v_\vep)\Big)\nrightarrow 0.\nn
\eea

\begin{dfn}[\textcolor{blue}{Left-complete relation, Right-complete relation, Complete relation}]
Let $A,B$ be sets and consider the mappings $p_1,p_2:A\times B\ra A\cup B$ given by $p_1(a,b)=a$, $p_2(a,b)=b$. A relation $R\subset A\times B$ is \ul{left-complete} if $p_1(R)=A$. $R$ is \ul{right-complete} if $p_2(R)=B$. $R$ is \ul{complete} if both left-complete and right-complete, i.e., $p_1(R)=A$ and $p_2(R)=B$.
\end{dfn}

\begin{prp}[\textcolor{OliveGreen}{Metrizabilty of $X(n)$}]\label{FSSmzty}
If a space $X$ is metrizable, then $X(n)$ is metrizable. Moreover, if the topology of $X$ is induced by a metric $d:X\times X\ra\Real$, then the topology of $X(n)$ is induced by the Hausdorff distance~
\bea
\textstyle d_H(\{x_i\}_i,\{y_i\}_i):=\max\big(\max_i\min_jd(x_i,y_j),\max_j\min_id(x_i,y_j)\big).\nn
\eea
\end{prp}
\begin{proof}
Let $d:X\times X\ra\Real$ be a metric that induces the topology of $X$. Then the topology of $X^n$ is induced by the metric $d_{\max}(x,y)=\max_{1\leq i\leq n}d(x_i,y_i)$, since $d_{\max}$-open balls satisfy $B^{d_{\max}}_\vep(x_1,...,x_n)=B^d_\vep(x_1)\times\cdots\times B^d_\vep(x_n)$,  for each $(x_1,...,x_n)\in X^n$,
and so are a common base for both the $d_{\max}$-topology and product topology of $X^n$. Consider the unordering map $q:X^n\ra X(n)$, $x=(x_1,...,x_n)\mapsto\{x_1,...,x_n\}$. Then for any $x,y\in X^n$,
\bea
d_H(q(x),q(y))=d_H\big(\{x_1,...,x_n\},\{y_1,...,y_n\}\big)\leq d_{\max}\big((x_1,...,x_n),(y_1,...,y_n)\big)=d_{\max}(x,y),\nn
\eea
and so $q:(X^n,d_{\max})\ra \big(X(n),d_H\big)$ is 1-Lipschitz (hence continuous). It follows that every $d_H$-open set in $X(n)$ is open in the quotient topology of $X(n)$, i.e., the quotient topology contains the $d_H$-topology.

It remains to show that the $d_H$-topology contains the quotient topology, i.e., every quotient-closed set is $d_H$-closed (which is equivalent to ``every quotient-open set is $d_H$-open''). Let $E\subset X(n)$ be quotient-closed. Let $x^k\in E$ be a $d_H$-convergent sequence, and let $x\in X(n)$ be its limit, i.e., $d_H(x^k,x)\ra 0$.

Our goal is to prove that $x\in E$. Take some vectors $z^k = (z_1^k, ... z_n^k)$ and $z = (z_1,..., z_n)$ in $X^n$ such that $q(z^k) = x^k$ and $q(z) = x$, i.e., $z^k\in q^{-1}(x^k)$ and $z\in q^{-1}(x)$. By the definition of $d_H$, for each $k$ we have a complete relation $R_k\subset\{1,\cdots,n\}\times\{1,\cdots,n\}$ such that $d(z_i^k, z_j)\le d_H\big(q(z^k),q(z)\big)=d_H(x^k,x)$ for each $(i,j)\in R_k$.

A key point here is that there are finitely many choices for $R_k$, and so by passing to a subsequence $\{k'\}$ we can assume $R_{k'} = R$ for all $k'$ (i.e., $R_{k'}$ is independent of $k'$).

Since $d_H(x^{k'}, x)\to 0$, we see that $d(z_i^{k'}, z_j) \to 0$ for each $(i,j)\in R$. Since $R$ is left-complete, we conclude that the sequence $z^{k'}$ converges in $X^n$. Let $w$ be its limit in $X^n$, i.e., $d_{\max}(z^{k'},w)\ra 0$. Then $w\in q^{-1}(E)$, since (i) $z^k\in q^{-1}(E)$, and (ii) $q^{-1}(E)$ is closed in $X^n$; recall that $q$ is quotient-continuous and $E$ is quotient-closed. That is, we have $z^{k'}\ra w\in q^{-1}(E)$.

Now, since $q :(X^n,d_{\max}) \to (X(n), d_H)$ is continuous, it follows that $x^{k'}=q(z^{k'})\sr{d_H}{\ral} q(w)\in q\big(q^{-1}(E)\big)=E$. Hence, $x=q(w)\in E$.
\end{proof}

\begin{rmk}[\textcolor{OliveGreen}{Alternative proof of Proposition \ref{FSSmzty} when $X$ is compact}]
(i) Recall that a continuous bijection $f:K\ra H$ from a compact space $K$ to a Hausdorff space $H$ must be a homeomorphism. (ii) By the continuity of $q:X^n\ra \big(X(n),d_H\big)$ and the universal property of quotient maps, there exists a unique continuous bijection $\wt{q}:X(n)\ra\big(X(n),d_H\big)$ satisfying $\wt{q}\circ q=q$ (i.e., $\wt{q}$ is the identity map), as in the following diagram.
\bc\bt
X^n\ar[d,"q"]\ar[rr,"q"]&& \big(X(n),d_H\big)\\
X(n)\ar[urr,dashed,"\exists!~\wt{q}=id"',"\cong"]
\et\ec
Thus, if $X$ is a compact metric space (so that $X(n)$ is compact), then $\wt{q}=id_{X(n)}$ is a homeomorphism by (i) and (ii) above, and so the quotient topology and the $d_H$-topology on $X(n)$ are equivalent.
\end{rmk}

\section{Finite subset spaces of metric spaces}\label{FSSfs}
Henceforth (due to Proposition \ref{FSSmzty} and Lemma \ref{MtzblUnion}) if $X$ is a metric space then, unless stated otherwise, we will assume {\small $X(n):=\big\{A\subset X:1\leq|A|\leq n\big\}=\big\{\{x_1,...,x_n\}:(x_1,...,x_n)\in X^n\big\}$} and {\small $FS(X):=\big\{A\subset X:1\leq|A|< \infty\big\}=\bigcup_{n\geq 1} X(n)$} are metric spaces with respect to the function
\begin{equation*}
\textstyle d_H\big(\{x_1,...,x_n\},\{x_1',...,x_n'\}\big):=\max\left\{\max_i\min_jd(x_i,x_j'),\max_i\min_jd(x_j,x_i')\right\}.
\end{equation*}

\begin{figure}[H]
\centering
\scalebox{1.5}{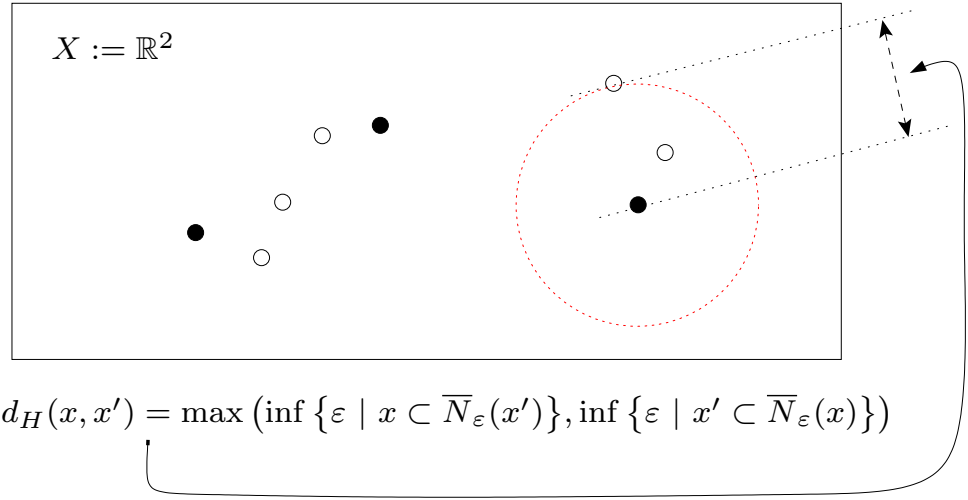} 
  \caption{~A visualization of Hausdorff distance in $\Real^2(n)$, for $n\geq 5$.}\label{dgdH}
\end{figure}

\begin{dfn}[\textcolor{blue}{\index{Minimum separation $\delta$}{Minimum separation} on finite subsets}]\label{MinSep}
Let $X$ be a metric space. The minimum separation is the function ~{\small$\delta_n:X(n)\ra[0,+\infty),~x=\{x_1,...,x_n\}\mapsto\min\limits_{i\neq j}d(x_i,x_j)$}.
\end{dfn}
Note that for any $x\in X(n)$, if $|x|<n$ then $\delta_n(x)=0$, i.e., $\delta_n|_{X(n-1)}=0$. Also note that $\delta_n:X(n)\ra \Real$ is related to the following \ul{total minimum separation}:
\bea
\label{TotMinSep}\textstyle \delta:FS(X)\ra\Real,~A\mapsto\min_{\substack{a,a'\in A\\a\neq a'}}d(a,a'),
\eea
which satisfies ~$\delta|_{X(n)\backslash X(n-1)}=\delta_n|_{X(n)\backslash X(n-1)}$, for all $n\geq 2$, i.e.,
\bea
\delta_n(x)=\left\{
              \begin{array}{ll}
                \delta(x), & x\in X(n)\backslash X(n-1) \\
                0, & x\in X(n-1)
              \end{array}
            \right\}\neq \delta(x),~~~~\txt{for all}~~x\in X(n).\nn
\eea

Nevertheless, unless it is specified otherwise, $\delta$ will always mean $\delta_n:X(n)\ra\Real$ for some $n$ under consideration.

\begin{lmm}\label{HausDistBound}
Let $X$ be a metric space. For all $x,y\in X(n)$, if $\delta_n(x)>2d_H(x,y)$ or $\delta_n(y)>2d_H(x,y)$, then there exists an enumeration $x_i,y_i$ of elements of $x,y$ such that
\begin{equation}
\label{KOVeq8}d(x_i,y_i)\leq d_H(x,y),~~~~~~\txt{for all}~~i=1,...,n.
\end{equation}
\end{lmm}
\begin{proof}
By symmetry it suffices to assume $\delta_n(x)>2d_H(x,y)$. Let $\rho:=d_H(x,y)$. Then by Remark \ref{HausDistRmk}, $x\subset \ol{y}_\rho=\bigcup_{i=1}^n\ol{N}_\rho(y_i)$ and $y\subset \ol{x}_\rho=\bigcup_{i=1}^n\ol{N}_\rho(x_i)$. Moreover, the balls $\ol{N}_\rho(x_i)$ are disjoint because $\delta_n(x)>2\rho$. Thus, each ball $\ol{N}_\rho(x_i)$ contains exactly one point of $y$, i.e., for each $x_i\in x$, there exists a \emph{unique} $y_j\in y$ such that $d(x_i,y_j)\leq\rho$. This means we can \emph{relabel} the points $x_i,y_i$ so that (\ref{KOVeq8}) holds.
\end{proof}

\begin{crl}\label{HausDistBound1}
Let $X$ be a metric space, and let $\vep>0$ be given. Then
\bea
\textstyle N^{d_H}_\vep(x):=\{y\in X(n)~|~d_H(x,y)<\vep\}=\Big(\bigcup_{a\in x}N^d_\vep(a)\Big)(n),~~~~\txt{for any}~~x\in X(n).\nn
\eea
Moreover, if $x=\{x_1,...,x_n\}\in X(n)$ has $\delta_n(x)\geq 2\vep$, then
\bea
N^{d_H}_\vep(x)=\big\{\{y_1,...,y_n\}\in X(n)~|~y_i\in N^d_\vep(x_i),~i=1,...,n\big\}=q\big(N^d_\vep(x_1)\times\cdots\times N^d_\vep(x_n)\big),\nn
\eea
where $q:X^n\ra X(n)$ is the unordering map.
\end{crl}

\begin{lmm}\label{DeltaCont}
Let $X$ be a metric space. The minimum separation $\delta_n:X(n)\ra\Real$ in Definition \ref{MinSep} is $2$-Lipschitz.
\end{lmm}
\begin{proof}
If $\delta_n(x)\leq 2d_H(x,y)$ and $\delta_n(y)\leq 2d_H(x,y)$, it is clear that
\begin{equation*}
|\delta_n(x)-\delta_n(y)|\leq\max\big(\delta_n(x),\delta_n(y)\big)\leq 2d_H(x,y).
\end{equation*}
So, assume $\delta_n(x)>2d_H(x,y)$ or $\delta_n(y)>2d_H(x,y)$. Then using Lemma \ref{HausDistBound},
\begin{equation*}
\textstyle|\delta_n(x)-\delta_n(y)|\leq\max\limits_{i\neq j}\Big|d(x_i,x_j)-d(y_i,y_j)\Big|\leq \max\limits_{i\neq j}\Big(d(x_i,y_i)+d(x_j,y_j)\Big)\sr{(\ref{KOVeq8})}{\leq} 2d_H(x,y).\qedhere
\end{equation*}
\end{proof}

\begin{lmm}\label{HdPPT2}
If $X$ is a metric space, the Hausdorff distance on $X(2)$ is given by
{\small\bea
d_H\Big(\{x_1,x_2\},\{y_1,y_2\}\Big)=\min\Big\{\max\big\{d(x_1,y_1),d(x_2,y_2)\big\},\max\big\{d(x_1,y_2),d(x_2,y_1)\big\}\Big\}.\nn
\eea}
\end{lmm}
\begin{proof}
By Lemma \ref{HausDistRmk}(ii),
{\small\begin{align}
&d_H\Big(\{x_1,x_2\},\{y_1,y_2\}\Big)=\inf\Big\{\vep:\{x_1,x_2\}\subset\ol{\{y_1,y_2\}}_\vep,~\{y_1,y_2\}\subset\ol{\{x_1,x_2\}}_\vep\Big\}\nn\\
&~~=\txt{\footnotesize $\inf\Big\{\vep:\dist(x_1,\{y_1,y_2\})\leq\vep,\dist(x_2,\{y_1,y_2\})\leq\vep,\dist(y_1,\{x_1,x_2\})\leq\vep,\dist(y_2,\{x_1,x_2\})\leq\vep\Big\}$}\nn\\
&~~=\inf\Big\{\vep:~\min\big\{d(x_1,y_1),d(x_1,y_2)\big\}\leq\vep,~\min\big\{d(x_2,y_1),d(x_2,y_2)\big\}\leq\vep,\nn\\
&~~~~~~~~~~~\min\big\{d(y_1,x_1),d(y_1,x_2)\big\}\leq\vep,~\min\big\{d(y_2,x_1),d(y_2,x_2)\big\}\leq\vep\Big\}\nn\\
&~~=\inf\Big\{\vep:~\min\big\{\max\{d(x_1,y_1),d(x_2,y_2)\},\!\!~\max\{d(x_1,y_2),d(y_1,x_2)\}\big\}\leq\vep\Big\}\nn\\
&~~=\min\Big\{\max\big\{d(x_1,y_1),d(x_2,y_2)\big\},\!\!~\max\big\{d(x_1,y_2),d(x_2,y_1)\big\}\Big\}.\nn\qedhere
\end{align}}
\end{proof}

\section{Some inherited properties of finite subset spaces}\label{FSSsi}
We have already seen that if a space $X$ is metrizable, then so is $X(n)$. Similarly, $X(n)$ inherits all properties of $X$ that are (i) passed onto $X^n$ and then (ii) preserved by continuity of the quotient map $q:X^n\ra X(n)$. Such properties include connectedness, path-connectedness, separability, compactness (where compactness of $X^n$ follows from Tychonoff's theorem), and more (see \cite[pp 877-878]{borsuk-ulam1}). We also have the following result on completeness.

\begin{prp}[\textcolor{OliveGreen}{\cite[Proposition 7.3.7]{BBI}}]\label{FSScmpl}
If $X$ is a complete metric space, so is $\big(\B_c^\ast(X),d_H\big)$. (Hence $X(n)$ is complete as a closed subspace of $\B_c^\ast(X)$.)
\end{prp}
\begin{proof}
Let $\{A_n\}$ be a Cauchy sequence in $\B_c^\ast(X)$, i.e., $d_H(A_n,A_m)\ra 0$. With $N(x)$ denoting a neighborhood of $x$, let $A$ $:=$ $\{$ $x\in X$: every $N(x)\cap A_n\neq\emptyset$ for infinitely many $n$$\}$. We will prove that $d_H(A,A_n)\ra 0$. Fix $\vep>0$ and let $n_0=n_0(\vep)$ be such that $d_H(A_n,A_m) < \vep$ for all $m,n\geq n_0$ (which is possible by Cauchyness of $\{A_n\}$). We will show that
\bea
\textstyle d_H(A,A_n):=\max\Big(\sup\limits_{x\in A}~\dist(x,A_n)~,~\sup\limits_{y\in A_n}~\dist(A,y)\Big) < 2\vep,~~\txt{for all}~~n\geq n_0,\nn
\eea
or equivalently, if we fix $n\geq n_0$, then (i) $\dist(x,A_n)< 2\vep$ for each $x\in A$, and (ii) $\dist(A,y)< 2\vep$ for each $y\in A_n$. So, let us fix $n\geq n_0$.

(i) \ul{$\dist(x, A_n) < 2\vep$ for each $x \in A$}: Fix $x\in A$ (along with $n\geq n_0$). Then by the definition of $A$, there exists $m\geq n_0$ such that $B_\vep(x)\cap A_m\neq\emptyset$, $\Ra$ there exists $y_m \in A_m$ such that $d(x,y_m) < \vep$. Therefore,
\bea
\textstyle\dist(x,A_n)\leq d(x,y_m) + \dist(y_m,A_n)\leq d(x,y) + d_H(A_m,A_n) < 2\vep.\nn
\eea

(ii) \ul{$\dist(A,y) < 2\vep$ for each $y\in A_n$}: Fix $y\in A_n$ (along with $n\geq n_0$). Let $n_1 := n$ and, for every integer
$k > 1$, choose $n_k=n_k(\vep)\sr{\txt{say}}{=} n_0(\vep/2^k)$ such that $n_k < n_{k+1}$ and $d_H(A_p, A_q) <\vep/2^k$ for
all $p,q\geq n_k$. The subsequence $n_k$ with this property is possible since $\{A_n\}$ is Cauchy wrt $d_H$. (In particular, with $p=n_k$, $q=n_{k+1}$, we have $d_H(A_{n_k},A_{n_{k+1}}) <\vep/2^k$ for all $k>1$). Define a sequence $\{y_k\}$, $y_k\in A_{n_k}$, as
follows. Let $y_1:= y$, and for all $k\geq 1$ pick $y_{k+1}\in A_{n_{k+1}}$ such that $d(y_k,y_{k+1})<\vep/2^k$ (which is possible since $d_H(A_{n_k},A_{n_{k+1}})<\vep/2^k$). Then, because {\small$\sum_{k=1}^\infty d(y_k,y_{k+1})<2\vep<\infty$}, the sequence $\{y_k\}$ is Cauchy and hence converges to a point $x\in X$. Since by construction, $x$ lies in infinitely many $A_m$, we know $x\in A$. Now, {\small$d(x,y)=\lim_k d(y_k,y)=\lim_k d(y_k,y_1)\leq\sum_kd(y_k,y_{k+1})<2\vep$}, and so {\small$\dist(A,y)\leq d(x,y)< 2\vep$}.
\end{proof}

\begin{crl}[\textcolor{OliveGreen}{\cite[Theorem 7.3.8]{BBI}}]\label{FSScmpc}
If $X$ is a compact metric space, so is $\big(\B_c^\ast(X),d_H)$.
\end{crl}
\begin{proof}
By Proposition \ref{FSScmpl}, $\B_c^\ast(X)$ is complete. Therefore it suffices (by Theorem  \ref{MetricCompact1}) to prove that $\B_c^\ast(X)$ is totally bounded. Fix $\vep>0$. Let $F\subset X$ be a finite set such that
\bea
X=F_\vep:=\{x\in X:\dist(x,F)<\vep\}.\nn
\eea
With the finite set $\F:=\P(F)\subset\B_c^\ast(X)$ of all subsets of $F$, we will prove that
\bea
\B_c^\ast(X)=\F_\vep:=\{A\in \B_c^\ast(X):\dist_H(A,\F)<\vep\}.\nn
\eea
Let $A\in\B_c^\ast(X)$. Consider the set {\small $F_A:=F\cap A_\vep=\{x\in F:\dist(x,A)<\vep\}\subset F$}.

Since $X=F_\vep$, for every $a\in A\subset F_\vep$ there exists an $x_a \in F$ such that
$d(x_a,a)<\vep$. Since $\dist(x_a,A)\leq d(x_a,a)<\vep$, the point $x_a\in F_A$. Therefore
$\dist(a,F_A)\leq d(a,x_a)<\vep$ for all $a\in A$. Since $\dist(x,A)<\vep$ for any $x \in F_A$ (by the
definition of $F_A$), it follows that
\bea
\textstyle d_H(A,F_A):=\max\Big(\sup\limits_{a\in A}~\dist(a,F_A)~,~\sup\limits_{x\in F_A}~\dist(A,x)\Big) < \vep.\nn
\eea
 This shows ~$\dist_H(A,\F)\leq d_H(A,F_A)<\vep$.
\end{proof}

{\flushleft \textbf{Snowflake properties of finite subset spaces}}
{\flushleft A} review of  section  \ref{PrlMET3} should be helpful with the following discussion. Fix $\ld\geq 1$ and $1<p\leq\infty$. Recall that a map of metric spaces $f:(X,d)\ra (Z,d_1)$ is $\ld$-biLipschitz if $d(x,y)/\ld\leq d_1(f(x),f(y))\leq \ld d(x,y)$ for all $x,y\in X$.

\begin{dfn}[\textcolor{blue}{Bi-Lipschitz equivalent metrics}]
Let $X$ be a set. Two metrics $d,d_1:X\times X\ra\Real$ are $\ld$-biLipschitz equivalent if the identity map $id_X:(X,d)\ra (X,d_1)$ is $\ld$-biLipschitz, i.e., $d(x,y)/\ld\leq d_1(x,y)\leq \ld d(x,y)$ for all $x,y\in X$.
\end{dfn}

\begin{dfn}[\textcolor{blue}{$L^p$-metric, $L^p$-metric space}]
Let $(X,d)$ be a metric space. The metric $d$ is an $L^p$-metric (making $(X,d)$ an $L^p$-metric space) if for all $x,y,z\in X$,
\bea
\label{SnowFlakeEq}d(x,y)\leq\left\{
            \begin{array}{ll}
              \big(d(x,z)^p+d(z,y)^p\big)^{1\over p}, &\txt{if}~~1<p<\infty \\
              \max\big(d(x,z),d(z,y)\big), &\txt{if}~~ p=\infty
            \end{array}
          \right\}.
\eea
\end{dfn}

For example, if $(X,d)$ is a metric space, then for any $1<p<\infty$, the metric $d_p(x,y):=d(x,y)^{1\over p}$ is an $L^p$-metric on $X$. Thus, $(X,d_p)$ is a snowflake (as defined below).

\begin{dfn}[\textcolor{blue}{\index{Snowflake metric space}{Snowflake metric space}}]
Fix $1<p\leq\infty$. A metric space $(X,d)$ is a $p$-snowflake if the metric $d$ is biLipschitz equivalent to an $L^p$-metric $d_1:X\times X\ra\Real$.
\end{dfn}

\begin{lmm}[\textcolor{OliveGreen}{\cite[Proposition 2.3, p.318]{tyson-wu2005}}]
If $1<p<\infty$, then a metric space $(X,d)$ is a $p$-snowflake $\iff$ there exists a constant $c>0$ such that for any finite set of points $x_0,x_1,..,x_N\in X$, we have ~$\sum_{i=1}^Nd(x_{i-1},x_i)^p\geq cd(x_0,x_N)^p$.
\end{lmm}
\begin{proof}
($\Ra$): Assume $(X,d)$ is a $p$-snowflake. Let $d_1$ be an $L^p$-metric on $X$ such that $d(x,y)/\ld\leq d_1(x,y)\leq\ld d(x,y)$. Then for any finite set of points $x_0,x_1,...,x_N\in X$,
\bea
\textstyle d(x_0,x_N)^p/\ld^p\leq d_1(x_0,x_N)^p\sr{(\ref{SnowFlakeEq})}{\leq} \sum_{i=1}^Nd_1(x_{i-1},x_i)^p\leq \ld^p\sum_{i=1}^Nd(x_{i-1},x_i)^p.\nn
\eea
($\La$): Conversely, assume there exists $c$ such that $\sum_{i=1}^Nd(x_{i-1},x_i)^p\geq cd(x_0,x_N)^p$ for every finite set of points $x_0,x_1,...,x_N\in X$. Let
{\small\begin{align}
\textstyle d_1(x,y):=\inf\big\{l(c)~|~c=\{x=x_0,x_1,...,x_N=y\}\big\},~~\txt{where}~~l(c):=\Big(\sum\limits_{i=1}^Nd(x_{i-1},x_i)^p\Big)^{1\over p}.\nn
\end{align}}
Then given $\vep>0$ and $z\in X$, we can pick two finite chains $c_{xz}$ (a finite chain from $x$ to $z$) and $c_{zy}$ (a finite chain from $z$ to $y$) such that with $c_{xy}:=c_{xz}\cdot c_{zy}$,
\bea
d_1(x,y)^p\leq l(c_{xy})^p=l(c_{xz})^p+l(c_{zy})^p\leq d'(x,z)^p+d_1(z,y)^p+2\vep.\nn
\eea
This shows $d_1$ is an $L^p$-metric on $X$ satisfying $c^{1/p} d\leq d_1\leq d$, and so $(X,d)$ is a $p$-snowflake.
\end{proof}

\begin{lmm}\label{SnwfNoRectCrv}
If $(X,d)$ is a $p$-snowflake, $1<p<\infty$, then every nonconstant path $\gamma:[0,1]\ra X$ is non-rectifiable (i.e., has infinite length).
\end{lmm}
\begin{proof}
Let $\gamma:[0,1]\ra X$ be a nonconstant path, i.e., there exist $t,t'\in[0,1]$, $t<t'$, such that $\gamma(t)\neq\gamma(t')$. Suppose $l(\gamma)<\infty$. Let $t=t_0,t_1,...,t_N=t'\in[t,t']$ be a finite sequence of points such that $l\left(\gamma|_{[t,t']}\right)/N\geq d(\gamma(t_{i-1}),\gamma(t_i))$ for all $i=1,...,N$. Then
\bea
\textstyle l\left(\gamma|_{[t,t']}\right)^p/N^{p-1}\geq\sum\limits_{i=1}^Nd(\gamma(t_{i-1}),\gamma(t_i))^p\geq cd(\gamma(t),\gamma(t'))^p.\nn
\eea
Taking the limit $N\ra\infty$, we get a contradiction.
\end{proof}

\begin{crl}\label{SnwfNoLipCrv}
If $X$ is a snowflake, then every Lipschitz path $\gamma:[0,1]\ra X$ is constant.
\end{crl}

\begin{crl}\label{SnwfNoLipCon}
If $Z$ is a geodesic space and $X$ is a snowflake, then every Lipschitz map $f:Z\ra X$ is constant.
\end{crl}

\begin{dfn}[\textcolor{blue}{\index{Snowflake curve}{Snowflake curve}}]
Let $X$ be a metric space. A curve $\gamma:[0,1]\ra X$ is a snowflake curve if the subspace $\gamma\big([0,1]\big)\subset X$ is a snowflake metric space.
\end{dfn}

\begin{example*}[\textcolor{blue}{von Koch snowflake curve, $K$: See \cite{tukia81,koskela94,mckemie1987}}]
In $\Real^2$, let {\small $K_0:=I=[(0,0),(1,0)]$. Next, let ~$K_1:=\left[\left(0,0\right),\left({1\over 3},0\right)\right]\cup\big[\left({1\over 3},0\right),\big({1\over 2},{1\over 3}{\sqrt{3}\over 2}\big)\big]\cup\big[\big({1\over 2},{1\over 3}{\sqrt{3}\over 2}\big),\left({2\over 3},0\right)\big]\cup\left[\left({2\over 3},0\right),\left(1,0\right)\right]$},~
a union of 4 line segments each of length $1/3$ formed from the line segment $K_0$. Similarly, we can apply the process to each of the 4 line segments of $K_1$ to obtain $K_2$ as a union of $4\times 4=4^2$ line segment each of length ${1\over 3}\times{1\over 3}$. Continuing this way, at the $j$th step, we get a set $K_j\subset\Real^2$ consisting of $4^j$ line segments each of length ${1\over 3^j}$. Note that the length of $K_j$ is $4^j\times{1\over 3^j}=\left(4\over 3\right)^j$.

In the limit $j\ra\infty$, we obtain a curve $K=\gamma(I)$, $\gamma:[0,1]\ra\Real^2$, of infinite length between $(0,0)$ and $(1,0)$. Moreover, between any two distinct points $\gamma(t),\gamma(t')\in K$, the segment $\gamma|_{[t,t']}$ also has infinite length.
\end{example*}

\begin{prp}[\textcolor{OliveGreen}{Snowflake subset spaces}]\label{SnwfSubSps}
Let $X$ be a space and $1<p<\infty$.
\bit[leftmargin=0.9cm]
\item[(i)] If two metrics $d_1,d_2:X\times X\ra\Real$ are $\ld$-biLipschitz equivalent, then so are the corresponding Hausdorff metrics $d_{1H},d_{2H}:X(n)\times X(n)\ra\Real$.
\item[(ii)] If $X$ is an $L^p$-metric space, then so is $X(n)$, for all $n\geq 1$.
\item[(iii)] If $X$ is a $p$-snowflake, then so is $X(n)$, for all $n\geq 1$.
\eit
\end{prp}
\begin{proof}
(i) This follows immediatly from the definitions. (ii) If $x,y\in X(n)$, then
$d_H(x,y)^p=\left(\max\left\{\max_i\min_jd(x_i,y_j),\max_j\min_id(x_i,y_j)\right\}\right)^p$ = $\max\big\{\max_i\min_jd(x_i,y_j)^p$, $\max_j\min_id(x_i,y_j)^p\big\}$ ,
where for any $z_k\in z\in X(n)$,
\begin{align}
&\textstyle \max_i\min_jd(x_i,y_j)^p\leq \max_i\min_j\left[d(x_i,z_k)^p+d(z_k,y_j)^p\right]= \max_id(x_i,z_k)^p+\min_jd(z_k,y_j)^p\nn\\
&\textstyle~~~~\leq\max_id(x_i,z_k)^p+\max_i\min_jd(z_i,y_j)^p,\nn\\
&\textstyle~~\Ra~~\max_i\min_jd(x_i,y_j)^p\leq \max_i\min_jd(x_i,z_j)^p+\max_i\min_jd(z_i,y_j)^p,\nn
\end{align}
from which it follows that ~$d_H(x,y)^p\leq d_H(x,z)^p+d_H(z,y)^p$.

(iii) This follows immediately from (i) and (ii).
\end{proof}

{\flushleft A further} discussion of snowflake concepts can be found in \cite{tyson-wu2005}.

\section{Gaps in subsets of a metric space: Gap reducing property}\label{FSSgi}
\begin{dfn}[\textcolor{blue}{\index{Gap of a set}{Gap of a set}}]\label{SubsetGap}
Let $X$ be a metric space and $A\subset X$. The gap $\rho(A)$ of $A$ is the largest distance between any two nonempty complementary subsets of $A$, i.e.,
\bea
\rho(A):=\sup_{\emptyset\neq A'\subsetneq A}\dist(A',A-A')=\sup_{\substack{A'\sqcup A''=A\\ A',A''\neq\emptyset}}\dist(A',A'').\nn
\eea
\end{dfn}

Note that if $X$ is a metric space and $A\subset X$, the minimum distance $\delta(A)$ between points of $A$, which is an extension of (\ref{TotMinSep}), satisfies
\bea
\delta(A):=\inf_{\substack{a,a'\in A\\ a\neq a'}}d(a,a')=\inf_{\substack{A'\sqcup A''=A\\ A',A''\neq\emptyset}}\dist(A',A'')=\inf_{\emptyset\neq A'\subsetneq A}\dist(A',A-A').\nn
\eea

\begin{lmm}[\textcolor{OliveGreen}{Continuity of $\rho$}]
Let $X$ be a metric space. The gap function $\rho:FS(X)\ra\Real$ is $2$-Lipschitz with respect to Hausdorff distance, i.e., $|\rho(A)-\rho(B)|\leq 2d_H(A,B)$.
\end{lmm}
\begin{proof}
By symmetry, it is enough to show that $\rho(B)\geq\rho(A)-2d_H(A,B)$. If $\rho(A)\leq 2d_H(A,B)$, then the result holds trivially. So, further assume $\rho(A)>2d_H(A,B)$.

Let $\rho(A)=\dist(A',A'')$. Observe that for any $b\in B$, we have
\bea
\dist(A',b)=\dist(A,b)\leq d_H(A,B)~~~\txt{or}~~~\dist(A'',b)=\dist(A,b)\leq d_H(A,B),\nn
\eea
since $\min\big(\dist(A',b),\dist(A'',b)\big)=\dist(A,b)\leq d_H(A,B)$. On the other hand, if there exists $b\in B$ such that $\dist(A',b)\leq d_H(A,B)$ and $\dist(A'',b)\leq d_H(A,B)$, then
\bea
\rho(A)=\dist(A',A'')\leq\dist(A',b)+\dist(b,A'')\leq 2d_H(A,B),~~~~\txt{(a contradiction)}.\nn
\eea
It follows that for any $b\in B$ exactly one of the following holds:
\bea
\txt{Either}~~~\dist(A',b)=\dist(A,b)\leq d_H(A,B)~~~\txt{or}~~~\dist(A'',b)=\dist(A,b)\leq d_H(A,B).\nn
\eea
Thus, with $B':=\{b\in B:\dist(A',b)\leq d_H(A,B)\}$, $B'':=\{b\in B:\dist(A'',b)\leq d_H(A,B)\}$, we get the disjoint union $B=B'\sqcup B''$. Note that $B',B''$ are nonempty because, by the definition of the Hausdorff distance $d_H(A,B)$, for any $a'\in A'$, $a''\in A''$ there exist $b',b''\in B$ such that $d(a',b')\leq d_H(A,B)$, $d(a'',b'')\leq d_H(A,B)$. Now, for any $b'\in B',b''\in B''$,
\begin{align}
&\rho(A)=\dist(A',A'')\leq \dist(A',b')+d(b',b'')+\dist(b'',A'')\leq \dist(b',b'')+2d_H(A,B),\nn\\
&~~\Ra~~\rho(B)\geq\dist(B',B'')\geq \rho(A)-2d_H(A,B).\nn\hfill\qedhere
\end{align}
\end{proof}

\begin{dfn}[\textcolor{blue}{$d_H$-gap}]
Let $X$ be a metric space and $A\subset X$. The $d_H$-gap of $A$ is
\bea
\ld(A):=\inf_{\emptyset\neq A'\subsetneq A}d_H(A',A-A')=\inf_{\substack{\emptyset\neq A',A''\subset A\\ A'\sqcup A''=A}}d_H(A',A'').\nn
\eea
\end{dfn}

Note that the $d_H$-gap function $\ld:FS(X)\ra\Real$ is not continuous. To see this, let $X=(X,\|\|)$ be a normed space and choose sets $A=\{a,b,c,d\}$, $B=\{a,b,c\}\in X(4)$ such that $\vep:=\|a-b\|=\|c-d\|$ is small and $L:=\dist(\{a,b\},\{c,d\})$ is large ($\gg\vep$). Then $d_H(A,B)=\|c-d\|=\vep$, but
\bea
\ld(A)=\vep,~~\ld(B)\geq\dist(\{a,b\},\{c,d\})=L,~~\Ra~~|\ld(A)-\ld(B)|\gg \vep.\nn
\eea

\begin{dfn}[\textcolor{blue}{\index{Gap reducing property (GRP)}{Gap reducing property (GRP)}, Gap reducing space}]\label{GapRedProp}
A metric space $X$ is gap reducing, or has the gap reducing property (GRP), if for every finite set $A\subset X$, the pair of subsets $A',A''\subset A$ satisfying $\rho(A)=\dist(A',A'')$ can be chosen such that $\rho(A'),\rho(A'')\leq\rho(A)$.
\end{dfn}

Before proving, in Proposition \ref{GRPprp}, that every metric space has the GRP, we will first obtain a few motivational results.

\begin{dfn}[\textcolor{blue}{Unique-gap sets}]\label{UniqGapSet}
Let $X$ be a metric space and $A\subset X$ a finite set. We say $A$ has a unique gap if the complementary pair of subsets $A',A''\subset A$ satisfying $\rho(A)=\dist(A',A'')$ is unique.
\end{dfn}

\begin{lmm}[\textcolor{OliveGreen}{Avoiding a large subgap}]\label{UnifGapProp}
Let $X$ be a metric space, $A\subset X$ a finite set, and $\rho(A)=\dist(A',A'')$. If $\rho(A')>\rho(A)$, then the following are true.
\bit[leftmargin=0.9cm]
\item[(i)] For any nonempty complementary subsets $C',D'\subset A'$, if $\rho(A')=\dist(C',D')$, then
\bea
\dist(C',A'')=\dist(D',A'')=\rho(A).\nn
\eea
\item[(ii)] The pair $A',A''\subset A$ satisfying $\rho(A)=\dist(A',A'')$ is not unique. (That is, there exists a different pair of complementary subsets $A'_{\txt{new}},A''_{\txt{new}}\subset A$ such that $\rho(A)=\dist(A'_{\txt{new}},A''_{\txt{new}})$.)
\eit
\end{lmm}
\begin{proof}
(i) Suppose on the contrary that $\rho(A)<\dist(D',A'')$ or $\rho(A)<\dist(C',A'')$. Then we can consider the two cases separately as follows.
\bit[leftmargin=0.9cm]
\item[(a)] {\small $\rho(A)<\dist(D',A'')$} and {\small $\rho(A)=\dist(A',A'')=\dist(C',A'')<\rho(A')=\dist(C',D')$}: In this case, let $A'_{\txt{new}}:=D'$ and $A''_{\txt{new}}:=C'\cup A''$. Then
{\small\bea
\rho(A)=\dist(C',A'')<\min\big(\dist(C',D'),\dist(D',A'')\big)\leq\dist(A'_{\txt{new}},A''_{\txt{new}}),\nn
\eea} which is a contradiction since $A=A'_{\txt{new}}\cup A''_{\txt{new}}$.
\item[(b)] {\small $\rho(A)<\dist(C',A'')$} and {\small $\rho(A)=\dist(A',A'')=\dist(D',A'')<\rho(A')=\dist(C',D')$}: We obtain a contradiction in the same way as in (a).
\eit

{\flushleft (ii)} From part (i), $\rho(A)=\dist(A',A'')=\dist(C',A'')=\dist(D',A'')<\rho(A')=\dist(C',D')$. Thus, with $A'_{\txt{new}}:=D'$ and $A''_{\txt{new}}:=C'\cup A''$, we have
\begin{align}
\rho(A)&=\dist(C',A'')\leq\min\big((\dist(C',D'),\dist(D',A'')\big)\leq\dist(A'_{\txt{new}},A''_{\txt{new}})=\rho(A),\nn\\
\Ra&~~\rho(A)=\dist(A',A'')=\dist(A'_{\txt{new}},A''_{\txt{new}}).\nn\hfill\qedhere
\end{align}
\end{proof}

\begin{crl}[\textcolor{OliveGreen}{A unique-gap set is gap reducing}]\label{UniqGapRed}
Let $X$ be a metric space and $A\subset X$ a finite set. If the complementary pair of subsets $A',A''\subset A$ satisfying $\rho(A)=\dist(A',A'')$ is unique (i.e., $A$ has a unique gap), then $\rho(A'),\rho(A'')\leq \rho(A)$.
\end{crl}

\begin{crl}[\textcolor{OliveGreen}{GRP from density of unique-gap sets}]\label{UniqGapRed2}
Let $X$ be a metric space. If the set of unique-gap elements of $X(n)$ is dense, then $X$ has the gap reducing property.
\end{crl}
\begin{proof}
Let $A\in X(n)$, and $A_k\in X(n)$ a sequence of unique-gap sets such that $A_k\sr{d_H}{\ral} A$, i.e., $d_H(A_k,A)\ra 0$. Let $\rho(A_k)=\dist(A_k',A_k'')$, where $\rho(A_k'),\rho(A_k'')\leq\rho(A_k)$. Let $\rho(A)=\dist(A',A'')$, and recall that $|\rho(A_k)-\rho(A)|\leq 2d_H(A_k,A)$. Given any $\vep<\rho(A)/4$, choose $k_\vep$ such that $d_H(A_{k_\vep},A)<\vep$, which implies
\bea
&&|\rho(A_{k_\vep})-\rho(A)|\leq 2d_H(A_{k_\vep},A)<2\vep,\nn\\
&&~~\Ra~~\rho(A_{k_\vep})>\rho(A)-2\vep>4\vep-2\vep=2\vep.\nn
\eea
Define {\small $A_\vep':=A\cap N_\vep(A_{k_\vep}')$}, {\small $A_\vep'':=A\cap N_\vep(A_{k_\vep}'')$}, where {\small $N_r(C):=\{x\in X:\dist(x,C)<r\}$}. Then
\bea
A=A\cap N_\vep(A_{k_\vep})=A\cap\left[N_\vep(A_{k_\vep}'\sqcup A_{k_\vep}'')\right]=A\cap\left[ N_\vep(A_{k_\vep}')\sqcup N_\vep(A_{k_\vep}'')\right]=A_\vep'\sqcup A_\vep'',\nn
\eea
where by direct computation we get the inequalities
\begin{align}
&d_H(A_\vep',A_{k_\vep}')\leq \max\left(d_H(A_\vep',A_{k_\vep}'),d_H(A_\vep'',A_{k_\vep}'')\right)=d_H(A_\vep,A_{k_\vep})<\vep,\nn\\
&d_H(A_\vep'',A_{k_\vep}'')\leq \max\left(d_H(A_\vep'',A_{k_\vep}''),d_H(A_\vep'',A_{k_\vep}'')\right)=d_H(A_\vep,A_{k_\vep})<\vep.\nn
\end{align}
Therefore,
{\small\begin{align}
& \left|\dist(A_\vep',A_\vep'')-\rho(A_{k_\vep})\right|=\left|\dist(A_\vep',A_\vep'')-\dist(A_{k_\vep}',A_{k_\vep}'')\right|\leq d_H(A_\vep',A_{k_\vep}')+d_H(A_\vep'',A_{k_\vep}'')<2\vep,\nn\\
&~~\Ra~~|\rho(A)-\dist(A_\vep',A_\vep'')|\leq \left|\dist(A_\vep',A_\vep'')-\rho(A_{k_\vep})\right|+|\rho(A)-\rho(A_{k_\vep})|<4\vep,\nn
\end{align}}
and
{\small\begin{align}
&\rho(A_\vep')\leq|\rho(A_\vep')-\rho(A_{k_\vep}')|+\rho(A_{k_\vep}')\leq d_H(A_\vep',A_{k_\vep}')+\rho(A_{k_\vep}')< 2\vep+\rho(A_{k_\vep})<4\vep+\rho(A),\nn\\
&\rho(A_\vep'')\leq|\rho(A_\vep'')-\rho(A_{k_\vep}'')|+\rho(A_{k_\vep}'')\leq d_H(A_\vep'',A_{k_\vep}'')+\rho(A_{k_\vep}'')< 2\vep+\rho(A_{k_\vep})<4\vep+\rho(A).\nn\hfill\qedhere
\end{align}}
\end{proof}

\begin{crl}
Every normed space has the GRP.
\end{crl}

\begin{dfn}[\textcolor{blue}{\index{Graph}{Graph}, Path, Loop, Connected Graph, Tree, Subgraph, Subtree, \index{Spanning tree}{Spanning tree}}]
Let $X$ be a metric space and $A\subset X$ a discrete set (i.e., discrete subspace). A \ul{graph} $G$ over $A$ consists of \ul{vertices} $V=A$ (i.e., every element of $A$ is a vertex of $G$), \ul{edges} $E\subset V\times V$ (i.e., a relation between elements of $A$) with an edge $e_{uv}$ between $u,v\in V$ iff $(u,v)\in E$. Briefly, we write $G=(V,E)=\big(V_G,E_G\big)$.

A \ul{path} in $G$ with endpoints $u,v\in V$ is a collection of vertices $\gamma=\{u=v_0,v_1,...,v_k=v\}$ such that $(v_{i-1},v_i)\in E$. We say $u$ and $v$ are connected by $\gamma$. A graph $G$ is a \ul{connected graph} if every two points $u,v$ in $G$ are connected by a path.

In a graph $G$, a path $\gamma=\{u=v_0,v_1,...,v_k=v\}$ is a \ul{loop} if its endpoints coincide, i.e., if $u=v$. A graph $T$ is a \ul{tree} if it is (i) connected and (ii) contains no loops, or equivalently, if every two points in $T$ are connected by a unique path.

Given two graphs $G=(V,E)$ and $G'=(V',E')$, we say $G$ is a \ul{subgraph} of $G'$ (written $G\subset G'$) if $V\subset V'$ and $E\subset E'$. A \ul{subtree} is a subgraph (of a tree) that is (itself) a tree.

Given a graph $G$, a \ul{spanning tree} of $G$ is a subgraph $T\subset G$ such that (i) $T$ is a tree and (ii) $V_T=V_G$.
\end{dfn}

\begin{lmm}
A graph is connected if and only if it has a spanning tree.
\end{lmm}
\begin{proof}Let $G$ be a graph. If $G$ has a spanning tree, then it is clear that $G$ is connected since a tree is connected. Conversely, assume $G$ is a connected graph.

If $G$ is finite, then by repeatedly eliminating loops (where a loop is eliminated by removing/excluding one edge from the loop -- a process that clearly preserves connectedness of the graph) we obtain a spanning tree of $G$. If $G$ is infinite, we can proceed as follows.

Consider the set $\P:=\{\txt{Trees}~~T\subseteq G\}$ as a poset under inclusion $\subseteq$. If $\{T_\ld\}_{\ld\in\Ld}$ is a chain in $\P$ (i.e., $T_{\ld_1}\subseteq T_{\ld_2}$ if $\ld_1\leq\ld_2$), then $T_1:=\bigcup T_\ld$ is also in $\P$: Indeed, $T_1\subseteq G$ and $T_1$ cannot contain a loop, otherwise the loop will lie in some $T_\ld$, and so $T_1$ is a tree. That is, every chain in $\P$ has an upper bound in $\P$, and so by Zorn's Lemma, $\P$ has a maximal element $T$. Suppose $T=(V_T,E_T)$ is not a spanning tree of $G$. Then $V_G\backslash V_T\neq\emptyset$. Since $G$ is connected, there is an edge $e$ connecting $V_T$ to a vertex $v\in V_G\backslash V_T$. It follows that $T'=\big(\{v\}\cup V_T,\{e\}\cup E_T\big)$ is a tree in $\P$ strictly containing $T$ (a contradiction).
\end{proof}

\begin{prp}\label{GRPprp}
Every metric space has the gap reducing property.
\end{prp}
\begin{proof}
Let $X$ be a metric space and $A\subset X$ a finite set. Let $\rho$ denote $\rho(A)$. Given $a,b\in A$, let $a\sim b$ if there exists a sequence $a=a_0,a_1,\cdots,a_k=b$ such that $d(a_{i-1},a_i)<\rho$ for all $i=1,...,k$. Then $A=A_1\sqcup\cdots\sqcup A_m$, where $A_i$ are the equivalence classes under the equivalence relation $\sim$.

Consider the graph $G$ whose vertices are the equivalence classes $A_1,...,A_m$, with an edge connecting $A_i,A_j$ if $\dist(A_i,A_j)=\rho$. Then it is clear that $G$ is a connected graph (otherwise, the gap of $A$ will be larger than $\rho$).

Let $T$ be a spanning tree of $G$. Since a tree is precisely a connected graph in which any two vertices are connected by a unique path, by removing an edge $e$ from $T$, we can write ~$T-e=T'\sqcup T''$,~ where $T'$ and $T''$ are nonempty disjoint subtrees of $T$.

If $A'$ denotes the union of all vertices of $T'$, and $A''$ the union of all vertices of $T''$, then ~$A=A'\sqcup A''$, where $\rho(A)=\dist(A',A'')$ and $\rho(A'),\rho(A'')\leq\rho$ by construction.
\end{proof}


Note that the gap reducing property does not hold for the $d_H$-gap function $\ld:FS(X)\ra\Real$. To see this, let $X$ be a normed space and choose a set $A=\{a,b,c,d\}\in X(4)$ such that $\vep:=d(a,b)=d(c,d)$ is small and $L:=\dist(\{a,b\},\{c,d\})$ is large ($\gg\vep$). Then $\ld(A)=\vep$ and there exist no $A',A''\subset A=\{a,b,c,d\}$ such that $\ld(A)=d_H(A',A'')$ and $\ld(A'),\ld(A'')\leq\ld(A)$.

The following two definitions are not explicitly used anywhere in the thesis. However Definition \ref{GapDecomDfn} is potentially useful in thinking about Lipschitz retractions $X(n)\ra X(k)$ using a possible alternative/generalization of two-cluster decomposition in Definition \ref{TwoClustDec}. Similarly, Definition \ref{GromoHausDist} is potentially useful in answering further questions about the finite subset retraction (FSR) property in section \ref{FSRP2fs}. Moreover, certain results in this thesis admit generalizations, as it is the case with version II of the compact image lemma in Remark \ref{ComImLmm2}, the applications of which would involve Gromov-Hausdorff distance (Definition \ref{GromoHausDist}) as a possible extension of Hausdorff distance.

\begin{dfn}[\textcolor{blue}{Gap-decomposition, $k$-polar decomposition, $k$-polar elements}]\label{GapDecomDfn}
Let $A\in X(n)$. The gap-decomposition of $A$ is the decomposition into equivalence classes $A=A_1\sqcup\cdots\sqcup A_m$ considered in the proof of Proposition \ref{GRPprp}.

The gap-decomposition of $A=A_1\sqcup\cdots\sqcup A_m$ is $k$-polar (where $1\leq k\leq n-1$) if for some indices $i_1,...,i_k\in\{1,...,m\}$ we have $|A_{i_1}|=\cdots=|A_{i_k}|=1$ and $|A_i|\geq 2$ for all $i\not\in\{i_1,...,i_k\}$. In this case, we will say $A$ is a $k$-polar element of $X(n)$.

We denote the set of $k$-polar elements of $X(n)$ by $pX^k(n)$, and write $pX(n):=\bigcup_{k=1}^{n-1}pX^k(n)$ for the set of all polar elements.
\end{dfn}

\begin{dfn}[\textcolor{blue}{\index{Gromov-Hausdorff distance}{Gromov-Hausdorff distance}, Gromov-Hausdorff space}]\label{GromoHausDist}
Let $\C$ be the collection of compact metric spaces and $X,X'\in \C$. The smallest Hausdorff distance between isometric imbeddings of $X,X'$, denoted by $d_{GH}(X,X')$, is called the Gromov-Hausdorff distance between $X,X'$. That is,
\bea
d_{GH}(X,X'):=\inf\Big\{d_H\big(f(X),f'(X')\big):\txt{for isometric imbeddings}~X\sr{f}{\ral}Y,~X'\sr{f'}{\ral}Y\Big\}.\nn
\eea
With the collection of isometry classes  $[\C]=\{[X]:X\in\C\}$ of compact metric spaces, the resulting metric space $([\C],d_{GH})$ is called the Gromov-Hausdorff space.
\end{dfn}


\chapter{Quasiconvexity and Structure of Paths}\label{GQC}
This chapter relies on notation/terminology from section \ref{PrelimsMET} (especially the discussion on geodesics and quasigeodesics in subsection \ref{PrlMET3}), where the characterization of geodesics in Lemma \ref{GeodCharLmm} is particularly useful. We prove (in section \ref{GQCiq}, Theorem \ref{QConvThm}) that finite subset spaces of a quasiconvex metric space are themselves quasiconvex, and give a lower bound on the quasiconvexity constant. This result improves and generalizes \cite[Theorem 4.1]{BorovEtal2010} in which it was established that $\Real(n)$ is $4^n$-quasiconvex. We also give (in sections \ref{GQCcg}, \ref{GQCcq}, Propositions \ref{GeodExist}, \ref{QGeodExist} respectively) a detailed description of (quasi)geodesics in $X(n)$ for any metric space $X$. Finally, in section \ref{GQCmb} (Lemmas \ref{NoBIPLmm1}, \ref{NoBIPLmm2}), we briefly discuss failure of metrical convexity and the binary intersection property for finite subset spaces of normed spaces.

\section{Inherited quasiconvexity and lower bound on the constant}\label{GQCiq}
This section is based on \cite[Section 6, pp 24-31]{akofor2019}.
\begin{dfn}[\textcolor{blue}{\index{Spaced points}{Spaced pairs} of points in a metric space, $k$-spaced finite set}]
Let $(X,d)$ be a metric space. Two points $x,y\in X$ are spaced (or form a spaced pair) if
\begin{equation*}
\ol{N}_r(x)\cap\ol{N}_r(y)=\emptyset~~~~\txt{for all}~~~~0<r<d(x,y).
\end{equation*}
\end{dfn}
Equivalently, $x,y\in X$ are spaced $\iff$ $d(x,y)\leq \max\big\{d(x,z),d(z,y)\big\}$ for all $z\in X$. Note that because $\max\big\{d(x,z),d(z,y)\big\}=d_H\big(\{x,y\},\{z\}\big)$ and $d(x,y)=\rho(\{x,y\})$, the above definition can be generalized as follows: A finite set $F\subset X$ is \ul{$k$-spaced} if the $\rho(F)$-neighborhood, $N_{\rho(F)}^{d_H}\big(F\big):=\big\{F'\in FS(X)~|~d_H(F,F')<\rho(F)\big\}$, of $F$ in $FS(X)$ has empty intersection with $X(k)$, i.e.,
\bea
N_{\rho(F)}^{d_H}\big(F\big)\cap X(k)=\big\{F'\in X(k)~|~d_H(F,F')<\rho(F)\big\}=\emptyset,\nn
\eea
where $\rho:\P^\ast(X)\ra\Real$ is the gap function. In particular $x,y\in X$ are spaced $\iff$ $\{x,y\}\subset X$ is $1$-spaced.

\begin{lmm}\label{SpacedPairLmm}
If $X$ is a normed space, then for $n\geq 3$ the metric space $X(n)$ contains spaced pairs.
\end{lmm}
\begin{proof}
Since every normed space contains an isometrically imbedded copy of $\Real$, it is enough to show that $\Real(n)$ contains spaced pairs.

Fix a number $m>3$. Let $x=\{x_1,...,x_n\}\in \Real(n)\backslash \Real(n-1)$ be given by $x_1=0$, $x_2=m-1$, $x_3=m+1$, and $x_i=(i-2)m+1$ for $i\geq 4$. Similarly, let $y=\{y_1,...,y_n\}\in \Real(n)\backslash \Real(n-1)$ be given by $y_1=-1$, $y_2=1$, $y_3=m$, and $y_i=(i-2)m+2$ for $i\geq 4$. That is,
\begin{align*}
x&=\{x_1,...,x_n\}=\{0,m-1,m+1,2m+1,3m+1,\cdots, (n-2)m+1\},\\
y&=\{y_1,...,y_n\}=\{-1,1,m,2m+2,3m+2,\cdots, (n-2)m+2\}.
\end{align*}

\begin{figure}[H]
\centering
\scalebox{1}{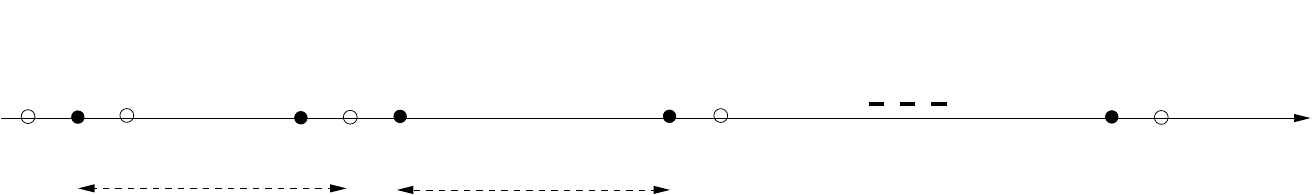} 
  \caption{~A spaced pair of points in $\Real(n)$.}\label{dg6}
\end{figure}

The points $x,y$ form a spaced pair because $\ol{N}_r(x)\cap \ol{N}_r(y)=\emptyset$ for all $0<r<1=d_H(x,y)$, where the proof is as follows.

Suppose on the contrary that $z=\{z_1,...,z_n\}\in\ol{N}_r(x)\cap\ol{N}_r(y)$ for some $0<r<1$. Consider the sets $A_1:=\{y_1,x_1,y_2\}=\{-1,0,1\}$, $A_2:=\{x_2,y_3,x_3\}=\{m-1,m,m+1\}$, and $A_k:=\{x_k,y_k\}=\{(k-1)m+1,(k-1)m+2\}$ for $3\leq k\leq n-1$. Then (by the definition of Hausdorff distance) $A_1,A_2$ each contains at least two elements of $z$, while $A_3,...,A_{n-1}$ each contains at least one element of $z$, i.e., at least $2+2+(n-3)=n+1$ elements of $z$ are required (a contradiction).
\end{proof}

\begin{rmk}\label{SpacedPairRmk}
Note that if $X$ is any geodesic space, then Lemma \ref{SpacedPairLmm} remains true. To prove this, let some geodesic segment in $X$ play the role of the real line in the proof of Lemma \ref{SpacedPairLmm}. Alternatively, observe that the existence of spaced pairs is an isometric property in the sense that spaced pairs are preserved by isometric homeomorphisms.
\end{rmk}

\begin{dfn}[\textcolor{blue}{\index{Complete relation}{Complete relation}, Incomplete relation, \index{Proximal relation}{Proximal relation}}]
Let $A,B$ be sets. A \ul{relation} on $A$ and $B$ is a set $R\subset A\times B$. Given a relation $R\subset A\times B$, let
\begin{align*}
&A_b(R):=\{a\in A:(a,b)\in R\},~~~~B_a(R):=\{b\in B:(a,b)\in R\},\nn\\
&A(R)\textstyle :=\bigcup_{b\in B}A_b(R)=\{a\in A:(a,b)\in R~\txt{for some}~b\in B\},~~~~\txt{(Left projection of $R$)}\\
&B(R)\textstyle :=\bigcup_{a\in A}B_a(R)=\{b\in B:(a,b)\in R~\txt{for some}~a\in A\},~~~~\txt{(Right projection of $R$)}.
\end{align*}
Then, we say $R$ is \ul{complete} if $A(R)=A$ and $B(R)=B$. Otherwise, $R$ is \ul{incomplete}.

If $X$ is a metric space and $A,B\subset X$, a relation $R\subset A\times B$ is \ul{proximal} if $d(a,b)\leq d_H(A,B)$ for all $(a,b)\in R$.
\end{dfn}

Note that by the definition of Hausdorff distance, for any $x,y\in X(n)$ there exists a proximal complete relation $R\subset x\times y$. This knowledge will be used in the proof of Corollary \ref{ProxSplitCrl}.

\begin{dfn}[\textcolor{blue}{\index{Proximal map}{Proximal map} between points of $X(n)$}]
Let $X$ be a metric space and $x,y\in X(n)$. A map $p : x\ra y$ is proximal if $d(a, p(a))\leq d_H(x,y)$ for all $a\in x$, i.e., the relation $R_p:=\big\{\big(a,p(a)\big):a\in x\big\}\subset x\times y$ is proximal.
\end{dfn}

\begin{dfn}[\textcolor{blue}{Orders of an element of a relation}]
Let $R\subset A\times B$ be relation and $(a,b)\in R$. The left and right orders of $(a,b)$ in $R$ are
\begin{equation*}
O_l(a,b):=|A_b(R)|=\big|(A\times\{b\})\cap R\big|,~~~~O_r(a,b):=|B_a(R)|=\big|(\{a\}\times B)\cap R\big|.
\end{equation*}
\end{dfn}

\begin{dfn}[\textcolor{blue}{Essential and Inessential elements of a complete relation}]
Let $R\subset A\times B$ be a complete relation. We say $(a,b)\in R$ is essential if $O_l(a,b)=1$ or $O_r(a,b)=1$. Otherwise, we say $(a,b)$ is inessential.
\end{dfn}
Note that if $R\subset A\times B$ is complete, then an element $(a,b)\in R$ is essential (resp. inessential) $\iff$ the relation $R\backslash\{(a,b)\}\subset A\times B$ is incomplete (resp. complete).

\begin{dfn}[\textcolor{blue}{\index{Reduced complete relation}{Reduced complete relation}}]
We say a complete relation $R\subset A\times B$ is reduced (or in reduced form) if every element of $R$ is essential.
\end{dfn}

\begin{prp}[\textcolor{OliveGreen}{Characterization of reduced complete relations}]\label{RedRelDecomp}
Let $R\subset A\times B$ be a complete relation. Then $R$ is reduced $\iff$ there exist disjoint union decompositions $A=A'\sqcup A_0\sqcup A''$, $B=B'\sqcup B_0\sqcup B''$ and maps $f:A'\ra B'$, $g:B''\ra A''$, $h:A_0\ra B_0$ such that $f,g$ are surjective, $h$ is bijective, and
\begin{equation*}
R=\Big\{\big(a,f(a)\big):a\in A'\Big\}\sqcup \Big\{\big(a,h(a)\big):a\in A_0\Big\}\sqcup\Big\{\big(g(b),b\big):b\in B''\Big\}.\nn
\end{equation*}
\end{prp}
\begin{proof}
($\Ra$): Assume that $R$ is reduced. Define sets $A_1',B_1''$ and maps $f_1:A_1'\ra B$, $g_1:B_1''\ra A$ as follows.
\begin{align*}
A_1'&:=\left\{a\in A:\big|B_a(R)\big|=1\right\},~~~~f_1(a):=\txt{the unique element $b\in B$ such that}~(a,b)\in R.\\
B_1''&:=\left\{b\in B:\big|A_b(R)\big|=1\right\},~~~~g_1(b):=\txt{the unique element $a\in A$ such that}~(a,b)\in R.
\end{align*}
Then, with $R_1:=\Big\{\big(a,f_1(a)\big):a\in A_1'\Big\}$, $R_2:=\Big\{\big(g_1(b),b\big):b\in B_1''\Big\}$, we have
\begin{equation*}
R=R_1\cup R_2=(R_1\backslash R_0)\sqcup R_0\sqcup(R_2\backslash R_0),
\end{equation*}
where $R_0:=R_1\cap R_2$ consists of elements $(u,v)\in R$ such that
\begin{align*}
(u,v)&=\big(a,f_1(a)\big)=\big(g_1(b),b\big)~~\Big[~\iff~~u=a=g_1(b),~~v=b=f_1(a)~\Big]\nn\\
     &~~~~\txt{for some}~~a\in A_1',~~b\in B_1''.
\end{align*}
Let $A_0:=A_1'\cap g_1(B_1'')$, $B_0:=f_1(A_1')\cap B_1''$. Then $f_1(A_0)=B_0$, $g_1(B_0)=A_0$, since
\begin{align*}
A_0=&\Big\{a\in A:~a=g_1(b)~\txt{for some}~b\in B,~\txt{with unique}~(a,f_1(a))\in R,~(g_1(b),b)\in R\Big\}\\
   =&\Big\{a\in A:~\txt{there exists}~b\in B~\txt{such that}~a=g_1(b),~b=f_1(a)\Big\}\\
   =&\Big\{a\in A:~\txt{there exists}~b\in B~\txt{such that}~(a,f_1(a))=(g_1(b),b)\in R\Big\},
\end{align*}
and similarly,
\begin{align*}
B_0=\Big\{b\in B:~\txt{there exists}~a\in A~\txt{such that}~(a,f_1(a))=(g_1(b),b)\in R\Big\}.
\end{align*}
It follows that $f_1|_{A_0}:A_0\ra B_0$, $g_1|_{B_0}:B_0\ra A_0$ are mutually inverse bijections, and
\begin{equation*}
R_0=\big\{(a,f_1(a)):a\in A_0\big\}=\big\{(g_1(b),b):b\in B_0\big\}.
\end{equation*}
Hence, we can set $h=f_1|_{A_0}$, $A'=A_1'\backslash A_0$, $f=f_1|_{A'}$, $B''=B_1''\backslash B_0$, $g=g_1|_{B''}$, $B'=f(A')$, $A''=g(B'')$.

($\La$): The converse is clear by the properties of the maps $f,g,h$.
\end{proof}

\begin{crl}[\textcolor{OliveGreen}{Characterization of reduced complete relations}]\label{RedRelDecomp2}
Let $R\subset A\times B$ be a complete relation. Then $R$ is reduced $\iff$ there exist disjoint unions $A=A'\sqcup A''$, $B=B'\sqcup B''$ and surjective maps $f:A'\ra B'$, $g:B''\ra A''$ such that
\bea
R=\Big\{\big(a,f(a)\big):a\in A'\Big\}\sqcup\Big\{\big(g(b),b\big):b\in B''\Big\}.\nn
\eea
\end{crl}

\begin{lmm}[\textcolor{OliveGreen}{Reduction of a finite complete relation}]\label{RelReductLmm}
Let $X$ be a metric space and $x,y\in X(n)$. Any (proximal) complete relation $R\subset x\times y$ contains a (proximal) reduced complete relation $R_{rc}\subset x\times y$, which means $R_{rc}\subset R$.
\end{lmm}
\begin{proof}
Since $R\subset x\times y$ is finite, we can obtain a reduced complete relation $R_{rc}\subset R\subset x\times y$ by repeatedly excluding inessential elements of $R$.
\end{proof}

\begin{crl}\label{ProxSplitCrl}
Let $X$ be a metric space and $x,y\in X(n)$. There exist proximal maps $f:x'\subset x\ra y$ and $g:y''\subset y\ra x$ such that $x=x'\sqcup g(y'')$ and $y=f(x')\sqcup y''$.
\end{crl}
\begin{proof}
By Lemmas \ref{RelReductLmm} and the definition of Hausdorff distance, a proximal reduced  complete relation $R\subset x\times y$ exists. Hence, by Corollary \ref{RedRelDecomp2}, the desired proximal maps exist.
\end{proof}

\begin{prp}[\textcolor{OliveGreen}{Sufficient condition for quasigeodesics in $X(n)$}]\label{QGeodExistSuff}
Let $X$ be a geodesic space and $x,y\in X(n)$. If there exists a complete relation $R\subset x\times y$  satisfying
\begin{equation}
\label{QGeodCond}|R|\leq n,~~~~d(a,b)\leq \ld d_H(x,y)~~\txt{for all}~~(a,b)\in R,
\end{equation}
then $x,y$ are connected by a $\ld$-quasigeodesic in $X(n)$.

\end{prp}
\begin{proof}
Assume some $R\subset x\times y$ satisfies (\ref{QGeodCond}). Then by Lemma \ref{GeodCharLmm}, the map $\gamma:[0,1]\ra X(n)$ given by
\begin{equation*}
\gamma(t):=\left\{\gamma_{(a,b)}(t):(a,b)\in R,~\gamma_{(a,b)}~\txt{a geodesic from $a$ to $b$}\right\}\nn
\end{equation*}
is a $\ld$-quasigeodesic from $x$ to $y$, since $\gamma(0)=x$, $\gamma(1)=y$, and for all $t,t'\in[0,1]$ we have
\begin{align*}
d_H(\gamma(t)&,\gamma(t'))=\max\left\{\!\max_{(a,b)\in R}\min_{(c,d)\in R}d\left(\gamma_{(a,b)}(t),\gamma_{(c,d)}(t')\right),\max_{(c,d)\in R}\min_{(a,b)\in R}d\left(\gamma_{(a,b)}(t),\gamma_{(c,d)}(t')\right)\!\right\}\\
   \leq&\max_{(a,b)\in R}d\left(\gamma_{(a,b)}(t),\gamma_{(a,b)}(t')\right)=|t-t'|\max_{(a,b)\in R}d(a,b)\leq\ld|t-t'|d_H(x,y).\qedhere
\end{align*}
\end{proof}

\begin{lmm}[\textcolor{OliveGreen}{Geodesics via proximal reduced complete relations}]\label{GeodPRCRLmm}
Let $X$ be a geodesic space. Then any two finite sets $x,y\subset X$ are connected by a geodesic  $\gamma:[0,1]\ra X(N)$, where $N:=\max(|x|,|y|,|x|+|y|-2)$. In particular, any two points $x,y\in X(n)$ are connected by a geodesic in $X\big(\max(n,2n-2)\big)$.
\end{lmm}
\begin{proof}
By Corollaries \ref{RedRelDecomp2} and \ref{ProxSplitCrl}, there exist proximal maps $f:x'\subset x\ra y$, $g:y''\subset y\ra x$ such that $x=x'\sqcup g(y'')$, $y=f(x')\sqcup y''$, and a proximal reduced complete relation $R\subset x\times y$ such that $R=\{(a,f(a)):a\in x'\}\sqcup\{(g(b),b):b\in y''\}$. Thus, we have the following 3 cases. (i) If $x'=\emptyset$, then $|R|=|y|$. (ii) If $y''=\emptyset$,  then $|R|=|x|$. (iii) If $x'\neq\emptyset$, $y''\neq\emptyset$, then
\begin{equation*}
|R|=|x'|+|y''|=|x|+|y|-\big(|f(x')|+|g(y'')|\big)\leq |x|+|y|-2.\nn
\end{equation*}
The conclusion now follows from Proposition \ref{QGeodExistSuff}.
\end{proof}

\begin{thm}[\textcolor{OliveGreen}{Quasiconvexity of $X(n)$}]\label{QConvThm}
If $X$ is a geodesic space, then $X(n)$ is 2-quasiconvex. Moreover, $X(2)$ is a geodesic space, and for $n\geq 3$, $\ld=2$ is the smallest quasiconvexity constant for $X(n)$. In particular, $X(n)$ for $n >2$ is never a geodesic space.
\end{thm}
\begin{proof}
Let $x,y\in X(n)$. If $x=y$, then $x,y$ are connected by the constant path. So, assume $x\neq y$. By Corollary \ref{ProxSplitCrl}, we have proximal maps $f:x'\subset x\ra y$, $g:y''\subset y\ra x$ such that, with $x'':=g(y'')$ and $y':=f(x')$, we have
\begin{equation*}
x=x'\sqcup g(y'')=x'\sqcup x'',~~~~y=f(x')\sqcup y''=y'\sqcup y''.
\end{equation*}
Let $z:=x''\cup y'$, for which we can verify that $d_H(x,z)\leq d_H(x',y')\leq d_H(x,y)$ and $d_H(z,y)\leq d_H(x'',y'')\leq d_H(x,y)$. Then $R_1=\{(a,f(a)):a\in x'\}\cup\{(c,c):c\in x''\}\subset x\times z$ and $R_2=\{(c,c):c\in y'\}\cup\{(g(b),b):b\in y''\}\subset z\times y$ are complete relations with the following properties.
\bit
\item $|R_1|\leq n$ and $d(u,v)\leq d_H(x,y)$ for each $(u,v)\in R_1$.
\item $|R_2|\leq n$ and $d(u,v)\leq d_H(x,y)$ for each $(u,v)\in R_2$.
\eit
If $z=x$ or $z=y$, then by Proposition \ref{QGeodExistSuff}, $x,y$ are connected by a quasigeodesic. So, assume $z\neq x$ and $z\neq y$. Then it follows again by Proposition \ref{QGeodExistSuff} that there exists a ${d_H(x,y)\over d_H(x,z)}$-quasigeodesic $\gamma_1:[0,1]\ra X(n)$ from $x$ to $z$, and there exists a ${d_H(x,y)\over d_H(z,y)}$-quasigeodesic $\gamma_2:[0,1]\ra X(n)$ from $z$ to $y$. Let {\small $\gamma=\gamma_1\cdot\gamma_2:[0,1]\ra X(n)$} be the path from $x$ to $y$ given by
\begin{equation*}
\gamma(t):=\left\{
             \begin{array}{ll}
               \gamma_1(2t), & \txt{if}~~t\in[0,1/2] \\
               \gamma_2(2t-1), & \txt{if}~~t\in[1/2,1]
             \end{array}
           \right\}.
\end{equation*}
Then
\begin{equation*}
\textstyle l(\gamma)=l(\gamma_1)+l(\gamma_2)\leq {d_H(x,y)\over d_H(x,z)}d_H(x,z)+{d_H(x,y)\over d_H(z,y)}d_H(z,y)=2d_H(x,y).
\end{equation*}
This shows that $X(n)$ is $2$-quasiconvex.

It follows from Lemma \ref{GeodPRCRLmm} that $X(2)$ is a geodesic space. Now let $n\geq 3$. To show $\ld=2$ is the smallest quasiconvexity constant for $X(n)$, let $x,y\in X(n)$ be a spaced pair (which exists by Remark \ref{SpacedPairRmk}). Then {\small $d_H(x,y)\leq\max\left\{d_H(x,z),d_H(z,y)\right\}$} for all $z\in X(n)$. If $\gamma:[0,1]\ra X(n)$ is a $\ld$-quasigeodesic from $x$ to $y$, i.e., {\small $d_H(\gamma(t),\gamma(t'))\leq\ld|t-t'|d_H(x,y)$}, then
\begin{align*}
d_H(x,&\gamma(1/2))\textstyle\leq{\ld\over 2}d_H(x,y),~~~~d_H(\gamma(1/2),y)\leq{\ld\over 2}d_H(x,y)\\
~~\Ra&\textstyle~~d_H(x,y)\leq\max\left\{d_H(x,\gamma(1/2)),d_H(\gamma(1/2),y)\right\}\leq{\ld\over 2}d_H(x,y),\\
~~\Ra&~~\ld\geq 2.\qedhere
\end{align*}
\end{proof}

\begin{crl}[\textcolor{OliveGreen}{Analog of Theorem \ref{QConvThm}}]\label{MetricQCQC}
If $X$ is an $\al$-quasiconvex metric space, then $X(n)$ is $2\al$-quasiconvex. Moreover, $X(2)$ is an $\al$-quasiconvex space, and for $n\geq 3$, $\ld=2$ is the smallest quasiconvexity constant for $X(n)$.
\end{crl}
\begin{proof}
Let $X$ be $\al$-quasiconvex. In the proof of Proposition \ref{QGeodExistSuff}, by replacing geodesics in $X$ with $\al$-quasigeodesics in $X$, we get the following fact.
\bit
\item \emph{Let $x,y\in X(n)$. If there exists a complete relation $R\subset x\times y$ satisfying
\bea
|R|\leq n,~~d(a,b)\leq \ld d_H(x,y)~~\txt{for all}~~(a,b)\in R,\nn
\eea
then $x,y$ are connected by an $\al\ld$-quasigeodesic in $X(n)$.}
\eit
Thus, in the proof of Theorem \ref{QConvThm}, by rescaling the Lipschitz constants of paths by $\al$, it follows that $X(n)$ is $2\al$-quasiconvex.

Using the obvious analog of Lemma \ref{GeodPRCRLmm} (with ``\emph{geodesic space}'' replaced by ``\emph{$\al$-quasiconvex space}''), we see that $X(2)$ is $\al$-quasiconvex.

Using Corollary \ref{ConstSpeedRP}, we can argue (as done in Remark \ref{SpacedPairRmk} when $X$ is a geodesic space) that $X(n)$ contains spaced pairs for $n\geq 3$. Thus, it follows as in the proof of Theorem \ref{QConvThm} that $2$ is the smallest quasiconvexity constant of $X(n)$ for $n\geq 3$.
\end{proof}

\section{Characterization of geodesics in finite subset spaces}\label{GQCcg}
Recall that if $X$ is a metric space, we write $FS(X)=\bigcup X(n)$ for the set of all finite subsets of $X$ as a metric space with respect to the Hausdorff distance $d_H$. Note that, even if $X$ is complete, $FS(X)$ is in general not complete.

\begin{note}\label{BndRelNote} 
If $X$ is a set and $x=\{x_1,..,x_n\},y=\{y_1,...,y_n\}\in X(n)$, then any index matching operation (not a map) of the form
\bea
\vphi:x\ra y,~x_i\mapsto y_{\sigma(i)},~~~~\txt{for a permutation}~~\sigma\in S_n,\nn
\eea defines a \emph{complete relation} $R:=\{(x_i,y_{\sigma(i)}):i=1,...,n\}\subset x\times y$ satisfying $|R|\leq n$.
\end{note}
Note however that the converse is not true, i.e., a complete relation $R\subset x\times y$ satisfying $|R|\leq n$ does not necessarily define an index matching operation of the form
\bea
\vphi:x\ra y,~x_i\mapsto y_{\sigma(i)},~~~~\txt{for a permutation}~~\sigma\in S_n.\nn
\eea
To see this, let $a,b,c,d,e\in X$ be distinct. In $X(4)$ let $x:=\{a,b,c,c\}$, $y:=\{d,d,e,e\}$, and $R:=\{(a,d),(a,e),(b,d),(c,e)\}$. Then $R$ is complete and $|R|=4$, but $x_1:=a$ cannot be mapped to a unique $y_j$ because ~$|\{j~|~(x_1,y_j)\in R\}|=2>1$.

\begin{lmm}[\textcolor{OliveGreen}{See the proof of Lemma \ref{GeodPRCRLmm}}]\label{GeodExistLmm1}
Let $(X,d)$ be a metric space. If $x,y\subset X$ are finite sets, then there exists a proximal complete relation $R\subset x\times y$ such that $|R|\leq \max(|x|,|y|,|x|+|y|-2)$.
\end{lmm}
The alternative proof (below) of the following special case of Lemma \ref{GeodExistLmm1} is interesting in its own right.
\begin{lmm}\label{GeodExistLmm2}
Let $(X,d)$ be a metric space. If $x,y\in FS(X)$, then there exists a complete relation $R\subset x\times y$ such that $\max\limits_{(u,v)\in R}d(u,v)=d_H(x,y)$, and $|R|\leq |x|+|y|$.

Moreover, if $|x|\geq 2$ and $|y|\geq 2$, then we can choose $R$ such that $|R|\leq|x|+|y|-2$.
\end{lmm}
\begin{proof}
Let $x\in X(n)\backslash X(n-1)$, $y\in X(m)\backslash X(m-1)$, and define
\bea
T_{nm}:=\big\{\txt{nonempty maps}~~\tau:\{1,...,n\}\ra\{1,...,m\}\big\}.\nn
\eea
Then by the fact that $x\cup y\subset \bigcup_{i,j}\ol{B}_{d_H(x,y)}(x_i)\cap \ol{B}_{d_H(x,y)}(y_j)$, we know there exist maps $x\sr{\vphi}{\ral}y$, $y\sr{\psi}{\ral}x$ (given by $\vphi(x_i)=y_{\tau(i)}$, $\psi(y_j)=x_{\eta(j)}$ for some $\tau\in T_{nm},\eta\in T_{mn}$) such that
\bea
&&d\left(x_i,y_{\tau(i)}\right)\leq d_H(x,y)~~~~\txt{for all}~~~~i=1,...,n,\nn\\
&&d\left(x_{\eta(j)},y_j\right)\leq d_H(x,y)~~~~\txt{for all}~~~~j=1,...,m.\nn
\eea
Let $R:=\left\{\left(x_i,y_{\tau(i)}\right):i=1,...,n\right\}\cup\left\{\left(x_{\eta(j)},y_j\right):j=1,...,m\right\}\subset x\times y$. We have the following two facts.
{\flushleft (i)} $|R|\leq|x|+|y|$, and we can choose the maps $\tau,\eta$ to satisfy ~$\max\limits_{(u,v)\in R}d(u,v)=d_H(x,y)$.
{\flushleft (ii)} Let $|x|,|y|\geq 2$. Then we have redundancies in the relation $R$ because the maps $\tau,\eta$ are nonempty. Let $R=R_1\cup R_2$, were
\bea
R_1:=\left\{\left(x_i,y_{\tau(i)}\right):i=1,...,n\right\},~~~~R_2:=\left\{\left(x_{\eta(j)},y_j\right):j=1,...,m\right\}.\nn
\eea
Since $n=|x|\geq 2$, the set $\{\tau(i):i=1,...,n\}$ is nonempty. Pick $j_0=\tau(i_0)\in\{\tau(i):i=1,...,n\}$ such that $d\left(x_{i_0},y_{j_0}\right)=\min_id\left(x_i,y_{\tau(i)}\right)$. Since $y_{j_0}$ is already paired as  $\left(x_{i_0},y_{j_0}\right)$ in $R_1$, we can remove $\left(x_{\eta(j_0)},y_{j_0}\right)$ from $R_2$ to get
\bea
R_1:=\left\{\left(x_i,y_{\tau(i)}\right):i=1,...,n\right\},~~~~\wt{R}_2:=\left\{\left(x_{\eta(j)},y_j\right):j\neq j_0\right\}.\nn
\eea
Similarly, because $m=|y|\geq 2$ implies the set $\{\eta(j):j\neq j_0\}$ is nonempty, we can pick $i_0'=\eta(j_0')\in\{\eta(j):j\neq j_0\}$ such that $d\left(x_{i_0'},y_{j_0'}\right)=\min_{j\neq j_0}d\left(x_{\eta(j)},y_j\right)$. Since $x_{i_0'}$ is already paired as $\left(x_{i_0'},y_{j_0'}\right)$ in $R_2$, we can remove $\left(x_{i_0'},y_{\tau(i_0')}\right)$ from $R_1$ to get
\bea
\wt{R}_1:=\left\{\left(x_i,y_{\tau(i)}\right):i\neq i_0'\right\},~~~~\wt{R}_2:=\left\{\left(x_{\eta(j)},y_j\right):j\neq j_0\right\}.\nn
\eea
Hence, we can replace $R$ with $\wt{R}:=\wt{R}_1\cup\wt{R}_2=R\backslash\left\{\left(x_{i_0'},y_{\tau(i_0')}\right),\left(x_{\eta(j_0)},y_{j_0}\right)\right\}\subset x\times y$.
\end{proof}

\begin{lmm}[\textcolor{OliveGreen}{Bounded image lemma}]\label{BndImgLmm}
Let $X$ be a metric space. If $K\subset X(n)$ is bounded, then ~$\bigcup K:=\bigcup_{x\in K}x$~ is bounded in $X$.
\end{lmm}
\begin{proof}
Let $R>0$. Fix $c\in X$. If $A\subset X$, then $B_R(A):=\{u\in X:\dist(u,A)<R\}$ denotes the $R$-neighborhood of $A$ in $X$. If $W\subset X(n)$, then $N_R(W):=\{x\in X(n):\dist_H(x,W)<R\}$ denotes the $R$-neighborhood of $W$ in $X(n)$. Note that
\bea
&&X_{R,c}(n):=\{x\in X(n):c\in x,~\diam x\leq R\}=\{x\in X(n):x\subset B_R(c)\}=N_R(\{c\})\nn\\
&&~~~~=\big(B_R(c)\big)(n).\nn
\eea
Since $K\subset X(n)$ is bounded, there exists $R>0$ such that $K\subset N_R(\{c\})=\big(B_R(c)\big)(n)$, i.e., $x\subset B_R(c)$ for all $x\in K$, and so $\bigcup_{x\in K}x\subset B_R(c)$.
\end{proof}

\begin{lmm}[\textcolor{OliveGreen}{\index{Compact image lemma}{Compact image lemma} I}]\label{ComImLmm}
Let $X$ be a metric space. If $\C\subset X(n)$ is compact, then $K:=\bigcup\C=\bigcup_{C\in \C}C$ is compact in $X$.
\end{lmm}

\begin{proof}
Let $\{x_k\}\subset K$ be a sequence. Then each $x_k\in C_k$ for some $C_k\in\C$. Since $\C$ is compact in $X(n)$, $\{C_k\}$ has a subsequence $\left\{C_{f(k)}\right\}$ that converges in $\C$. Let $C_{f(k)}\ra C_0\in\C$. Then because $C_0$ is compact (being a finite set), we have the following two facts.
\bit
\item[(i)] There exists a sequence $c_k\in C_0$ such that
\bea
d\left(x_{f(k)},c_k\right)=d\left(x_{f(k)},C_0\right)\leq d_H\left(C_{f(k)},C_0\right)\ra0.\nn
\eea
\item[(ii)] The sequence $c_k$ has a subsequence $c_{g(k)}$ that converges in $C_0$.
\eit
Let $c_{g(k)}\ra c_0\in C_0$. Then
\bea
d\left(x_{f\circ g(k)},c_0\right)\leq d\left(x_{f\circ g(k)},c_{g(k)}\right)+d\left(c_{g(k)},c_0\right)\ra0.\nn
\eea
That is, $\{x_k\}\subset K$ has the subsequence $x_{f\circ g(k)}$ that converges to $c_0\in C_0\subset K$.
\end{proof}

\begin{rmk}[\textcolor{OliveGreen}{Compact image lemma II}]\label{ComImLmm2}
By its proof, Lemma \ref{ComImLmm} holds in the following more general form: Let $X$ be a metric space and $\C^\ast(X)$ nonempty compact subsets of $X$. If $\C\subset\big(\C^\ast(X),d_H\big)$ is compact, then $K:=\bigcup\C=\bigcup_{C\in \C}C$ is a compact set in $X$.
\end{rmk}

\begin{lmm}[\textcolor{OliveGreen}{Pointwise convergent subsequence}]\label{PointWiseConv}
If $T$ is a countable set and $K$ a sequentially compact space (e.g., a compact metric space), then any given sequence of maps $f_k:T\ra K$ has a pointwise convergent subsequence $f_{s(k)}:T\ra K$.
\end{lmm}
\begin{proof}
Since $K$ is sequentially compact, for each $t\in T$, the sequence $\{f_k(t)\}$ has a convergent subsequence. Consider an enumeration $T=\{t_1,t_2,\cdots\}$. Then $\left\{f_k(t_1)\right\}$ has a convergent subsequence $\left\{f_{s_1(k)}(t_1)\right\}$, i.e., there exists a subsequence $\left\{f_{s_1(k)}\right\}\subset\{f_k\}$ such that $\left\{f_{s_1(k)}(t_1)\right\}$ converges. Similarly, because $\left\{f_{s_1(k)}(t_2)\right\}$ has a convergent subsequence, there exists a further subsequence $\left\{f_{s_2(k)}\right\}\subset \left\{f_{s_1(k)}\right\}\subset\{f_k\}$ such that $f_{s_2(k)}(t_2)$ converges. Continuing this way, we get subsequences $S_1\supset S_2\supset S_3\supset\cdots$ of $\{f_k\}$ which can be represented in an array as follows.
\bea
&& S_1=\left\{f_{s_1(k)}\right\}:~~~~~~f_{s_1(1)}~~~f_{s_1(2)}~~~f_{s_1(3)}~~~f_{s_1(4)}~~~\cdots~~~~\big(\txt{Converges on $\{t_1\}$}\big)\nn\\
&& S_2=\left\{f_{s_2(k)}\right\}:~~~~~~f_{s_2(1)}~~~f_{s_2(2)}~~~f_{s_2(3)}~~~f_{s_2(4)}~~~\cdots~~~~\big(\txt{Converges on $\{t_1,t_2\}$}\big)\nn\\
&& S_3=\left\{f_{s_3(k)}\right\}:~~~~~~f_{s_3(1)}~~~f_{s_3(2)}~~~f_{s_3(3)}~~~f_{s_3(4)}~~~\cdots~~~~\big(\txt{Converges on $\{t_1,t_2,t_3\}$}\big)\nn\\
&& S_4=\left\{f_{s_4(k)}\right\}:~~~~~~f_{s_4(1)}~~~f_{s_4(2)}~~~f_{s_4(3)}~~~f_{s_4(4)}~~~\cdots~~~~\big(\txt{Converges on $\{t_1,t_2,t_3,t_4\}$}\big)\nn\\
&&\vdots~~\hspace{2cm}~~\vdots~~\hspace{2cm}~~\vdots~~\hspace{2cm}~~\vdots~~\hspace{2cm}~~\vdots~~\hspace{2cm}~~\vdots\nn
\eea
Consider the diagonal sequence ~$S=\left\{f_{s_k(k)}\right\}:~~f_{s_1(1)}~~f_{s_2(2)}~~f_{s_3(3)}~~f_{s_4(4)}~~\cdots$~~.~~ Then for each $i=1,2,\cdots$, the sequence $S\cap S_i=S\backslash\left\{f_{s_1(1)},...,f_{s_{i-1}(i-1)}\right\}~\subset~S_i$~ (and hence $S$ also)~ converges on $\{t_1,...,t_i\}$. Thus, taking $i\ra\infty$, we see that $S=\bigcup_i(S\cap S_i)$ converges pointwise on $\{t_1,t_2,\cdots\}=T$. Define $f_{s(k)}:=f_{s_k(k)}$.
\end{proof}

\begin{lmm}[\textcolor{OliveGreen}{Constituent/component paths of a geodesic in $X(n)$}]\label{GeodExistLmm}
Let $X$ be a metric space, and $x,y\in X(n)$. If $\gamma:[0,1]\ra X(n)$ is a geodesic from $x$ to $y$, then for any $a\in x$, there exists a path $c:[0,1]\ra X$ satisfying the following.
\bit[leftmargin=0.9cm]
\item[(i)] $c(0)=a\in x$, ~~~$c(t)\in\gamma(t)$ ~for all ~$t\in[0,1]$. ~~(In particular, $c(1)\in\gamma(1)=y$.)
\item[(ii)] $d(c(t),c(t'))\leq|t-t'|d_H(x,y)$,~~ for all ~~$t,t'\in[0,1]$.
\eit
\end{lmm}
\begin{proof}
Let $\gamma:[0,1]\ra X(n)$ be a geodesic from $x$ to $y$, and let $\rho:=d_H(x,y)$. Then we have $\gamma(0)=x$, $\gamma(1)=y$, ~$|\gamma(t)|\leq n$ for all $t\in[0,1]$, ~and
\begin{align}
d_H(\gamma(t),\gamma(t'))=\max\left\{\max_{u\in\gamma(t)}\min_{u'\in\gamma(t')}d(u,u'),\max_{u'\in\gamma(t')}\min_{u\in\gamma(t)}d(u,u')\right\}=|t-t'|\rho,~~\forall~t,t'\in[0,1].\nn
\end{align}
For fixed $t,t'\in[0,1]$, this equation says for every $u\in\gamma(t)$, there exists $u'\in \gamma(t')$ such that
\bea
d(u,u')\leq|t-t'|\rho,~~~~\txt{(and vice versa)}.\nn
\eea
Let $D_k:=\{0=t_0<t_1<\cdots<t_k=1\}$, $k\geq 1$, be partitions of $[0,1]$ such that $D_k\subset D_{k+1}$ and $\bigcup D_k$ is dense in $[0,1]$ (e.g., $D_k=\{l/2^k:0\leq l\leq 2^k\}$). Fix $k\geq 1$. Then for each $a\in x$, we can define a map $g_k:D_k\ra X$ as follows. Let $g_k(t_0)=g_k(0):=a\in x=\gamma(0)$. Next, pick $g_k(t_1)\in\gamma(t_1)$ such that $d(g_k(t_0),g_k(t_1))\leq |t_0-t_1|\rho$. For the general step, given $g_k(t_i)$, pick $g_k(t_{i+1})\in\gamma(t_{i+1})$ such that $d(g_k(t_i),g_k(t_{i+1}))\leq|t_i-t_{i+1}|\rho$. This gives a map $g_k:D_k\ra X$ from $a\in x$ to some $b_k=g_k(1)\in y$ satisfying
\bea
d\left(g_k(t),g_k(t')\right)\leq|t-t'|\rho,~~~~\txt{for all}~~t,t'\in D_k.\nn
\eea

Consider the dense set $D:=\bigcup_{k=1}^\infty D_k$. For each $k$, let $f_k:D\ra X$ be an extension of $g_k:D_k\ra X$ such that $f_k(t)\in\gamma(t)$ for all $t\in D$. Then $f_k(D)\subset K:=\bigcup_{t\in[0,1]}\gamma(t)$. Since $K\subset X$ is compact (Lemma \ref{ComImLmm}) and $D$ is countable, it follows from Lemma \ref{PointWiseConv} that $\{f_k\}$ has a pointwise convergent subsequence $\left\{f_{s(k)}\right\}$, where we know $f_{s(k)}$ is $\rho$-Lipschitz on $D_{s(k)}$. Let $f_{s(k)}\ra f$. Then $f:D\ra X$ is $\rho$-Lipschitz and $f(t)\in\gamma(t)$, as follows: Indeed, given $t,t'\in D$, we can choose $N$ such that $t,t'\in D_{s(k)}$ for all $k\geq N$, and so
\bea
&&d\left(f_{s(k)}(t),f_{s(k)}(t')\right)\leq|t-t'|\rho,~~~~\txt{for all}~~k\geq N,\nn\\
&&\textstyle~~\Ra~~d(f(t),f(t'))\leq |t-t'|\rho,~~~~\txt{for all}~~~~t,t'\in D=\bigcup_kD_{s(k)},\nn\\
&&\dist(f(t),\gamma(t))\leq d\left(f(t),f_{s(k)}(t)\right)+\dist\left(f_{s(k)}(t),\gamma(t)\right)\nn\\
&&~~~~=d\left(f(t),f_{s(k)}(t)\right)\ra0,~~~~\txt{for all}~~t\in D,\nn\\
&&~~\Ra~~f(t)\in\gamma(t),~~~~\txt{for all}~~t\in D.\nn
\eea
Since $f$ is Lipschitz, and $D$ is dense in $[0,1]$, $f$ extends (by Lemma \ref{UniExtThm}) to a Lipschitz map $c:[0,1]\ra X$. It remains to show that $c(t)\in\gamma(t)$ for all $t\in[0,1]$, and that $c$ is $\rho$-Lipschitz.

Fix $t\in[0,1]$. Since $D$ is dense in $[0,1]$, pick $t_j\in D$ such that $t_j\ra t$. Then
\begin{align}
&\dist(c(t),\gamma(t))\leq d(c(t),c(t_j))+\dist(c(t_j),\gamma(t))=d(c(t),c(t_j))+\dist(f(t_j),\gamma(t))\nn\\
&~~~~\leq d(c(t),c(t_j))+d_H\big(\gamma(t_j),\gamma(t)\big)=d(c(t),c(t_j))+|t_j-t|\rho\ra 0,\nn\\
&~~\Ra~~c(t)\in\gamma(t).\nn
\end{align}
If there exist distinct $t_c,t_c'\in[0,1]$ such that $d(c(t_c),c(t_c'))>|t_c-t_c'|\rho$, then by continuity,
\bea
d(c(t),c(t'))>|t-t'|\rho~~\txt{for all}~~t\in[t_c-\vep,t_c+\vep],~t'\in[t_c'-\vep',t_c'+\vep'],~~\txt{for some}~~\vep,\vep'>0,\nn
\eea
which is a contradiction since $D$ is dense in $[0,1]$.
\end{proof}

\begin{figure}[H]
\centering
\scalebox{1.5}{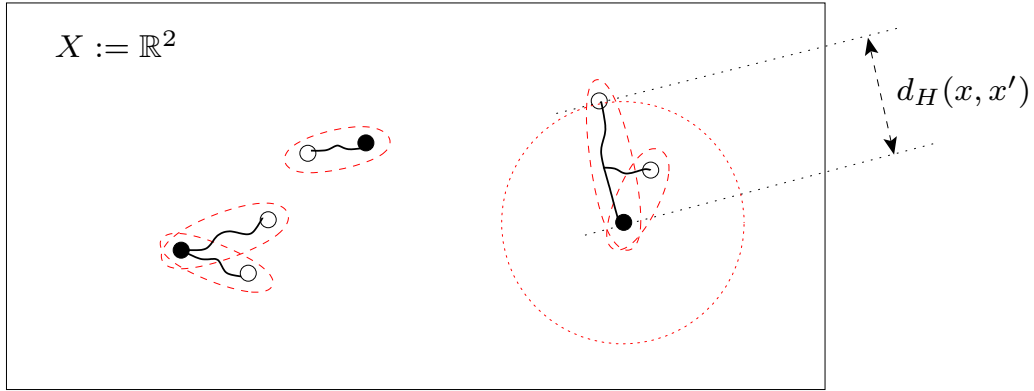} 
  \caption{~A typical path in $\Real^2(n)$, for $n\geq 5$, consisting of a (proximal) complete relation and component paths in $\Real^2$.}\label{dgdHR}
\end{figure}

\begin{prp}[\textcolor{OliveGreen}{Criterion for geodesics in $X(n)$}]\label{GeodExist}
Let $X$ be a metric space. For any $x,y\in X(n)$, a geodesic exists from $x$ to $y$ $\iff$ there exists a complete relation $R\subset x\times y$ and a collection of paths, $\left\{\gamma_{(a,b)}:(a,b)\in R\right\}$, $\gamma_{(a,b)}:[0,1]\ra X$ a path in $X$ from $a$ to $b$,  such that the following hold.
\bit[leftmargin=0.9cm]
\item[(i)] For each $(a,b)\in R$, we have ~~$d\left(\gamma_{(a,b)}(t),\gamma_{(a,b)}(t')\right)\leq|t-t'|d_H(x,y)~~\txt{for all}~~t,t'\in[0,1]$.

(Thus, $\max_{(a,b)\in R}d(a,b)\leq d_H(x,y)$.)

\item[(ii)] The path $\Gamma:[0,1]\ra FS(X)$ given by $\Gamma(t):=\left\{\gamma_{(a,b)}(t):(a,b)\in R\right\}$ lies in $X(n)$, i.e.,
\bea
|\Gamma(t)|\leq n~~~~\txt{for all}~~t\in[0,1].\nn
\eea
\eit
\end{prp}

\begin{proof}
Let $x,y\in X(n)$.

{\flushleft($\Ra$)}: Let $\gamma:[0,1]\ra X(n)$ be a geodesic from $x$ to $y$. Then by Lemma \ref{GeodExistLmm}, for each $\al\in x$, we have a path $c_\al:[0,1]\ra X$, $c_\al(t)\in\gamma(t)$, from $\al\in x$ to $c_\al(1)\in y$ satisfying
\bea
d\left(c_\al(t),c_\al(t')\right)\leq|t-t'|d_H(x,y),~~~~\txt{for all}~~t,t'\in[0,1].
\eea
Similarly, for each $\beta\in y$, we get a path $c^\beta:[0,1]\ra X$, $c^\beta(t)\in\gamma(1-t)$, from $\beta\in y$ to $c^\beta(1)\in x$ satisfying
\bea
d\left(c^\beta(t),c^\beta(t')\right)\leq|t-t'|d_H(x,y),~~~~\txt{for all}~~t,t'\in[0,1].
\eea
Let $R:=\left\{\big(c_\al(0),c_\al(1)\big):\al\in x\right\}\cup\left\{\big(\ol{c}^\beta(0),\ol{c}^\beta(1)\big):\beta\in y\right\}$, and consider the collection
\small\begin{align}
\left\{\gamma_{(a,b)}~\big|~(a,b)\in R\right\}~:=~\left\{\gamma_{\big(c_\al(0),c_\al(1)\big)}:=c_\al~\Big|~\al\in x\right\}\cup\left\{\gamma_{\big(\ol{c}^\beta(0),\ol{c}^\beta(1)\big)}:=\ol{c}^\beta~\Big|~\beta\in y\right\},\nn
\end{align}
where given a path $\sigma:[0,1]\ra Y$, we define $\ol{\sigma}:[0,1]\ra Y$ by $\ol{\sigma}(t):=\sigma(1-t)$. Then $R\subset x\times y$ is a complete relation, and we get a path ~$\Gamma:[0,1]\ra FS(X)$ satisfying $|\Gamma(t)|\leq n$ for all $t\in[0,1]$;
\bea
\Gamma(t)~:=~\left\{\gamma_{(a,b)}(t):(a,b)\in R\right\}=\{c_\al(t):\al\in x\}\cup\{c^\beta(1-t):\beta\in y\}~\subset~\gamma(t).
\eea

{\flushleft($\La$)}: Conversely, assume we have a complete relation $R\subset x\times y$ satisfying (i) and (ii). Then the map $\Gamma:[0,1]\ra X(n)$ given by $\Gamma(t):=\left\{\gamma_{(a,b)}(t):(a,b)\in R\right\}$ is a geodesic from $x$ to $y$, since $\Gamma(0)=x$, $\Gamma(1)=y$, and
{\footnotesize\begin{align}
&d_H\big(\Gamma(t),\Gamma(t')\big)=\max\left\{\max_{(a,b)\in R}\min_{(a',b')\in R}d\left(\gamma_{(a,b)}(t),\gamma_{(a',b')}(t')\right),\max_{(a',b')\in R}\min_{(a,b)\in R}d\left(\gamma_{(a,b)}(t),\gamma_{(a',b')}(t')\right)\right\}\nn\\
&~~~~\leq\max_{(a,b)\in R}d\left(\gamma_{(a,b)}(t),\gamma_{(a,b)}(t')\right)\leq|t-t'|d_H(x,y),~~~~\txt{for all}~~t,t'\in[0,1].\nn\qedhere
\end{align}}
\end{proof}

\begin{note}\label{GeodExistNt}
If $X$ is a geodesic space, then it is easy to see that for every $x\in X(1)$, $y\in X(n)$, there exists a geodesic $\gamma:[0,1]\ra X(n)$ between $x$ and $y$.
\end{note}

\begin{prp}[\textcolor{OliveGreen}{Sufficient condition for geodesics in $X(n)$}]\label{GeodExistSuff}
Let $X$ be a geodesic space and $x,y\in X(n)$. If there exists a complete relation $R\subset x\times y$  satisfying
\bea
\label{GeodCond}|R|\leq n,~~~~d(a,b)\leq d_H(x,y)~~\txt{for all}~~(a,b)\in R,
\eea
then $x,y$ are connected by a geodesic in $X(n)$.
\end{prp}
\begin{proof}
Set $\ld=1$ in Proposition \ref{QGeodExistSuff}.
\end{proof}

\begin{rmk}[\textcolor{OliveGreen}{Counterexample}]
The sufficient condition in Proposition \ref{GeodExistSuff} is not necessary (i.e., there exist geodesics that do not satisfy this condition).
\end{rmk}
\begin{proof}
In $\Real(3)$, consider the points $x=\{-1,1,m+3/2\}$ and $y=\{-3/2,m-1,m+1\}$, where $m>5$. Then $d_H(x,y)=5/2$. Observe that the sets $U=\{-3/2,-1,1\}$ and $V=\{m-1,m+1,m+3/2\}$ are separated by a distance of $m-2$, which is greater than $d_H(x,y)$. Thus, a complete relation $R\subset x\times y$ satisfying $d(a,b)\leq d_H(x,y)$ for all $(a,b)\in R$ can only pair elements within $U$ and within $V$. The only such complete relation is
\bea
R=\left\{(1,-3/2),(-1,-3/2),(m+3/2,m-1),(m+3/2,m+1)\right\}\subset x\times y.\nn
\eea
Since $|R|=4>3$, the sufficient condition (\ref{GeodCond}) cannot be satisfied.

\begin{figure}[H]
\centering
\scalebox{1}{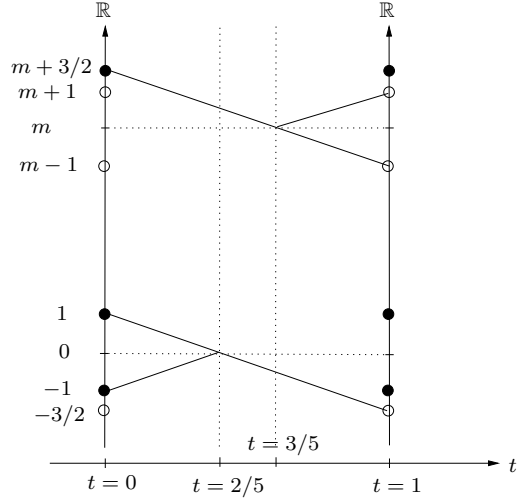} 
  \caption{~A geodesic in $\Real(3)$ with component paths that split up.}\label{dg1}
\end{figure}

However, the four paths $\gamma_{(1,-3/2)}, \gamma_{(-1,-3/2)}, \gamma_{(m+3/2,m-1)}, \gamma_{(m+3/2,m+1)}:[0,1]\ra\Real$ defined below satisfy the requirements of Proposition \ref{GeodExist}, and so there is a geodesic between $x,y$.
\begin{align}
&\textstyle\gamma_{(1,-3/2)}(t):=(1-t)1+t(-3/2)=1-{5\over 2}t,\nn\\
&\textstyle\gamma_{(-1,-3/2)}(t):=
\left\{
  \begin{array}{ll}
    (1-t)(-1)+t(3/2), & t\in[0,2/5] \\
    \gamma_{(1,-3/2)}(t), & t\in[2/5,1]
  \end{array}
\right\}=
\left\{
  \begin{array}{ll}
    {5\over 2}t-1, & t\in[0,2/5] \\
    \gamma_{(1,-3/2)}(t), & t\in[2/5,1]
  \end{array}
\right\},\nn\\
&\textstyle\gamma_{(m+3/2,m-1)}(t):=(1-t)(m+3/2)+t(m-1)=m+{3\over 2}-{5\over 2}t,\nn\\
&\textstyle\gamma_{(m+3/2,m+1)}(t):=\left\{
                            \begin{array}{ll}
                              \gamma_{(m+3/2,m-1)}(t), & t\in[0,3/5]\\
                             (1-t)(m-3/2)+t(m+1) , & t\in[3/5,1]
                            \end{array}
                          \right\}\nn\\
&\textstyle~~~~=\left\{
                            \begin{array}{ll}
                              \gamma_{(m+3/2,m-1)}(t), & t\in[0,3/5]\\
                             m-{3\over 2}+{5\over 2}t , & t\in[3/5,1]
                            \end{array}
                          \right\}.\nn\qedhere
\end{align}
\end{proof}

\begin{crl}[\textcolor{OliveGreen}{Self geodesics of a set in $X(n)$}]\label{SelfGeod}
Let $X$ be a geodesic space. If $x\in X(n)$ and $x'\subset x$, then a geodesic $\gamma:[0,1]\ra X(n)$ exists from $x'$ to $x$.
\end{crl}
\begin{proof}
We know that for any $a\in x\backslash x'$, there exists $a'_a\in x'$ such that $d\big(a'_a,a\big)\leq d_H(x',x)$. Thus, we have the complete relation $R=\{(a',a'):a'\in x'\}\cup\left\{\big(a'_a,a\big):a\in x\backslash x'\right\}\subset x'\times x$ satisfying (\ref{GeodCond}).
\end{proof}

\begin{crl}[\textcolor{OliveGreen}{Sufficient condition for geodesics in $X(n)\backslash X(n-1)$}]\label{GeodExistCor}
Let $X$ be a geodesic space and $x,y\in X(n)\backslash X(n-1)$. If there exists a bijection $\vphi:x\ra y$, $x_i\mapsto\vphi(x_i)=y_{\sigma(i)}$, for some $\sigma\in S_n$, such that ~$d\left(x_i,y_{\sigma(i)}\right)\leq d_H(x,y)$ for all $i$, then $x,y$ are connected by a geodesic in $X(n)$.
\end{crl}
\begin{proof}
Given points $x,y\in X(n)\backslash X(n-1)$, we have a complete relation $R\subset x\times y$ satisfying $|R|\leq n$ if and only if we have a bijection $\vphi:x\ra y$ (given by $\vphi(x_i)=y_{\sigma(i)}$ for some permutation $\sigma\in S_n$) such that $R=\left\{\left(x_i,y_{\sigma(i)}\right):i=1,...,n\right\}$.
Hence, the result follows by Proposition \ref{GeodExistSuff}.
\end{proof}

\begin{crl}
Let $X$ be a geodesic space and $x,y\in X(n)$. If $\delta_n(x)>2d_H(x,y)$ or $\delta_n(y)>2 d_H(x,y)$, then (by Corollary \ref{GeodExistCor} and Lemma \ref{HausDistBound}) a geodesic $\gamma:[0,1]\ra X(n)$ exists between $x$ and $y$.
\end{crl}

\begin{dfn}[\textcolor{blue}{\index{Antipodal points}{Antipodal (or opposite) points} in a geodesic space}]
Let $X$ be a geodesic space and $x,y,z\in X$. Then we say $x$ and $z$ are \ul{antipodal(ly located) about $y$} if there exists a geodesic $\gamma:[0,1]\ra X$ from $\gamma(0)=x$ through $\gamma(t_0)=y$ (for some $t_0\in(0,1)$) to $\gamma(1)=z$ such that $d(x,y)=d(y,z)$.
\end{dfn}

The following result was established in Theorem \ref{QConvThm}, but the alternative proof below is also interesting.
\begin{crl}\label{GeodExistCor1}
Let $X$ be a geodesic space. If $n\geq 3$, then $X(n)$ is not a geodesic space.
\end{crl}
\begin{proof}
For a fixed $\vep>0$ choose points $x=\{x_1,...,x_n\}$, $y=\{y_1,...,y_n\}$ in $X(n)\backslash X(n-1)$ such that the associated sets $A:=\{x_1,...,x_{n-1},y_n\}$, $B:=\{y_1,...,y_{n-1},x_n\}$ in $X$ satisfy
\bea
diam(A)\leq\vep,~~~~diam(B)\leq\vep,~~~~\dist(A,B):=\min\Big\{d(a,b):a\in A,~b\in B\Big\}>\vep.\nn
\eea
Then the only complete relation $R\subset x\times y$ satisfying condition (i) of Proposition \ref{GeodExist} is
\bea
R=\{(x_1,y_n),...,(x_{n-1},y_n)\}\cup\{(x_n,y_1),....,(x_n,y_{n-1})\}.\nn
\eea
Since $R$ contains $2n-2$ elements, it is easy to arrange the points in $A,B\subset X$ such that any associated collection of paths $\left\{\gamma_{(a,b)}:(a,b)\in R\right\}$ violates condition (ii) of Proposition \ref{GeodExist}. One way is to let the points $x_1,...,x_{n-1}\in A$ (respectively, $y_1,...,y_{n-1}\in B$) be distributed in a small neighborhood of the boundary of a ball of radius $\vep/2$ centered at $y_n\in A$ (respectively, $x_n\in B$), in such a way that at least two of the $x_i$'s (resp. $y_j$'s) are antipodally located about $y_n$ (resp. $x_n$). That is, with at least two antipodal points,
\bea
x_1,...,x_{n-1}\in N_\delta\left(\del\ol{B}_{\vep/2}(y_n)\right),~~~~y_1,...,y_{n-1}\in N_\delta\left(\del\ol{B}_{\vep/2}(x_n)\right),~~~~\delta<<\vep,\nn
\eea
where in the case with $\dim X\geq 2$, we can let $\delta\ra0$, and so have
\bea
x_1,...,x_{n-1}\in \del\ol{B}_{\vep/2}(y_n),~~~~y_1,...,y_{n-1}\in \del\ol{B}_{\vep/2}(x_n).\nn
\eea
Let the two $x_i$'s that are antipodal about $y_n$ in $A$ be $x_{i_1},x_{i_2}$, and the two $y_j$'s that are antipodal about $x_n$ in $B$ be $y_{j_1},y_{j_2}$. Then it is clear that these points satisfy
\bea
d(y_n,x_{i_1})=d(y_n,x_{i_2})=d_H(x,y)=\vep/2,~~~~d(x_n,y_{j_1})=d(x_n,y_{j_2})=d_H(x,y)=\vep/2.\nn
\eea
Therefore, at least four of the component paths, namely, $\gamma_{(x_{i_1},y_n)}$, $\gamma_{(x_{i_2},y_n)}$, $\gamma_{(x_n,y_{j_1})}$, $\gamma_{(x_n,y_{j_2})}$, are necessarily geodesics in $X$ that satisfy
\bea
&& d\left(\gamma_{(x_{i_1},y_n)}(t),\gamma_{(x_{i_1},y_n)}(t')\right)= d_H(x,y)|t-t'|,~~~~\txt{for all}~~t,t'\in[0,1],\nn\\
&& d\left(\gamma_{(x_{i_2},y_n)}(t),\gamma_{(x_{i_2},y_n)}(t')\right)= d_H(x,y)|t-t'|,~~~~\txt{for all}~~t,t'\in[0,1],\nn\\
&& d\left(\gamma_{(x_n,y_{j_1})}(t),\gamma_{x_n,y_{j_1}}(t')\right)= d_H(x,y)|t-t'|,~~~~\txt{for all}~~t,t'\in[0,1],\nn\\
&& d\left(\gamma_{(x_n,y_{j_2})}(t),\gamma_{x_n,y_{j_2}}(t')\right)= d_H(x,y)|t-t'|,~~~~\txt{for all}~~t,t'\in[0,1],\nn
\eea
and so also satisfy the inter-path distances
\bea
&&\left|d\left(\gamma_{(x_{i_1},y_n)}(t),\gamma_{(x_{i_1},y_n)}(t)\right)-d\left(\gamma_{(x_{i_2},y_n)}(t'),\gamma_{(x_{i_2},y_n)}(t')\right)\right|\leq 2d_H(x,y)|t-t'|,\nn\\
&&\left|d\left(\gamma_{(x_n,y_{j_1})}(t),\gamma_{x_n,y_{j_1}}(t)\right)-d\left(\gamma_{(x_n,y_{j_2})}(t'),\gamma_{x_n,y_{j_2}}(t')\right)\right|\leq 2d_H(x,y)|t-t'|.\nn
\eea
If two of the four paths, say $\gamma_1$ and $\gamma_2$ (denoting either $\gamma_{(x_{i_1},y_n)}$ and $\gamma_{(x_{i_2},y_n)}$ in $A$, or $\gamma_{(x_n,y_{j_1})}$ and $\gamma_{(x_n,y_{j_2})}$ in $B$), joint up at some time $t'\in(0,1)$, then
\bea
&&|d(\gamma_1(t),\gamma_2(t))-0|\leq 2d_H(x,y)|t-t'|,~~\Ra~~d(\gamma_1(0),\gamma_2(0))\leq 2d_H(x,y)t'\nn\\
&&~~\txt{and}~~d(\gamma_1(1),\gamma_2(1))\leq 2d_H(x,y)|1-t'|,\nn
\eea
which gives a contradiction since the antipodal location of the endpoints implies
\bea
d(\gamma_1(0),\gamma_2(0))>2d_H(x,y)t'~~\txt{or}~~d(\gamma_1(1),\gamma_2(1))>2d_H(x,y)|1-t'|,~~\txt{for any}~~t'\in(0,1).\nn
\eea
This shows it is impossible for $\gamma_{(x_{i_1},y_n)}$, $\gamma_{(x_{i_2},y_n)}$ (and also impossible for $\gamma_{(x_{i_1},y_n)}$, $\gamma_{(x_{i_2},y_n)}$) to join up into a single path. Hence, condition (ii) of Proposition \ref{GeodExist} cannot be satisfied.

\begin{figure}[H]
\centering
\scalebox{0.8}{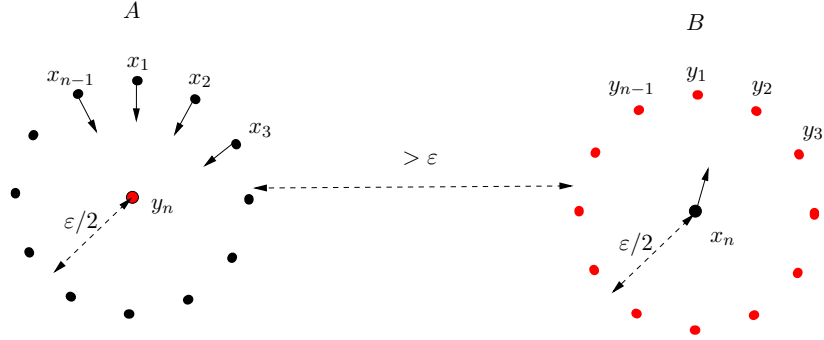} 
\caption{~Impossibility of a geodesic: All $n-1$ component paths from $x_1,..,x_{n-1}$ must joint up into a single path while the lone component path from $x_n$ must split up into $n-1$ paths, resulting in a collection of component paths that violates condition (ii) of Proposition \ref{GeodExist}.}\label{dg4}
\end{figure}

\end{proof}

\begin{rmk}[\textcolor{OliveGreen}{Non 1-ALR's}]\label{BanMetConvX(n)}
If $X$ is a Banach space, then for $n\geq 3$, $X(n)$ is not a 1-ALR. Indeed we know $X(n)$ is complete but not geodesic (hence not metrically convex due to Proposition \ref{ComMetricConv}), and so cannot be a 1-ALR by Proposition \ref{1ALRcrit}.
\end{rmk}

\begin{dfn}[\textcolor{blue}{\index{Geodesic! gap}{Geodesic gap of $X(n)$} in $FS(X)$}]
Let $X$ be a geodesic space. By Lemma \ref{GeodPRCRLmm}, for each pair of points $x,y\in X(n)$, there exists a smallest number $m_n(x,y)\geq n$ such that $X\big(m_n(x,y)\big)$ contains a geodesic between $x$ and $y$. That is,
\bea
m_n(x,y):=\min\big\{k\geq n~\big|~\txt{a geodesic $\gamma:[0,1]\ra X(k)$ exists with $\gamma(0)=x$, $\gamma(1)=y$}\big\}.\nn
\eea
(Note that if $m_n(x,y)>n$, then $X\big(m_n(x,y)-1\big)$, and hence each $X(k)$ for $n\leq k\leq m_n(x,y)-1$, does not contain a geodesic between $x$ and $y$.)

We define the \ul{geodesic gap} of $X(n)$ in $FS(X)$ to be the smallest number $m_n\geq n$ such that $X(m_n)$ contains a geodesic between all points $x,y\in X(n)$, i.e.,
\bea
\textstyle m_n:=\max_{x,y\in X(n)}m_n(x,y).\nn
\eea
\end{dfn}

\begin{prp}[\textcolor{OliveGreen}{Value of the geodesic gap}]
Let $X$ be a geodesic space. For $n\geq 2$, the geodesic gap of $X(n)$ in $FS(X)$ is $m_n=2n-2$.
\end{prp}
\begin{proof}
By the example used in the proof of Corollary \ref{GeodExistCor1}, we have the lower bound
\bea
m_n=\max_{x,y\in X(n)}m_n(x,y)\geq n+(n-2)=2n-2,~~~~\txt{for any}~~~~n\geq 2.\nn
\eea
Also, by Lemma \ref{GeodExistLmm2}, we have the upper bound ~$m_n\leq \max_{x,y\in X(n)}\big(|x|+|y|-2\big)= 2n-2$. Hence, $m_n=2n-2$.
\end{proof}

\begin{dfn}[\textcolor{blue}{\index{Distantly separated sets}{Distantly separated sets}}]
Let $X$ be a metric space. We say bounded sets $A_1,...,A_k\subset X$ are distantly separated if the largest diameter of the sets is less than the smallest distance between any two of them, i.e.,
\bea
\max_i~\diam(A_i)~<~\min_{i\neq j}~\dist(A_i,A_j):=\min_{i\neq j}\big\{d(a,b):a\in A_i,b\in A_j\big\}.\nn
\eea
\end{dfn}

\begin{prp}\label{GeodNumBound}
Let $X$ be a geodesic space. If $n\geq 3$, then for any given integer $n+1\leq m\leq 2n-2$, we can find two points $x,y\in X(n)$ such that the following hold.
\bit
\item[(i)] $x$ and $y$ are not connected by a geodesic in $X(m-1)$.
\item[(ii)] $x$ and $y$ are connected by a geodesic in $X(m)$.
\eit
\end{prp}
\begin{proof}
The main idea is to consider points $x,y\in X(n)\backslash X(n-1)$ for which the set $x\cup y$ in $X$ splits up into distantly separated sets in $X$, after which we then apply the procedure used in the proof of Corollary \ref{GeodExistCor1} for the \emph{nonexistence of a geodesic} (call this procedure $NP$), and Proposition \ref{GeodExistSuff} for the \emph{existence of a geodesic} (call this proposition $EP$). All conclusions below about existence/nonexistence of a geodesic are, without further mention (for brevity), based on the items $NP$ and $EP$ from the above discussion. We will proceed by induction on $n$.

For $n=3$, $m$ takes only one value, namely, $m=4$. For this case, choose $x=\{x_1,x_2,x_3\}$, $y=\{y_1,y_2,y_3\}$ in $X(3)\backslash X(2)$ such that the sets $A=\{x_1,x_2,y_3\}$ and $B=\{y_1,y_2,x_3\}$ are distantly separated in $X$. Then $x,y$ are connected by a geodesic in $X(4)$ but not in $X(3)$.

For $n=4$, $m\in\{5,6\}$. With $m=5$, choose $x=\{x_1,x_2,x_3,x_4\}$, $y=\{y_1,y_2,y_3,y_4\}$ in $X(4)\backslash X(3)$ such that the sets $A=\{x_1,x_2,y_3\}$, $B=\{y_1,y_2,x_3\}$, $C=\{x_4,y_4\}$ are distantly separated in $X$. Similarly, with $m=6$, choose $u=\{u_1,u_2,u_3,u_4\}$, $v=\{v_1,v_2,v_3,v_4\}$ in $X(4)\backslash X(3)$ such that the sets $D=\{u_1,u_2,u_3,v_4\}$, $E=\{v_1,v_2,v_3,u_4\}$ are distantly separated in $X$. Then $x,y$ are connected by a geodesic in $X(5)$ but not in $X(4)$, while $u,v$ are connected by a geodesic in $X(6)$ but not in $X(5)$.

For general $n$, $m\in\{n+1,\cdots,2n-2\}$. With $n+1\leq m\leq (2n-2)-1=2n-3$, choose $x=x'\sqcup\{x_n\}$, $y=y'\sqcup\{y_n\}$ in $X(n)\backslash X(n-1)$, where $x',y'\in X(n-1)\backslash X(n-2)$, such that the set $A=\{x_n,y_n\}$ and the subsets considered in the $AB$-splitting of $x'\cup y'$ (for $m':=m-1$) used in the previous induction stage (i.e., $(n-1)$th induction stage) are distantly separated in $X$. Similarly, with $m=2n-2$, choose $x=\{x'\}\sqcup\{x_n\}$, $y=y'\sqcup\{y_n\}$ in $X(n)\backslash X(n-1)$, where $x',y'\in X(n-1)\backslash X(n-2)$, such that the sets $A=x'\cup\{y_n\}$, $B=y'\cup\{x_n\}$ are distantly separated in $X$. Then for each $m\in \{n+1,\cdots,2n-2\}$, the conditions (i) and (ii) of the proposition are satisfied.
\end{proof}

\begin{crl}
Let $X$ be a geodesic space. For $n\geq 3$, there is no 1-Lipschitz retraction
\bea
X(m)\ra X(n),~~~~\txt{for all}~~~~n+1\leq m\leq 2n-2.\nn
\eea
\end{crl}
\begin{proof}
Fix $m\in\{n+1,\cdots,2n-2\}$. Then by Proposition \ref{GeodNumBound}, there are distinct points $x,y\in X(n)$ that are not connected by a geodesic in $X(m-1)$ but are connected by a geodesic $\gamma:[0,1]\ra X(m)$ such that $\gamma(0)=x$, $\gamma(1)=y$. Let $z\in\gamma([0,1])$ such that $z\neq x$ and $z\neq y$. Suppose there exists a $1$-Lipschitz retraction $r:X(m)\ra X(n)$. Then
\begin{align}
&d_H(x,z)+d_H(z,y)=d_H(x,y)=d_H(r(x),r(y))\sr{(s)}{<}d_H(r(x),r(z))+d_H(r(z),r(y))\nn\\
&~~=d_H(x,r(z))+d_H(r(z),y),~~\Ra~~d_H(x,z)<d_H(r(x),r(z))~\txt{or}~d_H(z,y)<d_H(r(z),r(y)),\nn
\end{align}
which implies $\Lip(r)\geq\max\left\{{d_H(r(x),r(z))\over d_H(x,z)},{d_H(r(z),r(y))\over d_H(z,y)}\right\}>1$ (a contradiction), where step (s) above holds because the path $r\circ\gamma:[0,1]\sr{\gamma}{\ral}X(m)\sr{r}{\ra}X(n)$ from $x$ to $y$ is not a geodesic.
\end{proof}
\begin{crl}
Let $X$ be a geodesic space. For $n\geq 3$, there is no 1-Lipschitz retraction
\bea
X(m)\ra X(n),~~~~\txt{for all}~~~~m>n.\nn
\eea
\end{crl}

\section{Characterization of quasigeodesics in finite subset spaces}\label{GQCcq}
We have already seen that when $X$ is a geodesic space, then Proposition \ref{QGeodExistSuff} gives a sufficient condition for the existence of quasigeodesics in $X(n)$. The following is a full characterization of quasigeodesics in $X(n)$, for any metric space $X$.

\begin{prp}[\textcolor{OliveGreen}{Criterion for quasigeodesics in $X(n)$}]\label{QGeodExist}
Let $X$ be a metric space. For any $x,y\in X(n)$, a $\ld$-quasigeodesic exists from $x$ to $y$ $\iff$ there exists a complete relation $R\subset x\times y$ and a collection of paths $\left\{\gamma_{(a,b)}:(a,b)\in R\right\}$, $\gamma_{(a,b)}:[0,1]\ra X$ a path in $X$ from $a$ to $b$,
such that the following hold.
\bit[leftmargin=0.9cm]
\item[(i)] For each $(a,b)\in R$, we have ~$d\left(\gamma_{(a,b)}(t),\gamma_{(a,b)}(t')\right)\leq \ld|t-t'|d_H(x,y)$ for all $t,t'\in[0,1]$.
(Thus, $\max_{(a,b)\in R}d(a,b)\leq \ld d_H(x,y)$.)

\item[(ii)] The path $\Gamma:[0,1]\ra FS(X)$ given by $\Gamma(t):=\left\{\gamma_{(a,b)}(t):(a,b)\in R\right\}$ lies in $X(n)$, i.e.,
\bea
|\Gamma(t)|\leq n~~~~\txt{for all}~~t\in[0,1].\nn
\eea
\eit
\end{prp}
\begin{proof}
This is the same as the proof of Proposition \ref{GeodExist}.
\end{proof}

As elaborated in Theorem \ref{QConvThm}, if $X$ is a geodesic space, then $X(n)$ is 2-quasiconvex. Below we give a simpler proof of 3-quasiconvexity of $X(n)$.
\begin{crl}[\textcolor{OliveGreen}{$3$-quasiconvexity of $X(n)$}]\label{SubSpaceQC}
If $X$ is a geodesic space, then each $X(n)$ is $3$-quasiconvex.
\end{crl}
\begin{proof}
Let $x,y\in X(n)$. Recall that for each $x_i\in x$, there exists $y(i)\in y$ such that $d(x_i,y(i))\leq d_H(x,y)$. Let $y':=\{y(i):i=1,...,n\}\subset y$, where wlog $x\neq y'$. Then
\begin{align}
&d_H(x,y')=\max\left\{\max_i\min_jd(x_i,y(j)),\max_j\min_id(x_i,y(j))\right\}\leq \max_id(x_i,y(i))\leq d_H(x,y),\nn\\
&d_H(y',y)\leq d_H(y',x)+d_H(x,y)\leq 2d_H(x,y).\nn
\end{align}
Thus, by Proposition \ref{QGeodExistSuff}, a ${d_H(x,y)\over d_H(x,y')}$-quasigeodesic $\gamma_1:[0,1]\ra X(n)$ exists from $x$ to $y'$, and we also know (by Corollary \ref{SelfGeod}) that a geodesic $\gamma_2:[0,1]\ra X(n)$ exists from $y'$ to $y$. Consider the path $\gamma:[0,1]\ra X(n)$ from $x$ to $y$ given by
\bea
\gamma(t):=(\gamma_1\cdot\gamma_2)(t):=\left\{
            \begin{array}{ll}
              \gamma_1(2t), & t\in[0,1/2] \\
              \gamma_2(2t-1), & t\in[1/2,1]
            \end{array}
          \right.\nn
\eea
Then {\small~$l(\gamma)=l(\gamma_1)+l(\gamma_2)\leq{{d_H(x,y)\over d_H(x,y')}}d_H(x,y')+d_H(y',y)\leq d_H(x,y)+2d_H(x,y)=3d_H(x,y)$}.
\end{proof}

\section{Metrical convexity and binary intersection property}\label{GQCmb}
We have seen that a complete metric space is metrically convex if and only if geodesic. We have also seen that if $X$ is geodesic, then $X(n)$ is not geodesic for $n\geq 3$. It follows that if $X$ is a Banach space, then $X(n)$ is not metrically convex for $n\geq 3$. We will also show that the binary intersection property (BIP) fails in $X(n)$, $n\geq 2$, if $X$ is a normed space of dimension two or more.

\begin{lmm}[\textcolor{OliveGreen}{Failure of BIP for $X(2)$}]\label{NoBIPLmm1}
Let $X$ be a normed space such that $\dim X\geq 2$. Then $X(2)$ does not have the BIP.
\end{lmm}
\begin{proof}
Since every normed space of dimension $\geq 2$ contains a copy of $\Real^2$, it suffices to prove the result for $\Real^2(2)$. Fix a point $c\in\Real^2$ (say $c:=(0,0)$). Let $x=\{x_1,x_2\}$, $y=\{y_1,y_2\}$, $z=\{z_1,z_2\}$ be points in $\Real^2(2)\backslash\Real^2(1)$ such that the elements of the set $x\cup y\cup z=\{x_1,x_2,y_1,y_2,z_1,z_2\}$ are uniformly distributed on the unit circle $\del\ol{B}_1(c)\subset\Real^2$ in the order $x_1,y_1,z_1,x_2,y_2,z_2$ counterclockwise, as in Figure \ref{dg2}. Thus, the points of $x\cup y\cup z$ are the vertices of a regular hexagon, of side length $1$, in $\Real^2$ centered at $c$. We have
\bea
&&d(x_1,y_1)=d(y_1,z_1)=d(z_1,x_2)=d(x_2,y_2)=d(y_2,z_2)=d(z_2,x_1)=1,\nn\\
&&\textstyle c={x_1+x_2\over 2}={y_1+y_2\over 2}={z_1+z_2\over 2}.\nn
\eea

\begin{figure}[H]
\centering
\scalebox{1}{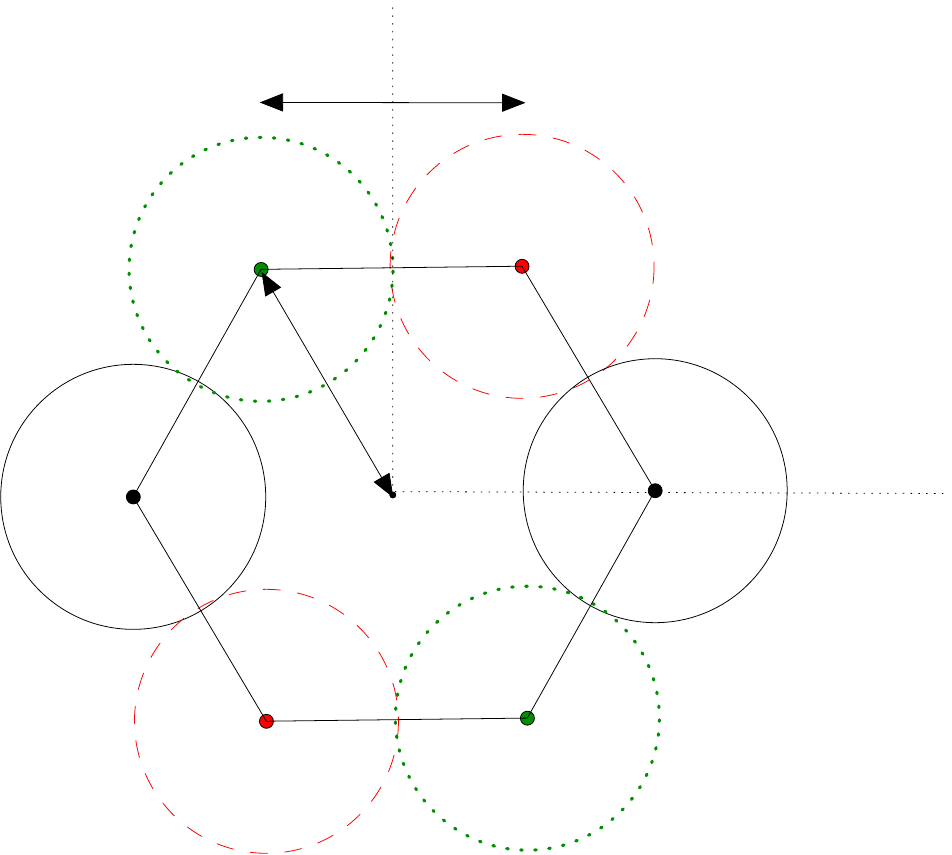} 
 \caption{~Pairwise intersecting balls in $\Real^2(2)$ with no common point.}\label{dg2}
\end{figure}

The collection of balls $\left\{\ol{B}_{1\over 2}(x),\ol{B}_{1\over 2}(y),\ol{B}_{1\over 2}(z)\right\}$ is pairwise intersecting in $\Real^2(2)$, and
\bea
&&\ol{B}_{1\over 2}(x)=\left\{u=\{u_1,u_2\}:u_1\in \ol{B}_{1\over 2}(x_1),~u_2\in \ol{B}_{1\over 2}(x_2)\right\},\nn\\
&&\ol{B}_{1\over 2}(y)=\left\{v=\{v_1,v_2\}:v_1\in \ol{B}_{1\over 2}(y_1),~v_2\in \ol{B}_{1\over 2}(y_2)\right\},\nn\\
&&\ol{B}_{1\over 2}(z)=\left\{w=\{w_1,w_2\}:w_1\in \ol{B}_{1\over 2}(z_1),~w_2\in \ol{B}_{1\over 2}(z_2)\right\},\nn
\eea
where the equalities follow by Corollary \ref{HausDistBound1}.

Suppose $e=\{e_1,e_2\}\in \ol{B}_{1\over 2}(x)\cap\ol{B}_{1\over 2}(y)\cap\ol{B}_{1\over 2}(z)$. Then by the properties of Hausdorff distance, $e_1$ (as well as $e_2$) must lie in the intersection of three balls from the collection of six balls {\small$\C=\left\{\ol{B}_{1\over 2}(c'):c'\in x\cup y\cup z\right\}$} in $\Real^2$. Since no three balls from the collection $\C$ have a common point in $\Real^2$, we have a contradiction.
\end{proof}

\begin{dfn}[\textcolor{blue}{Regular circular polygon}]
Points $x_1,...,x_n\in\Real^2$ form a regular circular polygon if they are uniformly placed (in that order) on a circle $\del\ol{B}_r(c)$ in $\Real^2$, where the radius $r$ is chosen such that~ $d(x_1,x_2)=d(x_2,x_3)=\cdots=d(x_{n-1},x_n)=d(x_n,x_1)=1$.
\end{dfn}

\begin{lmm}[\textcolor{OliveGreen}{Failure of BIP for $X(n)$}]\label{NoBIPLmm2}
Let $X$ is a normed space such that $\dim X\geq 2$. Then for all $n\geq 2$, $X(n)$ does not have the BIP.
\end{lmm}
\begin{proof}
Since every normed space of dimension $\geq 2$ contains a copy of $\Real^2$, it suffices to prove the result for $\Real^2(n)$. Let $x=\{x_1,...,x_n\}$, $y=\{y_1,...,y_n\}$, $z=\{z_1,...,z_n\}$ be points in $\Real^2(n)\backslash\Real^2(n-1)$ such that the following conditions hold.
\bit[leftmargin=0.7cm]
\item The sets $A_1=\{x_1,x_2,y_1,y_2,z_1,z_2\}$, $A_j=\{x_{j+1},y_{j+1},z_{j+1}\}$, $2\leq j\leq n-1$, are distantly separated in $\Real^2$.
\item For each $i$ the elements of $A_i$ lie on a circle and form a regular polygon in $\Real^2$, where the elements of $A_1$ in particular do so in the order $x_1,y_1,z_1,x_2,y_2,z_2$ counterclockwise.
\eit

\begin{figure}[H]
\centering
\scalebox{1.0}{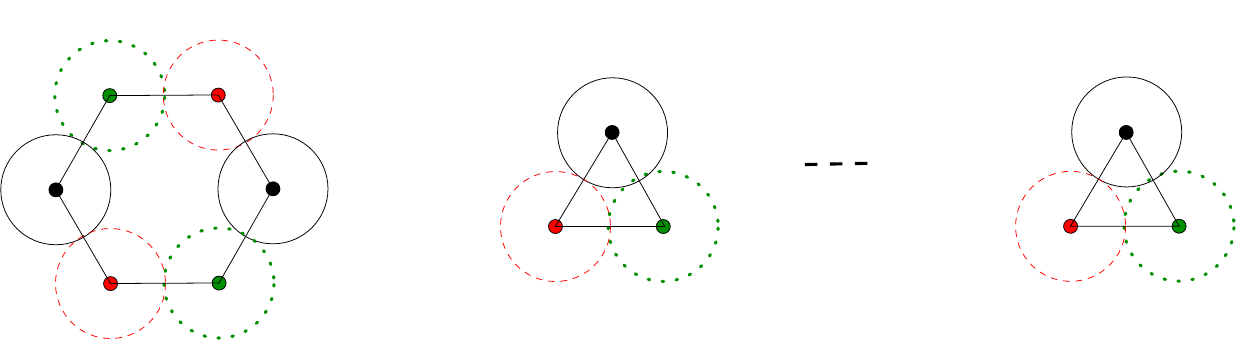} 
\caption{~Pairwise intersecting balls in $\Real^2(n)$ with no common point.}\label{dg3}
\end{figure}

Then the collection of balls {\small$\left\{\ol{B}_{1\over 2}(x),\ol{B}_{1\over 2}(y),\ol{B}_{1\over 2}(z)\right\}$} is pairwise intersecting in $\Real^2(n)$, since
\bea
&&\ol{B}_{1\over 2}(x)=\left\{u=\{u_1,...,u_n\}\in\Real^2(n):u_i\in\ol{B}_{1\over 2}(x_i),~i=1,...,n\right\},\nn\\
&&\ol{B}_{1\over 2}(y)=\left\{v=\{v_1,...,v_n\}\in\Real^2(n):v_i\in\ol{B}_{1\over 2}(y_i),~i=1,...,n\right\},\nn\\
&&\ol{B}_{1\over 2}(z)=\left\{w=\{w_1,...,w_n\}\in\Real^2(n):w_i\in\ol{B}_{1\over 2}(z_i),~i=1,...,n\right\},\nn
\eea
where the equalities follow by Corollary \ref{HausDistBound1}.

Suppose $e=\{e_1,...,e_n\}\in \ol{B}_{1\over 2}(x)\cap\ol{B}_{1\over 2}(y)\cap\ol{B}_{1\over 2}(z)$. Then by the properties of Hausdorff distance, at least one element of $e$ must lie in the intersection of three balls from the collection of six balls {\small$\C=\left\{\ol{B}_{1\over 2}(c):c\in\{x_1,x_2,y_1,y_2,z_3,y_3\}\right\}$} in $\Real^2$. Since no three balls from the collection $\C$ have a common point in $\Real^2$, we have a contradiction.
\end{proof}

\begin{rmk}[\textcolor{OliveGreen}{Non 1-ALR's}]
Let $X$ be a normed space such that $\dim X\geq 2$. Then for $n\geq 2$, $X(n)$ is not a 1-ALR. (This follows from Lemma \ref{NoBIPLmm2} and Proposition \ref{1ALRcrit}.)
\end{rmk}

\begin{dfn}[\textcolor{blue}{\index{BIP collection}{BIP collection}}]
Let $X$ be a metric space. A pairwise intersecting collection of balls $\left\{\ol{B}_{r_\al}(x_\al):\al\in\A\right\}$ in $X$ is a BIP collection if ~$\bigcap\ol{B}_{r_\al}(x_\al)\neq\emptyset$.
\end{dfn}
Recall that $X$ has the BIP $\iff$ every pairwise intersecting collection of balls in $X$ is a BIP collection.

\begin{dfn}[\textcolor{blue}{Compatible set in $X(n)$}]
Let $X$ be a metric space. A set $C=\{x_\al:\al\in\A\}\subset X(n)$ is compatible if there exist relations $\vphi_\beta^\al:x_\al\ra x_\beta$, $(x_\al)_i\mapsto (x_\beta)_{\sigma_\beta^\al(i)}$ (for permutations $\sigma_\beta^\al\in S_n$) such that
\bea
\vphi_\beta^\al=\vphi_\beta^\rho \vphi_\rho^\al:x_\al\sr{\vphi_\rho^\al}{\ral}x_\rho\sr{\vphi_\beta^\rho}{\ral}x_\beta,~~ \vphi_\al^\al=id,~~\txt{for all}~~\al,\beta,\rho\in \A.\nn
\eea
\end{dfn}

\begin{dfn}[\textcolor{blue}{Compatible collection of balls in $X(n)$}]
Let $X$ be a metric space. A collection of balls $\left\{\ol{B}_{r_\al}(x_\al):\al\in\A\right\}$ in $X(n)$ is compatible if the set of centers $\{x_\al:\al\in\A\}$ is compatible in $X(n)$.
\end{dfn}

\begin{dfn}[\textcolor{blue}{\index{CIP collection}{CPI collection} of balls in $X(n)$}]
Let $X$ be a metric space. A collection of balls $\left\{\ol{B}_{r_\al}(x_\al):\al\in\A\right\}$ in $X(n)$ is compatibly pairwise intersecting (or a CPI collection) if there exist relations
\bea
\label{CPI0}\vphi_\beta^\al:x_\al\ra x_\beta,~~(x_\al)_i\mapsto (x_\beta)_{\sigma_\beta^\al(i)},~~~~\sigma_\beta^\al\in S_n,
\eea
satisfying the following conditions.
\bea
\label{CPI1}&&d\left((x_\al)_i,(x_\beta)_{\sigma_\beta^\al(i)}\right)\leq d_H(x_\al,x_\beta),~~~~\txt{for all}~~~~i=1,...,n.\\
\label{CPI2}&& \vphi_\beta^\al=\vphi_\beta^\rho \vphi_\rho^\al:x_\al\sr{\vphi_\rho^\al}{\ral}x_\rho\sr{\vphi_\beta^\rho}{\ral}x_\beta,~~~~\vphi_\al^\al=id,~~~~\txt{for all}~~~~\al,\beta,\rho\in \A.
\eea
\end{dfn}
By Note \ref{BndRelNote}, (\ref{CPI0}) and (\ref{CPI1}) together imply the collection of balls $\{\ol{B}_{r_\al}(x_\al):\al\in\A\}$ is such that any two centers $x_\al,x_\beta$ admit a complete relation $R_{\al\beta}\subset x_\al\times x_\beta$ satisfying
\bea
|R_{\al\beta}|\leq n~~~~\txt{and}~~~~\max_{(a,b)\in R_{\al\beta}}\|a-b\|\leq d_H(x_\al,x_\beta).\nn
\eea

\begin{prp}\label{BIPexist}
If a Banach space $X$ has the BIP, then every CPI collection of balls in $X(n)$ is a BIP collection.
\end{prp}
\begin{proof}
We know that $X$ is both complete and geodesic, and so metrically convex. Thus, for any balls $\ol{B}_r(a),\ol{B}_{r'}(a')$ in $X$, we have $\ol{B}_r(a)\cap \ol{B}_{r'}(a')\neq\emptyset$ $\iff$ $\|a-a'\|\leq r+r'$.

Consider a CPI collection of balls $\left\{\ol{B}_{r_\al}(x_\al):\al\in\A\right\}$ in $X(n)$. By hypotheses, there exist permutations $\sigma_\beta^\al\in S_n$ giving bijections $\vphi_\beta^\al:x_\al\ra x_\beta$, $(x_\al)_i\mapsto (x_\beta)_{\sigma_\beta^\al(i)}$ satisfying (\ref{CPI1}) and (\ref{CPI2}). Thus,
\bea
\max_i\left\|(x_\al)_i-(x_\beta)_{\sigma_\beta^\al(i)}\right\|=d_H(x_\al,x_\beta)\leq d_H(x_\al,z)+d_H(z,x_\beta)\leq r_\al+r_\beta,\nn
\eea
for $z\in \ol{B}_{r_\al}(x_\al)\cap \ol{B}_{r_\beta}(x_\beta)$. Since $X$ is metrically convex, the above expression shows that for a fixed $\al_0\in \A$, each of the collections $\left\{\ol{B}_{r_\al}\left((x_\al)_{\sigma_\al^{\al_0}(i)}\right)\right\}_{\al\in \A}$, $i=1,...,n$,
intersect pairwise in $X$ and so have a common point of intersection $u_i\in\bigcap_{\al\in\A}\ol{B}_{r_\al}\left((x_\al)_{\sigma_\al^{\al_0}(i)}\right)$, i.e., $\big\|u_i-(x_\al)_{\sigma_\al^{\al_0}(i)}\big\|\leq r_\al$, $i=1,...,n$. The following computation shows that the point $u=\{u_i\}\in X(k)$ lies in the intersection $\bigcap_{\al\in\A}\ol{B}_{r_\al}(x_\al)$.
\begin{align}
 d_H(u,x_\al)&\textstyle=\max\left\{\max_i\min_j\|u_i-(x_\al)_j\|,\max_j\min_i\|u_i-(x_\al)_j\|\right\}\nn\\
&\textstyle=\max\left\{\max_i\min_j\left\|u_{(\sigma_\al^{\al_0})^{-1}(i)}-(x_\al)_j\right\|,\max_j\min_i\left\|u_{(\sigma_\al^{\al_0})^{-1}(i)}-(x_\al)_j\right\|\right\}\nn\\
&\textstyle\leq \max_i\left\|u_{(\sigma_\al^{\al_0})^{-1}(i)}-(x_\al)_i\right\|=\max_i\left\|u_i-(x_\al)_{\sigma_\al^{\al_0}(i)}\right\|\leq r_\al.\nn
\end{align}
\end{proof}


\chapter{Topological and Lipschitz $k$-Connectedness}\label{LipConn}
This chapter relies on notation/terminology from section \ref{PrelimsAT}. With the help of Lipschitz homotopy and Lipschitz contractibility (discussed in sections \ref{LipConnLH1}, \ref{LipConnLC1}), we discuss Lipschitz $k$- connectedness (in section \ref{LipConnLH2}, \ref{LipConnLC2}). A notion of conical Lipschitz contractibility introduced in section \ref{LipConnLC1} is used in section \ref{LipConnLC2} to prove in Theorem \ref{LipConSuff5} that finite subset spaces of a normed space are Lipschitz $k$-connected for all $k\geq 0$. This result will be used later (in section \ref{FSRP2lr}, Theorem \ref{FinDimLipRet}) to prove that if $X$ is a finite-dimensional normed space, then we have Lipschitz retractions $X(n)\ra X(n-1)$. Throughout, $I$ denotes the interval $[0,1]\subset\Real$.

\section{Lipschitz homotopy and contractibility}\label{LipConnLH1}
\begin{dfn}[\textcolor{blue}{\index{Lipschitz! homotopy}{Lipschitz homotopy}}]
Let $X,Y$ be metric spaces. Two continuous maps $f,g:X\ra Y$ are Lipschitz homotopic (written $f\simeq_Lg$) if there exists a Lipschitz homotopy $H:X\times I\ra Y$ between $f$ and $g$. (Note that if $f\simeq_Lg$, then $f$ and $g$ are necessarily Lipschitz, since a restriction of a Lipschitz map is Lipschitz. Also, recall that the metric on $X\times I$ is given by $d((x,t),(x',t'))=\max(d(x,x'),|t-t'|)$.)
\end{dfn}

By \cite[Theorem 1]{ferry-weinbg2013}, there exist nontrivial metric spaces $X,Y$ for which any two Lipschitz maps $f,g:X\ra Y$ can be connected by a Lipschitz homotopy. Also, as in \cite[Theorems 2,3]{ferry-weinbg2013}, we can find metric spaces $X,Y$ and Lipschitz maps $f,g:X\ra Y$ such that no Lipschitz homotopy $X\times I\ra Y$ exists between $f,g$, but $f,g$ can be connected by a homotopy $H:X\times I\ra Y$ with Lipschitz components $H_t:=H|_{X\times\{t\}}:X\ra Y$ for all $0\leq t\leq 1$.

\begin{example}
Let $X=\{u\}$ be a one-point metric space. Let $Y=[0,1]_p$ be the interval $[0,1]$ with the metric $d_p(t,t')=|t-t'|^{1\over p}$, $p>1$. Then the maps $f,g:X\ra Y$ given by $f(u)=0$, $g(u)=1$ are connected by the homotopy $H:X\times[0,1]\ra Y$ given by $H(u,t)=t$ (with Lipschitz components $H_t(u)=H(u,t)=t$), but there exists no Lipschitz homotopy $H:X\times[0,1]\ra Y$ between $f$ and $g$, since $Y$ is a snowflake.
\end{example}
\begin{example}
A normed space $X$ is contractible to any point $x_0\in X$ through the homotopy $H:X\times I\ra X$, $H(x,t)=(1-t)x+tx_0$, with Lipschitz components $H_t(x)=H(x,t)$. However, by Lemma \ref{LipConBound} below, a normed space (being unbounded) is not Lipschitz contractible.
\end{example}

\begin{dfn}[\textcolor{blue}{\index{Lipschitz! homotopy equivalence}{Lipschitz homotopy equivalence}, Lipschitz homotopy inverse}]
Let $X,Y$ be metric spaces. $h:X\ra Y$ is a Lipschitz homotopy equivalence (making $X,Y$ Lipschitz homotopy equivalent, written $X\simeq_LY$) if $f$ is Lipschitz and there exists a Lipschitz map $\wt{h}:Y\ra X$ (called Lipschitz homotopy inverse of $h$) such that $h\circ\wt{h}\simeq_Lid_Y$ and $\wt{h}\circ h\simeq_Lid_X$.
\end{dfn}

\begin{dfn}[\textcolor{blue}{\index{Lipschitz! contractibility}{Lipschitz contractibility}}]
A metric space $X$ is $\ld$-Lipschitz contractible if it is $\ld$-Lipschitz homotopy equivalent to a point. Equivalently, $X$ is $\ld$-Liscphitz contractible if the identity $id_X:X\ra X$ is $\ld$-Lipschitz homotopic to a constant map $X\ra X$.
\end{dfn}

Note that by Corollary \ref{SnwfNoLipCon}, a snowflake is not Lipschitz contractible. Thus, the contractible space $X=L\cup K$ (with $L$ a line and $K$ a snowflake) constructed in the proof of Proposition \ref{NoLipRetCon} is not Lipschitz contractible.

\begin{lmm}[\textcolor{OliveGreen}{Lipschitz contractibility bound}]\label{LipConBound}
A Lipschitz contractible metric space is bounded.
\end{lmm}
\begin{proof}
Let $X$ be a metric space and $H:X\times I\ra X$ a $\ld$-Lipschitz homotopy. Then
{\small\begin{align}
d\big(H(x,t),H(x',t')\big)\leq \ld\max\left\{d(x,x'),|t-t'|\right\},~~H|_{X\times\{0\}}=\txt{const.}=x_0\in X,~~H|_{X\times\{1\}}=id_X.\nn
\end{align}}
Therefore, for all $x,x'\in X$, we have
\bea
d(x,x')\leq d(x,x_0)+d(x_0,x')=d\big(H(x,1),H(x,0)\big)+d\big(H(x',0),H(x',1)\big)\leq 2\ld.\qedhere\nn
\eea
\end{proof}
This result shows that a weaker notion of Lipschitz contractibility, e.g., \emph{local Lipschitz contractibility}, would be more useful. The notion of \emph{bounded conical Lipschitz contractibility} in section \ref{LipConnLC2} below is closely related to the idea of local Lipschitz contractibility. Therefore, we make the following definition (by which a normed space in particular, although not Lipschitz-contractible, is locally Lipschitz contractible).

\begin{dfn}[\textcolor{blue}{\index{LL contractibility}{Local Lipschitz (LL) contractibility}}]
A metric space $X$ is locally $\ld$-Lipschitz contractible if every point $x\in X$ has a $\ld$-Lipschitz contractible neighborhood.
\end{dfn}

\begin{lmm}\label{RigRugLmm1}
Let $X$ be a (locally) Lipschitz contractible metric space (via a Lipschitz homotopy $H$). Then every Lipschitz retract $A\subset X$ is (locally) Lipschitz contractible, as shown in the diagram below (in which $r:X\ra A$ is a Lipschitz retraction).

\bc\bt
X\times I\ar[rrr,"H"] &&& X\ar[d,"r"]\\
A\times I\ar[u,hook,"i"']\ar[rrr,dashed,"G:=r\circ H\circ i"] &&& A
\et\ec
\end{lmm}

\begin{prp}\label{RigRugPrp}
If $A,B$ are metric spaces, then $A\times B$ is Lipschitz contractible $\iff$ $A$ and $B$ are Lipschitz contractible.
\end{prp}
\begin{proof}
If $A\times B$ is Lipschitz contractible, then $A,B\subset A\times B$ are Lipschitz contractible as Lipschitz retracts. Conversely, assume $A,B$ are Lipschitz contractible through $H_A:A\times I\ra A$, $H_B:B\times I\ra B$. To show $A\times B$ is also Lipschitz contractible, consider the map $H:A\times B\times I\ra A\times B$ given by  $H\big((a,b),t\big):=\big(H_A(a,t),H_B(b,t)\big)$. Then
\begin{align*}
d\Big(H&\big((a,b),t\big),H\big((a',b'),t'\big)\Big)=d\Big(\big(H_A(a,t),H_B(b,t)\big),\big(H_A(a',t'),H_B(b',t')\big)\Big)\nn\\
&=\max\Big(d\big(H_A(a,t),H_A(a',t')\big),d\big(H_B(b,t),H_B(b',t')\big)\Big)\nn\\
&\leq K\max\Big(d\big((a,t),(a',t')\big),d\big((b,t),(b',t')\big)\Big),~~~~~~~~K:=\max(\Lip H_A,\Lip H_B),\nn\\
&\leq K\max\Big(\max\big(d(a,a'),|t-t'|\big),\max\big(d(b,b'),|t-t'|\big)\Big)\nn\\
&\leq K\max\Big(\max\big(d(a,a'),d(b,b')\big),|t-t'|\Big)\leq K\max\Big(d\big((a,a'),(b,b')\big),|t-t'|\Big)\nn\\
&=K d\Big(\big((a,b),t\big),\big((a',b'),t'\big)\Big),\nn\qedhere
\end{align*}
\end{proof}

\begin{lmm}\label{RigRugLmm2}
Let $A,B$ be metric spaces and $X:=A\times B$. If $A$ contains no rectifiable curves, or $B$ contains no rectifiable curves, then $X$ is not Lipschitz contractible.
\end{lmm}
\begin{proof}
Without loss of generality, let $B$ contain no rectifiable curves. Suppose $A\times B$ is Lipschitz contractible, so that we have a Lipschitz homotopy $H:A\times B\times I\ra A\times B$. Since $H$ has the general form $H((a,b),t)=\big(H_1(a,b,t),H_2(a,b,t)\big)$, where $H_1:A\times B\times I\ra A$ and $H_2:A\times B\times I\ra B$, we have
\bea
&&d\big(H(a,b,t),H(a',b',t')\big)=\max\Big[d\big(H_1(a,b,t),H_1(a',b',t')\big),d\big(H_2(a,b,t),H_2(a',b',t')\big)\Big]\nn\\
&&~~~~\leq \Lip(H)\max\Big[\max\big(d(a,a'),d(b,b')\big),|t-t'|\Big],\nn\\
&&~~\Ra~~d\big(H_2(a,b,t),H_2(a,b,t')\big)\leq \Lip(H)|t-t'|,~~~~\txt{for any fixed}~~(a,b)\in A\times B,\nn
\eea
which is a contradiction since $B$ contains no rectifiable curves.
\end{proof}

\section{Homotopy equivalence and $k$-connectedness for finite subset spaces}\label{LipConnLH2}
\begin{lmm}[\textcolor{OliveGreen}{Lipschitz homotopy equivalence of finite subset spaces}]\label{LipHEfssRecall}
Let $X,Y$ be metric spaces. Then $X\simeq_L Y$ $\iff$ $X(n)\simeq_LY(n)$ for all $n\geq 1$. (Similarly, for any spaces $X,Y$ we have, topologically, $X\simeq Y$ $\iff$ $X(n)\simeq Y(n)$ for all $n\geq 1$.)
\end{lmm}
\begin{proof}
Let $f:X\ra Y$ be a Lipschitz homotopy equivalence with Lipschitz homotopy inverse $g:Y\ra X$, i.e., there exist Lipschitz homotopies {\small $H:X\times I\ra X$}, {\small $H':Y\times I\ra Y$} with
\bea
H(x,0)=g\circ f(x),~~H(x,1)=x,~~~~~~~~H'(y,0)=f\circ g(y),~~H'(y,1)=y.\nn
\eea
Then we get an induced Lipschitz homotopy equivalence $F:X(n)\ra Y(n)$, with Lipschitz homotopy inverse $G:Y(n)\ra X(n)$, given by
\bea
F\big(\{x_1,...,x_n\}\big):=\big\{f(x_1),...,f(x_n)\big\},~~~~G\big(\{y_1,...,y_n\}\big):=\big\{g(y_1),...,g(y_n)\big\}.\nn
\eea
Note that we have Lipschitz homotopies $G\circ F\sr{\wt{H}}{\simeq}id_{X(n)}$ and $F\circ G\sr{\wt{H}'}{\simeq}id_{Y(n)}$ given by
\bea
&&\wt{H}:X(n)\times I\ra X(n),~~~~\wt{H}\big(\{x_1,...,x_n\},t\big):=\big\{H(x_1,t),\cdots, H(x_n,t)\big\},\nn\\
&&\wt{H}':Y(n)\times I\ra Y(n),~~~~\wt{H}'\big(\{y_1,...,y_n\},t\big):=\big\{H'(y_1,t),\cdots, H'(y_n,t)\big\}.\nn
\eea
Indeed, we have the required relations
\bea
&&\wt{H}(\{x_1,...,x_n\},0)=\{H(x_1,0),...,H(x_n,0)\}=\{gf(x_1),...,gf(x_n)\}=GF(\{x_1,...,x_n\}),\nn\\
&&\wt{H}(\{x_1,...,x_n\},1)=\{H(x_1,1),...,H(x_n,1)\}=\{x_1,...,x_n\},\nn
\eea
and
\bea
&&\wt{H}'(\{y_1,...,y_n\},0)=\{H'(y_1,0),...,H'(y_n,0)\}=\{fg(y_1),...,fg(y_n)\}=FG(\{y_1,...,y_n\})\nn\\
&& \wt{H}'(\{y_1,...,y_n\},1)=\{H'(y_1,1),...,H'(y_n,1)\}=\{y_1,...,y_n\}.\nn\qedhere
\eea
\end{proof}

\begin{question}
Given a Lipschitz retraction $r_n:X(n)\ra X(n-1)$, we get the diagram
\bea\bt
X(n)\ar[rrrr,"r_n"] &&&& X(n-1)\ar[d,"F"]\\
Y(n)\ar[u,"G"]\ar[rrrr,dashed,"s_n~=~F\circ r_n\circ G"] &&&& Y(n-1)
\et\nn\eea
in which the induced Lipschitz map $s_n:Y(n)\ra Y(n-1)$ satisfies
\bea
s_n|_{Y(n-1)}=F\circ G|_{Y(n-1)}\simeq_Lid_{Y(n-1)}.\nn
\eea
This shows that the map $H_k(Y(n-1))\sr{i_\ast}{\ral}H_k(Y(n))$ induced by the inclusion $Y(n-1)\sr{i}{\hookrightarrow}Y(n)$ has a left inverse $H_k(Y(n))\sr{(s_n)_\ast}{\ral}H_k(Y(n-1))$. Can we conclude that the Lipschitz map $s_n:Y(n)\ra Y(n-1)$ is Lipschitz homotopic to a Lipschitz retraction $\wt{r}_n:Y(n)\ra Y(n-1)$?
\end{question}

\begin{dfn*}[\textcolor{blue}{Recall: Connected space}]
A space $X$ is connected if it cannot be written as a disjoint union of two nonempty open sets.
\end{dfn*}

\begin{question}\label{ConQsn}
We know that if $X$ is contractible, then so is $X(n)$ by Lemma \ref{LipHEfssRecall}. Is it also true that if $X$ is connected, then so is $X(n)$? \emph{The answer is yes}. Indeed, if $X(n)$ is not connected, then neither is $X^n$ (since $X(n)$ is a continuous image of $X^n$), and hence neither is $X$ (since $X$ connected implies $X^n$ connected [footnote\footnote{A product of connected spaces $X\times Y$ is connected: Indeed, suppose $X\times Y=U_1\cup U_2$ for nonempty disjoint open sets $U_1,U_2$. Then for any $y\in Y$, if the connected subspace $X\times\{y\}$ lies in $U_1$, then by the connectedness of $Y$, every subset of the form $\{x\times Y\}$, $x\in X$, (and hence all of $X\times Y$) also lies in $U_1$ (a contradiction).}]).
\end{question}

\begin{dfn}[\textcolor{blue}{\index{Topological! $k$-connectedness}{Topological $k$-connectedness}, Path-connected space, Simply-connected space}]\label{TopConnDfn}
A space $X$ is (topologically) \ul{$k$-connected} if for each $0\leq l\leq k$, every continuous map $f:S^l=\del B^{l+1}\ra X$ extends to a continuous map $F:B^{l+1}\ra X$. Here, $B^{l+1}\subset\Real^{l+1}$ is the unit ball. In this case, we say $X$ is (topologically) $k$-connected. A $0$-connected space is called a \ul{path-connected} space, while a $1$-connected space is called a \ul{simply-connected} space.
\end{dfn}

Criteria for topological $k$-connectedness were given in Lemma \ref{NoHoleCrit}. If $X$ is $k$-connected for $k\geq 1$, it does not follow that $X(n)$ is also $k$-connected (see the discussion in Question \ref{XXnLipCon} below). However, if $X$ is $0$-connected, then so is $X(n)$ by the following.

\begin{rmk}[\textcolor{OliveGreen}{Path-connectedness of $X(n)$}]\label{PathConRmk}
If $X$ is path-connected, so is $X(n)$. Indeed, if $X$ is path-connected, then so is $X^n$. Thus, given $x,y\in X(n)$, any path $\gamma:[0,1]\ra X^n$ from a point in $q^{-1}(x)$ to a point in $q^{-1}(y)$ gives a path $q\circ\gamma:[0,1]\sr{\gamma}{\ral}X^n\sr{q}{\ral}X(n)$ from $x$ to $y$.
\end{rmk}

\section{Lipschitz contractibility over a cone}\label{LipConnLC1}
Recall that if $X$ is a space, the \emph{cylinder} of $X$ is the product space ~$C(X):=X\times I$, and the \emph{cone} of $X$ is the quotient space $\Cone(X):={C(X)\over X\times\{0\}}={X\times I\over X\times\{0\}}$, where with the quotient map $q:X\times I\ra \Cone(X)$, $(x,t)\mapsto q(x,t)\eqv[(x,t)]$, the \emph{vertex} of the cone is the point $q(x,0)=X\times\{0\}\in \Cone(X)$ and its \emph{base} is the subset $q(X\times\{1\})\cong X\times\{1\}\cong X$. Note that $\Cone(X)$ is contractible through the homotopy $H:\Cone(X)\times I\ra \Cone(X)$ given by $H\big(q(x,s),t\big):=q(x,st)$.

\begin{lmm}\label{TunnelLmm}
Let $(X,d)$ be a metric space and $\vphi:[0,\infty)\ra[0,\infty)$ a concave nondecreasing function such that $\vphi(0)\geq0$. Then (i) $d_\vphi:=\vphi\circ d$ is a metric on $X$, and (ii) the Hausdorff metric $d_H$ on $X(n)$ satisfies ~$(d_H)_\vphi=(d_\vphi)_H$.
\end{lmm}
(Recall that $\vphi$ is concave iff ~$\al\vphi(t)+(1-\al)\vphi(t')\leq\vphi(\al t+(1-\al)t')$ for all $0\leq\al\leq 1$.)
\begin{proof}
(i) To show $d_\vphi$ is a metric, it suffices to prove the triangle inequality. Observe that for $0\leq \al\leq 1$, we have $\al\vphi(t)\leq \al\vphi(t)+(1-\al)\vphi(0)\leq\vphi(\al t+(1-\al)0)=\vphi(\al t)$,
\bea
\textstyle~~\Ra~~\vphi(t+t')={t\over t+t'}\vphi(t+t')+{t'\over t+t'}\vphi(t+t')\leq\vphi(t)+\vphi(t'),~~~~\txt{for all}~~t,t'\in[0,\infty).\nn
\eea
Given $x,y,z\in X$, let $d_1=d(x,y)$, $d_2=d(x,z)$, $d_3=d(y,z)$. Then because $d_1\leq d_2+d_3$ and $\vphi$ is nondecreasing, we have $\vphi(d_1)\leq\vphi(d_2+d_3)\leq\vphi(d_2)+\vphi(d_3)$, which shows $d_\vphi$ is a metric on $X$. (ii) Therefore,
\begin{align}
&(d_H)_\vphi(\{x_1,...,x_n\},\{y_1,...,y_n\})=\vphi\Big[d_H(\{x_1,...,x_n\},\{y_1,...,y_n\})\Big]\nn\\
&~~~~=\vphi\left[\max\left(\max_i\min_jd(x_i,y_j),\max_j\min_id(x_i,y_j)\right)\right]\nn\\
&~~~~\sr{(s)}{=}\max\left(\max_i\min_j\vphi\left[d(x_i,y_j)\right],\max_j\min_i\vphi\left[d(x_i,y_j)\right]\right)\nn\\
&~~~~=(d_\vphi)_H(\{x_1,...,x_n\},\{y_1,...,y_n\}),\nn
\end{align}
where step (s) holds because $\vphi$ is nondecreasing.
\end{proof}

\begin{dfn}[\textcolor{blue}{\index{Cone! metrics}{Cone metrics} on the cylinder/cone of a metric space}]\label{ConeMetrics1}
Let $(X,d)$ be a metric space. Then by Lemma \ref{TunnelLmm}(i), we know $d_{b1}(x,x'):=\min\big(\pi,d(x,x')\big)$ and $d_{b2}(x,x'):=\min\big(2,d(x,x')\big)$ are metrics on $X$.

Following \cite[Definition 3.6.12, p.91]{BBI}, a metric on $\Cone(X)$ is given by
\bea
&&\textstyle d_{c1}\big((x,t),(x',t')\big):=\sqrt{t^2+t'^2-2tt'\cos d_{b1}(x,x')}=\sqrt{|t-t'|^2+2tt'\big[1-\cos d_{b1}(x,x')\big]}\nn\\
&&\textstyle~~~~=\sqrt{|t-t'|^2+4tt'\sin^2{d_{b1}(x,x')\over 2}}~\sr{\sin|\theta|\leq|\theta|}{\leq}~ \sqrt{|t-t'|^2+tt'd_{b1}(x,x')^2}\nn\\
\label{Dc1Ineq}&&\textstyle~~~~\leq\sqrt{(1+tt')}\max\big(|t-t'|,d_{b1}(x,x')\big).
\eea
As shown in \cite[Lemma 2.2]{kovalev2015}, the above metric is biLipschitz equivalent to the metric
\bea
\label{Dc2Ineq}d_{c2}\big((x,t),(x',t')\big):=|t-t'|+\min(t,t')d_{b2}(x,x'),
\eea
where the bi-Lipschitz equivalence of $d_{c1}$ and $d_{c2}$ (according to \cite[Lemma 2.2]{kovalev2015}) is given by
\bea
\label{ConeBiLipEq}d_{c2}/12\leq d_{c1}\leq 10d_{c2}.
\eea
\end{dfn}

\begin{lmm}[\textcolor{OliveGreen}{\cite[Proposition 3.6.13, p.91]{BBI}}]
Let $(X,d)$ be a metric space. The map $d_{c1}:\Cone(X)\times \Cone(X)\ra\Real$ is a metric.
\end{lmm}
\begin{proof}
By construction, it suffices to assume $\diam X\leq\pi$, in which case,
\bea
\textstyle d_{c1}\big((x,t),(x',t')\big)=\sqrt{t^2+t'^2-2tt'\cos d(x,x')}.\nn
\eea
Nonnegativity and symmetry are clear. So we will prove the triangle inequality. Let $y=(x,t)$, $y'=(x',t')$, $y''=(x'',t'')$ be points in $X\times I$. Let $\ol{y}$, $\ol{y}'$, $\ol{y}''$ be points in $\Real^2$ such that
\bea
&&\|0-\bar{y}\|=\|\ol{y}\|=t,~~\|0-\bar{y}'\|=\|\ol{y}'\|=t',~~\|0-\bar{y}''\|=\|\ol{y}''\|=t'',\nn\\
&&\theta_{\ol{y},0,\ol{y}'}=d(x,x'),~~\theta_{\ol{y}',0,\ol{y}''}=d(x',x''),\nn
\eea
where for any vectors $a,b,c\in\Real^2$, $\theta_{a,b,c}$ is the angle between $a$ and $c$ at $b$. Then
\bea
&&\textstyle\|\ol{y}-\ol{y}'\|~=~\sqrt{\|\ol{y}\|^2+\|\ol{y}'\|^2-2\|\ol{y}\|\|\ol{y}'\|\cos\theta_{\ol{y},0,\ol{y}'}}~=~\sqrt{t^2+t'{}^2-2tt'\cos d(x,x')}\nn\\
&&\textstyle~~~~=d_{c1}\big((x,t),(x',t')\big),~~~~\txt{and}~~~~\|\ol{y}'-\ol{y}''\|~\sr{\txt{similarly}}{=}~d_{c1}\big((x',t'),(x'',t'')\big).\nn
\eea
Thus we have the following two cases.
\bit[leftmargin=0.4cm]
\item $\theta_{\ol{y},0,\ol{y}'}+\theta_{\ol{y}',0,\ol{y}''}\leq\pi$: In this case, $\theta_{\ol{y},0,\ol{y}''}=\theta_{\ol{y},0,\ol{y}'}+\theta_{\ol{y}',0,\ol{y}''}$, and so
\bea
&&\textstyle\|\ol{y}-\ol{y}''\|=\sqrt{\|\ol{y}\|^2+\|\ol{y}''\|^2-2\|\ol{y}\|\|\ol{y}''\|\cos\theta_{\ol{y},0,\ol{y}''}}\nn\\
&&\textstyle~~~~=\sqrt{\|\ol{y}\|^2+\|\ol{y}''\|^2-2\|\ol{y}\|\|\ol{y}''\|\cos\big(d(x,x')+d(x',x'')\big)}\nn\\
&&\textstyle~~~~\geq \sqrt{\|\ol{y}\|^2+\|\ol{y}''\|^2-2\|\ol{y}\|\|\ol{y}''\|\cos\big(d(x,x'')\big)}=d_{c1}\big((x,t),(x'',t'')\big),\nn
\eea
which implies $d_{c1}\big((x,t),(x'',t'')\big)\leq |\ol{y}-\ol{y}''\|\leq |\ol{y}-\ol{y}'\|+|\ol{y}'-\ol{y}''\|=d_{c1}\big((x,t),(x',t')\big)+d_{c1}\big((x',t'),(x'',t'')\big)$
\item $\theta_{\ol{y},0,\ol{y}'}+\theta_{\ol{y}',0,\ol{y}''}>\pi$: In this case, because $\theta_{\ol{y},0,\ol{y}'}=d(x,x')\leq\pi$ and $\theta_{\ol{y}',0,\ol{y}''}=d(x',x'')\leq\pi$,
{\small\bea
&&d_{c1}\big((x,t),(x'',t'')\big)=\sqrt{t^2+t''{}^2-2tt''\cos(d(x,x''))}\leq t+t''=\|0-\ol{y}\|+\|0-\ol{y}''\|\nn\\
&&~~~~\leq \|\ol{y}-\ol{y}'\|+\|\ol{y}''-\ol{y}'\|=d_{c1}\big((x,t),(x',t')\big)+d_{c1}\big((x',t'),(x'',t'')\big).\nn
\eea}where the last inequality holds because of the geometric setup, due to which triangle $(\ol{y},0,\ol{y}'')$ is contained in (hence smaller than) triangle  $(\ol{y},\ol{y}',\ol{y}'')$.\qedhere
\eit
\end{proof}

\begin{notation*}
In discussions where neither $d_{c1}$ nor $d_{c2}$ is preferred, we will simply write $d_c$.
\end{notation*}

\begin{rmk}[\textcolor{OliveGreen}{$d_c$-continuity of the quotient map of the cone}]
It follows from (\ref{Dc1Ineq}) and (\ref{Dc2Ineq}) that $q:X\times I\ra\big(\Cone(X),d_c\big)$ is $2$-Lipschitz, i.e., the quotient map $q:X\times I\ra\Cone(X)$ is $2$-Lipschitz, and hence continuous, with respect to the $d_c$-topology of $\Cone(X)$. Therefore, the quotient topology of $\Cone(X)$ contains the $d_c$-topology.
\end{rmk}

\begin{dfn}[\textcolor{blue}{Recall: Homotopy, Contracting homotopy, Contractible space}]
Let $X,Y$ be spaces. Recall that a \ul{homotopy} is a continuous map of the form $H:X\times I\ra Y$. A homotopy $H:X\times I\ra X$ is a \ul{contracting homotopy} if there is a point $x_0\in X$ such that $H(x,0)=x_0$, and $H(x,1)=x$, for all $x\in X$. A space $X$ is \ul{contractible} if a contracting homotopy $X\times I\ra X$ exists.
\end{dfn}

\begin{dfn}[\textcolor{blue}{\index{Conical homotopy}{Conical homotopy}, \index{Conical contractibility}{Conical contractibility}, Conical contracting homotopy}]\label{ConicHomContr}
Let $X,Y$ be spaces. A \ul{conical homotopy} is a continuous map of the form $H:\Cone(X)\ra Y$. A conical homotopy $H:\Cone(X)\ra X$ is a \ul{conical contracting homotopy} if $H(x,1)=x$ for all $x\in X$ (i.e., $H:\Cone(X)\ra X$ is a ``retraction''). A space $X$ is \ul{conically contractible} if a conical contracting homotopy $H:\Cone(X)\ra X$ exists.
\end{dfn}

\begin{dfn}[\textcolor{blue}{\index{CL contractibility}{Conical Lipschitz (CL) contractibility}}]
A metric space $X$ is conically $\ld$-Lipschitz contractible if there exists a contracting conical $\ld$-Lipschitz homotopy
\bea
H:(\Cone(X),d_c)\ra X\nn
\eea
(i.e., $H$ is a $\ld$-Lipschitz retraction in the sense that $H(x,1)=x$ for all $x\in X$).
\end{dfn}

\begin{dfn}[\textcolor{blue}{Modified-Lipschitz contractibility}]\label{LipContrSumm}
A metric space $X$ is $\ld$-Lipschitz contractible if every compact subset of $X$ is contained in a conically $\ld$-Lipschitz contractible set.
\end{dfn}
This definition is a weaker version of Definition \ref{ConBdLipContr} (later). Definition \ref{ConBdLipContr} is satisfied by normed spaces, but seems to be too strong for a general metric space -- especially if we are only interested in applying it to $n$-connectedness. With Definition \ref{LipContrSumm}, there is no loss of generality for Lipschitz $n$-connectedness, since it is based on continuous maps $S^k\ra X$, $0\leq k\leq n$, which have compact images.

\begin{prp}[\textcolor{OliveGreen}{Equivalence of contractibility and conical contractibility}]\label{ConicLipHom2}
A contracting homotopy $X\times I\ra X$ exists $\iff$ a conical contracting homotopy $\Cone(X)\ra X$ exists.
\end{prp}
\begin{proof}
($\Ra$): Assume a contracting homotopy $H:X\times I\ra X$ exists, i.e., $H$ is continuous and $H(x,0)=x_0$, $H(x,1)=x$. Then by the universal property of quotient maps, we get a continuous map $\wt{H}:\Cone(X)\ra X$ satisfying $\wt{H}\circ q=H$ (from which it follows that $\wt{H}$ is a conical contracting homotopy).

\bc\bt
X\times I\ar[d,"q"']\ar[rr,"H"]&& X\\
\Cone(X)\ar[urr,dashed,"\wt{H}"']
\et\ec

($\La$): Assume a conical contracting homotopy $H:\Cone(X)\ra X$ exists, i.e., $H$ is continuous and $H(x,1)=x$. Then $H\circ q:X\times I\ra X$ is a contracting homotopy, with
{\small\begin{align}
(H\circ q)(x,0)=H\big(q(x,0)\big)=x_0:=H(\txt{vertex}),~~(H\circ q)(x,1)=H\big(q(x,1)\big)=H(x,1)=x.\nn
\end{align}}
\end{proof}

\begin{question}[\textcolor{OliveGreen}{Conical Lipschitz contractibility is stronger than Lipschitz contractibility}]
Is Proposition \ref{ConicLipHom2} still true if we (i) replace ``homotopy'' with ``Lipschitz homotopy'' and (ii) replace $\Cone(X)$ with $\big(\Cone(X),d_c\big)$?

Answer = No. By Lemma \ref{ConicLipHom1} below, the direction ($\La$) works (i.e., a \ul{conical Lipschitz homotopy implies a Lipschitz homotopy}), but by Proposition \ref{MetrizQuotSp3}(ii) below, the direction ($\Ra$) does not work in general.
\end{question}

\begin{prp}\label{MetrizQuotSp3}
Let $X,Y$ be metric spaces. (i) Any Lipschitz homotopy $H:X\times I\ra Y$ that is constant on $X\times\{0\}$ (e.g., a contracting Lipschitz homotopy when $Y=X$) induces a conical homotopy $\wt{H}:\Cone(X)\ra Y$. (ii) Moreover, $\wt{H}$ is Lipschitz $\iff$ $H$ is Lipschitz with respect to $d_c$ on $(X\times I)\backslash(X\times\{0\})$, that is, $\iff$ $d\big(H(x,t),(x',t')\big)\leq K d_c\big((x,t),(x',t')\big)$, in addition to $d\big(H(x,t),(x',t')\big)\leq L d_{X\times I}\big((x,t),(x',t')\big)=L\max(|t-t'|,d(x,x'))$, outside the set $X\times\{0\}$.
\end{prp}
\begin{proof}
(i) Since $H$ is constant on $X\times\{0\}$ it follows by the universal property of quotient maps that $H$ induces a continuous map $\wt{H}:\Cone(X)\ra Y$ such that $\wt{H}\circ q=H$.
\bc\bt
X\times I\ar[d,"q"']\ar[rr,"H"]&& Y\\
\Cone(X)\ar[urr,dashed,"\wt{H}"']
\et\ec

(ii) Let $A:=(X\times I)\backslash(X\times\{0\})$. Then because $q|_A=id$, we have $\wt{H}=H$ on $A$. Thus, if $\wt{H}$ is Lipschitz, it is clear that $H$ is Lipschitz with respect to $d_c$ on $A$ (since $\wt{H}|_A$ is Lipschitz). Conversely, assume $H$ is Lipschitz with respect to $d_c$ on $A$. Then it is clear that $\wt{H}$ is Lipschitz on $A$. It remains to show $\wt{H}$ is Lipschitz near the vertex of $\Cone(X)$.

Observe that by the assumption above $\wt{H}$ is Lipschitz on $\Cone(X)\backslash\{v\}$, where $v$ is the vertex of $\Cone(X)$. Let $v_n\ra v$. Taking the limit $m\ra\infty$ in $d\big(\wt{H}(v_n),\wt{H}(v_m)\big)\leq Ld_c(v_n,v_m)$, we get $d\big(\wt{H}(v_n),\wt{H}(v)\big)\leq Ld_c(v_n,v)$.
\end{proof}

Note that the above proof involves a special case of the fact (from the ``uniform extension theorem'', Lemma \ref{UniExtThm}) that if a map between metric spaces $f:X\ra Y$ is Lipschitz on a dense subset of $X$, then $f$ is Lipschitz on all of $X$.

\begin{lmm}\label{ConicLipHom1}
Every conical Lipschitz homotopy $\big(\Cone(X),d_c\big)\ra Y$ gives a Lipschitz homotopy $X\times I\ra Y$.
\end{lmm}
\begin{proof}
Recall that $d_c\big((x,t),(x',t')\big)\sr{(\ref{Dc1Ineq}),(\ref{ConeBiLipEq})}{\leq} 12\sqrt{2}\max\big(d(x,x'),|t-t'|\big)$, which shows $q:X\times I\ra\big(\Cone(X),d_c\big)$, i.e., the quotient map $X\times I\ra\Cone(X)$ is $12\sqrt{2}$-Lipschitz with respect to the $d_c$-topolgy of $\Cone(X)$, as in the following diagram. Hence, a conical Lipschitz homotopy $H:(\Cone(X),d_c)\ra Y$ gives a Lipschitz homotopy $H\circ q:X\times I\sr{q}{\ral}\big(\Cone(X),d_c\big)\sr{H}{\ral}Y$.

\bc\bt
X\times I\ar[d,"q"']\ar[rrrr,"H\circ q"]\ar[drr,dashed,"q"] &&&& Y\\
\Cone(X) \ar[rr,"id"]&& \big(\Cone(X),d_c\big)\ar[urr,"H"]
\et\ec
\end{proof}

\begin{rmk}[\textcolor{OliveGreen}{Relevance of the cone for homotopy in general}]
The concept of a conical homotopy, as in Definition \ref{ConicHomContr}, seems relevant only when describing contractibility of spaces. It does not seem relevant in defining homotopy of maps (or homotopy equivalence of spaces) in general. This means that our earlier definition of Lipschitz homotopy (and of Lipschitz homotopy equivalence) might not be replaceable.
\end{rmk}

\section{Lipschitz $k$-connectedness of finite subset spaces of a normed space}\label{LipConnLC2}
\begin{dfn}[\textcolor{blue}{\index{Lipschitz! $k$-connectedness}{Lipschitz $k$-connectedness}}]\label{LipConnDfn}
Let $k\geq 0$. A metric space $X$ is Lipschitz $k$-connected if there exists a constant $\ld\geq 0$ such that for each $0\leq l\leq k$, every $c$-Lipschitz map $f:S^l=\del B^{l+1}\ra X$ extends to a $\ld c$-Lipschitz map $F:B^{l+1}\ra X$. Here, $B^{l+1}\subset\Real^{l+1}$ is the unit ball.

In this case, we say $X$ is $\ld$-Lipschitz $k$-connected.
\end{dfn}

\begin{question}\label{LipConQsn}
If $X$ is a normed space, is $X(n)$ Lipschitz $k$-connected for all $k\geq 0$? I.e., does there exist $\ld\geq 0$ such that every $c$-Lipschitz map $f:S^k\ra X(n)$ extends to a $\ld c$-Lipschitz map $F:B^{k+1}\ra X(n)$, for all $k\geq 0$? We will give a positive answer in Theorem \ref{LipConSuff5}.
\end{question}

\begin{dfn}[\textcolor{blue}{\index{Homogeneous metric}{Homogeneous metric} on a vector space}]
If $(X,d)$ is a metric vector space, we say $d$ is homogeneous if $d(\al x,\al y)\leq|\al|d(x,y)$ for any scalar $\al$. Note that $d$ is homogeneous $\iff$ $d(\al x,\al y)=|\al|d(x,y)$ for any scalar $\al$.
\end{dfn}

\begin{dfn*}[\textcolor{blue}{Recall: Lipschitz contractibility}]
A metric space $X$ is $\ld$-Lipschitz contractible if $X$ is $\ld$-Lipschitz homotopy equivalent to a point $x_0\in X$. Equivalently, $X$ is $\ld$-Lipschitz contractible if the identity $id:X\ra X$ is $\ld$-Lipschitz homotopic to a constant map $X\ra X$, i.e there exists a $\ld$-Lipschitz homotopy $H:X\times I\ra X$ such that $H(x,0)=x_0$ and $H(x,1)=x$ (for a point $x_0\in X$).
\end{dfn*}

\begin{dfn*}[\textcolor{blue}{Recall: Local Lipschitz contractibility}]
A metric space $X$ is $\ld$-locally Lipschitz contractible if every point $x\in X$ has a $\ld$-Lipschitz contractible neighborhood.
\end{dfn*}

\begin{dfn}[\textcolor{blue}{\index{BCL contractibility}{Bounded conical Lipschitz (BCL) contractibility}}]\label{ConBdLipContr}
A metric space $X$ is boundedly conically $\ld$-Lipschitz (BCL-) contractible if every bounded set $B\subset X$ is contained in a conically $\ld$-Lipschitz (CL-) contractible set.
\end{dfn}

Consider the following question. In a normed space, is bounded $d_{c2}$-conical $\ld$-Lipschitz contractibility scale-invariant? I.e., if $X$ is a normed space and $K\subset X$ is $d_{c2}$-conically $\ld$-Lipschitz contractible, does it follow that $\al K$ is also $d_{c2}$-conically $\ld$-Lipschitz contractible for any positive scalar $\al$? Equivalently, given a $\ld$-Lipschitz retraction $r:\Cone(K)\ra K$, does it follow that we have a $\ld$-Lipschitz retraction $R:\Cone(\al K)\ra\al K$? The answer is no, as the following remark shows.

\begin{rmk}[\textcolor{OliveGreen}{In a normed space, BCL contractibility is not scale-invariant.}]\label{ScaInvConLipCon2}
Let $X$ be a normed space. For $R>1$, if the ball $B_R\subset X$ is boundedly $d_{c2}$-conically $\ld$-Lipschitz contractible, then $\ld$ must grow with $R$, otherwise $X$ itself will be boundedly $d_{c2}$-conically $\ld$-Lipschitz contractible (a contradiction since $X$ is unbounded). This shows that the answer to the above question is no (for $\al>1$).
\end{rmk}

\begin{lmm}\label{LipHEfssRecCrl}
A metric space $X$ is ($d_{c2}$-conically) $\ld$-Lipschitz contractible $\iff$ $X(n)$ is ($d_{c2}$-conically) $\ld$-Lipschitz contractible for all $n\geq 1$. Moreover, a $\ld$-Lipschitz retraction $r:\Cone(X(n))\ra X(n)$ can be chosen to preserve cardinality, in the sense that $|r\big(x,t\big)|\leq|x|$ for each $x\in X(n)$.
\end{lmm}

\begin{proof}
The case with ``conically'' removed is an immediate corollary of Lemma \ref{LipHEfssRecall}. It remains to consider conical Lipschitz contractibility.
{\flushleft ($\Ra$)}: Assume $X$ is $d_{c2}$-conically $\ld$-Lipschitz contractible. First note that the map
\bea
f:\Cone(X(n))\ra \Cone(X)(n),~~~~\big(\{x_1,\cdots,x_n\},t\big)\mapsto\big\{(x_1,t),\cdots,(x_n,t)\big\}\nn
\eea
 is an isometry (hence 1-Lipschitz), since
\begin{align}
&(d_{c2})_H\Big(f\big(\{x_1,\cdots,x_n\},t\big),f\big(\{x'_1,\cdots,x'_n\},t'\big)\Big)\nn\\
&~~=(d_{c2})_H\Big(\big\{(x_1,t),\cdots,(x_n,t)\big\},\big\{(x'_1,t'),\cdots,(x'_n,t')\big\}\Big)\nn\\
&~~=\max\left\{\max_i\min_jd_{c2}\big((x_i,t),(x_j',t')\big),\max_j\min_id_{c2}\big((x_i,t),(x_j',t')\big)\right\}\nn\\
&~~=\max\left\{\max_i\min_j\big[|t-t'|+\min(t,t')d_{b2}(x_i,x_j')\big],\max_j\min_i\big[|t-t'|+\min(t,t')d_{b2}(x_i,x_j')\big]\right\}\nn\\
&~~=|t-t'|+\min(t,t')\max\left\{\max_i\min_jd_{b2}(x_i,x_j'),\max_j\min_id_{b2}(x_i,x_j')\right\}\nn\\
&~~= |t-t'|+\min(t,t')(d_{b2})_H\big(\{x_1,...,x_n\},\{x'_1,...,x'_n\}\big)\nn\\
&~~\sr{(s)}{=}|t-t'|+\min(t,t')(d_H)_{b2}\big(\{x_1,...,x_n\},\{x'_1,...,x'_n\}\big)\nn\\
&~~=(d_H)_{c2}\Big(\big(\{x_1,\cdots,x_n\},t\big),\big(\{x'_1,\cdots,x'_n\},t'\big)\Big),\nn
\end{align}
where step (s) holds by application of Lemma \ref{TunnelLmm} with $\vphi(\tau)=\min(2,\tau)$, and $d_{b2}$ is as in Definition \ref{ConeMetrics1}.

Let $r:\Cone(X)\ra X$ be a $\ld$-Lipschitz retraction, which induces a $\ld$-Lipschitz retraction
\bea
\wt{r}:\Cone(X)(n)\ra X(n),~~~~\wt{r}\Big(\big\{(x_1,t),...,(x_n,t)\big\}\Big):=\big\{r(x_1,t),...,r(x_n,t)\big\}.\nn
\eea
Then we get a $\ld$-Lipschitz retraction
\bea
\wt{r}\circ f:\Cone(X(n))\sr{f}{\ral}\Cone(X)(n)\sr{\wt{r}}{\ral}X(n).\nn
\eea
By construction, it is clear that the retraction $\wt{r}\circ f$ preserves cardinality, since
\bea
\wt{r}\circ f\big(\{x_1,...,x_n\},t\big)=\wt{r}\big\{(x_1,t),...,(x_n,t)\big\}=\big\{r(x_1,t),...,r(x_n,t)\big\}.\nn
\eea

{\flushleft ($\La$)}: If each $X(n)$ is conically $\ld$-Lipschitz contractible, it is clear that $X\sr{1-\txt{biLip}}{\cong} X(1)$ is conically $\ld$-Lipschitz contractible.
\end{proof}

\begin{lmm}\label{LipConSuff3}
If a metric space $X$ is boundedly ($d_{c2}$-conically) $\ld$-Lipschitz contractible, then $X(n)$ is boundedly ($d_{c2}$-conically) $\ld$-Lipschitz contractible.
\end{lmm}
\begin{proof}
Let $A\subset X(n)$ be bounded. Then $\bigcup A:=\bigcup_{a\in A}a\subset X$ is bounded (by the bounded image lemma, Lemma \ref{BndImgLmm}), and thus contained in a ($d_{c2}$-conically) $\ld$-Lipschitz contractible set $K\subset X$. Let $A':=\{x\in X(n):x\subset K\}=K(n)$. Then $\bigcup A\subset K$ implies $A\subset\big(\bigcup A\big)(n)\subset K(n)=A'$, and by Lemma \ref{LipHEfssRecCrl}, $A'$ is ($d_{c2}$-conically) $\ld$-Lipschitz contractible.
\end{proof}

\begin{lmm}\label{LipConSuff4}
Let $X$ be a normed space. The ball $B_R=B_R(0)\subset X$ is boundedly $d_{c2}$-conically $\max(1,R)$-Lipschitz contractible. That is, the following are true.
\bit[leftmargin=0.9cm]
\item[(i)] If $0<R\leqslant 1$, then $B_R$ is boundedly $d_{c2}$-conically $1$-Lipschitz contractible.
\item[(ii)] If $1<R<\infty$, then $B_R$ is boundedly $d_{c2}$-conically $R$-Lipschitz contractible.
\eit
\end{lmm}
\begin{proof}
{\flushleft (i)} \ul{$0<R\leqslant 1$}:  Define a map $r:\Cone(B)\ra B$ by $r(x,t):=tx$. Then, using the triangle inequality,
\bea
\|r(x,t)-r(x',t')\|=\|tx-t'x'\|\leq R|t-t'|+\min(t,t')\|x-x'\|\leq d_{c2}\big((x,t),(x',t')\big).\nn
\eea
{\flushleft (ii)} \ul{$1<R<\infty$}: With $r:\Cone(B_R)\ra B_R$ given by the same map as in (i), we again have
\begin{align}
&\textstyle\|r(x,t)-r(x',t')\|=\|tx-t'x'\|\leq R|t-t'|+\min(t,t')\|x-x'\|\nn\\
&\textstyle~~=R\left(|t-t'|+\min(t,t'){\|x-x'\|\over R}\right)\sr{(a)}{=}R\left(|t-t'|+\min(t,t')\min\Big(2,{\|x-x'\|\over R}\Big)\right)\nn\\
&\textstyle~~\sr{(b)}{\leq} R\Big(|t-t'|+\min(t,t')\min\big(2,\|x-x'\|\big)\Big)=Rd_{c2}\big((x,t),(x',t')\big),\nn
\end{align}
where step (a) holds because $\|x-x'\|\leq\diam B_R=2R$, and step (b) holds because $R>1$.
\end{proof}

\begin{crl}\label{LipConSuff4a}
The unit ball of a normed space is $24$-Lipschitz $k$-connected.
\end{crl}
\begin{proof}
Let $X$ be a normed space, $B_R:=B_R(0)$ the ball of radius $R>0$ centered at $0\in X$, and $B:=B_1(0)$ the unit ball. Let $f:S^k\ra B$ be Lipschitz. Note that ${\diam f(S^k)\over 2}={\diam f(S^k)\over\diam(S^k)}\leq\Lip f$, i.e., $\diam f(S^k)\leq2\Lip(f)$. (Without loss of generality) translate $f$ so that $0\in f(S^k)$. Then $f(S^k)\subset B_{\diam f(S^k)}\subset B_{2\Lip(f)}$, i.e.,
\bea
\textstyle {1\over 2\Lip f}f(S^k)\subset B.\nn
\eea
Define $\wt{f}:S^k\ra B$ by $\wt{f}:={1\over 2\Lip f} f$, where $\Lip(\wt{f})={1\over 2}$. Then, with a 1-Lipschitz retraction $\Cone(B)\sr{r}{\ral}B$, we get the extension $B^{k+1}\sr{F}{\ral}B$ of $f$ given by $F(ts):=(2\Lip f)~r\big(\wt{f}(s),t\big)$,
\bea\bt[row sep=tiny]
F:B^{k+1}\ar[r,"\cong"] & \Cone(S^k)={S^k\times I\over S^k\times\{0\}}\ar[r,"\wt{f}\times id"] & {B\times I\over B\times\{0\}}=\Cone(B)\ar[r,"r"] & B.\\
 st\ar[r,mapsto] & (s,t) & &
\et\nn
\eea
\begin{align}
\|F(ts)&-F(t's')\|=(2\Lip f)~\|r\big(\wt{f}(s),t\big)-r\big(\wt{f}(s'),t'\big)\|\leq(2\Lip f)~ d_{c2}\Big(\big(\wt{f}(s),t\big),\big(\wt{f}(s'),t'\big)\Big)\nn\\
&\leq(2\Lip f)~d_{c2}\big((s,t),(s',t')\big)\sr{(\ref{ConeBiLipEq})}{\leq}12(2\Lip f)~ d_{c1}\big((s,t),(s',t')\big)\nn\\
&\textstyle=(24\Lip f)~\sqrt{|t-t'|^2+4tt'\sin^2{\|s-s'\|\over 2}}~\sr{(a)}{\leq}~(24\Lip f)~\sqrt{|t-t'|^2+4tt'\sin^2{\theta_{s,0,s'}\over 2}}\nn\\
&\textstyle=(24\Lip f)~\sqrt{t^2+t'{}^2-2tt'\cos{\theta_{s,0,s'}}}=(24\Lip f)~\|ts-t's'\|,\nn
\end{align}
where step (a) holds because the distance $\|s-s'\|$ on $S^k$ is smaller than the length $\theta_{s,0,s'}$ of the shortest arc from $s$ to $s'$ on $S^k$.
\end{proof}

\begin{crl}\label{LipConSuff4b}
A normed space is $24$-Lipschitz $k$-connected.
\end{crl}
\begin{proof}
Let $X$ be a normed space, and consider a Lipschitz map $f:S^k\ra X$. Since $f$ is bounded, there is a scalar $\al\neq 0$ such that $\al f(S^k)\subset B$, where $B:=B_1(0)$ is the unit ball. Thus, by Corollary \ref{LipConSuff4a}, the map $\wt{f}=\al f:S^k\ra B$ has a $24\al\Lip(f)$-Lipschitz extension $\wt{F}:B^{k+1}\ra B$. Hence, $F={1\over \al}\wt{F}$ is a $24\Lip(f)$-Lipschitz extension of $f$.
\end{proof}

\begin{dfn*}[\textcolor{blue}{Recall: Gap of a set}]
Let $X$ be a metric space and $A\subset X$. The gap $\rho(A)$ of $A$ is the largest distance between any two nonempty complementary subsets of $A$, i.e.,
\bea
\rho(A):=\sup_{\emptyset\neq A'\subsetneq A}\dist(A',A\backslash A')=\sup_{\substack{A'\sqcup A''=A \\ A',A''\neq\emptyset}}\dist(A',A'').\nn
\eea
\end{dfn*}

\begin{lmm}[\textcolor{OliveGreen}{Bound on finite set diameter}]\label{SubsetGapLmm}
Let $X$ be a metric space and $A\subset X$ a finite set. Then for every $k\leq|A|$, there exists a set $A_k\subset A$ satisfying $|A_k|=k$ and ~$\diam A_k\leq (k-1)\rho(A)$. In particular, $\diam A\leq(|A|-1)\rho(A)$.
\end{lmm}
\begin{proof}
Let $n:=|A|$. We proceed by induction on $1\leq k\leq n$. The result is clear for $k=1$. For the induction step, suppose $k\leq n$ implies the existence of $A_k\subset A$ satisfying $\diam A_k\leq(k-1)\rho(A)$.

If $k+1\leq n$ (which implies $1\leq k<n$), then $\dist(A_k,A\backslash A_k)\leq\rho(A)$, and so there exists $a_{k+1}\in A\backslash A_k$ such that $\dist(A_k,a_{k+1})\leq\rho(A)$. Thus, with $A_{k+1}:=A_k\cup\{a_{k+1}\}$, we get
{\small\bea
\diam A_{k+1}\leq\diam A_k+\dist(A_k,a_{k+1})\leq (k-1)\rho(A)+\rho(A)=k\rho(A)=\big((k+1)-1\big)\rho(A).\nn\qedhere
\eea}
\end{proof}

\begin{thm}[\textcolor{OliveGreen}{Answer to Question \ref{LipConQsn}}]\label{LipConSuff5}
If $X$ is a normed space, then there exists a constant $\ld_n\geq 0$ such that $X(n)$ is $\ld_n$-Lipschitz $k$-connected for all $k\geq 0$.
\end{thm}
\begin{proof}
Let $R>0$. If $A\subset X$, then $B_R(A):=\{u\in X:\dist(u,A)<R\}$ denotes the $R$-neighborhood of $A$ in $X$. If $W\subset X(n)$, then $N_R(W):=\{x\in X(n):\dist_H(x,W)<R\}$ denotes the $R$-neighborhood of $W$ in $X(n)$. We also write $B=B_1=B_1(0)$. Note that
\bea
X_0(n):=\{x\in X(n):0\in x,~\diam x\leq 1\}\subset N_1(\{0\})=B(n).\nn
\eea
Let $f:S^k\ra X(n)$ be $c$-Lipschitz. If $c=0$, then $f$ is constant, and so extends to a constant map $F:B^{k+1}\ra X(n)$. So, assume $c>0$. Since ~$\diam f(S^k)\leq 2\Lip f=2c$,~ it follows that for any $s\in S^k$, we have
\begin{align}
\textstyle f(S^k)\subset N_{\diam f(S^k)}\big(f(s)\big)\subset N_{2c}\big(f(s)\big),~~\Ra~~{1\over 2c}f(S^k)\subset N_1\left({1\over 2c}f(s)\right).\nn
\end{align}

Fix any $s_0\in S^k$ such that $0\in f(s_0)$, which is possible by translating $f$. Then we have two cases.
\bit[leftmargin=0.4cm]
\item\ul{$\diam f(s_0)\leq Kc$}: With $z:={f(s_0)\over\diam f(s_0)}\in B(n)$, we have
\bea
&&f(S^k)\subset N_{2c}\big(f(s_0)\big)=\diam f(s_0)~N_{{2c\over\diam f(s_0)}}\big(z\big)\subset \diam f(s_0)~N_{{2c\over\diam f(s_0)}}\big(B(n)\big)\nn\\
&&~~~~\subset \diam f(s_0)~N_{{2c\over\diam f(s_0)}+1}\big(\{0\}\big)=\big[2c+\diam f(s_0)\big]~N_1\big(\{0\}\big),\nn\\
&&\textstyle~~\Ra~~f(S^k)\subset (2+K)c~B(n),~~\Ra~~{1\over(2+K)c}f(S^k)\subset B(n).\nn
\eea
By Lemmas \ref{LipHEfssRecCrl} and \ref{LipConSuff4}, there is a $1$-Lipschitz retraction $r:\Cone B(n)\ra B(n)$. Define an extension $B^{k+1}\sr{F}{\ral}X(n)$ of $f$ by $F(ts):=(2+K)c~r\big(\wt{f}(s),t\big)$, where $\wt{f}:={1\over(2+K)c}f:S^k\ra B(n)$. That is,
\bea\bt[row sep=tiny]
F:B^{k+1}\ar[r,"\cong"] & \Cone(S^k)={S^k\times I\over S^k\times\{0\}}\ar[r,"\wt{f}\times id"] & {B(n)\times I\over B(n)\times\{0\}}=\Cone\big(B(n)\big)\ar[r,"r"] & B(n).\\
 st\ar[r,mapsto] & (s,t) & &
\et\nn
\eea
Then
\begin{align}
d_H\Big(F(ts)&,F(t's')\Big)=(2+K)c~~d_H\Big(r\big(\wt{f}(s),t\big),r\big(\wt{f}(s'),t'\big)\Big)\nn\\
&\leq(2+K)c~~ (d_H)_{c2}\Big(\big(\wt{f}(s),t\big),\big(\wt{f}(s'),t'\big)\Big)\nn\\
&=(2+K)c~\bigg(|t-t'|+\min(t,t')\min\left[2,d_H\big(\wt{f}(s),\wt{f}(s')\big)\right]\bigg)\nn\\
&\leq(2+K)c~~d_{c2}\big((s,t),(s',t')\big)\sr{(\ref{ConeBiLipEq})}{\leq}12(2+K)c~~ d_{c1}\big((s,t),(s',t')\big)\nn\\
&\textstyle=\ld c~\sqrt{|t-t'|^2+4tt'\sin^2{\|s-s'\|\over 2}}~\sr{(a)}{\leq}~\ld c~\sqrt{|t-t'|^2+4tt'\sin^2{\theta_{s,0,s'}\over 2}}\nn\\
&\textstyle=\ld c~\sqrt{t^2+t'{}^2-2tt'\cos{\theta_{s,0,s'}}}=\ld c~\|ts-t's'\|,\nn
\end{align}
where $\ld:=12(2+K)$ and step (a) holds because the distance $\|s-s'\|$ on $S^k$ is smaller than the length $\theta_{s,0,s'}$ of the shortest arc from $s$ to $s'$ on $S^k$.

\item\ul{$\diam f(s_0)> Kc$}: There exists (by Lemma \ref{SubsetGapLmm}) a decomposition ~$f(s_0)=x'\cup x''$~ (with $x',x''$ nonempty) such that ~$\dist(x',x'')\geq{\diam f(s_0)\over n}$,~ and so $\dist(x',x'')>{K\over n}c$. We will choose $K\geq 5n$, so that
\bea
\dist(x',x'')> 5c,\nn
\eea
which implies the neighborhoods $B_{2c}(x')$ and $B_{2c}(x'')$ in $X$ are distantly separated. Since
{\small\bea
f(S^k)\subset N_{2c}\big(f(s_0)\big)=N_{2c}\big(x'\cup x''\big)~~\txt{i.e.,}~~f(s)\in N_{2c}\big(f(s_0)\big)=N_{2c}\big(x'\cup x''\big)~~\txt{for all}~~s\in S^k,\nn
\eea}
it follows (from the definition of Hausdorff distance -- see the footnote\footnote{
Observe that if $d_H\big(f(s),f(s_0)\big)=\max\left(\max_i\min_jd\big(f(s)_i,f(s_0)_j\big),\max_i\min_jd\big(f(s)_j,f(s_0)_i\big)\right)\leq 2c$, then $\max_i~\dist\big(f(s)_i,f(s_0)\big)\leq 2c$.
}) that for each $s\in S^k$, we have
\bea
f(s)\subset B_{2c}\big(f(s_0)\big)=B_{2c}(x'\cup x'')=B_{2c}(x')\cup B_{2c}(x'').\nn
\eea
Let $x'_s=f_1(s):=f(s)\cap B_{2c}(x')$ and $x''_s=f_2(s):=f(s)\cap B_{2c}(x'')$. Then we have the decomposition ~$f(s)=x'_s\cup x''_s=f_1(s)\cup f_2(s)$. Observe that, by direct calculation at step (a) below,
\begin{align}
d_H\big(f(s),&f(s_0)\big)=d_H\big(f_1(s)\cup f_2(s),x'\cup x''\big)\sr{(a)}{=}\max\Big(d_H\big(f_1(s),x'\big),d_H\big(f_2(s),x''\big)\Big)\leq 2c,\nn\\
&~~\Ra~~f_1(s)\in N_{2c}(x')~~~~\txt{and}~~~~f_2(s)\in N_{2c}(x'').\nn
\end{align}
The maps $f_1:S^k\ra X(n-1)$ and $f_2:S^k\ra X(n-1)$ are $c$-Lipschitz, since for $j=1,2$
\bea
&&d_H\Big(f_j(s),f_j(s')\Big)\leq \max\left[d_H\Big(f_1(s),f_1(s')\Big)~,~d_H\Big(f_2(s),f_2(s')\Big)\right]\nn\\
&&~~~~\sr{(a)}{=}d_H\Big(f_1(s)\cup f_2(s)~,~f_1(s')\cup f_2(s')\Big)=d_H\big(f(s),f(s')\big),\nn
\eea
where step (a) holds because $f_1(s)\cup f_1(s')$ and $f_2(s)\cup f_2(s')$ are distantly separated.

Thus, by induction on $n$ (where $n=1$ holds by Corollary \ref{LipConSuff4a}), $f_1,f_2$ have Lipschitz extensions
\bea
&&\textstyle F_1:B^{k+1}\ra X(n-1),~~~~F_1(ts):=K_1c~r_1\big(\wt{f}_1(s),t\big),~~~~\wt{f}_1={1\over K_1c}f_1:S^k\ra B(n-1),\nn\\
&&\textstyle F_2:B^{k+1}\ra X(n-1),~~~~F_2(ts):=K_2c~r_2\big(\wt{f}_2(s),t\big),~~~~\wt{f}_2={1\over K_2c}f_2:S^k\ra B(n-1),\nn
\eea
where $r_1:\Cone(B(n-1))\ra B(n-1)$ and $r_2:\Cone(B(n-1))\ra B(n-1)$ are $1$-Lipschitz retractions that exist by Lemmas \ref{LipHEfssRecCrl} and \ref{LipConSuff4}. Note that by Lemma \ref{LipHEfssRecCrl}, the above retractions each preserve cardinality in the sense that
\bea
|F_j(ts)|=\big|r_j\big(\wt{f}_j(s),t\big)\big|\leq \big|\wt{f}_j(s)\big|=|f_j(s)|.\nn
\eea
Define an extension $F:S^k\ra X(n)$ of $f$ by ~$F(ts):=F_1(ts)\cup F_2(ts)$. With $\ld:=12\max(K_1,K_2)$,
\begin{align}
&\textstyle d_H\big(F(ts),F(t's')\big)\leq{\ld\over 12}c~d_H\Big(r_1\big(\wt{f}_1(s),t\big)\cup r_2\big(\wt{f}_2(s),t\big)~,~r_1\big(\wt{f}_1(s'),t'\big)\cup r_2\big(\wt{f}_2(s'),t'\big)\Big)\nn\\
&~~~~\leq \ld c~\|ts-t's'\|.\nn\hspace{5cm}\qedhere
\end{align}
\eit
\end{proof}

\begin{dfn}[\textcolor{blue}{\index{Marginal $k$-connectedness}{Marginal $k$-connectedness}, \index{Marginal Lipschitz $k$-connectedness}{Marginal Lipschitz $k$-connectedness}}]
Let $X$ be a topological space (resp. metric space) , $\ld\geq 0$, and $k\geq 0$ an integer. Then $X$ is marginally topologically (resp. marginally $\ld$-Lipschitz) $k$-connected if every continuous (resp. $c$-Lipschitz) map $f:S^k\ra X$ extends to a continuous (resp. $\ld c$-Lipschitz) map $F:B^{k+1}\ra X$, where $B^{k+1}\subset\Real^{k+1}$ is the unit ball bounded by $S^k$.
\end{dfn}

It is clear that a space $X$ is (Lipschitz) $k$-connected for all $k\geq 0$ if and only if $X$ is marginally (Lipschitz) $k$-connected for all $k\geq 0$.

\begin{note}
If $X$ is marginally Lipschitz $k$-connected and $Z\subset X$ is a Lipschitz retract, then $Z$ is also marginally Lipschitz $k$-connected (as shown in the following diagram).
\bc\bt
S^k\ar[rrrr,bend left=25,"g"]\ar[d,hook]\ar[rr,near end,"f"]&& Z\ar[rr,hook] && X\ar[d,"r"]\\
B^{k+1}\ar[urrrr,near start,dashed,"G"]\ar[rrrr,dashed,near end,"F:=r\circ G"]&&&& Z
\et\ec
\end{note}

\begin{lmm}\label{S1LipCon}
The circle $S^1$ is marginally $\pi\over 2$-Lipschitz $k$-connected for $k\geq 2$.
\end{lmm}
\begin{proof}
Fix $k>1$. Consider the universal cover $p:\Real\ra S^1$. Let $f:S^k\ra S^1$ be $c$-Lipschitz. Since $\pi_1(S^k)$ is trivial, it follows from \cite[Proposition 1.33, p.61]{hatcher2001} that there exists a unique continuous lift $\wt{f}:S^k\ra\Real$ such that $p\circ\wt{f}=f$.
\bc\bt
B^{k+1}\ar[rrrr,dashed,"\wt{F}"] &&&& \Real\ar[d,"p(t)=e^{it}"]\\
S^k\ar[u,hook,"i"]\ar[urrrr,dashed,near start,"\wt{f}"]\ar[rrrr,near end,"f"] &&&& S^1
\et\ec
(Note that $\wt{f}$ is uniformly continuous since $S^k\subset\Real^k$ is a compact metric space). Since $p$ is locally isometric, it follows that $\wt{f}$ is locally $c$-Lipschitz, hence ${\pi\over 2}c$-Lipschitz, since $S^k$ is ${\pi\over 2}$-quasiconvex: Indeed, if $s,s'\in S^k$, then the path $\gamma:[0,1]\ra S^k$ from $s$ to $s'$ satisfying $\gamma(t)\cdot\gamma(t')=\cos\big((t-t')\theta_{s,0,s'}\big)$ is a ${\pi\over 2}$-quasigeodesic.

By the McShane-Whitney extension theorem, $\wt{f}$ has a ${\pi\over 2}c$-Lipschitz extension $\wt{F}:B^{k+1}\ra\Real$. Hence, $F=p\circ\wt{F}:B^{k+1}\ra S^1$ is a ${\pi\over 2}c$-Lipschitz extension of $f$.
\end{proof}

\begin{question}\label{XXnLipCon}
Let $X$ be a metric space. If $X$ is marginally $\ld$-Lipschitz $k$-connected, does it follow that $X(n)$ is marginally $\ld$-Lipschitz $k$-connected? Topologically, the answer is ``no'': Indeed, with $X=S^1$, we have $\pi_3(S^1)=\{0\}$, while $S^1(3)\cong S^3$ (by \cite[Corollary 5.3]{chinen-koyama2010} and \cite[Theorems 1, 3]{tuffley2002} ) implies $\pi_3\big(S^1(3)\big)=\Integer$.
\end{question}

\begin{prp}
If the homeomorphism $S^1(3)\cong S^3$ can be chosen to be Lipschitz (i.e., the answer to \cite[Question 4.20]{chinen2015} is positive), then the answer to Question \ref{XXnLipCon} is also ``no'' for the Lipschitz case.
\end{prp}
\begin{proof}
By Lemma \ref{S1LipCon}, $X=S^1$ is marginally ${\pi\over 2}$-Lipschitz $k$-connected for $k\geq 2$. On the other hand, if $S^1(3)\cong_{\Lip} S^3$, then because the homotopy class of the identity map $S^3\ra S^3$ is a nontrivial Lipschitz element of $\pi_3\big(S^1(3)\big)=\pi_3(S^3)=\Integer$, it follows that $S^1(3)$ is not marginally Lipschitz $3$-connected.
\end{proof}


\chapter{FSR Property I: Finite Subsets of Normed Spaces}\label{FSRP1}
The finite subset retraction (FSR) property concerns the existence of Lipschitz retractions $X(n)\ra X(n-1)$. In this chapter we introduce the FSR property and investigate some of its implications for normed spaces. Sections \ref{LipRetV3}, \ref{LipRetVn12}, \ref{Retractions} are entirely based on \cite{akofor2019}. It was asked in \cite[Question 3.4]{kovalev2016} whether Lipschitz retractions $X(n)\ra X(n-1)$ exist for all $n\geq2$ when $X$ is a Banach space. A related question in \cite[Remark 3.4]{bacac-kovalev2016} similarly asked whether strictly convex or uniformly convex Banach spaces admit such Lipschitz retractions. The results of this chapter, and those of the previous chapters, provide partial answers and tools of investigation towards answering the above questions.

Question \ref{MajQuest} below appeared (as Problem 1.4) in the collection of open problems ``\emph{AimPL: Mapping theory in metric spaces}'' published by the American Institute of Mathematics and available at \url{http://aimpl.org/mappingmetric} (see \cite{ProbList}).

\begin{question}\label{MajQuest}
If $X$ is an ALR, does it follow that $X(n)$ is an ALR? (For example, since $\ell^\infty$ is an ALR, is $\ell^\infty(n)$ also an ALR?)
\end{question}
Recall that the ALR property of $X(n)$ means that any larger metric space $Y\supset X(n)$ retracts onto $X(n)$ via a Lipschitz map $r:Y\ra X(n)$. A natural candidate for such a space $Y$ is $X(m)$, $m>n$, which thus provides us a tool with which to test Question \ref{MajQuest}. Since a composition of Lipschitz retractions is a Lipschitz retraction, it suffices to consider retractions with $m-n=1$. Thus, we focus on the existence of Lipschitz retractions $X(n) \to X(n-1)$.

\begin{dfn}[\textcolor{blue}{\index{FSR property}{Finite subset retraction (FSR) property}}]
We say a metric space $X$ has the \emph{finite subset retraction ($\FSR$) property} if there exist Lipschitz retractions $X(n)\ra X(n-1)$ for all $n\geq 2$. We also say $X$ has the \emph{$\FSR(k)$ property} if for each $1\leq l\leq k$, there exist Lipschitz retractions $X(n)\ra X(l)$ for all $n\geq l$. Equivalently, $X$ has the \emph{$\FSR(k)$ property} if there exist Lipschitz retractions (i) $X(n)\ra X(n-1)$ for all $2\leq n\leq k$,  and (ii) $X(n)\ra X(k)$ for all $n\geq k$.
\end{dfn}

We will see in sections \ref{LipRetV3} and \ref{LipRetVn12} (Theorems \ref{LrThmV3}, \ref{LrThmVn2}) that every normed space has the FSR(3) property, and in Remark \ref{HilbConRmk} of section \ref{Retractions}, which is based on Theorem \ref{MainThm}, that (as established in \cite{kovalev2016}) every Hilbert space has the FSR property. We will also see in section \ref{Retractions} (Theorem \ref{MainThm}) that if $X$ is a normed space, then we have retractions $X(n)\ra X(n-1)$ that are H\"{o}lder continuous on bounded sets. Section \ref{Selections} examines an interesting property shared by the constructed retractions and discusses some of its consequences.

\section{Concrete Lipschitz retractions $X(2)\ra X$ and $X(3)\ra X(2)$}\label{LipRetV3}
The Lipschitz retractions in this section, unlike those in Section \ref{LipRetVn12}, have concrete Lipschitz constants.

\begin{dfn}[\textcolor{blue}{Addition and scalar multiplication of sets}]\label{SumSclMult} Let $X$ be a vector space. If $A,B\subset X$ and $\ld$ is a scalar, we write $A+B:=\{a+b:a\in A,b\in B\}$ and $\ld A:=\{\ld a:a\in A\}$.
\end{dfn}

\begin{dfn}[\textcolor{blue}{Scale-invariant map, Translation-invariant map, \index{Affine! map}{Affine map}}]
Let $X$ be a vector space and $f:\E\subset\P^\ast(X)\ra\P^\ast(X)$ a map. We say $f$ is \ul{scale-invariant} (or commutes with scaling) if for any $t\in \Real$, we have $f(tA)=tf(A)$ for all $A\in\E$ such that $tA\in\E$. Similarly, $f$ is \ul{translation-invariant} (or commutes with translations) if for any $v\in X$, we have $f(A+v)=f(A)+v$ for all $A\in\E$ such that $A+v\in\E$.

We say $f$ is \ul{affine} if $f$ is both scale-invariant and translation-invariant, i.e., for any $t\in\Real$, $v\in X$, we have $f(tA+v)=tf(A)+v$ for all $A\in\E$ such that $tA+v\in\E$.
\end{dfn}

\begin{dfn}[\textcolor{blue}{Recall: Proximal map between points of $X(n)$}]
Let $X$ be a metric space and $x,y\in X(n)$. A map $p:x\ra y$ is proximal if $d(a,p(a))\leq d_H(x,y)$ for all $a\in x$, i.e., the relation $R_p:=\big\{\big(a,p(a)\big):a\in x\big\}\subset x\times y$ is proximal.
\end{dfn}

\begin{dfn}[\textcolor{blue}{\index{Normalized element}{Normalized element} in $X(n)$, Set of normalized elements}]
Let $X$ be a metric space and $x\in X(n)$. We say $x$ is \ul{normalized} if $\diam(x)=1$. We will write $N\big(X(n)\big)$ for all normalized elements of $X(n)$.
\end{dfn}

\begin{dfn}[\textcolor{blue}{\index{Central element}{Central element} in $X(n)$, Set of central elements, Set of normalized central elements}]
Let $X$ be a normed space and $x\in X(n)$. We say $x$ is \ul{central} if $0\in x$. We will write $X_0(n)=\{x\in X(n):0\in x\}$ for all central elements of $X(n)$. Accordingly, we will write $N(X_0(n))$ for the set of normalized central elements of $X(n)$.
\end{dfn}

Note that every element $x\in X(n)$ can be written (not uniquely) as
\begin{equation*}
x=t x_0+v,~~~~\txt{for some}~~t\in[0,+\infty),~~x_0\in N(X_0(n)),~~v\in X.
\end{equation*}

\begin{lmm}[\textcolor{OliveGreen}{Homogeneous Lipschitz Extension}]\label{HomExtLmm}
Let $X$ be a normed space and $1\leq k\leq n-1$. Any translation-invariant Lipschitz map $R:N(X_0(n))\ra X(k)$ satisfying $R(x)\subset\Conv(x)$ and $R|_{N(X_0(n))\cap X(k)} = id$ can be extended to an affine Lipschitz retraction $r:X(n)\ra X(k)$ with $\Lip(r)=6\Lip(R)+5$.
\end{lmm}
\begin{proof}
Let $R:N(X_0(n))\ra X(k)$ be a Lipschitz map such that $R(x)\subset\Conv(x)$ and $R|_{N(X_0(n))\cap X(k)}$ = $id$. Define a map $r:X(n)\ra X(k)$ by $r(tx+v):=tR(x)+v$ for all $x\in N(X_0(n))$, $t\in[0,+\infty)$, and $v\in X$. Then $r$ is well defined because $R$ is translation-invariant. For any $x,y\in N(X_0(n))$ and $t,s\in[0,+\infty)$, since $0\in x$ and diameter is $2$-Lipschitz with respect to $d_H$,
\begin{equation*}
d_H(tx,sx)\leq |t-s|=|\diam(tx)-\diam(s y)|\leq 2d_H(tx,sy).
\end{equation*}
 Thus, using the triangle inequality, $0\in y$, and $R(y)\subset \Conv(y)$, we have
\begin{align*}
d_H(&r(tx),r(sy))\leq d_H(r(tx),r(ty))+d_H(r(ty),r(sy))=d_H(tR(x),tR(y))+d_H(tR(y),sR(y))\\
    &\leq\Lip(R)d_H(tx,ty)+\diam(R(y))|t-s|\leq \Lip(R)\Big[d_H(tx,sy)+d_H(sy,ty)\Big]+|t-s|\\
    &\leq\Lip(R)d_H(tx,sy)+\Big(\Lip(R)+1\Big)|t-s|\leq \Big(3\Lip(R)+2\Big)d_H(tx,sy),
\end{align*}
which shows $r$ is Lipschitz on $X_0(n)$. Given $x,y\in X(n)$, let $u\in x$, $v\in y$ such that $\|u-v\|\leq d_H(x,y)$. Then
\begin{align*}
d_H\big(r(x),r(y)\big)&=d_H\big(r(x-u)+u,r(y-v)+v\big)\leq d_H\big(r(x-u),r(y-v)\big)+\|u-v\|\\
 &\leq\big(2\Lip(r|_{X_0(n)})+1\big)d_H(x,y).\qedhere
\end{align*}
\end{proof}

\begin{dfn}[\textcolor{blue}{\index{Thin sets}{Thin sets in $X(3)$}}]
Let $X$ be a normed space and let $x\in X(3)$. Then $x$ is called ``\ul{thin} in $X(3)$'' if $x$ is normalized and $0\leq\delta(x)\leq{1\over 3}$. We will denote all thin sets in $X(3)$ by $\Thin\big(X(3)\big)$.
\end{dfn}

\begin{notation}
If $x=\{x_1,x_2,x_3\}\in \Thin\big(X(3)\big)$, we will assume wlog that
\begin{equation*}
\delta(x)=d(x_1,x_2)\leq d(x_2,x_3)\leq d(x_1,x_3)=1,
\end{equation*}
where the triangle inequality implies $d(x_2,x_3)\geq 2/3$.
\end{notation}

\begin{dfn}[\textcolor{blue}{\index{Vertex map}{Vertex map} of thin elements of  $X(3)$}]\label{ExtrMinDfn}
This is the map $V:\Thin\big(X(3)\big)\ra X$ given by $V(\{x_1,x_2,x_3\}):=x_3$.
\end{dfn}

\begin{dfn}[\textcolor{blue}{\index{Average map}{Average map}}]
If $X$ is a normed space, the average $\Avg:X(n)\ra X$ is given by $\Avg(x)~:=~{1\over |x|}\sum_{a\in x}a$. In particular, if $x=\{x_1,...,x_n\}\in X(n)\backslash X(n-1)$, then we can write $\Avg(x)={1\over n}\sum_{i=1}^n x_i$.
\end{dfn}

Throughout the rest of this section, we will assume $X$ is a normed space.

\begin{lmm}[\textcolor{OliveGreen}{Continuity of the Average map}]\label{ContAvgLmm}
Let $X$ be a normed space and $x,y\in X(n)\backslash X(n-1)$.
\bit[leftmargin=0.9cm]
\item[(i)] If a proximal bijection $x\ra y$ exists, then ~~$d_H(\Avg(x),\Avg(y))\leq d_H(x,y)$.
\item[(ii)] If $\max\Big(\delta(x),\delta(y)\Big)\leq 2d_H(x,y)$, and $\al\diam(x)\leq\delta(x)$ or $\al\diam(y)\leq\delta(y)$ for constant $\al>0$, then
\begin{equation*}
\textstyle d_H(\Avg(x),\Avg(y))\leq \left(1+{2\over\al}\right)d_H(x,y).
\end{equation*}
\eit
\end{lmm}
\begin{proof}
(i) Let $x\ra y$, $x_i\mapsto y(i)$ be a proximal bijection. Then
\begin{align*}
\textstyle d_H(\Avg(x),\Avg(y))=\left\|{\sum x_i\over n}-{\sum y_i\over n}\right\|=\left\|{\sum x_i\over n}-{\sum y(i)\over n}\right\|\leq d_H(x,y).
\end{align*}

(ii) It suffices to assume $\al\diam(y)\leq\delta(y)$. Consider a proximal map $x\ra y$, $x_i\mapsto y(i)$. Then
\begin{align*}
d_H(&\Avg(x),\Avg(y))\textstyle=\left\|{\sum x_i\over n}-{\sum y_i\over n}\right\|\leq {\sum\|x_i-y(i)\|+\sum_i\|y(i)-y_i\|\over n}\\
  &\textstyle\leq d_H(x,y)+\diam(y)\leq d_H(x,y)+{\delta(y)\over\al}\leq\left(1+{2\over \al}\right)d_H(x,y).\qedhere
\end{align*}
\end{proof}

\begin{prp}\label{LrPrpV2}
Let $X$ be a normed space. There exists a $1$-Lipschitz retraction $X(2)\ra X$.
\end{prp}
\begin{proof}
Define $r:X(2)\ra X$ by $r(x)=\Avg(x)$. Then $r|_X=id$. If $x,y\in X(n)$, consider the following.
\bit[leftmargin=0.9cm]
\item[(i)] $x=\{x_1\}\in X$ and $y=\{y_1,y_2\}\in X(2)\backslash X$:~ In this case,
\begin{equation*}
\|r(x)-r(y)\|\leq \max(\|x_1-y_1\|,\|x_1-y_2\|)=d_H(x,y).
\end{equation*}
\item[(ii)] $x,y\in X(2)\backslash X$: A proximal bijection $x\ra y$ exists, and so by Lemma \ref{ContAvgLmm}(i), $d_H(r(x),r(y))\leq d_H(x,y)$.
\eit
Hence, $r$ is a $1$-Lipschitz retraction.
\end{proof}

\begin{lmm}\label{LipCenLmm1}
If $X$ is a normed space, the following map is $3$-Lipschitz.
\begin{equation*}
\textstyle f:\Thin\big(X(3)\big)\ra X(2),~~~~f(x):=\Big\{\Avg\big(x\backslash V(x)\big),V(x)\Big\}=\left\{{x_1+x_2\over 2},x_3\right\}.
\end{equation*}
\end{lmm}
\begin{proof}
Let $x,y\in X(3)$ be thin sets. Observe that $d_H(f(x),x)\leq{1\over 2}\delta(x)$, and so
\begin{equation*}
\textstyle d_H(f(x),f(y))\leq {1\over 2}\delta(x)+{1\over 2}\delta(y)+d_H(x,y).
\end{equation*}
Thus, if $\delta(x)\leq 2d_H(x,y)$ and $\delta(y)\leq 2d_H(x,y)$, then $d_H(f(x),f(y))\leq 3d_H(x,y)$. So, assume
\begin{equation*}
\textstyle d_H(x,y)<{1\over 2}\delta(x)~~~~\txt{or}~~~~ d_H(x,y)<{1\over 2}\delta(y),~~~~~~~~\left(\Ra~~d_H(x,y)<{1\over 2}{1\over 3}={1\over 6}\right).
\end{equation*}
Then by Lemma \ref{HausDistBound}, we have a proximal bijection $x_i\mapsto y(i)$ such that $\|x_i-y(i)\|\leq d_H(x,y)<1/6$ for all $i$. This bijection satisfies $\{x_1,x_2\}\mapsto\{x(1),x(2)\}=\{y_1,y_2\}$, i.e., $x_3\mapsto y(3)=y_3$, since
\begin{equation*}
\textstyle\|y(1)-y(2)\|\leq\|x_1-y(1)\|+\|x_1-x_2\|+\|x_2-y(2)\|<{2\over 3}\leq\|y_2-y_3\|.
\end{equation*}
Hence,
\begin{equation*}
\textstyle d_H(f(x),f(y))\leq \max\left(\left\|{x_1+x_2\over 2}-{y_1+y_2\over 2}\right\|,\|x_3-y_3\|\right)\leq d_H(x,y).\qedhere
\end{equation*}
\end{proof}

\begin{dfn}[\textcolor{blue}{\index{Lipschitz! partition of unity}{Lipschitz partition of unity}}]\label{LipPoU}
Consider the maps $\vphi_1,\vphi_2:\Real\ra\Real$ given by
\begin{equation*}
\vphi_1(t):=
\left\{
  \begin{array}{ll}
    1, & t\leq {1\over 5} \\
    -20t+5, & t\in \left[{1\over 5},{1\over 4}\right]\\
    0, & t\geq {1\over 4}
  \end{array}
\right\},~~~~
\vphi_2(t):=
\left\{
  \begin{array}{ll}
    0, & t\leq {1\over 5}\\
    20t-4, & t\in\left[{1\over 5},{1\over 4}\right]\\
    1, & t\geq {1\over 4}
  \end{array}
\right\}.
\end{equation*}
\end{dfn}
The above maps form a $20$-Lipschitz partition of unity.

\begin{lmm}[\textcolor{OliveGreen}{Gluing with strips}]\label{LSGluLmm}
Let $X$ be a metric space and $\vphi:X\ra\Real$ a Lipschitz function. Consider a finite collection of intervals $\{I_k=(a_k,b_k):k=1,...,N\}$ such that $a_k<a_{k+1}<b_k<b_{k+1}$ for all $k=1,...,N-1$ and $\Real=\bigcup_{k=1}^N(a_k,b_k)$.

If a map $g:E\subset X\ra X$ satisfies $\sup_{x\in E}d(x,g(x))\leq D<\infty$ and is Lipschitz on each of the sets $E_k:=\vphi^{-1}(a_k,b_k)=\{x\in E:a_k<\vphi(x)<b_k\}$, then $g$ is Lipschitz, and
\begin{equation*}
\textstyle\Lip(g)=\max\left\{\max\limits_{1\leq k\leq N}\Lip(g|_{E_k}),{\left(1+{2D\Lip(\vphi)\over\vep}\right)}\right\},
\end{equation*}
where ~$\vep:=\min\limits_k\diam(I_k\cap I_{k+1})=\min\limits_{1\leq k\leq N-1}|a_{k+1}-b_k|$.
\end{lmm}
\begin{proof}
Pick any $x,y\in E$. If $|\vphi(x)-\vphi(y)|<\vep$, then $x,y\in E_k=\vphi^{-1}(a_k,b_k)$ for some $k$, and so $d(g(x),g(y))\leq \Lip(g|_{E_k})d(x,y)$. On the other hand, if $|\vphi(x)-\vphi(y)|\geq\vep$, then by the triangle inequality, we get
\begin{equation*}
\textstyle d(g(x),g(y))\leq \left(1+{2D\Lip(\vphi)\over\vep}\right)d(x,y).\qedhere
\end{equation*}
\end{proof}

\begin{dfn}[\textcolor{blue}{\index{Interpolation map}{Interpolation map of $X(3)$}}]\label{IntMapDfnV3}
This is the map $R:N(X_0(3))\ra X(2)$,
\begin{equation*}
R(x):=\vphi_1\big(\delta(x)\big)R_1(x)+\vphi_2\big(\delta(x)\big)R_2(x),
\end{equation*}
where $\vphi_1,\vphi_2$ are as in Definition \ref{LipPoU},
\begin{equation*}
\textstyle R_1(x):=f(x)=\left\{{x_1+x_2\over 2},x_3\right\},~~~~R_2(x):=\Avg(x)={x_1+x_2+x_3\over 3},
\end{equation*}
and, we add and multiply sets by scalars as in Definition \ref{SumSclMult}.
\end{dfn}

\begin{lmm}\label{LipIntLmm}
The interpolation map $R:N(X_0(3))\ra X(2)$ is $121$-Lipschitz.
\end{lmm}
\begin{proof}
By Lemma \ref{LSGluLmm} with $\vphi=\delta$, $g=R$, and $\{(a_k,b_k)\}$ = $\{(-\infty,1/5)$, $(1/6,1/3)$, $(1/4,+\infty)\}$, it suffices to show that $R$ is Lipschitz on each of the following sets (within $N(X_0(3))$).
\begin{equation*}
\textstyle A:=\left\{\delta\leq {1\over 5}\right\},~~~~C:=\left\{{1\over 6}\leq\delta\leq{1\over 3}\right\} ,~~~~B:=\left\{\delta\geq{1\over 4}\right\}.
\end{equation*}
\bit[leftmargin=0.2cm]
\item[]\ul{$R$ is $3$-Lipschitz on $A$}: This follows from Lemma \ref{LipCenLmm1}.
\item[]\ul{$R$ is $9$-Lipschitz on $B$}: If $x,y\in B$, then $1\leq 4\delta(x),4\delta(y)$. Consider cases as follows. If $\delta(x)>2d_H(x,y)$ or $\delta(y)>2d_H(x,y)$, we have a proximal bijection $x\ra y$, and so by Lemma \ref{ContAvgLmm}(i),
\begin{equation*}
d_H(R(x),R(y))\leq d_H(x,y).
\end{equation*}
On the other hand, if $\delta(x),\delta(y)\leq 2d_H(x,y)$, then by Lemma \ref{ContAvgLmm}(ii),
\begin{equation*}
\textstyle d_H(R(x),R(y))\leq\left(1+{2\over 1/4}\right)d_H(x,y)=9d_H(x,y).
\end{equation*}
\item[]\ul{$R$ is $44$-Lipschitz on $C$}: This follows from the fact that $\vphi_i(\delta(x))$, $R_i(x)$, $i=1,2$ are bounded Lipschitz maps. With $\delta=\delta(x)$ and $\delta'=\delta(y)$, we have
{\small\begin{align*}
d_H(R(x)&,R(y))=d_H\Big(\vphi_1(\delta)R_1(x)+\vphi_2(\delta)R_2(x)~,~\vphi_1(\delta')R_1(y)+\vphi_2(\delta')R_2(y)\Big)\\
    &\leq d_H\Big(\vphi_1(\delta)R_1(x)~,~\vphi_1(\delta')R_1(y)\Big)+\Big\|\vphi_2(\delta)R_2(x)-\vphi_2(\delta')R_2(y)\Big\|\\
    &\leq \big[(20+3)+(20+1)\big]d_H(x,y)=44d_H(x,y).
\end{align*}}
The Lipschitz constant of $R$ can be calculated from Lemma \ref{LSGluLmm} as
\begin{equation*}
\textstyle \Lip(R)=\max\left\{44,1+{2\times 1\times 2\over 1/5-1/6}\right\}=121.\qedhere
\end{equation*}
\eit
\end{proof}

\begin{thm}\label{LrThmV3}
Let $X$ be a normed space. There exist Lipschitz retractions $X(n)\ra X(n-1)$ for $n=2,3$. From Lemmas \ref{HomExtLmm} and \ref{LipIntLmm}, the retraction $X(3)\ra X(2)$ has Lipschitz constant $6(121)+5=731$.
\end{thm}
\begin{proof}
The case of $n=2$ is given by Proposition \ref{LrPrpV2}. So, let $n=3$. Then by Lemma \ref{HomExtLmm}, it is enough to show that the interpolation map $R:N(X_0(3))\ra X(2)$ is Lipschitz, which follows from Lemma \ref{LipIntLmm}.
\end{proof}

\section{Lipschitz retractions $X(n)\ra X$ and $X(n)\ra X(2)$}\label{LipRetVn12}

If not stated, we will assume $X$ is a normed space. Let $\K_n(X)$:=$\{$convex compact subsets of $X$ of dimension $\leq n$$\}$. By Theorem 1.2 of \cite{shvartsman2004}, there exists an affine Lipschitz selector $S:\K_n(X)\ra X,~A\mapsto S(A)\in A$.
\begin{prp}\label{SvartLRVn1}
If $X$ is a normed space, there exist affine Lipschitz retractions $X(n)\ra X$ for all $n\geq 1$.
\end{prp}
\begin{proof}
Consider the map $s=S\circ\Conv:X(n)\sr{\Conv}{\ral}\K_n(X)\sr{S}{\ral}X$, where the convex hull operation is
\begin{equation*}
\textstyle \Conv(x):=\left\{\sum_{i=1}^n\al_ix_i:~\sum_{i=1}^n\al_i=1,~\al_i\geq0\right\}.
\end{equation*}
Given $\sum\al_ix_i\in \Conv(x)$, it follows from the definition of $d_H(x,y)$ that for each $x_i\in x$, there exists $y(i)\in y$ such that $\|x_i-y(i)\|\leq d_H(x,y)$. Since $\sum\al_iy(i)\in\Conv(y)$ and $\left\|\sum\al_i x_i-\sum\al_iy(i)\right\|\leq d_H(x,y)$, it follows by symmetry in the definition of Hausdorff distance that $d_H(\Conv(x),\Conv(y))\leq d_H(x,y)$. Hence,
\begin{align*}
d_H(s(x),s(y))&=\|s(x)-s(y)\|=\|S\circ \Conv(x)-S\circ \Conv(y)\|\\
 &\leq\Lip(S)d_H(x,y).\qedhere
\end{align*}
\end{proof}

\begin{dfn}[\textcolor{blue}{\index{Two-cluster decomposition}{Two-cluster decomposition} of an element of $X(n)$}]\label{TwoClustDec}
Let $x\in X(n)$ and consider numbers $\al,\beta>0$. A decomposition $x=x'\cup x''$ is an $(\al,\beta)$-decomposition if $x',x''$ are nonempty, $\diam(x')\leq\al$, $\diam(x'')\leq\al$, and $\dist(x',x'')\geq\beta$.
\end{dfn}

\begin{rmk}\label{UniqDecRmk}
Observe that an $(\al,\beta)$-decomposition $x=x'\cup x''$ is unique (up to permutation of the clusters $x',x''$) if $\al<\beta$. This is because if $x=\wt{x}'\cup\wt{x}''$ is any $(\al,\beta)$-decomposition, then neither $\wt{x}'$ nor $\wt{x}''$ can intersect both $x'$ and $x''$.

Moreover,  if $\alpha\le \alpha'<\beta' \le \beta$, then the $(\alpha', \beta')$-decomposition is the same as the $(\alpha,\beta)$-decomposition. This is because the $(\alpha,\beta)$-decomposition is also an $(\alpha', \beta')$-decomposition, which is unique. In particular, if $x$ has an $(\al,\beta)$-decomposition with $\al<\beta$, then for any number $0<c<{\beta-\al\over 2}$, any $(\al+c,\beta-c)$-decomposition is unique and equals the $(\al,\beta)$-decomposition.
\end{rmk}

Fix a number $\tau>6$.

\begin{dfn}[\textcolor{blue}{\index{Second order thin sets}{2nd order thin sets} in $X(n)$, Collection of thin sets}]
Let $X$ be a metric space and $x\in X(n)$. We say $x$ is a thin set of order $2$ (or $2$-thin set) if $x$ is normalized  (i.e., $\diam(x)=1$) and ~{\small $\dist_H(x,X(2)):=\inf\limits_{z\in X(2)}d_H(x,z)<{1\over\tau}$}. We will denote the collection of all 2-thin sets in $X(n)$ by $\Thin_2\big(X(n)\big)$.
\end{dfn}
\begin{lmm}[\textcolor{OliveGreen}{Cluster decomposition of a $2$-thin set}]
Let $X$ be a normed space and $x\in\Thin_2(X(n))$. Then $x$ admits a unique $\left({2\over\tau},1-{4\over\tau}\right)$-decomposition.
\end{lmm}
\begin{proof}
(i) \ul{Existence}: Since $\dist_H(x,X(2))<{1\over\tau}$, there exists $\{a,b\}\in X(2)$ such that $d_H(x,\{a,b\})<{1\over\tau}$. Since diameter is 2-Lipschitz with respect to $d_H$, we have $|\diam(x)-\diam(\{a,b\})|\leq 2d_H(x,\{a,b\})$. Thus,
\begin{equation*}
\textstyle \|a-b\|>1-{2\over\tau}.
\end{equation*}
Observe that for any $u\in x$, we have either $\|u-a\|<{1\over\tau}$ or $\|u-b\|<{1\over\tau}$ but not both: Otherwise, if  $\|u-a\|<{1\over\tau}$ and $\|u-b\|<{1\over\tau}$ then the triangle inequality gives $1-{2\over\tau}<\|a-b\|<{2\over\tau}$, which implies $\tau<4$ (a contradiction since $\tau>6$ by assumption). Let
\begin{align*}
x'&:=\{u\in x:\|u-a\|<1/\tau\}=x\cap N_{1/\tau}(a),\\
x''&:=\{u\in x:\|u-b\|<1/\tau\}=x\cap N_{1/\tau}(b).
\end{align*}
Note that for any $u\in x'$, $v\in x''$, we have $\big|\|u-v\|-\|a-b\|\big|<{2\over\tau}$, and so $\|u-v\|>1-{4\over\tau}$. Hence,
\begin{equation}
\label{2ThinDefEq}\textstyle \diam(x')<{2\over\tau},~~~~\diam(x'')<{2\over\tau},~~~~\dist(x',x'')\geq 1-{4\over\tau},
\end{equation}
where $x',x''$ are nonempty because $\diam(x)=1>{2\over\tau}$.

(ii) \ul{Uniqueness}: This follows from Remark \ref{UniqDecRmk} and the fact that $\tau>6$.
\end{proof}

\begin{dfn}[\textcolor{blue}{\index{Skeleton map}{Skeleton map}}]\label{SkelMapDfn}
This is the map $J:N\big(X(n)\big)\ra X(2)$ with
\begin{equation}
\label{SkelMapEq}J(x):=\left\{
       \begin{array}{ll}
         R_1(x):=\big\{s(x')~,~s(x'')\big\}, & x\in \Thin_2\big(X(n)\big) \\
       R_2(x):=\big\{s(x)\big\}, & x\in N\big(X(n)\big)\backslash \Thin_2\big(X(n)\big)\\
       \end{array}
     \right\},
\end{equation}
where $x\in\Thin_2\big(X(n)\big)$ decomposes as $x=x'\cup x''$, and $s:X(n)\ra X$ are the affine Lipschitz retractions from Proposition \ref{SvartLRVn1}.
\end{dfn}

\begin{lmm}\label{FirstLipLmm}
The map $R_1:\Thin_2\big(X(n)\big)\ra X(2)$ is Lipschitz.
\end{lmm}
\begin{proof}
If $x\in\Thin_2(X(n))$, let $x=x'\cup x''$ be the $(\al,\beta)=(2/\tau,1-4/\tau)$-decomposition of $x$. Pick a number $0<\rho<{\beta-\al\over 4}={1\over 4}\left(1-{6\over\tau}\right)$. Let $y\in\Thin_2(X(n))$ such that $d_H(x,y)\leq\rho$, and define
\begin{align*}
y'&:=\{u\in y:\dist(u,x')\leq\rho\}=y\cap N_\rho(x'),\\
y''&:=\{u\in y:\dist(u,x'')\leq\rho\}=y\cap N_\rho(x'').
\end{align*}
Observe that $y=y'\cup y''$ (by the definition of Hausdorff distance), $\diam(y')<\al+2\rho$, $\diam(y'')<\al+2\rho$, and $\dist(y',y'')>\beta-2\rho$. Thus, $y',y''$ give a unique $(\al+2\rho,\beta-2\rho)$-decomposition since $\al+2\rho<\beta-2\rho$ (where $\al+2\rho<{1\over 2}-{1\over\tau}<1$, and so $y',y''$ are nonempty). By Remark \ref{UniqDecRmk}, this $(\al+2\rho,\al-2\rho)$-decomposition of $y$ is the same as the $(\al,\beta)$-decomposition of $y$.

By construction, $x'\cup y'$ and $x''\cup y''$ are \emph{distantly separated} in the sense that
\begin{equation*}
\dist(x'\cup y',x''\cup y'')>\beta-2\rho>\al+2\rho>\max\left(\diam(x'\cup y'),\diam(x''\cup y'')\right).
\end{equation*}
It follows by direct calculation that $d_H(x,y)=\max\left(d_H(x',y'),d_H(x'',y'')\right)$. Hence, we have
\begin{align*}
d_H(R_1(x),R_1(y))&=d_H\big(\{s(x'),s(x'')\},\{s(y'),s(y'')\}\big)\leq \Lip(s)\max\Big(d_H(x',y'),d_H(x'',y'')\Big)\\
  &\leq \Lip(s)d_H(x,y),~~~~\txt{if}~~~~d_H(x,y)\leq\rho.
\end{align*}

On the other hand, since $R_1(x)\subset \Conv(x)$, we also have
\begin{align*}
d_H(R_1(x)&,R_1(y))\leq d_H(x,y)+d_H(R_1(x),x)+d_H(R_1(y),y)\\
  &\leq d_H(x,y)+\diam(x)+\diam(y)=d_H(x,y)+2\\
  &\textstyle\leq \left(1+{2\over\rho}\right)d_H(x,y),~~~~\txt{if}~~~~d_H(x,y)\geq\rho.\qedhere
\end{align*}
\end{proof}

\begin{dfn}[\textcolor{blue}{\index{Lipschitz! partition of unity}{Lipschitz partition of unity}}]\label{LipPoUVn2}
Fix $\tau>0$. The functions $\vphi_1,\vphi_2:\Real\ra\Real$ given by
\begin{equation*}
\vphi_1(t):=
\left\{
  \begin{array}{ll}
    1, & t\leq{1\over 3\tau} \\
    -(6\tau)t+3, & t\in\left[{1\over 3\tau},{1\over 2\tau}\right]\\
    0, & t\geq {1\over 2\tau}
  \end{array}
\right\},~~~~
\vphi_2(t):=
\left\{
  \begin{array}{ll}
    0, & t\leq{1\over 3\tau} \\
    (6\tau)t-2, & t\in\left[{1\over 3\tau},{1\over 2\tau}\right]\\
    1, & t\geq {1\over 2\tau}
  \end{array}
\right\}
\end{equation*}
form a Lipschitz partition of unity.
\end{dfn}

\begin{dfn}[\textcolor{blue}{Interpolation map}]
This is the map $R:N(X_0(n))\ra X(2)$ given by
\begin{equation*}
R(x):=\vphi_1(\delta)R_1(x)+\vphi_2(\delta)R_2(x),~~~~\delta:=\dist_H\big(x,X(2)\big),
\end{equation*}
where $R_1,R_2$ are as in (\ref{SkelMapEq}), $\vphi_1,\vphi_2$ are as in Definition \ref{LipPoUVn2} and, we add and multiply sets by scalars as in Definition \ref{SumSclMult}.
\end{dfn}

\begin{lmm}\label{LipIntLmmVn2}
The interpolation map $R:N(X_0(n))\ra X(2)$ is Lipschitz.
\end{lmm}
\begin{proof}
By Lemma \ref{LSGluLmm} with $\vphi=\delta$, $g=R$, $\{(a_k,b_k)\}$ = $\{(-\infty,1/(3\tau))$, $(1/(4\tau),1/\tau)$, $(1/(2\tau),+\infty)\}$, it suffices to show that $R$ is Lipschitz on each of the following sets (within $N(X_0(n))$).
\begin{equation*}
\textstyle A:=\left\{\delta\leq {1\over 3\tau}\right\},~~~~C:=\left\{{1\over 4\tau}\leq\delta\leq {1\over\tau}\right\},~~~~B:=\left\{\delta\geq{1\over 2\tau}\right\}.
\end{equation*}
$R$ is Lipschitz on $A$ by Lemma \ref{FirstLipLmm}, and Lipschitz on $B$ by the definition of $s$. On $C$, with $\delta=\dist_H\big(x,X(2)\big)$ and $\delta'=\dist_H\big(y,X(2)\big)$, we have
\begin{align*}
d_H(R(x)&,R(y))=d_H\Big(\vphi_1(\delta)R_1(x)+\vphi_2(\delta)R_2(x),\vphi_1(\delta')R_1(y)+\vphi_2(\delta')R_2(y)\Big)\\
    &\leq d_H\Big(\vphi_1(\delta)R_1(x),\vphi_1(\delta')R_1(y)\Big)+\Big\|\vphi_2(\delta)R_2(x)-\vphi_2(\delta')R_2(y)\Big\|.
\end{align*}
The result now follows because $\vphi_i(\delta)$, $R_i\big(x\big)$, $i=1,2$ are bounded Lipschitz maps.
\end{proof}

\begin{thm}\label{LrThmVn2}
Let $X$ be a normed space. There exist Lipschitz retractions $X(n)\ra X$ and $X(n)\ra X(2)$.
\end{thm}
\begin{proof}
The case of $X(n)\ra X$ is Proposition \ref{SvartLRVn1}. So, consider the case of $X(n)\ra X(2)$. By Lemma \ref{HomExtLmm}, it is enough to show that the interpolation map $R:N(X_0(n))\ra X(2)$ is Lipschitz, which follows from Lemma \ref{LipIntLmmVn2}.
\end{proof}

\section{Concrete H{\"o}lder retractions $X(n)\ra X(n-1)$}\label{Retractions}
In this section, unless stated otherwise, $X$ is a real normed space. As usual, we denote by $X^\ast$ the set of continuous \emph{linear functions/functionals} $x^\ast:X\ra\Real$ as a normed space (called \emph{dual space} of $X$) with norm $\|x^\ast\|:=\sup_{\|x\|\leq 1}|x^\ast(x)|=\sup_{\|x\|=1}|x^\ast(x)|$. If $x\in X$ and $x^\ast\in X^\ast$, we will sometimes write the number $x^\ast(x)$ as $\langle x,x^\ast\rangle$ for convenience.

\begin{dfn}[\textcolor{blue}{\index{Norming functional}{Norming functional}}] A linear functional $x^\ast\in X^\ast$ is a norming functional for $x_0\in X$ if $\|x^\ast\|=1$ and $x^\ast(x_0)=\|x_0\|$. If $x^\ast$ is a norming functional of $x_0$, we will also refer to $z^\ast:=\|x_0\|x^\ast$ as a norming functional of $x_0$. Note that $\|z^\ast\|=\|x_0\|$ and $z^\ast(x_0)=\|x_0\|^2$.
\end{dfn}
If $X$ is a normed space, then by the Hahn-Banach theorem, every $x_0\in X$ has a norming functional.

\begin{dfn}[\textcolor{blue}{\index{Fr{\'e}chet-G{\^a}teaux derivative}{Fr{\'e}chet-G{\^a}teaux derivative}}] Let $X,Y$ be normed spaces. A map $F:A\subset X\ra Y$ is (Fr{\'e}chet-) differentiable at $x\in A$ if there exists a linear map $dF_x:X\ra Y$ and a continuous map $o_x\in C(X,Y)$ such that
\begin{equation*}
\textstyle F(x+h)=F(x)+dF_xh+o_x(h)~~~~\txt{for all}~~h\in X,~~~~\txt{with}~~~~\lim\limits_{\|h\|\ra0}{\|o_x(h)\|\over\|h\|}=0.
\end{equation*}
The map $dF:X\ra L(X,Y)$, $x\mapsto  dF_x$ is called the (Fr{\'e}chet) derivative of $F$, and the linear map $dF_x:X\ra Y$ is called the (Fr{\'e}chet) derivative of $F$ at $x$.

    When the limit is required to exist only ``linearly'' (i.e., in one direction at a time), we get a weaker (G{\^a}teaux) version of the derivative: $F:A\subset X\ra Y$ is G{\^a}teaux-differentiable at $x\in A$ if there exists a linear map $DF_x:X\ra Y$ and a continuous map $o_x\in C(X,Y)$ such that for every $h\in X$ with $\|h\|=1$,
\begin{equation*}
\textstyle F(x+th)=F(x)+~t~DF_xh+o_x(th),~~~~\txt{for all}~~t\in \Real,~~~~\txt{with}~~\lim\limits_{t\ra0}{\|o_x(th)\|\over|t|}=0.
\end{equation*}
The map $DF:X\ra L(X,Y)$, $x\mapsto DF_x$ is called the G{\^a}teaux derivative of $F$, and the map $D_hF:X\ra Y$, $x\mapsto DF_xh$ is called the directional derivative of $F$ along $h$. Accordingly, the linear map $DF_x:X\ra Y$ is called the G{\^a}teaux derivative of $F$ at $x$, and the vector $DF_xh\in Y$ is called the directional derivative of $F$ at $x$ along $h$.
\end{dfn}

\begin{rmk*}[\textcolor{OliveGreen}{\index{Chain! rule}{Chain rule}}] Let ~$A\subset X\sr{F}{\ral}Y$~ and ~$B\subset Y\sr{G}{\ral}Z$~ be maps of normed spaces such that $F$ is differentiable at $x_0\in A$ and $G$ is differentiable at $y_0:=F(x_0)\in B$. Then $A\subset X\sr{G\circ F}{\ral}Z$ is differentiable at $x_0$, and~ $d(G\circ F)_{x_0}=dG_{F(x_0)}\circ dF_{x_0}:X\sr{dF_{x_0}}{\ral}Y\sr{dG_{F(x_0)}}{\ral}Z$.

In particular, if ~$[0,1]\sr{u}{\ral}X\sr{\|\cdot\|}{\ral}\Real$,~ where $(X,\|\cdot\|)$ is a normed space with a differentiable norm and $u$ is a $C^1$-smooth path, then~ ${d\|u(t)\|\over dt}=d\|\cdot\|_{u(t)}\left({du(t)\over dt}\right)$~ for all $t\in[0,1]$.
\end{rmk*}
\begin{proof}
See Remark \ref{ChainRule}.
\end{proof}

\begin{rmk*}[\textcolor{OliveGreen}{\index{Mean value theorem}{Mean value theorem}: \cite[Theorem 1.8, p.13]{ambro-prodi1993}}] If $F:O\subset X\ra Y$ is G{\^a}teaux-differentiable and $O$ is open, then for any $x_1,x_2\in X$ such that ~$[x_1,x_2]:=\big\{\eta(t)=(1-t)x_1+tx_2:t\in[0,1]\big\}\subset O$,~ we have
\begin{equation*}
\textstyle\|F(x_1)-F(x_2)\|\leq C(x_1,x_2)\|x_1-x_2\|,~~\txt{where}~~C(x_1,x_2):=\sup\limits_{x\in[x_1,x_2]}\left\|DF_x\right\|.
\end{equation*}
\end{rmk*}
\begin{proof}
See Remark \ref{MeanValThm}.
\end{proof}

\begin{dfn}[\textcolor{blue}{\index{Semi-inner products}{Semi-inner products}}]\label{SemiIPs} Let $X$ be a normed space and $x,y\in X$. We define
\begin{equation*}
\textstyle\langle x,y\rangle_-:=\inf\limits_{y^\ast\in\F y}\langle x,y^\ast\rangle,~~~~\langle x,y\rangle_+:=\sup\limits_{y^\ast\in\F y}\langle x,y^\ast\rangle,
\end{equation*}
where $\F:X\ra\P(X^\ast)$ is the duality map of $X$, given by the set of norming functionals
\begin{equation*}
\F x~:=~\left\{x^\ast\in X^\ast:\|x^\ast\|=\|x\|,~x^\ast(x)=\|x\|^2\right\},~~~~\txt{for all}~~x\in X.
\end{equation*}
\end{dfn}

\begin{prp}[\textcolor{OliveGreen}{Semi-inner products as derivatives: \cite[Proposition 12.3(d), p.115]{deimling}}]\label{SemiInnCh}
Let $X$ be a normed space. The semi-inner products in Definition \ref{SemiIPs} are determined by one-sided derivatives of the norm as follows.
\begin{equation}
\label{SemiInnEq}\textstyle\langle x,y\rangle_-=\|y\|\lim\limits_{t\uparrow 0}{\|y+tx\|-\|y\|\over t},~~~~\langle x,y\rangle_+=\|y\|\lim\limits_{t\downarrow 0}{\|y+tx\|-\|y\|\over t}.
\end{equation}
\end{prp}
\begin{proof}
See Lemma \ref{SemiInnAp}.
\end{proof}

\begin{lmm}[\textcolor{OliveGreen}{Derivative of the norm along trajectories}]\label{AbsCont}
If $X$ is a normed space and $\gamma:[0,1]\ra X$ is a $C^1$-smooth path, then the following are true.
\bit[leftmargin=0.9cm]
\item[(i)] The function $\vphi:[0,1]\ra\Real$, $\vphi(t)=\|\gamma(t)\|$ is absolutely continuous, i.e.,
\begin{equation}
\label{AbsContEq}\textstyle\vphi' ~~~~\txt{exists a.e.},~~~~\vphi'\in L([0,1]), ~~~~\txt{and}~~~~ \vphi(t)=\vphi(0)+\int_0^t\vphi'(s)ds.
\end{equation}
\item[(ii)] With $\gamma'(t)={d\over dt}\gamma(t):=d\gamma_t$, the derivative of $\vphi$ can be expressed in the following form:
\begin{equation}
\label{LinTrans}\textstyle\vphi'(t)=\lim\limits_{h\downarrow 0}{\|\gamma(t)+h\gamma'(t)\|-\|\gamma(t)\|\over h}~\sr{(\ref{SemiInnEq})}{=}~{1\over\|\gamma(t)\|}\big\langle\gamma'(t),\gamma(t)\big\rangle_+,~~~~\txt{for a.e.}~~t\in[0,1].
\end{equation}
\eit
\end{lmm}
\begin{proof}
(i) Let $C:=\sup_{[0,1]}\|\gamma'\|$, $\gamma'$ the derivative of $\gamma$. By the mean value theorem,
\begin{equation*}
\big|\vphi(a)-\vphi(b)\big|=\big|\|\gamma(a)\|-\|\gamma(b)\|\big|\leq\big\|\gamma(a)-\gamma(b)\big\|\leq C|a-b|,
\end{equation*}
which shows $\vphi$ is absolutely continuous.

(ii) If $\vphi$ is differentiable at $t\in[0,1]$, then
\begin{equation*}
\textstyle {d\over dt}\|\gamma(t)\|=\lim\limits_{h\ra0}{\|\gamma(t+h)\|-\|\gamma(t)\|\over h}~\sr{(s)}{=}~\lim\limits_{h\ra 0}{\|\gamma(t)+h\gamma'(t)\|-\|\gamma(t)\|\over h},
\end{equation*}
where step (s) holds because $\gamma(t+h)=\gamma(t)+h\gamma'(t)+o_t(h)$, and so
\begin{equation*}
\textstyle \lim\limits_{h\ra 0}\left|{\|\gamma(t+h)\|-\|\gamma(t)+h\gamma'(t)\|\over h}\right|\leq \lim\limits_{h\ra 0}{\|o_t(h)\|\over|h|}=0.\qedhere
\end{equation*}
\end{proof}

\begin{lmm}[\textcolor{OliveGreen}{Semi-monotonicity of the radial projection}]\label{accret-lmm}
If $X$ is a normed space, the map $X\backslash\{0\}\ra X$ given by $x\mapsto \hat{x}:={x\over\|x\|}$ satisfies
\begin{equation}
\label{AccretEq}\langle\hat{x}-\hat{y},x-y\rangle_-\geq0~~~~\txt{for all}~~~~x,y\in X\backslash\{0\}.
\end{equation}
\end{lmm}
\begin{proof}
Fix any two vectors $x,y\in X\backslash\{0\}$. If $\|x\|=\|y\|$, then $\langle\hat{x}-\hat{y},x-y\rangle_-=\|x\|^{-1}\langle x-y,x-y\rangle_-\geq0$. So, assume $\|x\|>\|y\|$. Consider the convex function
\begin{equation*}
\vphi(t)=\|x-y+t(\hat{x}-\hat{y})\|=\big\|(t+\|x\|)\hat{x}-(t+\|y\|)\hat{y}\big\|,~~~~\txt{for}~~t\in\Real.
\end{equation*}
By Proposition \ref{SemiInnCh}, it suffices to show the left-sided derivative of $\vphi$ is nonnegative at $t=0$.

Observe that $\vphi(-\|x\|)=\vphi(-\|y\|)= \|x\|-\|y\|$. Since $\vphi$ is convex and $\|x\|\neq \|y\|$, it follows that $\vphi$  attains its minimum on $\big[-\|x\|, -\|y\|\big]$. So, $\vphi$ is nondecreasing on $\big[-\|y\|,+\infty\big)$. Hence, both one-sided derivatives of $\vphi$ are nonnegative at $t=0$.
\end{proof}

\begin{thm}[\textcolor{OliveGreen}{Analog of Theorem 1.1 in \cite{kovalev2016}}]\label{MainThm}
If $X$ is a normed space, then for each $n\geq 2$ there exists a retraction ~$r:X(n)\ra X(n-1)$ that is H{\"o}lder-continuous on bounded subsets of $X(n)$.
\end{thm}
\begin{proof}
We will proceed in six steps to construct the retraction and prove its continuity.

{\flushleft \ul{1. Evolution equation and collision time}}:~ Equip $X^n$ with the metric {\small $d\big(x,y\big)=\sum_{i=1}^n\|x_i-y_i\|$}, which makes $X^n$ a normed space. Let $D:=\{x\in X^n:x_i=x_j~\txt{for some}~i\neq j\}$, and consider the map
\begin{align*}
&J\textstyle=(J_1,...,J_n):X^n\backslash D\ra X^n,~~~~J_i(x):=\sum_{j\neq i}{x_i-x_j\over\|x_i-x_j\|},\\
&\|J\|\textstyle=\sum_i\|J_i\|=\sum_i\left\|\sum_{j\neq i}{x_i-x_j\over\|x_i-x_j\|}\right\|\leq n(n-1).
\end{align*}
Note that ~$X(n)\backslash X(n-1)=\big\{x=\{x_1,...,x_n\}:~(x_1,...,x_n)\in X^n\backslash D\big\}$. Given $(x_1,...,x_n)\in X^n\backslash D$, consider the system of ordinary differential equations
\begin{align}
\label{KOVeq2}&\textstyle{du_i(t)\over dt}=-J_i\big(u(t)\big),~~~~u_i(0)=x_i,~~~~i=1,...,n,\\
\label{KOVeq3}&\textstyle\left\|{du_i(t)\over dt}\right\|=\|J_i\big(u(t)\big)\|\leq n-1,~~~~i=1,...,n.
\end{align}
Beginning with $u(0)=x\in X^n\backslash D$, by Picard's theorem, the system (\ref{KOVeq2}) continues to have a unique solution
\begin{equation*}
u(t)\in X^n\backslash D,~~~~\txt{with each}~~~~u_i(t)\in\txt{Span}\{x_1,...,x_n\}~~\txt{in}~~X,
\end{equation*}
until we reach the set $D$, a situation we will refer to as ``\emph{collision}''. Let ~$T(x):=\sup\{t:~t\geq 0,~u(t)\in X^n\backslash D\}$,~ i.e., $[0,T(x))$ is the maximal interval of existence of the solution of (\ref{KOVeq2}).
\begin{rmk}\label{CommColl}
For any $0<\tau<T(x)$, $u(t+\tau)$ is the unique solution of the system
\begin{equation*}
\textstyle {dw_i(t)\over dt}=-J_i\big(w(t)\big),~~~~w_i(0)=u_i\big(\tau\big).
\end{equation*}
This implies the point of collision for $w$ is the same as for $u$, i.e., {\small $u\big(T(x)\big)=w\Big(T\big(u(\tau)\big)\Big)=u\Big(T\big(u(\tau)\big)+\tau\Big)$}. Equivalently, we have
\begin{equation}
\label{TimeTrans}T(x)=T\big(u(\tau)\big)+\tau,~~\txt{or}~~T\big(u(\tau)\big)=T(x)-\tau,~~~~\txt{for all}~~~~0<\tau<T(x).
\end{equation}
\end{rmk}

{\flushleft \ul{2. Bounds on the collision time $T(x)$}}:~ With $\delta$ as in Definition \ref{MinSep},
\begin{align*}
\delta(x)&=\min_{i\neq j}\|x_i-x_j\|\leq \min_{i\neq j}\Big(\|x_i-u_i(T(x))\|+\|u_i(T(x))-u_j(T(x))\|+\|u_j(T(x))-x_j\|\Big)\\
  &\sr{(s)}{\leq} (n-1)T(x)+\min_{i\neq j}\|u_i(T(x))-u_j(T(x))\|+(n-1)T(x)\\
  &=(n-1)T(x)+0+(n-1)T(x)=2(n-1)T(x),
\end{align*}
where the mean value theorem is used at step (s). Therefore,
\begin{equation*}
\textstyle T(x)\geq{\delta(x)\over 2(n-1)}.
\end{equation*}

Renumbering the points $x_i$, we may assume $\delta(x)=\|x_1-x_2\|$. Let $\vphi(t):=\|u_1(t)-u_2(t)\|$. Then by Lemma \ref{AbsCont}, $\vphi$ is absolutely continuous and, for all $t$, its derivative satisfies
\begin{align*}
\vphi'(t)~&\textstyle\sr{(\ref{LinTrans})}{=}\left\langle{du_1\over dt}-{du_2\over dt},{u_1-u_2\over\|u_1-u_2\|}\right\rangle_+\sr{(\ref{KOVeq2})}{=}
-\left\langle J_1(u)-J_2(u),{u_1-u_2\over \|u_1-u_2\|}\right\rangle_-\nn\\
   &\textstyle=-\left\langle\sum_{j\neq 1}{u_1-u_j\over\|u_1-u_j\|}-\sum_{j\neq 2}{u_2-u_j\over\|u_2-u_j\|},{u_1-u_2\over\|u_1-u_2\|}\right\rangle_-\nn\\
   &\textstyle~=-2-\sum_{j= 3}^n{\left\langle{u_1-u_j\over\|u_1-u_j\|}-{u_2-u_j\over\|u_2-u_j\|},(u_1-u_j)-(u_2-u_j)\right\rangle_-\over\|u_1-u_2\|}\nn\\
   &\textstyle\sr{(\ref{AccretEq})}{\leq} -2,~~~~\txt{for almost all}~~~~0<t<T(x),
\end{align*}
where step (\ref{AccretEq}) refers to the property of the radial projection proved in Lemma \ref{accret-lmm}. Upon integration of the above inequality, we get $\vphi\big(T(x)\big)-\vphi(0)\leq -2T(x)$, which implies
\begin{equation}
\label{KOVeq4}\textstyle T(x)~\leq~{\delta(x)\over 2}.
\end{equation}

{\flushleft \ul{3. Definition of the retraction}}:~ Define $r:X(n)\ra X(n-1)$ as follows. If $x\in X(n)\backslash X(n-1)$, let $r(x)=r\big(\{x_i\}\big):=\left\{u_i\big(T(x)\big)\right\}=u\big(T(x)\big)$, which is a well defined map since the order of enumeration is unimportant. If $x\in X(n-1)$, let $r(x):=x$. Then $r|_{X(n-1)}=id_{X(n-1)}$. It remains to show that $r$ is continuous. Specifically, we will show that for all $x,y\in X(n)$,
\begin{equation}
\label{KOVeq5}d_H\big(r(x),r(y)\big)~\leq~n(2n-1)\diam\left(x\cup y\right)^{1-{1\over 2n-1}}d_H\left(x,y\right)^{1\over 2n-1}.
\end{equation}
For $x,y\in X^n$, let $u,v$ be the solutions of (\ref{KOVeq2}) with initial data $u(0)=x$, $v(0)=y$. Recall that for all $x,y\in X(n)$,
\begin{equation*}
d_H\big(x,y\big):=\max\left\{\max_{1\leq i\leq n}\min_{1\leq j\leq n}\|x_i-y_j\|,\max_{1\leq i\leq n}\min_{1\leq j\leq n}\|x_j-y_i\|\right\}.
\end{equation*}

\begin{figure}[H]
\centering
\scalebox{1}{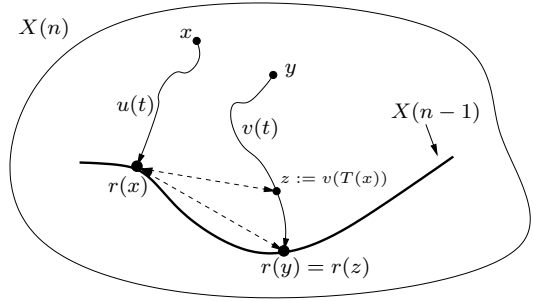} 
  \caption{~Evolution paths until collision occurs.}\label{dg5}
\end{figure}

{\flushleft \ul{4. Estimate of Hausdorff distance to the collision point}}:~ Using (\ref{KOVeq3}) and (\ref{KOVeq4}),
\begin{align*}
d_H&\big(r(x),x\big)=\max\left\{\max_i\min_j\|u_i(T(x))-x_j\|,\max_i\min_j\|u_j(T(x))-x_i\|\right\}\nn\\
   \leq&\max\left\{\max_i\min_j\Big((n-1)T(x)+\|x_i-x_j\|\Big),\max_i\min_j\Big((n-1)T(x)+\|x_j-x_i\|\Big)\right\}\nn\\
   =&\max\left\{\Big((n-1)T(x)+\max_i\min_j\|x_i-x_j\|\Big),\Big((n-1)T(x)+\max_i\min_j\|x_j-x_i\|\Big)\right\}\nn\\
   =&\textstyle\max\left\{\Big((n-1)T(x)+0\Big),\Big((n-1)T(x)+0\Big)\right\}=(n-1)T(x)\leq(n-1){\delta(x)\over 2},
\end{align*}
and a similar bound holds for $d(r(y),y)$. Therefore,

\begin{equation}
\label{KOVeq6}\textstyle d_H\big(r(x),x\big)\leq{n-1\over 2}\delta(x),~~~~d_H\big(r(y),y\big)\leq{n-1\over 2}\delta(y).
\end{equation}

{\flushleft \ul{5. Estimate of Hausdorff distance at first collision}}: By Lemma \ref{AbsCont}, the function $g(t):=\sum_{i=1}^n\|u_i(t)-v_i(t)\|$ is absolutely continuous and its derivative satisfies the following: For a.e. $0<t<T:=\min\Big(T(x),T(y)\Big)$,
{\small\begin{align*}
g'(t)\sr{(\ref{LinTrans})}{=}\textstyle &\textstyle\sum\limits_i\left\langle{du_i\over dt}-{dv_i\over dt},{u_i-v_i\over\|u_i-v_i\|}\right\rangle_+\sr{(\ref{KOVeq2})}{=}-\sum\limits_i\left\langle J_i(u)-J_i(v),{u_i-v_i\over\|u_i-v_i\|}\right\rangle_-,\\
|g'(t)|\leq&\textstyle \sum\limits_i\|J_i(u)-J_i(v)\|\leq \sum\limits_i\sum\limits_{j\neq i}\left\|{u_i-u_j\over\|u_i-u_j\|}-{v_i-v_j\over\|v_i-v_j\|}\right\|= 2\sum\limits_{i<j}\left\|{u_i-u_j\over\|u_i-u_j\|}-{v_i-v_j\over\|v_i-v_j\|}\right\|\\
\sr{(s)}{\leq}&\textstyle 4{\sum_{i<j}\left\|(u_i-u_j)-(v_i-v_j)\right\|\over\max(\|u_i-u_j\|,\|v_i-v_j\|)}
   \leq {4(n-1)\sum_i\|u_i-v_i\|\over\min\limits_{i<j}\Big(\max(\|u_i-u_j\|,\|v_i-v_j\|)\Big)}
   \leq {4(n-1)\over\max\Big(\delta\big(u(t)\big),\delta\big(v(t)\big)\Big)}g(t)\\
   \sr{(\ref{KOVeq4})}{\leq}&\textstyle {2(n-1)\over\max\Big(T\big(u(t)\big),T\big(v(t)\big)\Big)}g(t)
   \sr{(\ref{TimeTrans})}{=}{2(n-1)\over\max\Big(T(x),T(y)\Big)-t}g(t),\nn
\end{align*}}
where step (s) is due to the inequality $\left|{x\over\|x\|}-{y\over\|y\|}\right|\leq{2\|x-y\|\over\max(\|x\|,\|y\|)}$ from \cite{dunkl1964,thele1974}. It follows that
\begin{equation*}
\textstyle |g'(t)|\leq{2(n-1)\over T-t}g(t),~~~~\txt{for a.e.}~~0<t<T.
\end{equation*}
Since we also have ~$|g'(t)|\leq 2(n-1)$~ for a.e. $0<t<T$, it follows that ~$|g'(t)|\leq 2(n-1)\min\left({1\over T-t}g(t),1\right)$~ for a.e. $0<t<T$. We consider various cases as follows.
\bit[leftmargin=0.9cm]
\item[(i)]\ul{$g(\tau)=T-\tau$ for some $\tau\in(0,T)$}:~ The bound $g'(t)\leq{2(n-1)\over T-t}g(t)$ for a.e. implies
\begin{equation*}
\textstyle g(\tau)\leq \exp\left(\int_0^\tau{2(n-1)\over T-s}ds\right)g(0)=
\textstyle \left(\frac{T}{T-\tau}\right)^{2(n-1)} g(0)=\left(\frac{T}{g(\tau)}\right)^{2(n-1)} g(0),
\end{equation*}
which in turn implies
\begin{equation}
\label{gBound1} g(\tau)\leq T^{1-{1\over 2n-1}}g(0)^{1\over 2n-1}.
\end{equation}
Also, the bound $g'(t)\leq 2(n-1)$ a.e. implies ~$g(T)-g(\tau)\leq 2(n-1)(T-\tau) = 2(n-1)g(\tau)$,~ i.e.,
\begin{equation}
\label{gBound2}g(T)\leq (2n-1)g(\tau)\sr{(\ref{gBound1})}{\leq}(2n-1)T^{1-{1\over 2n-1}}g(0)^{1\over 2n-1}.
\end{equation}
\item[(ii)]\ul{$g(t)<T-t$ for all $t\in (0,T)$}:~ The bound $g'(t)\leq{2(n-1)\over T-t}g(t)$ a.e. shows (\ref{gBound1}) holds for all $t$, that is, for all $0<t<T$,
\begin{equation}
\label{gBound3}g(t)\leq T^{1-{1\over 2n-1}}g(0)^{1\over 2n-1}\leq (2n-1)T^{1-{1\over 2n-1}}g(0)^{1\over 2n-1}.
\end{equation}
\item[(iii)]\ul{$g(t)> T-t$ for all $t\in(0,T)$}:~  The bound $g'(t)\leq 2(n-1)$ a.e. implies ~$g(T)-g(t)\leq 2(n-1)(T-t)< 2(n-1)g(t)$,~ i.e., $g(T)< (2n-1)g(t)$, for all $0<t<T$. Therefore, in the limit $t\ra0$, we get
\begin{equation}
\label{gBound4}g(T)\leq (2n-1)g(0).
\end{equation}
\eit
Since $T\leq{\min\big(\delta(x),\delta(y)\big)\over 2}<\diam(x\cup y)$ and $\diam(x\cup y)\geq g(0)$, it follows that all three cases above imply
\begin{equation}
\label{gBound5}g(T)\leq(2n-1)~\diam(x\cup y)^{1-{1\over 2n-1}}~g(0)^{1\over 2n-1}~=~C_n(x,y)~g(0)^{1\over 2n-1},
\end{equation}
where $C_n(x,y):=(2n-1)~\diam(x\cup y)^{1-{1\over 2n-1}}$. Therefore,
{\small\begin{align*}
d_H\big(u(T)&,v(T)\big)\leq d\big(u(T),v(T)\big)=\max_i\|u_i(T)-v_i(T)\|\sr{(\ref{gBound5})}{\leq}~C_n(x,y)\max_i\|u_i(0)-v_i(0)\|^{1\over 2n-1},
\end{align*}}
which implies
\begin{equation}
\label{KOVeq7} d_H\big(u(T),v(T)\big)\leq~C_n(x,y)\max_i\|x_i-y_i\|^{1\over 2n-1}.
\end{equation}

{\flushleft \ul{6. Estimate of Hausdorff distance between collision points}}: Assume wlog that $T(x)\leq T(y)$. We consider two cases as follows, where $\rho:=d_H(x,y)$.
\bit[leftmargin=0.8cm]
\item \ul{Case 1}: $\delta(x)+\delta(y)\leq 4\rho$. In this case, we obtain (\ref{KOVeq5}) as follows.
\begin{align*}
d_H\big(&r(x),r(y)\big)\textstyle\leq d_H\big(r(x),x\big)+d_H(x,y)+d_H(y,r(y))\sr{(\ref{KOVeq6})}{\leq}~ {n-1\over 2}\delta(x)+d_H(x,y)+{n-1\over 2}\delta(y)\\
  &\leq \rho+2(n-1)\rho=(2n-1)\rho\leq n(2n-1)\diam\left(x\cup y\right)^{1-{1\over 2n-1}}d_H\left(x,y\right)^{1\over 2n-1}.
\end{align*}
\item \ul{Case 2}: $\delta(x)+\delta(y)>4\rho$. In this case, $\delta(x)>2\rho$ or $\delta(y)>2\rho$, and so (\ref{KOVeq8}) applies. By definition, $r(x)=u(T(x))$. Let $z:=v(T(x))$. Then $d_H\big(r(x),z\big)=d_H\Big(u\big(T(x)\big),v\big(T(x)\big)\Big)\sr{(\ref{KOVeq7})}{\leq}C_n(x,y)\max_i\|x_i-y_i\|^{1\over 2n-1}\sr{(\ref{KOVeq8})}{\leq} C_n(x,y)\rho^{1\over 2n-1}$. That is,
\begin{align}
\label{KOVeq9}d_H\big(r(x),z\big)\leq C_n(x,y)\rho^{1\over 2n-1}.
\end{align}
Since $\delta$ is 2-Lipschitz (Lemma \ref{DeltaCont}) and $\delta\big(r(x)\big)=0$, with $z_i:=v_i\big(T(x)\big)$,
\begin{equation*}
\delta(z)=|\delta\big(r(x)\big)-\delta(z)|\leq 2d_H\big(r(x),z\big)\sr{(\ref{KOVeq9})}{\leq}2C_n(x,y)\rho^{1\over 2n-1},
\end{equation*}
which together with (\ref{KOVeq6}) implies
\begin{equation}
\label{KOVeq10}\textstyle d_H\big(r(z),z\big)\sr{(\ref{KOVeq6})}{\leq}{n-1\over 2}\delta(z)\leq (n-1)C_n(x,y)\rho^{1\over 2n-1}.
\end{equation}
Using Remark \ref{CommColl} (which says $r(z)=r(y)$) at step (R\ref{CommColl}) below,
\begin{align*}
d_H\big(r(x),r(y)\big)&\sr{\txt{(R\ref{CommColl})}}{=}d_H\big(r(x),r(z)\big)\leq d_H\big(r(x),z\big)+d_H\big(r(z),z\big)\sr{(\ref{KOVeq9}),(\ref{KOVeq10})}{\leq} n C_n(x,y)\rho^{1\over 2n-1},
\end{align*}
which in turn implies (\ref{KOVeq5}), i.e.,
\begin{equation*}
d_H\big(r(x),r(y)\big)\leq n(2n-1)\diam\left(x\cup y\right)^{1-{1\over 2n-1}}d_H\left(x,y\right)^{1\over 2n-1}.
\end{equation*}
\eit
This completes the proof.
\end{proof}

\begin{rmk}[\textcolor{OliveGreen}{Connection with \cite{kovalev2016}}]\label{HilbConRmk}
When $X$ is a Hilbert space, the norm is differentiable. Thus, with the metric $d(x,y)=\left(\sum_i\|x_i-y_i\|^2\right)^{1/2}$ on $X^n$, the function $g(t)={1\over 2}\sum_i\|u_i(t)-v_i(t)\|^2$ satisfies the following (for all $0<t<T$):
{\small\begin{align*}
g'(t)&\textstyle =\sum\limits_i\left\langle{du_i\over dt}-{dv_i\over dt},u_i-v_i\right\rangle_+=-\sum_i\Big\langle J_i(u)-J_i(v),u_i-v_i\Big\rangle_-\\
  &\textstyle =-\sum\limits_{1\leq i<j\leq n}\left\langle{u_i-u_j\over\|u_i-u_j\|}-{v_i-v_j\over\|v_i-v_j\|},(u_i-v_i)^\ast-(u_j-v_j)^\ast\right\rangle\\
  &\textstyle \sr{(a)}{=}-\sum\limits_{1\leq i<j\leq n}\left\langle{u_i-u_j\over\|u_i-u_j\|}-{v_i-v_j\over\|v_i-v_j\|},(u_i-u_j)^\ast-(v_i-v_j)^\ast\right\rangle\sr{(b)}{\leq}0,
\end{align*}}where step (a) is due to linearity of the duality map (by Riesz representation theorem), and step (b) is due to monotonicity of the radial projection $F(x)=x/\|x\|$, $x\neq0$. This leads to $g(t)\leq g(0)$ for all $0<t<T$. Hence, as in \cite{kovalev2016}, the map $r:X(n)\ra X(n-1)$ is a Lipschitz retraction.
\end{rmk}

\begin{rmk}[\textcolor{blue}{Strict convexity and Lipschitz retractions}]
Let $X$ be a Banach space. The following are strong indications that Lipschitz retractions $X(n)\ra X(n-1)$, for all $n\geq 2$, may exist when $X$ is strictly convex (Definition \ref{StrictConvNS}).
\begin{enumerate}[leftmargin=0.8cm]
\item If $X$ is a geodesic space (resp. $\al$-quasiconvex space), then $X(n)$ is $2$-quasiconvex by Theorem \ref{QConvThm} (resp. $2\al$-quasiconvex by Corollary \ref{MetricQCQC}).
\item If $Y$ is a quasiconvex space, then a map of metric spaces $f:Y\ra Z$ is Lipschitz if and only if locally Lipschitz (Definition \ref{LocLipMap} and Lemma \ref{GluLmmII}). Thus, if $Y=X(n)$ and $Z=X(n-1)$, to obtain a Lipschitz retraction $X(n)\ra X(n-1)$, it is enough to find a locally Lipschtz retraction $X(n)\ra X(n-1)$.
\item A characterization in \cite{petryshyn1970} of strict convexity of Banach spaces in terms of the duality map indicates that if $X$ is strictly convex, the quantity $g(t):=\sum_i\|u_i(t)-v_i(t)\|$ (or a variant of it based on the particular norm considered on $X^n$) introduced in the proof of Theorem \ref{MainThm} may satisfy $g'(t)\leq0$ for all $t\in[0,1]$ provided the starting points $u(0)$ and $v(0)$ are close, i.e., $d_H\big(u(0),v(0)\big)$ is small. This behavior would yield the desired locally Lipschitz retraction $X(n)\ra X(n-1)$.
\end{enumerate}
\end{rmk}

\section{Invariance of the convex hull under the retractions}\label{Selections}
Let $X$ be a normed space. We have seen that there exist continuous retractions $r:X(n)\ra X(n-1)$ for each $n\geq 2$, where $r$ is an explicit map. Given $x\in X(n)$, let $\Conv(x)=\Conv\big(\{x_1,...,x_n\}\big)\subset X$ denote the convex hull of $x$ in $X$. We will show that $r(x)\subset\Conv(x)$ for all $x\in X(n)$. This will (as explained in Remark \ref{SvRmk}) lead to a result that is related to Shvartsman's theorem on Lipschitz selection in \cite[Theorem 1.2]{shvartsman2004}.

\begin{dfn}[\textcolor{blue}{Matrix exponential integral: \index{Peano-Baker series}{Peano-Baker series}}] If ~$M:\Real\ra Mat_{n\times n}(\Real)$, ~we define
\begin{align}\label{PBSeries}
\textstyle T e^{\int_0^tM(s)ds}&:=\textstyle I+\sum_{k=1}^\infty\int_{[0,t]^k}{TM^k(t_1,...,t_k)\over k!}dt_1\cdots dt_k\nn\\
 &\textstyle=I+\sum_{k=1}^\infty\int_{t>t_1>\cdots>t_k>0}M(t_1)\cdots M(t_k)dt_1\cdots dt_k,
\end{align}
where, with the characteristic function ~$\chi_{(0,\infty)}(t):=
\left\{
  \begin{array}{ll}
    0, & t\leq 0 \\
    1, & t>0
  \end{array}
\right\}
$,
\begin{align}\label{TOSeries}
&\textstyle TM^k(t_1,...,t_k):=\sum\limits_{\sigma\in S_k}~~\prod\limits_{j=1}^{k-1}\chi_{(0,\infty)}\left(t_{\sigma(j)}-t_{\sigma(j+1)}\right)~~M\left(t_{\sigma(1)}\right)M\left(t_{\sigma(2)}\right)\cdots M\left(t_{\sigma(k)}\right)\nn\\
&\textstyle~~~~=\sum\limits_{\sigma\in S_k}~\chi_{(0,\infty)}\left(t_{\sigma(1)}-t_{\sigma(2)}\right)\cdots\chi_{(0,\infty)}\left(t_{\sigma(k-1)}-t_{\sigma(k)}\right)~M\left(t_{\sigma(1)}\right)\cdots M\left(t_{\sigma(k)}\right).
\end{align}
(Note that nonzero terms of the sum in (\ref{TOSeries}) are those for which ~$t_{\sigma(1)}\geq t_{\sigma(2)}\geq\cdots\geq t_{\sigma(n)}$). In particular, $TM^2(t,t')=\chi_{(0,\infty)}(t-t')M(t)M(t')+\chi_{(0,\infty)}(t'-t)M(t')M(t)$.
\end{dfn}
The series (\ref{PBSeries}) is known as Peano-Baker series. By \cite[Theorem 1]{baake-schlaegel2011}, the series converges if the real-valued function $\vphi(t)=\|M(t)\|$ is locally integrable (i.e., integrable on bounded intervals).

In $X^n$, let $E:=\{x\in X^n:x_i\neq x_j~\txt{if}~i\neq j\}$. Consider the following system of equations in which $J_i:E\ra X$ is given by $J_i(x):=\sum_{j\neq i}{x_i-x_j\over\|x_i-x_j\|}$.
\bea
\label{EvEq1}\textstyle{du_i\over dt}=-J_i(u),~~~~u_i(0)=x_i,~~~~i=1,...,n.
\eea

\begin{prp}\label{ConvPrp}
The solution $u_i(t)$ of (\ref{EvEq1}) lies in the convex hull $\Conv(x)\subset X$. Moreover, if $t>0$, $u_i(t)$ lies in the interior of $\Conv(x)$.
\end{prp}
\begin{proof}
Observe that $\sum_iJ_i(x)=0$ for all $x\in X$, and so the solution of (\ref{EvEq1}) satisfies
\bea
\label{EvEq2}\textstyle\sum_{i=1}^nu_i(t)=\sum_{i=1}^nx_i,~~~~\txt{for all time}~~t.
\eea
Moreover, (\ref{EvEq1}) can be written in matrix form ~${du\over dt}=Mu$~ as follows. With ~$\al_{ij}:={1\over\|u_i-u_j\|}$,
{\tiny\begin{align}
{d\over dt}\left[\!\!
             \begin{array}{c}
               u_1 \\
               u_2 \\
               \vdots \\
               u_{i-1} \\
               u_i \\
               u_{i+1} \\
               \vdots\\
               u_{n-1} \\
               u_n\\
             \end{array}
           \!\!\right]=
\left[\!\!
  \begin{array}{ccccccccc}
    -\sum\limits_{j\neq 1}\al_{1j} & \al_{12} & \cdots &  & \al_{1i} & & \cdots & \al_{1~n-1} & \al_{1n} \\
    \al_{12} & -\sum\limits_{j\neq 2}\al_{2j} & \cdots &  & \al_{2i} & & \cdots &\al_{2~n-1} & \al_{2n}\\
    \vdots & \vdots &  &  & \vdots & & & \vdots & \vdots\\
     &  & &  & \al_{i-1~i} &  &  & & \\
     & & & & & & & &\\
    \al_{1i} & \al_{2i} & \cdots & \al_{i-1~i} & -\sum\limits_{j\neq i}\al_{ij} & \al_{i~i+1} & \cdots & \al_{i~n-1} & \al_{in}\\
    & & & & & & & &\\
    &  &  & & \al_{i~i+1} &  &  & & \\
    \vdots & \vdots &  &  & \vdots & & & \vdots & \vdots\\
    \al_{1~n-1} & \al_{2~n-1} & \cdots &  &  \al_{i~n-1}&  & \cdots & -\sum\limits_{j\neq n-1}\al_{n-1~j} & \al_{n-1~n}\\
    \al_{1n} & \al_{2n} & \cdots       &  & \al_{in} &  & \cdots & \al_{n-1~n} & -\sum\limits_{j\neq n}\al_{nj}\\
  \end{array}
\!\!\right]
           \left[\!\!
             \begin{array}{c}
               u_1 \\
               u_2 \\
               \vdots \\
               u_{i-1} \\
               u_i \\
               u_{i+1} \\
               \vdots\\
               u_{n-1} \\
               u_n\\
             \end{array}
           \!\!\right]\nn
\end{align}}Therefore, via \emph{iteration} of the integral $u(t)=x+\int_0^tM(s)u(s)ds$ (obtained from ${du\over dt}=Mu$), we can use the series (\ref{PBSeries}) to write the system (\ref{EvEq1}) in the integral form ~$u(t)=Te^{\int_0^tM(s)ds}x$,~ i.e.,
\bea
\label{EvEq3}u(t)=T e^{A(t)}x,~~~~A(t):=\int_0^tM(s)ds.
\eea
Thus, it is enough to show that the matrix $Te^A$ is \emph{right-stochastic} (i.e., has \emph{nonnegative entries} and its \emph{row sums} equal 1).

Indeed, the matrix $T e^A$ has \emph{nonnegative entries} because we can write it in the form
\bea
\label{EvEq4}T e^{A(t)}=Te^{\big(A(t)+ctI\big)-ctI}=T e^{A(t)+ctI}e^{-ctI},~~~~\txt{for any}~~~~c\in\Real,
\eea
which is due to the fact that for any $c\in\Real$ the equation ${du\over dt}=Mu$ can be written in the form
\bea
\label{EvEq5}\textstyle{d(e^{ctI}u)\over dt}=(M+c I)e^{ctI}u.
\eea
Also, the column sums (and also row sums by \emph{symmetry}) of $M$ are zero, i.e., $\sum_{i=1}^nM_{ij}(t)=0$ for all $j$, and so the same is true of ~$M^l(t_1,...,t_l):=M(t_1)M(t_2)\cdots M(t_l)$, ~$l\geq 1$,~ since
\bea
\label{EvEq6}\textstyle\sum\limits_{i=1}^nM^l_{ij}=\sum\limits_{i=1}^n\sum\limits_kM_{ik}M^{l-1}_{kj}=\sum\limits_k\left(\sum\limits_{i=1}^nM_{ik}\right)M^{l-1}_{kj}=\sum\limits_k0M^{l-1}_{kj}=0,
\eea
which implies the column sums of $T e^A$ (and also row sums of $Te^A$ by symmetry of $M$) equal $1$. This proves that $u_i(t)$ lies in the convex hull of $\{x_1,...,x_n\}$.

Since the entries of $Te^{A(t)}$ are positive (i.e., strictly nonnegative) for $t>0$, it follows that if $t>0$, then $u_i(t)$ lies in the interior of the convex hull of $\{x_1,...,x_n\}$.
\end{proof}

\begin{rmk}
For $n=2$, with $\beta(t):=\int_0^t{ds\over\|u_1(s)-u_2(s)\|}$ and $Q:=\left[
  \begin{array}{cc}
    -1 & 1 \\
    1 & -1 \\
  \end{array}
\right]$, we have
\bea
\textstyle Te^{A(t)}=e^{\beta(t)Q}=I+{1\over 2}\left(1-e^{-2\beta(t)}\right)Q
= {1\over 2}\left[
   \begin{array}{cc}
     1+e^{-2\beta(t)}& 1-e^{-2\beta(t)} \\
     1-e^{-2\beta(t)} & 1+e^{-2\beta(t)} \\
   \end{array}
 \right].\nn
\eea
Therefore, $J_i\big(u(t)\big)=J_i(x)$ at all times $t$, and so
\bea
\textstyle u_i(t)=x_i-J_i(x)t,~~~~~~~~J_1(x)=-J_2(x)={x_1-x_2\over\|x_1-x_2\|},\nn
\eea
meaning we always get a Lipschitz retraction $X(2)\ra X$ via the procedure used in the proof of Theorem \ref{MainThm}.
\end{rmk}

\begin{crl}
Let $\H$ be a Hilbert space and $X$ a normed space. (i) Every convex subset $C\subset\H$ admits Lipschitz retractions $C(n)\ra C(n-1)$. (ii) Every convex subset $K\subset X$ admits retractions $K(n)\ra K(n-1)$ that are H\"older continuous on bounded subsets of $K$.
\end{crl}

\begin{dfn}[\textcolor{blue}{\index{Selection of a map}{Selection of a set-valued map}}]
Let $X,Y$ be spaces. Given a set-valued map $F:X\ra 2^Y$, a selection of $F$ is any map $f:X\ra Y$ such that $f(x)\in F(x)$ for all $x\in X$.
\end{dfn}

\begin{crl}[\textcolor{OliveGreen}{Selection theorems}]\label{SelectThms}
Let $\H$ be a Hilbert space and $X$ a Banach space.
\begin{enumerate}[leftmargin=0.9cm]
\item[(a)] The map $F:\H(n)\ra 2^\H$, $F(x):=\Conv(x)$ has a Lipschitz selection.
\item[(b)] The map $F:X(n)\ra 2^X$, $F(x):=\Conv(x)$ has a selection that is H$\ddot{\txt{o}}$lder-continuous on bounded sets.
\end{enumerate}
\end{crl}
\begin{proof}
By Remark \ref{HilbConRmk}, we have Lipschitz retractions $r_n:\H(n)\ra\H(n-1)$, and by Theorem \ref{MainThm}, we have retractions $s_n:X(n)\ra X(n-1)$ that are H$\ddot{\txt{o}}$lder-continuous on bounded sets. Thus, by Proposition \ref{ConvPrp}, the following maps give the desired selections of ~$F(x):=\Conv(x)$.
\begin{align}
&f=r_2\circ r_3\circ\cdots \circ r_n:\H(n)\sr{r_n}{\ral}\H(n-1)\sr{r_{n-1}}{\ral}\cdots\sr{r_3}{\ral}\H(2)\sr{r_2}{\ral}\H,\nn\\
& g=s_2\circ s_3\circ\cdots \circ s_n:X(n)\sr{s_n}{\ral}X(n-1)\sr{s_{n-1}}{\ral}\cdots\sr{s_3}{\ral}X(2)\sr{s_2}{\ral}X.\nn\qedhere
\end{align}
\end{proof}

\begin{rmk}\label{SvRmk}
Corollary \ref{SelectThms} is related to a result in \cite[Theorem 1.2]{shvartsman2004} that asserts the existence (for each $n>0$) of an \ul{affine} Lipschitz selector $S:\K_n(X)\ra X$, $K\mapsto S(K)\in K$, where, as a metric space with respect to Hausdorff distance,
\bea
\K_n(X):=\{\txt{convex compact subsets $K\subset X$ of dimension $\dim K\leq n$}\},\nn
\eea
and hence the existence of an \ul{affine} Lipschitz selection of $F(x):=\Conv(x)$ given by
\bea
s=S\circ\Conv:X(n)\sr{\Conv}{\ral}\K_n(X)\sr{S}{\ral}X.\nn
\eea

The main advantages of Corollary \ref{SelectThms} are simplicity and concreteness. Note, however, that our Lipschitz selections in Corollary \ref{SelectThms} are not affine.
\end{rmk}

\chapter{FSR Property II: Further Considerations}\label{FSRP2}
This chapter relies on notation/terminology from sections \ref{PrelimsAT} and \ref{PrelimsMET}. Recall that the \emph{finite subset retraction ($\FSR$) property} concerns the existence of Lipschitz retractions $r_n:X(n)\ra X(n-1)$. In this chapter, we continue investigation on the FSR property for finite subsets of normed spaces, especially the existence of Lipschitz retractions $X(n)\ra X(n-1)$ for finite-dimensional normed spaces $X$ (in section \ref{FSRP2lr}, Theorem \ref{FinDimLipRet}) and the derivation of a lower bound on the growth rate of Lipschitz constants (in section \ref{GrwLipConst}, Theorem \ref{LipConstGr}).

Lemma \ref{FinSubDoubl} (which proves that if $X$ is a doubling metric space then so is $X(n)$) partly solves \cite[Problem 4.1]{BorovEtal2010} in which it was asked to show that $\Real(n)$ is doubling. Also, Theorem \ref{LipConstGr} provides a partial negative answer to \cite[Question 3.2]{kovalev2016} and \cite[Remark 3.4]{bacac-kovalev2016} which asked whether $\Lip(r_n)$ can be bounded above by a constant that is independent of $n$.

We also consider questions on the FSR property for metric spaces more general than normed spaces (in sections \ref{FSRP2fs}, \ref{FSRP2sn}). Section \ref{FSRP2cf} concerns a few counterexamples and facts that do not involve finite subset spaces, but which are potentially relevant in answering questions on the FSR property.

This chapter differs from the previous chapter in that it mostly poses questions (such as in sections \ref{FSRP2fs}, \ref{FSRP2sn}), while providing some answers that seriously rely on other major results in the literature (especially in sections \ref{FSRP2lr}, \ref{GrwLipConst}), and is therefore less self-contained in comparison to the previous chapter.

\section{The finite subset retraction (FSR) property}\label{FSRP2fs}
\begin{dfn*}[\textcolor{blue}{Recall: FSR property}]
We say a metric space $X$ has the \emph{finite subset retraction ($\FSR$) property} if there exist Lipschitz retractions $X(n)\ra X(n-1)$ for all $n\geq 2$. We also say $X$ has the \emph{$\FSR(k)$ property} if for each $1\leq l\leq k$, there exist Lipschitz retractions $X(n)\ra X(l)$ for all $n\geq l$. Equivalently, $X$ has the \emph{$\FSR(k)$ property} if there exist Lipschitz retractions (i) $X(n)\ra X(n-1)$ for all $2\leq n\leq k$,  and (ii) $X(n)\ra X(k)$ for all $n\geq k$.
\end{dfn*}

\begin{dfn}[\textcolor{blue}{\index{Disjoint union}{Disjoint union}, \index{Positively separated union}{Positively separated union}}]
Let $X=A\cup B$ be a metric space such that $A,B$ are nonempty. We will say $X=A\cup B$ is a \ul{disjoint union} (written $X=A\sqcup B$) if $A\cap B=\emptyset$. We will say a disjoint union $X=A\sqcup B$ is a \ul{positively separated union} if $\dist(A,B)>0$.
\end{dfn}

Among other things, we will see that for a metric space $X$ to possess the FSR property, \emph{quasiconvexity} is neither sufficient (e.g., $X=S^1$ in Lemma \ref{NoRetI}) nor necessary (e.g., $X=$ snowflake in Lemma \ref{SnFkFSRP}). Similar existence/nonexistence results, regarding sufficiency/insufficiency of \emph{(Lipschitz) contractibility} of $X$ for the FSR property, are given in Proposition \ref{NoLipRetCon} and Lemma \ref{SnFkFSRP}. In addition to the above, the following is a list of questions, on the FSR property, for which partial answers at least are given in this section.

\begin{questions}
Let $X$ be a metric space.
\begin{enumerate}[leftmargin=0.9cm]
\item If $A,B$ are metric spaces that have the FSR property, does it follow that $X:=A\times B$ has the FSR property? A partial answer is given in Proposition \ref{FDRPprp0}.
\item If $X=A\cup B$ and $A,B$ have the FSR property, does it follow that $X$ has the FSR property? Recall that $X=(A\backslash B)\sqcup(A\cap B)\sqcup(B\backslash A)$. (\emph{Answer = No}, by Lemma \ref{NoRetI}, with the example $S^1=A\cup B$, where $A,B\subset\Real^2$ are circular arcs.)
\item If $A,B$ are metric spaces that have the FSR property, does it follow that a positively separated union $X=A\sqcup B$ has the FSR property? A partial answer is given in Proposition \ref{FDRPprp1}.
\item Does a finite-dimensional normed space have the FSR property? A positive answer is given in Theorem \ref{FinDimLipRet}.
\item If $X$ has the FSR property, does it follow that $X/Y$ has the FSR property for all $Y$? (\emph{Answer = No}, by Lemmas \ref{NoRetI}, \ref{NoRetIII} and the example $[0,1]/\del[0,1]\cong S^1$.)
\item If $Z\subset X$ is a Lipschitz retract and $X$ has the FSR property, does it follow that $Z$ has the FSR property? (\emph{Answer = Yes}, by Corollary \ref{LipInductII})
\item If $X$ is a normed space (with the FSR property), can we find Lipschitz retractions $r_n:X(n)\ra X(n-1)$ such that $\Lip(r_n)$ is bounded above by a constant independent of $n$? If $\dim X\geq 2$, a negative answer is given in Theorem \ref{LipConstGr}.
\end{enumerate}
\end{questions}

\begin{notation*}[\textcolor{blue}{\index{Based Lipschitz retractions}{Based Lipschitz retractions}}]
Denote any Lipschitz retraction {\small$X(n)\ra X(k)$} by
\bea
\label{BasedLReq}r^X_{nk}:X(n)\ra X(k).
\eea
\end{notation*}

\begin{prp}\label{FDRPprp0}
Let $A,B$ be metric spaces and $X:=A\times B$. If $A,B$ have the $\FSR(1)$ property, then so does $X$.
\end{prp}
\begin{proof}
Given Lipschitz retractions $r_{n1}^A:A(n)\ra A$, $r_{n1}^B:B(n)\ra B$, define $r:(A\times B)(n)\ra A\times B$ by ~$r\big(\{(a_1,b_1),\cdots,(a_n,b_n)\}\big):=\left(r_{n1}^A(a_1,...,a_n),r_{n1}^B(b_1,...,b_n)\right)$. Then
\begin{align}
&d\Big(r\big(\{(a_1,b_1),\cdots,(a_n,b_n)\}\big),r\big(\{(a'_1,b'_1),\cdots,(a'_n,b'_n)\}\big)\Big)\nn\\
&~~=d\Big(\left(r_{n1}^A(a_1,...,a_n),r_{n1}^B(b_1,...,b_n)\right),\left(r_{n1}^A(a'_1,...,a'_n),r_{n1}^B(b'_1,...,b'_n)\right)\Big)\nn\\
&~~\leq\max\Big(d\left(r_{n1}^A(a_1,...,a_n),r_{n1}^A(a'_1,...,a'_n)\right),d\left(r_{n1}^B(b_1,...,b_n),r_{n1}^B(b'_1,...,b'_n)\right)\Big)\nn\\
&~~\leq\max\left(\Lip r_{n1}^A,\Lip r_{n1}^B\right)\max\Big[d_H\big(\{a_1,...,a_n\},\{a'_1,...,a'_n\}\big),d_H\big(\{b_1,...,b_n\},\{b'_1,...,b'_n\}\big)\Big]\nn\\
&~~\leq\max\left(\Lip r_{n1}^A,\Lip r_{n1}^B\right)\max\Big[d_H\big(\{a_i\},\{a'_i\}\big),d_H\big(\{b_i\},\{b'_i\}\big)\Big]\nn\\
&~~\leq \max\left(\Lip r_{n1}^A,\Lip r_{n1}^B\right)d_H\big(\{(a_1,b_1),\cdots,(a_n,b_n)\},\{(a'_1,b'_1),\cdots,(a'_n,b'_n)\}\big),\nn
\end{align}
where the last inequality is due to the following:
{\footnotesize\begin{align}
& d_H\big(\{(a_i,b_i)\},\{(a'_i,b'_i)\}\big)=\max\left(\max_i\min_jd\big((a_i,b_i),(a'_j,b'_j)\big),\max_j\min_id\big((a_i,b_i),(a'_j,b'_j)\big)\right)\nn\\
&~~=\max\left(\max_i\min_j\max\Big(d(a_i,a'_j),d(b_i,b'_j)\Big),\max_j\min_i\max\Big(d(a_i,a'_j),d(b_i,b'_j)\Big)\right)\nn\\
&~~\geq \max\left(\max\Big(\max_i\min_jd(a_i,a'_j),\max_i\min_jd(b_i,b'_j)\Big),\max\Big(\max_j\min_id(a_i,a'_j),\max_j\min_id(b_i,b'_j)\Big)\right)\nn\\
&~~= \max\left(\max\Big(\max_i\min_jd(a_i,a'_j),\max_j\min_id(a_i,a'_j)\Big),\max\Big(\max_i\min_jd(b_i,b'_j),\max_j\min_id(b_i,b'_j)\Big)\right)\nn\\
&~~=\max\Big[d_H\big(\{a_i\},\{a'_i\}\big),d_H\big(\{b_i\},\{b'_i\}\big)\Big].\nn
\end{align}}
\end{proof}

\begin{prp}\label{FDRPprp1}
Let $X$ be a \ul{bounded} metric space that can be written as a positively separated union $X=A\sqcup B$. If $A,B$ have the $\FSR(2)$ property, then so does $X$.
\end{prp}
\begin{proof}
Observe that with $A(k)\ast B(n-k):=\{a\cup b:a\in A(k),b\in B(n-k)\}$, we have
\bea
\textstyle X(n)=\big\{\{x_1,...,x_n\}:(x_1,...,x_n)\in X^n\}\big\}=A(n)\sqcup\left[\bigcup_{k=1}^{n-1}A(k)\ast B(n-k)\right]\sqcup B(n).\nn
\eea
Let $AB(n):=\bigcup_{k=1}^{n-1}A(k)\ast B(n-k)$, so that $X(n)=A(n)\sqcup AB(n)\sqcup B(n)$. Then explicitly,
{\footnotesize\begin{align}
&\textstyle X(1)=A(1)\sqcup B(1),\nn\\
&\textstyle X(2)=A(2)\sqcup\overbrace{\Big(A(1)\ast B(1)\Big)}^{AB(2)}\sqcup B(2),\nn\\
&\textstyle X(3)=A(3)\sqcup\overbrace{\Big(\big[A(2)\ast B(1)\big]\cup\big[A(1)\ast B(2)\big]\Big)}^{AB(3)}\sqcup B(3),\nn\\
&\textstyle X(4)=A(4)\sqcup\overbrace{\Big(\big[A(3)\ast B(1)\big]\cup\big[\textcolor{blue}{A(2)\ast B(2)}\big]\cup\big[A(1)\ast B(3)\big]\Big)}^{AB(4)}\sqcup B(4),\nn\\
&\textstyle X(5)=A(5)\sqcup\overbrace{\Big(\big[A(4)\ast B(1)\big]\cup\big[\textcolor{blue}{A(3)\ast B(2)}\big]\cup\big[\textcolor{blue}{A(2)\ast B(3)}\big]\cup\big[A(1)\ast B(4)\big]\Big)}^{AB(5)}\sqcup B(5),\nn\\
&\textstyle~~~~~~~~~~~~~\vdots~~~~~~~~~~~~~~~~~~~~~~~~~~~~~~~~\vdots~~~~~~~~~~~~~~~~~~~~~~~~~~~~~~~~\vdots~~~~\nn\\
&\textstyle X(n-1)=A(n-1)\sqcup\overbrace{\Big(\big[A(n-2)\ast B(1)\big]\cup\cdots\cup\big[A(1)\ast B(n-2)\big]\Big)}^{AB(n-1)}\sqcup B(n-1),\nn\\
&\textstyle X(n)=A(n)\sqcup\overbrace{\Big(\big[A(n-1)\ast B(1)\big]\cup\cdots\cup\big[A(1)\ast B(n-1)\big]\Big)}^{AB(n)}\sqcup B(n).\nn
\end{align}}Note that if either $x\in A(n)$ and $y\in AB(n)$, or $x\in AB(n)$ and $y\in B(n)$, then
\bea
d_H(x,y):=\max_{i,j}\max\big(\dist(x_i,y),\dist(x,y_j)\big)\geq\dist(A,B).\nn
\eea
Equivalently, if $d_H(x,y)<\dist(A,B)$ then either $x,y\in A(n)$, or $x,y\subset AB(n)$, or $x,y\in B(n)$. Define maps $R_{n1}:X(n)\ra X$ and $R_{n2}:X(n)\ra X(2)$ by
\begin{align*}
&R_{n1}(x):=\left\{
             \begin{array}{ll}
               r^A_{n1}(x\cap A), & \txt{if}~~x\cap A\neq\emptyset \\
               r^B_{n1}(x), & \txt{if}~~x\cap A=\emptyset
             \end{array}
           \right\},\nn\\
&R_{n2}(x):=\left\{
            \begin{array}{ll}
              r^A_{n2}(x), ~~& \txt{if}~~x\in A(n) \\
              G(x)=r_{n1}^A(x\cap A)~\cup~r_{n1}^B(x\cap B),~~ &\txt{if}~~ x\in AB(n) \\
              r^B_{n2}(x),~~ & \txt{if}~~x\in B(n)
            \end{array}
          \right\}.\nn
\end{align*}
Then for a fixed $0<\delta<\dist(A,B)$, and any $x,y\in X(n)$, we have
\bea
d_H\left(R_{n1}(x),R_{n1}(y)\right)\leq
\left\{
  \begin{array}{ll}
    {\diam X\over\delta}d_H(x,y), & ~~\txt{if}~~d_H(x,y)\geq\delta \\
    \max\left(\Lip(r^A_{n1}),\Lip(r^B_{n1})\right)d_H(x,y), & ~~\txt{if}~~d_H(x,y)\leq\delta\nn
  \end{array}
\right\},\nn
\eea
which shows $R_{n1}$ is Lipschitz. Similarly, with $d_H(x,y)\leq\delta<\dist(A,B)$ at step (s) below,
\bea
&&d_H\big(G(x),G(y)\big)=d_H\Big(r_{n1}^A(x\cap A)~\cup~r_{n1}^B(x\cap B),r_{n1}^A(y\cap A)~\cup~r_{n1}^B(y\cap B)\Big)\nn\\
&&~~~~\leq\max\left[d_H\Big(r_{n1}^A(x\cap A),r_{n1}^A(y\cap A)\Big),d_H\Big(r_{n1}^B(x\cap B),r_{n1}^B(y\cap B)\Big)\right]\nn\\
&&~~~~\leq \max\left(\Lip r^A_{n1},\Lip r^B_{n1}\right)\max\left[d_H(x\cap A,y\cap A),d_H(x\cap B,y\cap B)\right]\nn\\
&&~~~~\sr{(s)}{=}\max\left(\Lip r^A_{n1},\Lip r^B_{n1}\right)d_H(x,y).\nn
\eea
It follows that, like $R_{n1}$, the map $R_{n2}$ also satisfies
\bea
d_H\left(R_{n2}(x),R_{n2}(y)\right)\leq
\left\{
  \begin{array}{ll}
    {\diam X\over\delta}d_H(x,y), & ~~\txt{if}~~d_H(x,y)\geq\delta \\
    \max\left(\Lip(r^A_{n1}),\Lip(r^B_{n1})\right)d_H(x,y), & ~~\txt{if}~~d_H(x,y)\leq\delta\nn
  \end{array}
\right\}.\nn\qedhere
\eea
\end{proof}

\begin{crl}\label{FDRPprp2}
Let $X$ be a \ul{bounded} metric space that can be written as a positively separated union $X=A\sqcup B$. If $A,B$ have the $\FSR$ property, then the obvious Lipschitz retractions $A(n)\cup B(n)\ra A(n-1)\cup B(n-1)$ extend to Lipschitz maps $(A\cup B)(n)\ra (A\cup B)(n-1)$ that are not necessarily retractions.
\end{crl}
\begin{proof}
Let $r_n^A:A(n)\ra A(n-1)$, $r_n^B:B(n)\ra B(n-1)$, and $R_{n2}:X(n)\ra X(2)$ be Lipschitz retractions (where $R_{n2}$ exists as constructed in the proof if Proposition \ref{FDRPprp1}). Define $R_n:X(n)\ra X(n-1)$ by
{\small$
R_n(x):=\left\{\!\!\!
            \begin{array}{ll}
              r^A_n(x), ~~& \txt{if}~~x\in A(n) \\
              R_{n2}(x),~~ &\txt{if}~~ x\in AB(n) \\
              r^B_n(x),~~ & \txt{if}~~x\in B(n)
            \end{array}
          \!\!\!\right\}.\nn
$}
Then for a fix $0<\delta\leq\dist(A,B)$, and any $x,y\in X(n)$, we have
{\small\bea
d_H\big(R_n(x),R_n(y)\big)\leq
\left\{
  \begin{array}{ll}
    {\diam X\over\delta}d_H(x,y), & ~~\txt{if}~~d_H(x,y)\geq\delta \\
    \max\Big\{\Lip(r_n^A),\Lip(R_{n2}),\Lip(r_n^B)\Big\}d_H(x,y), & ~~\txt{if}~~d_H(x,y)<\delta
  \end{array}
\right\}.\nn\qedhere
\eea}
\end{proof}

The proof of the above corollary shows that to find a Lipschitz retraction $(A\times B)(n)\ra (A\sqcup B)(n-1)$, it suffices to find a Lipschitz retraction $AB(n)\ra (A\sqcup B)(n-1)$ replacing $R_{n2}$ above. Such a map may be defined piecewise on the disjoint pieces of $AB(n)$, which \emph{might not be easy} because the disjoint pieces of $AB(n)$ are \emph{not separated} the way $A(n)$, $AB(n)$, $B(n)$ are separated.

Note that with $\wt{X}(n):=X(n)\backslash X(n-1)$, and the observation that $\wt{X}(1)=X(1)$, we can further expand the middle component $AB(n)$ of $(A\sqcup B)(n)=A(n)\sqcup AB(n)\sqcup B(n)$ in the form $AB(n)=\wt{AB}(n)\sqcup AB(n-1)$, where
\bea
\textstyle\wt{AB}(n):=\bigsqcup\limits_{k=1}^{n-1}\wt{A}(k)\ast\wt{B}(n-k)\cong \bigsqcup\limits_{k=1}^{n-1}\wt{A}(k)\times\wt{B}(n-k).\nn
\eea
In particular,
\bea
&&AB(4)=\Big[\wt{A}(3)\ast B(1)~\sqcup~\wt{A}(2)\ast\wt{B}(2)~\sqcup~A(1)\ast\wt{B}(3)\Big]~\sqcup~AB(3)\nn\\
&&~~~~\cong \Big[\wt{A}(3)\times B(1)~\sqcup~\wt{A}(2)\times\wt{B}(2)~\sqcup~A(1)\times\wt{B}(3)\Big]~\sqcup~AB(3).\nn
\eea

\begin{example}[\textcolor{OliveGreen}{\index{Rickman's Rug}{Rickman's Rug}}]\label{RigRugEx}
Let $U:=[0,1)$, $V\subset\Real^2$ a snowflake curve, and $X:=U\times V$. It is clear that $X$ is not a snowflake metric space since $X$ contains rectifiable curves. Also, $X$ is not Lipschitz contractible since $V$ contains no rectifiable curves (Lemma \ref{RigRugLmm2}).
\end{example}

\begin{question}\label{RigRugQs}
Does Rickman's Rug $X=U\times V$ have the FSR property? (Note that $U$ and $V$ both have the FSR property, $U$ as a convex subset of a Hilbert space and $U$ by Lemma \ref{SnFkFSRP}.)
\end{question}

\section{Lipschitz retraction of finite subsets of finite-dimensional normed spaces}\label{FSRP2lr}
Recall that $X(n)$ is Lipschitz $k$-connectedness for all $k\geq 0$ when $X$ is a normed space. Using this fact, we will show that for a finite-dimensional normed space $X$, we have Lipschitz retractions $X(n)\ra X(n-1)$.

\begin{rmk}
Using the fact that any finite-dimensional normed space is isomorphic to a Euclidean space, one can deduce the existence of Lipschitz retractions $X(n) \to X(n-1)$ from the corresponding result for Hilbert spaces in \cite{kovalev2016}. However, this approach results in Lipschitz constants that grow to $\infty$ with the dimension of the space and therefore cannot lead to the solution for general normed spaces. Therefore, we look for other approaches to finite-dimensional spaces in the hope of finding retractions with dimension-independent Lipschitz constants.
\end{rmk}

\begin{dfn}[\textcolor{blue}{\index{Doubling! metric space}{Doubling metric space}, \index{Doubling! dimension}{Doubling dimension}}]
A metric space $X$ is \ul{doubling} if there exists an integer $N$ such that every ball $B_R(x)$ can be covered by $N$ half-radius balls $B_{R/2}(x_1),\cdots,B_{R/2}(x_N)$. The \ul{doubling dimension} $\dim_DX$ of $X$ is given by
{\small\bea
\textstyle 2^{\dim_D(X)}:=\inf\left\{N:B_R(x)\subset\bigcup\limits^N_{i=1}B_{R/2}(x_i)\right\}.\nn
\eea}That is, $\dim_D(X)$ is the smallest integer $m$ such that every ball $B_R(x)$ can be covered by $2^m$ half-radius balls $B_{R/ 2}(x_1),\cdots,B_{R/2}(x_{2^m})$.
\end{dfn}

\begin{facts}~
\begin{enumerate}[leftmargin=0.9cm]
\item If $X$ is doubling then so is $X^n$, where ~$d\big((x_1,...,x_n),(x_1',...,x_n')\big):=\max_id(x_i,x_i')$.
\begin{proof}
A ball $B_R(x_1,...,x_n)$ in $X^n$ can be written as $B_R(x_1,...,x_n)=B_R(x_1)\times\cdots\times B_R(x_n)$ for balls $B_R(x_i)$ in $X$. Thus, if a collection of half-radius balls $\C_i$ covers $B_R(x_i)$, then the collection $\left\{B_1\times \cdots\times B_n:B_i\in\C_i\right\}$ covers $B_R(x_1,...,x_n)$.
\end{proof}
\item Every finite-dimensional normed space is doubling.
\begin{proof}
Let $X$ be a normed space with $\dim X<\infty$. We know that all norms on $X$ are equivalent, and that equivalence of norms is a biLipschitz homeomorphism. Hence, by the following, $X$ is doubling.
\bit
\item[(i)] Because $\Real$ is clearly doubling, the normed space $\big(\Real^n,\max|\cdot|\big)$ is doubling.
\item[(ii)] A biLipschitz homeomorphism preserves doubling: \cite[Prop. 3.29(b), p.224]{brudnyis2012v1}\qedhere
\eit
\end{proof}

\item The Nagata dimension $\dim_NX$ of a doubling metric space $X$ is finite: \cite[Theorem 4.31, p.350]{brudnyis2012v1}
\item If $X$ is a finite dimensional normed space, then $\dim_NX(n)<\infty$. (This follows from the items above.)

\item $\big(\Real,\min(\|x-x'\|,1)\big)$ is not doubling because for $0<\vep<1$, the ball $B_{1+\vep}(0)=\Real$ cannot be covered by finitely many balls $B_{(1+\vep)/2}(x_i)$.
\item A Lipschitz map does not preserve doubling: Let $X=\big(\Real,\min(|x-x'|,1)\big)$. Then the map $f:\Real\ra X$, $f(x)=x$ is 1-Lipschitz, but $\Real$  is doubling while $X$ is not doubling.
\end{enumerate}
\end{facts}

\begin{dfn}[\textcolor{blue}{\index{Quasisymmetric map}{Quasisymmetric map}}]
A map of metric spaces $f:X\ra Z$ is quasisymmetric if (i) $f$ is an imbedding (i.e., a homeomorphism onto its image), and (ii) there exists a homeomorphism $\eta:[0,\infty)\ra[0,\infty)$ such that
\begin{align}
d(x,a)\leq td(x,b)~\Ra~d\big(f(x),f(a)\big)\leq\eta(t)d\big(f(x),f(b)\big),~~\txt{for all}~~x,a,b\in X,~t\in[0,\infty).\nn
\end{align}
In this case, we say $f$ is $\eta(t)$-quasisymmetric.
\end{dfn}
A quasisymmetric map is a generalization of a bi-Lipschitz map. Note that a $c$-biLipschitz map is $c^2t$-quasisymmetric. Also, for any $0<\vep<1$, the identity map $id:(X,d)\ra (X,d^\vep)$ is $t^\vep$-quasisymmetric.

\begin{thm}[\textcolor{OliveGreen}{\cite[Theorem 12.1, p.98]{heinonen2}}]\label{QSymmEmbedThm}
A metric space $X$ quasimmetrically imbeds into some $\Real^d$ $\iff$  $X$ is doubling.
\end{thm}

\begin{thm}[\textcolor{OliveGreen}{Assouad's Theorem: \cite[Theorem 12.2, p.98]{heinonen2}}]\label{AssouadThm}
Given a metric space $(X,d)$ and $0<\vep<1$, let $d^\vep$ be the snowflake metric on $X$ given by $d^\vep(x,x'):=d(x,x')^\vep$.

If $X$ is doubling, then there exists a bi-Lipschitz imbedding $(X,d^\vep)\ra\Real^N$ for some $N$.
\end{thm}

\begin{lmm}[\textcolor{OliveGreen}{Doubling property of finite subset spaces}]\label{FinSubDoubl}
If $X$ is a doubling metric space, then so is $X(n)$.
\end{lmm}
\begin{proof}
Since $X$ is doubling, it follows by Assouad's Theorem that we have a biLipschitz imbedding $(X,d^\vep)\ra \Real^N$. This induces a bi-Lipschitz imbedding $(X(n),d_H^\vep)\ra \Real^N(n)$.

By \cite[Theorem 1.2]{kovalev2015}, $\Real^N(n)$ admits a bi-Lipschitz imbedding into some $\Real^{N_n}$, i.e., we have a bi-Lipschitz imbedding $(X(n),d_H^\vep)\ra \Real^{N_n}$.
\bc\bt
(X,d)\ar[d]\ar[rr,"\txt{Quasisymm}"]&& (X,d^\vep)\ar[d]\ar[rr,hook,"\txt{BiLip}"] && \Real^N\ar[d] &&\\
(X(n),d_H)\ar[rr,"\txt{Quasisymm}"]&& \big(X(n),\big(d^\vep)_H\big)=\big(X(n),d^\vep_H\big)\ar[rr,hook,"\txt{BiLip}"] && \Real^N(n)\ar[rr,hook,"\txt{BiLip}"] &&\Real^{N_n}
\et\ec
Hence, it follows from Theorem \ref{QSymmEmbedThm} that $X(n)$ is doubling.

(Note: Since (i) every subspace of a doubling space is doubling and (ii) a bi-Lipschitz map preserves doubling, we see that $(X(n),d_H^\vep)$ is doubling.)
\end{proof}

\begin{thm}\label{FinDimLipRet}
If $X$ is a finite-dimensional normed space, then Lipschitz retractions $X(n)\ra X(n-1)$ exist, for all $n\geq 2$.
\end{thm}
\begin{proof}
Since $X(n)$ is doubling (Lemma \ref{FinSubDoubl}), the Nagata dimension of $X(n)$ is finite. Moreover, by Theorem \ref{LipConSuff5}, $X(n)$ is Lipschitz $k$-connected for all $k\geq 0$. Thus, by \cite[Theorem 6.26, p.20]{brudnyis2012v2}, there exists a constant $\ld_n=\ld_n(\dim X)$ such that every $c$-Lipschitz map $X(n-1)\ra X(n-1)$ extends to a $\ld_nc$-Lipschitz map $X(n)\ra X(n-1)$. In particular, $id:X(n-1)\ra X(n-1)$ extends to a $\ld_n$-Lipschitz retraction $r:X(n)\ra X(n-1)$.

\bc
\hspace{4cm}
\bt
X(n-1)\ar[d,hook]\ar[rr,"id"]& & X(n-1)\\
X(n)\ar[urr,dashed,"r"'] & &
\et
\qedhere
\ec
\end{proof}

\begin{dfn}[\textcolor{blue}{\index{Doubling! measure}{Doubling measure}}]
Let $X$ be a metric space. A measure $\mu$ on $X$ is doubling if there exists a constant $C$ such that all balls in $X$ satisfy
\bea
\mu\big(B_{2R}(x)\big)\leq C~\mu\big(B_R(x)\big).\nn
\eea
\end{dfn}
By \cite{Luu-Sak98} a complete doubling metric space admits a doubling measure, and by \cite{Luu98}, a metric space that admits a doubling measure is doubling.

\section{Growth of the Lipschitz constants}\label{GrwLipConst}
Let $\H$ be a Hilbert space or Hadamard space. From \cite{kovalev2016}, if $\H$ is a Hilbert space, the Lipschitz constants of retractions $r_n:\H(n)\ra\H(n-1)$ satisfy $\Lip(r_n)\leq\max(n^{3/2},2n-1)$. From \cite{bacac-kovalev2016}, if $\H$ is a Hadamard space, the Lipschitz constants of retractions $r_n:\H(n)\ra\H(n-1)$ satisfy $\Lip(r_n)\leq\max(4n^{3/2}+1,2n^2+n^{1/2})$. These upper bounds clearly grow with $n$. It was therefore asked in \cite[Question 3.2]{kovalev2016}, and later also asked in \cite[Remark 3.4]{bacac-kovalev2016}, whether $\Lip(r_n)$ can be bounded above by a constant that is independent of $n$. It will be shown in Theorem \ref{LipConstGr} that if $X$ is any normed space of dimension $\dim X\geq 2$, then every sequence of Lipschitz retractions $r_n:X(n)\ra X(n-1)$ satisfies $\Lip(r_n)\geq {n\over 4\pi\sqrt{2}}$. That is, Theorem \ref{LipConstGr} provides a partial negative answer to the above-mentioned question.
\begin{dfn}[\textcolor{blue}{\index{Bi-orthogonal system}{Bi-orthogonal system}, Unit bi-orthogonal system}]
If $X$ is a normed space, then $\{v_i\}_{i\in I}\subset X$ and $\{u_i\}_{i\in I}\in X^\ast$ form a bi-orthogonal system $(v_i,u_i)_{i\in I}$ in $X\times X^\ast$ if $\langle v_i,u_j\rangle=\delta_{ij}$ (where $\langle v_i,u_j\rangle:=u_j(v_i)$).

A bi-orthogonal system $(v_i,u_i)_{i\in I}$ is said to be unit if $\|v_i\|=\|u_i\|=1$ for all $i\in I$.
\end{dfn}

\begin{note}[\textcolor{OliveGreen}{Existence of bi-orthogonal systems}]
Let $X$ be a normed space. For any finite sets $\{v_1,...,v_n\}\subset X$ and $\{u_1,...,u_n\}\in X^\ast$, if we let $v_i^\ast:=\sum_k[\langle v_i,u_k\rangle]^{-1}u_k$ (whenever the matrix inverse exists), then $\langle v_i,v_j^\ast\rangle=\delta_{ij}$ i.e., $(v_i,v_i^\ast)_{i=1}^{i=n}$ is a bi-orthogonal system in $X\times X^\ast$.
\end{note}

\begin{lmm}[\textcolor{OliveGreen}{Auerbach's Lemma: Existence of unit bi-orthognal systems}]\label{AuerbLmm}
If $X$ is an $n$-dimensional normed space, there exists a unit bi-orthogonal system $(z_i,\vphi_i)_{i=1}^{i=n}$ in $X\times X^\ast$.
\end{lmm}
\begin{proof}
Let $(v_i,v_i^\ast)_{i=1}^{i=n}$ be any bi-orthogonal system in $X\times X^\ast$ (which we know always exists). Define a map $T:\ol{B}_1(0)^n\ra\Real$ by $T(x_1,...,x_n):=\det[x_{ij}]:=\det[\langle x_i,v_j^\ast\rangle]$, where $x_i=\sum_j\langle x_i,v_j^\ast\rangle v_j=\sum_jx_{ij}v_j$. Then $T$ is a continuous map on a compact set in $X^n$ and so attains it maximum (and minimum). Let $T$ attain its maximum at $(z_1,...,z_n)\in \ol{B}_1(0)^n$.

Since $T$ is multilinear (i.e., linear in each argument), the map $T_i:\ol{B}_1(0)\ra\Real$ given (for fixed $x_1,\cdots,x_{i-1},x_{i+1},\cdots,x_n\in X$) by
\bea
T_i(x):=T(x_1,\cdots,x_{i-1},x,x_{i+1},\cdots,x_n)\nn
\eea
extends to a linear functional $T_i:X\ra\Real$. Since convexity (footnote\footnote{The map $f(x):=|T_i(x)|$ is convex, and maxima of a convex function occur on the boundary of its domain.}) implies
\bea
&&\max_{x\in \ol{B}_1(0)}T_i(x)=\max_{x\in \ol{B}_1(0)}|T_i(x)| =\max_{\|x\|\leq 1}|T_i(x)|= \|T_i\|=\max_{\|x\|=1}|T_i(x)|,\nn\\
&&\max_{(x_1,...,x_n)\in \ol{B}_1(0)}T(x_1,...,x_n)=\max_{x_1\in \ol{B}_1(0)}\max_{x_2\in \ol{B}_1(0)}\cdots\max_{x_n\in \ol{B}_1(0)}T(x_1,...,x_n),\nn
\eea
it follows that the maximum point $(z_1,...,z_n)\in \ol{B}_1(0)$ of $T$ consists of unit vectors, i.e.,
\bea
\|z_i\|=1~~~~\txt{for all}~~~~i=1,...,n.\nn
\eea

For each $i$, define a linear functional $\vphi_i:X\ra\Real$ by
\bea
\textstyle\vphi_i(x):={T(z_1,\cdots,z_{i-1},x,z_{i+1},\cdots,z_n)\over T(z_1,\cdots,z_n)}={\det[\langle z_i(x),u_j\rangle]\over\det[\langle z_i,u_j\rangle]},\nn
\eea
where $\big(z_1(x),...,z_n(x)\big):=(z_1,\cdots,z_{i-1},x,z_{i+1},\cdots,z_n)$. By definition, we have $\vphi_i(z_j)=\delta_{ij}$ (for all $i,j$) and $\|\vphi_i\|\leq 1$ (for all $i$). Also, since $\vphi_i(z_i)=1$ implies $\|\vphi_i\|\geq 1$, it follows that $\|\vphi_i\|=1$.
\end{proof}

\begin{lmm}\label{AuLmmBnd}
If $X$ is a normed space and $\M$ an $n$-dimensional linear subspace, then there exists a linear projection $P:X\ra\M$ such that $\|P\|\leq n$.
\end{lmm}
\begin{proof}
By Auerbach's Lemma, there exists a unit bi-orthogonal system $(z_i,\vphi_i)_{i=1}^n$ in $X\times X^\ast$, where $\{z_1,..,z_n\}$ is a basis for $\M$. The linear projection $P:X\ra\M$ given by $P(x):=\sum_{i=1}^n\vphi_i(x)z_i$ satisfies $\|P(x)\|\leq n\|x\|$ for all $x\in X$, i.e., $\|P\|\leq n$.
\end{proof}

\begin{dfn}[\textcolor{blue}{Recall: Minimum separation, Total minimum separation}]
Let $X$ be a metric space and $x=\{x_1,...,x_n\}\in X(n)$. We define the minimum separation $\delta_n(x)$ and total minimum separation $\delta(x)$ by
\bea
\delta_n(x):=\min_{i\neq j}d(x_i,x_j),~~~~\delta(x):=\min\{d(a,b):a,b\in x,~a\neq b\}.\nn
\eea
\end{dfn}
\begin{lmm}\label{PointProxLmm}
If $X$ is a normed space, then the distance function $\dist_H(x,X(n-1))$ satisfies
\bea
\dist_H(x,X(n-1))={1\over 2}\delta_n(x)=
\left\{
  \begin{array}{ll}
    {1\over 2}\delta(x), & \txt{if}~~ x\in X(n)\backslash X(n-1)\\
    0, &\txt{if}~~ x\in X(n-1)
  \end{array}
\right\}.\nn
\eea
\end{lmm}
\begin{proof}
It is clear that $\dist_H(x,X(n-1))=0$ if $x\in X(n-1)$. So, assume $x\in X(n)\backslash X(n-1)$. Let $y\in X(n)$. If $d_H(x,y)<{1\over2}\delta(x)$, then $|y|=n$, i.e., $|y|<n$ implies $d_H(x,y)\geq{1\over2}\delta(x)$. Thus, $d_H(x,z)\geq{1\over 2}\delta(x)$ for all $z\in X(n-1)$, and so
\bea
\textstyle\dist_H(x,X(n-1))\geq{1\over 2}\delta(x).\nn
\eea
Moreover, since $\delta(x)=\|a_0-b_0\|$ for some $a_0,b_0\in x$, with {\small $z:=\big(x\backslash\{a_0,b_0\}\big)\cup\{(a_0+b_0)/2\}\in X(n-1)$},
\bea
\textstyle\dist_H(x,X(n-1))\leq d_H(x,z)\leq{1\over 2}\|a_0-b_0\|={1\over 2}\delta(x).\nn
\eea
\end{proof}

\begin{crl}\label{PointProxCrl}
If $X$ is a normed space, then for any finite set $x\subset X$, we have
\bea
\delta(x)=\delta_{|x|}(x)=2\dist_H\Big(x,X\big(|x|-1\big)\Big).\nn
\eea
\end{crl}

\begin{crl}[\textcolor{OliveGreen}{Upper bound on the displacement of a Lipschitz retraction}]\label{DispBound}
Let $X$ be a normed space. If $r:X(n)\ra X(n-1)$ is a Lipschitz retraction, then
\bea
\textstyle d_H(x,r(x))\leq{1+\Lip(r)\over 2}\delta(x),~~~~\txt{for any}~~x\in X(n).\nn
\eea
\end{crl}
\begin{proof}
Let $r:X(n)\ra X(n-1)$ be a Lipschitz retraction. Fix $x\in X(n)$. Then for any $y\in X(|x|-1)$,
\bea
&&\textstyle d_H(x,r(x))\sr{(s)}{\leq} d_H(x,y)+d_H(r(y),r(x))\leq\big(1+\Lip(r)\big)d_H(x,y),\nn\\
&&\textstyle~~\Ra~~d_H(x,r(x))\leq \big(1+\Lip(r)\big)\dist_H(x,X(|x|-1))={1+\Lip(r)\over 2}\delta(x),\nn
\eea
where step (s) uses the fact that $X(|x|-1)\subset X(n-1)$, and so $y\in X(n-1)$.
\end{proof}

\begin{lmm}[\textcolor{OliveGreen}{Some homotopy invariants: \cite[Theorem 2.10, p.111]{hatcher2001} for homology, and \cite[Exercise 2 of Section 4.1, p.358]{hatcher2001} generalizing \cite[Proposition 1.18, p.37]{hatcher2001} for homotopy}]\label{HtyInvLmm}
(i) If $X\simeq Y$, then $H_n(X)\cong H_n(Y)$ for all $n\geq 0$. (ii) Similarly, if $X\simeq Y$ and $X,Y$ are path-connected, then $\pi_n(X)\cong\pi_n(Y)$ for all $n\geq 0$.
\end{lmm}
\begin{proof}
Part (i) is Corollary \ref{HtyInvCrl1}. For the proof of part (ii), see the indicated references.
\end{proof}

Note that given Lemma \ref{HtyInvLmm}, either one of Lemmas \ref{HtyRetInc},\ref{HomRetInc} suffices for proving Lemma \ref{NoRetIII}.

\begin{lmm}[\textcolor{OliveGreen}{No retraction I}]\label{NoRetI}
There is no continuous retraction $S^1(2)\ra S^1$.
\end{lmm}
\begin{proof}
Suppose a continuous retraction $R:S^1(2)\ra S^1$ exists. Consider a path $h:[0,1]\ra [0,2\pi]$, $h(0)=0$, $h(1)=2\pi$. Fix the point $x_0=\{1\}\in S^1(2)$. Let $\gamma:[0,1]\ra S^1(2)$, $\gamma(t):=\{e^{if(t)},e^{ig(t)}\}$, be a loop based at $x_0$, where $f,g:[0,1]\ra [0,2\pi]$ are given by
\bea
&& f(t):=\left\{
          \begin{array}{ll}
            h(2t), & t\in[0,1/2] \\
            2\pi, & t\in[1/2,1]
          \end{array}
        \right\},~~~~g(t)=\left\{
          \begin{array}{ll}
            0, & t\in[0,1/2] \\
            h\big(2(t-1/2)\big), & t\in[1/2,1]
          \end{array}
        \right\}.\nn
\eea
That is, $\gamma$ winds twice over the loop $\eta:[0,1]\ra S^1$, $\eta(t)=e^{ih(t)}$. Thus, the \emph{winding number} of $R\circ\gamma:[0,1]\sr{\gamma}{\ral}S^1(2)\sr{R}{\ral} S^1$ (a loop in $S^1$ based at $x_0$) is even.

On the other hand, we know $f$ and $g$ are \emph{path-homotopic} (since $[0,2\pi]$ is \emph{contractible}), which implies $R\circ\gamma$ is \emph{path-homotopic} to a loop $[0,1]\ra S^1$ with winding number $1$ (a contradiction).
\end{proof}

\begin{note*} Due to \cite{mostovoy2004,tuffley2002}, Lemma \ref{NoRetI} also follows from the facts that (i) $S^1(2)\cong M$ (the M\"obius band), where the inclusion $S^1\hookrightarrow S^1(2)$ is given by the inclusion $\del M\hookrightarrow M$, and (ii) there exists no continuous retraction $M\ra\del M$.
\end{note*}
\begin{proof}[Proof] of (ii):~ If a continuous retraction $r:M\ra \del M$ exists, then the map $i_\ast:\pi_1(\del M)\cong\Integer\ra\pi_1(M)\cong\Integer,~z\mapsto 2z$ induced by the inclusion $i:\del M\hookrightarrow M$ has a left inverse $r_\ast:\pi_1(M)\ra\pi_1(\del M)$, which is a contradiction.
\end{proof}

\begin{lmm}[\textcolor{OliveGreen}{No retraction II}]\label{NoRetII}
For $1\leq k\leq n-1$, there exists no continuous retraction $S^n\ra S^{n-k}$.
\end{lmm}
\begin{proof}
By Lemmas \ref{HomRetInc} and \ref{HlgGpSph}, if $S^{n-k}\subset S^n$ ($1\leq k\leq n-1$) is a retract, then $\Integer\cong \wt{H}_{n-k}(S^{n-k})\ra \wt{H}_{n-k}(S^n)=0$ is injective (a contradiction).

Alternatively (with homotopy), we can use the fact that $A\subset X$ is a retract $\iff$ every continuous map $f:A\ra Z$ extends to a continuous map $F:X\ra Z$. (In particular, if $S^{n-1}\subset S^n$ is a retract, then by Lemma \ref{NoHoleCrit} and the facts
\bea
\textstyle S^{n-1}\cong\del I^n,~~~~S^n\cong\del I^{n+1}\cong{I^n\over\del I^n}\cong{I_1^n\sqcup I_2^n\over\{x\sim \vphi(x)\},\vphi:\del I_1^n\sr{\cong}{\ral}\del I_2^n},\nn
\eea
every space $Z$ has $\pi_n(Z)=0$, which is a contradiction.
\end{proof}

\begin{lmm}[\textcolor{OliveGreen}{\cite[Theorem 4]{tuffley2002}}]\label{HtyTypeThmI}
$S^1(2k-1)\simeq S^{2k-1}\simeq S^1(2k)$, where ``$\simeq$'' denotes homotopy equivalence.
\end{lmm}

\begin{lmm}[\textcolor{OliveGreen}{\cite[Theorem 5]{tuffley2002}}]\label{HtyTypeThmII}
The map $H_{2k-1}(S^1(2k-1))\cong\Integer\ra H_{2k-1}(S^1(2k))\cong\Integer$ induced by the inclusion $S^1(2k-1)\hookrightarrow S^1(2k)$ is multiplication by $2$.
\end{lmm}

Some of the results of \cite{tuffley2002} can also be found in \cite{wu-wen47}.

\begin{lmm}[\textcolor{OliveGreen}{No retraction III}]\label{NoRetIII}
If $1\leq k\leq n-1$, there is no continuous retraction $S^1(n)\ra S^1(n-k)$.
\end{lmm}
\begin{proof}
Recall that homotopy equivalent spaces have the same homology. We consider three cases follows.
\bit[leftmargin=0.4cm]
\item\ul{Case 1: $2\leqslant k\leqslant n-1$.} Suppose there exists a continuous retraction $r:S^1(n)\ra S^1(n-k)$. Then by taking the homology $H_n$ and applying Lemmas \ref{HlgGpSph},\ref{HtyTypeThmI}, we get $\Integer\cong 0$ (a contradiction).
\item\ul{Case 2: $k=1$ and $n=2m+1$, $m\geqslant 1$.} Suppose there exists a continuous retraction $r:S^1(2m+1)\ra S^1(2m)$. Then by taking the homology $H_{2m}$ and applying Lemmas \ref{HlgGpSph},\ref{HtyTypeThmI}, we get $\Integer\cong 0$ (a contradiction).
\item\ul{Case 3: $k=1$ and $n=2m$, $m\geqslant 1$.} Suppose there is a continuous retraction $r:S^1(2m)\ra S^1(2m-1)$. The maps $S^1(2m-1)\sr{i}{\ral} S^1(2m)\sr{r}{\hookrightarrow}S^1(2m-1)$ and the $H_{2m-1}$-induced maps $\Integer\sr{i_\ast}{\ral}\Integer\sr{r_\ast}{\ral}\Integer$ satisfy $r\circ i=id_{S^1(2m-1)}$, $r_\ast\circ i_\ast=id_{H_{2m-1}(S^1(2m-1))}$. By Lemma \ref{HtyTypeThmII}, $i_\ast$ is multiplication by $2$, which has no left inverse. Thus, the map $r_\ast$ does not exist, and so the retraction $r$ is not continuous (a contradiction).
\eit
\end{proof}

\begin{lmm}\label{LipInductI}
Let $X$ be a metric space. If $Z\subset X$ is a $c$-Lipschitz retract, then so is $Z(n)\subset X(n)$.
\end{lmm}
\begin{proof}
Given a Lipschitz retraction $r:X\ra Z$, let $R:X(n)\ra Z(n)$, $R(x)=\{r(x_1),...,r(x_n)\}$. Then
\begin{align}
&\textstyle d_H\big(R(x),R(y)\big)=\max\left\{\max_i\min_jd(r(x_i),r(y_j)),\max_j\min_id(r(x_i),r(y_j))\right\}\nn\\
&\textstyle~~~~\leq \Lip(r)\max\left\{\max_i\min_jd(x_i,y_j),\max_j\min_id(x_i,y_j)\right\}=\Lip(r)d_H(x,y).\nn
\end{align}
Also, it is clear that if $z\in Z(n)$, then $R(z)=z$.
\end{proof}

\begin{crl}\label{LipInductII}
If $Z\subset X$ is a $\mu$-Lipschitz retract, then any $c$-Lipschitz retraction {\small$r_n:X(n)\ra X(n-1)$} induces a $\mu c$-Lipschitz retraction $\wt{r}_n:Z(n)\ra Z(n-1)$ as in the diagram:
\bea\bt
X(n)\ar[rr,"r_n"] && X(n-1)\ar[d,"R"]\\
Z(n)\ar[u,hook,"i"]\ar[rr,dashed,"\wt{r}_n=R\circ r_n\circ i"] && Z(n-1)
\et\nn\eea
where $R$ is a $\mu$-Lipschitz retraction that exists by Lemma \ref{LipInductI}.
\end{crl}

\begin{rmk}\label{LipInductRmk}
Let $X$ be a locally convex space and $\M\subset X$ a finite-dimensional linear subspace. Then by \cite[Proposition 5.24, p.122]{conway}, $\M$ is complemented in $X$ (i.e., there exists a linear subspace $\N\subset X$ such that $\M\cap\N=\{0\}$, $X=\M+\N$, and the map $\M\times\N\ra X$, $(m,n)\mapsto m+n$ is a homeomorphism).

Thus, if $X$ is a normed space, then by a result of Kadets and Snobar (see \cite[p.335]{plichko-yost2000}) every $n$-dimensional linear subspace $\M\subset X$ is a $\sqrt{n}$-Lipschitz retract (an improvement of the bound in Lemma \ref{AuLmmBnd}). Moreover, it follows from \cite[Theorem 5.1]{phelps57} that every closed convex subset $C\subset E$ of an inner product space $E$ is a 1-Lipschitz retract through the nearest point projection $P:E\ra C$, $x\mapsto Px$, where $Px$ is given by $\|x-Px\|=d(x,C)$. In particular, if $\H$ is a Hilbert space, then every $n$-dimensional linear subspace $\M\subset \H$ is a $1$-Lipschitz retract.
\end{rmk}

\begin{thm}[\textcolor{OliveGreen}{Lower bound on growth rate of the Lipschitz constant}]\label{LipConstGr}
Let $X$ be a normed space such that $\dim X\geq 2$. If $r_n:X(n)\ra X(n-1)$ is a Lipschitz retraction, then $\Lip(r_n)\geq {n\over 4\pi\sqrt{2}}$ (or $\Lip(r_n)\geq {n\over 4\pi}$ if $X$ is a Hilbert space).
\end{thm}
\begin{proof}
Let $r_n:X(n)\ra X(n-1)$ be a Lipschitz retraction. Then by Corollary \ref{LipInductII} and Remark \ref{LipInductRmk}, we get a $\sqrt{2}\Lip(r_n)$-Lipschitz retraction $R_n:\Real^2(n)\ra \Real^2(n-1)$.

 Suppose $\Lip(r_m)<{1\over 4\pi\sqrt{2}}m$, i.e., $\Lip(R_m)<{1\over 4\pi}m$, for an integer $m\geq 1$. Let $x\in S^1(m)\backslash S^1(m-1)$, where $S^1\subset\Real^2$ is the unit circle centered at $0\in\Real^2$. Then
\bea
\textstyle\delta(x)\leq{2\pi\over|x|}={2\pi\over m}.\nn
\eea
Since $\dist_H\big(x,\Real^2(m-1)\big)={1\over 2}\delta(x)$ by Lemma \ref{PointProxLmm}, and $\Real^2(m-1)\subset\Real^2(m)$ is closed, there exists $y_x\in \Real^2(m-1)$ such that $d_H(x,y_x)={1\over 2}\delta(x)$. This implies for all $x\in S^1(m)\backslash S^1(m-1)$,
\begin{align}
&\textstyle d_H\big(R_m(x),y_x\big)=d_H\big(R_m(x),R_m(y_x)\big)\leq \Lip(R_m)d_H(x,y_x)={\Lip(R_m)\over 2}\delta(x),\nn\\
&\textstyle~~\Ra~~d_H\big(R_m(x),x\big)\leq[1+\Lip(R_m)]d_H(x,y_x)\leq 2\Lip(R_m){\delta(x)\over 2}<{1\over 2},\nn\\
&\textstyle~~\Ra~~R_m(x)\in A^1(m-1),~~~~A^1:=\big\{u\in\Real^2:1/2<\|u\|< 3/2\big\}\cong S^1\times[0,1].\nn
\end{align}
It follows that we have a Lipschitz map $R_m|_{S^1(m)}:S^1(m)\ra A^1(m-1)$. Composing this with the Lipschitz retraction $\pi_{m-1}:A^1(m-1)\ra S^1(m-1)$ given by $\pi_{m-1}(x):=\left\{{a\over\|a\|}:a\in x\right\}$, we get a Lipschitz retraction
\bea
\wt{R}_m=\pi_{m-1}\circ R_m|_{S^1(m)}:S^{1}(m)\sr{R_m}{\ral}A^1(m-1)\sr{\pi_{m-1}}{\ral}S^1(m-1).\nn
\eea
This contradicts the conclusion of Lemma \ref{NoRetIII}.
\end{proof}

Since homology only (and not homotopy) was used, Lemmas \ref{HtyRetInc},\ref{HtyGpSph} are not relevant. To use homotopy only (and not homology), we simply need to replace Lemmas \ref{HomRetInc},\ref{HlgGpSph} with Lemmas \ref{HtyRetInc},\ref{HtyGpSph}.

\section{Some non-existence results and questions on the FSR property}\label{FSRP2sn}
The following result is based on the fact from Lemma \ref{HausDistBound} and Corollary \ref{HausDistBound1} that locally, the subspace $DX(n):=X(n)\backslash X(n-1)$ of $X(n)$ looks like a subspace of $X^n$.

\begin{lmm}[\textcolor{OliveGreen}{Local Lipschitz decomposition}]\label{LocLipDecI}
Let $\vphi:Z\ra X(n)$ be a Lipschitz map and $z_0\in Z$. If $\vphi(z_0)\in X(n)\backslash X(n-1)$, then for some $\vep>0$, there exist Lipschitz maps $\vphi_1,...,\vphi_n:N_\vep(z_0)\subset Z\ra X$ such that $\Lip(\vphi_i)\leq\Lip(\vphi)$ for all $i=1,...,n$, and $\vphi(z)=\big\{\vphi_1(z),...,\vphi_n(z)\big\}$ for all $z\in N_\vep(z_0)$.
\end{lmm}
\begin{proof}
Fix a number $\tau>4$. Recall that for all $x,y\in X(n)$, if $\delta(x)>2d_H(x,y)$ or $\delta(y)>2d_H(x,y)$, then there exists an enumeration $x_i,y_i$ of elements of $x,y$ such that
\begin{equation*}
d(x_i,y_i)\leq d_H(x,y),~~~~~~\txt{for all}~~i=1,...,n.
\end{equation*}
This means that if $x_0=\big\{(x_0)_1,...,(x_0)_n\big\}\in X(n)\backslash X(n-1)$, then we have $1$-Lipschitz maps
{\small\begin{align}
f_i:N_{\delta(x_0)\over \tau}(x_0)\subset X(n)\ra X,~~f_i(x):=x_i~~(\txt{the unique element of $x$ with $d\big(x_i,(x_0)_i\big)\leq d_H(x,x_0)$}),\nn
\end{align}}satisfying $id_{X(n)}(x)=x=\{f_1(x),...,f_n(x)\}$ for all $x\in N_{\delta(x_0)\over \tau}(x_0)$.  That is, the identity $id_{X(n)}:X(n)\ra X(n)$ decomposes locally (on $X(n)\backslash X(n-1)$) into the $1$-Lipschitz maps $f_i$.

(Note that $f_i$ is Lipschitz, with $d(f_i(x),f_i(y))=d(x_i,y_i)\leq d_H(x,y)$, because we know there is $y(i)\in y$ satisfying $d(x_i,y(i))\leq d_H(x,y)$ for each $i$, and since the neighborhoods {\footnotesize$\left\{N_{\delta(x_0)\over \tau}\big((x_0)_i\big):i=1,...,n\right\}$} are distantly separated, it follows that $y(i)=y_i:=f_i(y)$.)

Thus, if $\vphi:Z\ra X(n)$ is a Lipschitz map, then with a fixed $z_0\in Z$ satisfying {\small$x_0:=\vphi(z_0)\in X(n)\backslash X(n-1)$}, we see that on some open set {\footnotesize $N_\vep(z_0)\subset \vphi^{-1}\left(N_{\delta(x_0)\over \tau}(x_0)\right)$},
{\small\bea
\vphi_i=f_i\circ\vphi:N_\vep(z_0)\subset Z\sr{\vphi}{\ral}N_{\delta(x_0)\over \tau}(x_0)\subset X(n)\sr{f_i}{\ral} X,\nn
\eea} are Lipschitz maps that (i) decompose $\vphi$ as $\vphi(z)=\{\vphi_1(z),...,\vphi_n(z)\}$ for all $z\in N_\vep(z_0)$, and (ii) satisfy ~$\Lip(\vphi_i)=\Lip(f_i\circ\vphi)\leq\Lip(f_i)\Lip(\vphi)\leq\Lip(\vphi_i)$.
\end{proof}

\begin{crl}[\textcolor{OliveGreen}{Local Lipschitz decomposition}]\label{LocLipDecII}
Let $\vphi:Z\ra X(n)$ be a Lipschitz map and $z_0\in Z$. Fix $1\leq k\leq n$ and let $Z_k:=\vphi^{-1}\big(X(k)\big)$. If $\vphi(z_0)\in X(k)\backslash X(k-1)$, then for some $\vep>0$, there exist Lipschitz maps $\vphi_1,...,\vphi_k:N_\vep(z_0)\cap Z_k\subset Z\ra X$ such that $\Lip(\vphi_i)\leq\Lip(\vphi)$ for all $i=1,...,k$, and $\vphi(z)=\big\{\vphi_1(z),...,\vphi_k(z)\big\}$ for all $z\in N_\vep(z_0)\cap Z_k$.

It follows that for any $z_0\in Z$, there exists $\vep=\vep(z_0)>0$ and Lipschitz maps
\bea
\vphi_1,...,\vphi_{|\vphi(z_0)|}:Z_{|\vphi(z_0)|}\ra X,~~~~Z_{|\vphi(z_0)|}:=N_\vep(z_0)\cap\vphi^{-1}\Big(X\big(|\vphi(z_0)|\big)\Big),\nn
\eea
satisfying $\Lip(\vphi_i)\leq\Lip(\vphi)$ for all $i=1,...,|\vphi(z_0)|$ and $\vphi(z)=\big\{\vphi_1(z),...,\vphi_{|\vphi(z_0)|}(z)\big\}$ for all $z\in Z_{|\vphi(z_0)|}$.
\end{crl}
\begin{proof}
Apply the result of Lemma \ref{LocLipDecI} to the map $\psi=f|_{Z_k}:Z_k\ra X(k)$.
\end{proof}

\begin{prp}\label{NoLipRetCon}
If $X$ is a (locally) contractible metric space, there need not exist Lipschitz retractions $r:X(n)\ra X(n-1)$ for all $n\geq 2$.
\end{prp}
\begin{proof}
Let $X\subset\Real^2$ be the union of the interval $E=[(0,0),(1,0)]$ and a snowflake curve (with no self intersections) $K\subset\Real^2$ from $(1,0)$ to $(2,0)$, i.e., $X=L\cup K$. Suppose $r:X(4)\ra X(3)$ is a Lipschitz retraction. With the origin $a=(0,0)$, a variable point $b(t)=(t,0)\in L$, $0\leq t\leq 1$, and distinct nearby points $c,d\in K\backslash\{(1,0),(2,0)\}$, consider the path $\gamma:[0,1]\ra X(4),~t\mapsto\{a,b(t),c,d\}$. Then we have a Lipschitz map
\bea
\vphi:=r\circ\gamma:[0,1]\sr{\gamma}{\ral}X(4)\sr{r}{\ral}X(3),~t\mapsto r\big(\{a,b(t),c,d\}\big).\nn
\eea
Note that $\vphi(0)=r\big(\{a,c,d\}\big)=\{a,c,d\}\in X(3)\backslash X(2)$. By Lemma \ref{LocLipDecI}, there exists $\vep>0$ such that for $0\leq t<\vep$, $\vphi$ decomposes into Lipschitz maps $\vphi_1,\vphi_2,\vphi_3:[0,1]\ra X$ as
\bea
\vphi(t)=\{\vphi_1(t),\vphi_2(t),\vphi_3(t)\},~~~~\vphi_1(0)=a,~~\vphi_2(0)=c,~~\vphi_3(0)=d.\nn
\eea
By Corollary \ref{SnwfNoLipCrv}, $\vphi_2,\vphi_3$ are constant maps, with $\vphi_2(t)\eqv c$, $\vphi_3(t)\eqv d$. It follows that
\bea
\vphi(t)=\{\vphi_1(t),c,d\},~~~~\txt{for all}~~t\in[0,1].\nn
\eea
With $c,d$ chosen such that $\dist\big(\{a,b(1)\},\{c,d\}\big)>\max\big(\|a-b(1)\|,\|c-d\|\big)$ (i.e., distant separation), it follows that when $x:=\vphi(1)=\{a,b(1),c,d\}\in X(4)$,
\bea
&&\label{DBeqn}\textstyle d_H\big(r(x),x\big)=d_H\Big(\{\vphi_1(1),c,d\},\{a,b(1),c,d\}\Big)=d_H\Big(\{a,b(1)\},\{\vphi_1(1)\}\Big)\nn\\
&&\textstyle~~~~=\max\Big\{\|\vphi_1(1)-a\|,\|\vphi_1(1)-b(1)\|\Big\}\geq{\|a-b(1)\|\over 2}={1\over 2},
\eea
where the last inequality is due to the triangle inequality.

On the other hand, by Corollary \ref{DispBound},
\bea
\textstyle d_H\big(r(x),x\big)\leq{1+\Lip(r)\over 2}\delta(x)\leq{1+\Lip(r)\over 2}\|c-d\|.\nn
\eea
This contradicts (\ref{DBeqn}), since we can choose $\|c-d\|$ to be as small as we wish.
\end{proof}

Note that in the proof of Proposition $\ref{NoLipRetCon}$, if only local contractibility is required, then the line $L\subset\Real^2$ and the snowflake $K\subset\Real^2$ in $X=L\cup K$ can be chosen to be disjoint (i.e., such that $L\cap K=\emptyset$) without altering the proof.

\begin{lmm}\label{SnFkFSRP}
Let $X$ be a metric space. (Local) Lipschitz contractibility of $X$ is not essential for the existence of Lipschitz retractions $X(n)\ra X(n-1)$.
\end{lmm}
\begin{proof}
Let $K$ be the classical von Koch snowflake curve based on the interval {\small$\big[(0,0),(1,0)\big]$} $\subset\Real^2$. Then by \cite{tukia81,koskela94,mckemie1987}, there exist a natural homeomorphism $\vphi:[0,1]\ra K$ and a constant $c$ such that
\bea
|t-s|^{1\over p}/c\leq\|\vphi(t)-\vphi(s)\|\leq c|t-s|^{1\over p},~~~~\txt{for all}~~s,t\in[0,1],\nn
\eea
where $p:={\log 4\over\log 3}$. Define $\phi:[0,1](n)\ra K(n)$ by $\phi(x)=\{\vphi(x_1),...,\vphi(x_n)\}$. If $r:[0,1](n)\ra[0,1](n-1)$ is a Lipschitz retraction, then we get a Lipschitz retraction
\bea
R=\phi\circ r\circ\phi^{-1}:K(n)\sr{\phi^{-1}}{\ral}[0,1](n)\sr{r}{\ral}[0,1](n-1)\sr{\phi}{\ral}K(n-1),\nn
\eea
where using {\footnotesize$d_H(R(x),R(y)):=\max\Big(\max\limits_i\min\limits_j\|R(x)_i-R(y)_j\|, \max\limits_j\min\limits_i\|R(x)_i-R(y)_j\|\Big)$}, we get
\bea
&&d_H(R(x),R(y))\leq cd_H\big(r\circ\phi^{-1}(x),r\circ\phi^{-1}(y)\big)^{1\over p}\leq c\Lip(r)^{1\over p}d_H\big(\phi^{-1}(x),\phi^{-1}(y)\big)^{1\over p}\nn\\
&&~~~~\leq c^2\Lip(r)^{1\over p}d_H(x,y).\nn
\eea
However, we know that a snowflake is not (locally) Lipschitz contractible.
\end{proof}

Note that the existence of Lipschitz retractions $[0,1](n)\ra[0,1](n-1)$ follows from Remark \ref{LipInductRmk}, since $[0,1]\subset\Real$ is a 1-Lipschitz retract. Alternatively, for a Hilbert space $\H$, we have Lipschitz retractions $r:\H(n)\ra\H(n-1)$ satisfying $r(x)\subset\Conv(x)$, and so every convex subset $K\subset\H$ inherits Lipschitz retractions $K(n)\ra K(n-1)$.

\begin{lmm}[\textcolor{OliveGreen}{Componentwise Lipschitzness and continuity}]\label{ComLipCont}
Let $(X,d_1),(Y,d_2)$ be metric spaces and $(X\times Y,d)$ the product metric space, i.e., $d\big((x,y),(x',y')\big):=\max\big(d_1(x,y),d_2(x',y')\big)$. If $f:X\times Y\ra Z$ is a map such that $f^y:X\ra Z$, $x\mapsto f(x,y)$ is $L^y$-Lipschitz, and $f_x:Y\ra Z$, $y\mapsto f(x,y)$ is $L_x$-Lipschitz, then $f$ is continuous. Moreover, $f$ is Lipschitz if $L^y$ and $L_x$ are constants.
\end{lmm}
\begin{proof} Observe that
~$d\big(f(x,y),f(x',y')\big)\leq d\big(f(x,y),f(x,y')\big)+d\big(f(x,y'),f(x',y')\big)$ $\leq$ $\big(L_x+L^y\big)d\big((x,y),(x',y')\big)$.
\end{proof}

\begin{example}[\textcolor{OliveGreen}{Geometry of a cusp}]\label{CuspLipCont}
Let $X:=\big\{(s,|s|^{1\over 2}):s\in[-1,1]\big\}\subset\Real^2$. Then $X$ is
\bit[leftmargin=0.8cm]
\item[(i)] \ul{Homeomorphic to a line segment}: The map $C:[-1,1]\ra X$, $s\mapsto\big(s,|s|^{1\over 2}\big)$ is a homeomorphism.
\item[(ii)] \ul{Lipschitz-contractible}: The map $H:X\times I\ra X$, {\small $H\big((s,|s|^{1\over 2}),t\big):=\big(t^2s,|t^2s|^{1\over 2}\big)=\big(t^2s,t|s|^{1\over 2}\big)$}, by Lemma \ref{ComLipCont}, is a $\big(1+5^{1\over 2}\big)$-Lipschitz homotopy from $id_X$ to the constant map $X\ra(0,0)$.
\item[(iii)] \ul{Not quasiconvex}: Given any $s\in[-1,1]$ and a curve $\gamma_s:[0,1]\ra X$ between $a_-:=\big(-s,|s|^{1\over 2}\big)$ and $a_+:=\big(s,|s|^{1\over 2}\big)$, we have
\begin{align}
\textstyle{l(\gamma_s)\over\|a_+-a_-\|}>{\|a_+-(0,0)\|+\|a_--(0,0)\|\over\|a_+-a_-\|}={2|s|^{1\over 2}(1+|s|)^{1\over 2}\over 2|s|}=\left({1+|s|\over|s|}\right)^{1\over 2}\ral\infty~~\txt{as}~~s\ra 0.\nn
\end{align}
This means the cusp is not bi-Lipschitz equivalent to a line segment since a bi-Lipschitz map preserves quasiconvexity.
\eit
\end{example}

\begin{question}
If $X$ is the cusp in Example \ref{CuspLipCont}, do we have Lipschitz retractions $X(n)\ra X(n-1)$?
\end{question}

\begin{question}
If $X$ is a locally Lipschitz-contractible quasiconvex space, do we have Lipschitz retractions $X(n)\ra X(n-1)$ for all $n\geq 2$?
\end{question}

\begin{question}
Let $X,Y$ be metric spaces such that $X\simeq_LY$. If Lipschitz retractions $X(n)\ra X(n-1)$ exist, do we get Lipschitz retractions $Y(n)\ra Y(n-1)$?

Given a Lipschitz retraction $r_n:X(n)\ra X(n-1)$, we get the diagram
\bea\bt
X(n)\ar[rrrr,"r_n"] &&&& X(n-1)\ar[d,"F"]\\
Y(n)\ar[u,"G"]\ar[rrrr,dashed,"s_n~=~F\circ r_n\circ G"] &&&& Y(n-1)
\et\nn\eea
in which the induced Lipschitz map $s_n:Y(n)\ra Y(n-1)$ satisfies
\bea
s_n|_{Y(n-1)}=F\circ G|_{Y(n-1)}\simeq_Lid_{Y(n-1)}.\nn
\eea
This shows that the map $H_k(Y(n-1))\sr{i_\ast}{\ral}H_k(Y(n))$ induced by the inclusion $Y(n-1)\sr{i}{\hookrightarrow}Y(n)$ has a left inverse $H_k(Y(n))\sr{(s_n)_\ast}{\ral}H_k(Y(n-1))$, i.e., $(s_n)_\ast\circ i_\ast=id_{H_k(Y(n-1))}$.

\end{question}

\section{Counterexamples and other facts}\label{FSRP2cf}
\begin{dfn*}[\textcolor{blue}{Recall: Connected space}]
A space $X$ is connected if it cannot be written as a disjoint union $X=A\sqcup B$ of nonempty open sets $A$ and $B$.
\end{dfn*}

\begin{rmk}\label{LipConSuffRmk}
If $X$ is Lipschitz $k$-connected for all $k\geq 0$, it does not follow that $X$ is boundedly conically $\ld$-Lipschitz contractible. For example, consider the (connected but \ul{not path-connected}) closed \emph{sine-curve} $S:=\{(0,0)\}\cup\big\{(u,\sin(\pi/u)):u\in(0,1]\big\}$, and let $X\subset\Real^2$ be the (locally connected but \ul{not locally path-connected}) \emph{Polish circle} given by
\bea
X:=S\cup[(1,0),(1,-2)]\cup[(1,-2),(0,-2)]\cup[(0,-2),(0,1)].\nn
\eea
Let $S_n:=\{(0,0)\}\cup\big\{(u,\sin(\pi/u)):u\in(0,1/n)\big\}$.

Topologically, $X$ is $k$-connected for all $k\geq 0$ but is not contractible (\cite[Example 5.1.40, p.164]{aguilar.etal2002}): Indeed, no nontrivial closed path exists in $X$. So, since $S^k$ is compact and locally path-connected, the image of a continuous map $f:S^k\ra X$ cannot contain the non-locally path-connected neighborhood $N_\vep\big((0,0)\big)$. Similarly, a contracting homotopy $H:X\times I\ra X$ with $H(x,0)=(0,0)$ and $H(x,1)=x$ will restrict to a contracting homotopy $H:N_\vep\big((0,0)\big)\times I\ra N_\vep\big((0,0)\big)$, which is a contraction since the sine-curve $S_n$ is not contractible.

Suppose on the contrary $X$ is boundedly conically $\ld$-Lipschitz contractible. Then because $X$ is bounded (hence conically $\ld$-Lipschitz contractible), it follows from Lemma \ref{ConicLipHom1} that there exists a $\ld$-Lipschitz homotopy $H:X\times I\ra X$, where $H(x,1)=x$ and $H(x,0)=x_0\in X$. Thus, for each $x_\tau:=(\tau,\sin(\pi/\tau))\in X$, $\tau\in(0,1]$, the curve $\gamma_{x_0,x_\tau}(t)=H(x_\tau,t)$ from $x_0$ to $x_\tau$ satisfies
\bea
\|\gamma_{x_0,x_\tau}(t)-\gamma_{x_0,x_\tau}(t')\|\leq \ld|t-t'|~~\big(~~\Ra~~l(\gamma_{x_0,x_\tau})\leq \ld~~\big)~~\txt{for all}~~\tau\in(0,1],\nn
\eea
and so has length $\leq \ld$ (a contradiction). Hence, $X$ is \ul{not conically Lipschitz contractible}.

It is also true that $X$ is \ul{not quasiconvex}, otherwise, we get a contradiction in the form
\bea
&&\|\gamma_{x_0,x_\tau}(t)-\gamma_{x_0,x_\tau}(t')\|\leq \ld \|\gamma_{x_0,x_\tau}(0)-\gamma_{x_0,x_\tau}(1)\||t-t'|=\ld \|x_0-x_\tau\||t-t'|\nn\\
&&~~~~\leq 4\ld |t-t'|~~\big(~~\Ra~~l(\gamma_{x_0,x_\tau})\leq 4\ld~~\big)~~\txt{for all}~~\tau\in(0,1].\nn
\eea
As shown in Lemma \ref{LipConSuffLmm} below, $X$ is Lipschitz $k$-connected for all $k\geq 0$.
\end{rmk}
\begin{lmm}\label{LipConSuffLmm}
The space $X$ in Remark \ref{LipConSuffRmk} is Lipschitz $k$-connected for all $k\geq 0$.
\end{lmm}
\begin{proof}
Let $f:S^k\ra X$ be Lipschitz. It is clear by quasiconvexity that $f:S^0\ra X$ has a Lipschitz extension $F:B^1\ra X$. So, assume $k\geq 1$. To show $f$ has a Lipschitz extension $F:B^{k+1}\ra X$, it suffices to show that $f(S^k)$ is an arc of length $\leq 2\pi\Lip(f)$ that is a locally path-connected subset of $X_\vep:=X\backslash N_\vep(\{(0,0)\})$, for some $\vep>0$. $f(S^k)$ is indeed locally path-connected (and so some $N_\vep(\{(0,0)\})\not\subset f(S^k)$) since $S^k$ is locally path-connected and $f$ is continuous. To show $f(S^k)\subset X_\vep$ (for some $\vep>0$), observe that if we fix two distinct points $s_0,s_1\in S^k$, then every other point $s\in S^k$ lies on a path $\gamma_s:[0,1]\ra S^k$ from $s_0$ to $s_1$ of length $l(\gamma_s)\leq\pi$. It follows that for every point $s\in S^k$, the point $f(s)\in X_\vep$ lies on the path $c_s=f\circ\gamma_s:[0,1]\sr{\gamma_s}{\ral}\gamma([0,1])\sr{f}{\ral}X_\vep$ from $f(s_0)$ to $f(s_1)$ of length $l(c_s)\leq\pi\Lip(f)$. Hence, $f(S^k)=\bigcup_{s\in S^k}f\big(\gamma([0,1])\big)\subset X_\vep$, and so $f(S^k)$ is an arc such that ~$\txt{Length}\big(f(S^k)\big)\leq\sup_{s\in S^k} l(c_s)\leq 2\pi\Lip(f)$.
\end{proof}

\begin{dfn}[\textcolor{blue}{\index{Simple curve}{Simple curve}}]
A subset $\Gamma\subset\Real^2$ that is homeomorphic to $[0,1]$, i.e., $\Gamma\subset\Real^2$ and $\Gamma\cong[0,1]$.
\end{dfn}

\begin{lmm}
Let $\Gamma$ be a simple curve. If $\Gamma$ is biLipschitz equivalent to $[0,1]$, then $\Gamma$ is Lipschitz contractible.
\end{lmm}
\begin{proof}
If $\vphi:[0,1]\ra\Gamma$ is a biLipschitz map and $H:[0,1]\times I\ra [0,1]$ is a contracting Lipschitz homotopy, then we get a contracting Lipschitz homotopy
\bea
\vphi^{-1}\circ H\circ(\vphi\times id):\Gamma\times I\sr{\vphi\times id}{\ral}[0,1]\times I\sr{H}{\ral}[0,1]\sr{\vphi^{-1}}{\ral}\Gamma.\nn\qedhere
\eea
\end{proof}
The converse of this result is not true. Recall that the cusp $\gamma:[-1,1]\ra\Real^2$, $\gamma(t)=\big(t,|t|^{1\over 2}\big)$ is Lipschitz-contractible but not quasiconvex, and so not biLipschitz equivalent to $[0,1]$, because a biLipschitz map preserves quasiconvexity.

\begin{lmm}[\textcolor{OliveGreen}{Conical Lipschitz contractibility of the cusp}]
The cusp $X=\big\{(s,|s|^{1\over 2}):s\in[-1,1]\big\}\subset\Real^2$ is (conically) Lipschitz contractible.
\end{lmm}
\begin{proof}
The case with ``conically'' removed was proved earlier in Example \ref{CuspLipCont}(ii). So, we will now prove $X$ is conically Lipschitz contractible. Consider the map
\bea
r:\Cone(X)\ra X,~~~~r\big((s,|s|^{1\over 2}),t\big):=\big(t^2s,|t^2s|^{1\over 2}\big)=\big(t^2s,t|s|^{1\over 2}\big).\nn
\eea
Let $x:=\big(s,|s|^{1\over 2}\big)$ and $x':=\big(s',|s'|^{1\over 2}\big)$. Then
\bea
\|r(x,t)-r(x,t')\|\leq 5^{1\over 2}|t-t'|~~~~\txt{and}~~~~\|r(x,t)-r(x',t)\|\leq t\|x-x'\|.\nn
\eea
Thus, using the triangle inequality, we get
\bea
\|r(x,t)-r(x',t')\|\leq 5^{1\over 2}|t-t'|+\min(t,t')\|x-x'\|\leq 5^{1\over 2}d_{c2}\big((x,t),(x',t')\big),\nn
\eea
where ~$d_{c2}\big((x,t),(x',t')\big)=|t-t'|+\min(t,t')\|x-x'\|$.
\end{proof}

\begin{dfn*}[\textcolor{blue}{Recall: Natural parametrization}]
Let $X$ be a metric space. A path $\gamma:[a,b]\ra X$ is a natural parametrization if $l(\gamma([t,t']))=|t-t'|$ for all $t,t'\in[a,b]$.
\end{dfn*}

\begin{dfn}[\textcolor{blue}{\index{Chord-arc condition}{Chord-arc condition}}]\label{ChordArcCond}
A curve $\Gamma\subset\Real^2$ is chord-arc if there exists a number $\ld\geq 0$ such that for any two points $a,b\in\Gamma$, the arc $\Gamma_{a,b}\subset\Gamma$ from $a$ to $b$ has length
\bea
l(\Gamma_{a,b})\leq \ld\|a-b\|.\nn
\eea
\end{dfn}

\begin{crl}[\textcolor{OliveGreen}{Recall of Corollary \ref{ConstSpeedRP}: Constant speed parametrization}]\label{ConSpdRPRcll}
Let $X$ be a metric space. If $\gamma:[a,b]\ra X$ is a rectifiable curve, there exists a curve $\eta:[a,b]\ra X$, such that $\eta([a,b])=\gamma\big([a,b]\big)$, $\eta(a)=\gamma(a)$, $\eta(b)=\gamma(b)$, and
\bea
\textstyle d(\eta(t),\eta(t'))\leq l(\eta([t,t']))={l(\gamma)\over b-a}|t-t'|,~~\txt{for all}~~t,t'\in[a,b].\nn
\eea
\end{crl}

\begin{lmm}
A curve $\Gamma\subset\Real^2$ is biLipschitz equivalent to $[0,1]$ $\iff$ it satisfies the chord-arc condition (Definition \ref{ChordArcCond}).
\end{lmm}
\begin{proof}
($\Ra$): Let $\vphi:[0,1]\ra\Gamma$ be biLipschitz, i.e., $|t-t'|/c\leq\|\vphi(t)-\vphi(t')\|\leq c|t-t'|$. If $a,b\in\Gamma$, let $a=\vphi(t_1)$, $b=\vphi(t_2)$. Then
\bea
l(\Gamma_{a,b})=l\big(\vphi([t_1,t_2])\big)\leq c|t_1-t_2|\leq c^2\|a-b\|.\nn
\eea

($\La$): Assume $\Gamma$ is chord-arc. Then there exists $\ld\geq 0$ such that $l(\Gamma_{a,b})\leq \ld$ for all $a,b\in\Gamma$. Thus, $\Gamma$ is rectifiable. Let $\vphi:[0,1]\ra\Gamma$ be a parametrization with constant speed $c$, i.e.,  $l(\vphi([t,t']))=c|t-t'|$ for all $t,t'\in[0,1]$. (This exists by Corollary \ref{ConSpdRPRcll}). Then,
\bea
\|\vphi(t)-\vphi(t')\|\leq c|t-t'|=l\big(\vphi([t,t'])\big)\sr{\txt{chord-arc}}{\leq}\ld\|\vphi(t)-\vphi(t')\|,~~~~\txt{for all}~~~~t,t'\in[0,1].\nn\qedhere
\eea
\end{proof}
Note that we have the following containments.
\bea
\left\{\substack{\txt{chord-arc}\\\txt{curves}}\right\}\subsetneq\left\{\substack{\txt{Lipschitz-Contractible}\\\txt{curves}}\right\}\subsetneq\left\{\substack{\txt{Rectifiable}\\\txt{curves}}\right\}.\nn
\eea

\begin{lmm}[\textcolor{OliveGreen}{Which simple curves are Lipschitz contractible?}]
If a simple curve $\Gamma\subset\Real^2$ is quasiconvex and rectifiable, then it is Lipschitz contractible.
\end{lmm}
\begin{proof}
It suffices to show that $\Gamma\sr{\txt{biLip}}{\cong}[0,1]$. Since $\Gamma$ is rectifiable, there exists a parametrization $\gamma:[0,1]\ra\Gamma$ with constant speed $c$. Thus,
\bea
\|\gamma(t)-\gamma(t')\|\leq l(\gamma([t,t']))=c|t-t'|~~~~\txt{for all}~~~~t,t'\in[0,1].\nn
\eea
On the other hand, since $\Gamma$ is quasiconvex, we have $l(\gamma([t,t']))\leq \ld\|\gamma(t)-\gamma(t')\|$, for some $\ld\geq 1$.
\end{proof}      

{\appendix \appendixpage \noappendicestocpagenum \addappheadtotoc
\chapter{Preliminaries}\label{Prelims} 
In this chapter we introduce basic mathematical concepts, as well as set up associated notation/terminology, that will be relevant in our subsequent discussion of finite subset spaces (to begin in chapter \ref{FSS}). In an effort to lighten notation, any refinements required in introducing the counterparts of these concepts for finite subset spaces will be discussed wherever they arise in later chapters. In order to make the discussion self contained, we will attempt to provide proof outlines for all essential/relevant results.     
\section{Topological Space Concepts}\label{PrelimsTSC} 

\subsection{Topological spaces, Standard examples, Space description tools}
\begin{dfn}[\textcolor{blue}{Powerset, \textcolor{magenta}{(Point-set) Topology}, \index{Topological! space}{(Topological) space}, \index{Open set}{Open set}, \index{Closed set}{Closed set}}]\label{PtSetTop}
Let $X$ be a set and $\P(X)$, or $2^X$, the collection of all subsets (called the \ul{powerset}) of $X$. A \ul{topology} $\T\subset\P(X)$ on $X$ (making {\small $X=(X,\T)$} a \ul{topological space} or simply a \ul{space}) is a family of subsets (called \ul{open sets} of the topology) with the following properties.
\begin{enumerate}[leftmargin=0.9cm]
\item $\T$ contains both the empty set $\emptyset$ and the whole set $X$: i.e., $\emptyset,X\in\T$.
\item $\T$ is closed under finite intersections: i.e., if $U,V\in \T$, then $U\cap V\in\T$.
\item $\T$ is closed under arbitrary unions: i.e., if $\{U_\al:\al\in A\}\subset\T$, then $\bigcup_{\al\in A}U_\al\in\T$.
\end{enumerate}
The complements $\{C=O^c:O\in\T\}$ are called \ul{closed sets} of the topology, where for any set $A\subset X$, its \ul{complement} $A^c=X\backslash A=X-A:=\{x\in X:x\not\in A\}$ consists of all elements of $X$ that are not in $A$.
\end{dfn}
There exist many equivalent ways of specifying the topology of a space, some of which involve distances or metrics (Definition \ref{MetricSp}). The description of a topology without reference to distances or metrics (such as in Definition \ref{PtSetTop}) is often called \ul{\emph{point-set topology}}.

\begin{example*}[\textcolor{blue}{Discrete topology, \index{Discrete space}{Discrete space}, Indiscrete topology, \index{Indiscrete space}{Indiscrete space}}]
The topology $\T$ of a space $X$ is \ul{discrete} (making {\small $(X,\T)$} a discrete space) if every point $x_0\in X$ is an open set $\{x_0\}\subset X$ [equivalently, every subset of $X$ is an open set, i.e., {\small $\T=\P(X)$} ]. The topology $\T$ of a space $X$ is \ul{indiscrete} (making {\small $(X,\T)$} an indiscrete space) if $\T=\{\emptyset,X\}$.
\end{example*}

\begin{lmm}[\textcolor{OliveGreen}{Open set criterion}]\label{OpenCrit}
In a space $X$, a set $A\subset X$ is open $\iff$ for every point $x\in A$, there is an open set $O$ such that $x\in O\subset A$.
\end{lmm}
\begin{proof}
If $A$ is open, then for every $x\in A$, we have $x\in O:=A\subset A$. Conversely, if for every $x\in A$, we have $x\in O_x\subset A$ for some open set $O_x$, then $A=\bigcup_{x\in A}O_x$ is open as a union of open sets.
\end{proof}

\begin{dfn}[\textcolor{blue}{\index{Metric}{Metric}, \index{Ball}{Ball}, \textcolor{magenta}{Metric topology}, \index{Metric space}{Metric space}}]\label{MetricSp}
Let $X$ be a set. A \ul{metric} on $X$ is a function $d:X\times X\ra\Real$ such that for all $x,y,z\in X$,
\bit
\item[] $d(x,y)=d(y,x)\leq d(x,z)+d(z,y)$, and $d(x,y)=0$ iff $x=y$.
\eit
Given a metric $d$ on $X$, the \ul{($d$-) ball} of radius $R>0$ centered at a point $x$ in $X$ is the set
\bea
B_R(x)=B_R^d(x):=\{x'\in X:d(x,x')<R\}.\nn
\eea
Given a metric $d$ on $X$, the \ul{metric ($d$-) topology} on $X$ (making $X=(X,d)$ a \ul{metric space}) is the topology on $X$ with open sets defined as follows: A (nonempty) set $A\subset X$ is $d$-open if for every $x\in A$, there exists a $d$-ball $B_R(x)\subset A$.
\end{dfn}
Note that a function $d:X\times X\ra\Real$ satisfying $d(x,y)=d(y,x)\leq d(x,z)+d(z,y)$ is called a \index{Pseudometric}{\ul{pseudometric}}, and becomes a metric only if it also satisfies $d(x,y)=0$ iff $x=y$.

\begin{rmk}[\textcolor{OliveGreen}{A metric space is a topological space}]
Observe that in a metric space $X=(X,d)$, any $d$-ball $B_R(x)\subset X$ is open because for each $y\in B_R(x)$, we have $B_{R-d(x,y)}(y)\subset B_R(x)$. It follows that a set $A\subset X$ is $d$-open $\iff$ it is a union of $d$-balls. Also, the intersection $B_R(x)\cap B_{R'}(x')$ of two $d$-balls is open because for each $y\in B_R(x)\cap B_{R'}(x')$,
\bit
\item[] $B_{R(x,x')}(y)\subset B_R(x)\cap B_{R'}(x')$, ~\txt{where} ~ $R(x,x'):=\min\big(R-d(x,y),R'-d(x',y)\big)$.
\eit
It follows that finite intersections of $d$-open sets are $d$-open sets, since
\bit
\item[]$\big(\bigcup B_{R_\al}(x_\al)\big)\cap\big(\bigcup B_{\beta}(x_\beta)\big)=\bigcup_{\al,\beta}~B_{R_\al}(x_\al)\cap B_{R_\beta}(x_\beta)$.
\eit
Hence, $d$-open sets satisfy the three properties in Definition \ref{PtSetTop}.
\end{rmk}

\begin{dfn}[\textcolor{blue}{\index{Distance between sets}{Distance between sets}: $\dist(A,B)$}]
Let $X$ be a metric space and $A,B\subset X$. The distance between $A$ and $B$ is ~$\dist(A,B):=\inf\{d(a,b):a\in A,b\in B\}$.
\end{dfn}

The proofs of the following useful identities are brief and follow directly from definitions.
\begin{lmm}[\textcolor{OliveGreen}{DeMorgan's laws}]\label{DemogLaws}
Given any collection of subsets $\{S_\al\}_{\al\in A}$ of a set $S$,
\bit
\item[] $\left(\bigcup_{\al\in A} S_\al\right)^c=\bigcap_{\al\in A} S_\al^c$, ~~ $\left(\bigcap_{\al\in A} S_\al\right)^c=\bigcup_{\al\in A} S_\al^c$.
\eit
\end{lmm}

\begin{dfn}[\textcolor{blue}{Subset, \textcolor{magenta}{\index{Subspace topology}{Subspace topology}}, Subspace}]
In a space $(X,\T)$, every subset $A\subset X$ inherits the topology $A\cap\T :=\{A\cap O~|~O\in\T\}$ on $A$, called \ul{subspace topology} on $A$, which makes $A=(A,A\cap\T)$ a \ul{subspace} of $X$.

The subspace-complements $\{O^{c_A}=A-(A\cap O)=A\cap O^c:O\in\T\}$ are the \ul{closed sets} of the subspace topology, where $O^c=X-O$ is a closed set in $X$.
\end{dfn}

Note that the subspace topology is also called the \emph{\ul{relative topology}}, and accordingly, open (closed) sets in the subspace topology are said to be \emph{relatively open} (\emph{relatively closed}).

\begin{lmm}
A relatively open subset of an open set is itself open (as an intersection of two open sets). Similarly, a relatively closed subset of a closed set is itself closed (as an intersection of closed sets).
\end{lmm}

\begin{dfn}[\textcolor{blue}{Cover, \index{Open cover}{Open cover}, \index{Compact space}{Compact space}, \index{Connected space}{Connected space}}]
Let $(X,\T)$ be a space. A collection of sets $\A\subset\P(X)$ covers (or is a cover for) $X$ if $X=\bigcup_{A\in\A}A$. A cover $\U$ of $X$ is an open cover if $\U\subset\T$. $X$ is \ul{compact} if every open cover of $X$ has a finite subcover. $X$ is called \ul{connected} if $X$ cannot be written as a disjoint union of two nonempty open sets.
\end{dfn}
\begin{convention}
In a space $X$, we say a set $A\subset X$ is compact/connected if it is compact/connected as a subspace, i.e., with respect to its subspace topology. More generally, if $(X,\T)$ is a space with a property $P$ \ul{that does not determine the topology $\T$ in any way}, then we say a set $A\subset X$ has property $P$ if the subspace $(A,A\cap\T)$ has property $P$.
\end{convention}

\begin{dfn}[\textcolor{blue}{\index{Lindel\"of space}{Lindel\"of space}}]
A space s.t. every open cover has a countable subcover.
\end{dfn}

\begin{dfn}[\textcolor{blue}{\index{Neighborhood!}{Neighborhood}, Open neighborhood, Closed neighborhood}]
A neighborhood (nbd) of a set $A\subset(X,\T)$ is a set $N\subset X$ such that $A\subset O\subset N$ for some open set $O\in\T$. An open (resp. closed) neighborhood is a nbd that is an open (resp. closed) set.
\end{dfn}
Observe (from Lemma \ref{OpenCrit}) that in a topological space, a set is open if and only if it is a nbd of each of its points. In most applications we will encounter, neighborhoods that are not open will not be required because they serve the same purpose as open neighborhoods. Thus, we will henceforth assume every nbd is open, unless it is specified otherwise.

\begin{dfn}[\textcolor{blue}{\index{Base}{Base} of a topology}] In a space $(X,\T)$, A collection of open sets $\B\subset\T$ is a base for $\T$ if every element of $\T$ is a union of elements of $\B$, i.e., if
\bit
\item[] $\T=\left\{\bigcup \A :\A\subset\B\right\}$, ~ where ~ $\bigcup \A:=\bigcup\{A:A\in\A\}$.
\eit
\end{dfn}

\begin{lmm}[\textcolor{OliveGreen}{Base criterion}]
In a space $(X,\T)$, a collection of open sets $\B\subset\T$ is a base for $\T$ $\iff$ for any point $x_0\in X$ and any open neighborhood $U\in\T$ of $x_0$, we have $x_0\in B\subset U$ for some $B\in\B$.
\end{lmm}
\begin{proof}
If $\B$ is a base for $\T$ and $U$ is an open neighborhood of $x_0\in X$, then with $U=\bigcup B_\al$ for $B_\al\in\B$, we have $x_0\in B_\al\subset U$ for some $\al$. Conversely, if every open neighborhood $U$ of every point $x_0\in X$ satisfies $x_0\in B_{x_0,U}\subset U$ for some $B_{x_0,U}\in\B$, then every open set $O\in\T$ satisfies $O=\bigcup_{x\in O}B_{x,O}$, and so $\B$ is a base for $\T$.
\end{proof}

\begin{dfn}[\textcolor{blue}{\index{Local topology}{Local topology}, \index{Neighborhood! base}{Neighborhood base}}] Let $A\subset X=(X,\T)$. All open sets containing $A$ together with the empty set, $\T_A=\{\emptyset\}\cup\{O\in\T:A\subset O\}$, form a topology on $X$ called \ul{local topology} at $A$. Any base $\B_A$ for $\T_A$ is called a \ul{neighborhood base} (at $A\subset X$). Note that $\T=\bigcup_{A\subset X}\T_A$, and $\B=\bigcup_{A\subset X}\B_A$ is a base for $\T$.
\end{dfn}

\begin{dfn}[\textcolor{blue}{\index{First countable space}{First countable space}, \index{Second countable space}{Second countable space}}] Let $(X,\T)$ be a space. $(X,\T)$ is first countable if every point $x_0\in X$ has a countable neighborhood base.\\
$(X,\T)$ is second countable if $\T$ has a countable base.
\end{dfn}

\begin{dfn}[\textcolor{blue}{\index{Closure}{Closure} of a set}]
If $A\subset (X,\T)$, the closure $\ol{A}$ or $\Cl A$ of $A$ is the smallest closed set containing $A$, i.e., {\small $\ol{A}=\Cl A:=\bigcap\{C\supset A:C~\txt{closed}\}$}. (It is clear that $A$ is closed iff $A=\ol{A}$.)
\end{dfn}

\begin{lmm}[\textcolor{OliveGreen}{Closure criterion}]\label{ClosureCrit}
The closure of a set $A\subset(X,\T)$ is given by
\bea
\ol{A}=\{x\in X:N(x)\cap A\neq\emptyset~\txt{for every nbd $N(x)$ of $x$}\}.\nn
\eea
\end{lmm}
\begin{proof}
Let $B:=\{x\in X:N(x)\cap A\neq\emptyset~\txt{for every nbd $N(x)$ of $x$}\}$, which is closed because $x\in B^c$ iff some $N(x)\cap A=\emptyset$, which in turn implies that for every $y\in N(x)$, some $N(y)\cap A=\emptyset$, which implies $N(x)\subset B^c$. It is also clear that $A\subset B$, and so $\ol{A}\subset\ol{B}=B$.

If $x\in B$, i.e., every $N(x)\cap A\neq\emptyset$, suppose $x\not\in\ol{A}$, i.e., there is a closed set $C\supset A$ with $x\not\in C$. Then $x\in C^c$ = some $N(x)$ and so $C^c\cap A\neq\emptyset$ (a contradiction since $A\subset C$). It follows that $B\subset \ol{A}$.
\end{proof}

\begin{dfn}[\textcolor{blue}{\index{Limit point}{Limit points} of a set}]
Let $X$ be a space and $A\subset X$. A point $x\in X$ is a limit point of $A$, written $x\in A'$, if $\big(N(x)-\{x\}\big)\cap A\neq\emptyset$ for every nbd $N(x)$ of $x$. (It follows from Lemma \ref{ClosureCrit} that $\ol{A}=A\cup A'$).
\end{dfn}

\begin{dfn}[\textcolor{blue}{\index{Interior}{Interior} of a set}] If $A\subset (X,\T)$, the interior $A^o$ or $\Int A$ of $A$ is the largest open set contained in $A$, i.e., ~$A^o=\Int A:=\bigcup\big\{O\subset A:O~\txt{open}\big\}$.  (It is clear that $A$ is open iff $A=A^o$.)
\end{dfn}

{\flushleft\ul{Note}}: Using DeMorgan's laws, we easily see from the definitions that
\bea
\label{ComplEqs}\left(\ol{A}\right)^c=(A^c)^o~~~~\txt{and}~~~~(A^o)^c=\ol{A^c},
\eea
which simply mean that the complement of the closure (resp. the interior) is the interior (resp. a closure) of the complement.

\begin{lmm}\label{UnionClIntLmm}
Let $X$ be a space and $A,B\subset X$. Then (i) $\ol{A\cup B}=\ol{A}\cup\ol{B}$ and (ii) $(A\cap B)^o=A^o\cap B^o$.
\end{lmm}
\begin{proof}
(i) $\ol{A\cup B}\subset\ol{\ol{A}\cup\ol{B}}=\ol{A}\cup\ol{B}$. Also, $x\in \ol{A}\cup\ol{B}$ implies $x\in\ol{A}$ or $x\in\ol{B}$, which implies every $N(x)\cap A\neq\emptyset$ or every $N(x)\cap B\neq\emptyset$, which implies every $N(x)\cap(A\cup B)\neq\emptyset$ (i.e., $\ol{A}\cup\ol{B}\subset\ol{A\cup B}$). (ii) Thus, $(A\cap B)^o=\left[\ol{(A\cap B)^c}\right]^c=\left[\ol{A^c\cup B^c}\right]^c=\left[\ol{A^c}\cup\ol{B^c}\right]^c=A^o\cap B^o$.
\end{proof}

\begin{crl}\label{UnionClIntCrl}
If $B\subset A\subset X$, then the \emph{subspace-closure} $\Cl_A(B)$ and \emph{subspace-interior} $\Int_A(B)$ of $B$ in $A$ (i.e., closure and interior in the subspace topology) are given by
\bea
\label{SubSpCl}&&\textstyle\Cl_A(B)=\bigcap\{C_1:~B\subset C_1,~\txt{$C_1\subset A$ closed in $A$}\}=\bigcap\{C\cap A:~B\subset C,~\txt{$C\subset X$ closed}\}\nn\\
&&~~~~=\Cl(B)\cap A,\\
\label{SubSpInt}&&\Int_A(B)\sr{(\ref{ComplEqs})}{=} \big[\Cl_A(B^{c_A})\big]^{c_A}=\big[\Cl(B^c\cap A)\cap A\big]^c\cap A=\big[\Cl(B^c\cap A)\big]^c\cap A\nn\\
&&~~~~=\Int(B\cup A^c)\cap A.
\eea
\end{crl}

\begin{dfn}[\textcolor{blue}{\index{Boundary}{Boundary} of a set, \index{Clopen set}{Clopen (or closed open) set}, \index{Dense set}{Dense set}}] Let $X$ be a space. The boundary of a set $A\subset X$ is $\del A:=\ol{A}-A^o=\ol{A}\cap(A^o)^c\sr{(\ref{ComplEqs})}{=}\del A^c$. A set $A\subset X$ is clopen (equivalently, both closed and open) if $\del A=\emptyset$. A set $A\subset X$ is dense if $\ol{A}=X$.
\end{dfn}
Note that $\del A^o,\del\ol{A}\subset\del A$ (where equalities hold if $A\subset X$ is a ball in a metric space $X$), but the example $X=\Real$, $A=\Rational\subset\Real$ shows no equality holds in general. Besides $\del A=\del A^c$ (implying \ul{$A$ is clopen iff $A^c$ is clopen}), another relation showing symmetry between $A,A^c$ is $\del A^o\cap\del\ol{A}\sr{(\ref{ComplEqs})}{=}\ol{A^o}\cap\ol{(A^c)^o}$.

\begin{lmm}\label{SubSpBdLmm}
If $B\subset A\subset X$ and $B_1\subset X$, the \emph{subspace-boundary} $\del_A$ in $A$ (i.e., boundary in the subspace topology of $A$) satisfies (i) $\del_AB\subset (\del B)\cap A$ and (ii) $\del_A(B_1\cap A)\subset (\del B_1)\cap A$. (Equality holds in (i) if $A$ is closed in $X$ and $B$ is open in $X$.)
\end{lmm}
\begin{proof}
(i) By direct calculation, we have
\begin{align}
&\del_AB=\Cl_A(B)-\Int_A(B)=[\ol{B}\cap A]\cap[(B\cup A^c)^o\cap A]^c=\ol{B}\cap A\cap\ol{B^c\cap A}\nn\\
&~~~~\subset \ol{B}\cap A\cap\ol{B^c}\cap\ol{A}=\ol{B}\cap(B^o)^c\cap A=(\del B)\cap A,\nn
\end{align}
where equality holds if $A$ is closed in $X$ and $B$ is open in $X$. (ii) Moreover, if $B=B_1\cap A$ for some set $B_1\subset X$, then
\begin{align}
&\del_A(B_1\cap A)=\del_AB=\ol{B}\cap A\cap\ol{B^c\cap A}=\ol{B_1\cap A}\cap A\cap\ol{B_1^c\cap A}\nn\\
&~~~~\subset\ol{B_1}\cap \ol{A}\cap A\cap\ol{B_1^c}=\ol{B_1}\cap A\cap\ol{B_1^c}=\ol{B_1}\cap A\cap(B_1^o)^c=(\del B_1)\cap A.\nn\qedhere
\end{align}
\end{proof}

For any $A,B\subset X$, we have $\del(A\cup B)=\ol{A\cup B}-(A\cup B)^o\subset \ol{A}\cup\ol{B}-(A^o\cup B^o)=(\del A\cap\ol{B^c})\cup(\ol{A^c}\cap\del B)$, where equality holds if $A,B$ are open. Similarly, $\del(A\cap B)=\del(A^c\cup B^c)\subset(\del A\cap\ol{B})\cup(\ol{A}\cap\del B)$, where equality holds if $A,B$ are closed.

\begin{dfn}[\textcolor{blue}{\index{Separable space}{Separable space}}] A space $X$ is separable if it contains a countable dense set, i.e., there exists a countable set $A\subset X$ such that $\ol{A}=X$.
\end{dfn}

\begin{dfn}[\textcolor{blue}{\index{Generated topology}{Topology generated by sets}}] The topology generated by a collection of sets $\A\subset\P(X)$ is the smallest topology $\langle\A\rangle$ on $X$ containing $\A$. Note that $\langle\A\rangle$ is the intersection of all topologies on $X$ containing $\A$.
\end{dfn}

\begin{dfn}[\textcolor{blue}{\index{Subbase}{Subbase} of a topology}] In a space $(X,\T)$, a collection of open sets $\B\subset\T$ is said to be a subbase for $\T$ if $\B$ generates $\T$, i.e., if $\T=\langle\B\rangle$.
\end{dfn}
Note that a base is a subbase, but a subbase might not be a base. Also, a base covers $X$, but a subbase might not (except as an extra requirement in the definition of a subbase).
\begin{lmm}[\textcolor{OliveGreen}{Subbase-to-base}]\label{SubbaseCrit}
If $X$ is a set, then for any $\B\subset\P(X)$ we have
\bea
\langle\B\rangle=\{\txt{unions of finite intersections of elements of}~\B\cup\{\emptyset,X\}\}.\nn
\eea
Equivalently, ~$\D:=\{\txt{finite intersections of elements of}~\B\cup\{\emptyset,X\}\}$~ is a base for $\langle\B\rangle$.
\end{lmm}
\begin{proof}
Observe that unions $\T_1=\bigcup\F\I(\B\cup\{\emptyset,X\}):=\left\{\bigcup_\al\bigcap_{i=1}^{n_\al}B_{i,\al}\right\}$ of finite intersections $\F\I(\B\cup\{\emptyset,X\}):=\{\bigcap_{i=1}^nB_i:B_i\in\B\cup\{\emptyset,X\}\}$ of elements of $\B\cup\{\emptyset,X\}$ form a topology. Indeed, the collection $\T_1$ contains $\emptyset,X$, is clearly closed under unions, and is closed under finite intersections because finite intersections of any unions can be written as unions of finite intersections as follows:
\bit
\item[] $\Big(\bigcup_\al~\bigcap_{i=1}^{n_\al}B_{i,\al}\Big)\cap\left(\bigcup_\beta~\bigcap_{i=1}^{n_\beta}B_{i,\beta}\right)=\bigcup_{\al,\beta}~\bigcap_{i=1,j=1}^{n_\al,n_\beta}B_{i,\al}\cap B_{j,\beta}$.
\eit
Moreover, it is clear that any topology containing $\B$ contains the finite intersections $\F\I(\B\cup\{\emptyset,X\})$, and so contains $\T_1$. Hence, $\T_1=\langle\B\rangle$.
\end{proof}
\begin{crl}
For a set $X$, a collection $\B\subset\P(X)$ is a base for some topology on $X$ iff
\bit[leftmargin=0.9cm]
\item[(i)] $\B$ contains $\emptyset$,
\item[(ii)] $\B$ covers $X$ (i.e., $X=\bigcup\B$), and
\item[(iii)] Finite intersections of elements of $\B$ are unions of elements of $\B$.
\eit
\end{crl}

\subsection{Continuous maps, Open maps, Quotient maps}
\begin{dfn}[\textcolor{blue}{\index{Continuous map}{Continuous map}, Continuity at a point}] A map of spaces $f:X\ra Y$ is continuous, also written $f\in \C(X,Y):=\{\txt{continuous}~f:X\ra Y\}$, if for any open set $V\subset Y$, the preimage $f^{-1}(V)\subset X$ is open. That is, $f:(X,\T_X)\ra (Y,\T_Y)$ is continuous if $f^{-1}(\T_Y)\subset\T_X$.

$f$ is continuous at $x\in X$ if for every open neighborhood $V$ of $f(x)$, the preimage $f^{-1}(V)$ is an open neighborhood of $x$ (Equivalently, for any open set $V\ni f(x)$, there exists an open set $U\ni x$ such that $f(U)\subset V$).
\end{dfn}

Note that a continues map takes a compact/connected set to a compact/connected set.
\begin{lmm}
A map of spaces $f:X\ra Y$ is continuous $\iff$ every set $A\subset X$ satisfies {\footnotesize $\ol{A}\subset f^{-1}\left(\ol{f(A)}\right)$}, i.e.,{\footnotesize $f(\ol{A})\subset\ol{f(A)}$}. Equivalently, by setting $A:=f^{-1}(B^c)$ in {\footnotesize $\ol{A}\subset f^{-1}\left(\ol{f(A)}\right)$} and applying (\ref{ComplEqs}), we see that $f:X\ra Y$ is continuous $\iff$ every set $B\subset Y$ satisfies $f^{-1}(B^o)\subset [f^{-1}(B)]^o$.
\end{lmm}
\begin{proof} If $f$ is continuous then {\footnotesize $f^{-1}\left(\ol{f(A)}\right)$} is closed, and so {\footnotesize $A\subset f^{-1}(f(A))\subset f^{-1}\left(\ol{f(A)}\right)$} implies {\footnotesize $\ol{A}\subset f^{-1}\left(\ol{f(A)}\right)$}. Conversely, if {\footnotesize $\ol{A}\subset f^{-1}\left(\ol{f(A)}\right)$} for any $A\subset X$, then for open $V\subset Y$, $A:=f^{-1}(V^c)$ satisfies {\footnotesize $\ol{f^{-1}(V^c)}\subset f^{-1}\left(\ol{f(f^{-1}(V^c))}\right)\subset f^{-1}\left(\ol{V^c}\right)=f^{-1}(V^c)\subset \ol{f^{-1}(V^c)}$}.
\end{proof}

\begin{lmm}[\textcolor{OliveGreen}{Continuity criterion}]\label{ContnyCrit1}
A map of spaces $f:X\ra Y$ is continuous $\iff$ continuous at each point $x\in X$.
\end{lmm}
\begin{proof}
If $f$ is continuous then it is clear that for every open set $V\ni f(x)$, $f^{-1}(V)$ is an open set containing $x$. Conversely, if $f$ is continuous at every point $x\in X$, then for any open set $V\subset Y$ (with $V_y\subset V$ denoting an open set containing $y$) we have the open set
\bit
\item[] $f^{-1}(V)=f^{-1}\left(\bigcup_{y\in V}V_y\right)=\bigcup_{y\in V}f^{-1}(V_y)$.\qedhere
\eit
\end{proof}

\begin{crl}[\textcolor{OliveGreen}{Metric space continuity}]\label{ContnyCrit2}
A map of metric spaces $f:X\ra Y$ is continuous $\iff$ for any $\vep>0$ and $x\in X$, there exists $\delta_x(\vep)>0$ such that $f\big(B_{\delta_x(\vep)}(x)\big)\subset B_\vep\big(f(x)\big)$, i.e., for every $x'\in X$, $d(x,x')<\delta_x(\vep)$ implies $d(f(x),f(x'))<\vep$.
\end{crl}

\begin{crl}\label{ContnyCrit3}
A map of metric spaces $f:X\ra Y$ is continuous $\iff$ for each $x\in X$ there exists a function $\omega_x:[0,\infty)\ra[0,\infty)$ such that (i) $\omega_x(t)$ is nondecreasing near $t=0$, (ii) $\lim_{t\ra 0}\omega_x(t)=0=\omega_x(0)$, and (iii) $d(f(x),f(x'))\leq\omega_x\big(d(x,x')\big)$ for all $x'\in X$.
\end{crl}
\begin{proof}
Assume $f$ is continuous, i.e., for any $\vep>0$ and $x\in X$, there exists $\delta_x(\vep)>0$ (where $\delta_x(\vep)\leq\vep$ wlog) such that $d(x,x')<\delta_x(\vep)$ implies $d(f(x),f(x'))<\vep$, for all $x'\in X$. Define the function
\bea
\omega_x(t):=\sup\left\{d(f(x),f(x')):d(x,x')\leq t,x'\in X\right\},\nn
\eea
which is nondecreasing, $\omega_x(0)=0$, and $d(f(x),f(x'))\leq\omega_x(d(x,x'))$ for all $x'\in X$. Then
\bea
d(x,x')<\delta_x(\vep)~~\Ra~~d(f(x),f(x'))\leq\omega_x(d(x,x'))\leq \omega_x(\delta_x(\vep))\leq\vep.\nn
\eea
Thus by taking $\vep\ra 0$, and noting $\delta_x(\vep)\leq\vep$, we see that $\omega_x(t)\ra 0$ as $t\ra 0$.

Conversely, assume there exists a function $\omega_x:[0,\infty)\ra[0,\infty)$ such that (i) $\omega_x(t)$ is nondecreasing near $t=0$, (ii) $\lim_{t\ra 0}\omega_x(t)=0=\omega_x(0)$, and (iii) $d(f(x),f(x'))\leq\omega_x\big(d(x,x')\big)$ for all $x'\in X$. Then for any $\vep>0$, because $\omega_x(t)\ra 0$ as $t\ra 0$, we can choose a $\delta=\delta_x(\vep)>0$ such that
\[
d(x,x')<\delta~~\Ra~~d(f(x),f(x'))\leq\omega_x\big(d(x,x')\big)\leq \omega_x\big(\delta\big)<\vep. \qedhere
\]
\end{proof}

\begin{dfn}[\textcolor{blue}{
\index{Uniform continuity}{Uniform continuity},
\index{Modulated uniform continuity}{Modulated uniform continuity},
\index{Lipschitz continuity}{Lipschitz continuity},
\index{H\"older continuity}{H\"older continuity}}]
Let $X,Y$ be metric spaces and $f:X\ra Y$ a continuous map. If the function $\delta_x$ in Corollary \ref{ContnyCrit2} is independent of $x$, we say $f$ is uniformly continuous (uc). If (in some neighborhood $[0,\al)$ of $0$ in $[0,\infty)$) the function $\omega_x$ in Corollary \ref{ContnyCrit3} is independent of $x$, we say $f$ is \ul{$\omega$-uniformly continuous} (or \ul{$\omega$-modulated uniformly continuous}). If $\omega(t)=ct$ for some constant $c$, we say $f$ is $c$-Lipschitz continuous. If $\omega(t)=ct^\al$ for some constants $c$ and $0<\al<1$, we say $f$ is $(c,\al)$-H\"older continuous.
\end{dfn}

\begin{dfn}[\textcolor{blue}{\index{Open map}{Open map}, \index{Closed map}{Closed map}, Openness at a point}] A map of spaces $f:X\ra Y$ is \ul{open} if for each open set $O\subset X$, the image $f(O)\subset Y$ is an open set. Similarly, $f:X\ra Y$ is \ul{closed} if for each closed set $C\subset X$, the image $f(C)\subset Y$ is closed in $Y$.

$f$ is open at $x\in X$ if for any open neighborhood $U$ of $x$, $f(U)$ is an open neighborhood of $f(x)$ (Equivalently, for any open set $U\ni x$, there exists an open set $V\ni f(x)$ such that $V\subset f(U)$).
\end{dfn}

\begin{lmm}[\textcolor{OliveGreen}{Open map criterion}]
A map of spaces $f:X\ra Y$ is open $\iff$ open at each point $x\in X$.
\end{lmm}
\begin{proof}
If $f$ is open, it is clear that $f$ is open at each $x\in X$. Conversely, if $f$ is open at every point $x\in X$, then for any open set $U\subset X$ (with $U_x\subset U$ denoting an open set containing $x$) we have the open set ~$f(U)=f\left(\bigcup_{x\in U}U_x\right)=\bigcup_{x\in U}f(U_x)$.
\end{proof}

\begin{crl}
A map of metric spaces $f:X\ra Y$ is open $\iff$ for any $\vep>0$ and $x\in X$, there exists $\delta_x(\vep)>0$ such that $B_{\delta_x(\vep)}\big(f(x)\big)\subset f\big(B_\vep(x)\big)$, i.e., for every $x'\in X$ there exists $x''\in f^{-1}(f(x'))$ such that $d\big(f(x),f(x')\big)<\delta_x(\vep)$ implies $d(x,x'')<\vep$.
\end{crl}

\begin{crl}\label{OpnssCrit3}
A map of metric spaces $f:X\ra Y$ is open $\iff$ for each $x\in X$ there exists a function $\omega_x:[0,\infty)\ra[0,\infty)$ such that (i) $\omega_x(t)$ is nondecreasing near $t=0$, (ii) $\lim_{\vep\ra 0}\omega_x(\vep)=0=\omega_x(0)$, and (iii) for any $x'\in X$ there exists $x''\in f^{-1}(f(x'))$ such that $d(x,x'')\leq \omega_x\Big(d\big(f(x),f(x')\big)\Big)$. In particular, we have
\bea
\dist\Big(x,f^{-1}\big(f(x')\big)\Big)\leq \omega_x\Big(d\big(f(x),f(x')\big)\Big),~~~~\txt{for all}~~~~x'\in X.\nn
\eea
(The proof is the same as that of Corollary \ref{ContnyCrit3}.)
\end{crl}

\begin{crl}\label{OpnssCrit4}
A \ul{continuous} map of metric spaces $f:X\ra Y$ is open $\iff$ for each $x\in X$ there exists a function $\omega_x:[0,\infty)\ra[0,\infty)$ such that (i) $\omega_x(t)$ is nondecreasing near $t=0$, (ii) $\lim_{\vep\ra 0}\omega_x(\vep)=0=\omega_x(0)$, and (iii) {\small$\dist\Big(x,f^{-1}\big(f(x')\big)\Big)\leq \omega_x\Big(d\big(f(x),f(x')\big)\Big)$}~ for all $x'\in X$.
(This is immediate from Corollary \ref{OpnssCrit3}, because $f^{-1}(f(x'))$ is a closed set.)
\end{crl}

\begin{lmm}[\textcolor{OliveGreen}{Continuity on compact sets}]
A continuous map of metric spaces $f:X\ra Y$ is uniformly continuous on compact subsets of $X$.
\end{lmm}
\begin{proof}
It suffices to assume $X$ is compact (since restrictions of $f$ are continuous). Fix $\vep>0$. Then for each $x\in X$, there is $\delta_x$ such that $d(x,x')<\delta_x$ implies $d(f(x),f(x'))<\vep/2$. Since $X$ is compact, let $\{B_{\delta_{x_i}/2}(x_i)\}_{i=1}^n$ cover $X$. Define $\delta:=\min_i\delta_{x_i}/2$. Pick any $u,v\in X$. Let $u\in B_{\delta_{x_j}/2}(x_j)$ for some $j$. Then $d(u,v)<\delta$ implies $u,v\in B_{\delta_{x_j}}(x_j)$, which in turn implies
\bea
d(f(u),f(v))\leq d(f(u),f(x_j))+d(f(x_j),f(v))<\vep/2+\vep/2=\vep.\nn\qedhere
\eea
\end{proof}
\begin{crl}
If a metric space $X$ is \ul{locally compact} (in the sense that every point of $X$ has a neighborhood whose closure is compact), then every continuous map of metric spaces $f:X\ra Y$ is \ul{locally uniformly continuous} (in the sense that every point of $X$ has a neighborhood in which $f$ is uniformly continuous).
\end{crl}

\begin{dfn}[\textcolor{blue}{\index{Homeomorphism}{Homeomorphism}, \index{Imbedding}{Imbedding}}]\label{HomeoImbed}
A homeomorphism is a bijective map of spaces $h:X\ra Y$ such that both $h,h^{-1}$ are continuous (in which case, we say $X,Y$ are homeomorphic, or $X\cong Y$). A map of spaces $f:X\ra Y$ is an imbedding (of $X$ into $Y$) if it is a homeomorphism onto its image (i.e., $f:X\ra f(X)$ is a homeomorphism).
\end{dfn}

\begin{dfn}[\textcolor{blue}{\index{BiLipschitz map}{BiLipschitz map}, \index{Isometric map}{Isometric map}, \index{BiLipschitz equivalence}{BiLipschitz equivalence}, \index{Isometry}{Isometry}}]
A map $f:(X,d_X)\ra(Y,d_Y)$ is $c$-biLipschitz (or a \ul{$c$-biLipschitz imbedding}) if $d_X(x,x')/c$ $\leq$ $d_Y\big(f(x),f(x')\big)$ $\leq$ $cd_X(x,x')$ for all $x,x'\in X$. A $1$-biLipschitz map is called an \ul{isometric map}. A surjective biLipschitz map is called a \ul{biLipschitz equivalence}. An \ul{isometry} is a surjective isometric map (i.e., an isometric homeomorphism).
\end{dfn}

\begin{dfn}[\textcolor{blue}{\index{Quotient! map}{Quotient map}, \textcolor{magenta}{\index{Quotient! topology}{Quotient topology}}, \index{Quotient! space}{Quotient space}}]
A map of spaces $f:X\ra Y$ is called a \ul{quotient map} (making $Y$ a \ul{quotient space} of $X$, and the topology of $Y$ a \ul{quotient topology}) if it is surjective and $B\subseteq Y$ is open $\iff$ $f^{-1}(B)$ is open in $X$. (Note that a quotient map is continuous, but does not have to be an open map).
\end{dfn}
\begin{examples*}
(1) Given a space $(X,\T)$ and any equivalence relation $\sim$ on $X$, the map $q:X\ra X_\sim$, $x\mapsto x_\sim$ onto the set of equivalence classes {\small$X_\sim={X\over\sim}=\left\{x_\sim:=\{x'\in X:x\sim x'\}~|~x\in X\right\}$} is a quotient map if $X_\sim$ is given the topology {\small$\T_\sim:=\big\{B\subset X_\sim:q^{-1}(B)\in\T\big\}$}.\\
(2) Given a space $X$, any map of sets $f:X\ra A$ induces an equivalence relation $\sim_f$ on $X$, where $x\sim_f x'$ $\iff$ $f(x)=f(x')$, with equivalence classes
\bit
\item[] {\small$X_{\sim_f}={X\over\sim_f}=\left\{x_{\sim_f}:=f^{-1}(f(x))~|~x\in X\right\}=\{f^{-1}(a):a\in f(X)\}$}.
\eit
(3) Given a space $X$ and any subspace $A\subset X$, we can define an equivalence relation $\sim$ on $X$ as follows: $x\sim x'$ if $x=x'$ or $x,x'\in A$. The space $X_\sim$ is denoted by ${X\over A}$ or $X/A$.\\
(4) Every quotient space of \ul{compact/connected} space is compact/connected (by continuity of the quotient map).
\end{examples*}

\begin{thm}[\textcolor{OliveGreen}{Universal property of quotient maps}]\label{UPQM}
Let $X$ be a space and $\sim$ an equivalence relation on $X$. For any continuous map $f:X\ra Y$ such that $x\sim x'$ $\Ra$ $f(x)=f(x')$, there exists a unique continuous map $\ol{f}:X_\sim\ra Y$ such that the following diagram commutes (i.e., $f$ factors through $q$ as $f=\ol{f}\circ q$).
\bc\bt
X\ar[d,"q"']\ar[rr,"f"]&& Y\\
X_\sim\ar[urr,dashed,"\ol{f}"']
\et\ec
\end{thm}
\begin{proof}
Let $X_\sim=\{x_\sim:x\in X\}$, where $x_\sim=\{x'\in X:x\sim x'\}$. \ul{\emph{Existence}}: Note that the hypotheses imply $\ol{f}:X_\sim\ra Y$, $x_\sim\mapsto f(x)$ is well defined since
\bea
x_\sim=x'_\sim~~\iff~~x\sim x'~~\Ra~~\ol{f}(x_\sim)=f(x)=f(x')=\ol{f}(x'_\sim),\nn
\eea
and moreover, $\ol{f}\circ q(x)=\ol{f}(x_\sim)=f(x)$, and so $\ol{f}\circ q=f$. \ul{\emph{Uniqueness}}: If $g:{X\over\sim}\ra Y$ is any continuous map such that $g\circ q=f$, then because $q$ is surjective, $f=g\circ q=\ol{f}\circ q$ ~$\Ra$~ $g=\ol{f}$. \ul{\emph{Continuity}}: If $B\subset Y$ is open, then because $f^{-1}(B)=q^{-1}\left(\ol{f}^{-1}(B)\right)$, $f$ is continuous, and $q$ is a quotient map, we see that $\ol{f}^{-1}(B)$ is open. Hence, $\ol{f}$ is continuous.
\end{proof}

\begin{crl}
Every continuous map of spaces $f:X\ra Y$ induces a unique injective continuous map $\ol{f}:X_{\sim_f}\ra Y$ such that $f=\ol{f}\circ q$, where $q:X\ra X_{\sim_f}$ is the quotient map.
\end{crl}

\begin{crl}[\textcolor{OliveGreen}{Continuous maps on a quotient space}]
If $Y=X_\sim$ is a quotient space of a space $X$, then continuous maps $\C(Y,Z)$ are precisely continuous maps $\C(X,Z)$ that are constant on equivalence classes $x_\sim:=\{x'\in X~|~x'\sim x\}$ for all $x\in X$.
\end{crl}
\begin{proof}
Let $q:X\ra Y$ be the quotient map. By the universal property of quotient maps, any continuous map $f:X\ra Z$ that is constant on equivalence classes (i.e., $x\sim x'$ $\Ra$ $f(x)=f(x')$~) gives a unique continuous map $\wt{f}:Y\ra Z$. Conversely, any continuous map $h:Y\ra Z$ also gives the unique continuous map $H=h\circ q:X\sr{q}{\ral}Y\sr{h}{\ral}Z$, which is by construction constant on equivalence classes (i.e., $x\sim x'$ $\Ra$ $H(x)=H(x')$~).
\end{proof}

\subsection{Sequential topology, Compactness in metric spaces}

\begin{dfn}[\textcolor{blue}{\index{Sequence}{Sequence}, \index{Convergent sequence}{Convergent sequence}, \index{Sequentially closed}{Sequentially closed}, \index{Sequentially continuous}{Sequentially continuous}, \index{Sequentially compact}{Sequentially compact}}]
Let $X$ be a space. A \ul{sequence} in $X$ is a map $s:\Natural\ra X$, $n\mapsto s_n$, which is often written as $\{s_n\}\subset X$ or simply as $s_n\in X$. A sequence $\{x_n\}\subset X$ \ul{converges} to $x$ (written $x_n\ra x$) if every neighborhood of $x$ contains all but finitely many $x_n$ (i.e., for any open set $U\ni x$, we have $\{x_N,x_{N+1},\cdots\}\subset U$ for some $N$). A set $A\subset X$ is \ul{sequentially closed} if for every sequence $x_n\in A$, $x_n\ra x\in X$ implies $x\in A$. A map of spaces $f:X\ra Y$ is \ul{sequentially continuous} if for every sequence $x_n\ra x$ in $X$, we have $f(x_n)\ra f(x)$ in $Y$. A set $K\subset X$ is \ul{sequentially compact} if every sequence $x_n\in K$ has a convergent subsequence $x_{n_k}\ra x\in K$.
\end{dfn}
Note that in a metric space $X=(X,d)$, a sequence $x_n\ra x$ if and only if $d(x,x_n)\ra 0$.

\begin{facts*}[\textcolor{OliveGreen}{Metric space peculiarities}]
(i) In a metric space $X$, a set $A\subset X$ is closed iff sequentially closed. (ii) A map of metric spaces $f:X\ra Y$ is continuous iff sequentially continuous. (iii) In a metric space $X$, a set $K\subset X$ is compact iff sequentially compact.
\end{facts*}
\begin{proof}
{\flushleft (i):} If $A$ is closed then by Lemma \ref{ClosureCrit}, every $x_n\in A$ satisfies ``$x_n\ra x\in X$ implies $x\in A$''. If $A$ is not closed, then with $x\in \ol{A}\backslash A$, we can pick a sequence $x_n\in B_{1/n}(x)\cap A$, which satisfies $1/n\geq d(x,x_n)\ra 0$ and so $x_n\ra x\not\in A$.
{\flushleft (ii):} If $f$ is continuous, it follows from Corollary \ref{ContnyCrit3} that ``$x_n\ra x$ implies $f(x_n)\ra f(x)$''. If $f$ is not continuous, then for some closed set $C\subset Y$, $f^{-1}(C)$ is not closed, i.e., by (i) there exists $x_n\in f^{-1}(C)$ such that $x_n\ra x\not\in f^{-1}(C)$, which implies $f(x_n)\not\!\!\ral f(x)$ since $C$ is closed and $f(x_n)\in C$ but $f(x)\not\in C$.
{\flushleft (iii):} This is Corollary \ref{MetricCompact2} below.
\end{proof}

\begin{dfn}[\textcolor{blue}{\index{Bounded set}{Bounded set}, \index{Totally bounded set}{Totally bounded set}, \index{Cauchy sequence}{Cauchy sequence}, \index{Complete set}{Complete set}}] Let $X$ be a metric space and $A\subset X$. The set $A$ is \ul{bounded} if $A\subset B_R(x)$ for some $R>0$, $x\in X$. The set $A$ is \ul{totally bounded} if for any $\vep>0$ there exist points $x_1,...,x_n\in X$ such that $A\subset\bigcup_{i=1}^nB_\vep(x_i)$ (i.e., $A$ can be covered with a finite number of balls of any radius). A sequence $x_n\in X$ is \ul{Cauchy} if $d(x_n,x_m)\ra 0$ as $\{m,n\}\ra\{\infty\}$. The set $A$ is \ul{complete} if every Cauchy sequence in $A$ converges (to a point) in $A$.
\end{dfn}

\begin{thm}[\textcolor{OliveGreen}{Compactness in a metric space}]\label{MetricCompact1}
In a metric space $X$, a set $A\subset X$ is compact $\iff$ complete and totally bounded.
\end{thm}
\begin{proof}
{\flushleft $\Ra$:} Assume $A$ is compact. It is clear that $A$ is totally bounded since every open cover of $A$ (using balls in particular) has a finite subcover. Also, $A$ is sequentially compact [[i.e., every sequence $\{a_n\}\subset A$ has a convergent subsequence: Indeed if $a_i\in A$ is an infinite sequence with no convergent subsequence, then $K=\{a_i\}$ is a closed set with no limit point, and so $A\subset(A-K)\cup\bigcup_i B_{\vep_i}(a_i)$ with $B_{\vep_i}(a_i)\cap K=a_i$. By compactness of $A$, $K$ is finite (a contradiction)]]. Hence $A$ is complete, since a Cauchy sequence converges iff it has a convergent subsequence.
{\flushleft $\La$:} Assume $A$ is complete and totally bounded. Then $A$ is sequentially compact (i.e., every sequence in $A$ has a subsequence that converges in $X$, hence converges in $A$ by completeness). Otherwise, if an infinite sequence $a_i\in A$ has no subsequence that converges in $X$, then $d(a_i,a_j)\not\ra 0$ since $A$ is complete. Thus, there is a subsequence $a_{i_k}$ and $\vep_0>0$ such that $d(a_{i_k},a_{i_{k'}})\geq 2\vep_0$ for all $k,k'$, and so we need an infinite number of balls $B_{\vep_0}(x_r)$ to cover $A$, i.e., $A\not\subset\bigcup_{r=1}^nB_{\vep_0}(x_r)$ for all $n\geq 1$ (a contradiction).

Let $A\subset\bigcup_{\al\in I}U_\al$, $U_\al\subset A$ be an infinite open cover of $A$. Let $A\subset \bigcup_{j=1}^mB_\vep(x_j^\vep)$, which is possible by total boundedness. Then by choosing $\vep$ small enough, we get $B_\vep(x_j^\vep)\cap A\subset U_{\al_j}$, and so $A\subset \bigcup_{j=1}^mU_{\al_j}$. Otherwise, if for every $\vep>0$, some $B_{\vep}(x_{j_\vep}^\vep)\cap A\not\subset U_\al$ for all $\al$, then for every $n$, some $a_n\in \left(B_{1/n}\left(x_n\right)\cap A\right)\backslash U_\al$ for all $\al$, where $x_n:=x_{j_{1/n}}^{1/n}$. By sequential compactness of $A$, a subsequence $a_{n_k}\ra a\in A$, and so
\bit
\item[] $d(x_{n_k},a)\leq d(x_{n_k},a_{n_k})+d(a_{n_k},a)\leq 1/n_k+d(a_{n_k},a)\ra 0~~\Ra~~x_{n_k}\ra a\in A.$
\eit
Since $\{U_\al\}$ covers $A$, some $U_{\al_0}\ni a$, i.e., some $B_\delta(a)\cap A\subset U_{\al_0}$. So, for sufficiently large $k$,
\bit
\item[] $a_{n_k}\in B_{1/n_k}\left(x_{n_k}\right)\cap A\subset B_\delta(a)\cap A\subset U_{\al_0}$, ~ (a contradiction).\qedhere
\eit
\end{proof}

\begin{crl}[\textcolor{OliveGreen}{Compactness in a metric space}]\label{MetricCompact2}
In a metric space $X$, a set $A$ is compact $\iff$ sequentially compact.
\end{crl}
\begin{proof}
If $A$ is compact, it is clear from the proof of Theorem \ref{MetricCompact1} that $A$ is sequentially compact. Conversely, assume $A$ is sequentially compact. Then $A$ is complete since every cauchy sequence in $A$ has a convergent subsequence, and so converges in $A$. To show $A$ is totally bounded, fix $\vep>0$. Then for any $a\in A$, the sequence $a_1:=a\in A$, $a_2\in A\backslash N_\vep(a_1)$, $a_3\in A\backslash\big(N_\vep(a_1)\cup N_\vep(a_2)\big)$, $\cdots$, $a_i\in A\backslash\big(N_\vep(a_1)\cup\cdots\cup N_\vep(a_{i-1})\big)$, $\cdots$
satisfies $d(a_i,a_j)\geq \vep$ for all $i,j$, and so terminates (otherwise we get an infinite sequence with no convergent subsequence). That is, $A\subset\bigcup_{i=1}^nN_\vep(a_i)$ for some $n$.
\end{proof}

\subsection{Product topology, Compactness of product spaces}

\begin{dfn}[\textcolor{blue}{\index{Box topology}{Box topology}}] Let $(X_\al,\T_\al)_{\al\in A}$ be a family of spaces. The box topology on {\footnotesize$X=\prod_{\al\in A}X_\al=\{(x_\al)_{\al\in A}:x_\al\in X_\al\}=\left\{A\sr{x}{\ral}\bigcup X_\al~|~x(\al)\in X_\al\right\}$} is the topology $\T_{\txt{box}}$ with base ~$\B_{\txt{box}}=\left\{\prod_{\al\in A}U_\al:~U_\al\in\T_\al\right\}$.
\end{dfn}

\begin{dfn}[\textcolor{blue}{Product, \textcolor{magenta}{Product topology}, \index{Product space}{Product space}}]\label{ProdTop}
Let $(X_\al,\T_\al)_{\al\in A}$ be a family of spaces. The product topology on { $X:=\prod_{\al\in A}X_\al=\{(x_\al)_{\al\in A}:x_\al\in X_\al\}$} = {\small $\left\{A\sr{x}{\ral}\bigcup X_\al~|~x(\al)\in X_\al\right\}$} is the topology $\T_{\txt{prod}}$ with base given by open rectangles
{\small\begin{align}
&\B_{\txt{prod}}=\left\{{\textstyle\prod_{\al\in A}}U_\al~\big|~\substack{U_\al\in\T_\al},~\substack{U_\al=X_\al~\txt{except for}\\\txt{finitely many}~\al~~}\right\}=\left\{U_F:={\textstyle\prod_{\al\in F}}U_\al\times{\textstyle\prod_{\al\in A\backslash F}}X_\al~\big|~\substack{\txt{finite}~F\subset A,\\ U_\al\in\T_\al~~~~}\right\},\nn
\end{align}}
where ~{\footnotesize $\prod_{\al\in A}U_\al:=\big\{(x_\al)_{\al\in A}:x_\al\in U_\al\big\}$}~ and ~{\footnotesize $U_F=\left\{x\in X:x_\al\in U_\al\in\T_\al,~\txt{for}~\al\in F\right\}$}.
\end{dfn}

\begin{notes*}
\begin{enumerate}[leftmargin=0.7cm]
\item A neighborhood base $\B_x$ for $x\in X$ consists of open rectangles of the form
{\small\begin{align}
\textstyle N_F(x)=\{y\in X:y_\al\in N(x_\al)\in\T_\al,~\txt{for}~\al\in F\}=\prod_{\al\in F}N(x_\al)\times\prod_{\beta\not\in F}X_\beta=\bigcap_{\al\in F}N_\al(x),\nn
\end{align}}where {\small $N_\al(x):=\left\{y\in X:y_\al\in N(x_\al)\in\T_\al\right\}=N(x_\al)\times\prod_{\beta\neq {\al}}X_\beta$} is a ``strip through $N(x_{\al})$''.
\item \ul{\emph{Notation (Product space)}}: Whenever $\prod_{\al\in A}X_\al$ is a space, we will always assume it is given the product topology, unless it is stated otherwise.
\item \ul{\emph{Cartesian power of a space}}: If $X_\al=X$ for all $\al\in A$ (for some space $X$), we write $X^A$ for $\prod_{\al\in A}X_\al=\prod_{\al\in A}X:=\big\{\txt{maps}~x:A\ra X\big\}=\big\{(x_\al)_{\al\in A}:x_\al\in X\big\}$. As usual, we write $X^n$ for $X^{\{0,1,2,\cdots,n-1\}}$. We will also write $X^\infty$ for $X^\Natural$.
\end{enumerate}
\end{notes*}

\begin{lmm}[\textcolor{OliveGreen}{\cite[Theorem 19.6, p 117]{munkres}}]\label{ComFunCont}
Let $\{X_\al\}_{\al\in A}$ be a family of spaces, and consider any map $f:Z\ra\prod_{\al\in A}X_\al$, $z\mapsto f(z)$, where it is clear that we necessarily have $f(z)=\left(f_\al(z)\right)_{\al\in A}$, for some maps $f_\al:Z\ra X_\al$. If $\prod_{\al\in A}X_\al$ is given the product topology, then $f$ is continuous $\iff$ each $f_\al$ is continuous.
\end{lmm}
\begin{proof}
If $f$ is continuous, then each $f_\al=\pi_\al\circ f$ is continuous, because the maps $\pi_\al:X\ra X_\al$, $x\mapsto x_\al$ are continuous due to $\pi_\al^{-1}(U_\al)=U_\al\times\prod_{\beta\neq\al}X_\beta$. Conversely, if each $f_\al$ is continuous, then so is $f$, since for each finite set $F\subset A$,
\begin{align}
&f^{-1}\left({\textstyle\prod_{\al\in F}}U_\al\times{\textstyle\prod_{\al\in A\backslash F}}X_\al\right)=\Big\{z\in Z:f_\al(z)\in U_\al~\txt{for}~\al\in F,~f_\al(z)\in X_\al~\txt{for}~\al\not\in F\Big\}\nn\\
&~~~~={\textstyle\bigcap_{\al\in F}}f_\al^{-1}(U_\al)~\cap~{\textstyle\bigcap_{\al\not\in F}}f_\al^{-1}(X_\al)={\textstyle\bigcap_{\al\in F}}f_\al^{-1}(U_\al)\nn
\end{align}
is open as a finite intersection of open sets.
\end{proof}
The proof above shows the converse of this result fails for the box topology on $\prod_{\al\in A} X_\al$ since an infinite intersection of open sets is not always open.

\begin{lmm}[\textcolor{OliveGreen}{Tube lemma}]
Let $X,K$ be spaces with $K$ compact. For any point $x\in X$ and any open neighborhood $N$ ( $\{x\}\times K\subset N\subset X\times K$) of $\{x\}\times K$, there is an open neighborhood $U$ ($x\in U\subset X$) of $x$ such that $U\times K \subset N$.
\end{lmm}
\begin{proof}
Since the rectangles $\{U\times V:~U~\txt{open in}~X,~V~\txt{open in}~K\}$ form a base for the topology of $X\times K$, let $N=\bigcup_{\al\in \A}U_\al\times V_\al$, where $U_\al\subset X$, $V_\al\subset K$ are open. Then
{\small\begin{align}
&\textstyle \{x\}\times K=\bigcup\limits_{\al\in \A}(U_\al\times V_\al)\cap(\{x\}\times K)=\bigcup\limits_{\al\in \A}(U_\al\cap\{x\})\times(V_\al\cap K)=\bigcup\limits_{\al\in \A}(U_\al\cap\{x\})\times V_\al\nn\\
&\textstyle ~~~~=\bigcup_{\al:x\in U_\al}\{x\}\times V_\al=\{x\}\times\bigcup_{\al:x\in U_\al}V_\al,~~~~\txt{(union over those $U_\al$ that contain $x$)}\nn\\
&\textstyle ~~\Ra~~K=\bigcup_{\al\in\A}V_\al=\bigcup_{\al:x\in U_\al}V_\al~\sr{\txt{K is compact}}{=}~\bigcup_{1\leq i\leq n:x\in U_{\al_i}}V_{\al_i},~~~~n<\infty.\nn
\end{align}}Hence, the finite intersection $U=\bigcap_{i:x\in U_{\al_i}}U_{\al_i}\subset X$ is an open set containing $x$, and
{\small\bea
\textstyle U\times K=U\times \bigcup\limits_{i:x\in U_{\al_i}}V_{\al_i}=\bigcup\limits_{i:x\in U_{\al_i}}(U\times V_{\al_i})\subset \bigcup\limits_{i:x\in U_{\al_i}}(U_{\al_i}\times V_{\al_i})\subset N.\nn\qedhere
\eea}
\end{proof}

\begin{crl}[\textcolor{OliveGreen}{Finite Tychonoff theorem}]
If $X,Y$ are compact spaces, then so is $X\times Y$.
\end{crl}
\begin{proof}
Let $\{E_\al\}$ be an open cover of $X\times Y$. For each $x\in X$, since $\{x\}\times Y$ is compact, let $\{x\}\times Y\subset\bigcup_{i=1}^{m_x} E_{\al^{x}_i}$. Define $N_x:=\bigcup_{i=1}^{m_x} E_{\al^{x}_i}$. Then by the tube lemma, there is a neighborhood $x\in U_x\subset X$ such that $U_x\times Y\subset N_x=\bigcup_{i=1}^{m_x} E_{\al^{x}_i}$. Thus, with compact $X\subset\bigcup_{j=1}^nU_{x_j}$, we have ~{\footnotesize$X\times Y=\bigcup_{j=1}^nU_{x_j}\times Y\subset \bigcup_{j=1}^nN_{x_j} =\bigcup_{j=1}^n\bigcup_{i=1}^{m_{x_j}}E_{\al^{x_j}_i}$}.
\end{proof}
The above result is true in the following general form but the proof is more involved.
\begin{thm*}[\textcolor{OliveGreen}{Tychonoff's theorem}]
For compact spaces $X_\al$, ~{\small $\prod_\al X_\al$} is compact.

(Proof: See Theorem \ref{TychThm})
\end{thm*}

\begin{lmm}\label{ProdClIntLmm}
If $X,Y$ are spaces and $A\subset X$, $B\subset Y$, then in $X\times Y$, (i) $\ol{A\times B}=\ol{A}\times\ol{B}$ and (ii) $(A\times B)^o=A^o\times B^o$. Hence, the boundary of a product of sets $A\times B\subset X\times Y$ is given by $\del(A\times B)=(\del A\times\ol{B})\cup(\ol{A}\times\del B)$.
\end{lmm}
\begin{proof}
(i) $\ol{A\times B}\subset \ol{\ol{A}\times\ol{B}}=\ol{A}\times\ol{B}$ (because $\ol{A}\times\ol{B}$ is closed). Also,
\bea
&& (a,b)\in \ol{A}\times\ol{B}~~\Ra~~a\in\ol{A},~~b\in\ol{B},\nn\\
&&~~\Ra~~\txt{every}~~[N(a)\times N(b)]\cap[A\times B]=[N(a)\cap B]\times[N(b)\cap B]\neq\emptyset\nn\\
&&~~\Ra~~\txt{every}~~N(a,b)\cap[A\times B]\neq\emptyset,~~\Ra~~(a,b)\in\ol{A\times B}.\nn
\eea
(ii) It follows that $(A\times B)^o=A^o\times B^o$ in $X\times Y$, since
\begin{align}
&(A\times B)^o=\left[\ol{(A\times B)^c}\right]^c=\left[\ol{A^c\times Y\cup X\times B^c}\right]^c=\left[\ol{A^c\times Y}\cup\ol{X\times B^c}\right]^c\nn\\
&~~~~=\left(\ol{A^c}\times Y\right)^c\cap\left(X\times\ol{B^c}\right)^c=\left([\ol{A^c}]^c\times Y\right)\cap\left(X\times[\ol{B^c}]^c\right)\nn\\
&~~~~=(A^o\times Y)\cap(X\times B^o)=A^o\times B^o.\nn\qedhere
\end{align}
\end{proof}

\subsection{Spaces with algebraic structures}

\begin{dfn}[\textcolor{blue}{\index{Topological! group}{Topological group}}]
A group $G$ with a topology $\T$ (on $G$) such that multiplication $G\times G\ra G,~(g,h)\mapsto gh$ and inversion $G\ra G,~g\mapsto g^{-1}$ are $\T$-continuous.
\end{dfn}

\begin{dfn}[\textcolor{blue}{\index{Topological! vector space (TVS)}{Topological vector space (TVS)}}]
A vector space $X$ (over $\mathbb{F}=\Real$ or $\Complex$) with a topology $\T$ such that addition $+:X\times X\ra X,~(x,y)\mapsto x+y$ and scalar multiplication $\mathbb{F}\times X\ra X,~(\al,x)\mapsto \al x$ are $\T$-continuous.
\end{dfn}

Examples of topological vector spaces are seminormed spaces in Definition \ref{VecSpExs} below.

\begin{dfn}[\textcolor{blue}{\index{Action}{Action of a topological group on a space}, \index{$G$-space}{$G$-space}}]
Let $G$ be a topological group (with identity $1_G$) and $X$ a space. An \ul{action} of $G$ on $X$ is a \ul{continuous map} $G\times X\ra X$, $(g,x)\mapsto gx$ such that ~(i) $g(g'x)=(gg')x$, and (ii) $1_Gx=x$, ~for all $g,g'\in G$, $x\in X$.

A \ul{$G$-space} is a space $X$ together with an action $G\times X\ra X$, $(g,x)\mapsto gx$.
\end{dfn}

\begin{dfn}[\textcolor{blue}{\index{Translation}{Translation} in a space, \index{Translation invariant space}{Translation invariant space}}]
If $X$ is a space and $V$ is a TVS, then a $V$-translation in $X$ is an action of $V$ on $X$ of the form
\bea
(V,+)\times X\ra X,~~(v,x)\mapsto x+v,~~~~(x+v_1)+v_2=x+(v_1+v_2)\eqv x+v_1+v_2.\nn
\eea
If the translation above exists, we say $X$ is a $V$-translation invariant space.
\end{dfn}

\begin{dfn}[\textcolor{blue}{\index{Scaling}{Scaling} in a space, \index{Scale invariant space}{Scale invariant space}}]
If $X$ is a space and $\mathbb{F}$ a field ($\Real$ or $\Complex$), then an $\mathbb{F}$-scaling in $X$ is an action of $\mathbb{F}$ on $X$ of the form
\bea
(\mathbb{F},\cdot)\times X\ra X,~~(v,x)\mapsto \al x,~~~~\al_1(\al_2 x)=(\al_1\al_2)x\eqv \al_1\al_2x.\nn
\eea
If the scaling above exists, we say $X$ is an $\mathbb{F}$-scale invariant space.
\end{dfn}

\begin{dfn}[\textcolor{blue}{
\index{Seminorm}{Seminorm},
\index{Seminormed space}{Seminormed space},
\index{Locally! convex space (LCS)}{Locally convex space (LCS)},
\index{Norm}{Norm},
\index{Normed space}{Normed space},
\index{Banach space}{Banach space}}]\label{VecSpExs}
Let $X$ be a vector space (over $\mathbb{F}=\Real$ or $\Complex$). A \ul{seminorm} on $X$ is a function $p:X\ra \Real$ such that $p(\al x+\beta y)\leq|\al|p(x)+|\beta|p(y)$ for all $x,y\in X$ and scalars $\al,\beta$. Associated with any seminorms $\{p_{a}\}_{{a}\in A}$ on $X$ is the topology (called \ul{seminorm topology}, making $X$ a \ul{seminormed space}) with \emph{subbase} (or generating set) consisting of all ``strips'' of the form $U_{\vep,{a}}(x)=\{y\in X:p_{a}(x-y)<\vep\}$, $\vep>0$, ${a}\in A$, $x\in X$.

A \ul{locally convex space} is a seminormed space $(X,\{p_{a}\}_{{a}\in A})$ such that for any $x\in X$, if $p_{a}(x)=0$ for all ${a}\in A$ then $x=0$. (In this case, we say $\{p_{a}\}_{{a}\in A}$ is a complete set of seminorms, which is the case iff the seminorm topology is Hausdorff).

A \ul{norm} on $X$ (making $X$ a \ul{normed space}) is a seminorm $p$ such that $p(x)=0$ implies $x=0$, in which case we write $p(x)$ as $\|x\|$. Every normed space $(X,\|\|)$ is a metric space with metric $d(x,y):=\|x-y\|$. A complete normed space is called a \ul{Banach space}. Examples of Banach spaces we will encounter include (for any set $A$ and a field $\mathbb{F}$ = $\Real$ or $\Complex$ -- depending on context)  the $\mathbb{F}$-vector space $\ell^\infty(A):=\left\{x:A\ra \mathbb{F}~\big|~\sup_{a\in A}|x(a)|<\infty\right\}$ with norm ~$\|x\|:=\sup_{a\in A}|x(a)|$.
\end{dfn}

\begin{rmks}
Let $X$ be an $\mathbb{F}$-vector space and $\{p_{a}\}_{{a}\in A}$ seminorms on $X$.\\
(i) By the subbase-to-base criterion in Lemma \ref{SubbaseCrit}, the seminorm topology on $X$ associated with $\{p_{a}\}_{{a}\in A}$ has a \emph{base} consisting of $\emptyset,X$ and \emph{finite intersections} of the ``strips''
\begin{align}
&\textstyle U_{\vep,{a}_1,...,{a}_n}(x)=\bigcap_{i=1}^nU_{\vep,{a}_i}(x),~~~~~~n\geq 1,~~\vep>0,~~a_i\in A,~~x\in X,\nn\\
&\textstyle~~~~=\{y\in X:p_{{a}_1,...,{a}_n}(x-y)<\vep\},~~~~~~~~p_{{a}_1,...,{a}_n}:=\max\{p_{{a}_1},...,p_{{a}_n}\}.\nn
\end{align}
(ii) If $p:X\ra\Real$ is any seminorm, the translation of $U_{\vep,p}(0):=\{p<\vep\}$ by a point $x\in X$ is
\bea
x+U_{\vep,p}(0)=x+\{z:p(z)<\vep\}=\{x+z:p(x+z-x)<\vep\}=\{y:p(y-x)<\vep\}=U_{\vep,p}(x).\nn
\eea
\end{rmks}

\begin{lmm}[\textcolor{OliveGreen}{\index{Minkowski function}{Minkowski function} of a convex set: \cite[Proposition 1.14, p.102]{conway}}]\label{MinkGauge}
Let $X$ be a locally convex space and $U\subset X$ a convex set that is ``\emph{symmetric}'' (so that $\ld U=U$ for any $\ld\in\mathbb{F}$ with $|\ld|=1$). Then $p:X\ra\Real$, $p(x):=\inf\left\{t> 0:{x\over t}\in U\right\}=\inf\left\{t> 0:x\in tU\right\}$ is a seminorm.
\end{lmm}
\begin{proof}
For any scalar $\al\in \mathbb{F}$, $p(\al x)=\inf\left\{t>0:{\al x\over t}\in U\right\}$ = $|\al|\inf\big\{{t\over|\al|}>0:{|\al|e^{i\arg\al}x\over t}\in e^{i\arg\al}U\big\}$ = $|\al|\inf\big\{s>0:{x\over s}\in U\big\}=|\al|p(x)$. Also, for any $t,t'>0$ with ${x\over t},{y\over t'}\in U$, we can choose $\ld\in (0,1)$ such that (by convexity of $U$) $\ld {x\over t}+(1-\ld){y\over t'}\in U$, along with ${\ld\over t}={1-\ld\over t'}$, imply ${1\over\ld}=1+{t'\over t}$ and ${\ld\over t}(x+y)={1\over t+t'}(x+y)\in U$, which in turn imply
\begin{align}
p(x)&+p(y)\textstyle=\inf\left\{t>0:{x\over t}\in U\right\}+\inf\left\{t>0:{y\over t}\in U\right\}\nn\\
&\textstyle\geq \inf\left\{t+t'>0:\ld{x\over t}+(1-\ld){y\over t'}\in U\right\}
 =  \inf\left\{t+t'>0:{x+y\over t+t'}\in U\right\}=p(x+y).\nn\qedhere
\end{align}
\end{proof}
Note (from the proof above) that if $U\subset X$ is any convex set, then (even without symmetry) the resulting Minkowski function $p:X\ra\Real$ is a \ul{sublinear functional} (See also \cite[Proposition 3.2,p.108]{conway}) in the sense it satisfies
\bea
p(\al x)=\al p(x)~~\txt{and}~~p(x+y)\leq p(x)+p(y),~~\txt{for all}~~\al\geq 0,~~x,y\in X.\nn
\eea

\begin{crl}\label{LCScrit}
A TVS is a LCS $\iff$ it has a base consisting of convex sets.
\end{crl}
\begin{proof}
($\Ra$) Seminorms produce convex neighborhoods: Seminorms are convex functions (so seminorm$^{-1}$($[0,\vep)$) is convex) and intersections of convex sets are convex sets.

($\La$) Conversely, convex neighborhoods come from seminorms: Note that the topology of a TVS is both translation-invariant and scale-invariant (i.e., translation and nonzero scaling preserve openness). Let $U$ be a convex set. From Lemma \ref{MinkGauge}, wlog, if $U$ is ``\emph{symmetric}'' (so that $\ld U=U$ for any $\ld\in\mathbb{F}$ with $|\ld|=1$), then its Minkowski function $p(x):=\inf\left\{t> 0:{x\over t}\in U\right\}=\inf\left\{t> 0:x\in tU\right\}$ is a seminorm. Again wlog, further choose the convex set $U$ to be ``\emph{absorbing}'' (so that $\ld U\subset U$ for any $\ld\in\mathbb{F}$ with $|\ld|<1$). (footnote\footnote{Note that if $U\subset X$ is a convex set with $0\in U$, then the absorbing property holds automatically.}). Observe that by definition, for any $t$, $x\in tU~\Ra~p(x)<t$, and so $tU\subset \{p<t\}$ for all $t$. In particular, $U\subset\{p<1\}$. On the other hand, suppose $x\in\{p<1\}$ but $x\not\in U$. Then $x\not\in tU\subset U$ for all $t$ such that $p(x)<t<1$. This contradicts the fact that $p(x)$ is the infimum of all $t$ such that $x\in tU$. Thus $\{p<1\}\subset U$, and so $U=\{p<1\}$.
\end{proof}

\begin{dfn}[\textcolor{blue}{Inner product, \index{Inner product space}{Inner product space}, \index{Hilbert space}{Hilbert space}}]
Let $X$ be a vector space (over $\mathbb{F}$ = $\Real$ or $\Complex$). An \ul{inner product} on $X$ (making $X$ an \ul{inner product space}) is a function $\langle\rangle:X\times X\ra\mathbb{F}$ such that for any $x,y,z\in X$ and scalar $\al$, (i) $\langle \al x,y\rangle=\al\langle x,y\rangle$ and $\langle x+y,z\rangle=\langle x,z\rangle+\langle y,z\rangle$, (ii) $\langle x,y\rangle=\overline{\langle y,x\rangle}$, (iii) $\langle x,x\rangle>0$ if $x\neq 0$. Every inner product space $(X,\langle\rangle)$ is a normed space with norm $\|x\|:=\sqrt{\langle x,x\rangle}$. A complete inner product space is called a \ul{Hilbert space}.
\end{dfn}
Note that a function $\langle\rangle:X\times X\ra\mathbb{F}$ that satisfies (i) and (iii) is called a \index{semi-inner product}{\ul{semi-inner product}}, and becomes an inner product only if it also satisfies the symmetry condition (ii).

\begin{dfn}[\textcolor{blue}{\index{Quotient! metric}{Quotient metric}, \index{Quotient! metric topology}{\textcolor{magenta}{Quotient metric topology}}, \index{Quotient! metric space}{Quotient metric space}}]
Let $(X,d)$ be a metric space and $\sim$ an equivalence relation on $X$. The quotient (pseudo)metric $\rho$ on the quotient space ${X\over\sim}=\{[x]=x_\sim:x\in X\}$, making $\big({X\over\sim},\rho\big)$ a quotient metric space, is
\bea
&&\rho([x],[y]):=\inf\Big\{l_\delta(c)~\big|~c:\{0,1,...,n\}\ra X,~c_0\in[x],~c_n\in[y],~n\in\Natural\Big\}\leq\dist([x],[y]),\nn\\
&&\textstyle l_\delta(c):=\sum\limits_{i=0}^n\delta(c_{i-1},c_i),
~~~~\delta(u,v):=\left\{
                     \begin{array}{ll}
                       0, & u\sim v  \\
                       d(u,v), & \txt{otherwise}
                     \end{array}
                   \right\},\nn
\eea
which is the smallest $\delta$-length $l_\delta(c)$ of finite chains $c=\{c_0,c_1,...,c_n\}$ of points in $X$ with endpoints in the equivalence classes $[x]$ and $[y]$, i.e., $c_0\in [x]$, $c_n\in[y]$.
\end{dfn}
It is clear that $\rho([x],[y])=\rho([y],[x])$. Also, $\rho([x],[y])\leq\rho([x],[z])+\rho([z],[y])$, since for any $\vep>0$, there exist two finite chains $c_{[x],[z]}$, $c_{[z],[y]}$ such that with $c_{[x],[y]}:=c_{[x],[z]}\cdot c_{[z],[y]}$ (another such finite chain by transitivity of $\sim$), we have
\bea
\rho([x],[y])\leq l_\delta\left(c_{[x],[y]}\right)\leq l_\delta\left(c_{[x],[z]}\right)+l_\delta\left(c_{[z],[y]}\right)<\rho([x],[z])+\rho([z],[y])+2\vep.\nn
\eea
For the special case where ${X\over\sim}={X\over A}$ for some $A\subset X$, we have
\bea
\rho([x],[y]) = \min\big\{d(x,y)~,~\dist(x,A)+\dist(y,A)\big\}.\nn
\eea

\begin{dfn}[\textcolor{blue}{Disjoint union, \textcolor{magenta}{Disjoint union topology}, \index{Disjoint union space}{Disjoint union space}}]
Let {\small $(X,\T_X)$, $(Y,\T_Y)$} be spaces. The \ul{disjoint union} {\small $X\sqcup Y$} of the sets {\small $X,Y$} is the set whose powerset {\small $\P(X\sqcup Y)=\{A\sqcup B:A\subset X,~B\subset Y\}$} is defined to be the  Cartesian pairing/product {\small $\P(X)\times\P(Y)=\{(A,B):A\subset X,~B\subset Y\}$} along with union and intersection given by
\bit[leftmargin=0.7cm]
\item[]\centering {\small $(A,B)\cup(A',B'):=(A\cup A',B\cup B')$, ~ $(A,B)\cap(A',B'):=(A\cap A',B\cap B')$}.
\eit
We write $A\sqcup B$ for {\small $(A,B)\in \P(X\sqcup Y)=\P(X)\times\P(Y)$}, so that {\small $\P(X\sqcup Y)=\{A\sqcup B:A\subset X,~B\subset Y\}$}, with the union and intersection rules
\bit[leftmargin=0.7cm]
\item[]\centering {\small $(A\sqcup B)\cup(A'\sqcup B'):=(A\cup A')\sqcup(B\cup B')$, ~ $(A\sqcup B)\cap(A'\sqcup B'):=(A\cap A')\sqcup(B\cap B')$}.
\eit
The \ul{disjoint union topology} $\T_{X\sqcup Y}$ on $X\sqcup Y$ (making {\small $X\sqcup Y=(X\sqcup Y,\T_{X\sqcup Y})$} a \ul{disjoint union space}) is given by {\small $\T_{X\sqcup Y}:=\T_X\sqcup\T_Y=\{U\sqcup V:U\in\T_X,V\in\T_Y\}$}.
\end{dfn}

\begin{dfn}[\textcolor{blue}{\index{Wedge sum}{Wedge sum} of spaces, \index{Basepoint}{Basepoint}}]
The wedge sum $\bigvee(X_\al,x_\al)$ of a collection $\{(X_\al,x_\al)\}$ of spaces $X_\al$ with basepoint $x_\al\in X_\al$ (a point selected on purpose) is
\bit
\item[] $\bigvee(X_\al,x_\al):=(X,x_0)$, ~ $X=\bigvee X_\al:={\bigsqcup X_\al\over\{x_\al\}}$, ~ $x_0:=x_\al~~\txt{for all}~~\al$.
\eit
\end{dfn}
\begin{lmm}[\textcolor{OliveGreen}{Gluing}]
If $f_\al:(X_\al,x_\al)\ra Y$ are continuous maps such that $f_\al(x_\al)=f_\beta(x_\beta)$ for all $\al,\beta$, then the map $F:\bigsqcup X_\al\ra Y$, $F(x)=f_\al(x)$ if $x\in X_\al$, is continuous. Moreover, there exists a unique continuous map $f:\bigvee(X_\al,x_\al)\ra Y$ such that $F=f\circ q$.
\end{lmm}
\begin{proof}
$F$ is continuous since $F^{-1}(V)=\bigsqcup f_\al^{-1}(V)$ for any (open) set $V\subset Y$. Moreover,
\bit
\item[] $x\sim x'$~~(i.e., $x=x'$ or $x,x'\in\{x_\al\}$) ~ $\Ra$ ~ $F(x)=F(x')$.
\eit
Thus, by Theorem \ref{UPQM} (universal property of quotient maps), we get the desired result.
{\small\bit\item[]
\bt
X_\al\ar[d,hook,"i_\al"']\ar[rr,"f_\al"]&& Y\\
\bigsqcup X_\al\ar[d,"q"]\ar[urr,"F"] &&\\
X=\bigvee X_\al\ar[uurr,dashed,"f"']
\et ~~~~
where~~$f|_{X_\al}:=f\circ q\circ i_\al= f_\al$ ~ for all ~~ $\al$. \qedhere
\eit}
\end{proof}

\subsection{Space types (by separation), Analytic and miscellaneous properties}

\begin{dfns}[\textcolor{blue}{Common types of spaces}] A space $X$ is called a
{\flushleft (a) \index{Kolmogorov space ($T_0$) space}{\ul{Kolmogorov space} (or \emph{$T_0$ space})}} if for any two distinct points $x,y\in X$, $x$ has a neighborhood that excludes $y$, or $y$ has a neighborhood that excludes $x$.
{\flushleft (b) \index{Fr$\acute{\txt{e}}$chet ($T_1$) space}{\ul{Fr$\acute{\txt{e}}$chet space}} (or \emph{$T_1$ space})} if for any two distinct points $x,y\in X$, $x$ has a neighborhood that excludes $y$, and $y$ also has a neighborhood that excludes $x$ (Equivalently, $X$ is a space in which every point is closed).
{\flushleft (c) \index{Hausdorff ($T_2$) space}{\ul{Hausdorff space}} (or \emph{$T_2$ space})} if distinct points $x,y\subset X$ have disjoint neighborhoods.
{\flushleft (d) \index{Regular space}{\ul{Regular space}}} if any closed set $C\subset X$ and any point $x\not\in C$ have disjoint neighborhoods (Equivalently, for every $x\in X$, every neighborhood $x\in U\subset X$ contains the closure of another neighborhood $x\in V\subset\ol{V}\subset U\subset X$). A regular Hausdorff space is called a \index{$T_3$ space}{\emph{$T_3$ space}}.

{\flushleft (e) \index{Normal space}{\ul{Normal space}}} if any two disjoint closed sets $C_1,C_2\subset X$, $C_1\cap C_2=\emptyset$, have disjoint neighborhoods. (Equivalently, for every closed set $C\subset X$, every neighborhood $C\subset U\subset X$ contains the closure of another neighborhood $C\subset V\subset\ol{V}\subset U\subset X$). A normal Hausdorff space (i.e., a normal $T_1$ space) is called a \index{$T_4$ space}{\emph{$T_4$ space}}.
\end{dfns}

\begin{lmm}[\textcolor{OliveGreen}{Bijection as homeomorphism}]\label{BijVsHomeo0}
A bijective continuous map $f:X\ra Y$ is a homeomorphism ($\iff$ it is open) $\iff$ it is closed. (Proof: If $U\subset X$ is open, then $f(U^c)=f(U)^c$ is closed, and so $f(U)$ is open, which shows $f^{-1}$ is continuous.)
\end{lmm}
Therefore, for any quotient map $q:X\ra Y$, if $q$ is injective on a closed (resp. open) set $A\subset X$, then $q|_A:A\ra Y$ is an imbedding, since if $B\subset A$ is closed (resp. open), then $B=q^{-1}\big(q(B)\big)$ implies $q(B)\subset q(A)$ is closed (resp. open).

\begin{lmm}[\textcolor{OliveGreen}{Closed image}]\label{ClosedImLmm}
Let $f :X \ra Y$ be a continuous map. If $X$ is compact and $Y$ is Hausdorff, then $f(X)$ is both compact and closed.
\end{lmm}
\begin{proof}
If $V_\al\subset Y$ is an open cover of $f(X)$, then $f^{-1}(V_\al)$ is an open cover of $X$, containing a finite subcover $f^{-1}(V_{\al_i})$, $i=1,...,n$. Thus $V_{\al_i}$ is a finite subcover of $f(X)$, i.e., $f(X)$ is compact. For closedness, pick $y\in Y\backslash f(X)$. Then for every $f(x)\in f(X)$, there exist disjoint open sets $A_x\ni y$ and $B_x\ni f(x)$. Since $f(X)$ is compact, its open cover $\{B_x:x\in X\}$ has a finite subcover $\{B_{x_i}\}_{i=1}^m$. Let $A=\bigcap A_{x_i}$ and $B=\bigcup B_{x_i}$. Then $A\cap f(X)=\emptyset$, since $f(X)\subset B$ and $A\cap B=\emptyset$. Hence $Y\backslash f(X)$ is open, i.e., $f(X)$ is closed.
\end{proof}
Note that the above lemma (by its proof) is an immediate corollary of two general facts, namely, (i) A continuous map takes compact sets to compact sets, (ii) A compact subset of a Hausdorff space is closed.

\begin{lmm}[\textcolor{OliveGreen}{Bijection as homeomorphism}]\label{BijVsHomeo}
Let $f:X\ra Y$ be a continuous map. If $X$ is compact and $Y$ is Hausdorff, then $f$ is a homeomorphism $\iff$ a bijection.
\end{lmm}
\begin{proof}
If $f$ is a homeomorphism, it is clear that $f$ is a bijection. Conversely, assume $f$ is a bijection. To show $f^{-1}:Y\ra X$ is continuous (i.e., that $f$ is open), let $A\subset X$ be open. Then the image ~$f(A)=(f^{-1})^{-1}\big((A^c)^c\big)=\big[(f^{-1})^{-1}(A^c)\big]^c=\big[f(A^c)\big]^c$~
is open since $A^c$ is compact (as a closed subset of a compact set) and so $f(A^c)$ is closed (by Lemma \ref{ClosedImLmm}).
\end{proof}

\begin{lmm}[\textcolor{OliveGreen}{Hausdorff quotient spaces}]\label{HausdQtSp}
Let $X$ be a regular Hausdorff space and $Y:=X_\sim=\{x_\sim:x\in X\}$ a quotient space of $X$. If the non-singleton equivalence classes $x^{ns}_\sim\subset X$\\
(i) are isolated closed sets (footnote\footnote{Isolated in the sense they have pairwise disjoint neighborhoods.}), and (ii) their union $\bigcup x^{ns}_\sim$ is closed in $X$, then $Y$ is a Hausdorff space. (In particular, if $X$ is a regular Hausdorff space and $C\subset X$ is closed, then ${X\over C}$ is Hausdorff.)
\end{lmm}
\begin{proof}
Let $q:X\ra Y$ be the quotient map, and let $y,y'\in Y$ be distinct. (a) If $y,y'\not\in\{x^{ns}_\sim:x\in X\}$ then by (ii) $y,y'$ have disjoint neighborhoods in $X$ on which $q$ is injective, and so have disjoint neighborhoods in $Y$. (b) If $y\not\in\{x^{ns}_\sim:x\in X\}$ and $y'\in\{x^{ns}_\sim:x\in X\}$, then by regularity, $y,y'$ have disjoint neighborhoods $O_y,O_{y'}$ in $X$ with $q$ injective on $O_y$ by (ii) and $O_{y'}$ disjoint from all other non-singleton equivalence classes by (i), and so have disjoint neighborhoods $q(O_y)$,$q(O_{y'})$ in $Y$. (c) If $y\in\{x^{ns}_\sim:x\in X\}$ and $y'\in\{x^{ns}_\sim:x\in X\}$, then by (ii), $y,y'$ have disjoint neighborhoods $O_y,O_{y'}$ in $X$ that are disjoint from all other non-singleton equivalence classes, and so have disjoint neighborhoods $q(O_y)$,$q(O_{y'})$ in $Y$.
\end{proof}

\begin{crl}
If $I=[0,1]\subset\Real$, then $S^n\cong{I^n\over\del I^n}$ and $I^{n+1}\cong{S^n\times I\over S^n\times\{1\}}$.
\end{crl}

\begin{lmm}[\textcolor{OliveGreen}{Non-Hausdorff quotient spaces}]\label{NonHausdQtSp}
If $X$ is a space and $O\subset X$ an open set such that $\del O\neq\emptyset$, then $Y:={X\over O}$ is not a $T_1$ space (hence a non-Hausdorff space).
\end{lmm}
\begin{proof}
Let $q:X\ra Y$ be the quotient map. Let $x\in\del O$. Since $O$ is open in $X$, we know $x\not\in O$, and so $x,O$ are distinct points of $Y$. If $B_x\subset Y$ is any open set containing $q(x)=x$, then because $q^{-1}(B_x)$ is open in $X$ and $x\in q^{-1}(B_x)$, we see that $q^{-1}(B_x)\cap O\neq\emptyset$. It follows that $B_x=q\big(q^{-1}(B_x)\big)\ni O$.
\end{proof}

\begin{dfn}[\textcolor{blue}{\index{Retraction}{Retraction}, Retract}]
Let $X$ be a space and $A\subset X$. Any continuous extension $r:X\ra A$ of the identity map $1_A:A\ra A$ (i.e., $r|_A=1_A$) is called a retraction of $X$ onto $A$. The subset $A\subset X$ is a retract (of $X$) if there exists a retraction $r:X\ra A$.
\end{dfn}

\begin{lmm}[\textcolor{OliveGreen}{Extension of continuous maps}]
A subspace $A\subset X$ is a retract $\iff$ every continuous map $f:A\ra Y$ extends to a continuous map $F:X\ra Y$.
\end{lmm}
\begin{proof}
If $A\subset X$ is a retract, then every continuous map $f:A\ra Y$ extends to a continuous map $F=f\circ r:X\sr{r}{\ral}A\sr{f}{\ral}Y$ (where $r$ is a retraction). Conversely, if every continuous map $f:A\ra Y$ extends to a continuous map $F:X\ra Y$, then in particular, $1_A:A\ra A$ extends to a retraction $r:X\ra A$.
\end{proof}

\begin{dfn}[\textcolor{blue}{\index{Uniform metric}{Uniform metric}, \textcolor{magenta}{Uniform metric topology}, \index{Function space}{Function space}}]
If $X$ is a space and $(M,d)$ a metric space (noting that $d_b(m,m'):=\min\big[d(m,m'),1\big]$ is another metric that induces the same topology on $M$ as $d$; see Lemma \ref{BdedMet}), the set of functions $\F(X,M)=\{\txt{maps}~f:X\ra M\}$ is a metric space (a \ul{function space}) with respect to the ``\ul{uniform metric}'' defined by
\bit
\item[] $d_u(f,g):=\sup\big\{d_b(f(x),g(x)):x\in X\big\}$.
\eit
If $\F(X,M)$ is given as a space, we will assume the $d_u$-topology, unless it is stated otherwise.
\end{dfn}
Note that if $M$ is bounded then (up to biLipschitz equivalence) we can simply define $d_u$ by $d_u(f,g):=\sup\big\{d(f(x),g(x)):x\in X\big\}$. In subsequent discussions involving $d_u$, we will for convenience not distinguish between $d_b=\min(d,1)$ and $d$. That is, $d$ in $d_u(f,g):=\sup\big\{d(f(x),g(x)):x\in X\big\}$ will be understood to be $d_b$ (especially when $M$ is unbounded).

\begin{lmm}[\textcolor{OliveGreen}{Completeness of the uniform metric}]\label{CompUniMet}
Let $X$ be a space and $(Y,d)$ a metric space. If $Y$ is complete, then the function space $\F(X,Y)$ is complete.
\end{lmm}
\begin{proof}
Assume $(Y,d)$ is complete, and let $\{f_n\}\subset\F(X,Y)$ be a Cauchy sequence. Then for any $\vep>0$, there exists $n_\vep$ such that ~$d_u(f_n,f_m)<\vep$ for all $n,m\geq n_\vep$. For every $x\in X$, $\{f_n(x)\}\subset Y$ is Cauchy (and so converges in $Y$) since ~$d\big(f_n(x),f_m(x)\big)\leq d_u(f_n,f_m)\ra 0$.

For each $x\in X$, let $f_n(x)\ra f(x)$ in $Y$. Then $\{f_n\}\subset\F(X,Y)$ converges to $f$. Indeed, we know ``$d_u(f_n,f_m)<\vep$ for $n,m\geq n_\vep$''. So, for fixed $n\geq n_\vep$, we have
\begin{align}
&d\left(f_{n+m}(x),f_n(x)\right)\leq d_u\left(f_{n+m},f_n\right)<\vep,~~~~\txt{for all}~~m\geq 1,~~x\in X,\nn\\
&~~\sr{(s)}{\Ra}~~d\left(f(x),f_n(x)\right)=\lim_{m\ra\infty}d\left(f_{n+m}(x),f_n(x)\right)\leq\vep,~~~~\txt{for all}~~x\in X,\nn\\
&~~\Ra~~d_u(f,f_n)\leq \vep,~~\Ra~~f_n\ra f~~\txt{in}~~\F(X,Y),\nn
\end{align}
where step (s) holds by {\footnotesize $\Big|d\big(f_{n+m}(x),f_n(x)\big)-d\big(f(x),f_n(x)\big)\Big|\leq d\big(f_{n+m}(x),f(x)\big)\ra0$} as $m\ra\infty$.
\end{proof}

\begin{lmm}[\textcolor{OliveGreen}{\index{Continuous limit theorem}{Continuous limit theorem}}]\label{ContLimThm}
Let $X$ be a space and $(Y,d)$ a metric space. (i) The space of continuous functions $C(X,Y)\subset\F(X,Y)$ is closed. (ii) If $X$ is also a metric space, the space of uniformly continuous functions $UC(X,Y)\subset\F(X,Y)$ is closed.
\end{lmm}
\begin{proof}
(i) Let $\{f_n\}\subset C(X,Y)$ be such that $f_n\ra f\in\F(X,Y)$. We need to show $f\in C(X,Y)$. That is, for any $x\in X$ and $\vep>0$, there exists an open set $O\ni x$ such that
\bea
f(O)\subset B_\vep(f(x)):=\{y\in Y:d(y,f(x))<\vep\}.\nn
\eea
Since $f_n\ra f$ uniformly, there exists $N=N_\vep$ such that
\bea
d\big(f_n(x),f(x)\big)\leq d_u(f_n,f)<\vep/3,~~~~\txt{for all}~~n\geq N,~~x\in X.\nn
\eea
Fix $x\in X$. Since $f_N$ is continuous at $x$, there exists a neighborhood $O_x$ of $x$ such that
\bea
d\big(f_N(x),f_N(x')\big)<\vep/3,~~~~\txt{for all}~~x'\in O_x.\nn
\eea
Therefore, $f$ is also continuous at $x$, since for all $x'\in O_x$
\bea
d\big(f(x),f(x')\big)\leq d\big(f(x),f_N(x)\big)+d\big(f_N(x),f_N(x')\big)+d\big(f_N(x'),f(x')\big)<\vep.\nn
\eea
{\flushleft (ii)} Now, assume $X$ is a metric space and, as before, let $\{g_n\}\subset UC(X,Y)$ be such that $g_n\ra g\in\F(X,Y)$. We need to show $g\in UC(X,Y)$. Since $g_n\ra g$ uniformly, there exists $N=N_\vep$ such that
\bea
d\big(g_n(x),g(x)\big)\leq d_u(g_n,g)<\vep/3,~~\txt{for all}~~n\geq N,~~x\in X.\nn
\eea
By uniform continuity of $g_N$, there exists $\delta=\delta_\vep>0$ such that
\bea
d(x,x')<\delta~~\Ra~~d\big(g_N(x),g_N(x')\big)<\vep/3,~~\txt{for all}~~x,x'\in X.\nn
\eea
Thus, $g$ is also uniformly continuous, since (for all $x,x'\in X$) the bound $d(x,x')<\delta$ implies
\bea
d\big(g(x),g(x')\big)\leq d\big(g(x),g_N(x)\big)+d\big(g_N(x),g_N(x')\big)+d\big(g_N(x'),g(x')\big)<\vep.\nn\qedhere
\eea
\end{proof}

\begin{crl}[\textcolor{OliveGreen}{Continuous limit theorem}]\label{ContLimCrl}
Let $X$ be a space and $(Y,d)$ a complete metric space. (i) The space of continuous functions $C(X,Y)\subset\F(X,Y)$ is complete. (ii) If $X$ is also a metric space, the space of uniformly continuous functions $UC(X,Y)\subset\F(X,Y)$ is complete.
\end{crl}

\begin{lmm}[\textcolor{OliveGreen}{\index{Uniform extension theorem}{Uniform extension theorem}}]\label{UniExtThm}
Let $X,Y$ be metric spaces, and $f:E\subset X\ra Y$. If (i) $f$ is uniformly continuous, (ii) $E$ is dense, and (iii) $Y$ is complete, then $f$ extends to a unique uniformly continuous map $F:X\ra Y$.

{\flushleft(Note: By its proof, the result remains true if "uniformly" is replaced with "Lipschitz".)}
\end{lmm}
\begin{proof}
Define a map $F:X\ra Y$ by $F(x):=\lim f(e_n)$ for any sequence $e_n\in E$ such that $e_n\ra x$. $F$ is \ul{\emph{well defined}} because if $e_n\ra x$, $e_n'\ra x$ and $f(e_n)\ra y$, $f(e_n')\ra y'$, then
\bea
d(y,y')\leq d(y,f(e_n))+d\left(f(e_n),f(e_n')\right)+d(f(e_n'),y')\ra 0,~~\Ra~~y=y',\nn
\eea
where $d\left(f(e_n),f(e_n')\right)\ra 0$ by uniform continuity of $f$ (because $d(e_n,e_n')\leq d(e_n,x)+d(x,e_n')\ra 0$). It is also clear that $F|_E=f$, i.e., $F$ is an \ul{\emph{extension}} of $f$.

To see that $F$ is \ul{\emph{uniformly continuous}} (hence \ul{\emph{unique}} since continuous maps that agree on a dense set agree everywhere), observe that
\bea
d(F(x),F(x'))=d\left(\lim f(e_n),\lim f(e_n')\right)=\lim_{k,l} d\left(f(e_k),f(e_l')\right)\sr{(s)}{\ral} 0~~\txt{uniformly as}~~x\ra x',\nn
\eea
where step (s) holds because there always exists $N_{x,x'}=N\big(d(x,x')\big)\geq 1$ such that
\bea
&& d(e_k,x)\leq d(x,x')~~\txt{and}~~d(e'_l,x')\leq d(x,x'),~~~~\txt{for all}~~k,l\geq N_{x,x'},\\
&&~~\Ra~~d(e_k,e_l')\leq d(e_k,x)+d(x,x')+d(x',e_l')\leq 3d(x,x'),~~~~\txt{for all}~~k,l\geq N_{x,x'}.\nn
\eea
\ul{Step (s)}: By uniform continuity of $f$, for any $\vep>0$, there exists $\delta=\delta_\vep$ such that for $k,l\geq 1$,
\bea
d(e_k,e'_l)<3\delta~~\Ra~~d(f(e_k),f(e'_l))<\vep/3.
\eea
Also, there exist {\small $0<\delta_x,\delta_{x'}\leq d(x,x')$} such that for any {\small $k,l\geq N_{x,x'}$, $d(x,e_k)<\delta_x$ $\Ra$ $d(F(x),f(e_k))<\vep/3$}, and {\small $d(x',e'_l)<\delta_{x'}$ $\Ra$ $d(F(x'),f(e'_l))<\vep/3$}. Hence, $d(x,x')<\delta$ implies
\[
d(F(x),F(x'))\leq d(F(x),f(e_k))+d(f(e_k),f(e'_l))+d(f(e'_l),F(x'))<\vep,~~\txt{if}~~k,l\geq N_{x,x'}. \qedhere
\]
\end{proof}

\begin{dfn}[\textcolor{blue}{\index{Dense map}{Dense map}}]
A map of spaces $f:X\ra Y$ is dense if its image is dense, i.e, $\ol{f(X)}=Y$.
\end{dfn}
\begin{dfn}[\textcolor{blue}{\index{Completion}{Completion} of a metric space}]
Let $X$ be a metric space. A metric space $Z$ is a completion of $X$ if (i) $Z$ is complete and (ii) there exists a dense isometric imbedding $\vphi:X\hookrightarrow Z$.
\end{dfn}
By the following theorem, every metric space has a completion that is unique up to isometry.

\begin{thm}[\textcolor{OliveGreen}{Existence and uniqueness of a completion}]\label{MetComThm}
For any metric space $X$, there is a complete metric space $\wt{X}$ and an isometric imbedding $i:X\hookrightarrow\wt{X}$ such that $i(X)$ is dense in $\wt{X}$. Moreover, $\wt{X}$ is unique up to isometry.
\end{thm}
\begin{proof}
\underline{\emph{Existence}}: Let $C=\{\txt{all Cauchy sequences in $X$}\}$. Define a relation $\sim$ on $C$ as:
\bea
(x_j)\sim(y_j)~~~~\txt{if}~~~~d(x_j,y_j)\ra 0,~~~~~~~~\txt{for all}~~~~(x_j),(y_j)\in C.\nn
\eea
Then $\sim$ is an equivalence relation due to the following.
\bit
\item[(a)] $d(x_j,x_j)=0\ra 0$ (So $\sim$ is reflexive).
\item[(b)] $d(x_j,y_j)\ra 0$ implies $d(y_j,x_j)=d(x_j,y_j)\ra 0$ (So $\sim$ is symmetric).
\item[(c)] $d(x_j,y_j)\ra 0$, $d(y_j,z_j)\ra 0$ implies $d(x_j,z_j)\leq d(x_j,y_j)+d(y_j,z_j)\ra 0$ (So $\sim$ is transitive).
\eit
Let $\big[(x_j)\big]:=\Big\{(y_j)\in C:(y_j)\sim(x_j)\Big\}$ be the equivalence class of $(x_j)$ in $C$. Consider the metric space $\wt{X}~:=~C/\sim~=~\Big\{\wt{x}=\big[(x_j)\big]~\Big|~(x_j)\in C\Big\}$, with metric
\bea
\wt{d}\left(\wt{x},\wt{y}\right):=\lim_{j\ra\infty} d(x_j,y_j),\nn
\eea
where $\wt{d}$ is well defined (i.e., exists and is independent of the choice of representatives of the equivalence classes) because $\big(d(x_j,y_j)\big)$ converges as a real Cauchy sequence, and if $(x'_j)\sim(x_j)$ and $(y'_j)\sim(y_j)$, then $|d(x'_j,y'_j)-d(x_j,y_j)|\leq d(x'_j,x_j)+d(y_j,y'_j)\ra 0$ implies
\bea
\lim_j d(x'_j,y'_j)=\lim_j d(x_j,y_j).\nn
\eea

\ul{Define the inclusion} $i:X\hookrightarrow \wt{X}$ by $i(x)=\big[(x,x,\cdots)\big]$, where $(x,x,\cdots)$ denotes the constant sequence. Then $i$ is isometric since~ $\wt{d}(i(x),i(y))=\lim d(x,y)=d(x,y)$. Note that every sequence $(x_j)$ such that $x_j\ra x$ belongs to the equivalence class $\big[(x,x,\cdots)\big]$, and so
\bea
i(x):=\big[(x,x,\cdots)\big]=\Big\{(x_j)\in C:~x_j\ra x\Big\}.\nn
\eea

\ul{To show that $i(X)$ is dense} in $\wt{X}$ we show that $\wt{X}\subset \overline{i(X)}=i(X)\cup i(X)'$. Let $\wt{x}=\big[(x_j)\big]\in\wt{X}$. If $(x_j)$ converges, then $\wt{x}\in i(X)\subset\overline{i(X)}$. If $(x_j)$ does not converge, let $\vep>0$ be given. We must show that $N_\vep(\wt{x})\backslash\{\wt{x}\}$ contains an element of $i(X)$. Recall that $\wt{y}=\big[(y_j)\big]\in N_\vep(\wt{x})$
\bea
~~\iff~~\wt{d}(\wt{x},\wt{y})=\lim_j d(x_j,y_j)<\vep.\nn
\eea
Since $(x_j)$ is Cauchy, there is $N=N(\vep)$ such that $d(x_j,x_{j'})<\vep$ for $j,j'\geq N$, and so
for any $k\geq N$,
\bea
&&\wt{d}\left(\wt{x},\big[i(x_k)\big]\right)=\lim_j d(x_j,x_k)<\vep,~~~~\Ra~~~~\big[i(x_k)\big]=\big[(x_k,x_k,\cdots)\big]\in N_\vep(\wt{x})\backslash\{\wt{x}\}.\nn
\eea
Therefore $i(X)$ is dense in $\wt{X}$.

\ul{Next we show that $\wt{X}$ is complete}. Let $(\wt{x}_k)$, where $\wt{x}_k=\big[(x_{jk})\big]$, be a Cauchy sequence in $\wt{X}$. We must show that $(\wt{x}_k)$ converges. Let $\vep>0$ be given. Then there is $N_1=N_1(\vep)$ such that
\bea
&&\wt{d}(\wt{x}_k,\wt{x}_{k'})=\lim_j~d(x_{jk},x_{jk'})<\vep~~~~\txt{for}~~~~k,k'\geq N_1.\nn
\eea
Since for each pair $k,k'\geq N_1$, $d(x_{jk},x_{jk'})\ra \wt{d}(\wt{x}_k,\wt{x}_{k'})$, it follows that there is $N_2(k,k')=N_2(\vep,k,k')$ with
\bea
&& |d(x_{jk},x_{jk'})-\wt{d}(\wt{x}_k,\wt{x}_{k'})|<\vep~~~~for~~~~j\geq N_2(k,k'),~~~~k,k'\geq N_1,\nn\\
&&~~\Ra~~d(x_{jk},x_{jk'})< \wt{d}(\wt{x}_k,\wt{x}_{k'})+\vep<2\vep~~~~for~~~~j\geq N_2(k,k'),~~~~k,k'\geq N_1.\nn
\eea
Consider the diagonal sequence $(x_{kk})$. Then because each representative of $\wt{x}_k$ is Cauchy,
\bea
d(x_{kk},x_{k'k'})\leq d(x_{kk},x_{jk})+d(x_{jk},x_{jk'})+d(x_{jk'},x_{kk'})<3(2\vep)=6\vep~~~~~~\txt{for}~~k,k'\geq N_1.\nn
\eea
Thus $(x_{kk})$ is a Cauchy sequence, i.e., $\big[(x_{kk})\big]\in \wt{X}$. Therefore, $(\wt{x}_k)$ converges in $\wt{X}$, since
\bea
\wt{d}\left(\wt{x}_k,\big[(x_{jj})\big]\right)=\lim_j~d(x_{jk},x_{jj})<6\vep~~~~\txt{for}~~k\geq N_1,~~~~\Ra~~~~\wt{x}_k\ra\big[(x_{jj})\big]\in\wt{X}.\nn
\eea

\underline{\emph{Uniqueness}}: Let $i_1:X\hookrightarrow\wt{X}_1$ and $i_2:X\hookrightarrow\wt{X}_2$ be isometric imbeddings, where $\wt{X}_1,\wt{X}_2$ are complete, $i_1(X)$ is dense in $\wt{X}_1$, and $i_2(X)$ is dense in $\wt{X}_2$. Then we get the map
\bea
j:i_1(X)\subset\wt{X}_1\ra \wt{X}_2,~i_1(x)\mapsto i_2(x)\nn
\eea
which is uniformly continuous on $i_1(X)$ as an isometric map, because
\bea
d(j(i_1(x)),j(i_1(y)))=d(i_2(x),i_2(y))=d(x,y)=d(i_1(x),i_1(y)),~~~~\txt{for all}~~x,y\in X.\nn
\eea
Since $i_1(X)$ is dense in $\wt{X}_1$, and $\wt{X}_2$ is complete, it follows by Theorem \ref{UniExtThm} that $j$ extends to a unique continuous map $\wt{X}_1\ra\wt{X}_2$. Hence we have an isometry $j:\wt{X}_1\ra \wt{X}_2$.
\end{proof}

\begin{thm}[\textcolor{OliveGreen}{\index{Baire category theorem}{Baire category theorem}}]
Let $(X,d)$ be a complete metric space. If $\{D_n\}_{n=1}^\infty$ is a countable collection of dense open (or dense $G_\delta$) sets, then $\bigcap_{n=1}^\infty D_n$ is dense. (\ul{Note}: A set is called a $G_\delta$ set if it is a countable intersection of open sets. The complements of $G_\delta$ sets, i.e., countable unions of closed sets, are called $F_\sigma$ sets.)
\end{thm}
\begin{proof}
We show every $x\in X$ lies in $\ol{\bigcap_{n=1}^\infty D_n}$, i.e., $B_r(x)\cap\bigcap_{n=1}^\infty D_n\neq\emptyset$ for all $r>0$. Observe that  $B_r(x)\cap\bigcap_{n=1}^\infty D_n=\bigcap_{n=1}^\infty B_r(x)\cap D_n=\bigcap_{n=1}^\infty E_n$,  where $E_n:=B_r(x)\cap D_n$ are dense open sets in the complete metric space $Y:=\ol{B_r(x)}=(\ol{B_r(x)},d)$.

Let $x_1\in E_1$. Since $E_1$ is open in $Y$, some $B_{2r_1}(x_1)\subset E_1$, and so $\ol{B_{r_1}(x_1)}\subset E_1$. Since $E_2$ is dense in $Y$, some $x_2\in E_2\cap B_{r_1}(x_1)\backslash\{x_1\}$, and so (as above) some $\ol{B_{r_2}(x_2)}\subset E_2\cap B_{r_1}(x_1)\backslash\{x_1\}$. Similarly, some $\ol{B_{r_3}(x_3)}\subset E_3\cap B_{r_2}(x_2)\backslash\{x_2\}$. Continuing this way, at the $n$th step, we get $\ol{B_{r_n}(x_n)}\subset E_n\cap B_{r_{n-1}}(x_{n-1})\backslash\{x_{n-1}\}$, resulting in a decreasing sequence of closed sets $\ol{B_{r_n}(x_n)}\supset\ol{B_{r_{n+1}}(x_{n+1})}$ of sizes $2^{-n}\geq r_n>r_{n+1}\ra 0$.

This means the sequence $\{x_n\}\subset Y$ is Cauchy, and so converges to a point $z\in Y$, since $Y$ is complete. Note that $z\in\bigcap\ol{B_{r_n}(x_n)}\subset \bigcap_{n=1}^\infty E_n$, otherwise, if some $\ol{B_{r_{n'}}(x_{n'})}\not\ni z$ then some neighborhood of $z$ misses a tail of $\{x_n\}$, which is a contradiction. Hence, $\bigcap E_n\neq\emptyset$.
\end{proof}

\section{Classification of Spaces: Algebraic Topology}\label{PrelimsAT} 
We will not be concerned with a detailed exposition of the algebraic topology of finite subset spaces. However, some questions concerning their metric geometry (which is our main concern) will require some basic knowledge of algebraic topology.

There are many equivalent ways of specifying the topology of a space, some of which we have already encountered. The topology of a space can be variously a point-set topology, variously a metric topology, variously a subspace topology, variously a quotient topology, variously a product topology, variously a disjoint union topology, and so on. In other words, a topological space is only unique up to homeomorphism. A property of a space is said to be \ul{\emph{topologically invariant}}, or called a \index{Topological! invariant}{\ul{\emph{topological invariant}}}, if it is invariant under homeomorphisms (i.e., if it does not depend on the way the topology is specified).

Since a topological space is only unique up to homeomorphism, in classifying spaces it is essential to consider \emph{classification methods that are topologically invariant} (in the sense they are based on topological invariants of some sort, and so do not depend on the way the topology of a space is specified). Algebraic topology is a study of \emph{topologically invariant classification methods} that are based on \emph{algebraic topological invariants} (i.e., topological invariants constructed from \emph{groups}, \emph{rings}, and \emph{modules}) defined over spaces. We will briefly describe two of such classification methods, namely, \emph{homotopy} and \emph{homology}.

\subsection{Homotopy}
{\flushleft Denote} the interval $[0,1]\subset\Real$ by $I$. On the unit $n$-sphere, $S^n=\del B_1(0)\subset\Real^{n+1}$, we will denote the north pole by $N:=(0,\cdots,0,1)$ and the south pole by $S:=(0,\cdots,0,-1)$. If $f:X\ra Y$ and $A\subset X$, $B\subset Y$ are such that $f(A)\subset B$, then we write $f:(X,A)\ra (Y,B)$.
\begin{dfn}[\textcolor{blue}{\index{Path}{Path}, \index{Reverse path}{Reverse path}, \index{Reverse map}{Reverse map}}]
If $X$ is a space, a continuous map of the form $\gamma:I\ra X$ is called a path in $X$ from $\gamma(0)$ to $\gamma(1)$. The reverse path of $\gamma$ is $\ol{\gamma}:I\ra X$, $\ol{\gamma}(t):=\gamma(1-t)$. Similarly, given a continuous map $f:I^n\ra X$, we define (up to homeomorphism of $I^n$) its reverse map $\ol{f}:I^n\ra X$ by $\ol{f}(t_1,...,t_n):=f(1-t_1,t_2,\cdots,t_n)$.
\end{dfn}
\begin{dfn}[\textcolor{blue}{\index{Homotopy!}{Homotopy}, \index{Invariant homotopy}{Invariant homotopy}}] Let $X,Y$ be spaces.
A family of maps $\{f_t:X\ra Y\}_{t\in I}$ is called a homotopy (or continuous family of continuous maps) if the joint map ~$F:X\times I\ra Y$, $(x,t)\mapsto f_t(x)$~ is continuous. The homotopy is invariant if $F(x,t)$ is constant in $t$, i.e., if the map $F:X\times I\ra Y$ is equivalent to a single map $f:X\ra Y$.
\end{dfn}
The joint map $F$ is itself called the homotopy, because any continuous map $F:X\times I\ra Y,~(x,t)\mapsto F(x,t)$ gives a homotopy $\big\{f_t:X\ra Y,~x\mapsto f_t(x):=F(x,t)\big\}_{t\in I}$ -- Recall that a continuous multivariate function is continuous in each variable separately. Note that a homotopy $F:X\times I\ra Y$ is equivalently a continuous family of paths $\gamma_X$ $:=$ $\left\{\gamma_x:I\ra Y\right\}_{x\in X}$, where $\gamma_x(t)=F(x,t)=f_t(x)$ is a path in $Y$ from $F(x,0)=f_0(x)$ to $F(x,1)=f_1(x)$. Thus, the homotopy is a \emph{\ul{continuous deformation}} in $Y$ of the image $F(X,0)=f_0(X)\subset Y$, through the images $F(X,t)=f_t(X)$ into the image $F(X,1)=f_1(X)\subset Y$, of $X$.

\begin{dfn}[\textcolor{blue}{\index{Homotopy! relative to a set}{Homotopy relative to a set}}]
If $A\subset X$, a homotopy $F:X\times I\ra Y$ is called a homotopy relative to $A$ if the subhomotopy $F|_{A\times I}:A\times I\ra Y$ is invariant. That is, the homotopy leaves the image $F(A,0)=F(A,t)=F(A,1)\subset Y$ of $A$ undeformed.
\end{dfn}

\begin{dfn}[\textcolor{blue}{\index{Homotopic maps}{Homotopic maps}}] Continuous maps $f,g:X\ra Y$ are homotopic (written $f\simeq g$) if there is a homotopy $\{f_t:X\ra Y\}_{t\in I}$ such that $f_0=f$, $f_1=g$ (i.e., such that the continuous map $F(x,t)=f_t(x)$ satisfies ~$F|_{X\times\{0\}}=f$,~ $F|_{X\times\{1\}}=g$). That is, the image $f(X)$ can be continuously deformed into the image $g(X)$.
\end{dfn}

\begin{dfn}[\textcolor{blue}{\index{Nullhomotopic map}{Nullhomotopic map}}] A continuous map $f:X\ra Y$ is null homotopic if it is homotopic to a constant map $c_{y_0}:X\ra Y,~x\mapsto y_0$ (for some $y_0\in Y$). That is, the image $f(X)\subset Y$ can be continuously deformed into a point.
\end{dfn}

\begin{dfn}[\textcolor{blue}{\index{Homotopy! equivalence}{Homotopy equivalence}, Homotopy equivalent spaces, Homotopy inverse}]
A continuous map $f:X\ra Y$ is a \ul{homotopy equivalence} if there exists a continuous map $g:Y\ra X$ (called \ul{homotopy inverse} of $f$) such that $f\circ g\simeq 1_Y$ and $g\circ f\simeq 1_X$. We say $X$ and $Y$ are \ul{homotopy equivalent} (or of the \ul{same homotopy type}), written $X\simeq Y$.
\end{dfn}

\begin{dfn}[\textcolor{blue}{\index{Contractible space}{Contractible space}}]
A space $X$ is contractible if $X\simeq\{x_0\}$ for a point $x_0\in X$ (or equivalently, the identity $1_X:X\ra X$ is nullhomotopic).
\end{dfn}

\begin{lmm}[\textcolor{OliveGreen}{Contractibility criterion}]  A space $X$ is contractible $\iff$ $1_X:X\ra X$ is nullhomotopic.
\end{lmm}
\begin{proof}
{\flushleft ($\Ra$):} Assume $X\simeq\{x_0\}$ for some $x_0\in X$. Let $f:X\ra\{x_0\}$ be a homotopy equivalence with a homotopy inverse $g:\{x_0\}\ra X$. Then $1_X\simeq g\circ f=c_{g(x_0)}:X\ra X$. ($\La$): Assume $1_X\simeq c_{x_0}:X\ra X$ for some $x_0\in X$. Then $c_{x_0}$ gives a map $p:X\ra \{x_0\},~x\mapsto c_{x_0}(x)=x_0$. Also, we always have the inclusion map $i:\{x_0\}\hookrightarrow X,~x_0\ra x_0$. Thus,
\bit
\item[] $i\circ p=c_{x_0}\simeq 1_X$, ~~ $p\circ i=1_{\{x_0\}}\simeq 1_{\{x_0\}}$,
\eit
which shows $ p:X\ra\{x_0\}$ is a homotopy equivalence.
\end{proof}

\begin{crl}[\textcolor{OliveGreen}{Contractibility criterion}]  A space $X$ is contractible $\iff$ every continuous map $f:Z\ra X$ is nullhomotopic $\iff$ every continuous map $f:X\ra Z$ is nullhomotopic.
\end{crl}

\begin{dfn}[\textcolor{blue}{Gluing of compatible continuous maps on $[0,1]^n$}]~ Let $I=[0,1]$. Given continuous maps $f,g:[0,1]\ra X$ such that $f(1)=g(0)$, we can glue them to obtain the continuous map
~{\small$f\cdot g:[0,1]\ra X$,~ $f\cdot g(t)=
\left\{
  \begin{array}{ll}
    f(2t), & t\in[0,1/2] \\
    g(2t-1), & t\in[1/2,1]
  \end{array}
\right\}$}.
Similarly, if continuous maps $f,g:I^n\ra X$ are such that for the two ``opposite'' faces $A=\{0\}\times I^{n-1}\subset\del I^n$ and $B=\{1\}\times I^{n-1}\subset\del I^n$ of $I^n$ we have (up to homeomorphism of $I^n$)
\begin{align}
f(1,t_2,\cdots,t_n)=g(0,t_2,\cdots,t_n),~~\txt{for all}~~(1,t_2,\cdots,t_n)\in B,~~(0,t_2,\cdots,t_n)\in A,\nn
\end{align}
then we can glue them together (along the coordinate $t_1$) to obtain the continuous map
\bit
\item[] $f\cdot g(t_1,...,t_n):=\left\{
  \begin{array}{ll}
    f(2t_1,t_2,\cdots,t_n), & t_1\in[0,1/2] \\
   g(2t_1-1,t_2,\cdots,t_n), & t_1\in[1/2,1]
  \end{array}
\right\}$.
\eit
\end{dfn}
\begin{rmk}[\textcolor{OliveGreen}{Eckmann-Hilton}]\label{ProdHtyRmk}
Let $f,g,f\cdot g,f\ast g:(I^2,\del I^2)\ra (X,x_0)$ satisfy
\begin{align}
f\cdot g(t,s):=\left\{
  \begin{array}{ll}
    f(2t,s), & t\in[0,1/2] \\
   g(2t-1,s), & t\in[1/2,1]
  \end{array}
\right\},~~~~f\ast g(t,s):=\left\{
  \begin{array}{ll}
    f(t,2s), & s\in[0,1/2] \\
   g(t,2s-1), & s\in[1/2,1]
  \end{array}
\right\}.\nn
\end{align}
Then for any continuous maps $f,g,h,k:(I^2,\del I^2)\ra (X,x_0)$, we have
{\footnotesize\begin{align}
&\Big((f\ast g)\cdot(h\ast k)\Big)(t,s)=\left\{
  \begin{array}{ll}
    f\ast g(2t,s), & t\in[0,1/2] \\
   h\ast k(2t-1,s), & t\in[1/2,1]
  \end{array}
\right\}\nn\\
&~~=\left\{
  \begin{array}{ll}
    f(2t,2s), & t\in[0,1/2],s\in[0,1/2] \\
        g(2t,2s-1), & t\in[0,1/2],s\in[1/2,1] \\\\
   h(2t-1,2s), & t\in[1/2,1],s\in[0,1/2]\nn\\
      k(2t-1,2s-1), & t\in[1/2,1],s\in[1/2,1]
  \end{array}
\right\}=\left\{
  \begin{array}{ll}
    f(2t,2s), & t\in[0,1/2],s\in[0,1/2] \\
        h(2t-1,2s), & t\in[1/2,1],s\in[0,1/2]\\\\
        g(2t,2s-1), & t\in[0,1/2],s\in[1/2,1] \\
      k(2t-1,2s-1), & t\in[1/2,1],s\in[1/2,1]
  \end{array}
\right\}\nn\\
&~~=\left\{
  \begin{array}{ll}
    f\cdot h(t,2s), & s\in[0,1/2] \\
   g\cdot k(t,2s-1), & s\in[1/2,1]
  \end{array}
\right\}=\Big((f\cdot h)\ast(g\cdot k)\Big)(t,s).\nn
\end{align}}
Note that this remains true if we consider maps $f,g,f\cdot g,f\ast g:(I^n,\del I^n)\ra (X,x_0)$ and suppress/fix all but any two of the $n$ arguments.
\end{rmk}

\begin{dfn}[\textcolor{blue}{\index{Homotopy! groups}{Homotopy groups} of a space, \index{Homotopy! class}{Homotopy class} of a map}]
Let $X$ be a space. Fix $x_0\in X$, and an integer $n\geq 0$. The $n$th \ul{homotopy group} of $X$ at $x_0$ is
{\small\begin{align}
&\pi_n(X,x_0):=\big\{\txt{continuous}~f:(I^n,\del I^n)\ra (X,x_0),~\txt{up to homotopy relative to $\del I^n$}\big\}\nn\\
&~~~~={\big\{\txt{continuous}~f:(I^n,\del I^n)\ra(X,x_0)\big\}\over\sim}=\big\{[f]:~\txt{for continuous}~f:(I^n,\del I^n)\ra (X,x_0)\big\}\nn\\
&~~~~=\left\{\txt{continuous}~{f'}:(S^n,N)\ra (X,x_0),~\txt{up to homotopy relative to $\{N\}$}\right\}\nn\\
&~~~~={\{\txt{continuous}~{f'}:(S^n,N)\ra(X,x_0)\}\over\sim}=\left\{[{f'}]:~\txt{for continuous}~{f'}:(S^n,N)\ra(X,x_0)\right\},\nn\\
&\txt{where}~~~~[f]~:=~\{g:g\sim f\},~~~~\txt{with}~~~~f\sim g~\iff~~f\simeq g~~\txt{relative to}~~\del I^n,\nn\\
&\txt{and}~~~~~~~[{f'}]~:=~\{{g'}:{g'}\sim{f'}\},~~~~\txt{with}~~~~{f'}\sim{g'}~\iff~~{f'}\simeq{g'}~~\txt{relative to}~~\{N\},\nn
\end{align}}
as a set with a product given by ~$[f][g]:=[f\cdot g]$, ~for all ~$[f],[g]~\in~\pi_n(X,x_0)$. The set of homotopic maps $[f]$ is called the \ul{homotopy class} of $f$.
\end{dfn}

\begin{rmks*}
(i) By the universal property of quotient maps, any continuous map $f:I^n\ra X$, such that $f|_{\del I^n}$ is constant, corresponds to a unique continuous map ${f'}:S^n\cong I^n/\del I^n\ra X$.
\bit
\item[]\bt
 &I^n\ar[d,"q"']\ar[rr,"f"]&& X\\
S^n\ar[r,draw=none,"\cong"description]& {I^n/\del I^n}\ar[urr,"{f'}"'] &&
\et ~~~~ where ~~~~ $f=f'\circ q$.
\eit
Thus, the elements of $\pi_n(X,x_0)$ are like ``wedge sums of continuous images of $S^n$ in $X$''.
{\flushleft (ii)} $\pi_n(X,x_0)$ is a \emph{\ul{group}} for every $n\geq 1$, with \emph{\ul{identity}} $e=[c_{x_0}]$ and \emph{\ul{inverse}} $[f]^{-1}=[\ol{f}]$.
{\flushleft (iii)} $\pi_n(X,x_0)$ is an \emph{\ul{abelian group}} for every $n\geq 2$. This is because the identity in Remark \ref{ProdHtyRmk} shows homotopy groups defined using the new product $\ast$ are isomorphic to the earlier homotopy groups defined using $\cdot$. Hence, by setting $f\simeq k\simeq c_{x_0}$ in the identity, we see that the homotopy groups are abelian for $n\geq 2$. Equivalently, we have $f\cdot g\simeq g\cdot f$. By the description of $\pi_n(X,x_0)$ above, this product resembles a ``wedge sum of two continuous images of $S^n$ in $X$'', i.e., a surjective continuous map $f'\cdot g':=(f\cdot g)':S^n\ra f'(S^n)\vee g'(S^n)$.
\end{rmks*}

\begin{dfn}[\textcolor{blue}{\index{Homotopy! $n$-chains}{Homotopy $n$-chains}}]
The homotopy $n$-chain of $X$ is $\C_n(X,x_0)$ $:=$ $\{\txt{continuous}~f:(I^n,\del I^n)\ra (X,x_0)\}$. Using $\C_n(X,x_0)$, we can write {\small $\pi_n(X,x_0)={\C_n(X,x_0)\over\sim}=\big\{[f]:f\in\C_n(X,x_0)\big\}$}, where $f\sim g$ if $f\simeq g$ relative to $\del I^n$.
\end{dfn}

\begin{lmm}[\textcolor{OliveGreen}{$n$-connectedness criteria: \cite[p 346]{hatcher2001}}]\label{NoHoleCrit}
For any space $X$, the following three conditions are equivalent (where $S^n\cong\del I^{n+1}$).
\bit[leftmargin=0.9cm]
\item[(a)] Every continuous map $S^n\ra X$ is nullhomotopic (i.e., homotopic to a constant map).
\item[(b)] Every continuous map $S^n\ra X$ extends to a continuous map $I^{n+1} \ra X$.
\item[(c)] $\pi_n(X,x_0) = 0$ for all $x_0 \in X$.
\eit
\end{lmm}
\begin{proof}
{\flushleft $\bullet$(a)$\Ra$(b):} Let $f:S^n\ra X$ be a continuous map. By (a), there is a point $x_0\in X$ and a continuous map $F:S^n\times I\ra X$ with $F(s,0)=f(s)$, $F(s,1)=x_0$ for all $s\in S$.
Since $F|_{S^n\times\{1\}}$ is constant, there exists (by the universal property of quotient maps) a unique continuous map $\wt{F}:(S^n\times I)/(S^n\times\{1\})\ra X$ such that $\wt{F}\circ q=F$.
\bit
\item[]\hspace{2cm}\bt
&S^n\times I\ar[d,"q"']\ar[rr,"F"]&& X\\
I^{n+1}\ar[r,"\cong"',"\vphi"]&(S^n\times I)/(S^n\times\{1\})\ar[urr,dashed,"\wt{F}"']
\et\eit
{\flushleft $\bullet$(b)$\Ra$(c):} Let $f:(I^n,\del I^n)\ra (X,x_0)$ be a continuous map. Since $f|_{\del I^n}$ is constant, there exists (by the universal property of quotient maps) a unique continuous map $\wt{f}:I^n/\del I^n\ra X$ such that the following diagram commutes, i.e., $\wt{f}\circ q=f$.
\bit
\item[]\bt
&I^n\ar[d,"q"']\ar[rr,"f"]&& X\\
(\del I^{n+1},N)\eqv (S^n,N)\ar[r,"\cong"',"\vphi"]&(I^n/\del I^n,\del I^n)\ar[urr,dashed,"\wt{f}"']
\et\eit
By (b), $\wt{f}\circ\vphi:(S^n,N)\ra X$ extends to a continuous map $F:I^{n+1}\ra X$. Since $I^{n+1}$ is contractible, there is a homotopy $h:(I^{n+1},N)\times I\ra (I^{n+1},N)$ such that $h|_{I^{n+1}\times\{0\}}=1_{I^{n+1}}$, $h|_{I^{n+1}\times\{1\}}$ is constant, and $h(N,t)=N$ (i.e., the homotopy fixes the base point). Thus, $H=F\circ h:I^{n+1}\times I\sr{h}{\ral}I^{n+1}\sr{F}{\ral} X$ is a homotopy between $F$ and a constant map (i.e., $F:I^{n+1}\ra X$ is nullhomotopic). Thus, $F|_{\del I^{n+1}}$ (and hence $f$) is homotopic to a constant map. Hence, $\pi_n(X,x_0)=0$.

{\flushleft $\bullet$(c)$\Ra$(a):} Let $f:S^n\ra X$ be a continuous map. The map {\small $g:=f\circ\vphi\circ q:I^n\sr{q}{\ral}(I^n/\del I^n,\del I^n)\sr{\vphi}{\ral}(S^n,N)\sr{f}{\ral}(X,x_0)$}, with $\vphi(\del I^n)=N$ and $x_0:=f(N)$, is in $\pi_n(X,x_0)$. Moreover, there exists (by universal property of quotient maps) a unique continuous map $\wt{g}:I^n/\del I^n\ra X$ giving a commutative diagram
\bit\item[]\hspace{3cm}
\bt
I^n\ar[d,"q"']\ar[rr,"g"]&& X\\
I^n/\del I^n\ar[urr,dashed,near start,"\wt{g}"]\ar[rr,"\vphi","\cong"'] && S^n\ar[u,"f"']
\et
\eit
By (c), $g$ is homotopic to the constant map $c_{x_0}:I^n\ra X,~s\mapsto x_0$ via a homotopy $\{g_t\}_{t\in I}$. Thus, for each $t\in I$, there exists (by universal property of quotient maps) a unique continuous map $\wt{g}_t:I^n/\del I^n\ra X$ giving a commutative diagram
\bit
\item[]\bt
I^n\ar[d,"q"']\ar[rr,"g_t"]&& X\\
I^n/\del I^n\ar[urr,dashed,near start,"\wt{g}_t"]\ar[rr,"\vphi","\cong"'] && S^n\ar[u,"f_t:=\wt{g}_t\circ\vphi^{-1}"']
\et
~~~~
\bt
I^n\times I\ar[d,"q\times 1_I"']\ar[rr,"G"]&& X\\
I^n/\del I^n\times I\ar[urr,dashed,near start,"\wt{G}"]\ar[rr,"\vphi\times 1_I","\cong"'] && S^n\times I\ar[u,"F:=\wt{G}\circ(\vphi\times 1_I)^{-1}"']
\et
\eit
which show $f$ is homotopic to a constant map via the homotopy $f_t:=\wt{g}_t\circ\vphi^{-1}$.
\end{proof}

\begin{lmm}[\textcolor{OliveGreen}{Homotopy groups of a retract are subgroups}]\label{HtyRetInc} If $A\subset X$ is a retract (and $x_0\in A$), the map $\pi_n(A,x_0) \ra \pi_n(X,x_0)$ induced by the inclusion $A \hookrightarrow X$ is injective.
\end{lmm}
\begin{proof}
By hypotheses, we have a continuous map $X\sr{r}{\ral}A$ (along with the inclusion $A\sr{i}{\ral}X$) such that $r\circ i=1_A:A\sr{i}{\hookrightarrow}X\sr{r}{\ral}A$. Thus, the induced homomorphisms
\begin{align}
&\C_n(A,x_0)\sr{i_\#}{\ral}\C_n(X,x_0),~f\mapsto i\circ f,~~~~\pi_n(A,x_0)\sr{i_\ast}{\ral}\pi_n(X,x_0),~[f]\mapsto [i_\#(f)]=[i\circ f],\nn\\
&\C_n(X,x_0)\sr{r_\#}{\ral}\C_n(A,x_0),~f\mapsto r\circ f,~~~~\pi_n(X,x_0)\sr{r_\ast}{\ral}\pi_n(A,x_0),~[f]\mapsto [r_\#(f)]=[r\circ f]\nn
\end{align}
satisfy $r_\#\circ i_\#=(r\circ i)_\#=(1_A)_\#=1_{\C_n(A,x_0)}$, and $r_\ast\circ i_\ast=(r\circ i)_\ast=(1_A)_\ast=1_{\pi_n(A,x_0)}$.
\end{proof}

\begin{lmm}[\textcolor{OliveGreen}{Homotopy groups of spheres: \cite[Cor. 4.9, p.349]{hatcher2001}}]\label{HtyGpSph}
{\small $\pi_i(S^n)\cong
\left\{
  \begin{array}{ll}
    \Integer, & i=n\\
    0, & i<n
  \end{array}
\right\}$}
\end{lmm}

\subsection{Homology}
{\flushleft This section} provides the same desired information as the homotopy section.
\begin{dfn}[\textcolor{blue}{\index{Chain! complex}{Chain complex}, \index{Chain! group}{Chain group}, \index{Boundary map}{Boundary map}}]
A chain complex $C_\#$ is a sequence of abelian groups (or modules) $C_n$, $n\in\Integer$, with homomorphisms ~$\del_n:C_n\ra C_{n-1}$~ such that ~$\del_n\circ\del_{n+1}=0$. The chain complex $C_\#=(C_\#,\del)$ is often expressed in the form
\begin{align}
C_\#:\cdots\ral C_{n+1}\sr{\del_{n+1}}{\ral}C_n\sr{\del_n}{\ral}C_{n-1}\sr{\del_{n-1}}{\ral}\cdots,\nn
\end{align}
where $C_n$ is called the $n$th \textcolor{blue}{chain group}, and $\del_n$ is called the \index{Differential}{\textcolor{blue}{differential}} (or \textcolor{blue}{boundary map}).
\end{dfn}

Note that $\del_n\circ\del_{n+1}=0$ $\iff$ $\im(\del_{n+1})\subset\ker(\del_n)$ in $C_n$.

\begin{dfn}[\textcolor{blue}{\index{Subcomplex}{Subcomplex}, \index{Quotient! complex}{Quotient complex}}]
Let $A_\#=(A_\#,\del^A)$ and $B_\#=(B_\#,\del^B)$ be complexes. Then $A_\#$ is a subcomplex of $B_\#$, written $A_\#\subset B_\#$, if for each $n$, we have (i) $A_n\subset B_n$ and (ii) $\del_n^B|_{A_n}=\del^A_n$. The quotient complex $B_\#/A_\#=(B_\#/A_\#,\ol{\del})$ of $B_\#$ by $A_\#$ is given by
\bea
(B_\#/A_\#)_n:=B_n/A_n,~~~~\ol{\del}_n(b+A_n):=\del^B_nb+A_{n-1}.\nn
\eea
\end{dfn}

\begin{dfn}[\textcolor{blue}{\index{Homology}{Homology} of a chain complex $C$, \index{Cycles}{Cycles}, \index{Boundaries}{Boundaries}, \index{Homology classes}{Homology classes}}]
The $n$th homology of $C_\#$ is the quotient
{\small\begin{align}
H_n(C_\#):={\ker(\del_n)\over\im(\del_{n+1})}={Z_n(C_\#)\over B_n(C_\#)},~~~
\left.
   \begin{array}{l}
    Z_n(C_\#):=\ker\del_n=\{\txt{\textcolor{blue}{$n$-cycles}}\}=\del_n^{-1}(0)\\
    B_n(C_\#):=\im\del_{n+1}=\{\txt{\textcolor{blue}{$n$-boundaries}}\}=\del_{n+1}(C_{n+1})
   \end{array}
 \right.\nn
\end{align}}
where the elements $[z]=z+\im\del_{n+1}$ of $H_n(C_\#)$ are called $n$th \index{Homology classes}{\ul{homology classes}} of $C_\#$.
\end{dfn}

\begin{dfn}[\textcolor{blue}{\index{Chain! map}{Chain map}, \index{Chain! isomorphism}{Chain isomorphism}}]
A \ul{chain map} (or map of chain complexes) $f_\#:(C_\#,\del)\ra(C_\#',\del')$ is a family of homomorphisms $\left\{f_n:C_n\ra C_n'\right\}_{n\in\Integer}$ such that $\del_n'\circ f_n=f_{n-1}\circ\del_n$ for each $n$ (written as $\del'\circ f_\#=f_\#\circ\del$).
\vspace{-0.1cm}
\bc\bt
\cdots\ar[r]& C_n\ar[d,"f_n"]\ar[r,"\del_n"] & C_{n-1}\ar[d,"f_{n-1}"]\ar[r,"\del_{n-1}"] &\cdots\\
\cdots\ar[r]& C_n'\ar[r,"\del_n'"] & C_{n-1}'\ar[r,"\del_{n-1}'"] &\cdots\\
\et\ec
\vspace{-0.7cm}
A chain map $f_\#:C_\#\ra C_\#'$ is a \ul{chain isomorphism} (making $C_\#,C_\#'$ \ul{isomorphic}, written $C_\#\cong C_\#'$) if each $f_n:C_n\ra C_n'$ is an isomorphism. In this case, we have the \ul{inverse chain map} $f_\#^{-1}:C_\#'\ra C_\#$ given by $(f_\#^{-1})_n:=(f_n)^{-1}$.
\end{dfn}
Note that if $f_\#,g_\#:C_\#\ra C_\#'$ are chain maps, then so are the sum $f_\#+g_\#:=\{f_n+g_n\}:C_\#\ra C_\#'$ and the scalar multiple $\ld f_\#:=\{\ld f_n\}:C_\#\ra C_\#'$.

\begin{rmk*}[\textcolor{OliveGreen}{\index{Induced homomorphism}{Induced homomorphism}}]
A chain map $f_\#:C_\#\ra C_\#'$ induces a homomorphism $f_\ast:H_n(C_\#)\ra H_n(C_\#')$, $f_\ast([z]):=[f_\#(z)]$, which is well defined because $f_\#(Z_n(C_\#))\subset Z_n(C_\#')$ and $f_\#(B_n(C_\#))\subset B_n(C_\#')$. The induction operation $\ast$ satisfies (i) $(1_{C_\#})_\ast=1_{H_n(C_\#)}$ and (ii) $(f\circ g)_\ast=f_\ast\circ g_\ast$.
\end{rmk*}

\begin{dfn}[\textcolor{blue}{\index{Quasi-isomorphism}{Quasi-isomorphism}}]
A chain map $f_\#:C_\#\ra C_\#'$ is called a quasi-isomorphism if the induced map $f_\ast:H_n(C_\#)\ra H_n(C_\#')$ is an isomorphism for each $n$.
\end{dfn}

\begin{rmk}[\textcolor{OliveGreen}{Delta chain map}]
For any chain complexes $C_\#,C_\#'$ and any family of homomorphisms {\footnotesize$h:=\{h_n:C_n\ra C_{n+1}'\}_{n\in\Integer}$}, the map {\footnotesize$\delta_\#(h):=h\circ\del+\del'\circ h:C_\#\ra C_\#'$} given by
    \bea
    \delta_n(h):=h_{n-1}\circ\del_n+\del_{n+1}'\circ h_n:~C_n\ra C_n',~~~~n\in\Integer,\nn
    \eea
is a chain map, because $\del'\circ\delta_\#(h)=\del'h\circ\del=\delta_\#(h)\circ\del$.
\bc\bt
\cdots\ar[r] &C_{n+1}\ar[rr]&& C_n\ar[ddll,dashed,"h_n"']\ar[dd,near end,"\delta_n(h)"']\ar[rr,"\del_n"] && C_{n-1}\ar[ddll,dashed,"h_{n-1}"']\ar[r,"{\del_{n-1}}"] &\cdots\\
 & && && &\\
\cdots\ar[r]&C_{n+1}'\ar[rr,"\del_{n+1}'"]&& C_n'\ar[rr,"\del_n'"] && C_{n-1}'\ar[r,"{\del_{n-1}'}"] &\cdots
\et\ec
Note that ~$\delta_\#(h_1+h_2)=\delta_\# h_1+\delta_\# h_2$.
\end{rmk}
\begin{dfn}[\textcolor{blue}{\index{Chain! homotopy}{Chain homotopy}, \index{Homotopic chain maps}{Homotopic chain maps}}]
A \ul{chain homotopy} between chain maps $f_\#,g_\#:C_\#\ra C_\#'$ is a family of homomorphisms $h:=\{h_n:C_n\ra C_{n+1}'\}_{n\in\Integer}$ such that $f_n-g_n=\delta_n(h)$ for each $n$, also written as $f_\#-g_\#=\delta_\#(h)$.
\bc\bt
\cdots\ar[r] &C_{n+1}\ar[rr]&& C_n\ar[ddll,dashed,"h_n"']\ar[dd,near end,"f_n-g_n"']\ar[rr,"\del_n"] && C_{n-1}\ar[ddll,dashed,"h_{n-1}"']\ar[r,"{\del_{n-1}}"] &\cdots\\
 & && && &\\
\cdots\ar[r]&C_{n+1}'\ar[rr,"\del_{n+1}'"]&& C_n'\ar[rr,"\del_n'"] && C_{n-1}'\ar[r,"{\del_{n-1}'}"] &\cdots
\et\ec
In this case, we say $f_\#,g_\#$ are \ul{$h$-homotopic}, and write $f_\#\sr{h}{\simeq}g_\#$.
\end{dfn}

\begin{rmks*}[\textcolor{OliveGreen}{Properties of chain homotopy}] (1) Chain homotopy is an equivalence relation: If $f_\#,f_\#',f_\#'':C_\#\ra C_\#'$ are chain maps, then (a) $f_\#\sr{0}{\simeq}f_\#$, (b) $f_\#\sr{h}{\simeq}g_\#$ $\Ra$ $g_\#\sr{-h}{\simeq}f_\#$, and (c) $f_\#\sr{h}{\simeq}f_\#'$, $f_\#'\sr{h'}{\simeq}f_\#''$ $\Ra$ $f_\#\sr{h+h'}{\simeq}f_\#''$. (2) Moreover, if $f_\#\simeq g_\#$, then $f_\ast=g_\ast$.
\end{rmks*}

\begin{dfn}[\textcolor{blue}{\index{Chain! homotopy equivalence}{Chain homotopy equivalence}, Chain homotopy inverse, Chain homotopy equivalent complexes}] A chain map $f_\#:C_\#\ra C_\#'$ is a \ul{chain homotopy equivalence} if there exists a chain map $g_\#:C_\#'\ra C_\#$ such that $f_\#\circ g_\#\simeq 1_{C_\#'}$ and $g_\#\circ f_\#\simeq 1_{C_\#}$. ($g_\#$ is called a \ul{homotopy inverse} of $f_\#$). In this case, we say $C_\#$ and $C_\#'$ are \ul{chain homotopy equivalent}, written $C_\#\simeq C_\#'$.
\end{dfn}

\begin{rmk*}~ $\left\{\substack{\txt{chain}\\ \txt{isomorphisms}}\right\}$ ~$\subset$~ $\left\{\substack{\txt{chain homotopy}\\ \txt{equivalences}}\right\}$ ~$\subset$~ $\left\{\substack{\txt{quasi-}\\ \txt{isomorphisms}}\right\}$
\end{rmk*}

\begin{dfn}[\textcolor{blue}{\index{Nullhomotopic chain map}{Nullhomotopic chain map}}]
A chain map $f_\#:C_\#\ra C_\#'$ such that $f_\#\simeq 0$.
\end{dfn}

\begin{dfn}[\textcolor{blue}{\index{Contractible complex}{Contractible complex}}]
A complex $C_\#$ such that $C_\#\simeq\{0\}$.
\begin{prp}
A complex $C_\#$ is contractible if and only if~ $1_{C_\#}\simeq 0:C_\#\ra C_\#$.
\end{prp}
\begin{proof}
($\Ra$): Assume $C_\#\simeq\{0\}$, and let $f_\#:C_\#\ra\{0\}$ be a homotopy equivalence with homotopy inverse $g_\#:\{0\}\ra C_\#$. Then $f_\#\circ g_\#\simeq 1_{\{0\}}$ and $0=g_\#\circ f_\#\simeq 1_{C_\#}$.\\
($\La$): Assume $1_{C_\#}\simeq 0$. Let $f_\#:C_\#\ra\{0\},~c\mapsto 0$ and $g_\#=i_{\{0\}}:\{0\}\hookrightarrow C_\#$. Then $f_\#\circ g_\#=1_{\{0\}}$, $g_\#\circ f_\#=0\simeq 1_{C_\#}$, and so $f_\#:C_\#\ra\{0\}$ is a homotopy equivalence.
\end{proof}

\end{dfn}

\begin{dfn}[\textcolor{blue}{\index{Exact (or acyclic) complex}{Exact (or acyclic) complex}}]
A complex $C_\#$ is exact (or acyclic) if $H_n(C_\#)=\{0\}$ for all $n\in\Integer$.
\end{dfn}

\begin{rmks*}
(1) A contractible complex is acyclic. Indeed, if $C_\#\simeq\{0\}$, then $1_{C_\#}\simeq 0$, and so
\bea
1_{H_n(C_\#)}=0:H_n(C_\#)\ra H_n(C_\#),~~\Ra~~H_n(C_\#)=\{0\}.\nn
\eea
(2) If a complex $C_\#$ is acyclic, then ~$i_{\{0\}}:\{0\}\hookrightarrow C_\#$~ is a quasi-isomorphism, since
    \bea
    \{0\}\sr{i_{\{0\}}}{\hal} C_\#~~\Ra~~\{0\}=H_n(\{0\})\sr{H_n(i_{\{0\}})}{\ral}H_n(C_\#)=\{0\}.\nn
    \eea
\end{rmks*}

{\flushleft\dotfill}

\begin{dfn}[\textcolor{blue}{\index{Points in general position}{Points in general position}}]
Points $v_0,...,v_n\in\Real^m$ ($m\geq n+1$) are in general position if they are not all contained in any affine hyperplane of dimension $<n$, i.e., if $v_1-v_0$, $v_2-v_0$, $\cdots$, $v_n-v_0$ are linearly independent vectors.
\end{dfn}

\begin{dfn}[\textcolor{blue}{\index{Affine! simplex}{Affine $n$-simplex}, Vertices, Edges, Orientation, Faces}]

An affine $n$-dimensional simplex (or affine $n$-simplex) spanned by $n+1$ points $v_0,...,v_n$ in general position (called \ul{vertices} of the simplex) in $\Real^m$ is (i) the $n$-dimensional set $[v_0,...,v_n]$ $:=$ the convex hull of $\{v_0,...,v_n\}$ = the smallest convex subset of $\Real^m$ containing $v_0,...,v_n$,
(ii) together with a strict \ul{ordering} of the vertices ~$v_0\prec v_1\prec\cdots\prec v_n$ that induces an \ul{orientation} on each \ul{edge} $[v_i,v_j]$ denoted by an arrow
$\bt v_i\ar[r]& v_j,\et$ if $i<j$. As a set,
{\small\begin{align}
[v_0,\cdots,v_n]=
\left\{\substack{
\txt{convex combinations}\\
\txt{of $v_0,...,v_n$}
}\right\}=\left\{{\textstyle\sum\limits_{i=0}^n}t_iv_i:{\textstyle\sum\limits_{i=0}^n}t_i=1,~t_i\geq 0,~~\txt{for all}~~i\right\}.\nn
\end{align}}Each subset of vertices $\{v_{i_0},...,v_{i_k}\}\subset\{v_0,...,v_n\}$, $k\leq n-1$, determines a $k$-dimensional subsimplex $[v_{i_0},...,v_{i_k}]\subset[v_0,...,v_k]$ called a \ul{$k$-face} of $[v_1,...,v_n]$, where the orientation of each face is inherited from that of the $n$-simplex (via the ordering $v_0\prec v_1\prec\cdots\prec v_n$). Naturally, the \ul{$0$-faces} are called vertices, and the \ul{$1$-faces} are called edges.
\end{dfn}

\begin{dfn}[\textcolor{blue}{\index{Baryocentric coordinates}{Baryocentric coordinates}}]
Given a point ~$p=\sum_{i=0}^nt_iv_i\in[v_0,...,v_n]$,~ the numbers $t_0,...,t_n\in[0,1]$ satisfying $\sum_{i=0}^n t_i=1$ are called the baryocentric coordinates of $p$.
\end{dfn}

\begin{dfn}[\textcolor{blue}{\index{Standard simplex}{Standard $n$-simplex}}]
This is the affine $n$-simplex ~$\Delta^n:=[e_0,...,e_n]\subset\Real^{n+1}$,~ where $e_j=(0,...,0,x_j=1,0,...,0)$ are the standard unit vectors, and as a set,
\begin{align}
\Delta^n=\big\{{\textstyle\sum}_{i=0}^nt_ie_i:{\textstyle\sum}_{i=0}^nt_i=1,~t_i\in[0,1]\big\}=\left\{(t_0,...,t_n)\in[0,1]^{n+1}:{\textstyle\sum}_{i=0}^nt_i=1\right\}.\nn
\end{align}
\end{dfn}

\begin{dfn}[\textcolor{blue}{Canonical homeomorphism of an affine simplex}]
For each affine $n$-simplex $[v_0,...,v_n]\subset\Real^m$, we have the canonical homeomorphism
\bea
\textstyle f_{[v_0,...,v_n]}:\Delta^n\ra [v_0,...,v_n],~~(t_0,...,t_n)\mapsto\sum_{i=0}^nt_iv_i,~~~~\txt{(i.e., $e_i\mapsto v_i$)}.\nn
\eea
(Thus, for each dimension $n$, there is only one affine $n$-simplex $=\Delta^n$ up to homeomorphism.)
\end{dfn}

\begin{dfn}[\textcolor{blue}{Boundary faces of affine simplex}]
The $i$th boundary face of $[v_0,...,v_n]$ is ~$[v_0,...,\wh{v}_i,...,v_n]:=[v_0,...,v_{i-1},v_{i+1},...,v_n]\subset[v_0,...,v_n]$.
\end{dfn}

\begin{dfn}[\textcolor{blue}{\index{Singular simplex}{Singular $n$-simplex} in a space $X$, Set of singular simplices}]
This is a continuous map $\sigma:\Delta^n\ra X$. We write ~$S_n(X):=\C(\Delta^n,X)=\{\txt{singular $n$-simplices in $X$}\}$.
\end{dfn}
{\flushleft As} a continuous map $\Delta^n\ra\Real^m$, $m\geq n+1$, the affine $n$-simplex $[v_0,...,v_n]$ is given by $[v_0,...,v_n](e_i):=v_i$, extended by linearity, i.e., $[v_0,...,v_n]\left(\sum t_ie_i\right):=\sum t_iv_i$. Equivalently,
\begin{align}
[v_0,...,v_n]=i_{[v_0,...,v_n]}\circ f_{[v_0,...,v_n]}:\Delta^n\sr{f_{[...]}}{\ral}[v_0,...,v_n]\sr{i_{[...]}}{\hal}\Real^m.\nn
\end{align}
In particular, $[v_0,...,v_n]\big|_{[e_0,...,\wh{e}_i,...,e_n]}=[v_0,...,\wh{v}_i,...,v_n]$.

\begin{dfn}[\textcolor{blue}{Boundary face maps of standard simplex, Boundary face (map) of an affine simplex}]
The $i$th \ul{boundary face map} ~$F_i^n:\Delta^{n-1}\hookrightarrow\Delta^n$~ of $\Delta^n$ is the composition

\bt
F_i^n:\Delta^{n-1}\ar[r,"\cong"]&{[e_0,...,\wh{e}_i,...,e_n]}\ar[r,hook]&\Delta^n.
\et
\vspace{-0.2cm}
{\flushleft For} $\sigma\in S_n(X)$, we define the $i$th \ul{boundary face (map)} of $\sigma$ as
\bea
\sigma|_{[e_0,...,\wh{e}_i,...,e_n]}:=\sigma\circ F_i^n:\Delta^{n-1}\sr{F_i^n}{\hal}\Delta^n\sr{\sigma}{\ral}X,~~~~(\txt{an element of $S_{n-1}(X)$}).\nn
\eea
\end{dfn}

\begin{dfn}[\textcolor{blue}{\index{Singular chain complex}{Singular chain complex} of a space $X$}]
This is the chain complex
\bea
\bt[row sep = tiny,column sep = small]
  & n=2& n=1& n=0& n=-1 &n=-2 &  \\
C_\#(X):~~\cdots\ar[r,"\del_3"]&C_2(X)\ar[r,"\del_2"] & C_1(X)\ar[r,"\del_1"]&C_0(X)\ar[r,"0"]&0\ar[r,"0"]& 0\ar[r,"0"]& \cdots,
\et\nn
\eea
where the chain group $C_n(X)$ is the free abelian group generated by $S_n(X)$, or equivalently, the free $\Integer$-module spanned by $S_n(X)$, i.e.,
\begin{align}
C_n(X):=\txt{Span}_\Integer S_n(X)={\textstyle\bigoplus\limits_{\sigma\in S_n(X)}}\Integer\sigma=\left\{
     \substack{\txt{formal finite sums}\\\txt{(singular $n$-chains in $X$)}}
~{\textstyle\sum} n_i\sigma_i~\big|~n_i\in\Integer,\sigma_i\in S_n(X)\right\},\nn
\end{align}
and the boundary map $C_n(X)\sr{\del_n}{\ral} C_{n-1}(X)$ is given by ~$\del_n\sum n_i\sigma_i:=\sum n_i\del_n\sigma_i$, where
\begin{align}
\del_n\sigma={\textstyle\sum\limits_{i=0}^n}(-1)^i\sigma|_{[e_0,...,\wh{e}_i,...,e_n]}:={\textstyle\sum\limits_{i=0}^n}(-1)^i\sigma\circ F_i^n,~~~~\txt{for}~~n\geq 0,~~~~~~\txt{for all}~~\sigma\in S_n(X).\nn
\end{align}
\end{dfn}
\begin{lmm}
The boundary map of the singular chain complex satisfies $\del_{n+1}\circ\del_n=0$.
\end{lmm}
\begin{proof}
It suffices by linearity to verify the desired property on $\sigma\in S_n(X)$.
\begin{align}
&\textstyle(\del_n\circ\del_{n+1})(\sigma)=\del_n\sum\limits_{i=0}^{n+1}(-1)^i\sigma|_{[e_0,...,\wh{e}_i,...,e_{n+1}]}=\sum\limits_{i=0}^{n+1}(-1)^i\del_n\sigma|_{[e_0,...,\wh{e}_i,...,e_{n+1}]}\nn\\
&\textstyle~~~~=\sum\limits_{i=0}^{n+1}(-1)^i\left(\sum\limits_{j=0}^{i-1}(-1)^j\sigma|_{[e_0,...,\wh{e}_j,...,\wh{e}_i,...,e_{n+1}]}+\sum\limits_{j=i+1}^{n+1}(-1)^{j-1}\sigma|_{[e_0,...,\wh{e}_i,...,\wh{e}_j,...,e_{n+1}]}\right)\nn\\
&\textstyle~~~~\sr{(a)}{=}\sum\limits_{j=0}^{n+1}(-1)^j\sum\limits_{i=0}^{j-1}(-1)^i\sigma|_{[e_0,...,\wh{e}_i,...,\wh{e}_j,...,e_{n+1}]}+\sum\limits_{i=0}^{n+1}(-1)^i\sum\limits_{j=i+1}^{n+1}(-1)^{j-1}\sigma|_{[e_0,...,\wh{e}_i,...,\wh{e}_j,...,e_{n+1}]}\nn\\
&\textstyle~~~~=\sum\limits_{i,j=0}^{n+1}\Big[\chi_{\{0,...,j-1\}}(i)-\chi_{\{i+1,...,n+1\}}(j)\Big](-1)^{i+j}\sigma|_{[e_0,...,\wh{e}_i,...,\wh{e}_j,...,e_{n+1}]}\sr{(b)}{=}0,\nn
\end{align}
where at step (a) we swap $i$ and $j$, and step (b) holds because ~$\chi_{\{i+1,...,n+1\}}(j)=\chi_{\{0,...,j-1\}}(i)$.
\end{proof}

\begin{dfn}[\textcolor{blue}{\index{Singular homology}{Singular homology} of a space}]
Let $X$ be a space. The $n$th singular homology of $X$ is the homology ~~{\small$H_n(X):=H_n\big(C_\#(X)\big)={Z_n\big(C_\#(X)\big)\over B_n\big(C_\#(X)\big)}={Z_n(X)\over B_n(X)}$}.
\end{dfn}

\begin{rmk}[\textcolor{OliveGreen}{Interpretation of the homology of a space $X$}]
{\small\begin{align}
&H_0(X)={\ker\del_0\over\im\del_1}={C_0(X)\over\im\del_1}\cong{\txt{Span}_\Integer\{\txt{points in $X$}\}\over\txt{Span}_\Integer\{\gamma(1)-\gamma(0)~:~\txt{$\gamma$ a path in $X$}\}}
={\sum_{x\in X}\Integer x\over\sum_{x\sim_py}\Integer(x-y)}\nn\\
&~~~~={\textstyle\sum}_{x\in X}\Integer\ol{x},~~~~\txt{with}~~~~\ol{x}:=x+{\textstyle\sum}\big\{\Integer(x-y):x\sim_py~\txt{in $X$}\big\},\nn\\
&~~~~\cong \txt{Span}_\Integer\{\txt{path components of $X$}\},\nn
\end{align}}
where $x\sim_py$ iff there is a path between $x$ and $y$, iff ~$\Integer \ol{x}=\Integer\ol{y}$.
\begin{align}
&H_1(X)={\ker\del_1\over\im\del_2}\cong{\txt{Span}_\Integer\{\gamma:~\txt{$\gamma$ a loop in $X$}\}\over\txt{Span}_\Integer\{\gamma:~\txt{$\gamma$ a nullhomotopic loop in $X$}\}}\sr{(s)}{\cong}\txt{Abelianization of $\pi_1(X)$},\nn\\
&H_n(X)={\ker\del_n\over\im\del_{n+1}}\cong{\txt{Span}_\Integer\{f:S^n\ra X~|~\txt{$f$ a continuous map}\}\over\txt{Span}_\Integer\{f:S^n\ra X~|~\txt{$f$ a nullhomotopic continuous map}\}}\nn\\
&~~~~~~~~\cong\txt{Span}_\Integer\{\txt{certain $n$-dimensional holes in $X$}\},\nn
\end{align}
where ~$S^0:=\{x\in\Real:x^2=1\}=\{\pm1\}$ ~~(i.e., a two-point set), and the \index{Abelianization}{\ul{Abelianization}} of a group $G$ (with identity $e$) is the quotient group $G_{ab}:=G/[G,G]$, where
\begin{align}
[G,G]:=\big\langle [a,b]:a,b\in G\big\rangle\eqv\big\langle aba^{-1}b^{-1}:a,b\in G\big\rangle,\nn
\end{align}
is the normal subgroup of $G$ (called the \index{Commutator subgroup}{\ul{commutator subgroup}} of $G$) generated by elements of the form $[a,b]:=(ab)(ba)^{-1}=aba^{-1}b^{-1}$ called \index{Commutator}{\ul{commutators}}. For step (s), see \cite[p.166]{hatcher2001}.
\end{rmk}

\begin{dfn}[\textcolor{blue}{\index{Augmented singular chain complex}{Augmented singular chain complex} of a space $X$}]
The chain complex
\bea
\bt[row sep = tiny, column sep = small]
  & n=2& n=1& n=0& n=-1 &n=-2 &  \\
\wt{C}_\#(X):~~\cdots\ar[r,"\del_3"]&C_2(X)\ar[r,"\del_2"] & C_1(X)\ar[r,"\del_1"]&C_0(X)\ar[r,"\vep"]&\Integer\ar[r,"0"]& 0\ar[r,"0"]& \cdots,
\et\nn
\eea
where the map $\vep:C_0(X)\ra\Integer$ is given by ~$\vep\left(\sum n_i\sigma_i\right):=\sum n_i$, with $\sigma_i\in S_0(X)$.
\end{dfn}

\begin{dfn}[\textcolor{blue}{\index{Reduced singular homology}{Reduced singular homology} of a space $X$}]
The $n$th reduced singular homology of $X$ is the homology ~$\wt{H}_n(X):=H_n\big(\wt{C}_\#(X)\big)\subset H_n(X)$, where we can verify that
\bea
H_0(X)\cong\wt{H}_0(X)\oplus\Integer~~~~\txt{and}~~~~H_n(X)=\wt{H}_n(X)~~~~\txt{for all}~~n\geq 1.\nn
\eea
\end{dfn}

\begin{dfn}[\textcolor{blue}{Induced homomorphisms}]
Given a continuous map $f:X\ra Y$, we get homomorphisms $f_\#:C_n(X)\ra C_n(Y)$, $f_\#(\sigma):=f\circ\sigma$ (extended by linearity) that give a chain map $f_\#:C_\#(X)\ra C_\#(Y)$ because $f_\#\circ\del^{C_\#(X)}=\del^{C_\#(Y)}\circ f_\#$, and thus also give homomorphisms $f_\ast:H_n(X)\ra H_n(Y)$, $f_\ast([z]):=[f_\#(z)]=[f\circ z]$. Note that we have the immediate relations
\begin{align}
(f\circ g)_\#=f_\#\circ g_\#,~~(1_X)_\#=1_{C_n(X)},~~~~(f\circ g)_\ast=f_\ast\circ g_\ast,~~(1_X)_\ast=1_{H_n(X)}.\nn
\end{align}
\end{dfn}

\begin{lmm}[\textcolor{OliveGreen}{\cite[Theorem 2.10, p.111]{hatcher2001}}]
Let $f,g:X\ra Y$ be continuous maps. If $f\simeq g$, then $f_\#\simeq g_\#$, and so $f_\ast=g_\ast:H_n(X)\ra H_n(Y)$.
\end{lmm}
\begin{proof}
Let $\Delta^n=[e_0,...,e_n]\subset\Real^{n+1}$ be the standard $n$-simplex. Let $v_i=e_i\times\{0\}$ and $w_i=e_i\times\{1\}$, so that $[v_0,...,v_n]$ (resp. $[w_0,...,w_n]$) is the lower (resp. upper) face of $\Delta^n\times I$. Moreover, for each $0\leq i\leq n$, the $(n+1)$-simplex $[v_0,...,v_i,w_i,...,w_n]\subset\Delta^n\times I$ has $[v_0,...,v_i,w_{i+1},...,w_n]$ as ``lower face'' and $[v_0,...,v_{i-1},w_i,...,w_n]$ as ``upper face'', and $\Delta^n\times I=\bigcup_{i=0}^n [v_0,...,v_i,w_i,...,w_n]$. Let $F:X\times I\ra Y$ be a homotopy from $f$ to $g$. If $\sigma\in S_n(X)$, then based on the above simplex-decomposition of the solid cylinder $\Delta^n\times I$, and hence of its image $F\circ(\sigma\times 1_I):\Delta^n\times I\sr{\sigma\times 1_I}{\ral}X\times I\sr{F}{\ral}Y$, we can define a map $h:C_n(X)\ra C_{n+1}(Y)$ by ~$h(\sigma):=\sum_{i=0}^n(-1)^iF\circ(\sigma\times 1_I)|_{[v_0,...,v_i,w_i,...,w_n]}$. We have
{\small\begin{align}
&\textstyle\del h(\sigma)=\sum\limits_{i\geq j}(-1)^{i+j}F\circ(\sigma\times 1_I)|_{[v_0,...,\wh{v}_j,...,v_i,w_i,...,w_n]}+\sum\limits_{i\leq j}(-1)^{i+j+1}F\circ(\sigma\times 1_I)|_{[v_0,...,v_i,w_i,...,\wh{w}_j,...,w_n]}\nn\\
&\textstyle~~=\sum_{i=0}^n(-1)^{2i}F\circ(\sigma\times 1_I)|_{[v_0,...,\wh{v}_i,w_i,...,w_n]}+\sum_{i=0}^n(-1)^{2i+1}F\circ(\sigma\times 1_I)|_{[v_0,...,v_i,\wh{w}_i,...,w_n]}\nn\\
&\textstyle~~~+\sum\limits_{i>j}(-1)^{i+j}F\circ(\sigma\times 1_I)|_{[v_0,...,\wh{v}_j,...,v_i,w_i,...,w_n]}+\sum\limits_{i<j}(-1)^{i+j+1}F\circ(\sigma\times 1_I)|_{[v_0,...,v_i,w_i,...,\wh{w}_j,...,w_n]}\nn\\
&\textstyle~~=F\circ(\sigma\times 1_I)|_{[w_0,...,w_n]}-F\circ(\sigma\times 1_I)|_{[v_0,...,v_n]}-h(\del\sigma)=g\circ\sigma-f\circ\sigma-h(\del\sigma).\nn
\end{align}}Hence, $g_\#-f_\#=\del\circ h+h\circ\del$.
\end{proof}

\begin{crl}\label{HtyInvCrl1}
If spaces $X\simeq Y$, then $C_\#(X)\simeq C_\#(Y)$, and so $H_\ast(X)\cong H_\ast(Y)$.
\end{crl}

\begin{lmm}[\textcolor{OliveGreen}{Homology groups of a retract are subgroups: \cite[Sec. 2.1, Ex. 11]{hatcher2001}}]\label{HomRetInc}
If $A\subset X$ is a retract, the map $H_n(A)\ra H_n(X)$ induced by the inclusion $A \ha X$ is injective.
\end{lmm}
\begin{proof}
By hypotheses, we have a continuous map $X\sr{r}{\ral}A$ (along with the inclusion $A\sr{i}{\hal}X$) such that $r\circ i=1_A:A\sr{i}{\hal}X\sr{r}{\ral}A$. Thus, the induced homomorphisms
\begin{align}
&C_n(A)\sr{i_\#}{\ral}C_n(X),~\sigma\mapsto i\circ\sigma,~~~~H_n(A)\sr{i_\ast}{\ral}H_n(X),~[z]\mapsto [i_\#(z)]=[i\circ z],\nn\\
&C_n(X)\sr{r_\#}{\ral}C_n(A),~\sigma\mapsto r\circ\sigma,~~~~H_n(X)\sr{r_\ast}{\ral}H_n(A),~[z]\mapsto [r_\#(z)]=[r\circ z]\nn
\end{align}
satisfy $r_\#\circ i_\#=(r\circ i)_\#=(1_A)_\#=1_{C_n(A)}$, and $r_\ast\circ i_\ast=(r\circ i)_\ast=(1_A)_\ast=1_{H_n(A)}$.
\end{proof}

\begin{dfn}[\textcolor{blue}{Pairs of spaces, \index{Map of pairs}{Map of pairs}, \index{Homotopy! of pairs}{Homotopy of pairs}, Homotopic maps of pairs, Homotopy equivalence of pairs, Relative singular chain complex, \index{Relative singular homology}{Relative singular homology}, \index{Local singular homology}{Local singular homology}}]
A \ul{pair of spaces} is an ordered pair $(X,A)$ in which $X$ is a space and $A\subset X$. A \ul{map of pairs} $f:(X,A)\ra(Y,B)$ is a map $f:X\ra Y$ such that $f(A)\subset B$ (i.e., $f|_A:A\ra B$).

Two maps of pairs $f,g:(X,A)\ra (Y,B)$ are \ul{homotopic as maps of pairs} (written $f\simeq_{op} g$) if they are homotopic through a \ul{homotopy of pairs} $H:(X,A)\times I\ra(Y,B)$, defined to be a homotopy $H:X\times I\ra Y$ such that $H|_{X\times\{t\}}:(X,A)\ra(Y,B)$ is a map of pairs for each $t\in I$.

Two pairs $(X,A)$ and $(Y,B)$ are homotopy equivalent, written $(X,A)\simeq(Y,B)$, if there exists a \ul{homotopy equivalence of pairs} $f:(X,A)\ra(Y,B)$, defined to be homotopy equivalence $f:X\ra Y$ with homotopy inverse a map of pairs $g:(Y,B)\ra(X,A)$ such that $f\circ g\simeq_{op}1_Y$ and $g\circ f\simeq_{op} 1_X$.

Since the simplices $S_n(A)\subset S_n(X)$ for each $n$, and the inclusion $i:A\hookrightarrow X$ is injective, it induces an inclusion $i_\#:C_\#(A)\hookrightarrow C_\#(X)$, i.e., $C_\#(A)\subset C_\#(X)$, and so we have a short exact sequence (SES) of singular chain complexes $0\ra C_\#(A)\hookrightarrow C_\#(X)\sr{\pi}{\twoheadrightarrow} C_\#(X,A)\ra 0$, where $C_\#(X,A):=C_\#(X)/C_\#(A)$ is called the \ul{relative singular chain complex} of $(X,A)$. The \ul{relative singular homology} of $(X,A)$ is defined as
\bea
H_n(X,A):=H_n\big(C_\#(X,A)\big),\nn
\eea
which is (equivalently) generated/spanned by classes $c:=[z]\in H_n(X)$, of $n$-cycles $z\in Z_n(X)\subset C_\#(X)$, that are not contained in $A$ in the sense that
\bea
\textstyle z=\sum\limits_{\al=1}^m n_\al\sigma_\al^n~~\txt{for}~~\sigma_\al^n\in S_n(X)~~\txt{such that}~~\left(\bigcup\limits_\al\sigma_\al^n\right)\cap(X-A)\neq\emptyset.\nn
\eea
We also define the \ul{local singular homology} of $(X,A)$, or of $X$ at $A\subset X$, by
\bea
H_n(X|A):=H_n\big(C_\#(X|A)\big),~~~~C_\#(X|A):=C_\#(X,X-A),\nn
\eea
which is (equivalently) generated/spanned by classes $c:=[z]\in H_n(X)$, of $n$-cycles $z\in Z_n(X)\subset C_\#(X)$, that intersect $A$ in the sense that
\bea
\textstyle z=\sum\limits_{\al=1}^m n_\al\sigma_\al^n~~\txt{for}~~\sigma_\al^n\in S_n(X)~~\txt{such that}~~\left(\bigcup\limits_\al\sigma_\al^n\right)\cap A\neq\emptyset.\nn
\eea
\end{dfn}

Note that associated with any triple of spaces $(X,A,B)$, i.e., $X\supset A\supset B$, is a SES (short exact sequence) of relative singular chain complexes given by
\bea
\textstyle 0\ra C_\#(A,B)\hookrightarrow C_\#(X,B)\sr{\pi}{\twoheadrightarrow} {C_\#(X,B)\over C_\#(A,B)}\cong C_\#(X,A)\ra 0.\nn
\eea
Similarly, further associated with the triple of spaces $(X,X-B,X-A)$, i.e., $X\supset X-B\supset X-A$, is a SES of local singular chain complexes given by
\bea
\textstyle 0\ra C_\#(X-B|A)\hookrightarrow C_\#(X|A)\sr{\pi}{\twoheadrightarrow} {C_\#(X|A)\over C_\#(X-B|A)}\cong C_\#(X|B)\ra 0.\nn
\eea

\begin{lmm}[\textcolor{OliveGreen}{Homology groups of spheres: \cite[Cor. 2.14, p.114]{hatcher2001}}]\label{HlgGpSph}
{\small $\wt{H}_i(S^n)\cong
\left\{
  \begin{array}{ll}
   \Integer, & i=n \\
    0, & i\neq n
  \end{array}
\right.$}
\end{lmm}

\begin{thm}[\textcolor{OliveGreen}{No retraction theorem}]\label{NoRetOrig}
There exists no continuous retraction $r:D^n\ra S^{n-1}\cong\del D^n$, where $D^n\cong I^n$ is the $n$-dimensional closed disc.
\end{thm}
\begin{proof}
The inclusion $i:S^{n-1}\hookrightarrow D^n$ induces the map $i_\ast:\wt{H}_{n-1}(S^{n-1})=\Integer\ra \wt{H}_{n-1}(D^n)=0$, which is not injective.
\end{proof}

\begin{thm}[\textcolor{OliveGreen}{\index{Brouwer fixed point theorem}{Brouwer fixed point theorem}}]
Every continuous map $f:D^n\ra D^n$ ($D^n\cong I^n\cong\ol{B}_1(0)\subset\Real^n$) has a fixed point (a point $x$ such that $f(x)=x$).
\end{thm}
\begin{proof}
Suppose $f$ has no fixed point, i.e., $f(x)\neq x$ for all $x$. Define a map
$r:D^n\ra S^{n-1}=\del D^n$ by $r(x):=\del D^n\cap L[x,f(x)]$, with
$L[x,f(x)]$ the line through $x$ and $f(x)$ in $\Real^n$.

\begin{figure}[H]
\centering
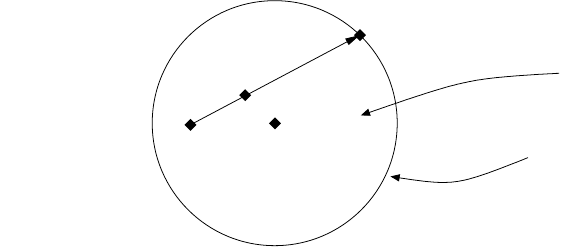
\end{figure}
{\flushleft For} an explicit formula, note that $r(x)=x+t(x-f(x))$ for $t\geq 0$, and $\|r(x)\|=1$. Thus,
\begin{align}
&1=\|r(x)\|=\|x\|^2+2t\langle x,x-f(x)\rangle+t^2\|x-f(x)\|^2,\nn\\
&~~\Ra~~\|x-f(x)\|^2t^2+2\langle x,x-f(x)\rangle t+\|x\|^2-1=0,\nn\\
&\textstyle~~\Ra~~t(x)={-\langle x,x-f(x)\rangle+\sqrt{(\langle x,x-f(x)\rangle)^2-\|x-f(x)\|^2(\|x\|^2-1)}\over\|x-f(x)\|^2},\nn\\
&~~\Ra~~r(x)=x+t(x)(x-f(x)),\nn
\end{align}
which shows $r$ is continuous as a composition of continuous maps, and $r(x)=x$ if $\|x\|=1$ (i.e., $r$ is a retraction). This contradicts the ``No retraction theorem'' (Theorem \ref{NoRetOrig}).
\end{proof}

\section{Topological Extension Theorems: Absolute Retracts}\label{PrelimsTET}  
This section is not strictly essential to our main purpose and serves only to introduce related terminology and research questions. The main results of interest are (1) extendability of a continues map from a closed subset of a normal Hausdorff space to $\Real^n$ in Corollary \ref{TietExThmCrl} of Tietze's extension theorem, (2) Paracompactness of metric spaces in Corollary \ref{MetricParac}, and (3) the absolute retraction property of convex subsets of locally convex spaces in Corollary \ref{DT-theoremCrl} of Dugundji's generalization of Tietze's extension theorem.

\subsection{Urysohn's lemma and Tietze's extension theorem}
\begin{dfn}[\textcolor{blue}{\index{Separating function}{Separating function}}] Let $X$ be a space and $A,B\subset X$. A separating function of $A,B$ is a continuous function $f:X\ra[0,1]\subset\Real$ such that $f|_A=0$ and $f|_B=1$.
\end{dfn}

\begin{dfn}[\textcolor{blue}{\index{Completely regular space}{Completely regular space}}] A space $X$ is completely regular if (i) $X$ is a $T_1$ space and, (ii) any closed set $C\subset X$ and any point $x_0\not\in C$ have a separating function.
\end{dfn}

\begin{thm}[\textcolor{OliveGreen}{\index{Urysohn's lemma}{Urysohn's lemma}: \cite[Theorem 4.1, p.146]{dugundji1966}}]\label{UrysThm}
Let $X$ be a Hausdorff space. Then $X$ is normal $\iff$ every two disjoint closed sets have a separating function.
\end{thm}
\begin{proof}
($\Ra$): Assume $X$ is normal. Let $C_0,C_1\subset X$ be disjoint closed sets. By normality of $X$, there is an open set $O\left(\sfrac{1}{2}\right)$ such that $C_0\subset O\left(\sfrac{1}{2}\right)\subset\ol{O\left(\sfrac{1}{2}\right)}\subset C_1^c$. Repeating this, we get two more open sets $O\left(\sfrac{1}{2^2}\right)$ and $O\left(\sfrac{3}{2^2}\right)$, in addition to $O\left(\sfrac{1}{2}\right)=O\left(\sfrac{2}{2^2}\right)$, such that
\bea
C_0\subset O\left(\tfrac{1}{2^2}\right)\subset\ol{O\left(\tfrac{1}{2^2}\right)}\subset O\left(\tfrac{2}{2^2}\right)\subset\ol{O\left(\tfrac{2}{2^2}\right)}\subset O\left(\tfrac{3}{2^2}\right)\subset\ol{O\left(\tfrac{3}{2^2}\right)}\subset C_1^c.\nn
\eea
Continuing this way, for each $n>0$ we get open sets $\left\{O\left(k\over 2^n\right):1\leq k\leq 2^n-1\right\}$ such that
\bea
C_0\subset O\left(\tfrac{1}{2^n}\right)\subset O\left(\tfrac{2}{2^n}\right)\subset O\left(\tfrac{3}{2^n}\right)\subset\cdots\subset O\left(\tfrac{2^n-1}{2^n}\right)\subset C_1^c,~~~~\ol{O\left(\tfrac{k-1}{2^n}\right)}\subset O\left(\tfrac{k}{2^n}\right).\nn
\eea\
Let ~$D:=\bigcup_{n\geq1}\left\{0={0\over 2^n},{1\over 2^n},{3\over 2^n},\cdots,{2^n-1\over 2^n},{2^n\over 2^n}=1\right\}\subset[0,1]$. Then $\ol{O(a)}\subset O(b)$ for all $a,b\in D\backslash\{0,1\}$ such that $a<b$. Let $d(x):=\inf\{d:x\in O(d)\}=\inf\{d_n(x):n\geq 1\}$, where $d_n$ is determined by $O\big(d_n(x)\big)=\bigcap\left\{O(\tfrac{k}{2^n}):x\in O(\tfrac{k}{2^n}),~0\leq k\leq 2^n\right\}$ due to the observation that $d\leq d'$ $\iff$ $O(d)\subset O(d')$. Consider the map ~$f:X\ra [0,1]$~ given by
\bit
\item[] \centering{$f(x):=d(x)\chi_{C_1^c}(x)+\chi_{C_1}(x)$, ~~ where ~~
{\footnotesize $\chi_A(x):=
\left\{
  \begin{array}{ll}
    1, & x\in A \\
    0, & x\not\in A
  \end{array}
\right\}
$}.}
\eit
Then $f|_{C_0}=\inf(D)=0$ and $f|_{C_1}=1$. It remains to show $f$ is continuous. That is, for any $x\in X$, given any $N_\vep(f(x))$, there is $N_\delta(x)$ such that $f\big(N_\delta(x)\big)\subset N_\vep(f(x))$. Observe that due to the ordering of the sets $\{O(d):d\in D\}$, we have
\bit[leftmargin=0cm]
\item[]\centering {\small $O(d)\cap O(d')
=\left\{
   \begin{array}{ll}
     O(d), & O(d)\subset O(d') \\
     O(d'), & O(d)\supset O(d')
   \end{array}
 \right\}$, ~~
$O(d)\backslash\ol{O(d')}
=\left\{
   \begin{array}{ll}
     \neq\emptyset, & O(d)\supset O(d') \\
     \emptyset, & O(d)\subsetneq O(d')
   \end{array}
 \right\}$.}
\eit
Thus, for any fixed number $d,d'\in D$, we have
\begin{align}
& f|_{O(d')}(x)=\inf_{x\in O(r)\cap O(d')}r=\min\left(d'~,~\inf_{x\in O(r)\subseteq O(d')}r\right)~\leq~d',\nn\\
& f|_{X\backslash\ol{O(d')}}(x)\geq\inf_{x\in O(r)\cap\left[X\backslash\ol{O(d')}\right]}r=\inf_{x\in O(r)\backslash\ol{O(d')}}r  =\inf_{x\in O(r)\supseteq O(d')}r~\geq~d',\nn\\
& f|_{O(d)\backslash\ol{O(d')}}(x)=\inf_{x\in O(r)\cap\left[O(d)\backslash\ol{O(d')}\right]}r=\inf_{\substack{x\in O(r)\\ O(d)\supseteq O(r)\supseteq O(d')}}r.\nn
\end{align}
That is, for any fixed numbers $d,d'\in D$, we have $0~\leq~f|_{O(d)}\leq d\leq f|_{X\backslash\ol{O(d)}}\leq 1$, and, if $d<d'$ then $d'\leq f|_{O(d)\backslash\ol{O(d')}}\leq d$.
Equivalently, for any fixed numbers $d,d'\in D$, we have
\bit[leftmargin=0cm]
\item[]\centering $f\big(O(d)\big)\subset[0,d)$, $f\left(X\backslash\ol{O(d)}\right)\subseteq(d,1]$, and, if $d<d'$ then $f\left(O(d)\backslash\ol{O(d')}\right)\subseteq[d',d]$.
\eit
Fix $x\in X$ and any $1>\vep>0$. Then we have the following cases.
\bit[leftmargin=1cm]
\item $f(x)=0$: Choose $d$ such that  $0=f(x)<d<\vep$. Then
\bea
x\in C_0\subset O(d)~~~~\txt{and}~~~~f\big(O(d)\big)\subseteq[0,d)\subset [0,\vep).\nn
\eea
\item $f(x)=1$: Choose $d$ such that $1-\vep<d<f(x)=1$. Then
\bea
x\in C_1\subset X\backslash\ol{O(d)}~~~~\txt{and}~~~~f\left(X\backslash\ol{O(d)}\right)\subseteq(d,1]\subset (1-\vep,1].\nn
\eea
\item $0<f(x)<1$: Choose $d,d'$ such that $0<f(x)-\vep<d<f(x)<d'<f(x)+\vep$. Then
\bea
x\in O(d')\backslash\ol{O(d)}~~~~\txt{and}~~~~f\left(O(d')\backslash\ol{O(d)}\right)\subseteq(d,d')\subset\Big(f(x)-\vep,f(x)+\vep\Big).\nn
\eea
\eit
{\flushleft ($\La$):}  Assume any two disjoint closed sets have a separating function. Let $C_0,C_1$ be disjoint closed sets in $X$, and let $f:X\ra [0,1]$ be their separating function. Then we get disjoint open neighborhoods $U:=f^{-1}\Big([0,1/2)\Big)\supset C_0$ and $V:=f^{-1}\Big((1/2,1]\Big)\supset C_1$, where $1/2$ can be replaced with any other number $0<a<1$.
\end{proof}

\begin{thm}[\textcolor{OliveGreen}{\index{Tietze's extension theorem}{Tietze's extension theorem}: \cite[Theorem 5.1, p.149]{dugundji1966}}]\label{TietzeThm}
Let $X$ be a Hausdorff space. Then (i) $X$ is normal $\iff$ every continuous function on a closed set $f:C\subset X\ra\Real$ extends to a continuous function $F:X\ra\Real$. Moreover, (ii) if $f$ is bounded as $|f|\leq c$, then $F$ is also bounded as $|F|\leq c$.
\end{thm}
\begin{proof}
($\Ra$): Assume $X$ is normal. Let $C\subset X$ be closed and $f:C\ra\Real$ continuous.
\bit[leftmargin=0cm]
\item[]\ul{\emph{Case 1 ($f$ bounded)}}: Let $c:=\sup_{x\in C}|f(x)|$. Then $f(X)\subset[-c,c]$. Consider the closed sets
\bea
\textstyle A_0:=f^{-1}\left(\left[-c,-{c\over 3}\right]\right),~~~~B_0:=f^{-1}\left(\left[{c\over 3},c\right]\right).\nn
\eea
Let ~$g_0={2c\over 3}\left(u_0-{1\over 2}\right):X\ra\left[-{c\over 3},{c\over 3}\right]$,~ where $u_0$ is the separating function of $A_0$ and $B_0$. Next, let $f_1=f-g_0|_C:C\ra\left[-{2c\over 3},{2c\over 3}\right]$, and consider the closed sets
\bea
\textstyle A_1:=f_1^{-1}\left(\left[-{2c\over 3},-{2c\over 9}\right]\right),~~~~B_1:=f_1^{-1}\left(\left[{2c\over 9},{2c\over 3}\right]\right).\nn
\eea
As before, let $g_1={4 c\over 9}\left(u_1-{1\over 2}\right):X\ra\left[-{2c\over 9},{2c\over 9}\right]$, where $u_1$ is the separating function of $A_1$ and $B_1$. Next, let $f_2=f_1-g_1|_C:C\ra\left[-{4c\over 9},{4c\over 9}\right]$, and consider the closed sets
\bea
\textstyle A_2:=f_2^{-1}\left(\left[-{4c\over 9},-{4c\over 27}\right]\right),~~~~B_2:=f_2^{-1}\left(\left[{4c\over 27},{4c\over 9}\right]\right).\nn
\eea
Continuing this way, for each $n\geq 1$, we have cont. functions $g_n:X\ra\left[-{2^nc\over 3^{n+1}},{2^nc\over 3^{n+1}}\right]$ and
\bea
\textstyle f_{n+1}=f-(g_0+g_1+\cdots+g_n)|_C:C\ral\left[-{2^{n+1}c\over 3^{n+1}},{2^{n+1}c\over 3^{n+1}}\right].\nn
\eea
Let $F_n=g_0+g_1+\cdots+g_n:X\ra\Real$ for $n\geq 1$. Then for $n<m$,
{\small\bea
\textstyle |F_n(x)-F_m(x)|\leq\sum\limits_{k=n+1}^m|g_k(x)|\leq\sum\limits_{k=n+1}^m{2^kc\over 3^{k+1}}\leq{c\over 3}\sum\limits_{k=n+1}^\infty\left(2\over 3\right)^k={c\over 3}{\left(2\over 3\right)^{n+1}\over 1-{2\over 3}}=\left(2\over 3\right)^{n+1}c.\nn
\eea}That is, $\{F_n\}\subset\F(X,\Real)$ is Cauchy, and so converges to a continuous function $F:X\ra\Real$,
\bea
\textstyle |F(x)|=\left|\sum\limits_{k=0}^ng_k(x)\right|\leq\sum\limits_{k=0}^n|g_k(x)|\leq \sum\limits_{k=0}^n{2^kc\over 3^{k+1}}\leq{c\over 3}\sum\limits_{k=0}^\infty\left({2\over 3}\right)^k=c.\nn
\eea
It is clear that $F|_C=f$, since $f_{n+1}:=f-F_n|_C\ra0$ uniformly (i.e., in $\F(C,\Real)$).

\item[]\ul{\emph{Case 2 ($f$ unbounded)}}: Let $h:\Real\ra(-1,1)$ be a homeomorphism. Then $h\circ f:C\sr{f}{\ral}\Real\sr{h}{\ral}(-1,1)$ is bounded and so extends to a continuous function $G:X\ra(-1,1)$. Thus, $f$ extends to the continuous function $F=h^{-1}\circ G:X\sr{G}{\ral}(-1,1)\sr{h^{-1}}{\ral}\Real$.
\eit

($\La$): Assume that for each closed set $C\subset X$, any continuous function $f:C\ra\Real$ extends to a continuous function $F:X\ra\Real$. Let $A,B\subset X$ be disjoint closed sets. Consider the function $f:A\cup B\ra\Real$ such that $f|_A=0$ and $f|_B=1$. Then $f$ is continuous, and so extends to a continuous function $F:X\ra\Real$. Thus, with any $0<\vep<1/2$, the sets $F^{-1}\big((-\vep,\vep)\big)$ and $F^{-1}\big((1-\vep,1+\vep)\big)$ are disjoint neighborhoods of $A$ and $B$.
\end{proof}

\begin{crl}\label{TietExThmCrl}
Let $X$ be a normal Hausdorff space. Then every continuous map on a closed set $f:C\subset X\ra\Real^n$ extends to a continuous map $F:X\ra \Real^n$.
\end{crl}
\begin{proof}
Let $f:C\subset X\ra\Real^n$ be continuous and $C$ closed. Then $f(x)$ = $(f_1(x)$,$\cdots$,$f_n(x))$, where each $f_i:C\subset X\ra\Real$ is continuous, and so extends to a continuous function $F_i:X\ra\Real$. Hence $f$ extends to the continuous map $F:X\ra\Real^n$, $F(x):=\big(F_1(x),\cdots,F_n(x)\big)$.
\end{proof}

\subsection{Paracompactness and full normality}
\begin{dfn}[\textcolor{blue}{Star of a set wrt a cover}]
Let $X$ be a space, $A\subset X$, and $\U$ a cover of $X$. The $\U$-star of $A$ is the union of all elements of $\U$ that intersect $A$, i.e.,
\bea
\textstyle \U\ast A=Star_\U(A):=\bigcup\{U\in\U:A\cap U\neq\emptyset\}\subseteq X.\nn
\eea
\end{dfn}

\begin{dfn}[\textcolor{blue}{\index{Star of a cover}{Star of a cover}}]
Let $X$ be a space $X$ and $\U$ a cover of $X$. Then
\bea
\textstyle\U_\ast=\U_{star}:=\left\{\U\ast U:U\in\U\right\},~~~~(\txt{a cover of $X$}).\nn
\eea
\end{dfn}

\begin{dfn}[\textcolor{blue}{\index{Locally! finite collection}{Locally finite collection}}]
Let $X$ be a space and $\A\subset\P(X)$ some sets. Then $\A$ is locally finite if each $x\in X$ has an open neighborhood $O_x\ni x$ that intersects at most finitely many $A\in\A$.
\end{dfn}

\begin{dfn}[\textcolor{blue}{\index{Refinement}{Refinement} of a cover}] Let $X$ be a space and $\U$, $\V$ two covers of $X$. Then  $\U$ is a refinement of $\V$, written $\U\leq\V$, if every element of $\U$ is contained in some element of $\V$ (i.e., for every $U\in\U$ there is $V\in\V$ such that $U\subset V$).
\end{dfn}

\begin{dfn}[\textcolor{blue}{\index{Paracompact space}{(Para)compact space}}] A space $X$ is (para)compact if every open cover of $X$ has a (locally) finite open refinement. (Note: Unlike in certain literature, regularity, or the Hausdorff property, is not part of our definition.)
\end{dfn}

\begin{lmm}
If $X$ is a (para)compact space, then so is every closed subspace $C\subset X$.
\end{lmm}
\begin{proof}
Let $\U\subset\P(C)$ be an open cover of $C$. For each $U\in\U$, let $\wt{U}$ be an open set in $X$ such that $U=C\cap\wt{U}$. Then $\V:=\{X-C\}\cup\{\wt{U}:U\in\U\}$ is an open cover of $X$. Let $\V_1$ be a (locally) finite open refinement of $\V$. We get a (locally) finite open refinement of $\U$, as
\bit
\item[]\centering $\U_1=C\cap \V:=\{C\cap V_1:V_1\in\V_1\}$. \qedhere
\eit
\end{proof}

\begin{lmm}[\textcolor{OliveGreen}{\cite[Lemma 39.1]{munkres}}]\label{NSMTlmm}
Let $X$ be a space. If $\A\subset\P(X)$ is locally finite, then (a) $\C=\left\{\ol{A}:A\in\A\right\}$ is locally finite and (b) {\small $\ol{\bigcup_{A\in\A}A}=\bigcup_{A\in\A}\ol{A}$}.
\end{lmm}
\begin{proof}
(a) By the closure criterion, an open set intersects $A$ $\iff$ it intersects $\ol{A}$.

{\flushleft (b)} Let $Z:=\bigcup_{A\in\A}A$. Since each $A\subset Z$, we have $\bigcup_{A\in\A}\ol{A}\subset\ol{Z}$. On the other hand, if $x\in\ol{Z}$, let $N(x)$ be a neighborhood of $x$ intersecting only some $A_1,...,A_n\in\A$ (equivalently, intersecting only $\ol{A}_1,...,\ol{A}_n\in\C$). Then {\small $N(x)\cap\ol{Z}\subset \ol{\bigcup_{i=1}^nA_i}=\bigcup_{i=1}^n\ol{A}_i$}, and so $x\in\bigcup_{A\in\A}\ol{A}$.
\end{proof}

\begin{lmm}[\textcolor{OliveGreen}{\cite[Theorem 41.1, p.253]{munkres}}]\label{PaHauNorm}
A (para)compact Hausdorff space is normal.
\end{lmm}
\begin{proof}
Let $X$ be a paracompact Hausdorff space. First we will prove that $X$ is regular. Let $x_0\in X$ and $C\subset X$ a closed set not containing $x_0$. Since $X$ is Hausdorff, each point $c\in C$ has a neighborhood $U_c$ disjoint from a neighborhood $V_c$ of $x_0$, and so $x_0\not\in \ol{U}_c$ (since $U_c\subset(V_c)^c$ $\Ra$ $\ol{U}_c\subset (V_c)^c$). Thus, we get an open cover of $X$ given by
\bea
\U:=\{C^c\}\cup\{U_c:c\in C\}.\nn
\eea
Since $X$ is paracompact, $\U$ has a locally finite open refinement $\V\leq\U$. Let
\bea
\D=\V|_C:=\{V\in\V:V\cap C\neq\emptyset\}=\{\txt{elements of $\V$ that meet $C$}\}.\nn
\eea
Observe that $x_0\not\in D$ for all $D\in\D$ (since $D\not\subset C^c$ and so $D\subset U_c$ for some $c\in C$, where $x_0\not\in\ol{U}_c$). It follows that for all $D\in\D$, we have $x_0\not\in \ol{D}$ (since $\ol{D}\subset\ol{U}_c$).
{\flushleft Let} $W:=\bigcup_{D\in\D}D$, an open set containing $C$. Since $\D$ is locally finite, Lemma \ref{NSMTlmm} implies $\ol{W}=\bigcup_{D\in\D}\ol{D}$, which shows $x_0\not\in\ol{W}$. Hence, $\ol{W}^c$ is a neighborhood of $x_0$ disjoint from $W$ (a neighborhood of $C$). This proves $X$ is regular.

To prove normality, we repeat (below) the same argument above with $x_0$ replaced by another closed set $C_0$ and the Hausdorff property replaced by regularity.
{\flushleft Let} $C_0,C\subset X$ be disjoint closed sets. Since $X$ is regular, each point $c\in C$ has a neighborhood $U_c$ disjoint from a neighborhood $V_c$ of $C_0$, and so $C_0\cap\ol{U}_c=\emptyset$ (since $U_c\subset (V_c)^c$ $\Ra$ $\ol{U}_c\subset (V_c)^c$). Thus, we get an open cover of $X$ given by
\bea
\U:=\{C^c\}\cup\{U_c:c\in C\}.\nn
\eea
Since $X$ is paracompact, $\U$ has a locally finite open refinement $\V\leq\U$. Let
\bea
\D=\V|_C:=\{V\in\V:V\cap C\neq\emptyset\}=\{\txt{elements of $\V$ that meet $C$}\}.\nn
\eea
Observe that $C_0\cap D=\emptyset$ for all $D\in\D$ (since $D\not\subset C^c$ and so $D\subset U_c$ for some $c\in C$, where $C_0\cap\ol{U}_c=\emptyset$). It follows that for all $D\in\D$, we have $C_0\cap \ol{D}=\emptyset$ (since $\ol{D}\subset\ol{U}_c$).
{\flushleft Let} $W:=\bigcup_{D\in\D}D$, an open set containing $C$. Since $\D$ is locally finite, Lemma \ref{NSMTlmm} implies $\ol{W}=\bigcup_{D\in\D}\ol{D}$ which shows $C_0\cap\ol{W}=\emptyset$. Hence, $\ol{W}^c$ is a neighborhood of $C_0$ disjoint from $W$ (a neighborhood of $C$). This proves $X$ is normal.
\end{proof}

\begin{dfn}[\textcolor{blue}{\index{Partition of unity}{Partition of unity}}] Let $X$ be a space. A partition of unity on $X$ is a collection of continuous functions $\big\{f_\al:X\ra[0,1]\big\}_\al$ such that
\bit[leftmargin=1cm]
\item[(i)] Each $x\in X$ has a neighborhood $V\ni x$ such that $f_\al|_V=0$ except for finitely many $\al$.
\item[(ii)] Each $x\in X$ satisfies $\sum\limits_\al f_\al(x)=1$. (The sum is well defined everywhere by (i).)
\eit
If $\U$ is an open cover of $X$, a partition of unity {\footnotesize $\P=\big\{f_\al:X\ra[0,1]\big\}_\al$} is \ul{subordinate to} (or \ul{dominated by}) $\U$ if for each $\al$, we have {\footnotesize $Supp(f_\al):=\ol{\{x\in X:f_\al(x)\neq 0\}}\subset U$} for some $U\in\U$.
\end{dfn}

\begin{lmm}[]\label{PaHauCov}
If $X$ is a paracompact Hausdorff space, every locally finite open cover $\{U_\al\}$ has a locally finite open refinement $\{V_\al\}\leq \{U_\al\}$ such that $\ol{V_\al}\subset U_\al$ for each $\al$.
\end{lmm}
\begin{proof}
Given $x\in X$, we know $x\in U$ for some $U\in\U:=\{U_\al\}$. Let a neighborhood $N(x)\subset U$ intersect only finitely many elements of $\U$. Since $X$ is a normal $T_1$ space, we can choose $N(x)$ to be disjoint from some neighborhood $N(U^c)$, and so
\bea
~~\Ra~~N(x)\subset N(U^c)^c~~\Ra~~\ol{N(x)}\subset N(U^c)^c\subset U.\nn
\eea
Thus, $\O$ $:=$ $\{$open $O\subset X:$ $O$ meets only finitely many $U_\al\in\U$, and $\ol{O}\subset U$ for some $U\in\U$$\}$ is an open refinement of $\U$. Since, by paracompactness, $\O$ has a locally finite open refinement with the same properties as $\O$ (i.e., that each element intersects only finitely many $U\in\U$ and has closure lying in some $U\in\U$), we can assume $\O$ is a locally finite open refinement of $\U$. Define another locally finite open refinement $\V=\{V_\al\}$ of $\U=\{U_\al\}$ by
\bea
\textstyle V_\al:=\bigcup\{O\in\O:~\ol{O}\subset U_\al\}~\subset~ U_\al.\nn
\eea
Note that $\V$ is a cover of $X$ because it is refined by $\O$ (i.e., every $O\in\O$ is contained in some $V_\al$), and $\V$ is locally finite because $\U$ is locally finite. We will now show that $\ol{V_\al}\subset U_\al$. Let $x\not\in U_\al$. Choose a neighborhood $N(x)$ that intersects only finitely many of the $O\subset V_\al$, say $O_1,...,O_n\subset V_\al$. Then $V_\al-(O_1\cup\cdots\cup O_n)\subset X-N(x)$, and so
\bea
&&V_\al=(O_1\cup\cdots\cup O_n)\cup\left[V_\al-(O_1\cup\cdots\cup O_n)\right]\subset \left(\ol{O_1}\cup\cdots\cup\ol{O_n}\right)\cup\big(X-N(x)\big),\nn\\
&&~~\Ra~~\ol{V_\al}\subset\ol{\left(\ol{O_1}\cup\cdots\cup\ol{O_n}\right)\cup\big(X-N(x)\big)}=\ol{O_1}\cup\cdots\cup\ol{O_n}\cup\big(X-N(x)\big).\nn
\eea
Since $\ol{O_i}\subset U_\al$, we have $x\not\in\ol{O_i}$, and since $x\not\in(X-N(x)\big)$, it follows that $x\not\in \ol{V_\al}$. That is, $X-U_\al\subset X-\ol{V_\al}$, and so $\ol{V_\al}\subset U_\al$.
\end{proof}

\begin{crl}[\textcolor{OliveGreen}{Existence of a partition of unity: \cite{michael1953}}]
A paracompact Hausdorff space admits a partition of unity.
\end{crl}
\begin{proof}
Let $X$ be a paracompact Hausdorff space and $\U$ an open cover of $X$. By paracompactness, let $\V$ be a locally finite open refinement of $\U$. Then by Lemma \ref{PaHauCov}, $\V$ has a locally finite open refinement $\W=\{W_V:V\in\V\}\leq\V$ such that $\ol{W_V}\subset V$ for all $V\in \V$.
Thus, we have a locally finite closed cover (i.e., a cover by closed sets) of $X$ given by
\bea
\C=\left\{C_V=\ol{W_V}:V\in \V\right\},~~\txt{such that}~~C_V\subset V~~\txt{for all}~~V\in\V.\nn
\eea
Since $C_V$ and $X-V$ are disjoint closed sets, we have a separating function $f_V:X\ra [0,1]$ such that $f_V|_{C_V}=1$ and $f_V|_{X-V}=0$. Thus, with the continuous function
\bea
\textstyle f:X\ra(0,\infty),~~~~f(x):=\sum_{V\in\V}f_V(x),\nn
\eea
we get a partition of unity on $X$ subordinate to $\V$ (and hence also subordinate to $\U$),
\bea
\textstyle \Phi_V:X\ra[0,1],~~~~\phi_V(x):={f_V(x)\over f(x)}={f_V(x)\over\sum_{V'\in\V}f_{V'}(x)},~~~~V\in\V.\nn
\eea
(Note that $f$ is well defined because for each $x$, $f_V(x)=0$ except for finitely many $V\in\V$. Recall that every $x\in X$ lies in finitely many $V\in\V$, since $\V$ is locally finite.)
\end{proof}

\begin{dfn}[\textcolor{blue}{\index{Star-refinement}{Star-refinement} of a cover}]
Let $X$ be a space and $\U$, $\V$ covers of $X$. Then $\U$ is a star-refinement of $\V$, written $\U\leq_\ast\V$, if $\U_\ast$ is a refinement of $\V$ (i.e., $\U_\ast\leq\V$).
\end{dfn}

\begin{dfn}[\textcolor{blue}{\index{Fully normal space}{Fully normal space}}] A space $X$ if fully normal if every open cover of $X$ has an open star-refinement.
\end{dfn}
\begin{lmm}
A fully normal space is a normal space.
\end{lmm}
\begin{proof}
Let $X$ be fully normal, and let $C,C'\subset X$ be disjoint closed sets. Then $\U:=\{X\backslash C,X\backslash C'\}$ is an open cover of $X$. Let $\O=\{O_\al\}$ be an open star-refinement of $\U$. Let
\bit[leftmargin=0cm]
\item[]\centering $O:=\O\ast C=\bigcup\{O_\al\in\O:O_\al\cap C\neq\emptyset\}$, ~ $O':=\O\ast C'=\bigcup\{O_\al\in\O:O_\al\cap C'\neq\emptyset\}$.
\eit
Let $O_\al\subset O$ (i.e., $O_\al\cap C\neq\emptyset$) and $O_{\al'}\subset O'$ (i.e., $O_{\al'}\cap C'\neq\emptyset$). If $O_\al\cap O_{\al'}\neq\emptyset$, then $O_\al\cup O_{\al'}\subset \O\ast O_\al$ lies in $X\backslash C$ or in $X\backslash C'$ (a contradiction). Hence, $O,O'$ are disjoint open neighborhoods of $C$ and $C'$.
\end{proof}

\begin{lmm}[\textcolor{OliveGreen}{Paracompactness and Full Normality: \cite[Theorems 1,2]{stone1948}}]\label{ParacFullnorm}
A space $X$ is paracompact and Hausdorff $\iff$ fully normal and $T_1$.
\end{lmm}
\begin{proof}
($\Ra$): Let $X$ be a paracompact Hausdorff space. Since $X$ is normal (Lemma \ref{PaHauNorm}), it follows that $X$ is $T_1$. Thus, it remains to prove that every open cover of $X$ has an open star refinement. Let $\U=\{U_\al\}_{\al\in A}$ be an open cover of $X$. Note that if $\V\leq\U$ is any open cover of $X$, then an open star refinement of $\V$ is also an open star refinement of $\U$. Thus, it suffices (by paracompactness of $X$) to assume $\U$ is locally finite. By Lemma \ref{PaHauCov}, there are open sets $\{X_\al\}_{\al\in A}$ such that $\ol{X}_\al\subset U_\al$ and $\bigcup X_\al=X$.

Since $\U$ is locally finite, each $x\in X$ has a neighborhood $V(x)$ that meets only finitely many $U_\al$, say $\{U_\al:\al\in A_x\}$ for a finite set $A_x:=\{\al\in A:V(x)\cap U_\al\neq\emptyset\}\subset A$. Let
\begin{align}
&B_x:=\{\al\in A_x:x\in U_\al\}=\{\al\in A:x\in V(x)\cap U_\al\},\nn\\
&C_x:=\left\{\al\in A_x:x\not\in\ol{X_\al}\right\}=\left\{\al\in A:x\in X-\ol{X_\al},~V(x)\cap U_\al\neq\emptyset\right\}.\nn
\end{align}
Then it is clear that $A_x\subset B_x\cup C_x$, and so $A_x=B_x\cup C_x$, where $B_x\cap C_x$ = $\{\al\in A:x\in V(x)\cap(U_\al-\ol{X_\al})\}$. We have a neighborhood of $x$ given by
\begin{align}
&\textstyle W(x):=V(x)\cap\bigcap\{U_\al:\al\in B_x\}\cap\bigcap\{X-\ol{X_\al}:\al\in C_x\}\nn\\
&\textstyle~~~~=\bigcap\left\{V(x)\cap[U_\al-\ol{X}_{\al'}]:\al\in B_x,\al'\in C_x\right\}\nn\\
&\textstyle~~~~=\bigcap\left\{V(x)\cap[U_\al-\ol{X}_{\al'}]:x\in[U_\al-\ol{X}_{\al'}],~\al,\al'\in A_x\right\}\nn
\end{align}
Thus, the sets $\W:=\{W(x):x\in X\}$ form an open cover of $X$. To show $\W$ is a star refinement of $\U$, it suffices (since each $W(x)\subset X_{\beta_x}$ for some $\beta_x\in A$) to show that for any $y\in X$, if $y\in X_\beta$ then $\W\ast y:=\bigcup\{W(x)\in\W:y\in W(x)\}\subset U_\beta$. So, let $y\in X$, and fix $\beta\in A$ such that $y\in X_\beta$. Pick any $x\in X$ such that $y\in W(x)$, i.e., $W(x)\subset\W\ast y$. Then
{\small\begin{align}
&W(x)\cap \ol{X_\beta}\neq\emptyset~~\txt{and}~~W(x)\subset V(x)~~\Ra~~V(x)\cap \ol{X_\beta}\neq\emptyset~~\Ra~~V(x)\cap U_\beta\neq\emptyset~~\Ra~~\beta\in A_x\nn\\
&~~\txt{and}~~\beta\not\in C_x~~~~[\txt{otherwise, if $\beta\in C_x$, then $W(x)\subset V(x)\cap \ol{X}_\beta\cap(X-\ol{X}_\beta)=\emptyset$, a contradiction}],\nn\\
&~~\Ra~~\beta\in B_x,~~\Ra~~x\in U_\beta,~~\Ra~~W(x)\subset U_\beta~~~~\txt{(by the definition of $W(x)$)}.\nn
\end{align}}
{\flushleft ($\La$)}:
Let $X$ be a fully normal $T_1$ space. Then we know $X$ is normal and $T_1$, and so Hausdorff. It remains to prove paracompactness. Let $\U=\{U_\al\}$ be an open cover of $X$. We will construct a locally finite open refinement of $\U$. By full normality, we have a sequence of open covers $\U\geq_\ast\U^1\geq_\ast\U^2\geq_\ast\cdots$ (i.e., $\U^n$ is star-refined by $\U^{n+1}$). For $A\subset X$, let
\begin{align}
\label{stone1}&\textstyle A_n:=\U^n\ast A=\bigcup_{x\in A}x_n=\bigcup_{x\in A}\U^n\ast x=\bigcup\{U^n\in\U^n:U^n\cap A\neq\emptyset\},\\
\label{stone2}&\textstyle A_{-n}:=X-(X-A)_n=X-\U^n\ast(X-A).
\end{align}
Note that $A_{-n}$ is a closed set. Since $X-A\subset (X-A)_n$ and $A\subset A_n$, we have $A_{-n}\subset A\subset A_n$. Moreover, by definition, $x\in A_{-n}=X-(X-A)_n$ $\iff$ $x$ is not in any element of $\U^n$ that meets $X-A$, $\iff$ any element of $\U^n$ that contains $x$ lies in $A$. Therefore,
\begin{align}
\label{stone3}&\textstyle A_{-n}=\left\{x\in A~|~x_n=\U^n\ast x\subset A\right\},
\end{align}
\begin{align}
\label{stone4}&\textstyle~~\Ra~~(A_{-n})_n=\bigcup_{x\in A_{-n}}x_n~\subset~ A.
\end{align}
Since $\U^n\geq\U_\ast^{n+1}$ and {\small $(A_{n+1})_{n+1}=\U^{n+1}\ast A_{n+1}=\bigcup_{U^{n+1}\cap A\neq\emptyset}\U^{n+1}\ast U^{n+1}$} is a union of $\U^{n+1}$-stars of those elements of $\U^{n+1}$ that meet $A$, it follows that
\bea
\label{stone5}\textstyle(A_{n+1})_{n+1}=\bigcup_{U^{n+1}\cap A\neq\emptyset}\U^{n+1}\ast U^{n+1}~\subset~ \bigcup_{U^n\cap A\neq\emptyset}U^n~=~A_n.
\eea
We also have the following:
\begin{enumerate}[leftmargin=0.7cm]
\item \ul{$A\subset B$ ~$\Ra$~ $A_n\subset B_n$ and $A_{-n}\subset B_{-n}$}: Because $A\subset B$ implies $X-B\subset X-A$, which implies $(X-B)_n\subset (X-A)_n$, which implies $X-(X-A)_n\subset X-(X-B)_n$.
\item \ul{$n\leq m$ ~$\Ra$~ $A_n\supset A_m$ and $A_{-n}\subset A_{-m}$}: Because, $A_{n+1}\subset (A_{n+1})_{n+1}\subset A_n$ and $(X-A)_{n+1}\subset (X-A)_n$ implies $A_{-n}\subset A_{-(n+1)}$.
\item \ul{$\ol{A}~\subset~A_n$}: Indeed, if $x\in\ol{A}-A_n$, then some element of $\U^n$ contains $x$ but does not meet $A$, which is a contradiction since $\ol{A}=\{x\in X:\txt{every open}~Nbd(x)\cap A\neq\emptyset\}$.
\item \ul{$x\in y_n$ $\iff$ $x$ and $y$ both lie in some $U^n\in\U^n$ $\iff$ $y\in x_n$.} (\textbf{Note}: This implies the open sets $x_n=\U^n\ast x$, for all $x\in X$, form a partition of $X$).
\end{enumerate}
{\flushleft For} each $\al$, let
\begin{align}
& V_\al^1:=(U_\al)_{-1}=X-(X-U_\al)_1=\{x\in U_\al:~x_1\subset U_\al\}\subset U_\al,~~~~(\txt{\textbf{Note}: $V_\al^1$ is closed}),\nn\\
& V_\al^2:=(V_\al^1)_2=\U^2\ast V_\al^1,\nn\\
& V_\al^3:=(V_\al^2)_3=\U^3\ast V_\al^2,\nn\\
& \cdots\hspace{1cm}\cdots\hspace{1cm}\cdots\nn\\
\label{stone6}& V_\al^n:=(V_\al^{n-1})_n=\U^n\ast V_\al^{n-1}.
\end{align}
Since $V_\al^1\subset V_\al^2\subset V_\al^3\subset\cdots$ and $V_\al^n$ is open for $n\geq 2$, it follows that $V_\al:=\bigcup_nV_\al^n$ is open. Also, by induction using (\ref{stone4}) and (\ref{stone5}), we have
\bea
&&V_\al^n\subset(V_\al^n)_n=\left(\left(V_\al^{n-1}\right)_n\right)_n\sr{(\ref{stone5})}{\subset} \left(V_\al^{n-1}\right)_{n-1}\sr{(\ref{stone5})}{\subset}\cdots\sr{(\ref{stone5})}{\subset}(V^1_\al)_1=\big((U_\al)_{-1}\big)_1\sr{(\ref{stone4})}{\subset}U_\al,\nn\\
\label{stone7}&&\textstyle~~\Ra~~V_\al:=\bigcup_n V_\al^n\subset U_\al.
\eea
For each $x\in X$, we have $x_1\subset (U^1)_1=\U^1\ast U^1\subset U_\al$ for some $U^1$ in $\U^1$ and some $U_\al$ in $\U$ (since $\U^1$ star-refines $\U$), and so by (\ref{stone3}), $x\in (U_\al)_{-1}=V_\al^1\subset V_\al$. Thus, $\{V_\al\}$ covers $X$, i.e.,
\bea
\label{stone8}\textstyle\bigcup_\al V_\al=X. ~~~~\txt{(So $\V:=\{V_\al\}$ is an open refinement of $\U=\{U_\al\}$)}
\eea
If $x\in V_\al$, then we know $x\in V_\al^{n-1}$ for some $n$, and so $x_n\subset(V_\al^{n-1})_n=V_\al^n\subset V_\al$, i.e.,
\bea
\label{stone9} x\in V_\al~~\Ra~~x_n\subset V_\al~~~~\txt{for some}~~n\geq 1.
\eea

For simplicity, assume (without loss of generality) that the indices $\{\al\}$ are a well-ordered set, i.e., totally ordered and each subset has a least element or infimum. Recursively define a sequence of closed sets as follows.
{\small\begin{align}
\label{stone10}&\textstyle C_{n1}:=(V_1)_{-n}=X-(X-V_1)_n=\{x\in X:x_n\subset V_1\},\nn\\
&\textstyle C_{n\al}:=\left(V_\al-\bigcup_{\beta<\al}C_{n\beta}\right)_{-n}=X-\left(X-\Big(V_\al-\bigcup_{\beta<\al}C_{n\beta}\Big)\right)_n=X-\left((X-V_\al)\cup\bigcup_{\beta<\al}C_{n\beta}\right)_n\nn\\
&\textstyle ~~=\left\{x\in X:~x_n\subset V_\al-\bigcup_{\beta<\al}C_{n\beta}\right\}\subset V_\al.
\end{align}} These sets, $\{C_{n\al}\}$, have the following (partitioning) property.
\bea
\label{stone11}\txt{If $\al\neq\gamma$, no $U^n\in\U^n$ can meet both $C_{n\al}$ and $C_{n\gamma}$}.
\eea
[[Indeed, assuming wlog $\gamma<\al$, if $U^n\in\U^n$ meets $C_{n\al}$, i.e., some $x\in U^n\cap C_{n\al}$, then by (\ref{stone3}) and (\ref{stone10}), $U^n\subset x_n\subset V_\al-\bigcup_{\beta<\al}C_{n\beta}$, which shows $U^n$ is disconnected from $C_{n\gamma}$.]]

The sets $\C:=\{C_{n\al}\}$ cover $X$, i.e.,
\bea
\label{stone12}\textstyle \bigcup_{n,\al}C_{n\al}=X.
\eea
[[Indeed, given $x\in X$, by (\ref{stone8}), there is a first (i.e., a smallest) $\al$ such that $x\in V_\al$. By (\ref{stone9}), $x_n\subset V_\al$ for some $n$. If $x\not\in C_{n\al}=\left(V_\al-\bigcup_{\beta<\al}C_{n\beta}\right)_{-n}$, then by (\ref{stone3}), $x_n\not\subset V_\al-\bigcup_{\beta<\al}C_{n\beta}$, i.e., some $y\in x_n-\left(V_\al-\bigcup_{\beta<\al}C_{n\beta}\right)=x_n\cap \bigcup_{\beta<\al}C_{n\beta}$. This implies $y\in x_n\cap C_{n\beta}$ for some $\beta<\al$, and so because $y\in x_n$ $\iff$ $x\in y_n$, we have
\bea
\textstyle x\in y_n\subset (C_{n\beta})_n=\left(\left(V_\beta-\bigcup_{\beta'<\beta}C_{n\beta'}\right)_{-n}\right)_n\sr{(\ref{stone4})}{\subset} V_\beta-\bigcup_{\beta'<\beta}C_{n\beta'}\subset V_\beta,\nn
\eea
which contradicts the minimality of $\al$.]]

Next, consider the open covers (open refinements of $\U$)
\bea
\label{stone13} P_{n\al}:=(C_{n\al})_{n+3},~~~~Q_{n\al}:=(C_{n\al})_{n+2},
\eea
which satisfy ~$C_{n\al}\subset P_{n\al}\subset\ol{P_{n\al}}\subset(P_{n\al})_{n+3}\subset Q_{n\al}=\left(\left(V_\al-\bigcup_{\beta<\al}C_{n\beta}\right)_{-n}\right)_{n+2}\subset V_\al$.
{\flushleft The} open cover $\{Q_{n\al}\}$ also has the partitioning property (\ref{stone11}) wrt $\U^{n+2}$, i.e.,
\bea
\label{stone14}\txt{If $\al\neq\gamma$, no $U^{n+2}\in\U^{n+2}$ can meet both $Q_{n\al}$ and $Q_{n\gamma}$}.
\eea
[[Indeed $Q_{n\al}=\bigcup\left\{U^{n+2}:~U^{n+2}\cap C_{n\al}\neq\emptyset\right\}$ and $Q_{n\gamma}=\bigcup\left\{U^{n+2}:~U^{n+2}\cap C_{n\gamma}\neq\emptyset\right\}$. But each $U^{n+2}$ lies in some $U^n$, and by (\ref{stone11}) no $U^n$ can meet both $C_{n\al}$ and $C_{n\gamma}$. Thus, no $U^{n+2}$ can meet both $C_{n\al}$ and $C_{n\gamma}$, and hence no $U^{n+2}$ can meet both $Q_{n\al}$ and $Q_{n\gamma}$.]]
{\flushleft Similarly}, the open cover $\{P_{n\al}\}$ has the partitioning property (\ref{stone11}) with respect to $\U^{n+3}$.

The sets $C_n:=\bigcup_\al\ol{P_{n\al}}$ are closed. [[Indeed, if $x\in \ol{C_n}$, then because $\{P_{n\al}\}$ covers $X$, every neighborhood $N(x)$ meets some $P_{n\al}=(C_{n\al})_{n+3}$. Choose $N(x)$ small enough so that $N(x)\subset U^{n+3}\subset x_{n+3}$ for some $U^{n+3}\in\U^{n+3}$. Then by property (\ref{stone14}) for the cover $\{P_{n\al}\}$, $N(x)$ can meet at most one $P_{n\al}$, and so $x\in \ol{P_{n\al}}\subset C_n$.]]

Now define the following open sets (the desired locally finite open refinement of $\U$).
\bea
\textstyle \O:~~O_{1\al}=Q_{1\al},~~~~O_{n\al}=Q_{n\al}-\bigcup_{j=1}^{n-1}C_j~~~~~~~~\big(\subset Q_{n\al}\subset V_\al\subset U_\al\big).\nn
\eea
These sets cover $X$, i.e., $\bigcup_{n,\al}O_{n\al}=X$. [[To see this, let $x\in X$. Then by (\ref{stone12}), $x$ lies in some $C_{n\al}\subset\ol{P_{n\al}}\subset Q_{n\al}$. Let $m$ be the smallest integer such that $x\in\ol{P_{m\beta}}\subset Q_{m\beta}$ for some $\beta$ (i.e., $x\not\in \ol{P_{n\al}}$ for all $n<m$ and all $\al$, equivalently, $x\not\in C_1,...,C_n$ for all $n<m$). Then $x\in Q_{m\beta}-\bigcup_{j=1}^{m-1}C_j=O_{m\beta}$.]]

It remains to show that $\O=\{\O_{n\al}\}$ is locally finite. [[Let $x\in X$. Then $x\in C_{n\al}$ for some $(\al,n)$. Thus, $x_{n+3}\subset (C_{n\al})_{n+3}=P_{n\al}\subset C_n$ is disjoint from $O_{k\beta}=Q_{k\beta}-\bigcup_{j=1}^{k-1}C_j$ for $k>n$, since $O_{n+1~\beta}=Q_{n+1~\beta}-\bigcup_{j=1}^nC_j$, $O_{n+2~\beta}=Q_{n+2~\beta}-\bigcup_{j=1}^{n+1}C_j$, $\cdots$. On the other hand, for $k\leq n$, we have $x_{n+3}\subset U^{k+2}$ for some $U^{k+2}\in\U^{k+2}$. Thus, by (\ref{stone14}), $x_{n+3}$ can meet
\bea
\textstyle O_{k\beta}=Q_{k\beta}-\bigcup_{j=1}^{k-1}C_j,~~~~1\leq k\leq n,\nn
\eea
only for one value of $\beta$. Hence, $x_{n+3}$ meets at most $n$ of the sets in $\O$.]]
\end{proof}

\begin{crl}[\textcolor{OliveGreen}{\cite[Theorem 1.2, p.13]{howes1995}}]\label{MetricParac}
Every metric space is paracompact.
\end{crl}
\begin{proof}
Let $X$ be a metric space (hence $T_1$). We need to show $X$ is fully normal. Let $\U=\{U_\al\}$ be an open cover of $X$. For each $x\in X$, there exists $0<\vep(x)<1$ such that $B_{4\vep(x)}(x)\subset U_{\al(x)}$ for some $\al=\al(x)$. So, we have an open cover of $X$,
\bea
\V=\left\{B_{\vep(x)}(x):x\in X\right\},~~~~B_{\vep(x)}(x)\subset B_{4\vep(x)}(x)\subset U_{\al(x)}.\nn
\eea
We will show that $\V$ is a star-refinement of $\U$. Given $y\in X$, define
\begin{align}
&{y^\ast}:=\left\{x\in X:~y\in B_{\vep(x)}(x)\right\}=\{x\in X:d(x,y)<\vep(x)\},\nn\\
&\textstyle~\Ra~\V\ast y:=\bigcup\left\{B_{\vep(x)}(x):y\in B_{\vep(x)}(x)\right\}=\bigcup\left\{B_{\vep(x)}(x):d(x,y)<\vep(x)\right\}=\bigcup_{x\in{y^\ast}}B_{\vep(x)}(x).\nn
\end{align}
For any $x',z\in{y^\ast}$ (i.e., $d(x',y)<\vep(x')$, $d(z,y)<\vep(z)$),
\begin{align}
& B_{\vep(x')}(x')=\{u:d(u,x')<\vep(x')\}\subset \{u:d(u,y)\leq d(u,x')+d(x',y)<2\vep(x')\}= B_{2\vep(x')}(y)\nn\\
&~~~~\subset \{u:d(u,z)\leq d(u,y)+d(y,z)<2\vep(x')+\vep(z)\}=B_{2\vep(x')+\vep(z)}(z)\nn\\
&\textstyle~~~~\subset B_{4\vep(z)}(z),~~~~\txt{if we choose $z\in{y^\ast}$ such that}~~\vep(z)>{2\over 3}\vep(x'),\nn\\
&\textstyle~~\Ra~~B_{\vep(x')}(x')\subset B_{4\vep(z)}(z)\subset U_{\al(z)},~~\txt{for $x',z\in y^\ast$ such that}~~\vep(z)>{2\over 3}\vep(x').\nn
\end{align}
Thus, if we choose $z_y\in{y^\ast}$ such that $\vep(z_y)>{2\over 3}\sup\{\vep(x):x\in{y^\ast}\}$, then
\begin{align}
B_{\vep(x')}(x')\subset B_{4\vep(z_y)}(z_y)\subset U_{\al(z_y)}~~\txt{for all}~~x'\in{y^\ast},~~\Ra~~\V\ast y\subset B_{4\vep(z_y)}(z_y)\subset U_{\al(z_y)}.\nn
\end{align}
Thus, $\V$ is a star refinement of $\U$, because for any $x\in X$,
\bea
\textstyle\V\ast B_{\vep(x)}(x)=\bigcup_{y\in B_{\vep(x)}(x)}\V\ast y\sr{(s)}{=}\V\ast y\subset U_{\al(z_y)},~~~~\txt{for any}~~y\in B_{\vep(x)}(x),\nn
\eea
where step (s) is a special property of the selected cover $\V$ of $X$ by balls.
\end{proof}

\subsection{Dugundji's extension of Tietze's theorem}~\\\vspace{-0.2cm}\\
Some useful facts on the topology of cellular complexes can be found in \cite[p.519]{hatcher2001}. Let $X$ be a space and $\sigma:\Delta^n\ra X$ a singular $n$-simplex (i.e., a continuous map). We introduce the \ul{unoriented-boundary} of a simplex as follows: Denoting the union of all faces of $\Delta^n$ of dimension $\leq n-1$ by $\del\Delta^n:=\bigcup\{\txt{$k$-faces of}~\Delta^n~|~k\leq n-1\}$, we define
\bea
\delta\sigma\eqv\delta\sigma(\Delta^n):=\sigma(\Delta^n)-\sigma(\Delta^n-\del\Delta^n).\nn
\eea
Note that the discussion above and that to follow carries through unchanged (up to homeomorphism) if we replace the affine $n$-simplex $\Delta^n$ with the unit $n$-disc $D^n\cong[0,1]^n$.

\begin{dfn}[\textcolor{blue}{\index{Cell}{$n$-cell} (a type of singular $n$-simplex) in a space $X$, \index{Characteristic map}{Characteristic map} for a cell, \index{Delta cell}{$\Delta$-cell}}]
An \ul{$n$-cell} in $X$ is a subspace $e^n\subset X$ with a Hausdorff closure $\ol{e}\subset X$ and a continuous map $f:\Delta^n\ra X$ (called the \ul{characteristic map} for $e^n$) such that:
\bit[leftmargin=1cm]
\item[(i)] $f(\Delta^n)=\ol{e}^n$ (the closure of $e^n$)~ and ~$f(\Delta^n-\del\Delta^n)=e^n$.
\item[(ii)] $f|_{\Delta^n-\del\Delta^n}:\Delta^n-\del\Delta^n\ral X$ is a homeomorphism onto $e^n$. (Note that for constructions requiring greater generality, this condition may be dropped, or modified accordingly.)
\eit
A \ul{$\Delta$-cell} in $X$ is a cell $e\subset X$ whose characteristic map $f:\Delta^{\dim e}\ra X$ is injective on $\del\Delta^{\dim e}$ (i.e., $f|_{\del\Delta^{\dim e}}:\Delta^{\dim e}\ra X$ is injective).
\end{dfn}
Note that (i) and (ii) imply that $f(\del\Delta^n)=\delta e^n:=\ol{e}^n-e^n$ ( = the unoriented-boundary of $f$ as a singular simplex in $X$; a \ul{non-topological} notion of boundary in the space $X$), since
\bea
f^{-1}(\delta e^n):=f^{-1}(\ol{e}^n-e^n)=f^{-1}(\ol{e}^n)-f^{-1}(e^n)=\Delta^n-(\Delta^n-\del\Delta^n)=\del\Delta^n.\nn
\eea
Moreover, for any cell $e\subset X$, the characteristic map $f:\Delta^{\dim e}\ra X$ is a \ul{quotient map} $\Delta^{\dim e}\ra\ol{e}$ (since $f$ is closed, as a continuous map from a compact space to a Hausdorff space). Also, as a compact Hausdorff space, $\ol{e}=f(\Delta^{\dim e})\subset X$ is a normal subspace.

\begin{dfn}[\textcolor{blue}{\index{Affine! cell}{Affine $n$-cell} in a vector space $X$}]
An $n$-cell $e^n\subset X$ with an affine characteristic map $f:\Delta^n\ra X$. (Here, $f$ is called affine if $f$ admits a linear extension {\small $F:\Span_\Real(\Delta^n)\subset\Real^{n+1}\ra X$, $F(\ld x+y)=\ld F(x)+F(y)$ for $x,y\in\Span_\Real(\Delta^n)$} and scalar $\ld$.)
\end{dfn}
Note that if $e^n$ is an affine $n$-cell, then $\ol{e}^n$ is an affine $n$-simplex.
\begin{dfn}[\textcolor{blue}{Faces of an affine $n$-cell $e^n$ in a vector space $X$}]
The faces of $\ol{e}^n$ (or of $e^n$) are the images of the faces of $\Delta^n$ under the characteristic map $f:\Delta^n\ra X$ of $e^n$.
\end{dfn}

\begin{dfn}[\textcolor{blue}{\index{Cellular complex}{Cellular complex} (\cite[p.221]{whitehead1949}), \index{Skeleton}{$n$-skeleton}, \index{Delta complex}{$\Delta$-complex}}]
A \ul{cellular complex} is a space {\small $K=(K,\C)$} together with a collection of cells {\small$\C=\C_K$} with the following properties.
\bit[leftmargin=0.9cm]
\item[(i)] The cells $\C$ are disjoint, and cover $K$, i.e., $K=\bigsqcup\C$.
\item[(ii)] For each $n$-cell $e^n\in\C$ with characteristic map $f:\Delta^n\ra K$, the image $f\big(\del\Delta^n\big)=\delta e^n\subset K$ lies in the union of cells in $\C$ of dimension $\leq n-1$.
\eit
That is, if we write {\small $K=\bigcup_{n\geq0}K^n$}, where {\small$K^n:=\bigcup\{e^k\in\C:~k\leq n\}$} = union of all cells in $\C$ of dimension $\leq n$, then for each $n$-cell $e^n\in\C$, we have {\small$f\big(\del\Delta^n\big)=\delta e^n\subset K^{n-1}\subset K^n$}, where
\bea
\textstyle K^n~=~K^{n-1}~\sqcup~\bigsqcup_{\al\in \Gamma_n}e_\al^n,~~~~\txt{for $n$-cells ~$e_\al^n\in\C$}~~~~(\txt{\ul{$n$-skeleton} of $K$}).\nn
\eea

A cellular complex $(K,\C)$ is called a \ul{$\Delta$-complex} if the following hold:
\bit[leftmargin=0.9cm]
\item[(a)] Every cell $e\in\C$ is a $\Delta$-cell (i.e., a cell whose characteristic map $f:\Delta^{\dim e}\ra X$ is injective on $\del\Delta^{\dim e}$).
\item[(b)] For each $n$-cell $e^n\in\C$ with characteristic map $f:\Delta^n\ra K$, each face image
    \bea
    f\big([e_0,...,\widehat{e}_i,...,e_n]\big)\subset K,~~~~i=0,...,n,\nn
    \eea
    is an $(n-1)$-cell in $\C$.
\eit
\end{dfn}
{\flushleft Note} that in our definition above, unlike \cite{whitehead1949}, we do not require the space $K$ to be Hausdorff. Any space $X$ has at least one cell complex structure, namely, as a complex with all of its points as 0-cells (and no higher-dimensional cells), i.e., $X=X^0$. Also, not necessarily all cells in a cellular complex $K=(K,\C)$ are open (as subsets of $K$).

Observe that by construction, for any two cells $e_1^n,e_2^m\in\C$, if $n\leq m$ then $\ol{e}_1^n\cap e_2^m=\emptyset$. In particular, any two $n$-cells $e_1^n,e_2^n\in\C$ are \ul{separated} in the sense that $\ol{e}_1^n\cap e_2^n=e_1^n\cap\ol{e}_2^n=\emptyset$.

\begin{dfn}[\textcolor{blue}{\index{Subcomplex}{Subcomplex}}] Given a complex $K=(K,\C)$, a subset $L\subset K$ is a subcomplex if $L=\bigcup\C_L$ for some subcollection of cells $\C_L\subset \C$ such that $e\in\C_L$ $\Ra$ $\ol{e}\subset L$. (Note that each $K^n$ is a finite-dimensional subcomplex of $K$. $K^n$ is a subcomplex because for any cell $e\in K^n$, we have $\ol{e}=e\sqcup\delta e\subset K^n\cup K^{n-1}=K^n$.)
\end{dfn}

\begin{dfn}[\textcolor{blue}{\index{Finite complex}{Finite complex}}] A complex $K=(K,\C)$ such that $\C$ is a finite collection.
\end{dfn}

\begin{dfn}[\textcolor{blue}{\index{Locally! finite complex}{Locally-finite complex}}] A complex $K=(K,\C)$ such that each point $x\in K$ has a neighborhood that meets only finitely many cells in $\C$. Equivalently, each point $x\in K$ lies in the interior of a finite subscomplex (i.e., for each $x\in K$, there is a finite subcomplex $L\subset K$ with $x\in L^o=\txt{Int}(L)$).
\end{dfn}

\begin{dfn}[\textcolor{blue}{\index{Generated subcomplex}{Generated subcomplex}}] In a complex $K=(K,\C)$, the subcomplex generated by $A\subset K$ is {\footnotesize $K(A):=\bigcap\{L:A\subset L\subset K,\txt{$L$ a subcomplex}\}$}, i.e., the smallest subcomplex containing $A$. (Note that for any $e\in \C$, we have {\small $K(e)=K(\ol{e})=K(x)$} for each $x\in e$.)
\end{dfn}

\begin{dfn}[\textcolor{blue}{\index{Closure-finite complex}{Closure-finite complex}}] A complex $K=(K,\C)$ such that for each $e\in \C$, $\ol{e}$ intersects only finitely many cells in $\C$. Equivalently, $K=(K,\C)$ is closure-finite if the subcomplex $K(e)$ is finite for each cell $e\in\C$ (i.e., $K(x)$ is finite for each point $x\in K$).
\end{dfn}

\begin{dfn}[\textcolor{blue}{\index{Weak complex}{Weak complex}}] A complex $K=(K,\C)$ whose topology (called \ul{weak topology}) is such that $A\subset X$ is closed (open) $\iff$ $A\cap\ol{e}$ is closed (relatively open) for each $e\in\C$. (See an alternative/equivalent in the footnote\footnote{Equivalently, $K$ is a weak complex if $A\subset X$ is closed (open) $\iff$ $A\cap K^n$ is closed (relatively open) for each $n\geq 0$. To prove this, recall that each $K^n$ is a subcomplex, and use (i) induction on $n$ based on  the decomposition $K^n~=~K^{n-1}~\sqcup~\bigsqcup_{\al\in \Gamma_n}e_\al^n=K^{n-1}~\cup~\bigcup_{\al\in \Gamma_n}\ol{e}_\al^n$, for $n$-cells $e_\al^n\in\C$, and (ii) the fact that for each $n$-cell $e^n\in\C$, we have $A\cap\ol{e}^n=[A\cap K^n]\cap\ol{e}^n$ with $\ol{e}^n\subset K^n$ a subspace.}.)
\end{dfn}

Note that every weak complex $(K,\C)$ is automatically a $T_1$ space, since for any $x\in K$ and $e\in\C$, $\ol{e}\cap\{x\}$ is a closed subset of $\ol{e}$ (because $\ol{e}$ is $T_1$) and hence a closed subset of $K$.

\begin{lmm}
Every weak complex is a Hausdorff space.
\end{lmm}
\begin{proof}
Let $K=\bigcup_{n\geq 0} K^n=(K,\C)$ be a  weak complex. Recall that $K$ is a $T_1$ space (i.e., points are closed in $K$) by the definition of the weak topology. Let $\C_n\subset\C$ denote the cell decomposition of the \ul{weak subcomplex} $K^n=(K^n,\C_n)=\bigcup_{0\leq r\leq n}K^r$. Observe that each $n$-cell $e^n\in\C$ is open in $K^n$, because for each $0\leq r\leq n$, we have either $e^n\cap K^r=\emptyset$ (if $r\leq n-1$) or $e^n\cap K^r=e^n$ (if $r=n$). It follows that $K^0\subset K$ is a discrete (hence Hausdorff) subspace. By induction on $n$ based on the decomposition $K^n~=~K^{n-1}~\sqcup~\bigsqcup_{\al\in \Gamma_n}e_\al^n$, for $n$-cells $e_\al^n\in\C$, we see that each $K^n\subset K$ is a Hausdorff subspace as follows: Let $x,x'\in K^n$. Then $x\in\ol{e}_{\al}^n$, $x'\in \ol{e}_{\al'}^n$ for some $\al,\al'\in\Gamma_n$. Consider the following cases.
\bit[leftmargin=0.6cm]
\item $x\in e_{\al}^n$ or $x'\in e_{\al'}^n$: Then it is clear that $x$ and $x'$ have disjoint neighborhoods in $\ol{e}_\al^n$ and in $\ol{e}_{\al'}^n$ (hence in $K^n$).
\item $x,x'\in\delta e_{\al}^n\cap\delta e_{\al'}^n\subset K^{n-1}$: Then by induction on $n$, these points have disjoint neighborhoods $Nb(x),Nb(x')\subset K^{n-1}$, and so have disjoint neighborhoods in $K^n$.
\eit
Let $N_n(x)$, $N_n(x')$ be the disjoint neighborhoods in $K^n$, chosen such that $N_n(x)\subset N_{n+1}(x)$, $N_n(x')\subset N_{n+1}(x')$, say by replacing $N_n(x)$ with $N_n(x)\cap N_{n+1}(x)$ and $N_n(x')$ with $N_n(x')\cap N_{n+1}(x')$. Then $N(x):=\bigcup N_n(x)$, $N(x'):=\bigcup N_n(x')$ are disjoint nbds in $K$, because
\bea
\textstyle N(x)\cap N(x')=\big(\bigcup N_n(x)\big)\cap\big(\bigcup N_n(x')\big)=\bigcup\limits_{n,n'}\big(N_n(x)\cap N_{n'}(x')\big)=\bigcup\limits_{n,n'}\emptyset=\emptyset.\nn
\eea
\end{proof}

\begin{lmm}\label{ClFnWknss}
Let $K=(K,\C)$ be a closure-finite complex. Then $K$ is weak $\iff$ for any $A\subset K$, ``$A\cap L$ closed for every finite subcomplex $L\subset K$'' implies ``$A$ is closed''.
\end{lmm}
\begin{proof}
($\Ra$) Assume $K$ is weak. If $A\cap L$ is closed for every finite subcomplex $L$, then for each $e\in\C$, $A\cap\ol{e}=[A\cap K(e)]\cap\ol{e}$ is closed. Hence $A$ is closed since $K$ is weak.

($\La$): Conversely, assume $A'\subset K$ is closed whenever $A'\cap L$ is closed for every finite subcomplex $L\subset K$. Let $A\subset K$. To show $K$ is weak, we need to show $A$ is closed iff $A\cap\ol{e}$ is closed for all $e\in\C$. If $A$ is closed, it is clear that $A\cap\ol{e}$ is closed for all $e\in\C$. On the other hand, if $A\cap\ol{e}$ is closed for every $e\in\C$, then because each finite complex $L$ is a finite union of cells of the type $\ol{e}$, it follows that $A\cap L$ is closed for each finite subcomplex $L$ (and hence $A$ is closed).
\end{proof}

\begin{dfn}[\textcolor{blue}{\index{CW-complex}{CW-complex} (\textcolor{blue}{\cite{whitehead1949}, p.223})}]
A complex $K=(K,\C)$ that is both closure-finite and weak.
\end{dfn}
By \cite[Proposition A.3, p.522]{hatcher2001}, a CW-complex is a normal space (a fact that is not required in the following discussion).

\begin{lmm}\label{ClFnWknss2}
Every (locally-) finite complex is a CW-complex.
\end{lmm}
\begin{proof}
(i) \ul{A finite complex $K=(K,\C)$ is CW}: $K$ is trivially both closure-finite and weak since $K=\bigcup_{e\in\C}e=\bigcup_{e\in\C}\ol{e}$ is a finite union of closed sets. (ii) \ul{A locally-finite complex $K=(K,\C)$ is CW}: For each $x\in K$, $K(x)$ is finite since $K$ is locally-finite, and so $K$ is closure-finite. Let $A\subset K$. To show $K$ is weak, it suffices by Lemma \ref{ClFnWknss} to show that ``$A\cap L$ closed for every finite subcomplex $L$'' implies ``$A$ closed''. So, assume $A\cap L$ is closed for every finite subcomplex $L$. For any $x\in int(K-A)$, let $L_x$ be a finite subcomplex such that $x\in int(L_x)$. Since $A\cap L_x$ is closed, $x\in int(L_x-A)=int(L_x-A\cap L_x)\subset K-A$, and so $A$ is closed (since $K-A$ is open).
\end{proof}

\begin{crl}
A space is a (locally) finite complex $\iff$ it is a (locally) compact CW-complex.
\end{crl}

\begin{lmm}\label{ClFnWknss3}
Let $(K,\C)$ be a CW complex. For any space $Y$, a map $f:K\ra Y$ is continuous $\iff$ $f|_{\ol{e}}:\ol{e}\subset K\ra Y$ is continuous for each $e\in\C$.
\end{lmm}
\begin{proof}
If $f$ is continuous, it is clear that $f|_{\ol{e}}$ is continuous for each $e\in\C$. Conversely, if $f|_{\ol{e}}$ is continuous for each $e\in\C$, then given a closed set $C\subset Y$, the set $f^{-1}(C)$ is closed in $K$ since {\footnotesize $f^{-1}(C)\cap\ol{e}=(f|_{\ol{e}})^{-1}(C)$} is closed in $\ol{e}$ for each $e\in\C$, and so $f$ is continuous.
\end{proof}

\begin{dfn}[\textcolor{blue}{\index{Affine! delta complex}{Affine $\Delta$-complex} in a vector space $X$}]
A $\Delta$-complex $A=(A,\C)\subset X$ whose cells $\C$ are affine cells. Note that an affine cellular complex (by construction) has the following defining properties.
\bit[leftmargin=0.9cm]
\item[(i)] The cells $\C$ are disjoint and cover $A$, i.e., $A=\bigsqcup\C$.
\item[(ii)] If $e\in\C$, then every face of $\ol{e}$ is also in $\C$.
\eit
\end{dfn}

\begin{dfn}[\textcolor{blue}{\index{Polytope}{Polytope} (or \index{Simplicial complex}{Simplicial complex}) in a vector space $X$}]
An affine cellular complex $P=(P,\C)\subset X$ in which the intersection of two closed cells in $\C$ is also a cell in $\C$. Thus, a polytope (by construction) has the following defining properties.
\bit[leftmargin=0.9cm]
\item[(i)] The cells $\C$ are disjoint and cover $P$, i.e., $P=\bigsqcup\C$.
\item[(ii)] If $e\in\C$, then every face of $\ol{e}$ is also in $\C$.
\item[(iii)] For all $e,e'\in\C$, the intersection $\ol{e}\cap\ol{e}'$ is a face of both $\ol{e}$ and $\ol{e}'$.
\eit
\end{dfn}

\begin{dfn}[\textcolor{blue}{\index{Triangulated space}{Triangulated space}}]
A space that is homeomorphic to a polytope.
\end{dfn}

\begin{dfn}[\textcolor{blue}{Star of a polytope cell}]
Let $(P,\C)$ be a polytope. For any $e\in\C$, let
\begin{align}
P\ast e=Star_P(e):={\textstyle\bigcup}\big\{\text{open cells}~e_0\in\C~\txt{with $e$ as a face of $\ol{e_0}$}\big\}.\nn
\end{align}
\end{dfn}

\begin{dfn}[\textcolor{blue}{\index{CW-polytope}{CW-polytope}: \cite[p.254]{dugundji1951}}]
A polytope $P=(P,\C)$ with the weak topology, i.e., $A\subset P$ is closed (resp. open) $\iff$ $A\cap\ol{e}$ is closed (resp. relatively open) for every cell $e\in\C$.
\end{dfn}

\begin{lmm}\label{ClFnWknss4}
(i) In a CW-polytope $P=(P,\C)$, for any cell $e\in\C$,
$P\ast e=Star_P(e)$ is an open set. (ii) For any space $Y$, a map $f:P\ra Y$ is continuous $\iff$ $f|_{\ol{e}}:\ol{e}\subset P\ra Y$ is continuous for each $e\in\C$.
\end{lmm}
\begin{proof}
{\flushleft (i)} $P\ast e:=\bigcup\{\txt{open cells $\wt{e}$ having $e$ as a face}\}$ is open as a union of open sets.
{\flushleft (ii)} This follows from Lemma \ref{ClFnWknss3}.
\end{proof}

\begin{dfn}[\textcolor{blue}{\index{Nerve of a cover}{Nerve of a cover}}]
Let $X$ be a space and $\U$ an open cover of $X$. Consider the real vector space $R=Span_\Real\B$ spanned by a basis $\B=\{p_U:U\in\U\}$ in 1-1 correspondence with the elements of $\U$. We will say $n+1$ points $p_{U_0},...,p_{U_n}\in\B$ determine an $n$-cell $e=e\left(p_{U_0},\cdots,p_{U_n}\right)$, i.e., an ``open'' $n$-simplex with vertices $\{p_{U_i}\}$, in $R$ if $U_0\cap \cdots\cap U_n\neq\emptyset$. (Note the intersection $U_{i_0}\cap\cdots\cap U_{i_k}$ is nonempty for any subset $\left\{U_{i_0},...,U_{i_k}\right\}\subset\{U_0,...,U_n\}$, which thus determines a $k$-face of the $n$-cell). Let $P=(P,\C)$ be the CW-polytope in $R$ determined by the collection $\C$ of all such cells.

This CW-polytope is denoted by $N(\U)$ and called the \ul{nerve} of the cover $\U$.
\end{dfn}

{\flushleft\dotfill}

\begin{dfn}[\textcolor{blue}{AR($\S$): \index{Absolute! retract}{Absolute retract} for a class of spaces $\S$ (\cite{michael1953J})}]
A space $R\in AR(\S)$ if it is a retract in every closed set inclusion $R\subset X\in\S$.

Equivalently, a space $R\in AR(\S)$ if every continuous map on a closed set $f:R\subset X\ra Y$, $X\in\S$, extends to a continuous map $F:X\ra Y$.
\bc\bt
R\ar[d,hook]\ar[rrr,"f"] &&& Y\\
X\in\S\ar[urrr,dashed,"F"']
\et\ec

If $R$ is an AR for all spaces, then we say simply that $R$ is an AR (or $R\in AR$).
\end{dfn}

\begin{dfn}[\textcolor{blue}{AE($\S$): \index{Absolute! extensor}{Absolute extensor} for a class of spaces $\S$ (\cite{michael1953J})}]
A space $E\in AE(\S)$ if every continuous map on a closed set $f:C\subset X\ra E$, $X\in\S$, extends to a continuous map $F:X\ra E$.
\bc\bt
C\ar[d,hook]\ar[rrr,"f"] &&& E\\
X\in\S\ar[urrr,dashed,"F"']
\et\ec

If $E$ is an AE for all spaces, then we say simply that $E$ is an AE (or $E\in AE$).
\end{dfn}
\begin{note}[\textcolor{OliveGreen}{Every AE is an AR}]
To see this, recall that every continuous map on a closed set $f:C\subset X\ra E$ extends to a continuous map $F:X\ra E$. Thus, if $E\subset X$ is closed, then $id_E:E\ra E$ extends to a retraction $r:X\ra E$.
\end{note}

\begin{dfn}[\textcolor{blue}{\index{Neighborhood! retract}{Neighborhood retract} (NR)}]
A subspace $R\subset X$ is a neighborhood retract of $X$ if $R\subset O$ is a retract for some open set $O\subset X$.
\end{dfn}

\begin{dfn}[\textcolor{blue}{\index{Absolute! neighborhood retract}{Absolute neighborhood retract} (ANR)}]
A space $R$ that is a neighborhood retract in every closed set inclusion $R\subset X$.
\end{dfn}

{\flushleft\dotfill}

\begin{lmm}[\textcolor{OliveGreen}{\index{Canonical cover}{Canonical cover}: \cite[Lemma 2.1]{dugundji1951}}]\label{CanCov}
Let $X$ be a metric space and $C\subset X$ a closed set. There exists a locally finite cover $\U$ of $X-C$ with the following properties (thus called a canonical cover of $X-C$).
\begin{enumerate}[leftmargin=0.7cm]
\item[(a)] Any neighborhood $N(b)$ of a boundary point $b\in\del C$ contains infinitely many $U\in \U$.
\item[(b)] Any neighborhood $N(c)$ of $c\in C$ contains some nbd $N'(c)\subset N(c)$ of $c$ such that
\bea
\U\ast N'(c)\subset N(c),~~~~\txt{i.e.,}~~~~\txt{for all}~~U\in\U,~~U\cap N'(c)\neq\emptyset~~\Ra~~U\subset N(c).\nn
\eea
\end{enumerate}
\end{lmm}
\begin{proof}
Consider the collection of neighborhoods $\big\{N(e):e\in X-C,~\diam N(e)<{1\over 2}\dist(e,C)\big\}$, which is an open cover of $X-C$. Since a metric space is paracompact, this cover has a locally finite open refinement $\U$, which has the desired properties.
\end{proof}

\begin{lmm}[\textcolor{OliveGreen}{\index{Polytope map}{Polytope map} of a cover: \cite[Theorem 2.31]{dugundji1951}}]\label{PolytMap}
Let $X$ be a metric space and $\U$ a locally finite cover of $X$. Then there exists a continuous map ~$g:X\ra N(\U)$~ (call it polytope map of $\U$) such that $g^{-1}\big(N(\U)\ast p_U\big)\subset U$, for any $U\in\U$, where
\bit[leftmargin=0.7cm]
\item[]\centering $N(\U)\ast p_U:=\bigcup\{\txt{open cells $e$ of $N(\U)$ with $p_U$ as a vertex of $\ol{e}$}\}$.
\eit
That is, $g(U)$ contains every open cell of $N(\U)$ for which $p_U$ is a vertex.
\end{lmm}
\begin{proof}
For each $U\in\U$, define a map ~$\ld_U:X\ra \Real$~ by ~$\ld_U(x):={\dist(x,X-U)\over\sum_U\dist(x,X-U)}$, where the sum is well defined since each $x\in X$ lies in finitely many $U\in\U$, and $\ld_U$ is continuous (being locally a ratio of finite sums of continuous functions).

Next, consider the mapping ~$g:X\ra N(\U)$~ given by ~$g(x):=\sum_U\ld_U(x)p_U$. We have
{\small\begin{align}
&\textstyle g^{-1}\big(N(\U)\ast p_U\big)=\left\{x\in X:g(x)\in N(\U)\ast p_U\right\}=\left\{x\in X:\sum_{U'}\ld_{U'}(x)p_{U'}\in N(\U)\ast p_U\right\}\nn\\
&\textstyle~~~~=\big\{x\in X:\sum_{U'}\ld_{U'}(x)p_{U'}\in e,~p_U\in\txt{Vertex}(\ol{e}),~\txt{for an open cell}~e(p_{U_0},...,p_{U_n})~\txt{in}~N(\U)\big\},\nn\\
&\textstyle~~~~\sr{(s)}{=}\left\{x\in X:g(x)=\sum_{i=0}^n\ld_{U_i}(x)p_{U_i},~U\cap U_0\cap\cdots\cap U_n\neq\emptyset\right\}\subset U,\nn
\end{align}}where step (s) is due to the following: Let $U_0,...,U_n$ be the only $U$'s containing $x$. Then, ``$\sum_{U'}\ld_{U'}(x)=1$'' and ``$\ld_{U'}(x)\neq0$ $\iff$ $x\in U'$'' imply that $\sum_{i=0}^n\ld_{U_i}(x)=1$ and so
\bea
\textstyle g(x)=\sum_{i=0}^n\ld_{U_i}(x)p_{U_i}~\in~e(p_{U_0},\cdots,p_{U_n}).\nn
\eea
Since we also have $p_U\in Vertex\big(e(p_{U_0},\cdots,p_{U_n})\big)$, it follows that $U\cap U_0\cap\cdots\cap U_n\neq\emptyset$ and so $x\in U$ (because $U=U_i$ for some $i$).

It remains to show $g$ is continuous. Given $x\in X$, let $x\in U_0,...,U_n$ only. Then $g(x)=\sum_{i=0}^n\ld_{U_i}(x)p_{U_i}$ is in the interior of $\ol{e}=\ol{e(p_{U_0},...,p_{U_n})}$. Let $V\subset N(\U)$ be any open set containing $g(x)$.

[[[Note that in $O:=U_0\cap\cdots\cap U_n\ni x$, $g$ can be decomposed into continuous maps as
\bea
g|_O:O\subset X\sr{(\ld_{U_0},...,\ld_{U_n})}{\ral}\Real^n\sr{L}{\ral}N(\U)\supseteq V\cap\ol{e},\nn
\eea
where $L$ is linear and given by $L(\al_0,...,\al_n)=\sum_{i=0}^n\al_ip_{U_i}$ (the continuity of which requires compatibility of the topology of $N(\U)$ with the Euclidean topology of $\Real^n$).]]]

Then $V\cap\ol{e}$ is open in (the Euclidean topology of) $\ol{e}$, and so the continuity of each $\ld_{U_i}$ shows there exists an open set $W\ni x$ such that $g(W)\subset V\cap\ol{e}\subset V$.
\end{proof}

\begin{lmm}[\textcolor{OliveGreen}{\index{Replacement map}{Replacement map} of the complement of a closed set: \cite[Theorem 3.1]{dugundji1951}}]\label{RepMap}
Let $X$ be a metric space and $C\subset X$ a closed set. Then there exists a space $Y$ and a continuous map ~$h:X\ra Y$~ (a replacement map of $X-C$) with the following properties.
\bit
\item[(a)] $h|_C:C\ra Y$ is a homeomorphism, and $h(C)\subset Y$ is closed.
\item[(b)] $Y-h(C)$ is an infinite CW-polytope, and $h(X-C)\subset Y-h(C)$.
\item[(c)] Each neighborhood of any point $b\in \del h(C)$ contains infinitely many cells of $Y-h(C)$.
\eit
\end{lmm}
\begin{proof}
Let $\U$ be a canonical cover of $X-C$ (as in Lemma \ref{CanCov}), and $N(\U)$ the nerve of the cover. Let $Y=C\sqcup N(\U)$ with the following topology.
\bit[leftmargin=0.7cm]
\item[(i)] In $Y$, $N(\U)$ has its defining topology, as a CW-polytope in $R=Span_\Real\{p_U:U\in\U\}$ over cells $e=e(p_{U_1},\cdots,p_{U_n})\in\C$ in $N(\U)$ determined by points $p_{U_1},\cdots,p_{U_n}\in\B=\{p_U:U\in\U\}$ such that $U_1\cap\cdots\cap U_n\neq\emptyset$. Thus, $N(\U)=\bigcup_{e\in\C}e$, with the weak topology.
\item[(ii)] In $Y$, the topology of $C$ is generated by (a subbase with) all sets in $Y$ of the form
\begin{align}
O_Y(c):=&\big(O_X(c)\cap C\big)~\sqcup~p_\U\big(O_X(c)\big),~~\txt{for all nbds $O_X(c)\subset X$, of every $c\in C$},\nn\\
 p_\U\big(O_X(c)\big):=&\textstyle \bigcup\big\{N(\U)\ast p_U:U\subset O_X(c)\big\}=\bigcup\big\{N(\U)\ast p_U:U\subset O_X(c)-C\big\},\nn
\end{align}
where $N(\U)\ast p_U:=\bigcup\left\{\txt{open cells}~e\in\C:p_U\in\txt{Vertex}(\ol{e})\right\}$, and $p_\U\big(O_X(c)\big)$ is called the \emph{canonical polytope over $X-C$ determined by $O_X(c)\subset X$}.
\eit
Note that $Y$ is Hausdorff, and the subspaces $C,N(\U)\subset Y$ maintain their original topologies.

Using the polytope map of $\U$, ~$g:X-C~\ra~N(\U)$,~ define the desired map $h:X\ra Y$ by ~$h(x):=
\left\{
  \begin{array}{ll}
    x, & x\in C \\
    g(x), & x\in X-C
  \end{array}
\right\}$.~ By construction, it suffices to prove continuity of $h$ at points of $\del C=C~\cap~\ol{X-C}$ only. Let $b\in\del C$. Pick a neighborhood of $h(b)=b$ in $Y$,
\bea
\textstyle O_Y(b):=\big(O_X(b)\cap C\big)~\sqcup~\bigcup\big\{N(\U)\ast p_U:U\subset O_X(b)\big\}.\nn
\eea
Since $\U$ is canonical, by Lemma \ref{CanCov}b, there is a neighborhood $O_X'(b)\subset O_X(b)$ such that for any $U\in\U$, ~$U\cap O_X'(b)\neq\emptyset$ implies $U\subset O_X(b)$. Also, it is clear that $\{U\in\U:~U\subset O_X'(b)\cap (X-C)\}\neq\emptyset$. We will now prove that $h\big(O_X'(b)\big)\subset O_Y(b)$.

First, since $O_X'(b)\subset O_X(b)$, we have $h\big(O_X'(b)\cap C\big)\subset O_X'(b)\cap C\subset O_Y(b)$. Next, if $x\in O_X'(b)\cap (X-C)$, let $x\in U_1\cap\cdots\cap U_n$ only. Then $g(x)\in int~e(p_{U_1},\cdots,p_{U_n})\subset \bigcup_iN(\U)\ast p_{U_i}$, and so $g(x)\in N(\U)\ast p_{U_i}$ for some $i$. But $x\in O_X'(b)\cap U_i$ implies $U_i\subset O_X(b)$, and so $g(x)\in O_Y(b)$. This shows ~$h\big(O_X'(b)\cap(X-C)\big)=g\big(O_X'(b)\cap(X-C)\big)\subset O_Y(b)$.

Hence, $h$ is continuous. It is now clear that $h$ has the desired properties (a),(b),(c).
\end{proof}

\begin{thm}[\textcolor{OliveGreen}{\index{Dugundji-Tietze extension theorem}{Dugundji-Tietze extension theorem}: \cite[Theorem 4.1]{dugundji1951}}]\label{DT-theorem}
Let $X$ be a metric space and $L$ a locally convex space. Then every continuous map on a closed set, $f:C\subset X\ra L$, extends to a continuous map $F:X\ra L$. Moreover, ~$F(X)\subset \Conv\big(f(C)\big)$.
\end{thm}
\begin{proof}
Let Let $\U$ be a canonical covering of $X-C$ and let $h:X\ra C\sqcup N(\U)$ be a replacement map of $X-C$, from Lemma \ref{RepMap}. It suffices to prove that every continuous map $f:C\subset C\sqcup N(\U)\ra L$ extends to a continuous map $\ol{F}:C\sqcup N(\U)\ra L$ (because we then obtain the desired extension $F=\ol{F}\circ h:X\ra L$ as in the diagram below).

\bc\bt
             && C\ar[dll,hook,bend right=15]\ar[d,hook]\ar[rrrr,"f"]         &&&& L \\
X \ar[rr,"h"]\ar[drr,bend right=10,"g"] && Y:=C\sqcup N(\U)\ar[urrrr,dashed,"\ol{F}"'] &&&& \\
             && N(\U)\ar[u,hook]
\et\ec
Let $N(\U)_0=\{p_U:U\in\U\}\subset N(\U)$ be the collection of all vertices of $N(\U)$. Define an extension $f_0:C\sqcup N(\U)_0\ra L$ of $f$ as follows. For each $U\in\U$, pick $x_U\in U$, and then choose $c_U\in C$ such that $d(x_U,c_U)<2\dist(x_U,C)$, i.e., $c_U\in B_{2\dist(x_U,C)}(x_U)\cap C$. Set
\bea
&&f_0(c):=f(c)~~\txt{for all}~~c\in C,\nn\\
&&f_0(p_U):=f(c_U),~~~~c_U\in B_{2\dist(x_U,C)}(x_U)\cap C,~~~~x_U\in U.\nn
\eea
\ul{Proof of continuity of $f_0$}: $f_0$ is continuous on $N(\U)_0$, since the subspace $N(\U)_0\subset N(\U)$ is discrete (the star of a vertex excludes other vertices). We now prove $f_0$ is continuous on $C$.

Let $V\ni f_0(c)=f(c)$ be an open neighborhood. Since $f$ is continuous, we have $\delta>0$ such that $f\left(B_\delta(c)\right)\subset V$, i.e., $d(c,c')<\delta$ $\Ra$ $f(c')\in V$. Let $O_X(c):=B_{\delta/3}(c)$ in $X$. Then for any $U\in\U$ such that $U\subset O_X(c)$, we have
\bea
&& d(c,c_U)\leq d(c,x_U)+d(x_U,c_U)< \delta/3+2d(x_U,C)<\delta/3+2(\delta/3)=\delta,\nn\\
&&~~\Ra~~c_U\in B_\delta(c),~~\Ra~~f_0(p_U)=f(c_U)\in V.\nn
\eea
That is, with $Y:=C\sqcup N(\U)$, all vertices {\small$p_U\in O_Y(c):=\big(O_X(c)\cap C\big)~\sqcup~\bigcup\big\{N(\U)\ast p_U:U\subset O_X(c)\big\}$}
satisfy $f_0(p_U)=f(c_U)\in V$, i.e., $f_0\big(N(\U)_0\cap O_Y(c)\big)\subset V$. So $f_0$ is continuous since
\bea
&&f_0\Big(\big[C\sqcup N(\U)_0\big]\cap O_Y(c)\Big)=f_0\Big(\left[C\cap O_Y(c)\right]\cup\left[N(\U)_0\cap O_Y(c)\right]\Big)\nn\\
&&~~~~=f_0\Big(\left[O_X(c)\cap C\right]\cup\left[N(\U)_0\cap O_Y(c)\right]\Big)=f_0\big(O_X(c)\cap C\big)\cup f_0\Big(N(\U)_0\cap O_Y(c)\Big)\nn\\
&&~~~~=f\big(O_X(c)\cap C\big)\cup f_0\Big(N(\U)_0\cap O_Y(c)\Big)\subset V.\nn
\eea
Now extend $f_0$ linearly (over the vertices $N(\U)_0$) to obtain a map $\ol{F}:C\sqcup N(\U)\ra L$. By linearity and the fact that ``a map on a CW-polytope is continuous $\iff$ it is continuous on the closure of each cell'', $\ol{F}$ is continuous on $N(\U)$. It remains to prove $\ol{F}$ is continuous on $C$.

Let $V\ni \ol{F}(c)=f_0(c)=f(c)$ be a convex neighborhood of $f(c)$ in $L$. Since $f_0$ is continuous at $c$, there is a neighborhood {\small $O_Y(c)=\big(O_X(c)\cap C\big)~\sqcup~\bigcup\big\{N(\U)\ast p_U:U\subset O_X(c)\big\}$} of $c$ in $Y=C\sqcup N(\U)$ such that $f_0\Big(O_Y(c)\cap[C\sqcup N(\U)_0]\Big)\subset V$. Let $O_X'(c)\subset O_X(c)$ be a neighborhood such that for any $U\in\U$, ~~$U\cap O_X'(c)\neq\emptyset$~ $\Ra$~ $U\subset O_X(c)$. Then vertices $p_U\in\{p_U:U\cap O_X'(c)\neq\emptyset\}\subset\{p_U:U\subset O_X(c)\}\subset O_Y(c)$, corresponding to sets in the neighborhood $O_X'(c)$, have images $\ol{F}(p_U)=f_0(p_U)\in V$. Also, given a vertex $p_{U'}\in\{p_{U'}:U'\subset O_X'(c)\}\subset \{p_U:U\cap O_X'(c)\neq\emptyset\}$, if a vertex {\small $p_U\in \ol{N(\U)\ast p_{U'}}$}, then
\begin{align}
U\cap U'\neq\emptyset~~\txt{and}~~U'\subset O_X'(c)~~\Ra~~ U\cap O_X'(c)\neq\emptyset,~~\Ra~~p_U\in O_Y(c),~~\Ra~~\ol{F}(p_U)\in V.\nn
\end{align}
Thus, for any cell $e(p_{U_1},...,p_{U_n})\subset\ol{N(\U)\ast p_{U'}}$ for any $p_{U'}\in\{p_{U'}:U'\subset O_X'(c)\}$, the images of its vertices $p_{U_1},...,p_{U_n}$ under $\ol{F}$ lie in the convex set $V$. It follows that the images (under $\ol{F}$) of linear extensions over such cells lie in $V$ due to convexity of $V$, and so
\bea
\textstyle \ol{F}\big(O'_Y(c)\big)\subset V,~~~~O'_Y(c):=\big(O'_X(c)\cap C\big)~\sqcup~\bigcup\big\{N(\U)\ast p_U:U\subset O'_X(c)\big\}.\nn
\eea
Because $L$ is locally convex (hence has a base for its topology consisting of convex sets), it follows that $\ol{F}$ is continuous.

From the definition of $f_0$ (mapping $N(\U)_0$ into $f(C)$) and the fact that every interior point of $N(\U)$ is a convex linear combination of some vertices in $N(\U)_0$ (and so under extension of $f_0$ by linearity, all points of $N(\U)$ are mapped into the convex hull of $f(C)$), it follows that $\ol{F}(C\sqcup N(\U))\subset \Conv\big(f(C)\big)$. Hence, $F(X)\subset \Conv\big(f(C)\big)$.
\end{proof}

\begin{rmk}
The extension $F:X\ra L$ of $f:C\subset X\ra L$ constructed in the theorem can be written in the following form (where ~$c_U\in B_{2\dist(x_U,C)}(x_U)\cap C$, ~~$x_U\in U$).
\bea
\textstyle F(x)=\left\{
       \begin{array}{ll}
         f(x), & x\in C \\
         \sum_U\ld_U(x)f(c_U), & x\in X-C
       \end{array}
     \right\},~~~~\ld_U(x)={\dist(x,X-U)\over\sum_{U'}\dist(x,X-U')}.\nn
\eea
\end{rmk}

\begin{crl}[\textcolor{OliveGreen}{\cite[Corollary 4.2]{dugundji1951}}]\label{DT-theoremCrl}
Every convex subset $K\subset L$ of a locally convex space $L$ is an absolute extensor (AE) for metric spaces.
\end{crl}
\begin{proof}
By the theorem, $L$ is an AE for metric spaces. Thus, if $X$ is a metric space, then any continuous map on a closed set $f:C\subset X\ra K\subset L$ extends to a continuous map $F:X\ra L$ with $F(X)\subset \Conv(f(C))\subset \Conv(K)=K$ (since $K$ is convex).
\end{proof}  
\section{Metrizability of Spaces}\label{PrelimsMS}   
As for the previous section, this section is not strictly essential to our main purpose, but serves to introduce related terminology and research questions. The main results of interest are (1) metrizability of product spaces in Remark \ref{UMTrmk1}, (2) metrizability of $n$-manifolds in Remarks \ref{CoMfoEmRmk1} and \ref{CoMfoEmRmk2}, (3) metrizability of an arbitrary space in Corollary \ref{NSMTcrl}, and (4) metrizability of the cone of a metrizable space in Proposition \ref{MetrizQuotSp4}.

\subsection{Regularity and the Urysohn metrization theorem}
\begin{dfn}[\textcolor{blue}{Recall: Imbedding of a space}]
A continuous map $f:X\ra Y$ is an imbedding of $X$ into $Y$ (written $f:X\hookrightarrow Y$) if it is a homeomorphism onto its image (i.e., the associated surjection $f:X\ra f(X)\subset Y$ is a homeomorphism).
\end{dfn}

\begin{dfn}[\textcolor{blue}{\index{Metrizable space}{Metrizable space}}]
A space $(X,\T)$ is metrizable if there exists a metric $d:X\times X\ra\Real$ such that $(X,\T)=(X,d)$, i.e., $\T$ is (equal to) a metric topology.
\end{dfn}
Note that given any set $X$, the discrete topology on $X$ is metrizable, and induced by the metric
{\small $d(x,y):=\left\{
          \begin{array}{ll}
            0, & x=y \\
            1, & x\neq y
          \end{array}
        \right\}$}. On the other hand, every non-Hausdorff space is non-metrizable, since a metric space is Hausdorff.

\begin{lmm}[\textcolor{OliveGreen}{\index{Bounded metric}{Bounded metric}}]\label{BdedMet}
Let $(X,d)$ be a metric space. (i) $d_b(x,y):=\min\big(d(x,y),1\big)$, for $x,y\in X$, is a metric on $X$. (ii) $(X,d)=(X,d_b)$ topologically, i.e., $(X,d)\cong(X,d_b)$.
\end{lmm}
\begin{proof}
{\flushleft (i)} It is clear that $d_b(x,y)=d_b(y,x)=0$ iff $x=y$. So we prove the triangle inequality. Observe that $d(x,y)\leq d(x,z)+d(z,y)$ implies
\bea
\min\big(d(x,y),1\big)\leq\min\big(d(x,z)+d(z,y),1\big).\nn
\eea
If $d(x,z)+d(z,y)\leq 1$, then $d(x,z)\leq 1$, $d(z,y)\leq 1$ also, and so
\bea
\min\big(d(x,y),1\big)\leq d(x,z)+d(z,y)=\min\big(d(x,z),1\big)+\min\big(d(z,y),1\big).
\eea
On the other hand, if $d(x,z)+d(z,y)>1$ ( = $t+1-t$ for any $t\in[0,1]$), then it is clear there exists $\al\in[0,1]$ such that $d(x,z)\geq\al$, $d(z,y)\geq1-\al$. This implies $\min\big(d(x,z),1\big)\geq\al$, $\min\big(d(z,y),1\big)\geq1-\al$, which in turn implies $\min\big(d(x,z),1\big)+\min\big(d(z,y),1\big)\geq 1$. Thus,
\bea
\min\big(d(x,y),1\big)\leq 1\leq \min\big(d(x,z),1\big)+\min\big(d(z,y),1\big).
\eea
{\flushleft (ii)} $d_b$ induces the same topology as $d$, since $d_b=d$ on any ball of radius $<1$, and balls of radii $<1$ form a common base for both topologies.
\end{proof}

\begin{lmm}[\textcolor{OliveGreen}{\cite[Theorem 20.5, p.125]{munkres}}]\label{UMTlmm1}
The space $\Real^\infty$ is metrizable.
\end{lmm}
\begin{proof}
Consider the bounded metric $d_b(x,y)=\min(|x-y|,1)$ on $\Real$, which induces the same topology as $|x-y|$, by Lemma \ref{BdedMet}. The function $d(x,y):=\sup_n{d_b(x_n,y_n)\over n}$, for $x,y\in\Real^\infty$, is a metric (as follows). (Note that the restriction of $d$ to $[0,1]^\infty$ is $d|_{[0,1]^\infty}(x,y)=\sup_n{|x_n-y_n|\over n}$).
Let $z\in\Real^\infty$ and $r>0$. Then for any finite set $F\subset\Natural$,
\begin{align}
\textstyle B_r(z)=\left\{x:\sup_n{d_b(x_n,z_n)\over n}<r\right\}\subset\left\{x:\sup_{n\in F}{d_b(x_n,z_n)\over n}<r\right\}=\prod_{n\in F}B'_{nr}(z_n)\times\prod_{n\not\in F}\Real,\nn
\end{align}
which shows the metric topology contains the product topology. On the other hand, fix $\vep>0$ ($\vep<r$) and choose $N$ such that ${1\over N}<\vep<r$. Let $F_\vep:=\{1,...,N\}$. Then
\begin{align}
\textstyle\prod_{n\in F_\vep}B'_{nr}(z_n)\times\prod_{n\not\in F_\vep}\Real=\left\{x:\sup_{n\in F_\vep}{d_b(x_n,z_n)\over n}<r\right\}\subset \left\{x:\sup_n{d_b(x_n,z_n)\over n}<r\right\}=B_r(z),\nn
\end{align}
which shows the product topology contains the metric topology.
\end{proof}

\begin{rmk}\label{UMTrmk1}
Lemma \ref{UMTlmm1} (by its proof) is true in the following more general form:
\bit[leftmargin=0.1cm]
\item[] If $\{X_n\}_{n\in\Natural}=\left\{(X_n,d_{X_n}\right\}_{n\geq 1}$ is any countable family of metrizable spaces, then the Cartesian product space $\prod_{n\in\Natural}X_n:=\big\{\txt{maps}~x:\Natural\ra X\big\}=\big\{(x_n)_{n\in \Natural}:x_n\in X\big\}$ is metrizable, and its topology is induced by the metric $d(x,y):=\sup_n{\min(d_{X_n}(x_n,y_n),1)\over n}$.
\eit
\end{rmk}

\begin{convention}[\textcolor{blue}{Finite product of metric spaces}]
If $(X_1,d_1),...,(X_n,d_n)$ are metric spaces, then (unless stated otherwise) the space $Y=\prod_{i=1}^nX_i$ will be viewed as the metric space $(Y,d)$ with
$d\big((x_1,...,x_n),(x'_1,...,x'_n)\big):=\max_id_i(x_i,x'_i)$.
\end{convention}

\begin{dfn}[\textcolor{blue}{Recall: Regular space, Normal space}]
Let $X$ be a space. Then $X$ is regular if any point $x\in X$ and any closed set $C\not\ni x$ have disjoint neighborhoods. $X$ is normal if every two disjoint closed sets $C_1,C_2\subset X$ have disjoint neighborhoods.
\end{dfn}
Note that in \cite[p.195]{munkres},  a ``regular space'' (resp. a ``normal space'') is assumed to be $T_1$. In our definition above a regular space (resp. a normal space) is not required to be $T_1$.

\begin{dfn}[\textcolor{blue}{Recall: Completely regular space}]
A $T_1$ space $X$ such that for any point $x_0\subset X$ and any closed set $C\subset X$ not containing $x_0$, there is a continuous function $f:X\ra[0,1]$ satisfying $f(x_0)=1$ and $f|_C=0$.
\end{dfn}

\begin{lmm}\label{Reg2ndNorm}
A regular space with a countable base is normal. (Converse: A normal Hausdorff space is completely regular, by the Urysohn lemma.)
\end{lmm}
\begin{proof}
Let $X$ be a regular space with a countable base $\B=\{B_1,B_2,\cdots\}$. Let $C_1,C_2\subset X$ be disjoint closed sets. Let $x\in C_1$. Then by regularity, $x$ has a neighborhood $O_x$ such that $\ol{O}_x\subset X-C_2$. Pick a base element $B_x\in \B$ such that $x\in B_x\subset O_x$. Since $\B$ is countable,
\bea
\B_1:=\{B_x\in\B:x\in C_1\}=\{U_1,U_2,\cdots\},~~~~\ol{U}_i\subset X-C_2~~\txt{for all}~~i,\nn
\eea
is a countable cover of $C_1$. Similarly, we can choose a countable cover of $C_2$,
\bea
\B_2:=\{B_y\in\B:y\in C_2\}=\{V_1,V_2,\cdots\},~~~~\ol{V}_i\subset X-C_1~~\txt{for all}~~i.\nn
\eea
Let $\B'_1=\{U'_1,U'_2,\cdots\}$ and $\B'_2=\{V'_1,V'_2,\cdots\}$ be the collections of open sets given by
{\small\begin{align}
\textstyle U'_n=U_n-\bigcup\limits_{i=1}^n\ol{V}_i=U_n\cap\left(\bigcup\limits_{i=1}^n\ol{V}_i\right)^c\subset U_n,~~~~V'_m=V_m-\bigcup\limits_{i=1}^m\ol{U}_i=V_m\cap\left(\bigcup\limits_{i=1}^m\ol{U}_i\right)^c\subset \bigcup\limits_{i=1}^m\ol{U}_i^c.\nn
\end{align}}Then $U_n'\cap V_m'=\emptyset$ for $n\leq m$, and by symmetry, $U_n'\cap V_m'=\emptyset$ if $n\geq m$, i.e., $U_n'\cap V_m'=\emptyset$ for all $n,m$. Hence, $U=\bigcup U'_i$, $V=\bigcup V'_i$ are disjoint neighborhoods of $C_1$, $C_2$.
\end{proof}

\begin{lmm}\label{UMTlmm2}
Let $X$ be a regular space with a countable base. There exists a countable collection of continuous
functions $\{f_n : X \ra [0,1]\}_{n\geq 1}$ with the following property.
\bit[leftmargin=0.7cm]
\item Given any point $x_0\in X$ and any neighborhood $U$ of $x_0$, there is an index $n'=n'(x_0,U)$ for which $f_{n'}(x_0)>0$ and $f_{n'}|_{U^c}=0$.
\eit
\end{lmm}
\begin{proof}
Since $X$ is normal, by Urysohn's lemma, there exists a (possibly uncountable) family $\F$ of functions having the said property, each of which $f_{(x_0,U)}\in\F$ corresponds to a pair $(x_0,U)$ consisting of a point $x_0\in X$ and a neighborhood $U$ of $x_0$. It remains to pick a countable subfamily of $\F$ that has the same property.

Let $\{U_n\}$ be a countable base for $X$. Given $x_0\in X$ and a neighborhood $U$ of $x_0$, pick an index $m$ such that $x\in U_m\subset U$. Since $X$ is regular, we can choose disjoint open sets $V\ni x$ and $W\supset U_m^c$, i.e., $x\in V\subset W^c\subset U_m\subset U$. Thus, we can choose another base element $U_n$ such that $x\in U_n\subset V\subset W^c\subset U_m\subset U$. That is, $x\in U_n\subset U_m\subset U$ and $\ol{U}_n\subset U_m$. Since $\ol{U}_n\subset U_m$, it follows from Urysohn's lemma that there exists a continuous function $f_{nm}:X\ra [0,1]$ with $f_{nm}|_{\ol{U}_n}=1$ and $f_{nm}|_{U_m^c}=0$ (hence $f_{nm}(x_0)>0$ and $f_{nm}|_{U^c}=0$).

Hence, the countable collection $\{f_{mn}: \ol{U}_n\subset U_m\}\subset\F$ has the desired property.
\end{proof}

\begin{thm}[\textcolor{OliveGreen}{\index{Urysohn metrization theorem}{Urysohn metrization theorem}: \cite[Theorem 34.1, p.215]{munkres}}]\label{UMT}
A regular Hausdorff space with a countable base (i.e., a second countable $T_3$ space) is metrizable.
\end{thm}
\begin{proof}
The strategy is to show that $X$ imbeds into a metrizable space (namely, $\Real^\infty$). Consider the countable collection of cont. functions $\{f_n : X \ra [0,1]\}_{n\geq 1}$ from Lemma \ref{UMTlmm2}.
{\flushleft \ul{Method I}}: Define a map $F:X\ra\Real^\infty$ by $F(x):=\big(f_n(x)\big)_{n=1}^\infty$. Then the following hold.
\vspace{-0.2cm}
{\flushleft (i) \ul{\emph{$F$ is continuous}}}: Since $\Real^\infty$ has the product topology and each $f_n$ is continuous.
\vspace{-0.2cm}
{\flushleft (ii) \ul{\emph{$F$ is injective}}}: For any $x,y\in X$, if $x\neq y$, then some neighborhood $U$ of $x$ excludes $y$, and so there is an index $n$ for which $f_n(x)>0$ and $f_n(y)=0$, i.e., $F(x)\neq F(y)$.
\vspace{-0.2cm}
{\flushleft (iii) \ul{\emph{$F$ is a homeomorphism onto $Z=F(X)$}}}: We need to show $F^{-1}$ is continuous (i.e., $F$ is open). Let $U\subset X$ be open. Let $z_0\in F(U)$. We need to find an open set $W$ such that $z_0\in W\subset F(U)$. Let $x_0=F^{-1}(z_0)\in U$ (i.e., $F(x_0)=z_0$). Pick an index $N$ for which $f_N(x_0)>0$ and $f_N|_{U^c}=0$. Consider the open set in $\Real^\infty$ given by
    \bea
    \textstyle V:=\pi_N^{-1}\big((0,\infty)\big)=\{x\in\Real^\infty:x_N\in(0,\infty)\}=(0,\infty)\times\prod_{i\neq N}\Real,\nn
    \eea
    where $\pi_N:\Real^\infty\ra\Real$ is the continuous map given by $\pi_N(x)=x_N$. Define
    \bea
    \textstyle W:=V\cap Z=\left((0,\infty)\times\prod_{i\neq N}\Real\right)\cap F(X),\nn
    \eea
    which is open in the subspace topology of $Z=F(X)$ in $\Real^\infty$. Now, $z_0\in W\subset F(U)$, since $\pi_N(z_0)=\pi_N(F(x_0))=f_N(x_0)>0$, which implies $z_0\in\pi_N^{-1}\big((0,\infty)\big)\cap Z=W$. Hence,
\bea
&&z\in W~~\Ra~~z=F(x)~~\txt{for some}~~x\in X,~~\txt{and}~~\pi_N(z)=\pi_N(F(x))=f_N(x)\in (0,\infty),\nn\\
&&~~\Ra~~x\in U~~\txt{since}~~f_N|_{U^c}=0,~~\Ra~~z=F(x)\in F(U),~~\Ra~~W\subset F(U).\nn
\eea

Hence, $F$ is an imbedding of $X$ into $\Real^\infty$ (with $F(X)\subset [0,1]^\infty$).

{\flushleft \ul{Method II}}: Consider the collection $\{f_n : X \ra [0,1]\}_{n\geq 1}$ from Lemma \ref{UMTlmm2}, now re-scaled such that $f_n(x)\leq{1\over n}$ for all $x$. Define a function $F:X\ra[0,1]^\infty$ by $F(x)=\big(f_n(x)\big)_{n=1}^\infty$. The bounded uniform metric on $\Real^\infty$, given by $\ol{\rho}(x,y):=\sup_n{\min(|x_n-y_n|,1)}$, when restricted to $[0,1]^\infty$ becomes the uniform metric $\rho(x,y)=\sup_n|x_n-y_n|=\ol{\rho}|_{[0,1]^\infty}(x,y)$.

We will show $F$ is an imbedding with respect (i.e., relative) to the metric $\rho$ on $[0,1]^\infty$.
\vspace{-0.2cm}
{\flushleft (i) \ul{\emph{$F$ is injective}}}: We know from Method I above that $F$ is injective.
\vspace{-0.2cm}
{\flushleft (ii) \ul{\emph{$F^{-1}$ is continuous}}}: Recall that the product topology on $\Real^\infty$ is induced by the metric
    \bea
    \textstyle d(x,y)=\sup_n\min(|x_n-y_n|,1)/n\leq\sup_n\min(|x_n-y_n|,1)=\ol{\rho}(x,y),\nn
    \eea
    which shows that in $\Real^\infty$ (hence in $[0,1]^\infty$ also) the uniform topology contains the product topology ($id:(\Real^\infty,\ol{\rho})\ra (\Real^\infty,d)$ is 1-Lipschitz). From Method I above, $F$ is open with respect to the product topology on $[0,1]^\infty$, and thus is also open with respect to the uniform topology on $[0,1]^\infty$.
\vspace{-0.2cm}
{\flushleft (iii) \ul{\emph{$F$ is continuous}}}: Let $x\in X$. We need to show for any $\vep>0$, there is a neighborhood $U$ of $x$ such that $F(U)\subset B_\vep(F(x))$, i.e., $\rho\big(F(x),F(y)\big)=\sup_n|f_n(x)-f_n(y)|<\vep$ for all $y\in U$. Let $x_0\in X$, and let $\vep>0$. Choose $N$ large enough so that $1/N\leq\vep/2$. For each $n=1,...,N$, by continuity of $f_n$, we have a neighborhood $U_n$ of $x_0$ such that
\bea
|f_n(x)-f_n(x_0)|\leq\vep/2~~~~\txt{for all}~~x\in U.\nn
\eea
Let $U:=U_1\cap\cdots\cap U_N$. Then for any $x\in U$,
{\footnotesize
$\textstyle |f_n(x)-f_n(x_0)|\left\{\!\!\!
                      \begin{array}{ll}
                        \sr{(s1)}{\leq}\vep/2, & \txt{if}~~~~n\leq N \\
                        \vspace{-0.4cm}\\
                        \sr{(s2)}{<}1/N\leq\vep/2, &\txt{if}~~~~n>N
                      \end{array}
                    \!\!\!\right\}
$},
where step (s1) is due to the choice of $U$ and (s2) holds because $f_n(X)\subset[0,1/n]$. Thus, for all $x\in U$, we have $\rho\big(F(x),F(x_0)\big)=\sup_n|f_n(x)-f_n(x_0)|\leq\vep/2<\vep$. \qedhere
\end{proof}

\begin{thm}[\textcolor{OliveGreen}{Imbedding theorem: \cite[Theorem 34.2, p.217]{munkres}}]\label{ImbedThm}
Let $X$ be a $T_1$ space, and $\{f_\al:X\ra\Real\}_{\al\in A}$ a family of continuous functions with the following property.
\bit[leftmargin=0.7cm]
\item For any point $x_0\in X$ and any neighborhood $U$ of $x_0$, there is an index $\al'=\al'(x_0,U)\in A$ such that $f_{\al'}(x_0)>0$ and $f_{\al'}|_{U^c}=0$.
\eit
Then the map $F:X\ra\Real^A$, $F(x):=\big(f_\al(x)\big)_{x\in A}$, is an imbedding of $X$ into $\Real^A$. Moreover, if each $f_\al(X)\subset[0,1]$, then $F$ is an imbedding of $X$ into $[0,1]^A$.
\end{thm}
\begin{proof}
This is similar to Method I of the proof of Theorem \ref{UMT}.
\vspace{-0.2cm}
{\flushleft (i) \ul{\emph{$F$ is continuous}}}: Because $\Real^A$ has the product topology and each $f_\al$ is continuous.
\vspace{-0.2cm}
{\flushleft (ii) \ul{\emph{$F$ is injective}}}: For any $x,y\in X$, if $x\neq y$, then some neighborhood $U$ of $x$ excludes $y$, and so there is an index $\al$ for which $f_{\al}(x)>0$ and $f_\al(y)=0$, i.e., $F(x)\neq F(y)$.
\vspace{-0.2cm}
{\flushleft (iii) \ul{\emph{$F$ is a homeomorphism onto $Z=F(X)$}}}: We need to show $F^{-1}$ is continuous (i.e., $F$ is open). Let $U\subset X$ be open. Let $z_0\in F(U)$. We need to find an open set $W$ such that $z_0\in W\subset F(U)$. Let $x_0=F^{-1}(z_0)$ (i.e., $F(x_0)=z_0$). Pick an index $\beta\in A$ for which $f_\beta(x_0)>0$ and $f|_{U^c}=0$. Consider the open set in $\Real^\infty$ given by
\bea
\textstyle V:=\pi_\beta^{-1}\big((0,\infty)\big)=\{x\in\Real^\infty:x_\beta\in(0,\infty)\}=(0,+\infty)\times\prod_{\al\in A\backslash\beta}\Real,\nn
\eea
where $\pi_\beta:\Real^A\ra\Real$ is the continuous map given by $\pi_\beta(x)=x_\beta$. Define
\bea
\textstyle W:=V\cap Z=\left((0,\infty)\times\prod_{\al\in A\backslash\beta}\Real\right)\cap F(X),\nn
\eea
which is open in the subspace topology of $Z=F(X)$ in $\Real^A$. Now, $z_0\in W\subset F(U)$, since $\pi_\beta(z_0)=\pi_\beta(F(x_0))=f_\beta(x_0)>0$, which implies $z_0\in \pi_\beta^{-1}\big((0,\infty)\big)\cap Z=W$. Hence,
\begin{align}
&z\in W~~\Ra~~z=F(x),~~\txt{for some}~~x\in X,~~\txt{and}~~\pi_\beta(z)=\pi_\beta(F(x))=f_\beta(x)\in (0,\infty),\nn\\
&~~\Ra~~x\in U~~\txt{since}~~f_\beta|_{U^c}=0,~~\Ra~~z=F(x)\in F(U),~~\Ra~~W\subset F(U).\nn
\end{align}

Hence, $F$ is an imbedding of $X$ into $\Real^A$. If each $f_\al(X)\subset[0,1]$, then it is clear that $F$ is an imbedding of $X$ into $[0,1]^A$.
\end{proof}

\begin{thm}[\textcolor{OliveGreen}{Properties of completely regular spaces}]\label{ComRegProps}
(i) A subspace of a completely regular space is completely regular. (ii) A product of completely regular spaces is completely regular.
\end{thm}
\begin{proof}
(i) Let $X$ be a completely regular space and $Z\subset X$. Let $x_0\in Z$ and $C\subset Z$ be a closed set not containing $x_0$. With $\ol{C}$ denoting the closure of $C$ in $X$, we have $x_0\not\in C=\ol{C}\cap Z$, and so $x_0\not\in \ol{C}$. Since $X$ is completely regular, we have a continuous function $f:X\ra [0,1]$ such that $f(x_0)=1$, $f|_{\ol{C}}=0$. Thus, the continuous function $f|_Z\ra[0,1]$ satisfies $f|_Z(x_0)=f(x_0)=1$, $(f|_Z)|_C=f|_C=0$. Hence, $Z$ is completely regular.
{\flushleft (ii)} Let $\{X_\al\}_{\al\in A}$ be a family of completely regular spaces, and $X=\prod_{\al\in A}X_\al$ the product space. Let $z=(z_\al)_{\al\in A}\in X$ and $C\subset X$ a closed set not containing $z$. Choose a base element (i.e., an open rectangle) $U=\prod_{\al\in A}U_\al$, where $U_\al=X_\al$ except for finitely many indices $\al_1,...,\al_n\in A$, such that $z\in U\subset C^c$, i.e., $U\cap C=\emptyset$. For each $i=1,...,n$, choose a continuous function $f_i:X_{\al_i}\ra [0,1]$ such that $f_i(z_{\al_i})=1$, $f_i|_{X_{\al_i}-U_{\al_i}}=0$. Then for each $i$, the map $g_i=f_i\circ\pi_{\al_i}:X\sr{\pi_{\al_i}}{\ral}X_{\al_i}\sr{f_i}{\ral}[0,1]$ is continuous and satisfies
\bea
g_i|_{\pi_{\al_i}^{-1}(z_{\al_i})}=1,~~~~g_i|_{X-\pi_{\al_i}^{-1}(U_{\al_i})}=g_i|_{\pi_{\al_i}^{-1}(X_{\al_i}-U_{\al_i})}=0.\nn
\eea
It follows that the continuous function $f:X\ra[0,1]$, $f(x):=g_1(x)g_2(x)\cdots g_n(x)$ satisfies
\bea
f(z)=1,~~~~f|_{X\backslash U}=0~~(~\Ra~f|_C=0).\nn
\eea
Hence $X$ is completely regular.
\end{proof}

\begin{crl}[\textcolor{OliveGreen}{Complete regularity criterion: \cite[Theorem 34.3, p.218]{munkres}}]\label{ComRegCrit}
A space $X$ is completely regular $\iff$ it is homeomorphic to a subspace of $[0,1]^A$ for some $A$.
\end{crl}
\begin{proof}
$(\Ra)$: If $X$ is completely regular, it follows from Theorem \ref{ImbedThm} that $X$ is homeomorphic to a subspace of $[0,1]^A$ for some $A$.

$(\La)$: If $X$ is homeomorphic to a subspace of $[0,1]^A$ for some $A$, then $X$ is completely regular as a subspace of a completely regular space (Theorem \ref{ComRegProps}).
\end{proof}

\subsection{Imbedding of manifolds}
{\flushleft The} discussion in this section is based directly on \cite[p.224]{munkres}.
\begin{dfn}[\textcolor{blue}{\index{Relatively compact}{Relatively compact} set, \index{Locally! compact}{Locally compact} space}]\label{LocCompSp}
Let $X$ be a space. A set $A\subset X$ is \ul{relatively compact} if its closure $\ol{A}$ is compact. The space $X$ is \ul{locally compact} if each point $x\in X$ has a relatively compact neighborhood.
\end{dfn}
\begin{dfn}[\textcolor{blue}{\index{Compactification}{Compactification} of a space, \index{One-point compactification}{One-point compactification}}]
Let $X$ be a space. A space $Y$ is called a \ul{compactification} of $X$ if (i) $Y$ is a compact Hausdorff space and (ii) there exists an imbedding $f:X\hookrightarrow Y$ with a proper dense image, i.e., $f(X)\neq Y$ and $\ol{f(X)}=Y$. If $Y\backslash f(X)$ is a singleton (i.e., consists of a single point), then $Y$ is a \ul{one-point compactification} of $X$ (which is unique up to homeomorphism by Theorem \ref{LCompChar}).
\end{dfn}

\begin{thm}[\textcolor{OliveGreen}{\cite[Theorem 29.1, p.183]{munkres}}]\label{LCompChar}
Let $X$ be a space. Then $X$ is a locally compact Hausdorff space $\iff$ $X$ has a one-point compactification, in the sense there exists a space $Y$ satisfying the following conditions: (i) $X$ is a subspace of $Y$. (ii) The set $Y-X$ consists of a single point. (iii) $Y$ is a compact Hausdorff space.

Moreover, if $Y$ and $Y'$ are two spaces satisfying these conditions, then there is a homeomorphism $h:Y\ra Y'$ such that $h|_X=id_X:X\ra X$.
\end{thm}
\begin{proof}
{\flushleft \ul{Step 1 (Uniqueness of $Y$)}}: Let $Y,Y'$ be spaces satisfying (i),(ii),(iii). Then $Y=X\sqcup\{p\}$, $Y'=X\sqcup\{p'\}$. Define a map $h:Y\ra Y'$ by $h(p)=p'$, $h|_X=id_X$. We will show $h$ is open (which implies $h$ is a homeomorphism by symmetry). Let $U\subset Y$ be open. If $p\not\in U$, then $h(U)=U\subset X\subset Y'$ is open in $Y'$ since $X$ is open in $Y'$. So, assume $p\in U$. Then $C:=Y-U\subset X\subset Y$ is compact (as a closed subset of the compact space $Y$). Thus $C$ is a compact subspace of $X\subset Y'$, and hence a compact subspace of $Y'$. Since $Y'$ is Hausdorff, $C$ is closed in $Y'$, and so $h(U)=h(Y-C)=h([Y-C-p]\cup p)=[Y-C-p]\cup p'=Y'-C$ is open in $Y'$.

{\flushleft \ul{Step 2 ($\Ra$) (Existence of $Y$)}}:
Assume $X$ is a locally compact Hausdorff space. Let $p$ be any object that is not an element of $X$, and define $Y:=X\cup\{p\}$. (Note that $p$ is often denoted by the symbol $\infty$). Let us give $Y$ the following topology (call it $\T$). $A\subset Y$ is open iff (i) $A=O$ is an open subset of $X$ or (ii) $A=Y-K$ for a compact set $K\subset X$.
\bit[leftmargin=0.7cm]
\item[(a)] It is clear that $\emptyset,Y\in\T$ since $\emptyset$ is both compact and open in $X$.
\item[(b)] $\T$ is closed under finite intersections because if $O_1,O_2$ are open in $X$ and $K_1,K_2$ are compact in $X$, then $O_1\cap O_2$, $(Y-K_1)\cap(Y-K_2)=Y-(K_1\cup K_2)$, and $O_1\cap(Y-K_1)=O_1\cap(X-K_1)$ are all in $\T$.
\item[(c)] Similarly, $\T$ is closed under arbitrary unions because if $\{O_\al\}$ are open in $X$ and $\{K_\beta\}$ are compact in $X$, then $\bigcup O_\al$, $\bigcup(Y-K_\beta)=Y-\bigcap K_\beta$, and $\left[\bigcup O_\al\right]\cup \left[\bigcup(Y-K_\beta)\right]=Y-\left(\bigcap K_\beta-\bigcup O_\al\right)$ are all in $\T$.
\eit

\ul{$X$ is a subspace of $Y$}: If $O$ is open in $X$ and $K$ is compact in $X$, then $O\cap X=O$ and $(Y-K)\cap X=X-K$ are both open in $X$. Moreover, $X\cap\T=\{A\cap X:A\in\T\}$ exhausts all open subsets of $X$. Hence, the topology of $X$ is the subspace topology induced by $\T$.

\ul{$Y$ is compact}: Let $\U$ be an open cover of $Y$. Then $\U$ must contain a set of the type $Y-K$ where $K\subset X$ is compact (otherwise, open subsets of $X$ alone cannot cover the point $p\in Y$). The sets $\V=\{U\cap X:U\in\U,~U\neq Y-K\}$ form an open cover of $K$ in $X$. Since $K$ is compact, finitely many of these sets $\V'=\{U_1\cap X,...,U_n\cap X\}$ cover $K$. Hence, $\U'=\{U_1,...,U_n\}\cup\{Y-K\}\subset\U$ is a finite subcover.

\ul{$Y$ is Hausdorff}: Let $u,v\in Y$ be distinct points. If $u,v\in X$, it is clear that $u,v$ have disjoint open neighborhoods in $Y$. If $u\in X$ and $v=p$, then because $X$ is locally compact, we can choose a compact set $K\subset X$ containing an open neighborhood $U\subset X$ of $u$. It follows that $U\ni u$ and $Y-K\ni v$ are disjoint open neighborhoods of $u,v$.

{\flushleft \ul{Step 3 ($\La$) (Converse)}}: Assume $X\subset Y$, where $Y=X\cup\{p\}$ satisfies (i),(ii),(iii). Then $X$ is Hausdorff (as a subspace of a Hausdorff space). Let $x\in X$. Let $U\ni x$ and $V\ni p$ be disjoint open sets in $Y$. Then $C:=Y-V$ is closed, hence compact, in $Y$. It follows that $C$ is a compact subspace of $X$ containing the neighborhood $U$ of $x$ (and so $\ol{U}\subset C$ is compact).
\end{proof}

\begin{dfn}[\textcolor{blue}{\index{Manifold}{Manifold}, Manifold without boundary, Internal point, Boundary point, \index{Manifold boundary}{Boundary of a manifold}, $n$-Manifold, \index{Curve}{Curve}, \index{Surface}{Surface}}]
Let $X,Y$ be spaces and $\C$ a subcollection of closed subsets of $X$. $Y$ is an \ul{$(X,\C)$-manifold} (or a manifold over $(X,\C)$) if $Y$ is (i) Hausdorff, (ii) second countable, and (iii) each point of $Y$ has a neighborhood that is homeomorphic to an \ul{open} subset of $X$, or homeomorphic to a \ul{closed} subset of $X$ in $\C$. If $Y$ is an $(X,\emptyset)$-manifold, we say $Y$ is an \ul{$X$-manifold without boundary}. If $Y$ is an $(X,\C)$-manifold, a point $y\in Y$ is an \ul{$(X,\C)$-internal point} of $Y$ if it has a neighborhood that is homeomorphic to an open subset of $X$, otherwise $y$ is an \ul{$(X,\C)$-boundary point} of $Y$. The \ul{$(X,\C)$-boundary} of $Y$ is the collection $\del_{(X,\C)} Y:=$ $\{$\txt{all $(X,\C)$-boundary points of $Y$}$\}$.

An $\big(\Real^n,\{\Real^{n-1}\times[0,\infty)\}\big)$-manifold is called an \ul{$n$-manifold} (or an $n$-dimensional manifold). A $1$-manifold is called a curve, and a $2$-manifold is called a surface.
\end{dfn}

Henceforth, ``\emph{manifold}'' we will mean ``\emph{$n$-manifold without boundary (for some $n\geq 1$)}'', unless specified otherwise.

\begin{lmm}\label{NormCoverLmm}
If $X$ is a normal space, then every finite open cover $\{U_1,...,U_n\}$ has a finite open refinement $\{V_1,...,V_n\}$ such that $\ol{V}_i\subset U_i$ for all $i=1,...,n$.
\end{lmm}
\begin{proof}
We will proceed by induction. Observe that the set $C=X-(U_2\cup\cdots\cup U_n)$ is a closed set contained in $U_1$. Since $X$ is normal, there are disjoint open sets $V_1\supset C$, $W_1\supset U_1^c$, and so $C\subset V_1\subset W_1^c\subset U_1$. Thus, $C\subset V_1\subset\ol{V}_1\subset U_1$, since $W_1^c$ is closed.

The collection $\{V_1,U_2,\cdots,U_n\}$ covers $X$. As before, we have the closed set
\bea
C=X-\big(V_1\cup U_3\cup\cdots\cup U_n\big)\subset U_2,\nn
\eea
and so there exists an open set $V_2\supseteq C$ such that $\ol{V}_2\subset U_2$, where we again see that $\{V_1,V_2,U_3,\cdots,U_n\}$ covers $X$. Continuing this way, at the $k$th step, we have the closed set
\bea
C=X-\big(V_1\cup\cdots\cup V_{k-1}\cup U_{k+1}\cup\cdots\cup U_n\big)\subset U_k,\nn
\eea
and so there exists an open set $V_k\supseteq C$ such that $\ol{V}_k\subset U_k$, where we as usual see that $\{V_1,\cdots,V_k,U_k,\cdots,U_n\}$ covers $X$. At step $k=n$, we get the desired cover $\{V_1,\cdots,V_n\}$.
\end{proof}

\begin{thm}[\textcolor{OliveGreen}{\cite[Theorem 36.1, p.225]{munkres}}]
Let $X$ be a normal space and $\{U_1,\cdots, U_n\}$ a finite open cover of $X$. Then there exists a finite partition of unity dominated by $\{U_i\}$, i.e., continuous functions $\left\{f_i:X\ra[0,1]\right\}_{i=1}^n$ such that
\bea
\textstyle f_i|_{X\backslash U_i}=0,~~\txt{for each}~~i=1,...,n,~~~~\txt{and}~~~~\sum_{i=1}^nf_i(x)=1~~\txt{for each}~~x\in X.\nn
\eea
\end{thm}
\begin{proof}
By Lemma \ref{NormCoverLmm}, there exists an open cover $\{V_1,...,V_n\}$ such that $\ol{V}_i\subset U_i$ for each $i$. Similarly, there exists a further open cover $\{W_1,...,W_n\}$ such that $\ol{W}_i\subset V_i$  for each $i$. By Urysohn's lemma, for each $i$, there exists a continuous function $g_i:X\ra[0,1]$ such that $g_i|_{\ol{W}_i}=1$ and $g_i|_{V_i^c}=0$ (which implies $g_i|_{U_i^c}=0$ as well). Since the collection $\{W_1,...,W_n\}$ covers $X$, the sum $\sum_i g_i(x)>0$ for each $x\in X$, and so the desired functions $f_i:X\ra[0,1]$ are given by ~$f_i(x):={g_i(x)\over\sum_jg_j(x)}$.
\end{proof}

\begin{thm}[\textcolor{OliveGreen}{\cite[Theorem 36.2, p.226]{munkres}}]\label{ComMfoldEmb}
Every compact $n$-manifold $X$ imbeds into $\Real^N$ for some $N\geq n$.
\end{thm}
\begin{proof}
Since $X$ is compact, we can choose a finite open cover $\{U_1,...,U_m\}$, along with imbeddings $g_i:U_i\ra\Real^n$ for each $i$. Note that $X$ is normal (as a compact Hausdorff space), and so there exists a finite partition of unity $\left\{\vphi_i:X\ra[0,1]\right\}_{i=1}^m$ dominated by $\{U_1,...,U_m\}$.

Let $A_i:=\txt{Supp}(\vphi_i):=\txt{closure}{\{\vphi_i>0\}}\subset U_i$. For each $i=1,...,m$, define a map $h_i:X\ra\Real^n$ by~
{\footnotesize $h_i(x):=\left\{
          \begin{array}{ll}
            \vphi_i(x)g_i(x), & x\in U_i \\
            0, & x\in A_i^c=X\backslash A_i .
          \end{array}
        \right\}$},~
which is well defined and continuous, because its restrictions to the open sets $U_i,A_i^c$ are continuous. Define the continuous map
\bea
F:X\ra\Real^{m+mn},~~~~F(x):=\big(\vphi_1(x),\cdots,\vphi_m(x),h_1(x),\cdots,h_m(x)\big).\nn
\eea
Since $X$ is compact, to show $F$ is an imbedding, it suffices to show $F$ is injective (footnote\footnote{Recall that a continuous bijection from a compact space to a Hausdorff space is a homeomorphism.}). Let $x,y\in X$ be such that $F(x)=F(y)$. Then $\vphi_i(x)=\vphi(y)$, $h_i(x)=h_i(y)$ for all $i$. We know $\phi_j(x)>0$ for some index $j\in\{1,...,m\}$ (since $\sum_i\vphi_i(x)=1>0$), which implies $\vphi_j(y)>0$ also. Thus, $x,y\in A_j\subset U_j$, and so $x=y$ (by injectivity of $g_j$) since
\bea
\vphi_j(x)=\vphi_j(y),~~\vphi_j(x)g_j(x)=h_j(x)=h_j(y)=\vphi_j(y)g_j(y)~~\Ra~~g_j(x)=g_j(y).\nn\qedhere
\eea
\end{proof}

\begin{lmm}
A locally compact Hausdorff space is completely regular.
\end{lmm}
\begin{proof}
Let $X$ be a locally compact Hausdorff space. Then $X$ has a one-point compactification, and is thus a subset of a compact Hausdorff space (Theorem \ref{LCompChar}). But a compact Hausdorff space is normal, hence completely regular, and a subspace of a complete regular space is completely regular (Theorem \ref{ComRegProps}). It follows that $X$ is completely regular.
\end{proof}

\begin{rmk}[\textcolor{OliveGreen}{Metrizability of $n$-manifolds}]\label{CoMfoEmRmk1}
Since an $n$-manifold is locally $\Real^n$, it is a locally compact Hausdorff space. Thus, because a locally compact Hausdorff space is completely regular, it follows by Theorem \ref{UMT} that $n$-manifolds are metrizable.
\end{rmk}

\begin{rmk}\label{CoMfoEmRmk2}
Theorem \ref{ComMfoldEmb} holds for separable $n$-manifolds even without compactness (namely, every separable $n$-manifold imbeds into $\Real^{2n+1}$), but the proof is more involved. See
Theorem \ref{ImbedThm2}, and related results in \cite{PWZ1961}. Note that separability is strictly necessary here, because \ul{every subspace of a separable metric space (such as $\Real^{2n+1}$) is separable}, and the proof is as follows: If $X$ is a separable metric space, $A\subset X$ any subset, and $\{x_1,x_2,\cdots\}\subset X$ a countable dense subset, then because $X=\bigcup_{i\geq 1}\left\{B_q(x_i):q\in\Rational\cap(0,\infty)\right\}$, if we pick $a_{q,i}\in B_q(x_i)\cap A$ whenever the intersection is nonempty, then the set $A':=\left\{a_{q,i}:q\in\Rational\cap(0,\infty),i\geq 1\right\}$ is a countable dense subset of $A$. To see this, let $a\in A$, $\vep>0$, and take a sequence $y_k:=x_{i(k)}\ra a$.
Pick $k$ and a positive rational $q_k$ such that $d(a,y_k)\leq q_k<\vep/2$. Then $B_{q_k}(y_k)\cap A\neq\emptyset$, which in turn implies $B_{q_k}(y_k)\cap A'\neq\emptyset$ (by the definition of $A'$). Thus, with $a'\in B_{q_k}(y_k)\cap A'$, we get
\bea
d(a,a')\leq d(a,y_k)+d(y_k,a')\leq 2q_k<\vep.\nn
\eea
\end{rmk}

\subsection{The Nagata-Smirnov metrizability criterion}
\begin{dfn}[\textcolor{blue}{\index{$G_\delta$ set}{$G_\delta$ set} in a space}]
A countable intersection of open sets (of the space).
\end{dfn}
Note that complements of $G_\delta$ sets (i.e., countable unions of closed sets) are called \index{$F_\sigma$ set}{$F_\sigma$ sets}.

\begin{dfn*}[\textcolor{blue}{Recall: Locally finite family}]
If $X$ is a space, a family of subsets $\{A_\al\}_{\al\in\Gamma}$ is locally finite if each $x\in X$ has a neighborhood that intersects only finitely many $A_\al$.
\end{dfn*}

\begin{dfn}[\textcolor{blue}{\index{Countably locally finite}{Countably locally finite} family}]
If $X$ is a space, a family of subsets $\A=\{A_\al\}$ is countably locally finite if it is a countable union of locally finite families, i.e., $\A=\bigcup_{n\in\Natural}\A_n$, where each $\A_n$ is a locally finite collection.
\end{dfn}

\begin{lmm}[\textcolor{OliveGreen}{Recall: \cite[Lemma 39.1]{munkres}}]\label{NSMTlmmRcll}
Let $X$ be a space. If $\A\subset\P(X)$ is locally finite, then (a) $\C=\left\{\ol{A}:A\in\A\right\}$ is locally finite and (b) $\ol{\bigcup_{A\in\A}A}=\bigcup_{A\in\A}\ol{A}$.
\end{lmm}
\begin{proof}
See Lemma \ref{NSMTlmm}.
\end{proof}

\begin{lmm}[\textcolor{OliveGreen}{\cite[Lemma 39.2]{munkres}}]\label{NSMTlmm0}
If $X$ is a metrizable space, then every open cover $\U$ of $X$ has
an open refinement that is countably locally finite.
\end{lmm}
\begin{proof}
Let $X$ be a metrizable space and $\U$ an open covering. Since a metric space is paracompact, $\U$ has an open refinement $\V$ that is locally finite, and hence countably locally finite. (For a proof without paracompactness, see the proof of \cite[Lemma 39.2]{munkres})
\end{proof}

\begin{lmm}[\textcolor{OliveGreen}{\cite[Lemma 40.1]{munkres}}]\label{NSMTlmm1}
Let $X$ be a regular space with a countably locally finite base $\B$. Then (i) $X$ is normal, and (ii) every closed set $C\subset X$ is a $G_\delta$ set in $X$.
\end{lmm}
\begin{proof}
\ul{Step I}: Let $W\subset X$ be open. There is a countable collection of open sets $\{U_1,U_2,\cdots\}$ such that ~$W=\bigcup U_n=\bigcup\ol{U}_n$ (proved as follows). Let $\B=\bigcup\B_n$, where each collection $\B_n$ is locally finite. Let $\C_n:=\{B\in\B_n:~\ol{B}\subset W\}$, which is locally finite (since $\B_n$ is locally finite). Define the open sets $U_n:=\bigcup_{B\in\C_n}B$. Then by Lemma \ref{NSMTlmmRcll},
\bea
\textstyle \ol{U}_n=\bigcup_{B\in\C_n}\ol{B}\subset W,~~~~\Ra~~~~\bigcup U_n\subset\bigcup\ol{U}_n\subset W.\nn
\eea
Let $x\in W$. Then because $X$ is regular, there is $B\in\B$ such that $x\in B$ and $\ol{B}\subset W$ (footnote\footnote{\label{RegFnt1} To check this, take disjoint neighborhoods $N(x)$ and $N(W^c)$, and then choose $B\in \B$ such that $x\in B\subset N(x)\subset N(W^c)^c\subset W$. Since $N(W^c)^c$ is closed, the result follows.}). Now, $B\in \B_m$ for some index $m$, which implies $B\in\C_m$, and so $x\in U_m$. Hence $W\subset\bigcup U_n$, and so $W=\bigcup U_n=\bigcup\ol{U}_n$.

{\flushleft \ul{Step II}}: Let $C\subset X$ be closed. Then $C$ is a $G_\delta$ set (proved as follows). From Step I, $C^c=\bigcup\ol{U}_n$, and so $C=\bigcap\ol{U}_n^c$, (a countable intersection of open sets).

{\flushleft \ul{Step III}}: $X$ is normal (proved as follows). Let $C,D\subset X$ be disjoint closed sets. By Step I, there are countable collections of open sets $\{U_1,U_2,\cdots\}$ and $\{V_1,V_2,\cdots\}$ such that
\begin{align}
&\textstyle D^c=\bigcup U_n=\bigcup\ol{U}_n,~~\txt{where $\{U_n\}$ covers $C$ and each $\ol{U}_n$ is disjoint from $D$},\nn\\
&\textstyle C^c=\bigcup V_n=\bigcup\ol{V}_n,~~\txt{where $\{V_n\}$ covers $D$ and each $\ol{V}_n$ is disjoint from $C$}.\nn
\end{align}
As in the proof of Lemma \ref{Reg2ndNorm}, if we define the sets (which satisfy $U'_n\cap V'_m=\emptyset$ for all $m,n$)
\bea
\textstyle U_n':=U_n-\bigcup_{i=1}^n\ol{V}_i,~~~~V_n':=V_n-\bigcup_{i=1}^n\ol{U}_i,\nn
\eea
then we get disjoint neighborhoods of $C$ and $D$ given by $U':=\bigcup_{n=1}^\infty U_n'$, $V':=\bigcup_{n=1}^\infty V_n'$.
\end{proof}

\begin{lmm}[\textcolor{OliveGreen}{\cite[Lemma 40.2]{munkres}}]\label{NSMTlmm2}
If $X$ is a normal $T_1$ space, then for every closed $G_\delta$ set $D\subset X$, there exists a continuous function $f:X\ra [0,1]$ such that $f|_D=0$ and $f|_{D^c}>0$.
\end{lmm}
\begin{proof}
Let $D=\bigcap U_n$, where $U_n$ are open sets. For each $n$, let $f_n:X\ra[0,1]$ be a continuous function (by Urysohn's lemma) such that $f_n|_D=0$ and $f_n|_{U_n^c}=1$. Define the function $f:X\ra\Real$ given for each $x\in X$ by the series $f(x):=\sum_{i=1}^\infty{f_n(x)\over 2^n}$
which converges uniformly, since $0\leq f(x)\leq\sum{1\over 2^n}$, and so $f$ is continuous. Moreover,
\bea
\textstyle f|_D=\sum\limits_{i=1}^\infty{f_n|_D\over 2^n}=0,~~~~f|_{D^c}=f|_{\bigcup U_n^c}=\sum\limits_{i=1}^\infty{f_n|_{\bigcup U_m^c}\over 2^n}\geq \sum\limits_{i=1}^\infty{f_n|_{U_n^c}\over 2^n}=\sum\limits_{i=1}^\infty{1\over 2^n}>0.\qedhere\nn
\eea
\end{proof}

\begin{thm}[\textcolor{OliveGreen}{Nagata-Smirnov: \cite[Theorem 40.3, p.250]{munkres}}]\label{NSMT}
Let $X$ be space. Then $X$ is metrizable $\iff$ $X$ is a regular Hausdorff space with a countably locally finite base $\B$.
\end{thm}
\begin{proof}
$\bullet$($\Ra$): Assume $X$ is metrizable. Then we know $X$ is a normal Hausdorff (hence regular Hausdorff) space. It remains to find a countably locally finite base for $X$. Choose a metric $d$ that induces the topology of $X$. Then for each $n\in \Natural$, we have the cover
\bea
\A_n:=\left\{B_{1/n}(x):x\in X\right\},\nn
\eea
which has a countably locally finite open refinement $\B_n$ by Lemma \ref{NSMTlmm0}. Thus, the collection $\B=\bigcup_{n=1}^\infty\B_n$ is also countably locally finite. We will show $\B$ is a base for $X$ (i.e., for any $x_0\in X$ and a neighborhood $U$ of $x_0$, there exists $B\in\B$ such that $x_0\in B\subset U$). Let $x\in X$, $\vep>0$, and consider the neighborhood $U=B_\vep(x)$ of $x$. Choose $n$ so that $1/n<\vep/2$. Since $\B_n$ covers $X$, we have $x\in B$ for some $B\in\B_n$. It follows that for any $y\in B$, we have
\bea
d(x,y)\leq \diam B\leq 2/n\leq 2(\vep/2)=\vep,~~\Ra~~y\in B_\vep(x),~~\Ra~~x\in B\subset B_\vep(x).\nn
\eea

{\flushleft $\bullet$($\La$):} Suppose $X$ is a regular Hausdorff space with a countably locally finite base $\B$. Then by Lemma \ref{NSMTlmm1}, $X$ is normal and every closed set in $X$ is a $G_\delta$ set in $X$. We will prove that $X$ imbeds into $(\Real^J,\ol{\rho})$ for some set $J$, where $\ol{\rho}$ is the uniform metric on $\Real^J$ given by
\bea
\textstyle\ol{\rho}(x,y):=\sup_{\al\in J}\min(|x_\al-y_\al|,1).\nn
\eea
Let $\B=\bigcup_{n\in\Natural}\B_n$, where each $\B_n$ is locally finite. For each $n\geq 1$, and each $B\in\B_n$, there exists (by Lemma \ref{NSMTlmm2} with $D=B^c$, and rescaling) a continuous function
\bea
f_{n,B}:X\ra [0,1/n]~~\txt{such that}~~f_{n,B}|_B>0~~\txt{and}~~f_{n,B}|_{B^c}=0.\nn
\eea
The collection $\{f_{n,B}:n\in\Natural,B\in\B_n\}$ \ul{\emph{separates points from closed sets}} in $X$, in the sense that given a point $x_0\in X$ and a closed set $C$ not containing $x_0$ (or equivalently, a neighborhood $U=C^c$ of $x_0$), then with a $B\in \B_n$ (for some $n$) satisfying $x_0\in B\subset C^c$, we have
\bea
f_{n,B}(x_0)>0~~(\txt{since}~~f_{n,B}|_B>0)~~\txt{and}~~f_{n,B}|_C=f_{n,B}|_{U^c}=0~~(\txt{since}~~f_{n,B}|_{B^c}=0,~C\subset B^c).\nn
\eea
Let $J:=\{(n,B):n\in\Natural, B\in\B_n\}=\bigcup_{n=1}^\infty\big(\{n\}\times\B_n\big)\subset\Natural\times\B$. Define the map
\bea
\textstyle F:X\ra[0,1]^J,~~~~F(x):=\big(f_{n,B}(x)\big)_{(n,B)\in J}.\nn
\eea
Then by Theorem \ref{ImbedThm}, $F$ imbeds $X$ into $[0,1]^J$ (with respect to the usual, = product, topology on $[0,1]^J$). In what follows, we will show that $F$ is also an imbedding with respect to the uniform topology $([0,1]^J,\rho)=([0,1]^J,\ol{\rho})\subset (\Real^J,\ol{\rho})$, where
\bea
\textstyle \rho(x,y)=\sup_{\al\in J}|x_\al-y_\al|.\nn
\eea
We know that the uniform topology contains the product topology, and so with respect to the uniform topology, $F$ is injective and maps open sets to open sets (i.e., $F^{-1}$ exists and is continuous). We will now show $F$ is also continuous with respect to the uniform topology.

Let $x_0\in X$ and fix $\vep>0$. We need to find a neighborhood $W$ of $x_0$ such that
\bea
\rho(F(x_0),F(x))<\vep~~~~\txt{for all}~~x\in W.\nn
\eea
Fix $n\geq 1$. Let $U_n$ be a neighborhood of $x_0$ that intersects finitely many $B\in \B_n$. Then
\bea
f_{n,B}|_{U_n}=0~~~~\txt{for all but finitely many}~~B\in\B_n.\nn
\eea
Thus, because $f_{n,B}$ is continuous, we can choose another neighborhood $V_n\subset U_n$ of $x_0$ such that for each of the (remaining) functions $f_{n,B}$, $B\in\B_n$, satisfying $f_{n,B}|_{V_n}>0$, we have
\bea
\textstyle \sup_{x,y\in V_n}|f_{n,B}(x)-f_{n,B}(y)|\leq\vep/2~~~~(\txt{which holds for all}~~B\in\B_n).\nn
\eea
Consider all the neighborhoods $\{V_n\}_{n\in\Natural}$ of $x_0$. Choose $N$ such that $1/N\leq\vep/2$, and define $W:=V_1\cap\cdots\cap V_N$. To show $W$ is the desired neighborhood of $x_0$, let $x\in W$. Then
{\small\bea
|f_{n,B}(x_0)-f_{n,B}(x)|\left\{
                             \begin{array}{ll}
                               \sr{(s1)}{\leq}\vep/2, & \txt{if}~~n\leq N \\
                               \vspace{-0.4cm}\\
                               \sr{(s2)}{\leq}1/n<\vep/2, & \txt{if}~~n>N
                             \end{array}
                           \right\},~~~~\txt{for all}~~~~(n,B)\in J,\nn
\eea}
where step (s1) is due to the choice of $V_1,...,V_n$, and step (s2) holds because $f_{n,B}(X)\subset[0,1/n]$. Hence, ~$\rho(F(x_0),F(x))=\sup_{(n,B)\in J}|f_{n,B}(x_0)-f_{n,B}(x)|\leq\vep/2<\vep$.
\end{proof}

\begin{crl}\label{NSMTcrl}
A space $X$ is metrizable $\iff$ $X$ is homeomorphic to a subspace of $[0,1]^J$ equipped with the uniform metric, for some set $J$.
\end{crl}
\begin{proof}
($\Ra$) If $X$ is metrizable, then by the proof of Theorem \ref{NSMT}, $X$ imbeds into some $[0,1]^J$, as a metric space equipped with the uniform metric.
($\La$) If $X$ imbeds into some $[0,1]^J$, as a metric space equipped with the uniform metric, then $X$ is metrizable as a subspace of a metric space.
\end{proof}

\subsection{Metrizability of quotient spaces}  
{\flushleft We} have already seen that the ``quotient space operation'' can destroy the Hausdorff property which is essential for metrizability. Let $I$ denote the unit interval $[0,1]\subset\Real$.

\begin{dfn}[\textcolor{blue}{\index{Local homeomorphism}{Local homeomorphism}}]
A continuous map $f:X\ra Y$ is a local homeomorphism if every point $x\in X$ has an open neighborhood $U\subset X$ such that (i) $f(U)\subset Y$ is open, and (ii) $f|_U:U\ra f(U)$ is a homeomorphism.
\end{dfn}

\begin{lmm}\label{MetrizQuotSp1}
A bijective local homeomorphism is a homeomorphism.
\end{lmm}
\begin{proof}
Let $f:X\ra Y$ be a bijective local homeomorphism. We need to show $f^{-1}:Y\ra X$ is continuous (i.e., $f=(f^{-1})^{-1}$ is open). Let $U\subset X$ be open. For each $x\in U$, let $U_x$ be a neighborhood of $x$ such that $f|_{U_x}:U_x\ra f(U_x)$ is a homeomorphism (where $f(U_x)$ is open). Then $f(U)=f\left(\bigcup_{x\in U}U\cap U_x\right)=\bigcup_{x\in U}f|_{U_x}(U\cap U_x)$ is open as a union of open sets (since each $f|_{U_x}$ is a homeomorphism, hence an open map).
\end{proof}

\begin{prp}\label{MetrizQuotSp2}
Let $X$ be a metric space and $q:X\ra Y={X\over\sim}$ a quotient of $X$. If
\bit[leftmargin=0.9cm]
\item[(i)] $d_Y$ is a metric on $Y$ such that $q:X\ra (Y,d_Y)$ is continuous, and
\item[(ii)] $X$ has a cover by compact sets $\{K_\al\}$ such that the interiors $\left\{\big(q(K_\al)\big)^o\right\}$ cover $Y$,
\eit
then $Y$ is metrizable, and the topology of $Y$ is induced by $d_Y$, i.e., $Y\cong (Y,d_Y)$.
\end{prp}
\begin{proof}
Since $X\sr{q}{\ral}(Y,d_Y)$ is continuous, by the universal property of quotient maps, we have a continuous map $f:Y\ra (Y,d_Y)$ such that $f\circ q=q$ (which implies $f$ is injective, and moreover, $f=id$ since $q$ is surjective).
\bea\bt
X\ar[d,"q"']\ar[rr,"q"] && (Y,d_Y)\\
Y\ar[rru,dashed,"f"']
\et
\hspace{2cm}
\bt
K_\al\ar[d,"q|_{K_\al}"']\ar[rr,"q|_{K_\al}"] && (Y,d_Y)\\
q(K_\al)\ar[rru,dashed,"f|_{q(K_\al)}"']
\et\nn\eea
Because each $q(K_\al)$ is compact, each $f|_{q(K_\al)}:q(K_\al)\ra f(q(K_\al))$ is a homeomorphism (as a continuous bijection from a compact space onto a Hausdorff space). Thus, $f$ is a local homeomorphism since each $f|_{q(K_\al)^o}:q(K_\al)^o\ra f(q(K_\al)^o)=[f(q(K_\al))]^o$ is a homeomorphism.

It follows that $f$ is a homeomorphism (since $f$ is an injective local homeomorphism).
\end{proof}

\begin{dfn}[\textcolor{blue}{\index{Cylinder of a space}{Cylinder of a space} $X$}]
The product space ~$C(X):=X\times I$.
\end{dfn}

\begin{dfn}[\textcolor{blue}{\index{Cone!}{Cone of a space} $X$, Vertex, Base, \index{Cone! metric $d_c$}{Cone metric $d_c$}}]
The \ul{cone} of $X$ is the quotient space $\Cone(X):={C(X)\over X\times\{0\}}={X\times I\over X\times\{0\}}$. With the quotient map $q:X\times I\ra \Cone(X)$, $(x,t)\mapsto q(x,t)\eqv[(x,t)]$, the \ul{vertex} of the cone is the point $v:=q(x,0)=X\times\{0\}\in \Cone(X)$ and its \ul{base} is the subset $q(X\times\{1\})\cong X\times\{1\}\cong X$.

A \ul{cone metric} $d_c:\Cone(X)\times\Cone(X)\ra[0,\infty)$ on $\Cone(X)$ is the map \ul{induced} (when possible) by any function $d^c:(X\times I)^2\ra[0,\infty)$ that satisfies the following three conditions.
\bit[leftmargin=0.9cm]
\item[(i)] $d^c\big((x,t),(x',t')\big)=d^c\big((x',t'),(x,t)\big)\leq d^c\big((x',t'),(x'',t'')\big)+d^c\big((x'',t''),(x',t')\big)$.
\item[(ii)] $d^c|_{[(X\times I)-(X\times\{0\})]^2}$ is a metric on $(X\times I)-(X\times\{0\})$.
\item[(iii)] $d^c\big((x,t),(x',0)\big)=c(x,t):=t$ ~for all~ $t\in[0,1]$, $x,x'\in X$.
\eit
\end{dfn}

\begin{prp}[\textcolor{OliveGreen}{$d_c$-Metrizability of $\Cone(X)$}]\label{MetrizQuotSp4}
Let $X$ be a metric space. $\Cone(X)\cong\big(Cone(X),d_c\big)$ $\iff$ $X$ is compact. (Assume $q:X\times I\ra(\Cone(X),d_c)$ is continuous.)
\end{prp}
\begin{proof}
($\Ra$): Assume $\Cone(X)\cong\big(Cone(X),d_c\big)$. Suppose $X$ is not compact. Then $X\times I$ is not compact (otherwise, $X$ will be compact as a closed subspace of a compact space). Let $x_n\in X$ be a sequence with no convergent subsequence. Then for any convergent sequence $0<t_n\ra 0$ in $I$, $(x_n,t_n)$ has no convergent subsequence in $X\times I$, since for any $x\in X$, $d\big((x_{n_k},t_{n_k}),(x,0)\big)=\max\big(d(x_{n_k},x),|t_{n_k}-0|\big)\ra 0$ iff $d(x_{n_k},x)\ra 0$. It follows that the set $A=\{(x_n,t_n)\}\subset X\times I$ is closed for any sequence $0<t_n\ra 0$.

Thus, $A$ is also closed in $\Cone(X)$ since the quotient map $q:X\times I\ra \Cone(X)$ is injective on $\big((X\times I)\backslash (X\times\{0\})\big)\supset A$, and so $A=q^{-1}(q(A))$ is closed $\iff$ $q(A)$ is closed. However, $A$ is not closed in $\big(\Cone(X),d_c\big)$ because $(x_n,t_n)\sr{d_c}{\ra}(x,0)$, since $d_c\big((x_n,t_n),(x,0)\big)=t_n\ra 0$. This contradicts the assumption $\Cone(X)\cong\big(Cone(X),d_c\big)$.
{\flushleft ($\La$):} Assume $X$ is compact. Then it follows from Proposition \ref{MetrizQuotSp2} that $\Cone(X)\cong \big(\Cone(X),d_c\big)$.
\end{proof}
\section{Metrical Extension Theorems: Absolute Lipschitz Retracts}\label{PrelimsMET}  
Unlike the previous two sections, this section is closely related to our subsequent discussion. The most relevant concepts are those of absolute Lipschitz retracts, biLipschitz equivalence, quasiconvexity, metrical convexity, and binary intersection property. Any extension theorems associated with these concepts are equally important. The main results of interest include the following. (i) Lipschitz extendability of maps from a subset of a metric space to $\Real^n$ in Corollary \ref{MWETlmm}. (ii) Bi-Lipschitz equivalence of the simplex and the cube in Corollary \ref{BiLipHomC}. (iii) The characterization of geodesics in Lemma \ref{GeodCharLmm}. (iv) Gluing of locally Lipschitz maps in Lemmas \ref{GluLmmI} and \ref{GluLmmII}. (v) Proof that $\ell^\infty(A)$ is an AUR in Lemma \ref{ExtToLinf}.
\subsection{The McShane-Whitney extension theorem}\label{PrlMET1}
(Based on \cite{heinonen}, Lemma 2.1, Theorem 2.3, Corollary 2.4)
\begin{lmm}\label{MWetLmm}
Let $X$ be a metric space. If $f_{\al}:X\ra\Real$ are $c$-Lipschitz maps, then so are $\inf_{\al}f_{\al}$ and $\sup_{\al} f_{\al}$.
\end{lmm}
\begin{proof}
By symmetry, it suffices to prove $\inf_{\al}f_{\al}$ is $c$-Lipschitz. The result follows from
\begin{align}
&\textstyle\left|\inf_{\al}f_{\al}(x)-\inf_{\al}f_{\al}(y)\right|=\max\left\{\inf_{\al}\sup_{\beta}[f_{\al}(x)-f_{\beta}(y)],-\inf_{\al}\sup_{\beta}[f_{\al}(x)-f_{\beta}(y)]\right\}\nn\\
&\textstyle~~~~\leq \max\left\{\inf_{\al}\sup_{\beta}|f_{\al}(x)-f_{\beta}(y)|,\inf_{\beta}\sup_{\al}|f_{\al}(x)-f_{\beta}(y)|\right\}=d_H\big(\{f_{\al}(x)\},\{f_{\al}(y)\}\big)\nn\\
&~~~~\leq\sup_{\al}\left|f_{\al}(x)-f_{\al}(y)\right|.\nn\qedhere
\end{align}
\end{proof}

\begin{thm}[\textcolor{OliveGreen}{McShane-Whitney extension theorem}]\label{MWET}
Let $X$ be a metric space. Every $c$-Lipschitz function $f:A\subset X\ra\Real$ extends to a $c$-Lipschitz function $F:X\ra\Real$.
\end{thm}
\begin{proof}
Observe the functions $\left\{f_a^{\pm}:X\ra\Real\right\}_{a\in A}$ given by $f_a^{\pm}(x)=f(a)\pm cd(x,a)$ are $c$-Lipschitz. Thus, by Lemma \ref{MWetLmm}, the functions $F(x)=\sup_{a\in A}f_a^-:X\ra\Real$ and $G(x)=\inf_{a\in A}f_a^+:X\ra\Real$ are $c$-Lipschitz extensions of $f$.
\end{proof}

\begin{crl}\label{MWETlmm}
Let $X$ be a metric space. Every $c$-Lipschitz map $f:A\subset X\ra\Real^n$ extends to a Lipschitz map $F:X\ra\Real^n$.
\end{crl}
\begin{proof}
Consider the $p$-norm on $\Real^n$, $1\leq p\leq\infty$. We have $f(x)=\big(f_1(x),...,f_n(x)\big)$, where $f_k:A\ra\Real$ has a $c$-Lipschitz extension $F_k:X\ra\Real$ by Theorem \ref{MWET}. Thus $F(x)=\big(F_1(x),...,F_n(x)\big)$ is a Lipschitz extension of $f$, with
\begin{align}
&\textstyle \|F(x)-F(y)\|=\big(\sum_k|F_k(x)-F_k(y)|^p\big)^{1/p}\leq n^{1\over p}cd(x,y),~~~~\txt{if}~~~~1\leq p<\infty,\nn\\
&\textstyle \|F(x)-F(y)\|=\max_k|F_k(x)-F_k(y)|\leq cd(x,y),~~~~\txt{if}~~~~p=\infty.\nn\qedhere
\end{align}
\end{proof}
{\flushleft Notice} from the proof of Corollary \ref{MWETlmm} that we can replace the target space $\Real$ in Theorem \ref{MWET} by any $\ell^\infty(\Gamma)$. By definition, the most general target space (replacing $\Real$) in Theorem \ref{MWET} is a 1-ALR (also called ``\ul{an injective metric space}'', or ``\ul{a hyperconvex metric space}'': See \cite{espinola-khamsi}).

\subsection{Bi-Lipschitz equivalence of simplices and cubes}\label{PrlMET2}
\begin{thm}[\textcolor{OliveGreen}{Brouwer's invariance of domain: \cite{tao2011}}]\label{BrInvDom}
Let $f:E\subset\Real^n\ra\Real^n$ be an injective continuous map. If $E$ is open, then so is $f(E)$.
\end{thm}
\begin{proof}
Let $E\subset\Real^n$ be open and consider an injective continuous map $f:E\subset\Real^n\ra\Real^n$. Since $E$ is open, for any $e\in E$, there exists an open neighborhood $N(e)$ of $e$ such that $N(e)\subset E$. Thus, $E=\bigcup_{e\in E}N(e)$, which implies $f(E)=\bigcup_{e\in E}f\big(N(e)\big)$. Since the union of open sets is open, it suffices to show that for any $e\in E$, $f\big(N(e)\big)$ is open. Moreover, because openness is invariant under rescaling and translations, it is enough to assume $E=B_1(0)$ and prove there exists an open neighborhood $N(0)$ of $0$ such that $f\big(N(0)\big)$ is open.

Since $f$ is injective, if there is an open neighborhood $N(f(0))$ of $f(0)$ such that $N(f(0))\subset f(B_1(0))$ (i.e., $f(0)$ is in the interior of $f(B_1(0))$), then $N(0):=f^{-1}\left[N(f(0))\right]$ is an open neighborhood of $0$ (by continuity of $f$) satisfying $f\big(N(0)\big)=N\big(f(0)\big)\subset f\big(B_1(0)\big)$. Thus, it suffices to prove $f(0)$ is an interior point of $f(B_1(0))$.

By Lemma \ref{BijVsHomeo}, the continuous bijection $f:B_1(0)\ra f(B_1(0))$ has a continuous inverse. By  Tietze's extension theorem, the continuous map $f^{-1}:f(B_1(0))\ra B_1(0)$ has a continuous extension $F':\Real^n\ra\Real^n$. Note that $F'(f(0))=0$, i.e., $f(0)$ is a zero of $F'$. For any map $F'':f(B_1(0))\ra\Real^n$, let
\bea
\|F'-F''\|:=\sup\limits_{y\in f(B_1(0))}\|F'(y)-F''(y)\|=\sup\limits_{x\in B_1(0)}\|F'(f(x))-F''(f(x))\|.\nn
\eea
Consider the map
\bea
h:B_1(0)\ra\Real^n,~~h(x):=x-F''(f(x))=F'(f(x))-F''(f(x)).\nn
\eea
If $\|F'-F''\|\leq 1$, then by Brouwer fixed point theorem, $h$ has a fixed point, and so $F''$ has a zero (i.e., there exists $x_0\in B_1(0)$ such that $F''(f(x_0))=0$).

Suppose on the contrary that $f(0)$ is not an interior point of $f(B_1(0))$. We will obtain a contradiction by finding a perturbation $F''$ of $F'$ such that $\|F'-F''\|\leq 1$ but $F''$ has no fixed point. Since $F'$ is continuous at $f(0)$, there is $\vep>0$ such that $F'\Big(B_{2\vep}\big(f(0)\big)\Big)\subset B_{0.1}(0)$, i.e.,
\bea
\|F'(y)\|=\|F'(y)-F'(f(0))\|\leq 0.1~~~~\txt{for}~~~~\|y-f(0)\|\leq 2\vep.\nn
\eea
Also, since $f(0)$ is not an interior point of $f(B_1(0))$, we know that $B_\vep(f(0))\cap f(B_1(0))^c\neq\emptyset$, i.e., there exists $c\in f(B_1(0))^c$ such that $\|c-f(0)\|<\vep$. Note that $F'(B_\vep(c))\subset B_{0.1}(0)$, since $\|y-c\|\leq\vep$ implies $\|y-f(0)\|\leq\|y-c\|+\|c-f(0)\|\leq 2\vep$, which in turn implies
\bea
\|F'(y)\|=\|F'(y)-F'(f(0))\|\leq 0.1~~~~\txt{for}~~~~\|y-c\|\leq \vep.\nn
\eea

For convenience, translate $f$ (i.e., replace $f(x)$ with $f(x)-c$) so that $c=0$. Then $f(B_1(0))$ avoids $0$ (formally $c$), $\|f(0)\|<\vep$, and $F'(B_\vep(0))\subset B_{0.1}(0)$, i.e.,
\bea
\label{BIDeq1}\|F'(y)\|=\|F'(y)-F'(f(0))\|\leq 0.1~~~~\txt{for}~~~~\|y\|\leq\vep.
\eea
Let $\Sigma:=\Sigma_1\cup\Sigma_2$, where $\Sigma_1:=f(B_1(0))\backslash N_\vep(0)=\{y\in f(B_1(0)):\|y\|\geq\vep\}$ and $\Sigma_2:=\del N_\vep(0)=\{y\in \Real^2:\|y\|=\vep\}$. Since $0=c\not\in\Sigma$, we get a continuous map $\Phi:f(B_1(0))\ra\Sigma$,
\bea
\label{BIDeq2}\textstyle\Phi(y):=\max\left({\vep\over\|y\|},1\right)y=
\left\{
  \begin{array}{ll}
    \vep{y\over\|y\|}, & y\in \Sigma_1=f\big(B_1(0)\big)\backslash N_\vep(0) \\
    y, & y\in f\big(B_1(0)\big)\cap N_\vep(0)
  \end{array}
\right\}.
\eea
Since $\Sigma_1$ is compact and $f(0)\not\in\Sigma_1$ implies $F'(y)\neq 0$ for all $y\in\Sigma_1$, there exists $0<\delta<0.1$ such that $\delta\leq \|F'(y)\|$ for all $y\in \Sigma_1$. By Weierstrass' theorem (on approximation of continuous functions with polynomials) there exists a polynomial $P:\Real^n\ra\Real^n$ such that
\bea
\label{BIDeq3}\|P(y)-F'(y)\|<\delta~~~~\txt{for all}~~~~y\in\Sigma,
\eea
where $P(y)\neq 0$ for all $y\in\Sigma_1$ (because $\delta\leq\|F'(y)\|$ for all $y\in \Sigma_1$). Since $\Sigma_2$ has \ul{measure zero}, we can choose $P$ such that $P(y)\neq 0$ for all $y\in\Sigma_2$ as well. [\ul{Reason-footnote}\footnote{Since $P$ is smooth and $\Sigma_2$ is a compact smooth surface of measure zero, the image $P(\Sigma_2)$ is also a compact smooth surface of measure zero (so, $P(\Sigma_2)$ has empty interior). Thus, if $0\in P(\Sigma_2)$, then by replacing $P(x)$ with $P_{\vep'}(x):=P(x)+\vep'\widehat{v}$ for an arbitrarily small $\vep'>0$ and a unit vector $\widehat{v}$, it follows that $0\not\in P_{\vep'}(\Sigma_2)$. We can perform this process and still have $0\not\in P_{\vep'}(\Sigma_1)$ since $P(\Sigma_1)$ is closed. Note that the set of zeros $Z(P):=P^{-1}(0)$ of $P$ is an algebraic set with empty interior (otherwise, if $\Int Z(P)\neq\emptyset$, then $P\eqv 0$).}]. Define the \ul{nonvanishing} function $F''=P\circ\Phi:f(B_1(0))\sr{\Phi}{\ral}\Sigma\sr{P}{\ral}\Real^n$.
By (\ref{BIDeq2}) and (\ref{BIDeq3}), we have
\bea
\|F'(y)-F''(y)\|<\delta~~~~\txt{for all}~~~~y\in f(B_1(0))~~\txt{such that}~~\|y\|>\vep.\nn
\eea
By (\ref{BIDeq1}) and (\ref{BIDeq2}), $\|F'(y)\|,\|F'(\phi(y))\|\leq0.1$ for $\|y\|\leq\vep$, and so by (\ref{BIDeq3}),
\bea
&&\|F'(y)-F''(y)\|=\|F'(y)-P(\Phi(y))\|\leq \|F'(y)-F'(\Phi(y))\|+\|F'(\Phi(y))-P(\Phi(y))\|\nn\\
&&~~~~<0.2+\delta<0.3<1,~~~~\txt{if}~~~~\|y\|\leq\vep.\nn
\eea
Hence, we have a contradiction, since $\|F'(y)-F''(y)\|<1$ for all $y\in f(B_1(0))$,
although ~$F''(y)=P(\Phi(y))\neq 0$~ for all ~$y\in f(B_1(0))$.
\end{proof}

\begin{crl}\label{BIDcrl1}
Let $K\subset\Real^n$ be compact. If $f:K\ra\Real^n$ is an injective continuous map, then (i) $f:K\ra f(K)$ is a homeomorphism and (ii) interior (resp. boundary) points of $K$ map to interior (resp. boundary) points of $f(K)$.
\end{crl}
\begin{proof}
By Lemma \ref{BijVsHomeo}, $f:K\ra f(K)$ is a homeomorphism. Hence, the correspondences interior$\sr{f}{\longleftrightarrow}$interior and boundary$\sr{f}{\longleftrightarrow}$boundary follow from Theorem \ref{BrInvDom}.
\end{proof}

\begin{prp}[\textcolor{OliveGreen}{Bi-Lipschitz equivalence of convex set and ball}]\label{BiLipHomP}
Every \ul{compact convex} set $K\subset\Real^n$ with \ul{nonempty interior} ($K^o\neq\emptyset$) is bi-Lipschitz equivalent to the closed unit ball $\ol{B}_1(0)\subset\Real^n$.
\end{prp}
\begin{proof}
Let $K\subset\Real^n$ be a compact convex set with $\Int(K)\neq\emptyset$. Since a translation is a bi-Lipschitz equivalence, we can assume $0\in K^o$. Define a map $F:\Real^n\ra \ol{B}_1(0)$ by
\bea
\textstyle F(x):=p_K(x){x\over\|x\|}=p_K(x)\hat{x}=p_K(\hat{x})x,~~~~\hat{x}:={x\over\|x\|},\nn
\eea
where $p_K(x)$ is the Minkowski function of $K$ (Lemma \ref{MinkGauge}) defined by
\bea
\textstyle p_K(x):=\inf\left\{c>0:{x\over c}\in K\right\}=\inf\left\{c>0:x\in cK\right\}.\nn
\eea
Note that $p_K(x)\leq 1$ for all $x\in K$. Also, $p_{\ol{B}_1(0)}(x)=\|x\|$. Recall that $p_K$ is a sublinear functional, i.e., $p_K(x+y)\leq p_K(x)+p_K(y)$ and $p_K(\ld x)=\ld p_K(x)$ for any scalar $\ld\geq 0$. Thus,
\bea
\textstyle\big|p_K(x)-p_K(y)\big|\leq p_K(x-y)\leq\|p_K\|\|x-y\|,~~~~\|p_K\|:=\sup_{\|x\|=1}p_K(x).\nn
\eea
We need to show $F|_K:K\ra\ol{B}_1(0)$ is a bi-Lipschitz homeomorphism. $F$ is a continuous bijection with inverse $F^{-1}(y)={1\over p_K(y)}\|y\|y={y\over p_K(\hat{y})}$, and so $F|_K:K\ra\ol{B}_1(0)$ is a homeomorphism by Corollary \ref{BIDcrl1}. Also, $F|_K$ is Lipschitz as a product of bounded Lipschitz maps, and $F|_K^{-1}$ is Lipschitz for the same reason.
\end{proof}

\begin{crl}[\textcolor{OliveGreen}{Bi-Lipschitz equivalence of simplex and cube}]\label{BiLipHomC}
The $n$-simplex
\bea
\Delta^n=\Conv(e_0,e_1,...,e_n)\subset\Real^{n+1}\nn
\eea
(viewed, under a suitable affine transformation, as a subset of $\Real^n$) is bi-Lipschitz homeomorphic to the $n$-cube $I^n=[0,1]^n\subset\Real^n$.

Moreover, if $\phi:\Delta^n\ra I^n$ is a bi-Lipschitz homeomorphism, then we have the correspondences $\Int(\Delta^n)\sr{\phi}{\longleftrightarrow}\Int(I^n)$ and $\del\Delta^n\sr{\phi}{\longleftrightarrow}\del I^n$.
\end{crl}

\subsection{Paths, quasiconvexity, and gluing lemmas}\label{PrlMET3}
\begin{dfn}[\textcolor{blue}{\index{Path}{Path}, \index{Parametrization}{Parametrization}, \index{Length}{Length} of a path, \index{Rectifiable path}{Rectifiable path}, \index{ML parametrization}{Minimum length (ML) parametrization}, \index{Natural parametrization}{Natural parametrization}, \index{Constant speed path}{Constant speed path}, \index{Subconstant speed path}{Subconstant speed path}}]\label{PathLenDfn}
A \ul{path} in a space $X$ is a continuous map $\gamma:[0,1]\ra X$. Given paths $\gamma,\eta:[0,1]\ra X$, $\eta$ is a \ul{(re)parametrization} of $\gamma$, written $\eta\in\txt{Para}(\gamma)$, if
\bea
\im\eta=\im\gamma,~~\eta(0)=\gamma(0),~~\eta(1)=\gamma(1).\nn
\eea

Let $X$ be a metric space and $\gamma:[0,1]\ra X$ a path (or any map). The \ul{length} of $\gamma$ is
\begin{equation*}
l(\gamma)~:=~\sup\big\{l_P(\gamma):P\subset [0,1]~~\txt{a finite set with $0,1\in P$}\big\},
\end{equation*}
where $l_P(\gamma):=\sum_{i=1}^kd(\gamma(t_{i-1}),\gamma(t_i))$ is the length of $\gamma$ over $P:=\{0=t_0<t_1<\cdots<t_k=1\}$. If $l(\gamma)<\infty$, we say $\gamma$ is \ul{rectifiable}. A parametrization $\eta$ of $\gamma$ is a \ul{minimum length (ML) parametrization} if
\bea
l(\eta)\leq l(\eta')~~~~\txt{for all}~~\eta'\in \txt{Para}(\gamma).\nn
\eea
A continuous map $\gamma:[a,b]\ra X$ (for $a\leq b$ in $\Real$) is a \ul{natural parametrization} if
\bea
l(\gamma([t,t']))=|t-t'|~~~~\txt{for all}~~t,t'\in[a,b].\nn
\eea
A path $\gamma:[0,1]\ra X$ has \ul{constant speed} $c\geq 0$ if
\bea
l(\gamma([t,t']))=c|t-t'|~~~~\txt{for all}~~t,t'\in[0,1].\nn
\eea
A path $\gamma:[0,1]\ra X$ has \ul{subconstant speed} $c\geq0$ if
\bea
l(\gamma([t,t']))\leq c|t-t'|~~~~\txt{for all}~~t,t'\in[0,1].\nn
\eea
\end{dfn}

In the above definition, $l(\gamma)$ depends on the way $\gamma$ is parameterized, i.e., if $\gamma,\eta:[0,1]\ra X$ are parametrizations of the same path, then we can have $l(\gamma)\neq l(\eta)$.

\begin{lmm}[\textcolor{OliveGreen}{\cite[Proposition 2.5.9, p.46]{BBI}}]\label{ConstSpeedPar}
Let $X$ be a metric space. If $\gamma:[0,1]\ra X$ is a rectifiable path, then with $l:=l(\gamma)$, there exists a nondecreasing continuous map $\vphi:[0,1]\ra [0,l]$ and a natural parametrization $\ol{\gamma}:[0,l]\ra X$ such that $\gamma=\ol{\gamma}\circ\vphi:[0,1]\sr{\vphi}{\ral}[0,l]\sr{\ol{\gamma}}{\ral}X$.
\end{lmm}
\begin{proof}
Define ~$\vphi:[0,1]\ra[0,l],~t\mapsto l\big(\gamma([0,t])\big)$,~ which is continuous and nondecreasing. Also define
\bea
\ol{\gamma}:[0,l]\ra X,~{s}\mapsto\gamma(t)\in\gamma\left(\vphi^{-1}({s})\right),~~\txt{i.e.,}~~t\in\vphi^{-1}({s}),~~\txt{or}~~\vphi(t)={s}.\nn
\eea
$\ol{\gamma}$ is well defined because for any $t,t'\in\vphi^{-1}({s})$, $t\leq t'$,
\bea
&&d(\gamma(t),\gamma(t'))\leq l(\gamma([t,t']))=l(\gamma([0,t']))-l(\gamma([0,t]))=\vphi(t')-\vphi(t)={s}-{s}=0,\nn\\
&&~~\Ra~~d(\gamma(t),\gamma(t'))=0,~~\Ra~~\gamma(t')=\gamma(t).\nn
\eea
Also, for any $t\in[0,1]$, we have $\gamma(t)=\ol{\gamma}(\vphi(t))$ since $t\in\vphi^{-1}(\vphi(t))$. Finally,
\bea
&&l(\ol{\gamma}([{s},{s}']))=l(\gamma([t_{s},t_{s}']))=\big|l(\gamma([0,t_{s}]))-l(\gamma([0,t_{{s}'}]))\big|\nn\\
&&~~~~=\big|\vphi(t_{s})-\vphi(t_{{s}'})\big|=|{s}-{s}'|,\nn
\eea
where ~$t_s\in\vphi^{-1}(s)$, $t_{s'}\in\vphi^{-1}(s')$.
\end{proof}

\begin{crl}[\textcolor{OliveGreen}{Constant speed parametrization of a rectifiable path}]\label{ConstSpeedRP}
Let $X$ be a metric space. If $\gamma:[0,1]\ra X$ is a rectifiable path, there exists a path $\eta:[0,1]\ra X$ such that $\eta([0,1])=\gamma\big([0,1]\big)$, $\eta(0)=\gamma(0)$, $\eta(1)=\gamma(1)$, and
\bea
\textstyle d(\eta(t),\eta(t'))\leq l(\eta([t,t']))=l(\gamma)|t-t'|,~~\txt{for all}~~t,t'\in[0,1].\nn
\eea
\end{crl}
\begin{proof}
By Lemma \ref{ConstSpeedPar}, $\gamma=\ol{\gamma}\circ\vphi:[0,1]\sr{\vphi}{\ral}[0,l(\gamma)]\sr{\ol{\gamma}}{\ral}X$, for a nondecreasing continuous map $\vphi:[0,1]\ra [0,l(\gamma)]$ and a natural parametrization $\ol{\gamma}:[0,l(\gamma)]\ra X$. Thus, with {\small$\psi:[0,1]\ra[0,c]$, $t\mapsto ct$}, where $c:=l(\gamma)$, we get the new parametrization {\small$\eta:=\ol{\gamma}\circ\psi:[0,1]\sr{\psi}{\ral}[0,c]\sr{\ol{\gamma}}{\ral}X$}, which satisfies ~$l(\eta([t,t']))=c|t-t'|$ for all $t,t'\in[0,1]$.
\end{proof}

\begin{dfn}[\textcolor{blue}{\index{Quasigeodesic}{Quasigeodesic}, \index{Standard parametrization}{Standard parametrization}, \index{Quasiconvex space}{Quasiconvex space}, \index{Geodesic!}{Geodesic}, \index{Geodesic! space}{Geodesic space}, \index{Uniquely geodesic space}{Uniquely geodesic space}}]\label{QsiGeoDfn}
Let $X$ be a metric space and $\ld\geq1$. A path {\footnotesize $\gamma:[0,1]\ra X$} is a \ul{$\ld$-quasigeodesic} if $\gamma$ has subconstant speed $c=\ld d(\gamma(0),\gamma(1))$, i.e.,
\bea
\textstyle l(\gamma([t,t']))\leq\ld d(\gamma(0),\gamma(1))|t-t'|,~~~~\txt{for all}~~~~t,t'\in[0,1].\nn
\eea

We say $X$ is a \ul{$\ld$-quasiconvex space} if for every $x,y\in X$, there exists a $\ld$-quasigeodesic $\gamma:[0,1]\ra X$ from $x$ to $y$, i.e., such that $\gamma(0)=x$, $\gamma(1)=y$.

A $1$-quasigeodesic is called a \ul{geodesic}, and similarly, a $1$-quasiconvex space is called a \ul{geodesic space}. A geodesic space $X$ is a \ul{uniquely geodesic space} if for every $x,y\in X$, there exists only one geodesic from $x$ to $y$.
\end{dfn}

Note that a $\ld$-quasigeodesic is also called a \ul{$\ld$-quasiconvex path}, \cite[p.205]{hakobyan-herron2008}. In \cite[p.317]{tyson-wu2005}, a quasigeodesic is differently defined to be a path that is a bi-Lipschitz imbedding. Injectivity of the path is not required in our definition. An alternative definition of a geodesic in terms of paths that are parameterized by arc length can be found in \cite[Definition 2.2.1, p.56]{papado2014}.

\begin{lmm}[\textcolor{OliveGreen}{Characterization and Sufficient condition for quasigeodesics}]\label{QgeodCharLmm}
Let $X$ be a metric space, $\gamma:[0,1]\ra X$ a path, and $\ld,\ld_1,...,\ld_k\geq 1$. Then the following are true:
\bit[leftmargin=0.8cm]
\item[(i)] $\gamma$ is a $\ld$-quasigeodesic $\iff$ $d(\gamma(t),\gamma(t'))\leq\ld d(\gamma(0),\gamma(1))|t-t'|$, for all $t,t'\in[0,1]$.
\item[(ii)] If $[0,1]=\bigcup_{i=1}^k[a_i,b_i]$ such that $b_i=a_{i+1}$, and $\gamma_i:=\gamma|_{[a_i,b_i]}$ is a $\ld_i$-quasigeodesic (for all $i=1,...,k$), then $\gamma$ is a $\ld$-quasigeodesic if ~$\ld:=\max\limits_i\ld_i{d(\gamma(a_i),\gamma(b_i))\over d(\gamma(0),\gamma(1))}$.
\eit
\end{lmm}

\begin{lmm}[\textcolor{OliveGreen}{Characterization of geodesics: See also \cite[Sec. 2.2, pp 56-60]{papado2014}}]\label{GeodCharLmm}
Let $X$ be a metric space and $\gamma:[0,1]\ra X$ a path. Then (i) $\gamma$ is a geodesic $\iff$ (ii) $d(\gamma(t),\gamma(t'))\leq d(\gamma(0),\gamma(1))|t-t'|$ for all $t,t'\in[0,1]$, $\iff$ (iii) up to a reparametrization, $d(\gamma(t),\gamma(t'))=d(\gamma(0),\gamma(1))|t-t'|$ for all $t,t'\in[0,1]$.
\end{lmm}
\begin{proof}
{\flushleft \ul{(i)$\Ra$(ii)}:} This is immediate, since ~$d(\gamma(t),\gamma(t'))\leq l(\gamma|_{[t,t']})$, for all $t,t'\in[0,1]$.
{\flushleft \ul{(ii)$\Ra$(i),(iii)}:} (ii) implies $d(\gamma(t),\gamma(t'))\leq d(\gamma(0),\gamma(1))|t-t'|$ for all $t,t'\in[0,1]$. By the definition of $l(\gamma)$, $l(\gamma)\leq d(\gamma(0),\gamma(1))\leq l(\gamma)$, i.e., $l(\gamma)=d(\gamma(0),\gamma(1))$. Also, observe that $d(\gamma(0),\gamma(1))\leq l_P(\gamma)\leq l(\gamma)=d(\gamma(0),\gamma(1))$ for any finite partition $P\subset[0,1]$. That is, $l_P(\gamma)=l_{P'}(\gamma)$ for any two finite partitions $P,P'\subset [0,1]$. In particular, if $P:=\{0,t,t',1\}$ and $P':=P\cup Q=\{0\}\cup Q\cup\{1\}$ for any finite partition $Q\subset[t,t']$, then
\bea
l_P(\gamma)=l_{P'}(\gamma)~~\Ra~~d(\gamma(t),\gamma(t'))=l_Q(\gamma|_{[t,t']})=l(\gamma|_{[t,t']}).\nn
\eea
Hence, up to a reparametrization (by Corollary \ref{ConstSpeedRP}), ~$d(\gamma(t),\gamma(t'))=d(\gamma(0),\gamma(1))|t-t'|$ ~for all ~$t,t'\in[0,1]$.
{\flushleft \ul{(iii)$\Ra$(ii)}:} This is again immediate.
\end{proof}

\begin{dfn}[\textcolor{blue}{\index{Minimum-length path}{Minimum-length path}}]
Let $X$ be a metric space, $x,y\in X$, and $\P_{x,y}(X):=\{\txt{paths}~\gamma:[0,1]\ra X,~\gamma(0)=x,\gamma(1)=y\}\subset \C\big([0,1],X\big)$. Then a path $\gamma:[0,1]\ra X$ is a \ul{minimum-length path} (or path of minimum length) if
\bea
l(\gamma)=\inf\big\{l(\eta)~|~\eta\in\P_{\gamma(0),\gamma(1)}(X)\big\}.\nn
\eea
\end{dfn}
It follows immediately from Lemma \ref{GeodCharLmm} that every geodesic is a minimum-length path.

\begin{dfn}[\textcolor{blue}{Strictly convex normed space}]\label{StrictConvNS}
A normed space $X$ is strictly convex if for any $x,y\in X$, we have $\|x+y\|=\|x\|+\|y\|$ $\iff$ $y=\ld x$ for some scalar $\ld\geq 0$.
\end{dfn}
\begin{lmm}
A normed space $X$ is uniquely geodesic $\iff$ strictly convex.
\end{lmm}
\begin{proof}
($\Ra$) Assume $X$ is uniquely geodesic. Then for any $a,b\in X$, the only geodesic from $a$ to $b$ is $\gamma_{a,b}:[0,1]\ra X$, $\gamma_{a,b}(t):=(1-t)a+tb$. Thus, for any $x,y\in X$, if $\|x+y\|=\|x-(-y)\|=\|x-0\|+\|0-y\|$, then for some $t\in[0,1]$ we have $\gamma_{x,-y}(t)=0=(1-t)x+t(-y)$, i.e., $y={1-t\over t}x$.

($\La$): Assume $X$ is strictly convex. Let $\gamma:[0,1]\ra X$ be a geodesic from $x$ to $y$. Then for any $z\in\gamma([0,1])$, we have $\|x-y\|=\|x-z+z-y\|=\|x-z\|+\|z-y\|$. Thus, by strict convexity, we get $z-y=\ld(x-z)$ for some $\ld\geq 0$, i.e., $z={\ld\over 1+\ld}x+{1\over 1+\ld}y=(1-t)x+ty$, where $0\leq t:=1/(1+\ld)\leq 1$. This shows $\gamma$ is uniquely given by $\gamma(t)=(1-t)x+ty$.
\end{proof}

\begin{lmm}[\textcolor{OliveGreen}{\index{Gluing lemma I}{Gluing lemma I}}]\label{GluLmmI}
Let $X$ be a $\lambda$-quasiconvex space, and $\{C_j:j=1,...,n\}$ a finite cover of $X$ by closed sets $C_j$. If $f:X\ra Y$ is a map such that each restriction $f_j=f|_{C_j}:C_j\ra Y$ is $c$-Lipschitz, then $f$ is $\lambda c$-Lipschitz.
\end{lmm}
\begin{proof}
We know $f$ is continuous because for any closed set $C'\subset Y$, we have the finite union of closed sets $f^{-1}(C')=\bigcup C_j\cap f^{-1}(C')=\bigcup f_j^{-1}(C')$.

Given $x,x'\in X$, let $\gamma:[0,1]\ra X$ be a $\ld$-quasigeodesic from $x$ to $x'$. Then, for any partition $P=\{0=t_0<t_1<\cdots<t_k=1\}$ of $[0,1]$, we have
\begin{align}
&\textstyle d\big(f(x),f(x')\big)\leq \sum_{i=1}^kd\big(f(\gamma(t_{i-1})),f(\gamma(t_i))\big)=\sum_{i=1}^kd\big(f_{j_{i-1}}(\gamma(t_{i-1})),f_{j_i}(\gamma(t_i))\big)\nn\\
&\textstyle~~~~\sr{(s)}{\lesssim}c\sum_{i=1}^{k'}d\big(\gamma(t'_{i-1}),\gamma(t'_i)\big)\leq c~\txt{length}(\gamma)\leq c\ld d(x,x'),\nn
\end{align}
where at step (s), we use a sufficiently fine partition $P=\{0=t_0<t_1<\cdots<t_k=1\}$, and  $P'=\{0=t'_0<t'_1<\cdots<t'_{k'}=1\}$ is some refinement of $P$ based on the following discussion. (\textcolor{blue}{footnote}\footnote{To avoid/remove the approximation sign at step (s), simply choose a partition $P\subset[0,1]$ containing the intersection points $\bigcup_j\gamma\cap \del C_j$ in the sense that $\bigcup_j\gamma\cap \del C_j\subset\gamma(P)$. Recall that (i) the cover $\{C_j\}_{j=1}^n$ is finite and that (ii) we can choose $\gamma$ to have only finitely many intersection points $\bigcup_j\gamma\cap \del C_j$ with the boundaries of the cover.})

To justify step (s), fix $i\in\{1,...,k\}$. Observe that if $\gamma\cap C_{j_{i-1}}\cap C_{j_i}=\emptyset$, then we can approximate $\gamma$ with a disjoint union of segments of $\gamma$ as ~$\gamma\approx\gamma|_{[0,t_{i-1}]}\sqcup\gamma|_{[t_i,1]}$. (Note that since $k$ is finite, by induction, we can approximate $\gamma$ with a disjoint union of at most $k+1$ segments of $\gamma$). On the other hand, if $z\in\gamma\cap C_{j_{i-1}}\cap C_{j_i}$, then
\begin{align}
&d\big(f_{j_{i-1}}(\gamma(t_{i-1})),f_{j_i}(\gamma(t_i))\big)\leq d\big(f_{j_{i-1}}(\gamma(t_{i-1})),f(z)\big)+d\big(f(z),f_{j_i}(\gamma(t_i))\big)\nn\\
&~~~~\leq c\big[d\big(\gamma(t_{i-1}),z\big)+d\big(z,\gamma(t_i)\big)\big].\nn\qedhere
\end{align}
\end{proof}

\begin{dfn}[\textcolor{blue}{\index{Locally! Lipschitz map}{Locally Lipschitz map}}]\label{LocLipMap}
Let $c\geq 0$. A map of metric spaces $f:X\ra Y$ is locally $c$-Lipschitz if for each $x\in X$, there exists a ball $B_{r_x}(x)$, $r_x=r_{x,f}>0$, such that
\bea
d(f(x),f(z))\leq cd(x,z)~~~~\txt{for all}~~~~z\in B_{r_x}(x).\nn
\eea
\end{dfn}

\begin{lmm}[\textcolor{OliveGreen}{\index{Gluing lemma II}{Gluing lemma II}}]\label{GluLmmII}
Let $X$ be a $\ld$-quasigeodesic space. If a map $f:X\ra Y$ is locally $c$-Lipschitz, then it is $\ld c$-Lipschitz. (Conversely, we know a $c$-Lipschitz map is locally $c$-Lipschitz.)
\end{lmm}
\begin{proof}
Fix $x,x'\in X$. Let $\gamma:[0,1]\ra X$ be a $\ld$-quasigeodesic from $x$ to $x'$. Define $C:=cd(x,x')$. Then $f\circ\gamma:[0,1]\ra Y$ is locally $\ld C$-Lipschitz on $[0,1]$, since for any $t\in[0,1]$, there exists a ball $B_{r_{\gamma(t)}}\big(\gamma(t)\big)$, $r_{\gamma(t)}>0$, such that
\bea
d\big(f(\gamma(t)),f(z)\big)\leq cd\big(\gamma(t),z\big)~~~~\txt{for all}~~~~z\in B_{r_{\gamma(t)}}\big(\gamma(t)\big),\nn
\eea
and so we get the neighborhood $U:=\gamma^{-1}\big(B_{r_{\gamma(t)}}\big(\gamma(t)\big)\big)$ of $t$ in $[0,1]$ satisfying
\bea
d\big(f(\gamma(t)),f(\gamma(s))\big)\leq cd\big(\gamma(t),\gamma(s)\big)\leq C\ld|t-s|,~~~~\txt{for all}~~s\in U.\nn
\eea

Since $[0,1]$ is compact, for any $t,t'\in [0,1]$ we can choose a  partition $P=\{t=t_0<t_1<\cdots<t_k=t'\}$ of $[t,t']$ such that for some refinement $Q=\{t=s_0<s_1<\cdots<s_{l}=t'\}$ of $P$,
\vspace{-0.2cm}
\begin{align}
\textstyle d\big(f(\gamma(t)),f(\gamma(t'))\big)\leq \sum_{i=1}^kd\big(f(\gamma(t_{i-1})),f(\gamma(t_i))\big)\sr{(a)}{\leq} \ld C\sum_{j=1}^l|s_{j-1}-s_j|=\ld C|t-t'|,\nn
\end{align}
where at step (a), we choose $P$ fine enough so that for each $i\in\{1,...,k\}$ there exists $s\in[t_{i-1},t_i]$ satisfying
\bea
d\big(f(\gamma(t_{i-1})),f(\gamma(s))\big)\leq\ld C|t_{i-1}-s|,~~~~d\big(f(\gamma(s)),f(\gamma(t_i))\big)\leq\ld C|s-t_i|.\nn
\eea
This shows $f\circ\gamma$ is $\ld C$-Lipschitz. Thus, $d\big(f(x),f(x')\big)\leq \ld C=\ld cd(x,x')$.
\end{proof}

The following concepts are based on \cite{sims,espinola.etal}.

\begin{dfn}[\textcolor{blue}{Standard constant curvature spaces}]
A standard 2-dimensional space of constant (Riemannian) curvature $\kappa$ is a metric space (based on $\Real^3$-subsets) of the form
{\footnotesize
\bea
M_\kappa^2:=\left\{
                \begin{array}{ll}
                  \big(S^2,d_+\big),~~d_+(x,y)={1\over\sqrt{\kappa}}\cos^{-1}\big(\langle x,y\rangle_+\big),~~\langle x,y\rangle_+:=x_1y_1+x_2y_2+x_3y_3, ~&~ \txt{if}~~\kappa>0 \\\\
\big(\Real^2,d_0\big),~~d_0(x,y)=\sqrt{\langle x-y,x-y\rangle},~~\langle x,y\rangle_0:=x_1y_1+x_2y_2+x_3y_3, ~&~ \txt{if}~~\kappa=0 \\\\
\big(H^2,d_-\big),~~d_-(x,y)={1\over\sqrt{-\kappa}}\cosh^{-1}\big(-\langle x,y\rangle_-\big),~~\langle x,y\rangle_-:=x_1y_1+x_2y_2-x_3y_3, ~&~ \txt{if}~~\kappa<0
             \end{array}
           \right\},\nn
\eea}
where
\bea
&&S^2:=\left\{x\in\Real^3~|~\langle x,x\rangle_+:=x_1^2+x_2^2+x_3^2=1\right\},\nn\\
&&H^2:=\left\{x\in\Real^3~|~\langle x,x\rangle_-:=x_1^2+x_2^2-x_3^2=-1,~~x_3>0\right\}.\nn
\eea
\end{dfn}

\begin{dfn}[\textcolor{blue}{\index{Geodesic! triangle}{Geodesic triangle} (in a geodesic space $X$), \index{Isometry of triangles}{Isometry of triangles}, Isometric geodesic triangles}]
Let $x,y,z\in X$. We denote by $\gamma[x,y]$ a \ul{geodesic (segment)} $\gamma:[0,1]\ra X$ from $x$ to $y$ (i.e., with $\gamma(0)=x$, $\gamma(1)=y$). Consider any geodesic segments $\gamma_1[x,y]$, $\gamma_2[y,z]$, $\gamma_3[z,x]$. A loop $\Delta:[0,1]\ra X$, $\Delta(0)=\Delta(1)$ is called a \ul{geodesic triangle} in $X$ with \ul{vertices} $x,y,z$ and \ul{sides} $\gamma_1[x,y]$, $\gamma_2[y,z]$, $\gamma_3[z,x]$ if $\Delta=\gamma_1\cdot\gamma_2\cdot\gamma_3:[0,1]\ra X$, where
\bea
&&(\gamma_1\cdot\gamma_2\cdot\gamma_3)|_{[0,1/3]}(t):=\gamma_1(3t-0),\nn\\
&&(\gamma_1\cdot\gamma_2\cdot\gamma_3)|_{[1/3,2/3]}(t):=\gamma_2(3t-1),\nn\\
&&(\gamma_1\cdot\gamma_2\cdot\gamma_3)|_{[2/3,1]}(t):=\gamma_3(3t-2).\nn
\eea
That is, $\Delta(0)=\Delta(1)=x$, $\Delta(1/3)=y$, $\Delta(2/3)=z$. We will write the triangle as a union
\bea
\Delta[x,y,z]:=\gamma_1[x,y]\cup\gamma_2[y,z]\cup\gamma_3[z,x].\nn
\eea

Given geodesic spaces $X,Y$, and two geodesic triangles $\Delta_X[x_1,x_2,x_3]\subset X$ and $\Delta_Y[y_1,y_2,y_3]\subset Y$, a map $\vphi:\Delta_X[x_1,x_2,x_3]\ra\Delta_Y[y_1,y_2,y_3]$ is an \ul{isometry of triangles}, making $\Delta_X[x_1,x_2,x_3]$ and $\Delta_Y[y_1,y_2,y_3]$ \ul{isometric geodesic triangles} ( written $\Delta_X[x_1,x_2,x_3]\sr{\vphi}{\cong}\Delta_Y[y_1,y_2,y_3]$ ) if it is an isometry satisfying $\vphi(x_i)=y_i$ for all $i=1,2,3$.
\end{dfn}
\begin{note*}
If $X$ is a uniquely geodesic space, then it is clear that the triangle $\Delta[x,y,z]$ is completely determined by its vertices $x,y,z$.
\end{note*}

\begin{dfn}[\textcolor{blue}{\index{Comparison triangle}{Comparison triangle} for a geodesic triangle, \index{Comparison point}{Comparison point}}]
A \ul{comparison triangle} for a geodesic triangle $\Delta[x,y,z]\subset X$ is a geodesic triangle $\ol{\Delta}[\ol{x},\ol{y},\ol{x}]\subset M_\kappa^2$ for which there exists an isometry of triangles $\ol{\Delta}[\ol{x},\ol{y},\ol{x}]\sr{\vphi}{\cong}\Delta[x,y,z]$.

The \ul{$\vphi$-comparison point} for a point $p\in\Delta[x,y,z]$ is its preimage {\small$\ol{p}:=\vphi^{-1}(p)\in\ol{\Delta}[\ol{x},\ol{y},\ol{x}]$}. Equivalently, the $\vphi$-comparison point $\ol{p}\in\ol{\Delta}[\ol{x},\ol{y},\ol{x}]$ for $p\in\Delta[x,y,z]$ is defined as follows:
\bit[leftmargin=0.6cm]
\item If {\small$p\in\gamma_1[x,y]$}, then {\small$\ol{p}\in\ol{\gamma}_1[\ol{x},\ol{y}]$} is given by {\small$d(\ol{x},\ol{p})=d(x,p)$}.
\item If {\small$p\in\gamma_2[y,z]$}, then {\small$\ol{p}\in\ol{\gamma}_2[\ol{y},\ol{z}]$} is given by {\small$d(\ol{y},\ol{p})=d(y,p)$}.
\item If {\small$p\in\gamma_3[z,x]$}, then {\small$\ol{p}\in \ol{\gamma}_3[\ol{z},\ol{x}]$} is given by {\small$d(\ol{z},\ol{p})=d(z,p)$}.
\eit
\end{dfn}

By construction, the comparison triangle (when it exists) is unique up to isometry. Also, because $\Real^2$ and $H^2$ are unbounded, comparison triangles always exist in $M_\kappa^2$ for $\kappa\leq0$. For $\kappa>0$, because $S^2$ is bounded, a comparison triangle for $\Delta[x,y,z]\subset X$ exists only if the perimeter $l(\Delta[x,y,z]):=d(x,y)+d(y,z)+d(z,x)$ is sufficiently small.

\begin{dfn}[\textcolor{blue}{\index{CAT($\kappa$) inequality}{CAT($\kappa$) inequality}}]
A geodesic triangle $\Delta\subset X$ satisfies the $CAT(\kappa)$ inequality if for all points $p,q\in\Delta$, we have
\bea
d(p,q)\leq d\left(\ol{p},\ol{q}\right),~~~~\txt{for all comparison points}~~\ol{p},\ol{q}\in\ol{\Delta}\subset M_\kappa^2.\nn
\eea
\end{dfn}

\begin{dfn}[\textcolor{blue}{\index{CAT($\kappa$) space}{CAT($\kappa$) space}, \index{Hadamard space}{Hadamard space}}]
A \ul{CAT($\kappa$) space} is a geodesic space in which every geodesic triangle satisfies the $\CAT(\kappa)$ inequality. A complete $\CAT(0)$ space is called a \ul{Hadamard space}.
\end{dfn}

\begin{dfn}[\textcolor{blue}{\index{Geodesic! convexity}{Geodesic convexity}}]
A subset $K\subset X$ of a metric space $X$ is geodesically convex if given any $x,y\in K$, every geodesic segment $\gamma[x,y]$ between $x$ and $y$ lies in $K$.
\end{dfn}

\begin{dfn}[\textcolor{blue}{\index{Convex! metric}{Convex metric}}] If $(X,d)$ is a geodesic space, the metric $d:X\times X\ra\Real$ is convex if for every $x\in X$ and any geodesic segment $\gamma:[0,1]\ra X$, the function $f:[0,1]\ra\Real$ given by $f(t):=d(x,\gamma(t))$ is convex.~ Equivalently,
\bea
d(x,\gamma(t))=d(x,\gamma((1-t)0+t1))\leq (1-t)d(x,\gamma(0))+td(x,\gamma(1)),~~~~\txt{for all}~~~~0\leq t\leq1.\nn
\eea
(For the equivalence, recall that every subsegment of a geodesic segment is itself a geodesic segment).
\end{dfn}

\begin{dfn}[\textcolor{blue}{\index{Busemann convexity}{Busemann convexity}}]
A geodesic space $(X,d)$ is Busemann convex if for for every pair of geodesic segments $\gamma_1:[0,1]\ra X$, $\gamma_2:[0,1]\ra X$, the function $g:[0,1]\ra\Real$ given by $g(t)=d(\gamma_1(t),\gamma_2(t))$ is convex. Equivalently (since every subsegment of a geodesic segment is itself a geodesic segment),
\bea
d(\gamma_1(t),\gamma_2(t))\leq (1-t)d(\gamma_1(0),\gamma_2(0))+td(\gamma_1(1),\gamma_2(1)),~~~~\txt{for all}~~~~0\leq t\leq1.\nn
\eea
\end{dfn}

\begin{note}
If a geodesic space $X$ is Busemann-convex, then it is a uniquely geodesic space and its metric is convex. Indeed, if $\gamma_{xy},\gamma_{xy}':[0,1]\ra X$ are two geodesics from $x$ to $y$, then the nonnegative convex function $f(t)=d\big(\gamma_{xy}(t),\gamma_{xy}'(t)\big)$ satisfies $f(0)=f(1)=0$, and so we have $f(t)=0$ for all $t$. To see that the metric of $X$ is convex, observe that if $\gamma:[0,1]\ra X$ is a geodesic and $x\in X$, then we have a convex function $f(t):=d(x,\gamma(t))=f(\gamma_x(t),\gamma(t))$, where $\gamma_x:[0,1]\ra X,~t\mapsto x$ is the constant geodesic.
\end{note}

\begin{dfn}[\textcolor{blue}{\index{Convex! bicombing}{Convex bicombing}}]
Let $X$ be a geodesic space. Consider the set-valued map $\Gamma:X\times X\ra 2^{C([0,1],X)}$ given by
\bea
\Gamma(x,y):=\big\{\txt{geodesics $\gamma[x,y]$ from $x$ to $y$}~\big|~ x,y\in X\big\}.\nn
\eea
A convex bicombing on $X$ is a selection of $\Gamma$,
\bea
\gamma:X\times X\ra C([0,1],X),~~(x,y)\mapsto \gamma_{xy}=\gamma[x,y]~\in~\Gamma(x,y),\nn
\eea
such that for all $x,x',y,y'\in X$, the following are true.
\bit[leftmargin=0.8cm]
\item[(i)] $\gamma_{xy}(t)=\gamma_{yx}(1-t)$ ~for all ~$t\in[0,1]$.
\item[(ii)] $d(\gamma_{xy}(t),\gamma_{x'y'}(t))\leq(1-t)d(x,x')+td(y,y')$ ~for all ~$t\in[0,1]$.
\eit
\end{dfn}

\begin{note*}
(a) This is a selection of geodesics satisfying the Busemann convexity condition.\\ (b) We know that every normed space $X$ is geodesic since for every $x,y\in X$, we have the special geodesic $\gamma_{xy}=[x,y]:[0,1]\ra X$,~ $\gamma_{xy}(t)=(1-t)x+ty=\gamma_{yx}(1-t)$. This clearly gives a bicombing on $X$, since for every $x,y,x',y'\in X$, we have the convex function
\bea
d(\gamma_{xy}(t),\gamma_{x'y'}(t))=\|\gamma_{xy}(t)-\gamma_{x'y'}(t)\|=\left\|x-x'+t\big(y-y'-(x-x')\big)\right\|.\nn
\eea
\end{note*}

\begin{dfn}[\textcolor{blue}{Weakly \index{Convex bicombing}{convex bicombing}}] Let $X$ be a geodesic space. A weakly convex bicombing on $X$ is a selection of geodesics
\bea
\gamma:X\times X\ra C([0,1],X),~~~~(x,y)\mapsto \gamma_{xy}=\gamma[x,y]~\in~\Gamma(x,y),\nn
\eea
such that for all $x,y,y'\in X$, the following are true.
\bit[leftmargin=0.8cm]
\item[(i)] $\gamma_{xy}(t)=\gamma_{yx}(1-t)$ ~~ for all ~~ $t\in[0,1]$.
\item[(ii)] There exists $C\geq 1$ such that ~~ $d(\gamma_{xy}(t),\gamma_{xy'}(t))\leq Ctd(y,y')$ ~~ for all ~~ $t\in[0,1]$.
\eit
\end{dfn}

\subsection{Lipschitz retractions and extensions}\label{PrlMET4}
(Based on \cite{ben-lind2000})
\begin{dfn}[\textcolor{blue}{\index{Intermediate points}{Intermediate points}}]
Let $X$ be a metric space and $x,y\in X$. A point $z\in X$ is intermediate to $x$ and $y$ (or lies between $x$ and $y$) if $d(x,y)=d(x,z)+d(z,y)$. Let
\bea
[x,y]_0:=\{z\in X:d(x,y)=d(x,z)+d(z,y)\}=\{\txt{points between $x$ and $y$}\}.\nn
\eea
\end{dfn}
{\flushleft Note} that if $X=(X,\|\cdot\|)$ is a normed space, then for any $x,y\in X$, we have
\bea
[x,y]:=\{z_\al=(1-\al)x+\al y:~0\leq\al\leq 1\}\subset [x,y]_0,\nn
\eea
where the inclusion $\subset$ can be strict, as with the case of $X=\Real^2_\infty=(\Real^2,\max|\cdot|)$. Also,
\bea
\|x-z_\al\|=\al\|x-y\|,~~~~\|z_\al-y\|=(1-\al)\|x-y\|.\nn
\eea

\begin{dfn}[\textcolor{blue}{\index{Metrical convexity}{Metrical convexity}}]
A metric space $X$ is metrically convex if for all $x,y\in X$, $0<\al<1$, there exists a point $z_\al\in X$ such that
    \bea
    d(x,z_\al)=\al d(x,y)~~~~\txt{and}~~~~d(z_\al,y)=(1-\al)d(x,y).\nn
    \eea
(Such a point is precisely an intermediate point $z_\al\in[x,y]_0$ such that $d(x,z_\al)=\al d(x,y)$.)
\end{dfn}
{\flushleft Note} that every convex subset $K$ of a normed space $X=(X,\|\cdot\|)$ is metrically convex (with respect to the metric $d(x,y)=\|x-y\|$), since for all $x,y\in K$ we have $[x,y]\subset K\cap[x,y]_0$.

\begin{lmm}
$X$ is metrically convex $\iff$ balls in $X$ intersect as follows.
\bea
\ol{B}_r(x)\cap \ol{B}_s(y)\neq\emptyset~~\txt{iff}~~d(x,y)\leq r+s,~~~~\txt{for all}~~~~x,y\in X.\nn
\eea
\end{lmm}
\begin{proof}
($\Ra$): Assume $X$ is metrically convex. If $z\in \ol{B}_r(x)\cap \ol{B}_s(y)$, then $d(x,y)\leq d(x,z)+d(z,y)\leq r+s$. Conversely, if $d(x,y)\leq r+s$, let $\al:={r\over r+s}$, i.e., $1-\al={s\over r+s}$. Then
\begin{align}
&\textstyle d(x,z_\al)=\al d(x,y)={r\over r+s}d(x,y)\leq r,~~\Ra~~z_\al\in \ol{B}_r(x),\nn\\
&\textstyle d(z_\al,y)=(1-\al)d(x,y)={s\over r+s}d(x,y)\leq s,~~\Ra~~z_\al\in \ol{B}_s(y).\nn
\end{align}
($\La$): Assume $\ol{B}_r(x)\cap \ol{B}_s(y)\neq\emptyset$ iff $d(x,y)\leq r+s$, for any $x,y\in X$. Given $x,y\in X$, $\al\in (0,1)$, let $r=\al d(x,y)$ and $s=(1-\al)d(x,y)$. Then because $d(x,y)=r+s$, there exists some $z_\al\in \ol{B}_r(x)\cap\ol{B}_s(y)$, and so
\bea
&&d(x,y)\leq d(x,z_\al)+d(z_\al,y)\leq r+s=d(x,y),~~\Ra~~d(x,z_\al)+d(z_\al,y)=r+s=d(x,y),\nn\\
&&~~\Ra~~d(x,z_\al)=r=\al d(x,y),~~~~d(z_\al,y)=s=(1-\al)d(x,y).\nn\qedhere
\eea
\end{proof}

\begin{dfn}[\textcolor{blue}{\index{Binary intersection property (BIP)}{Binary intersection property (BIP)}}]
A metric space $X$ has the binary intersection property if for any collection of closed balls $\left\{\ol{B}_{r_\al}(x_\al)\right\}_{\al\in\A}$,
    \bea
    \textstyle\ol{B}_{r_\al}(x_\al)\cap \ol{B}_{r_\beta}(x_\beta)\neq\emptyset~~~~\txt{for all}~~~~\al,\beta\in\A~~~~~~~~\Ra~~~~~~~~\bigcap_{\al\in\A}\ol{B}_{r_\al}(x_\al)\neq\emptyset.\nn
    \eea
\end{dfn}
\begin{rmks}
(i) $\Real=(\Real,|\cdot|)$ has the BIP. (ii) Thus, because balls in $\ell^\infty(A):=\big\{s:A\ra\Real~\big|~\|s\|:=\sup|s(A)|<\infty\big\}$ are ``squares'' (i.e., Cartesian products of real intervals),
\bea
\textstyle\ol{B}_r(s)=\Big\{t\in\ell^\infty(A):|s(a)-t(a)|\leq r,~\txt{for all}~a\in A\Big\}=\prod_{a\in A}\big[s(a)-r,s(a)+r\big],\nn
\eea
it follows that $\ell^\infty(A)$ also has the BIP: Recall that given sets $C_a,D_a$ for each $a\in A$,
\bea
\textstyle\left(\prod_{a\in A}C_a\right)\cap\left(\prod_{a\in A}D_a\right)=\big\{(c_a)_{a\in A}=(d_a)_{a\in A}\big\}=\prod_{a\in A}C_a\cap D_a,\nn
\eea
and so given a collection of balls $\{B_{r_\al}(s_\al)\}$ in $\ell^\infty(A)$, we have
{\footnotesize\bea
\textstyle B_{r_\al}(s_\al)\cap B_{r_\beta}(s_\beta)=\prod_{a\in A}B_{r_\al}(s_\al(a))\cap B_{r_\beta}(s_\beta(a))\neq\emptyset,~~~~B_{r_\al}(s_\al(a)):=\big[s_\al(a)-r_\al,s_\al(a)+r_\al\big],\nn
\eea} for all $(\al,\beta)$ $\iff$ for each $a\in A$ we have $B_{r_\al}(s_\al(a))\cap B_{r_\beta}(s_\beta(a))\neq\emptyset$ for all $(\al,\beta)$.

(iii) \ul{BIP given metrical convexity}: If $X$ is metrically convex, then $X$ has the binary intersection property $\iff$ for all balls $\{\ol{B}_{r_\al}(x_\al)\}$, we have ~$\bigcap\ol{B}_{r_\al}(x_\al)\neq\emptyset$~ iff
\bea
\textstyle d(x_\al,x_\beta)\leq r_\al+r_\beta~~\big(~\txt{iff}~~\ol{B}_{r_\al}(x_\al)\cap \ol{B}_{r_\beta}(x_\beta)\neq\emptyset\big)~~\txt{for all}~~\al,\beta.\nn
\eea
\end{rmks}

\begin{lmm}[\textcolor{OliveGreen}{Modulation of uniform continuity}]\label{UnifContLmm}
A map of metric spaces $f:X\ra Y$ is uniformly continuous $\iff$ for any $\vep>0$, there exists $\delta(\vep)>0$ ($\delta(\vep)\leq\vep$) such that one of the following is true: For all $x,x'\in X$,
\bit[leftmargin=0.7cm]
\item[(i)] $f\big(B_{\delta(\vep)}(x)\big)\subset B_\vep\big(f(x)\big)$. ~~(ii) $d(x,x')<\delta(\vep)$ $\Ra$ $d\big(f(x),f(x')\big)<\vep$.
\item[(iii)] There is a function $u=u_f:[0,\infty)\ra[0,\infty)$, with $u(t)\ra 0$ as $t\ra0$, such that ~$d(x,x')\leq t$ $\Ra$ $d\big(f(x),f(x')\big)\leq u(t)$,~ (for sufficiently small $t$).
\eit
\end{lmm}
\begin{proof} (i) $\Lra$ (ii) is clear from definitions. So we will prove (ii) $\Lra$ (iii).
{\flushleft ($\Ra$):} The function {\small $\omega(t):=\sup\big\{d\big(f(x),f(x')\big):d(x,x')\leq t\big\}$} is nondecreasing, with $\omega(0)=0$, and $d(x,x')\leq t$ $\Ra$ $d\big(f(x),f(x')\big)\leq\omega(t)$, for all $x,x'\in X$. So, if $f$ is uniformly continuous, then for any $\vep>0$, there exists $\delta(\vep)>0$ such that
\bea
d(x,x')\leq\delta(\vep)~~\Ra~~d\big(f(x),f(x')\big)\leq\omega\big(\delta(\vep)\big)\leq\vep,~~~~\txt{for all}~~~~x,x'\in X.\nn
\eea
Thus, for all possible choices of $\delta$, we have $\omega\big(\delta(\vep)\big)\ra0$ as $\vep\ra 0$. Because we can choose $\delta$ such that $\delta(\vep)\ra0$ as $\vep\ra 0$ (which is always possible, say in the form $\vep\geq \delta(\vep)>0$, wlog in the definition of continuity), we see that $\omega(t)\ra0$ as $t\ra0$ (since $\omega(t)$ is nondecreasing and $\omega(0)=0$). Set $u(t):=\omega(t)$.

{\flushleft ($\La$):} Conversely, if there exists a function $u(t)$ such that $u(t)\ra 0$ as $t\ra0$, and (for sufficiently small $t$) $d(x,x')\leq t$ $\Ra$ $d\big(f(x),f(x')\big)\leq u(t)$, $x,x'\in X$, then for any $\vep>0$, we can choose $t=\delta(\vep)>0$ such that $u(t)=u\big(\delta(\vep)\big)\leq\vep$, and so $f$ is uniformly continuous since
\bea
d(x,x')\leq t=\delta(\vep)~~\Ra~~d\big(f(x),f(x')\big)\leq u(t)\leq\vep,~~\txt{for all}~~ x,x'\in X.\nn\qedhere
\eea
\end{proof}
Recall that in case (iii) of Lemma \ref{UnifContLmm}, we say $f$ is \ul{$u$-continuous}. If $u(t)=ct$, we say $f$ is \ul{$c$-Lipschitz continuous}, and if $u(t)=ct^\al$, $\al>0$, we say $f$ is \ul{$(c,\al)$-H\"{o}lder continuous}. Note we can choose $t=d(x,x')$, and so if $f$ is $u$-continuous, then $d\big(f(x),f(x')\big)\leq u\big(d(x,x')\big)$ for sufficiently small $d(x,x')$.

\begin{dfn}[\textcolor{blue}{\index{Modulus of uniform continuity}{Modulus of uniform continuity}}]
If $f:X\ra Y$ is a map of metric spaces, the modulus of uniform continuity of $f$ is $\omega_f:\Real\ra\Real$, $\omega_f(t):=\sup_{d(x,x')\leq t}d\big(f(x),f(x')\big)$.
\end{dfn}
\begin{rmks}\label{ModUContRmk}
Let $f:X\ra Y$ be a map of metric spaces.\\
(i) $f$ is uniformly continuous $\iff$ $\omega_f(t)$ is continuous at $t=0$ ($\iff$ $\omega_f(t)\ra0$ as $t\ra0$), $\iff$ $f$ is $\omega_f$-continuous. Consequently, for any uniformly continuous maps ~$g\circ f:X\sr{f}{\ral}Y\sr{g}{\ral}Z$, it follows from the definition of $\omega_{g\circ f}(t)$ as a supremum that
\begin{align}
\omega_{g\circ f}(t)\leq \omega_g\circ\omega_f(t),~~\txt{since}~~d\big(g\circ f(x),g\circ f(x')\big)\leq \omega_g\big(d\big(f(x),f(x')\big)\big)\leq \omega_g\big(\omega_f\big(d(x,x')\big)\big).\nn
\end{align}
(ii) $f$ is $u$-continuous $\iff$ $\omega_f(t)\leq u(t)$. In particular, $f$ is Lipschitz-continuous $\iff$ $\omega_f(t)\leq ct$, and H$\ddot{\txt{o}}$lder-continuous $\iff$ $\omega_f(t)\leq ct^\al$, $\al>0$.
{\flushleft (iii)} When $X$ is metrically convex, so that $d(x,x')\leq r+s$ $\iff$ $B_r(x)\cap B_s(x')\neq\emptyset$ ($\iff$ $d(x,z)\leq r$ and $d(z,x')\leq s$ for some $z=z_{xx'}\in X$), then $\omega_f$ is subadditive, i.e.,
\begin{align}
&\textstyle \omega_f(r)+\omega_f(s)=\sup_{d(x,z)\leq r}d(f(x),f(z))+\sup_{d(z,x')\leq s}d(f(z),f(x'))\nn\\
&\textstyle~~~~\geq \sup_{d(x,x')\leq r+s}d(f(x),f(z))+\sup_{d(x,x')\leq r+s}d(f(z),f(x'))\nn\\
&\textstyle~~~~\geq \sup_{d(x,x')\leq r+s}d(f(x),f(x'))=\omega_f(r+s).\nn
\end{align}
\end{rmks}

\begin{dfn}[\textcolor{blue}{\index{Lipschitz! retraction}{Lipschitz retraction}, \index{Lipschitz! retract}{Lipschitz retract}}]
A Lipschitz retraction is a Lipschitz-continuous retraction. A metric subspace $Z\subset X$ is a $\ld$-Lipschitz retract if there exists a $\ld$-Lipschitz retraction $r:X\ra Z$. (Equivalently, every $c$-Lipschitz map $f:Z\subset X\ra Y$ extends to a $\ld c$-Lipschitz map $F:X\ra Y$.)
\end{dfn}

\begin{dfn}[\textcolor{blue}{ALR($\S$): \index{Absolute! Lipschitz retract}{Absolute Lipschitz retract} for a class of metric spaces $\S$}]
A metric space $Z$ is a $\ld$-ALR for metric spaces $\S$ (written $Z\in\txt{ALR}_\ld(\S)$) if it is a $\ld$-Lipschitz retract in every closed set inclusion $Z\subset X\in\S$. Equivalently, every $c$-Lipschitz map on a closed set $f:Z\subset X\ra Y$, $X\in\S$, extends to a $\ld c$-Lipschitz map $F:X\ra Y$.
\bc\bt
Z\ar[d,hook]\ar[rrr,"f"] &&& Y\\
X\in\S\ar[urrr,dashed,"F"']
\et\ec

If $Z$ is a $\ld$-ALR for all metric spaces, then we say simply that $Z$ is a $\ld$-ALR (or $Z\in ALR_\ld$).
\end{dfn}

\begin{dfn}[\textcolor{blue}{Subadditively modulated (SM) map}]
A map of metric spaces $f:X\ra Y$ is a uniformly modulated map (or SM map) map if $\omega_f(t)\leq\omega(t)$ for some subadditive nondecreasing function $\omega:[0,\infty)\ra[0,\infty)$ such that $\omega(t)\ra 0$ as $t\ra0$. In this case, we say $f$ is a $\omega$-SM map.
\end{dfn}

\begin{rmk*}
Every SM map is uniformly continuous. Moreover, by Remark \ref{ModUContRmk}(iii), if $X$ is metrically convex, then $f:X\ra Y$ is uniformly continuous $\iff$ $f$ is a SM map.
\end{rmk*}

\begin{dfn}[\textcolor{blue}{Uniform retraction, \index{Uniform retract}{Uniform retract}}]
A uniform retraction is a uniformly continuous retraction. A metric subspace $A\subset X$ is a uniform retract if there exists a uniform retraction $r:X\ra A$. Equivalently, every uniformly continuous map $f:A\subset X\ra Y$ extends to a uniformly continuous map $F:X\ra Y$. (Note that in this case,
\begin{align}
\textstyle\omega_f\leq\omega_F\leq\omega_f\circ\omega_r,~~\txt{for a nondecreasing}~~\omega_r:[0,\infty)\ra[0,\infty)~~\txt{s.t.}~~\lim_{t\ra 0}\omega_r(t)=0,\nn
\end{align}
in which case, we say $A\subset X$ is a $\omega_r$-uniform retract.)
\end{dfn}

\begin{dfn}[\textcolor{blue}{AUR($\S$): \index{Absolute! uniform retract}{Absolute uniform retract} for a class of metric spaces $\S$}]
A metric space $Z$ is an AUR for metric spaces $\S$ (written $Z$ $\in$ AUR($\S$)) if it is a uniform retract in every closed set inclusion $Z\subset X\in\S$. Equivalently, every uniformly continuous map on a closed set $f:Z\subset X\ra Y$, $X\in\S$, extends to a uniformly continuous map $F:X\ra Y$.
\bc\bt
Z\ar[d,hook]\ar[rrr,"f"] &&& Y\\
X\in\S\ar[urrr,dashed,"F"']
\et\ec
(Note that in the latter we can also consider the special case where $\omega_f\leq\omega_F\leq\omega_f\circ\omega$ for a nondecreasing $\omega:[0,\infty)\ra[0,\infty)$ such that $\lim_{t\ra 0}\omega(t)=0$, in which case, we say $Z$ is a $\omega$-AUR for the spaces $\S$, and write $Z\in\txt{AUR}_\omega(\S)$.)

If $Z$ is an AUR for all metric spaces, then we say simply that $Z$ is an AUR (or $Z\in AUR$).
\end{dfn}

\begin{lmm}[\textcolor{OliveGreen}{$\ell^\infty(A)$ is an AUR: \cite[Lemma 1.1]{ben-lind2000}}]\label{ExtToLinf}
Let $X$ be any metric space. (i) $X$ imbeds isometrically into some $\ell^\infty(A)$. (ii) Every $\omega$-SM map $f:Y\subset X\ra \ell^\infty(A)$ extends to a $\omega$-SM map $F:X\ra\ell^\infty(A)$. (In particular, if $f$ is $c$-Lipschitz, then so is $F$.)
\end{lmm}
\begin{proof}
(i) Fix $x_0\in X$, and consider the map $\vphi:X\ra\ell^\infty(X),~x\mapsto\vphi_x$, where $\vphi_x:X\ra\Real$ is given by $\vphi_x(a):=d(x,a)-d(x_0,a)$, for all $a\in A:=X$, which is well defined since $\|\vphi_x\|=\sup_a|d(x,a)-d(x_0,a)|\leq d(x,x_0)$. Then $\|\vphi_x-\vphi_{x'}\|=d(x,x')$, since
\begin{align}
\textstyle\|\vphi_x-\vphi_{x'}\|=\sup_a|\vphi_x(a)-\vphi_{x'}(a)|\leq d(x,x')=|\vphi_x(x)-\vphi_{x'}(x)|\leq\|\vphi_x-\vphi_{x'}\|.\nn
\end{align}
(ii) Let $F:X\ra \ell^\infty(A),~x\mapsto F_x$ be given by $F_x(a):=\inf\limits_{y\in Y}~\Big(f_y(a)+\omega\big(d(x,y)\big)\Big)$ for each $a\in A$. $F_x(a)$ is finite because ~$|f_y(a)-f_{y'}(a)|\leq\|f_y-f_{y'}\|\leq \omega_f\big(d(y,y')\big)\leq \omega\big(d(y,y')\big)$,
\begin{align}
&~~\Ra~f_y(a)-f_{y'}(a)\leq\omega\big(d(y,y')\big)\leq\omega\big(d(y,x)+d(x,y')\big)\leq\omega\big(d(y,x)\big)+\omega\big(d(x,y')\big)\nn\\
&~~\Ra~~f_y(a)-\omega\big(d(y,x)\big)\leq f_{y'}(a)+\omega\big(d(x,y')\big),~~\txt{for all}~~a\in A,~x\in X,~y,y'\in Y,\nn\\
&~~\Ra~~f_y(a)-\omega\big(d(y,x)\big)\leq \inf_{y'\in Y}\Big(f_{y'}(a)+\omega\big(d(x,y')\big)\Big)= F_x(a),~~a\in A,~x\in X,~y\in Y.\nn
\end{align}
Now, by setting $x=y\in Y$ above, we get $f_y(a)-0\leq F_y(a)\leq f_y(a)$, and so $F_y(a)=f_y(a)$ for all $y\in Y$ (i.e., $F$ is an extension of $f$). Finally, for all $a\in A$, $x,x'\in X$, $y\in Y$,
\begin{align}
&F_x(a)\leq f_y(a)+\omega\big(d(x,y)\big)\leq f_y(a)+\omega\big(d(x,x')\big)+\omega\big(d(x',y)\big),\nn\\
&~~\Ra~~F_x(a)\leq F_{x'}(a)+\omega\big(d(x,x')\big),~~\txt{and similarly},~~F_{x'}(a)\leq F_x(a)+\omega\big(d(x,x')\big),\nn\\
&~~\Ra~~|F_x(a)-F_{x'}(a)|\leq\omega\big(d(x,x')\big),~~\txt{for all}~~a\in A,~~\Ra~~\|F_x-F_{x'}\|\leq \omega\big(d(x,x')\big),\nn\\
&\textstyle~~\Ra~~\omega_F(t)=\sup_{d(x,x')\leq t}\|F_x-F_{x'}\|\leq \sup_{d(x,x')\leq t}\omega\big(d(x,x')\big)\leq \omega(t).\nn\qedhere
\end{align}
\end{proof}

\begin{crl}
If $X$ is metrically convex, then every uniformly cont. map $f:Y\subset X\ra\ell^\infty(A)$ extends to a uniformly cont. map $F:X\ra\ell^\infty(A)$, such that $\omega_F(t)=\omega_f(t)$.
\end{crl}

\begin{rmk}
Let $X$ be a metrically convex space, $Z$ any metric space, and $f:Y\subset X\ra Z$ uniformly continuous. If there exists no nondecreasing subadditive function $\omega(t)$ such that $\omega_f(t)\leq\omega(t)$, then $f$ cannot be extended to a uniformly continuous map $F:Y\subset X\ra Z$. Otherwise, if $F:X\ra Z$ is a uniformly continuous extension of $f$, then $\omega_f(t)\leq\omega_F(t)$, where $\omega_F(t)$ is nondecreasing and subadditive due to metrical convexity of $X$.
\end{rmk}

\begin{rmk}
It follows from Lemma \ref{ExtToLinf}(ii) that $\ell^\infty(A)$ is  a 1-ALR.
\end{rmk}

\begin{prp}[\textcolor{OliveGreen}{ALR criterion: \cite[Proposition 1.2]{ben-lind2000}}]\label{ALRcrit}
If $X$ is a metric space, then the following are equivalent. (i) $X$ is an absolute $\ld$-Lipschitz retract ($\ld$-ALR).
\bit[leftmargin=0.7cm]
\item[(ii)] Any $c$-Lipschitz map $f:Z\subset Y\ra X$ extends to a $\ld c$-Lipschitz map $F:Y\ra X$.
\item[(iii)] Any $c$-Lipschitz map $f:X\subset Y\ra Z$ extends to a $\ld c$-Lipschitz map $F:X\ra Z$.
\eit
\end{prp}
\begin{proof}
\ul{(i)$\Ra$(ii)}: By Lemma \ref{ExtToLinf}(i), $X$ imbeds isometrically into some $\ell^\infty(A)$ as $i_X:X\ra\ell^\infty(A)$. So, because $X$ is a $\ld$-ALR, we get a $\ld$-Lipschitz retraction $r:\ell^\infty(A)\ra X$. Also, by Lemma \ref{ExtToLinf}(ii), $g=i_X\circ f:Z\ra\ell^\infty(A)$ extends to a $c$-Lipschitz map $G:Y\ra \ell^\infty(A)$. Hence, $f:Z\ra X$ extends to $\ld c$-Lipschitz $F=r\circ G:Y\ra X$.

\bc\bt
Z\ar[d,hook]\ar[drr,"g=i_X\circ f"'description]\ar[rr,"f"] && X\ar[d,hook,"i_X"]\ar[rr,"id_X"] && X\\
Y\ar[rr,dashed,"G"]&& \ell^\infty(A)\ar[urr,dashed,"r"] &&
\et\ec
{\flushleft\ul{(ii)$\Ra$(iii)}:} By $(ii)$, $id_X:X\ra X$ extends to a $\ld$-Lipschitz retraction $r:Y\ra X$. So $f:X\ra Z$ extends to a $\ld c$-Lipschitz map $F=f\circ r:Y\ra Z$.

\bc\bt
 X\ar[d,hook]\ar[rr,"id_X"]  && X\ar[rrr,"f"] &&& Z\\
 Y \ar[urr,dashed,"r"]\ar[urrrrr,dashed,"F=f\circ r"'] && &&&
\et\ec
{\flushleft\ul{(iii)$\Ra$(i)}:} {\small With $Z=X$ and $f=id_X$, $X$ is a $\ld$-Lipschitz retract in every inclusion $X\subset Y$.}
\end{proof}

\begin{crl}
If $X$ is a $\ld$-ALR, then every $\mu$-Lipschitz retract $Z\subset X$ is a $\mu\ld$-ALR. (To prove this, use Proposition \ref{ALRcrit}(ii))
\end{crl}

\begin{prp}[\textcolor{OliveGreen}{AUR criterion}]\label{AURcrit}
If $X$ is a metric space, then the following are equivalent. (i) $X$ is an absolute uniform retract (AUR).
\bit[leftmargin=0.7cm]
\item[(ii)] Any uniformly cont. $f:Z\subset Y\ra X$ extends to a uniformly cont. $F:Y\ra X$.
\item[(iii)] Any uniformly cont. $f:X\subset Y\ra Z$ extends to a uniformly cont. $F:X\ra Z$.
\eit
\end{prp}
\begin{proof}
{\flushleft\ul{(i)$\Ra$(ii)}:} By Lemma \ref{ExtToLinf}(i), $X$ imbeds isometrically into some $\ell^\infty(A)$ as $i_X:X\ra\ell^\infty(A)$. So, because $X$ is an AUR, we get a uniform retraction $r:\ell^\infty(A)\ra X$ (an SM retraction since $\ell^\infty(A)$ is metrically convex). Also, by Lemma \ref{ExtToLinf}(ii), if $f:Z\ra X$ is a uc map, then the uc map $g=i_X\circ f:Z\ra\ell^\infty(A)$ extends to a uc map $G:Y\ra\ell^\infty(A)$. Hence, a uc map $f:Z\ra X$ extends to a uc map $F=r\circ G:Y\ra X$.
\bc\bt
Z\ar[d,hook]\ar[drr,"g=i_X\circ f"'description]\ar[rr,"f"] && X\ar[d,hook,"i_X"]\ar[rr,"id_X"] && X\\
Y\ar[rr,dashed,"G"] && \ell^\infty(A)\ar[urr,dashed,"r"] &&
\et\ec
{\flushleft\ul{(ii)$\Ra$(iii)}:} By $(ii)$, the uc map $id_X:X\ra X$ extends to a uniform retraction $r:Y\ra X$. So a uc map $f:X\ra Z$ extends to the uc map $F=f\circ r$.
\bc\bt
 X\ar[d,hook]\ar[rr,"id_X"]  && X\ar[rrr,"f"] &&& Z\\
 Y \ar[urr,dashed,"r"]\ar[urrrrr,dashed,"F=f\circ r"'] && &&&
\et\ec
{\flushleft\ul{(iii)$\Ra$(i)}:} With $Z=X$ and $f=id_X$, $X$ is a uniform retract in every inclusion $X\subset Y$.
\end{proof}

\begin{crl}
If $X$ is an AUR, then every uniform retract $Z\subset X$ is itself an AUR.
\end{crl}

\begin{prp}[\textcolor{OliveGreen}{1-ALR criterion: \cite[Proposition 1.4]{ben-lind2000}}]\label{1ALRcrit}
 A metric space $X$ is $1$-ALR $\iff$ metrically convex and BIP. (Equivalently, $X$ is $1$-ALR $\iff$ geodesic and BIP.)
\end{prp}
\begin{proof}
($\Ra$): Assume $X$ is $1$-ALR. Given $x,y\in X$, the $d(x,y)$-Lipschitz map {\small $f:\{0,1\}\subset[0,1]\ra X$, $f(0)=x$, $f(1)=y$}, extends to a $d(x,y)$-Lipschitz map $F:[0,1]\ra X$. It follows that for any $0<\al<1$, $z_\al:=F(\al)$ satisfies
\bea
&& d(x,z_\al)=d(F(0),F(\al))\leq d(x,y)\al,~~~~d(z_\al,y)=d(F(\al),F(1))\leq d(x,y)(1-\al),\nn\\
&&~~\Ra~~d(x,y)\leq d(x,z_\al)+d(z_\al,y)\leq d(x,y)\al+d(x,y)(1-\al)=d(x,y),\nn\\
&&~~\Ra~~d(x,z_\al)= \al d(x,y),~~~~d(z_\al,y)=(1-\al)d(x,y),\nn
\eea
which shows $X$ is metrically convex. (Note $F$ is a geodesic, and so $X$ is a geodesic space.)

To see that $X$ is a BIP space, let $\{\ol{B}_{r_\al}(x_\al)\}_\al$ be a collection of balls in $X$ such that $\ol{B}_{r_\al}(x_\al)\cap \ol{B}_{r_\beta}(x_\beta)\neq\emptyset$ for all $\al,\beta$. Since $X$ imbeds isometrically into some $\ell^\infty(A)$ (thus preserving radii and pairwise intersection of balls) we have maps $X\sr{\vphi}{\ral}\vphi(X)\sr{i}{\hookrightarrow}\ell^\infty(A)\sr{r}{\ral}\vphi(X)\sr{\vphi^{-1}}{\ral}X$, where $\vphi$ is an isometry (i.e., $i\circ\vphi$ is an isometric imbedding) and $r$ is a $1$-Lipschitz retraction. Because $\ell^\infty(A)$ has the binary intersection property, the images of the balls under $i\circ\vphi$ have a common intersection point $u\in\bigcap \ol{B}_{r_\al}\big(i\circ\vphi(x_\al)\big)=\bigcap i\circ\vphi\left(\ol{B}_{r_\al}(x_\al)\right)$ in $\ell^\infty(A)$, i.e., $d(u,i\circ\vphi(x_\al))\leq r_\al$ for all $\al$. Thus, the balls also have a common intersection point $\vphi^{-1}\circ r(u)$ in $X$, since
\bea
d\left(\vphi^{-1}\circ r(u),x_\al\right)=d\Big((\vphi^{-1}\circ r)(u),(\vphi^{-1}\circ r)\circ(i\circ\vphi)(x_\al)\Big)\leq 1~d\Big(u,i\circ\vphi(x_\al)\Big)\leq r_\al.\nn
\eea

($\La$): Assume $X$ is metrically convex and BIP. For a metric space $Y$, consider a $c$-Lipschitz map $f:Z\subset Y\ra X$. Let $y\in Y\backslash Z$, and set $\wt{Z}:=Z\cup\{y\}$. The balls $\left\{\ol{B}_{cd(z,y)}\big(f(z)\big)\right\}_{z\in Z}$ intersect pairwise in $X$, since metrical convexity of $X$ implies that for any $z,z'\in Z$,
\begin{align}
d\big(f(z),f(z')\big)\leq cd(z,z')\leq cd(z,y)+cd(z',y)~~\Ra~~\ol{B}_{cd(z,y)}\big(f(z)\big)\cap \ol{B}_{cd(z',y)}\big(f(z')\big)\neq\emptyset.\nn
\end{align}
Thus, by the BIP of $X$, there exists a point ~$x\in\bigcap\limits_{z\in Z}\ol{B}_{cd(z,y)}\big(f(z)\big)$. The map $\wt{f}:\wt{Z}\ra X$ given by $\wt{f}|_Z=f$ and $\wt{f}(y):=x\in\bigcap\limits_{z\in Z}\ol{B}_{cd(z,y)}\big(f(z)\big)$, i.e., $d\big(f(z),x\big)\leq cd(z,y)$ for all $z\in Z$, is a $c$-Lipschitz extension of $f$, since for any $z\in Z$, $d\big(\wt{f}(z),\wt{f}(y)\big)=d\big(f(z),x\big)\leq cd(z,y)$.
Thus, we get a nonempty poset $\P=\big\{\wt{f}:\wt{Z}\subset Y\ra X~\big|~\wt{f}~\txt{a $c$-Lipschitz extension of $f$}\big\}$ in which every nonempty chain $\big\{\wt{f}_\ld:\wt{Z}_\ld\subset Y\ra X\big\}_\ld$ has an upper bound $\bigcup\wt{f}_\ld:\bigcup\wt{Z}_\ld$ $\subset$ $Y\ra X$. It follows by Zorn's lemma that $\P$ has a maximal element $F:\wt{Z}_F\subset Y\ra X$. If $\wt{Z}_F\neq Y$, then $F$ fails to be maximal (since, as before, any $y\in Y\backslash\wt{Z}_F$ leads to a strictly larger $c$-Lipschitz extension $\wt{F}:\wt{Z}_F\cup\{y\}\subset Y\ra X$ of $f$). Hence, we have a $c$-Lipschitz extension $F:Y\ra X$ of $f$.
\end{proof}

\begin{prp}[\textcolor{OliveGreen}{\cite[Theorem 2.4.16]{BBI}}]\label{ComMetricConv}
Every geodesic space is metrically convex. Moreover, if $X$ is a complete metric space, then $X$ is geodesic $\iff$ metrically convex.
\end{prp}
\begin{proof}
If $X$ is geodesic, it is clear that $X$ is metrically convex (where completeness is not required). Conversely, assume $X$ is complete and metrically convex, and let $x,y\in X$. To construct a geodesic $\gamma:[0,1]\ra X$ from $x$ to $y$, let $c(0)=x$, $c(1)=y$. Next, by metrical convexity, there exists a point $c(1/2)$ such that
\bea
\textstyle d(c(1),c(1/2))=d(c(1/2),c(1))={1\over 2}d(c(0),c(1))={1\over 2}d(x,y).\nn
\eea
Similarly, there exist points $c((0+1/2)/2)=c(1/4)$ and $c((1/2+1)/2)=c(3/4)$ such that
\begin{align}
&\textstyle d(c(0),c(1/4))=d(c(1/4),c(1/2))={1\over 2}d(c(0),c(1/2))={1\over 4}d(x,y),\nn\\
&\textstyle d(c(1/2),c(3/4))=d(c(3/4),c(1))={1\over 2}d(c(1/2),c(1))={1\over 4}d(x,y),\nn
\end{align}
 from which we obtain
\begin{align}
&\textstyle d(c(0),c(t))=td(x,y),~~~~d(c(t),d(c(1)))=(1-t)d(x,y),~~~~\txt{for all}~~~~t\in\left\{0,{1\over 4},{2\over 4},{3\over 4},1\right\},\nn\\
&\textstyle~~\Ra~~ d(c(t),c(t'))\leq|t-t'|d(x,y),~~~~\txt{for all}~~~~t,t'\in\left\{0,{1\over 4},{2\over 4},{3\over 4},1\right\}.\nn
\end{align}
Continuing this way, for each $n\geq 1$ (and with $D_n:=\left\{{k\over 2^n}:0\leq k\leq 2^n\right\}$) we obtain
\begin{align}
&\textstyle d(c(0),c(t))=td(x,y),~~~d(c(t),d(c(1)))=(1-t)d(x,y),~~~\txt{for all}~~~t\in D_n,\nn\\
&\textstyle ~~\Ra~~d(c(t),c(t'))\leq|t-t'|d(x,y)~~~~\txt{for all}~~~~t,t'\in D_n.\nn
\end{align}
Since $D_n\subset D_{n+1}$, it follows that with $D:=\bigcup D_n$, we have
\begin{align}
d(c(t),c(t'))\leq|t-t'|d(x,y)~~~~\txt{for all}~~~~t,t'\in D.\nn
\end{align}
Since $D\subset[0,1]$ is dense, i.e., $\ol{D}=[0,1]$, and $X$ is complete, $c:D\subset [0,1]\ra X$ extends to a map $\gamma:[0,1]\ra X$ such that $d(\gamma(t),\gamma(t'))\leq|t-t'|d(x,y)$ for all $t,t'\in[0,1]$.
\end{proof}

\chapter{Digressions}\label{Digress}
\section{Partial ordering and Zorn's lemma}\label{ParOrdZorn}
Recall that our preferred set containment operations are $\subset$ for non-strict containment and $\subsetneq$ for strict containment. In a poset we will (for convenience) sometimes \ul{also} use $\subseteq$ (in place of $\subset$) for non-strict containment. That is $\subset,\subseteq$ each stand for ``non-strict containment'' while $\subsetneq$ stands for ``strict containment''.

\begin{dfn}[\textcolor{blue}{\index{Binary relation}{Binary relation}, Related element, Intermediate element}]
A binary relation on a set $A$ is any set $R\subset A\times A:=\{(a,b):a,b\in A\}$. If $a,b\in A$, we say ``$a$ is related to $b$ (through $R$)'', written $aR b$, if $(a,b)\in R$. For any $a,b,c\in A$, we say ``$c$ lies between $a$ and $b$ (with respect to $R$)'', or ``$c$ is intermediate to $a$ and $b$ (with respect to $R$)'', written $aRcRb$, if $aRc$ and $cRb$.
\end{dfn}
\begin{dfn}[\textcolor{blue}{\index{Reflexivity}{Reflexivity}, \index{Symmetry}{Symmetry}, \index{Transitivity}{Transitivity}, \index{Antisymmetry}{Antisymmetry}}]
Let $A$ be a set and $R\subset A\times A$. The binary relation $R$ is (i) reflexive if $aRa$ for all $a\in A$, (ii) symmetric if $aRb$ implies $b R a$ for all $a,b\in A$, (iii) antisymmetric if $aRbRa$ implies $a=b$ for all $a,b\in A$, (iv) transitive if $aRbRc$ implies $aRc$ for all $a,b,c\in A$.
\end{dfn}

Note that (unless a relation is symmetric) if $a$ is related to $b$, it does not follow that $b$ is related to $a$. Similarly (unless a relation, if transitive, is symmetric), if $c$ lies between $a$ and $b$, it does not follow that $c$ lies between $b$ and $a$.

\begin{dfn}[\textcolor{blue}{\index{Equivalence relation}{Equivalence relation}}]
A binary relation that is reflexive, symmetric, and transitive. That is, a binary relation $\sim~\subset A\times A$ such that for all $a,b,c\in A$,
\bea
(i)~a\sim a,~~~~(ii)~a\sim b~~\Ra~~b\sim a,~~~~(iii)~a\sim b\sim c~~\Ra~~a\sim c.\nn
\eea
\end{dfn}

\begin{dfn}[\textcolor{blue}{\index{Partial order}{Partial order}, \index{Poset}{Partially ordered set (poset)}}]
A partial order is a binary relation that is reflexive, antisymmetric, and transitive, i.e., a binary relation $\leq~\subset P\times P$ (making $P=(P,\leq)$ a partially ordered set or poset) such that for all $a,b,c\in P$,
\bea
(i)~a\leq a,~~~~(ii)~a\leq b\leq a~~\Ra~~a=b,~~~~(iii)~a\leq b\leq c~~\Ra~~a\leq c.\nn
\eea
If $a\leq b$ (or $b\geq a$) we say ``\ul{$a$ is smaller than $b$}'' (or ``\ul{$b$ is larger than $a$}''). Also, if $a\leq b$ and $a\neq b$, we write $a<b$ (a situation called \ul{strict ordering}).
\end{dfn}
Note that in a poset $P=(P,\leq)$, every subset $S\subset P$ is itself a poset $(S,\leq)$.

\begin{dfn}[\textcolor{blue}{\index{Comparable elements}{Comparable elements}, Linearly/totally ordered set (\index{Chain!}{chain}) in a poset}]
In a poset $(P,\leq)$, two elements $a,b\in P$ are comparable if $a\leq b$ or $b\leq a$.

A set $C\subset P$ is linearly ordered (or a chain) if every two elements in $C$ are comparable, i.e., $C\subset P$ is linearly ordered if for any $a,b\in C$, we have $a\leq b$ or $b\leq a$.
\end{dfn}

\begin{dfn}[\textcolor{blue}{\index{Upper bound}{Upper bound}, \index{Lower bound}{Lower bound}, \index{Directed set}{Directed set}}]
Let $(P,\leq)$ be a poset. An element $u\in P$ is an upper bound of a set $A\subset P$ (written $A\leq u$) if $a\leq u$ for all $a\in A$. An element $l\in P$ is a lower bound of a set $A\subset P$ (written $l\leq A$) if $l\leq a$ for all $a\in A$.

A directed set is a poset in which any two elements have an upper bound, i.e., $(P,\leq)$ is directed if for any $a,b\in P$, there exists $c\in P$ such that $a,b\leq c$ (i.e., $a\leq c$ and $b\leq c$).
\end{dfn}

\begin{dfn}[\textcolor{blue}{Strict upper bound, Strict lower bound}]
If $A\leq u$ and $u\not\in A$, we write $A<u$. Similarly, if $v\leq A$ and $v\not\in A$, we write $v<A$.
\end{dfn}

\begin{dfn}[\textcolor{blue}{\index{Minimal element}{Minimal element}, \index{Maximal element}{Maximal element}, \index{Least element}{Least element}, \index{Greatest element}{Greatest element}}]
Let $(P,\leq)$ be a poset. An element $m\in P$ is \ul{minimal} if for all $x\in P$,~ $x\leq m$ implies $m=x$. (That is, $m$ is an element that cannot be bigger than any element that is comparable with it.) An element $m\in P$ is \ul{maximal} if for all $x\in P$,~ $m\leq x$ implies $m=x$. (That is, $m$ is an element that cannot be smaller than any element that is comparable with it.)

An element $y\in P$ is a \ul{least} (resp. \ul{largest}) element if $y\leq x$ (resp. $x\leq y$) for all $x\in P$.
\end{dfn}

\begin{dfn}[\textcolor{blue}{\index{Well-ordered (po)set}{Well-ordered (po)set}}]
A (po)set $(P,\leq)$ is well ordered if (i) it is a chain and (ii) every nonempty subset $\emptyset\neq A\subset P$ has a minimal element.
\end{dfn}

\begin{dfn}[\textcolor{blue}{Greatest lower bound (glb or \index{Infimum}{infimum}), glb property}]
Let $(P,\leq)$ be a poset, and $A\subset P$ a nonempty subset with a lower bound. A \ul{glb or infimum} (greatest lower bound) of $A$, written $glb(A)$, is a greatest element in the set of all lower bounds of $A$, i.e.,
\bea
L(A):=\{l\in P:l\leq A\}\leq\txt{glb}(A)\leq A.\nn
\eea

A poset $(P,\leq)$ has the \ul{glb property} if every nonempty subset of $P$ with a lower bound has a greatest lower bound.
\end{dfn}

\begin{dfn}[\textcolor{blue}{Least upper bound (lub or \index{Supremum}{supremum}), lub property}]
Let $(P,\leq)$ be a poset, and $A\subset P$ a nonempty subset with an upper bound. A \ul{lub or supremum} (least upper bound) of $A$, written $lub(A)$, is a least element in the set of all upper bounds of $A$, i.e.,
\bea
A\leq\txt{lub}(A)\leq U(A):=\{u\in P:A\leq u\}.\nn
\eea

A poset $(P,\leq)$ has the \ul{lub property} if every nonempty subset of $P$ with an upper bound has a least upper bound.
\end{dfn}

\begin{dfn}[\textcolor{blue}{\index{Initial segment}{Initial segment}}]
Given chains $C,D$ in a poset $(P,\leq)$, we say $C$ is an initial segment of $D$ (written $C\preceq D$ or $C\leq_{\txt{seg}}D$) if $C\subseteq D$ and $C<D\backslash C$ (i.e., $C<d$ for all $d\in D\backslash C$). If $C\preceq D$ and $C\neq D$ (i.e., $C\subsetneq D$ and $C<D\backslash C$), we write $C\prec D$.
\end{dfn}

Note that the empty chain $\emptyset$ is (naturally) an initial segment of every nonempty chain. Also, the collection of chains $\C$ in $P$ gives posets $(\C,\subseteq)$ and $(\C,\preceq)$.

\begin{dfn}[\textcolor{blue}{Lower section (at a point)}]
Let $(P,\leq)$ be a poset. The lower section at $x\in P$ is the set of lower bounds ~$L_x:=\{y\in P:y\leq x\}$~ of~ $x$.
\end{dfn}
\begin{dfn}[\textcolor{blue}{\index{Regular poset}{Regular poset}}]
A (nonempty) poset in which every (nonempty) chain has an upper bound.
\end{dfn}

\begin{lmm}\label{ZrnPreLmm}
Let $(P,\leq)$ be a (nonempty) poset. The collections $\C=\C(P):=\{\txt{chains}~C\subset P\}$ and $\L=\L(P):=\{L_x~|~x\in P\}$ are posets $(\C,\subseteq)$, $(\L,\subseteq)$ [being subsets of the poset $(\P(P),\subseteq)$]. Moreover, the following hold. For all $x,y\in P$,
\bit[leftmargin=0.7cm]
\item[(a)] $x\leq y$ $\iff$ $L_x\subseteq L_y$ (and ~$x<y$ $\iff$ $L_x\subsetneq L_y$ ).
\item[(b)] $x$ is maximal in $P$ $\iff$ $L_x$ is maximal in $\L$.
\item[(c)] Assume $P$ is regular. If there exists a maximal chain $C\in\C$, then $P$ has a maximal element (an upper bound of $C$ in $P$).
\item[(d)] $(\C,\subseteq)$ is a regular poset. Also, $(P,\leq)$ is regular $\iff$ $(\L,\subseteq)$ is regular.
\eit
\end{lmm}
\begin{proof}
(a),(b),(c) are immediate consequences of the definitions. For (d), unions give upper bounds.
\end{proof}

The following proof is based on Pierre-Yves Gaillard's version of \cite{lewin1991} at \url{http://vixra.org/abs/1207.0064}.

\begin{thm}[\textcolor{OliveGreen}{\index{Zorn's lemma}{Zorn's lemma}: \cite{lewin1991}}]\label{ZornLemma}
A regular poset has a maximal element.
\end{thm}
\begin{proof}
Let $(P,\leq)$ be a regular poset (i.e., every nonempty chain in $P$ has an upper bound). For any sets $A\subset B$ in $P$, let ${}_{A<}B=U^\ast_A(B):=\{b\in B:A<b\}$ be the set of strict upper bounds of $A$ in $B$, and $B_{<A}=L^\ast_A(B):=\{b\in B:b<A\}$  be the set of strict lower bounds of $A$ in $B$. That is, ~$B_{<A}~<~A~<~{}_{A<}B$.

Suppose $P$ has no maximal element. Then, by Lemma \ref{ZrnPreLmm}(c), every chain $C\subset P$ is not maximal in $(\C(P),\subseteq)$, and so its set of strict upper bounds ${}_{C<}P$ is nonempty. Thus, we can define a strict upper bound selection (or choice) function
\bea
u:\{\txt{nonempty chains in}~P\}\ra P,~~C\mapsto u(C)\in{}_{C<}P=U^\ast_C(P)\nn
\eea
along with the associated \emph{nonempty} poset $(\W,\subseteq)$ of those well-ordered chains in $P$ given by
\bea
\W:=\left\{W\in \C(P):~W~\txt{is well-ordered, and}~u(W_{<w})=w~\txt{for all}~w\in W\right\}.\nn
\eea
Note that for all $W\in\W$, we have $W\cup\{u(W)\}\in \W$ as well, because ~{\small$u\left((W\cup\{u(W)\})_{<u(W)}\right)=u\left(W_{<u(W)}\right)=u(W)$} implies
~{\small $u\big((W\cup\{u(W)\})_{<w}\big)=
\left\{
  \begin{array}{ll}
   u(W) , & \txt{if}~~w=u(W) \\
    u\left(W_{<w}\right)=w, & \txt{if}~~w\in W
  \end{array}
\right\}=w$}.

For any $W_1,W_2\in \W$, let $L:=\bigcup\{C\in \W:C\preceq W_1,W_2\}\in\W$ be the union of all common \emph{initial segments} of $W_1$ and $W_2$. Then we see that $L\preceq W_1,W_2$, i.e., $L\in\W$, and that $L$ is maximal in $\W$ with respect to this property. If $L\prec W_1,W_2$, then there exist $w_1\in W_1$, $w_2\in W_2$ such that $L=(W_1){}_{<w_1}=(W_2){}_{<w_2}$, $L\cup\{w_1\}\preceq W_1$, and $L\cup\{w_2\}\preceq W_2$,
\bea
&&~~\Ra~~w_1=u\left((W_1){}_{<w_1}\right)=u(L)=u\left((W_2){}_{<w_2}\right)=w_2,\nn\\
&&~~\Ra~~L\cup\{u(L)\}\preceq W_1,W_2,~~\Ra~~L\subsetneq L\cup\{u(L)\}\in\W,\nn
\eea
contradicting maximality (wrt $\preceq W_1,W_2$) of $L$ in $\W$. Thus, $L=W_1$ or $L=W_2$, i.e.,
\bea
W_1\preceq W_2~~\txt{or}~~W_2\preceq W_1,~~~~\Ra~~~~W_1\subseteq W_2~~\txt{or}~~W_2\subseteq W_1,\nn
\eea
which shows $\W$ is totally ordered, hence a chain, in $(\C(P),\subseteq)$.

Let $W_0:=\bigcup\{W\in\W\}$. Then $W_0\in \W$, and so $W_0\subsetneq W_0\cup\{u(W_0)\}\in \W$ (a contradiction).
\end{proof}

\begin{rmk}
Instead of ``all chains have an upper bound,'' the above proof only requires that ``every well-ordered subset has an upper bound''.
\end{rmk}

\begin{crl}[\textcolor{OliveGreen}{Hausdorff's \index{Maximality principle}{maximality principle}}]
Every (nonempty) poset $(P,\leq)$ has a maximal chain (by Lemma \ref{ZrnPreLmm}d).
\end{crl}

\begin{crl}
Let $(P,\leq)$ be a poset. If $x\in P$ is a maximal element, then there exists a maximal chain $C_x\in\C(P)$ such that $C_x\leq x\in C_x$ (i.e., $x$ is a maximal element of $C_x$).
\end{crl}

\begin{crl}[\textcolor{OliveGreen}{Zermelo's \index{Well-ordering theorem}{well-ordering theorem}}]
Every set $S$ can be well-ordered.
\end{crl}
\begin{proof}
Let $(P,\leq)$ be the poset with $P:=\{\txt{well-orderings}~(A,\leq_A)~|~A\subset S\}$ consisting of all possible well-orderings of subsets of $S$, and
{\small\bea
(A,\leq_A)\leq (B,\leq_B)~~\txt{if}~~\leq_B\big|_A~\txt{equals}~\leq_A~~\txt{and}~~A\preceq B~~(\txt{i.e., $A$ is an initial segment of $B$}).\nn
\eea}Then $(P,\leq)$ is a regular poset, and so has a maximal element $(M,\leq_M)\in P$. Suppose $M\neq S$, i.e., there exists $e\in S\backslash M$. Then $M_e:=M\cup\{e\}$ can be well-ordered as $(M_e,\leq_{M_e})$ by declaring that (i) ``$\leq_{M_e}\big|_M$ equals $\leq_M$'' and (ii) $m\leq_{M_e}e$ (i.e., $m<_{M_e}e$) for all $m\in M$, which contradicts the maximality of $(M,\leq_M)$ in $P$.
\end{proof}

This corollary implies the \ul{axiom of choice for sets}, i.e., if $I$ is a set and $\{s_i~|~i\in I\}$ a collection of nonempty (disjoint) sets, then there exists a map
\bea
\textstyle f:I\ra\bigcup_{i\in I}s_i,~~i\mapsto f(i):=\min s_i~\in~s_i,\nn
\eea
where $\min s_i$ is based on some well-ordering of the set $\bigcup_{i\in I}s_i$. Since the axiom of choice for sets in turn leads to Zorn's lemma, it follows that Zorn's lemma is equivalent to the axiom of choice for sets.   
\section{Nets and Tychonoff's theorem}\label{NetsTych}
\begin{dfn}[\textcolor{blue}{\index{Net}{Net}, Tail of a net, Limit of a net, \index{Convergent net}{Convergent net}, \index{Cluster point}{Cluster point} of a net}]
Let $X$ be a space and $I$ a \ul{directed set}. A net in $X$ is a map of the form $f:I\ra X$.

A \ul{tail} of a net $I\sr{f}{\ra} X$ is the image $f\big(U_{i_0}\big)=\{f(i)\}_{i\geq i_0}$ of a set of the form  $U_{i_0}=\{i\in I:i\geq i_0\}\eqv[i_0\leq I]$ for some $i_0\in I$. (Note we can have two tails, none containing the other. However, any finite collection of tails $f([i_1\leq I])$, ..., $f([i_n\leq I])$ has a common intersection point, since every finite set in $I$, hence $\{i_1,...,i_n\}$, has an upperbound.)

A point $x_0\in X$ is a \ul{limit} of a net $I\sr{f}{\ra}X$ (written $x_0\in\lim f$) if every neighborhood of $x_0$ contains a tail of $f$. If $x_0$  is a limit of $f$, we also say $f$ \ul{converges to} $x_0$ (written $f\ra x_0$).

If a net $f:I\ra X$ converges to at least one point (i.e., $\lim f\neq\emptyset$) we say $f$ is \ul{convergent}. If every limit of $f:I\ra X$ lies in $A\subset X$ (i.e., $\lim f\subset A\subset X$) we say $f$ \ul{converges in} $A$.

A point $x_0\in X$ is a \ul{cluster point} of a net $I\sr{f}{\ra}X$ if each neighborhood of $x_0$ intersects every tail of $f$. (Note that a limit of a net is a cluster point of the net).
\end{dfn}
\begin{rmks*}
(i) If $X$ is Hausdorff and $f:I\ra X$ converges to a point $x_0\in X$, then the limit $x_0$ is \ul{unique}: otherwise, if $x_1\neq x_0$ is another limit, then $x_0$ and $x_1$ have disjoint neighborhoods which must, however, intersect because the two tails of $f$ contained in these neighborhoods must intersect (a contradiction).

(ii) On the other hand, if $X$ is a nontrivial indiscrete space, in the sense that $|X|\geq 2$ and $\T_X=\{\emptyset,X\}$, then all points have common neighborhoods and so any net in $X$ converges to more than one point
\end{rmks*}

Note that directed index sets suffice for characterizing closedness (hence topology) because the poset $\big(\P(X),\leq\big):=\big(\P(X),\supseteq\big)$ is a directed set, and in particular, for each $x\in X$, the neighborhood base $\B_x$ at $x$ is indexed as a directed set $(\B_x,\leq):=(\B_x,\supseteq)$.

\begin{lmm}
Let $X$ be a space and $A\subset X$. Then $A$ is closed $\iff$ every convergent net $f:I\ra A$ converges in $A$.
\end{lmm}
\begin{proof}
If $A$ is closed and $f:I\ra A$ converges to a point $x_0\in X$, then $x_0\in A$, otherwise a neighborhood of $x_0$ will miss a tail of $f$. Conversely, assume every convergent net $f:I\ra A$ converges in $A$. Suppose $A$ is not closed. Then there is a point $x_1\in\ol{A}-A$. So, $N(x_1)\cap A\neq\emptyset$ for every neighborhood $N(x_1)$ of $x_1$. Let $I$ consist of the neighborhoods $N(x_1)$ ordered by containment $\supseteq$ (i.e., {\small $N_1(x_1)\leq N_2(x_1)$ iff $N_1(x_1)\supseteq N_2(x_1)$}). Then we can pick a net {\small $f:I\ra A$, $N\mapsto f(N)\in N\cap A\backslash\{x_1\}$} that clearly converges to $x_1\not\in A$ (a contradiction).
\end{proof}

\begin{dfn}[\textcolor{blue}{\index{Eventual containment}{Eventual containment}, \index{Frequent containment}{Frequent containment}, \index{Cofinal set}{Cofinal set}}]
Let $f:I\ra X$ be a net and $A\subset X$. $f$ is \ul{eventually} in $A$ if $A$ contains a tail of $f$. (Recall that if $f$ is eventually in each $A_1,...,A_n\subset X$ then $A_1\cap \cdots\cap A_n\neq\emptyset$, since any finite collection of tails $f([i_1\leq I])$,...,$f([i_n\leq I])$ has a common intersection point, namely, $f(i)$ for a common upper bound $i\geq\{i_1,...,i_n\}$ which exists because $I$ is directed.)

$f$ is \ul{frequently} in $A$ if $f$ is not eventually in $X-A$, i.e., if every tail of $f$ intersects $A$  (equivalently, for each $i\in I$, there is $j\geq i$ such that $f(j)\in A$).

A set $K\subset I$ is \ul{cofinal} (or \ul{supremal}) if for each $i\in I$, there exists $k\in K$ such that $k\geq i$. For \ul{example}, if $f:I\ra X$ is frequently in $A\subset X$, then $C(A):=f^{-1}(A)\subset I$ is cofinal.
\end{dfn}
It may be convenient to think of a cofinal set $K\subset I$ as an exhaustive set of upper bounds in $I$ (i.e., $K$ is cofinal $\iff$ every element of $I$ has an upper bound in $K$).

\begin{rmks*}
(i) A cofinal set $K\subset I$ is directed. Indeed, any $k_1,k_2\in K\subset I$ have an upper bound $i\in I$, while there exists $k\in K$ such that $k\geq i$, and so $\{k_1,k_2\}\leq i\leq k$.
{\flushleft (ii)} $f:I\ra X$ is frequently in $A\subset X$ $\iff$ $f(K)\subset A$ for some cofinal set $K\subset I$ (i.e., $A$ contains the image of a cofinal set). Indeed, if $f$ is frequently in $A$, then we can set $K=f^{-1}(A)$. Conversely, if $f(K)\subset A$ for some cofinal set $K\subset I$, then for any $i\in I$, there is $k\in K$ such that $k\geq i$ and $f(k)\in A$.
{\flushleft (iii)} A net $I\sr{f}{\ral}X$ converges to $x_0\in X$ $\iff$ $f$ is eventually in every neighborhood of $x_0$.
{\flushleft (iv)} A point $x_0\in X$ is a cluster point of a net $I\sr{f}{\ral}X$ $\iff$ $f$ is frequently in every neighborhood of $x_0$.
\end{rmks*}

\begin{dfn}[\textcolor{blue}{\index{Cofinal map}{Cofinal map}, \index{Subnet}{Subnet}}]
A map of directed sets $\phi:I\ra J$ is \ul{cofinal} if for any $j\in J$, there exists $i=i_j\in I$ such that
\bea
i'\geq i~~\Ra~~\phi(i')\geq j,~~~~\txt{for all}~~i'\in I ~~\Big(\txt{i.e.,}~\phi([i\leq I])\subset[j\leq J]~\Big).\nn
\eea
(That is, a cofinal map $\phi:I\ra J$ is a map that ensures (i) each upper set in $J$ contains the image of an upper set from $I$, and (ii) the image $\phi(I)$ is a \ul{cofinal}, and therefore \ul{directed}, set in $J$.)

{\flushleft Let} $I\sr{f}{\ral}X$, $J\sr{g}{\ral}X$ be nets. Then $f$ is a \ul{subnet} of $g$ (written $f\subset g$) if $f$ factors through $g$ by a cofinal map, in the sense {\footnotesize $f=g\circ\phi:I\sr{\phi}{\ral}J\sr{g}{\ral}X$} for some cofinal map $\phi:I\ra J$.\hspace{1cm}
\bea\bt
I\ar[d,dashed,"\phi"']\ar[drr,"f=g\circ \phi"] &&\\
J\ar[rr,"g"] && X
\et\nn\eea
\end{dfn}

\begin{lmm}
(i) A net $g:J\ra X$ is eventually in $A\subset X$ $\iff$ every subnet {\footnotesize $f=g\circ\phi:I\sr{\phi}{\ral}J\sr{g}{\ral}X$} is eventually in $A$. (ii) A net $g:J\ra X$ is frequently in $A\subset X$ $\iff$ some subnet {\footnotesize $f_A=g\circ\phi_A:I\sr{\phi_A}{\ral}J\sr{g}{\ral}X$} is eventually in $A$.
\end{lmm}
\begin{proof}
(i) By the construction of a subnet, if $g:J\ra X$ is eventually in $A\subset X$, then every subnet $f=g\circ\phi\subset g$ is also eventually in $A$. The converse is also clear since $g\subset g$.

(ii) By construction of a subnet, if some subnet $f=g\circ\phi\subset g$ is eventually in $A\subset X$, then $g$ is frequently in $A$ (because the image of $\phi$ is a cofinal set). Conversely, if $g:J\ra X$ is frequently in $A\subset X$, then the subnet {\small $f_A:=g|_{I=g^{-1}(A)}:I\ra X$} of $g$ is eventually in $A$ (recall that if $g$ is frequently in $A$, then {\small $I=g^{-1}(A)\subset J$} is directed, and so $g|_I$ is a subnet of $g$).
\end{proof}

\begin{lmm}[\textcolor{OliveGreen}{\cite[Lemma 5, p.70]{kelley1975}}]
Let $f:I\ra X$ be a net and $\A\subset\P(X)$ a family of subsets such that (i) $f$ is frequently in each $A\in\A$, and (ii) $A_1,A_2\in\A$ implies $A_1\cap A_2\in\A$. Then there is a subnet {\footnotesize $g=f\circ\phi:J\sr{\phi}{\ral}I\sr{f}{\ral}X$} that is eventually in each $A\in\A$.
\end{lmm}
\begin{proof}
By hypotheses, $(\A,\supset)$ is a directed set. Let $K:=\{(i,A):i\in I,~A\in\A,~f(i)\in A\}\subset I\times\A$, where $I\times\A$ is a directed set under the \ul{product ordering}: $(i,A)\leq (j,B)$ iff $i\leq j$ and $A\supset B$. $I\times\A$ is directed because if $(i,A),(j,B)\in I\times\A$, then some $C\subset A\cap B$, $k\geq \{i,j\}$, and so $(k,C)\geq\{(i,A),(j,B)\}$.

Also $K\subset I\times\A$ is directed, because if $(i,A),(j,B)\in K$, then $f(i)\in A$, $f(j)\in B$, and some $C\subset A\cap B\in\A$ implies there is $k\geq\{i,j\}$ such that $f(k)\in C$, and so $(k,C)\geq \{(i,A),(j,B)\}$.

Define $\phi:K\ra I$ by $\phi(i,A):=i$. Then $\phi$ is cofinal, since for any $i\in I$, we can pick any $(i,A)\in K$, and for any $(i',A')\in K$ satisfying $(i',A')\geq (i,A)$, i.e., $i'\geq i$ and $A'\subset A$, we have $\phi(i',A')=i'\geq i$. Thus, {\small $f\circ\phi:K\sr{\phi}{\ral}I\sr{f}{\ral}X$} is a subnet of $f$.

To show $f\circ \phi$ is eventually in every member of $\A$, fix $A\in\A$. Pick $i\in I$ such that $f(i)\in A$, i.e., $(i,A)\in K$. Then for any $(j,B)\in K$ satisfying $(j,B)\geq (i,A)$, i.e., $j\geq i$ and $B\subset A$, we have {\small $f\circ\phi(j,B)=f(j)\in B\subset A$}. That is, $f\circ\phi\big([(i,A)\leq K]\big)\subset A$.
\end{proof}

\begin{thm}[\textcolor{OliveGreen}{\cite[Theorem 6, p.71]{kelley1975}}]
A point $x_0\in X$ is a cluster point of a net $f:I\ra X$ $\iff$ a subnet of $f$ converges to $x_0$.
\end{thm}
\begin{proof}
($\Ra$): Assume $x_0\in X$ is a cluster point of $f$. Then the set of neighborhoods $\A$ of $x_0$ satisfies the hypothesis of the above lemma, and so a subnet of $f$ converges to $x_0$.

($\La$): Assume $x_0\in X$ is not a cluster point of $f$. Then there is a neighborhood $N(x_0)$ of $x_0$ such that $f$ is not frequently in $N(x_0)$. This means $f$ is eventually in $X-N(x_0)$. It follows that every subnet of $f$ is eventually in $X-N(x_0)$, and so cannot converge to $x_0$.
\end{proof}

\begin{crl}\label{NetCompactCrl}
A space $X$ is compact $\iff$ every net in $X$ has a convergent subnet (equivalently, every net in $X$ has a cluster point).
\end{crl}
\begin{proof}
($\Ra$): Assume $X$ is compact. Suppose $f:I\ra X$ is a net with no cluster point, i.e., for each $x\in X$, there exists a neighborhood $N(x)$ of $x$ such that $f$ is not frequently in $N(x)$, i.e., $f$ is eventually in $X-N(x)$. Since $X$ is compact, we have a finite cover $X\subset\bigcup_{i=1}^nN_i(x_i)$. This means $\bigcap_{i=1}^n(X-N_i(x_i)\big)=\emptyset$. But, since $f$ is eventually in each $X-N_i(x_i)$, we also have $\bigcap_{i=1}^n(X-N_i(x_i)\big)\neq\emptyset$ (a contradiction).

($\La$): Assume $X$ is not compact. Then $X$ has an open cover $\{U_\al\}_{\al\in A}$ with no finite subcover. Let $I:=\{\txt{finite}~F\subset A\}$, which is a directed set under inclusion $\subset$ (i.e., $F_1\leq F_2$ iff $F_1\subset F_2$). Let $f:I\ra X$, $F\mapsto f(F)\in X\backslash\bigcup_{\al\in F}U_\al$ (which is well defined by hypotheses). Suppose $f$ has a cluster point $x_0\in U_{\al_0}$ (for some $\al_0\in A$). Then for any $F\in I$ there is $F'\geq F$ (i.e., $F'\supset F$) such that $f(F')\in U_{\al_0}\cap\big(X\backslash\bigcup_{\al\in F'}U_\al\big)$. In particular, with $F=\{\al_0\}$, any choice $F'\supset F=\{\al_0\}$ gives $f(F')\in\emptyset$ (a contradiction).
\end{proof}

\begin{thm}[\textcolor{OliveGreen}{\index{Tychonoff's theorem}{Tychonoff's theorem}: \cite{chernoff-1992}}]\label{TychThm}
If $X_\al$ are compact spaces, then so is $\prod X_\al$.
\end{thm}
\begin{proof}
Let $\{(X_{\al},\T_{\al})\}_{\al\in A}$ be a family of nonempty compact spaces, and for any $B\subset A$, let $X_B:=\prod_{\al\in B}X_{\al}=\{(x_{\al})_{\al\in B}:x_{\al}\in X_{\al}\}$ be given the product topology. Recall that a base-neighborhood of a point $x\in X_B$ has the form $N_F(x)$ for some finite set $F\subset B$, with
{\footnotesize\bea
\label{nbd-supp-eqn}\textstyle N_F(x)=\{y\in X_B:y_{\al}\in N(x_{\al})\in\T_{\al}~~\txt{for}~\al\in F\}=\prod_{\al\in F}N(x_{\al})\times\prod_{\beta\not\in F}X_\beta=\bigcap_{\al\in F}N_{\al}(x),
\eea}
where {\small $N_{\al}(x):=\left\{y\in X:y_{\al}\in N(x_{\al})\in\T_{\al}\right\}=N(x_{\al})\times\prod_{\beta\neq \al}X_\beta$} is a ``strip through $N(x_{\al})$''.

Let $f:I\ra X_A,~\al\mapsto f_\al$ be a net. We need to show $f$ has a cluster point. A \ul{partial cluster point} of $f$ is a cluster point of the net {\footnotesize $f_B:I\sr{f}{\ral}X_A\sr{p^B}{\twoheadrightarrow} X_B$}, for some $B\subseteq A$ (where $p^B(x)=p^B\big((x_\al)_{\al\in A}\big)=(x_\al)_{\al\in B}=x|_B$ is the projection). Let $\P=\{\txt{all partial cluster points of $f$}\}$, which is a poset under the ordering
\bea
\textstyle x_B\in X_B\leq x_{B'}\in X_{B'}~~~\txt{if}~~~B\subseteq B'~~\txt{and}~~x_{B'}\big|_B=x_B~~~~(\txt{Recall that}~~ B\sr{x_B}{\ral}\bigcup_{\al\in B}X_{\al}).\nn
\eea
$\P$ is nonempty because for any $\al\in A$, with $B:=\{\al\}$, the net \bt[column sep=small] f_{\{\al\}}:I\ar[r,"f"]& X_A\ar[r,two heads,"p^{\{\al\}}"]&X_{\{\al\}}\et has a cluster point since $X_\al\cong X_{\{\al\}}$ is compact (Corollary \ref{NetCompactCrl}).

Let {\footnotesize$L=\big\{x^{(\ld)}_{B_\ld}\in X_{B_\ld}:\ld\in\Ld\big\}$} be a chain (linearly ordered set) in $\P$, where $x^{(\ld)}_{B_\ld}$ is a cluster point of {\footnotesize $f_{B_\ld}:I\sr{f}{\ral}X_A\sr{p^{B_\ld}}{\twoheadrightarrow} X_{B_\ld}$} (i.e., $f_{B_\ld}$ is frequently in every neighborhood of $x^{(\ld)}_{B_\ld}$).

Define {\footnotesize$z:=\bigcup_{\ld\in\Ld}x^{(\ld)}_{B_\ld}:=\bigcup_{\ld\in \Ld}\big(x^{(\ld)}_{\al}:\al\in B_\ld\big)\in X_{\bigcup B_\ld}$} and {\footnotesize$g:=f_{\bigcup B_\ld}:I\sr{f}{\ral}X_A\sr{p^{\bigcup B_\ld}}{\twoheadrightarrow} X_{\bigcup B_\ld}$}, where
(i) {\footnotesize
$f_{\bigcup B_\ld}(i)=p^{\bigcup B_\ld}\circ f(i)=\bigcup_\ld p^{B_\ld}\circ f(i)=\bigcup_\ld f_{B_\ld}(i)~~\Ra~~f_{\bigcup B_\ld}=\bigcup f_{B_\ld}$}, (ii) $x\in X_{\bigcup B_\ld}$ $\iff$ $x|_{B_\ld}\in X_{B_\ld}$ for all $\ld$, and (iii) with the subspace ``imbedding'' $X_{B_\ld}\subset X_{\bigcup B_\ld}$, for any $x\in X_{\bigcup B_\ld}$ each base-neighborhood {\footnotesize $N(x)=\bigcup_\ld\big(N(x)\cap X_{B_\ld}\big)$} takes the form {\footnotesize $N(x)=\prod_{\al\in F}N(x_{\al})\times\prod_{\beta\in\left(\bigcup B_{\ld'}\right)\backslash F}X_\beta= \big(N(x)\cap X_{B_\ld}\big)\times\prod_{\beta\in\left(\bigcup B_{\ld'}\right)\backslash B_\ld}X_\beta$} for some $\ld$ such that $F\subset B_\ld$. It follows that $z$ is a cluster point of $g$, otherwise, if $f_{\bigcup B_\ld}$ is not frequently in some base-nbd {\footnotesize $N\big(\bigcup x^{(\ld)}_{B_\ld}\big)=\prod_{\al\in F}N(z_{\al})\times\prod_{\beta\in\left(\bigcup B_\ld\right)\backslash F}X_\beta$} of $z=\bigcup x^{(\ld)}_{B_\ld}$,
\bit[leftmargin=0.5cm]
\item i.e., $f_{\bigcup B_\ld}$ is eventually in {\footnotesize $X_{\bigcup B_\ld}-N\big(\bigcup x^{(\ld)}_{B_\ld}\big)=\left[X_F-\prod_{\al\in F}N(z_{\al})\right]\times\prod_{\beta\in\left(\bigcup B_\ld\right)\backslash F}X_\beta$},
\item i.e., {\footnotesize $f_{\bigcup B_\ld}\left([i_0\leq I]\right)\subset X_{\bigcup B_\ld}-N\big(\bigcup x^{(\ld)}_{B_\ld}\big)=\left[X_F-\prod_{\al\in F}N(z_{\al})\right]\times\prod_{\beta\in\left(\bigcup B_\ld\right)\backslash F}X_\beta$},
\eit
then ~{\footnotesize $f_{B_\ld}\left([i_0\leq I]\right)\subset \left[X_F-\prod_{\al\in F}N(z_{\al})\right]\times\prod_{\beta\in B_\ld\backslash F}X_\beta=X_{B_\ld}-N\big(x^{(\ld)}_{B_\ld}\big)$}~ for some $\ld$ such that $F\subset B_\ld$. That is, $f_{B_\ld}$ is not frequently in some nbd $N\big(x^{(\ld)}_{B_\ld}\big)$ of $x^{(\ld)}_{B_\ld}$, which is a contradiction.

(Note that it is enough to use base-neighborhoods as we have done, because a net that is not frequently in a given neighborhood $N(x)=\bigcup N_b(x)$ is also not frequently in some base-neighborhood $N_b(x)$.)

Thus, $L$ has an upper bound {\small$\bigcup_{\ld\in\Ld}x^{(\ld)}_{B_\ld}\in \P$}, and so by Zorn's lemma, $\P$ has a maximal element $x_B\in X_B$ for some $B\subseteq A$. It remains to show that $B=A$ (so that $x_B$ is a cluster point of $I\sr{f}{\ral}X_A$).

Suppose $B\subsetneq A$. Then some $c\in A\backslash B$. We know that $x_B$ is a cluster point of {\footnotesize $f_B:I\sr{f}{\ral}X_A\sr{p^B}{\twoheadrightarrow} X_B$}. Also, since $X_c$ is compact and nonempty, the net ~{\footnotesize $f_{\{c\}}:I\sr{f}{\ral}X_A\sr{p^{\{c\}}}{\twoheadrightarrow} X_c$}~ has a cluster point $x_c\in X_c$. Consider the point $x_{B\cup\{c\}}\in X_{B\cup\{c\}}$ given by
\bea
x_{B\cup\{c\}}=(x_B,x_c):B\cup\{c\}\ra X_{B\cup\{c\}},~~~~x_{B\cup\{c\}}(\al)=\left\{
                                                                    \begin{array}{ll}
                                                                      x_B(\al), & \al\in B \\
                                                                      x_c, & \al=c
                                                                    \end{array}
                                                                  \right\},\nn
\eea
which is is a cluster point of {\footnotesize\bt f_{B\cup\{c\}}:I\ar[r,"f"]&X_A\ar[r,"{p^{B\cup\{c\}}}"]&X_{B\cup\{c\}}\et}, because if {\footnotesize $f_B([i\leq I])\cap N(x_B)\neq\emptyset$} and {\footnotesize $f_{\{c\}}([i\leq I])\cap N(x_c)\neq\emptyset$}, then with $N(x_B\cup\{c\}):=N(x_B)\times N(x_c)$, we also have {\scriptsize $f_{B\cup\{c\}}([i\leq I])\cap N(x_{B\cup\{c\}})=\big(f_B([i\leq I])\times f_{\{c\}}([i\leq I])\big)\cap\big(N(x_B)\times N(x_c)\big)=\big(f_B([i\leq I])\cap N(x_B)\big)\times\big(f_{\{c\}}([i\leq I])\cap N(x_c)\big)\neq\emptyset$}. That is, $x_B<x_{B\cup\{c\}}\in\P$, which is a contradiction.
\end{proof}

\section{Subdifferential calculus on normed spaces}\label{SubdifCalc}
Our notation and some results such as Lemmas \ref{DiffConvFun}, \ref{SDLin}, \ref{NormSubdiff}, and \ref{SemiInnAp} are partly due to \cite{deimling,conway}. Unless stated otherwise, $X$ is a real normed space. As usual, we denote by $X^\ast$ the set of continuous \emph{linear functions/functionals} $x^\ast:X\ra\Real$ as a normed space (called \emph{dual space} of $X$) with norm $\|x^\ast\|:=\sup_{\|x\|\leq 1}|x^\ast(x)|=\sup_{\|x\|=1}|x^\ast(x)|$. If $x\in X$ and $x^\ast\in X^\ast$, we will sometimes write the number $x^\ast(x)$ as $\langle x,x^\ast\rangle$ for convenience.

\begin{dfn}[\textcolor{blue}{Fr{\'e}chet-G{\^a}teaux derivative}] A map of normed spaces $F:X\ra Y$ is (Fr{\'e}chet-) differentiable at $x\in X$ if there exists a linear map $dF_x:X\ra Y$ and a continuous map $\vep_x\in C(X,Y)$ such that
\begin{equation*}
\textstyle F(x+h)=F(x)+dF_xh+\vep_x(h)~~~~\txt{for all}~~h\in X,~~~~\txt{with}~~~~\lim\limits_{\|h\|\ra0}{\|\vep_x(h)\|\over\|h\|}=0.
\end{equation*}
The map $dF:X\ra L(X,Y)$, $x\mapsto  dF_x$ is called the (Fr{\'e}chet) derivative of $F$, and the linear map $dF_x:X\ra Y$ is called the (Fr{\'e}chet) derivative of $F$ at $x$.

    When the limit is required to exist only ``linearly'' (i.e., in one direction at a time), we get a weaker (G{\^a}teaux) version of the derivative: $F$ is G{\^a}teaux-differentiable at $x\in X$ if there exists a linear map $DF_x:X\ra Y$ and a continuous map $\vep_x\in C(X,Y)$ such that for every $h\in X$ with $\|h\|=1$,
\begin{equation*}
\textstyle F(x+th)=F(x)+~t~DF_xh+\vep_x(th),~~~~\txt{for all}~~t\in \Real,~~~~\txt{with}~~\lim\limits_{t\ra0}{\|\vep_x(th)\|\over|t|}=0.
\end{equation*}
The map $DF:X\ra L(X,Y)$, $x\mapsto DF_x$ is called the G{\^a}teaux derivative of $F$, and the map $D_hF:X\ra Y$, $x\mapsto DF_xh$ is called the directional derivative of $F$ along $h$. Accordingly, the linear map $DF_x:X\ra Y$ is called the G{\^a}teaux derivative of $F$ at $x$, and the vector $DF_xh\in Y$ is called the directional derivative of $F$ at $x$ along $h$.
\end{dfn}

\begin{rmk}[\textcolor{OliveGreen}{Chain rule}]\label{ChainRule} Let ~$O\subset X\sr{F}{\ral}Y$~ and ~$O'\subset Y\sr{G}{\ral}Z$~ be maps of normed spaces such that $F$ is differentiable at $x_0\in O$ and $G$ is differentiable at $y_0:=F(x_0)\in O'$. Then $O\subset X\sr{G\circ F}{\ral}Z$ is differentiable at $x_0$, and~ $d(G\circ F)_{x_0}=dG_{F(x_0)}\circ dF_{x_0}:X\sr{dF_{x_0}}{\ral}Y\sr{dG_{F(x_0)}}{\ral}Z$.

In particular, if ~$[0,1]\sr{u}{\ral}X\sr{\|\cdot\|}{\ral}\Real$,~ where $(X,\|\cdot\|)$ is a normed space with a differentiable norm and $u$ is a $C^1$-smooth path, then~ ${d\|u(t)\|\over dt}=d\|\cdot\|_{u(t)}\left({du(t)\over dt}\right)$~ for all $t\in[0,1]$.
\end{rmk}
\begin{proof}
For all $h\in X$,
\begin{align}
&G\circ F(x_0+h)=G\Big(F(x_0)+dF_{x_0}h+\vep_{x_0}^F(h)\Big)=G(y_0+h'),~~~~h':=dF_{x_0}h+\vep_{x_0}^F(h),\nn\\
&~~~~=G(y_0)+dG_{y_0}h'+\vep_{y_0}^G(h')\nn\\
&~~~~=G\circ F(x_0)+dG_{F(x_0)}\circ dF_{x_0}h+\vep_{x_0}(h),~~~~\vep_{x_0}(h):=dG_{y_0}\circ \vep_{x_0}^F(h)+\vep_{y_0}^G(h'),\nn\\
&~~~~~=G\circ F(x_0)+d(G\circ F)_{x_0}h+\vep_{x_0}(h),~~~~d(G\circ F)_{x_0}:=dG_{F(x_0)}\circ dF_{x_0}.\nn\qedhere
\end{align}
\end{proof}

Also, note that if $u:[0,1]\ra O\subset X$ is a differentiable path and $F:O\subset X\ra Y$ is differentiable, then
\bea
\textstyle\|F(u(t+\vep))\|-\|F(u(t))\|\leq\big\|F(u(t+\vep))-F(u(t))\big\| ~~\Ra~~{d\over dt}\big\|F\big(u(t)\big)\big\|\leq \left\|{d\over dt}F\big(u(t)\big)\right\|.\nn
\eea

\begin{rmk}[\textcolor{OliveGreen}{Mean value theorem: \cite[Theorem 1.8, p.13]{ambro-prodi1993}}]\label{MeanValThm} If $F:O\subset X\ra Y$ is G{\^a}teaux-differentiable and $O$ is open, then for any $x_1,x_2\in X$ such that ~$[x_1,x_2]:=\big\{\eta(t)=(1-t)x_1+tx_2:t\in[0,1]\big\}\subset O$,~ we have
\begin{equation*}
\textstyle\|F(x_1)-F(x_2)\|\leq C(x_1,x_2)\|x_1-x_2\|,~~\txt{where}~~C(x_1,x_2):=\sup\limits_{x\in[x_1,x_2]}\left\|DF_x\right\|.
\end{equation*}
\end{rmk}
\begin{proof}
Consider the path $\gamma:[0,1]\sr{\eta}{\ral}O\sr{F}{\ral} Y,~t\mapsto F(\eta(t))$. Let $y^\ast\in Y^\ast$ be a norming functional of $\gamma(0)-\gamma(1)=F(x_1)-F(x_2)$, i.e., $\|y^\ast\|=1$ and $y^\ast\big(\gamma(0)-\gamma(1)\big)=\|\gamma(0)-\gamma(1)\|=\|F(x_1)-F(x_2)\|$. Applying the usual mean value theorem to ~$y^\ast\circ\gamma:[0,1]\sr{\gamma}{\ral}Y\sr{y^\ast}{\ral}\Complex\cong\Real^2$,~ we get
\begin{align}
&\|F(x_1)-F(x_2)\|=|y^\ast\circ\gamma(0)-y^\ast\circ\gamma(1)|\leq|0-1|\sup_{0\leq t\leq 1}|(y^\ast\circ\gamma)'(t)|=\sup_{0\leq t\leq 1}\left|y^\ast\big(\gamma'(t)\big)\right|\nn\\
&~~~~\leq \|y^\ast\|\sup_{0\leq t\leq 1}\|\gamma'(t)\|=\sup_{0\leq t\leq 1}\left\|DF_{x(t)}(x_2-x_1)\right\|\leq\|x_1-x_2\|\sup_{0\leq t\leq 1}\left\|DF_{x(t)}\right\|\nn\\
&~~~~=\|x_1-x_2\|\sup_{x\in[x_1,x_2]}\left\|DF_x\right\|.\nn\qedhere
\end{align}
\end{proof}

\begin{dfn}[\textcolor{blue}{\index{Convex! functional}{Convex functional}}] A map ~$\vphi:O\subset X\ra\Real$~ such that $O$ is an open convex set and $\vphi(t x+(1-t)y)\leq t\vphi(x)+(1-t)\vphi(y)$ for all $x,y\in O,~~t\in [0,1]$.
\end{dfn}
We will be dealing with continuous convex functionals only.

\begin{dfn}[\textcolor{blue}{Norming functional}] A linear functional $x^\ast\in X^\ast$ is a norming functional for $x_0\in X$ if $\|x^\ast\|=1$ and $x^\ast(x_0)=\|x_0\|$. If $x^\ast$ is a norming functional of $x_0$, we will also refer to $z^\ast:=\|x_0\|x^\ast$ as a (un-normilized) norming functional of $x_0$. Note that $\|z^\ast\|=\|x_0\|$ and $z^\ast(x_0)=\|x_0\|^2$.
\end{dfn}
If $X$ is a normed space, then by Hahn-Banach theorem, every $x_0\in X$ has a norming functional. If $\H=\big(\H,{\scriptstyle\langle,\rangle}\big)$ is a Hilbert space, then by Riesz representation theorem, every $x_0\in\H$ has a unique norming functional $x_0^\ast\in \H^\ast\cong\H$ given by
~$x_0^\ast(x)=\left\langle x,x_0\right\rangle$,~ for all $x\in X$.

\begin{dfn}[\textcolor{blue}{Set-valued maps}] A set-valued map of normed spaces $X,Y$ is a map of the form
$F:X\ra \P(Y)=2^Y$, $x\mapsto Fx\subset Y$.
\end{dfn}
Note that if $A\subset X$, then ~$F(A)=\bigcup_{a\in A}Fa$. In particular, $F(X)=\bigcup_{x\in X}Fx$. Also, a set-valued map $F:X\ra 2^Y$ has ``inverse'' ~$F^{-1}:F(X)\ra 2^X$, $y\mapsto F^{-1}y:=\{x\in X:Fx\ni y\}$.

\begin{dfn}[\textcolor{blue}{\index{Duality map}{Duality map}}]\label{Duality}
 If $X$ is a normed vector space, the duality map of $X$ is the set-valued map $\F:X\ra\P(X^\ast)=2^{X^\ast}$ given by the set of \ul{un-normalized} norming functionals
~$\F x=\F_x:=\left\{x^\ast\in X^\ast:\|x^\ast\|=\|x\|,~x^\ast(x)=\|x\|^2\right\}$.
\end{dfn}
\begin{dfn}[\textcolor{blue}{Norming duality map}]\label{NormDuality} If $X$ is a normed vector space, the norming duality map of $X$ is the set-valued map $\widehat{\F}:X\ra\P(X^\ast)=2^{X^\ast}$ given by the set of norming functionals ~$\widehat{\F}x=\widehat{\F}_x:=\left\{x^\ast\in X^\ast:\|x^\ast\|=1,~x^\ast(x)=\|x\|\right\}$.
\end{dfn}
Note that ~$\widehat{\F}=\F\circ F:X\sr{F}{\ral}X\sr{\F}{\ral}2^{X^\ast}$,~ where $F(x)=\hat{x}:=x/\|x\|$ for $x\neq 0$.

\begin{dfn}[\textcolor{blue}{Semi-inner products on $X$}]\label{SemiIPs} For $x,y\in X$, we define
\bea
\langle x,y\rangle_-:=\inf\limits_{y^\ast\in\F y}\langle x,y^\ast\rangle,~~~~\langle x,y\rangle_+:=\sup\limits_{y^\ast\in\F y}\langle x,y^\ast\rangle.\nn
\eea
\end{dfn}

\begin{dfn}[\textcolor{blue}{Addition and scalar multiplication of sets}] If $A,B\subset X$ (a vector space) and $\ld$ is a scalar, we write ~$\ld A:=\{\ld a:a\in A\}$~ and ~$A+B:=\{a+b:a\in A,b\in B\}$.
\end{dfn}

\begin{dfn}[\textcolor{blue}{Strong ordering on subsets of real numbers}] If $A,B\subset\Real$, we write $A\leq B$ if $a\leq b$ for all $a\in A$, $b\in B$. (In particular, we will write $0\leq A$ if $\{0\}\leq A$~). We write $A\leq B$ if $-B\leq -A$, where $-A:=\{-a:a\in A\}$ and $-B:=\{-b:b\in B\}$.
\end{dfn}

\begin{dfn}[\textcolor{blue}{\index{Subgradient}{Subgradient} of a functional at a point}] Let $f:D\subset X\ra\Real$ and let $x_0\in D$. A linear functional $x^\ast\in X^\ast$ is a subgradient of $f$ at $x_0$ if
\bea
f(x)\geq f(x_0)+x^\ast(x-x_0)~~~~\txt{for all}~~~~x\in D.\nn
\eea
\end{dfn}

\begin{dfn}[\textcolor{blue}{\index{Subderivative}{Subderivative} of a functional}] The subderivative of a functional $f:D\subset X\ra\Real$ is the set-valued map $\del f:D\subset X\ra 2^{X^\ast}$,~ $\del f(x)=\del f_x:=\big\{x^\ast\in X^\ast:~x^\ast~\txt{is a subgradient of $f$ at $x$}\big\}$, i.e.,
\bea
\del f(x)=\del f_x:=\left\{x^\ast\in X^\ast:~f(y)\geq f(x)+x^\ast(y-x)~\txt{for all}~y\in D\right\}.\nn
\eea
\end{dfn}
Note that this definition implies ~$\del f_x\subset\{x^\ast\in X:\nabla^-f_x\leq x^\ast\leq\nabla^+f_x\}$,~ where
\bea
\textstyle\nabla^-f_x(y):=\lim\limits_{t\uparrow 0}{f(x+ty)-f(x)\over t},~~~~\nabla^+f_x(y):=\lim\limits_{t\downarrow 0}{f(x+ty)-f(x)\over t}.\nn
\eea

\begin{dfn}[\textcolor{blue}{\index{Subdifferentiable functional}{Subdifferentiable functional}}] A functional $f:D\subset X\ra\Real$ is subdifferentiable at $x_0\in D$ if $\del f_{x_0}\neq\emptyset$. $f$ is subdifferentiable if $f$ is subdifferentiable at every point $x_0\in D$.
\end{dfn}

\begin{lmm}[\textcolor{OliveGreen}{Subdifferentiability of convex functionals}]\label{DiffConvFun}
Every convex functional $\vphi:X\ra\Real$ is subdifferentiable, and its subderivative satisfies ~$\del\vphi_x=\{x^\ast\in X^\ast:\nabla^-\vphi_x\leq x^\ast\leq\nabla^+\vphi_x\}$.
\end{lmm}
\begin{proof}
Since the function $t\in\Real\mapsto\vphi(x+th)\in\Real$ is convex for all $x,h\in X$, the (difference quotient) function $\vphi_x(h,t):={\vphi(x+th)-\vphi(x)\over t}$ is increasing in $t$, and the one-sided directional derivatives
\begin{align}
&\textstyle\nabla^-\vphi_x(h):=\lim_{t\uparrow 0}{\vphi(x+th)-\vphi(x)\over t}~~~~\txt{(Left derivative at $x$ in the direction of $h$)}\nn\\
&\textstyle~~~~=\sup_{t<0}{\vphi(x+th)-\vphi(x)\over t}=\lim_{t\ra 0_-}{\vphi(x+th)-\vphi(x)\over t},\nn\\
&\textstyle\nabla^+\vphi_x(h):=\lim_{t\downarrow 0}{\vphi(x+th)-\vphi(x)\over t}~~~~\txt{(Right derivative at $x$ in the direction of $h$)}\nn\\
&\textstyle~~~~=\inf_{t>0}{\vphi(x+th)-\vphi(x)\over t}=\lim_{t\ra 0_+}{\vphi(x+th)-\vphi(x)\over t},\nn
\end{align}
exist and have the following properties.
\bit[leftmargin=0.9cm]
\item[(i)] $\nabla^-\vphi_x(h)\leq \nabla^+\vphi_x(h)$~ and~ $\nabla^-\vphi_x(h)=-\nabla^+\vphi_x(-h)$. (These follow from definitions). Also,
    $\nabla^\pm\vphi_x(h)=\nabla^{\pm}\ol{\vphi}_{-x}(-h)=-\nabla^\mp\ol{\vphi}_{-x}(h)$, where $\ol{\vphi}(x):=\vphi(-x)$.
\item[(ii)] $\nabla^+\vphi_x(\ld h)
=\left\{
   \begin{array}{ll}
     \ld \nabla^+\vphi_x(h), & \ld\geq0 \\
     \ld\nabla^-\vphi_x(h), & \ld<0
   \end{array}
 \right\}\geq \ld \nabla^+\vphi_x(h).
$
\item[(iii)] $\nabla^+\vphi_x(h+h')\leq \nabla^+\vphi_x(h)+\nabla^+\vphi_x(h')$, ~because for all~ $s\geq t>0$,
{\small\begin{align}
&\textstyle\nabla^+\vphi_x(h+h')\leq{\vphi\big(x+t(h+h')\big)-\vphi(x)\over t}\sr{\txt{convexity}}{\leq} {\vphi\big(x+2th)-\vphi(x)\over 2t}+{\vphi\big(x+2th')-\vphi(x)\over 2t}\nn\\
&\textstyle~~~~\sr{\txt{monotonicity}}{\leq}{\vphi\big(x+2th)-\vphi(x)\over 2t}+{\vphi\big(x+2sh')-\vphi(x)\over 2s}.\nn
\end{align}}Thus, the map ~$\nabla^+\vphi_x:X\ra \Real$~ is a sublinear functional.
\item[(iv)] If $\vphi=\|\cdot\|$, then $|\nabla^\pm\vphi_x(h)|\leq |\nabla^\pm\vphi_0(h)|=\|h\|$.
\eit
Now, for a fixed $h_0\in X$, the linear functional $x_0^\ast:\Real h_0\ra\Real,~\ld h_0\ra\ld\nabla^+\vphi_x(h_0)$ satisfies $x^\ast_0\leq \nabla^+\vphi_x$ on $\Real h_0$ and so by Hahn-Banach Theorem (for locally convex spaces), there exists a linear functional $x^\ast:X\ra\Real$ such that $x^\ast\leq\nabla^+\vphi_x$ on $X$. Moreover, $\nabla^-\vphi_x(h)=-\nabla^+\vphi_x(-h)\leq -x^\ast(-h)=x^\ast(h)$,
\bea
~~~~\Ra~~~~\nabla^-\vphi_x\leq x^\ast\leq \nabla^+\vphi_x~~~~\txt{on}~~X.~~~~~~~~~~~~~~~~(\ast)\nn
\eea
Also, $x^\ast\in X^\ast$, because for $t<0<s$, we have
\bea
&&\textstyle{\vphi\big(x+th\big)-\vphi(x)\over t}\leq\nabla^-\vphi_x(h)\leq x^\ast(h)\leq\nabla^+\vphi_x(h)\leq{\vphi\big(x+sh\big)-\vphi(x)\over s},\nn\\
&&~~\Ra~~x^\ast(h)\ra 0~~\txt{as}~~h\ra 0~~\txt{(since $\vphi$ is continuous)}.\nn
\eea
By setting $h=y-x$ and $s=1$, in the right inequality above, we see that for all $x\in D$, there exists a linear functional $x^\ast=\vphi^\ast[x]\in X^\ast$ (depending on $x$) such that
\bea
\vphi(y)\geq \vphi(x)+\vphi^\ast[x](y-x)~~~~\txt{for all}~~y\in D.~~~~~~~~~~~~~~~~(\ast\ast)\nn
\eea
This proves that $\vphi$ is subdifferentiable. Moreover, it is clear from the above discussion (due to monotonicty of the difference quotient $\vphi_x(h,t)={\vphi(x+th)-\vphi(x)\over t}$ in $t$) that an $x^\ast\in X^\ast$ satisfies $(\ast)$ $\iff$ it satisfies $(\ast\ast)$.
\end{proof}

Note that $(\ast\ast)$ above says that for all $x,y\in D$, there are functionals $\vphi^\ast[x],\vphi^\ast[y]\in X^\ast$ such that
\bea
&&\nabla^-\vphi_y(x-y)\leq \vphi^\ast[y](x-y)\leq\vphi(x)-\vphi(y)\leq \vphi^\ast[x](x-y)\leq\nabla^+\vphi_x(x-y),\nn\\
&&~~\Ra~~\left(\nabla^+\vphi_x-\nabla^-\vphi_y\right)(x-y)\geq 0,~~~~\left(\vphi^\ast[x]-\vphi^\ast[y]\right)(x-y)\geq 0.\nn
\eea

\begin{lmm}[\textcolor{OliveGreen}{Partial linearity of convex subderivative}]\label{SDLin}
If $\vphi,\psi:X\ra\Real$ are convex functionals, then ~$\del(\ld \vphi)_x=\ld~\del \vphi_x$ and $\del(\vphi+\psi)_x=\del\vphi_x+\del\psi_x$,~ for all $x\in X,~~\ld\geq0$.

Moreover, ~$\del\vphi_x=-\del\ol{\vphi}_{-x}$,~ where~ $\ol{\vphi}(x):=\vphi(-x)$.
\end{lmm}
\begin{proof}
The cases $\del(\ld \vphi)_x=\ld~\del \vphi_x$ and $\del(\vphi+\psi)_x\supseteq\del\vphi_x+\del\psi_x$ follow from the definitions and the fact that ~{\footnotesize $\del\vphi_x=\{x^\ast\in X^\ast:\nabla^-\vphi_x\leq x^\ast\leq\nabla^+\vphi_x\}$, $\del\psi_x=\{x^\ast\in X^\ast:\nabla^-\psi_x\leq x^\ast\leq\nabla^+\psi_x\}$}.

So, we will prove $\del(\vphi+\psi)_x\subset\del\vphi_x+\del\psi_x$. Let $x^\ast\in \del(\vphi+\psi)_x$. Then
\bea
&&(\vphi+\psi)(y)\geq(\vphi+\psi)(x)+x^\ast(y-x)~~~~\txt{for all}~~~~y\in X,\nn\\
&&~~\Ra~~\psi(y)-\psi(x)\geq\vphi(x)-\vphi(y)+x^\ast(y-x)~~~~\txt{for all}~~~~y\in X.\nn
\eea
Thus, the convex sets {\footnotesize $A=\{(y,t):\psi(y)-\psi(x)\geq t\}$} and {\footnotesize $B=\{(y,t):t\geq\vphi(x)-\vphi(y)+x^\ast(y-x)\}$} can be separated by a hyperplane in $X\times\Real$ (See \cite[Theorem 3.7, p.110]{conway}), i.e., there exist $f\in X^\ast$ and $\al\in\Real$ such that the linear map {\footnotesize $F(y,t)=-f(y)+t:X\times\Real\ra\Real$} satisfies
\bea
F(A)\geq\al\geq F(B),~~~~\txt{i.e.,}~~A\subset\{F\geq\al\}~~\txt{and}~~B\subset\{\al\geq F\}.\nn
\eea
The condition $A\subset\{F\geq\al\}$ means for all {\footnotesize $(y,t)$,~ $\psi(y)-\psi(x)\geq t~\Ra~-f(y)+t\geq\al$}, and so
\bea
\label{SepEq1}\psi(y)-\psi(x)\geq f(y)+\al~~~~\txt{for all}~~~~y\in X.
\eea
The condition $B\subset\{\al\geq F\}$ means for all {\footnotesize$(y,t)$,~ $t\geq\vphi(x)-\vphi(y)+x^\ast(y-x)~\Ra~\al\geq-f(y)+t$},
\bea
\label{SepEq2}~~\txt{and so}~~~~f(y)+\al\geq\vphi(x)-\vphi(y)+x^\ast(y-x)~~~~\txt{for all}~~~~y\in X.
\eea
It follows from (\ref{SepEq1}) and (\ref{SepEq2}) that
\bea
\psi(y)-\psi(x)\geq f(y)+\al\geq\vphi(x)-\vphi(y)+x^\ast(y-x)~~~~\txt{for all}~~~~y\in X.\nn
\eea
At the point $y=x$, we see that $\al=- f(x)$, and so $f\in\del\psi_x$ and $x^\ast-f\in\del\vphi_x$, since
\bea
\psi(y)-\psi(x)\geq f(y-x)\geq\vphi(x)-\vphi(y)+x^\ast(y-x)~~~~\txt{for all}~~~~y\in X.\nn
\eea
Hence, $x^\ast=(x^\ast-f)+f\in\del\vphi_x+\del\psi_x$, that is, $\del(\vphi+\psi)_x\subset \del\vphi_x+\del\psi_x$.
\end{proof}

\begin{lmm}[\textcolor{OliveGreen}{Subderivative of the norm}]\label{NormSubdiff}
The subderivative of the norm is the norming duality map (Definition \ref{NormDuality}), i.e., $\del\|\cdot\|_x=\widehat{\F}_x$, for all $x\in X$.
\end{lmm}
\begin{proof}
Let $\vphi:X\ra \Real$ be the convex functional given by $\vphi(x)=\|x\|$. Then
{\small
\bea
&&\del\vphi_x=\left\{x^\ast\in X^\ast~:~\|y\|\geq\|x\|+x^\ast(y-x)~~\txt{for all}~~y\in X\right\}\nn\\
&&~~~~=\left\{x^\ast\in X^\ast~:~x^\ast(y)\leq\|x+y\|-\|x\|~~\txt{for all}~~y\in X\right\}\nn\\
&&~~~~=\left\{x^\ast\in X^\ast~:~x^\ast(y)\geq\|x\|-\|x-y\|~~\txt{for all}~~y\in X\right\}\sr{(s)}{=}\widehat{\F}_x,\nn
\eea}where the proof of step (s) is as follows.
\bit[leftmargin=0.7cm]
\item $\sr{(s)}{\subset}$: Let $x^\ast\in\del\vphi_x$. Then with $y=0$, we get $x^\ast(x)\geq\|x\|$, and so $\|x^\ast\|\geq 1$. On the other hand, since $x^\ast(y)\leq\|y\|$ for all $y$, setting $y=\txt{sign}\big(x^\ast(z)\big)z$, we get $|x^\ast(z)|\leq\|z\|$ for all $z\in X$, and so $\|x^\ast\|\leq1$. It follows that $\|x^\ast\|=1$ and $x^\ast(x)=\|x\|$, i.e., $x^\ast\in\widehat{\F}_x$.
\item $\sr{(s)}{\supseteq}$: If $x^\ast\in\widehat{\F}_x$, then~ {\small$\|x\|+x^\ast(y-x)=\|x\|+(x^\ast(y)-\|x\|)=x^\ast(y)\leq\|y\|$},~ and so $x^\ast\in\del\vphi_x$.\qedhere
\eit
\end{proof}

\begin{lmm}[\textcolor{OliveGreen}{Derivative characterization of semi-inner products}]\label{SemiInnAp}
The semi-inner products (Definition \ref{SemiIPs}) are determined by one-sided derivatives of the norm as follows.
\bea
&&\textstyle\langle x,y\rangle_-=\|y\|\lim\limits_{t\uparrow 0}{\|y+tx\|-\|y\|\over t}=\|y\|~\nabla^-\|\cdot\|_y(x),\nn\\
&&\textstyle\langle x,y\rangle_+=\|y\|\lim\limits_{t\downarrow 0}{\|y+tx\|-\|y\|\over t}=\|y\|~\nabla^+\|\cdot\|_y(x).\nn
\eea
\end{lmm}
\begin{proof}
Note that $\langle x,y\rangle_-=-\langle -x,y\rangle_+$. By Lemmas \ref{DiffConvFun} and \ref{NormSubdiff}, we have
\bea
\widehat{\F}_x:=\left\{x^\ast\in X^\ast:\|x^\ast\|=1,~x^\ast(x)=\|x\|\right\}=\left\{x^\ast\in X^\ast:\nabla^-\|\cdot\|_x\leq x^\ast\leq\nabla^+\|\cdot\|_x\right\}.\nn
\eea
Now, $\widehat{\F}_x\subset B_{X^\ast}$ is $\sigma(X^\ast,X)$-closed: Indeed, if $x_\al^\ast\in\widehat{\F}_x$ is a net such that $x_\al^\ast(y)\ra x^\ast(y)$ for all $y$, then
\bea
&& \|x\|=x_\al^\ast(x)\ra x^\ast(x)~~\Ra~~x^\ast(x)=\|x\|,~~\Ra~~\|x^\ast\|\geq 1,\nn\\
&& ~~\Ra~~x^\ast(y)\leq|x^\ast(y)-x_\al^\ast(y)|+|x_\al^\ast(y)|\leq|x^\ast(y)-x_\al^\ast(y)|+\|y\|\ra \|y\|,\nn\\
&&~~\Ra~~x^\ast(y)\leq\|y\|~~~~\txt{for all}~~y,~~\Ra~~\|x^\ast\|\leq 1,\nn\\
&&~~\Ra~~\|x^\ast\|=1,~~~~x^\ast(x)=\|x\|,~~~~\Ra~~~~x^\ast\in \widehat{\F}_x.\nn
\eea
Thus, $\widehat{\F}_x$ is $\sigma(X^\ast,X)$-compact since $B_{X^\ast}$ is $\sigma(X^\ast,X)$-compact by Banach-Alaoglu theorem. This shows
\bea
&&\langle x,y\rangle_-=\|y\|\inf_{y^\ast\in\widehat{\F}y}y^\ast(x)\sr{(a)}{=}\|y\|~y_-^\ast(x)\sr{(b)}{=}\|y\|~\nabla^-\|\cdot\|_y(x),~~~~\txt{for some}~~y_-^\ast\in \widehat{\F}_y,\nn\\
&&\langle x,y\rangle_+=\|y\|\sup_{y^\ast\in\widehat{\F}y}y^\ast(x)\sr{(a)}{=}\|y\|~y_+^\ast(x)\sr{(b)}{=}\|y\|~\nabla^+\|\cdot\|_y(x),~~~~\txt{for some}~~y_+^\ast\in \widehat{\F}_y,\nn
\eea
where step (a) holds because the map $\widehat{\F}_y\ra\Real,~y^\ast\mapsto y^\ast(x)$ is continuous and so attains its maximum and minimum on $\widehat{\F}_y$, and step (b) holds because (as we already know from the proof of Lemma \ref{DiffConvFun}) we can construct a linear functional $y_\pm^\ast$ in {\small$\big\{y^\ast\in X^\ast:\nabla^-\|\cdot\|_y\leq y^\ast\leq\nabla^+\|\cdot\|_y\big\}$} that agrees with $\nabla^\pm\|\cdot\|_y$ at $x$.
\end{proof}

\begin{lmm}[\textcolor{OliveGreen}{Absolute continuity of the norm along trajectories, Chain rule II}]\label{AbsCont}~
\bit[leftmargin=0.8cm]
\item[(i)] If $X$ is a normed space and $[0,1]\sr{\gamma}{\ral}X$ is a $C^1$-smooth path, then the function $[0,1]\sr{\vphi}{\ral}\Real$, $\vphi(t)=\|\gamma(t)\|$, is absolutely continuous, i.e.,
\bea
\txt{$\vphi'$ exists a.e.,~~~~ $\vphi'\in L([0,1])$, ~~~~and ~~~~ $\vphi(t)=\vphi(0)+\int_0^t\vphi'(s)ds$.}\nn
\eea
\item[(ii)] The subderivative of $\vphi$ satisfies the chain rule ~$\del\vphi_t\subset\del\|\cdot\|_{\gamma(t)}\big(\gamma'(t)\big)$ for all $t$. Moreover,
    \bea
    \del\vphi_t=\del\|\cdot\|_{\gamma(t)}\big(\gamma'(t)\big),~~~~\txt{for almost every}~~t\in[0,1].\nn
    \eea
\eit
\end{lmm}
\begin{proof}
\bit[leftmargin=0.8cm]
\item[(i)] Let $C:=\sup_{[0,1]}\|\gamma'\|$, where $\gamma'$ is the derivative of $\gamma$. Then the mean value theorem gives
\bea
\big|\vphi(a)-\vphi(b)\big|=\big|\|\gamma(a)\|-\|\gamma(b)\|\big|\leq\big\|\gamma(a)-\gamma(b)\big\|\leq C|a-b|,\nn
\eea
which shows $\vphi$ is absolutely continuous.
\item[(ii)] At all $t$, $\vphi$ has left and right derivatives satisfying $\vphi_-'(t)\leq\vphi_+'(t)$, and
{\footnotesize\begin{align}
&\del\vphi_t\subset \Delta\vphi_t:=\left\{l\in\Real^\ast:\vphi_-'(t)\leq l\leq\vphi'_+(t)\right\}=\left\{l\in\Real^\ast:\nabla^-\|\cdot\|_{\gamma(t)}\big(\gamma'(t)\big)\leq l\leq\nabla^+\|\cdot\|_{\gamma(t)}\big(\gamma'(t)\big)\right\}\nn\\
&~~~~\supseteq\del\|\cdot\|_{\gamma(t)}\big(\gamma'(t)\big).\nn
\end{align}}
By part (i), ~$\del\vphi_t=\Delta\vphi_t=\del\|\cdot\|_{\gamma(t)}\big(\gamma'(t)\big)$~ a.e. In general, we have $\del\vphi_t\subset\del\|\cdot\|_{\gamma(t)}\big(\gamma'(t)\big)$, and the proof is as follows: Given any $l\in\del\vphi_t$, the map
\bea
x^\ast_0:\Real\gamma(t)+\Real\gamma'(t)\ra\Real,~~\al\gamma(t)+\beta\gamma'(t)\mapsto \al\|\gamma(t)\|+\beta l\nn
\eea
extends (by Hahn-Banach theorem) to an $x^\ast\in\del\|\cdot\|_{\gamma(t)}$. Indeed, $x^\ast_0\left({1\over\|\gamma(t)\|}\gamma(t)+0\gamma'(t)\right)=1$ implies $\|x_0^\ast\|\geq 1$. Meanwhile, because ~$\vphi_t(s)={\|\gamma(t)-s\gamma'(t)\|-\|\gamma(t)\|\over s}$ is increasing, we see that in the set
\bea
\del\vphi_t=\{l\in\Real^\ast:l(s)\leq\vphi(t+s)-\vphi(t)~~\txt{for all}~~s\in\Real\},\nn
\eea
the condition $l(s)\leq\|\gamma(t+s)\|-\|\gamma(t)\|$ for all $s$ is stricter than (and so implies) the condition
\bea
&&\textstyle{\|\gamma(t)-s\gamma'(t)\|-\|\gamma(t)\|\over-s}\leq\vphi_-'(t)\leq l\leq\vphi_+'(t)\leq {\|\gamma(t)+s\gamma'(t)\|-\|\gamma(t)\|\over s}~~~~\txt{for all}~~s>0,\nn\\
&&\textstyle~~\txt{or}~~\|\gamma(t)\|-\|\gamma(t)-s\gamma'(t)\|\leq l(s)\leq \|\gamma(t)+s\gamma'(t)\|-\|\gamma(t)\|~~\txt{for all}~~s>0,\nn\\
&&\textstyle~~\txt{or}~~\|\gamma(t)\|+l(s)\leq \|\gamma(t)+s\gamma'(t)\|~~\txt{for all}~~s\in\Real,\nn
\eea
which implies~ $\left|x_0^\ast\big(\al\gamma(t)+\beta\gamma'(t)\big)\right|~=~\big|~\al\|\gamma(t)\|+\beta l~\big|\leq \big\|\al\gamma(t)+\beta\gamma'(t)\big\|$,~ and so $\|x_0^\ast\|\leq 1$. \qedhere
\eit
\end{proof}

\begin{lmm}[\textcolor{OliveGreen}{Chain rule III}]
If $\vphi:X\ra\Real$ is \ul{convex} and $f:\Real\ra\Real$ is \ul{convex}, \ul{differentiable}, and $f'$ is \ul{nonnegative} on the range of $\vphi$, then ~$\del(f\circ\vphi)_x=f'_{\vphi(x)}\del\vphi_x$. (In particular, if $\vphi(x)={\|x\|^2\over 2}$, then $\del\vphi_x=\F_x$.)
\end{lmm}
\begin{proof}
Observe that for any $x^\ast\in X^\ast$, {\small$\nabla^-\vphi_x\leq x^\ast\leq \nabla^+\vphi_x$} $\iff$ {\small$\nabla^-(f\circ\vphi)_x\leq f'_{\vphi(x)}x^\ast\leq \nabla^+(f\circ\vphi)_x$}, where {\small$\nabla^\pm(f\circ\vphi)_x=f'_{\vphi(x)}\nabla^\pm\vphi_x$}.
\end{proof}

\begin{lmm}[\textcolor{OliveGreen}{Subderivative of max}]\label{SubdMax}
Given functionals $f_1,...,f_m:X\ra\Real$, define a functional $f:X\ra\Real$ by $f(x)=\max_if_i(x)$. Then~ $\txt{Conv}\left(\bigcup_{i\in I(x)}\del f_i(x)\right)\subset \del f(x)\subset\txt{Conv}\Big(\bigcup_{i=1}^m\del f_i(x)\Big)$,~ where ``Conv'' denotes ``convex hull'' and $I(x):=\big\{i\in\{1,...,m\}:f(x)=f_i(x)\big\}$ are the indices of those $f_i$ that attain the value of $f$ at $x$.
\end{lmm}
\begin{proof}
By definition, $\del f_i(x)=\{x^\ast\in X^\ast:f_i(y)\geq f_i(x)+x^\ast(y-x)~\txt{for all}~y\in X\}$, and
{\small\begin{align}
&\textstyle\txt{Conv}\left(\bigcup\limits_{i\in I(x)}\del f_i(x)\right)\sr{(a)}{\subset}\del f(x):=\left\{x^\ast\in X^\ast:\max\limits_jf_j(y)\geq \max\limits_jf_j(x)+x^\ast(y-x)~\txt{for all}~y\in X\right\}\nn\\
&\textstyle~~~~\sr{(b)}{\subset} \txt{Conv}\left(\bigcup\limits_{i=1}^m\del f_i(x)\right),\nn
\end{align}}where step (a) holds because by definition $\del f_i(x)\subset\del f(x)$ for each $i\in I(x)$, and $\del f(x)$ is a convex set. Also, step (b) holds because
{\small\begin{align}
&\textstyle\txt{Conv}\left(\bigcup\limits_{i=1}^m\del f_i(x)\right)=\txt{Conv}\left(\bigcup\limits_{i=1}^m\{x^\ast\in X^\ast:f_i(y)\geq f_i(x)+x^\ast(y-x)~\txt{for all}~y\in X\}\right)\nn\\
&\textstyle~~~~\supseteq\txt{Conv}\left(\bigcup\limits_{i=1}^m\bigcap\limits_{j=1}^m\{x^\ast\in X^\ast:f_i(y)\geq f_j(x)+x^\ast(y-x)~\txt{for all}~y\in X\}\right)\nn\\
&\textstyle~~~~=\txt{Conv}\left\{x^\ast\in X^\ast:\max\limits_if_i(y)\geq \max\limits_jf_j(x)+x^\ast(y-x)~\txt{for all}~y\in X\right\}=\del f(x).\nn \qedhere
\end{align}}
\end{proof}

\begin{dfn}[\textcolor{blue}{\index{Accretive map}{Accretive map} $X\ra X$}] A map $F:D\subset X\ra X$ is accretive if for all $x,y\in D$ there exists a norming functional $(x-y)^\ast\in \F_{x-y}$ of $x-y$ such that ~$\big\langle Fx-Fy,(x-y)^\ast\big\rangle\geq 0$. That is, $F:D\subset X\ra X$ is accretive if ~$\big\langle Fx-Fy,x-y\big\rangle_+\geq 0$ for all $x,y\in D$.
\end{dfn}

\begin{dfn}[\textcolor{blue}{\index{Strongly accretive map}{Strongly accretive map} $X\ra X$}] A map $F:D\subset X\ra X$ is strongly accretive if for any $x,y\in D$ and for any norming functional $(x-y)^\ast\in\F_{x-y}$ of $x-y$, we have ~$\big\langle Fx-Fy,(x-y)^\ast\big\rangle\geq 0$. That is, $F:D\subset X\ra X$ is strongly accretive if ~$\big\langle Fx-Fy,x-y\big\rangle_-\geq 0$ for all $x,y\in D$.
\end{dfn}

\begin{dfn}[\textcolor{blue}{\index{Monotone map}{Monotone map} $X\ra X^\ast$: \cite{deimling}}]
A map $F:X\ra X^\ast$ is monotone if
\bea
\langle x-y,Fx-Fy\rangle:=(Fx-Fy)(x-y)=ev_{x-y}(Fx-Fy)\geq 0~~~~\txt{for all}~~~~x,y\in X.\nn
\eea
\end{dfn}


\chapter{Imbedding of finite-dimensional spaces}\label{DimTh}
This chapter is based on \cite{HureWall41}. Wherever information about the spaces in question is insufficient for the desired conclusion(s), assume the spaces are separable metric spaces.
\begin{dfn}[\textcolor{blue}{\index{Topological! dimension}{(Topological) dimension} of the empty space}]
We define the (topological) dimension of the empty set (in any space) to be $\dim\emptyset:=-1$.
\end{dfn}

\begin{dfn}[\textcolor{blue}{(Topological) dimension of a nonempty space}]
Let $X$ be a (nonempty) space and $n\geq 0$ an integer.
\bit[leftmargin=0.7cm]
\item We say $X$ has dimension at most $n$ at $x\in X$ (written $\dim_xX\leq n$) if every neighborhood $U\ni x$ contains a neighborhood $V\ni x$ such that $\dim\del V\leq n-1$.
\item We say $X$ has dimension equal to $n$ at $x\in X$ (i.e., $\dim_xX=n$) if (i) $\dim_xX\leq n$ and (ii) $\dim_xX\not\leq n-1$, i.e., there exists a neighborhood $U_x\ni x$ such that $\dim_x\del V\geq n-1$ for every neighborhood $x\in V\subset U_x$.
\item We say $X$ has dimension at most $n$ (written $\dim(X)\leq n$) if $\dim_xX\leq n$ for all $x\in X$.
\item We say $\dim X=n$ if (i) $\dim X\leq n$ and (ii) $\dim X\not\leq n-1$, i.e., there exists a point $x_0\in X$ such that $\dim_{x_0}(X)=n$.
\item We say $X$ is infinite-dimensional (written $\dim X=\infty$) if $\dim X\not\leq n$ for each $n\geq 0$, i.e., for each $n\geq 0$ there exists $x=x_n\in X$ such that $\dim_xX>n$.
\eit
\end{dfn}

Since the definition of (topological) dimension depends on topology only, it is preserved by homeomorphisms (i.e., it is a topological invariant).

\begin{dfn}[\textcolor{blue}{\index{Separated sets}{Separated sets}}]
Let $X$ be a space. Sets $A,B\subset X$ are separated if $\ol{A}\cap B=\emptyset=A\cap\ol{B}$ (i.e., $A$ and $B$ are disjoint from the closure of each other).
\end{dfn}

\begin{lmm}
Let $X$ be a space. Then $A,B\subset X$ are separated (i.e., $\ol{A}\cap B=\emptyset=A\cap\ol{B}$) $\iff$ there exist open sets $U\supset A$ and $V\supset B$ such that $\ol{A}\cap V=\emptyset=U\cap\ol{B}$.
\end{lmm}
\begin{proof}
($\Ra$): If $A,B$ are separated, let $U:=\big(\ol{B}\big)^c$ and $V:=\big(\ol{A}\big)^c$. ($\La$): The converse is clear.
\end{proof}

\vspace{1.5cm}
(This space is due to the shortening of a previous version of the above proof.)
\vspace{1.5cm}

\begin{dfn}[\textcolor{blue}{\index{Completely normal space}{Completely normal space}}]
A space in which every two separated sets have disjoint neighborhoods.
\end{dfn}

\begin{lmm}
A space is completely normal $\iff$ every subspace is normal. (In particular, a metric space is completely normal, since every subspace of a metric space is a metric space.)
\end{lmm}
\begin{proof}
($\Ra$): Assume $X$ is completely normal, and let $A\subset X$. We need to show $A$ is normal. Let $C,C'\subset A$ be disjoint sets that are closed in $A$. Then {\small $\ol{C}\cap C'=\emptyset=C\cap\ol{C'}$}, since
\bea
\ol{C}\cap C'=\ol{C}\cap A\cap C'=\Cl_A(C)\cap C'=C\cap C'=\emptyset.\nn
\eea
Thus (by complete normality of $X$) there are disjoint open sets $U_1\supset C$, $V_1\supset C'$, which give disjoint sets $U=U_1\cap A\supset C$, $V=V_1\cap A\supset C'$ that are open in $A$.

($\La$): Assume every subset of a space $X$ is normal. Let $A,B\subset X$ be separated. We need to show $A,B$ have disjoint neighborhoods. Let $Z=X-(\ol{A}\cap\ol{B})=(\ol{A})^c\cup(\ol{B})^c$. Then $\ol{A}\cap Z$ and $\ol{B}\cap Z$ are closed in $Z$ and so have disjoint sets $U\supset \ol{A}\cap Z\supset A$, $V\supset\ol{B}\cap Z\supset B$ that are open in $Z$, and hence also open in $X$ since $Z$ is open in $X$.
\end{proof}

\section{Spaces of dimension $0$}

\begin{lmm}\label{IFSlmm1}
A space is $0$-dimensional $\iff$ it has a base consisting of clopen sets.
\end{lmm}
\begin{proof}
($\Ra$): Assume $X$ is a $0$-dimensional space. Let $x\in X$, and $U\ni x$ any open set. Then there is an open set $x\in V_x\subset U$ such that $\del V_x=\emptyset$, i.e., $V_x$ is clopen. $\B=\{V_x:x\in X\}$ is the desired base for $X$.

($\La$): Assume $X$ has a base $\B$ consisting of clopen sets. Then for any $x\in X$ and any neighborhood $U\ni x$, we know there exists $V\in \B$  such that $x\in V\subset U$. Hence $X$ is $0$-dimensional.
\end{proof}

\begin{thm}[\textcolor{OliveGreen}{\cite[Theorem II.1]{HureWall41}}]\label{IFSthmII1}
A nonempty subset of a $0$-dimensional space is $0$-dimensional.
\end{thm}
\begin{proof}
Let $X$ be a $0$-dimensional space and $A\subset X$ nonempty. For any $a\in A$, let $U'\subset A$ be open. Then $U'=U\cap A$ for an open set $U\subset X$. Since $X$ is $0$-dimensional, there is a clopen set $V\subset X$ such that $a\in V\subset U$. Let $V':=V\cap A$. Then it is clear that $V'\subset U'$ is open in $A$. Also, $V'$ is closed in $A$ since $A-V'=A\cap(V\cap A)^c=V^c\cap A$. Hence $A$ is $0$-dimensional.
\end{proof}

\begin{dfn}[\textcolor{blue}{\index{Separation of sets}{Separation of sets} by a set}]
Let $X$ be a space and $A,B,C\subset X$ pairwise disjoint sets. We say $A$ and $B$ are \ul{separated in $X$ by $C$} if $X-C=U\sqcup V$ for disjoint relatively open sets $U,V\subset X-C$ (i.e., $U,V$ are disjoint and open in $X-V$) such that $A\subset U$, $B\subset V$. In other words, \ul{$C$ separates $A,B$} if $A,B$ lie in disjoint relatively open sets that partition $X-C$.
\end{dfn}

\begin{lmm}\label{IFSlmm2}
Let $X$ be a space. $A,B\subset X$ are separated by $\emptyset$ $\iff$ there exists a clopen set $V\subset X$ such that $A\subset V$ and $V\cap B=\emptyset$.
\end{lmm}
\begin{proof}
This is clear by the fact that a set is clopen $\iff$ its complement is clopen.
\end{proof}

\begin{lmm}\label{IFSlmm3}
A regular space $X$ is $0$-dimensional $\iff$ every point $x\in X$ and every closed set $C\subset X$ not containing $x$ can be separated by $\emptyset$.
\end{lmm}
\begin{proof}
($\Ra$): Assume $\dim X=0$. Let $x\in X$ and $C\subset X$ a closed set not containing $x$. Then by regularity there is a neighborhood $U\ni x$ such that $\ol{U}$ (hence $U$) excludes $C$. Let $x\in V\subset U$ be a neighborhood such that $\del V=\emptyset$. Then $x$ and $C$ are separated by $\emptyset$ since $x\in V$ and $V\cap C=\emptyset$.

($\La$): Assume every point can be separated by $\emptyset$ from any closed set not containing. Let $x\in X$ and $U\ni x$ an open set. Then $x$ and $U^c$ are separated by $\emptyset$, i.e., there exists a clopen set $V\ni x$ such that $U^c\cap V=\emptyset$. It follows that $x\in V\subset U$.
\end{proof}

\begin{crl}\label{IFSlmm4}
Let $X$ be a $0$-dimensional regular Hausdorff space. (i) $X$ is connected $\iff$ $X$ consists of one point only. (ii) Every nonempty connected subset of $X$ consists of one point only.
\end{crl}

\begin{lmm}\label{IFSlmm5}
A second countable regular Hausdorff space $X$ is $0$-dimensional $\iff$ any two disjoint closed sets $C_1,C_2\subset X$ can be separated by $\emptyset$.
\end{lmm}
\begin{proof}
($\Ra$): Assume $\dim X=0$. Let $C_1,C_2\subset X$ be closed sets. Since $\emptyset$ can separate any point from any closed set not containing it, for each $x\in X$, there exists a clopen neighborhood $U(x)$ such that \ul{either} $U(x)\cap C_1=\emptyset$ \ul{or} $U(x)\cap C_2=\emptyset$. Since $X$ has a countable base $\B=\{B_1,B_2,\cdots\}$, with $x_i\in B_i\subset U_i:=U(x_i)$, we get a countable cover $U_1,U_2,\cdots$ of $X$. Define a new sequence of sets by
\bea
\textstyle V_1:=U_1,~~~~V_i:=U_i-\bigcup_{k=1}^{i-1}U_k=U_i\cap(X-\bigcup_{k=1}^{i-1}U_k),~~~~i=2,3,\cdots\nn
\eea
Then it is clear that (1) $\{V_i\}$ covers $X$, (2) $V_i\cap V_j=\emptyset$ for $i\neq j$, (3) each $V_i$ is open, and (4) for each $i$, \ul{either} $V_i\cap C_1=\emptyset$ \ul{or} $V_i\cap C_2=\emptyset$. Let $O_1=\bigcup\{V_i:V_i\cap C_2=\emptyset\}$, $O_2=\bigcup\{V_i:V_i\cap C_1=\emptyset\}$. Then it is clear that $X=O_1\cup O_2$, $O_1\cap O_2=\emptyset$, $C_1\subset O_1$, $C_2\subset O_2$.

($\La$): If $x\in X$ and $U\ni x$ is open, then $x$ and $U^c$ (being closed sets) can be separated by $\emptyset$ and so there is a clopen set $V\ni x$ such that $V\cap U^c=\emptyset$, i.e., $x\in V\subset U$.
\end{proof}

\begin{lmm}\label{IFSlmm5b}
Let $X$ be a second countable completely normal Hausdorff space and $A\subset X$ a $0$-dimensional subset (of the same type, i.e., second countable normal Hausdorff). For any two disjoint closed sets $C_1,C_2\subset X$, there exists a closed set $B$ separating $C_1,C_2$ such that $A\cap B=\emptyset$.
\end{lmm}
\begin{proof}
By normality, $C_1,C_2$ we have neighborhoods (i) $U_1\supset C_1$, $U_2\supset C_2$ such that $\ol{U}_1\cap\ol{U}_2=\emptyset$. Thus, $\ol{U}_1\cap A$ and $\ol{U}_2\cap A$ are two disjoint closed sets in $A$, and so can be separated by $\emptyset$ in $A$ (Lemma \ref{IFSlmm5}). Thus, we have disjoint sets $C_1',C_2'\subset A$, that are clopen in $A$, such that
\bea
(ii)~~A=C_1'\cup C_2',~~~~\ol{U}_1\cap A\subset C_1',~~~~\ol{U}_2\cap A\subset C_2'.\nn
\eea
It follows that
\bea
&&(iii)~~C_1'\cap\ol{U}_2~\cup~C_2'\cap\ol{U}_1=\emptyset,~~~~\ol{C_1'}\cap U_2~\cup~\ol{C_2'}\cap U_1=\emptyset,\nn\\
&&(iv)~~C_1'\cap\ol{C_2}~\cup~C_2'\cap\ol{C_1}=\emptyset,~~~~\ol{C_1'}\cap C_2~\cup~\ol{C_2'}\cap C_1=\emptyset,~~~~C_1'\cap\ol{C_2'}~\cup~\ol{C_1'}\cap C_2'=\emptyset.\nn
\eea
It follows from (iv), and $\ol{C}_1\cap\ol{C}_2=C_1\cap C_2=\emptyset$, that $C_1\cup C_1'$ and $C_2\cup C_2'$ are separated. Thus, by complete normality of $X$, there exists an open set $W$ such that
\bea
C_1\cup C_1'\subset W~~~~\txt{and}~~~~\ol{W}\cap(C_2\cup C_2')=\emptyset.\nn
\eea
The boundary $B:=\del W$ separates $C_1,C_2$ and is disjoint from $A=C_1'\cup C_2'$.
\end{proof}

\begin{thm}[\textcolor{OliveGreen}{Union theorem for dimension $0$: \cite[Theorem II.2]{HureWall41}}]\label{IFSthmII2}
A space that is a \ul{countable} union of $0$-dimensional \ul{closed} (or $F_\sigma$) sets is $0$-dimensional. (Assume all sets, including their unions, are normal Hausdorff spaces.)
\end{thm}
\begin{proof}
Let $X=C_1\cup C_2\cup\cdots$, where each $C_i$ is closed and $0$-dimensional. Let $K,L\subset X$ be disjoint closed sets. We need to show $K,L$ are separated by $\emptyset$. Since $K\cap C_1,L\cap C_1\subset C_1$ are disjoint closed sets and $C_1$ is $0$-dimensional, there are clopen subsets $A_1,A_2$ of $C_1$ such that
\bea
C_1=A_1\cup B_1,~~A_1\cap B_1=\emptyset,~~~~K\cap C_1\subset A_1,~~~~L\cap C_1\subset B_1.\nn
\eea
The sets $K\cup A_1$, $L\cup B_1$ are closed and disjoint in $X$, and so by normality of $X$ there are open sets $G_1,H_1$ such that $K\cup A_1\subset G_1$, $L\cup B_1\subset H_1$, $\ol{G}_1\cap\ol{H}_1=\emptyset$. Thus,
\bea
C_1\subset G_1\cup H_1,~~~~K\subset G_1,~~L\subset H_1,~~~~\ol{G}_1\cap\ol{H}_1=\emptyset.\nn
\eea
If we now repeat the above process with $K,L,C_1$ replaced by $\ol{G}_1,\ol{H}_1,C_2$ respectively, we get open sets $G_2,H_2$ such that
\bea
C_2\subset G_2\cup H_2,~~~~\ol{G}_1\subset G_2,~~\ol{H}_1\subset H_2,~~~~\ol{G}_2\cap\ol{H}_2=\emptyset.\nn
\eea
Continuing this way, we get sequences of sets $\{G_i\}$, $\{H_i\}$ such that
\bea
C_i\subset G_i\cup H_i,~~~~\ol{G}_{i-1}\subset G_i,~~\ol{H}_{i-1}\subset H_i,~~~~\ol{G}_i\cap\ol{H}_i=\emptyset.\nn
\eea
It follows that $G:=\bigcup G_i$, $H:=\bigcup H_i$ are disjoint open sets satisfying $K\subset G$, $L\subset H$, and $X=G\cup H$ (since $X=\bigcup C_i\subset G\cup H$).
\end{proof}

\begin{crl}\label{IFScrl1}
A space that is a countable union of $0$-dimensional $F_\sigma$ sets is $0$-dimensional.
\end{crl}

\begin{crl}\label{IFScrl2}
The union of two $0$-dimensional subsets of a space $X$, at least one of which is closed, is $0$-dimensional. (Assume the spaces are separable metric spaces.)
\end{crl}
\begin{proof}
Let $A,B\subset X$ be $0$-dimensional, and assume wlog that $B$ is closed. Then $A=(A\cup B)-B$ is open in $A\cup B$. Thus, $A$ is $F_\sigma$ (as an open set in a metric space -- footnote\footnote{Let $X$ be a metric space and $O\subset X$ an open set. Then $O=\bigcup_n\{x\in X:\dist(x,X-O)\geq 2^{-n}\}$ is $F_\sigma$, since each $C_n:=\{x\in X:\dist(x,X-O)\geq 2^{-n}\}$ is closed as the preimage of a closed set under a continuous function.}). It follows that $A\cup B$ is a countable union of $0$-dimensional $F_\sigma$ sets.
\end{proof}

\begin{crl}\label{IFScrl3}
$0$-dimensional space remains $0$-dimensional after the adjunction of a single point. (Assume the spaces are separable metric spaces.)
\end{crl}

\section{Spaces of finite dimension $n$}
\begin{lmm}\label{IFSlmm6}
A space has dimension $\leq n$ $\iff$ it has a base consisting of sets whose boundaries have dimension $\leq n-1$.
\end{lmm}
\begin{proof}
($\Ra$): Assume $X$ is a space of dimension $\leq n$. Let $x\in X$, and $U\ni x$ any open set. Then there is an open set $x\in V_x\subset U$ such that $\dim\del V_x\leq n-1$. $\B=\{V_x:x\in X\}$ is the desired base for $X$.

($\La$): Assume $X$ has a base $\B$ consisting of sets with boundaries of dimension $\leq n-1$. Then for any $x\in X$ and any neighborhood $U\ni x$, we know there exists $V\in \B$  such that $x\in V\subset U$. Hence $X$ has dimension $\leq n$.
\end{proof}

\begin{thm}[\textcolor{OliveGreen}{\cite[Theorem III.1]{HureWall41}}]\label{IFSthmIII1}
If a space has dimension $\leq n$, then every subspace has dimension $\leq n$.
\end{thm}
\begin{proof}
We proceed by induction on $n$. The case of $n=-1$ is clear. Suppose the case of $n-1$ is true as well. Let $X$ be a space with $\dim X\leq n$. Let $A\subset X$, $a\in A$, and let $a\in U\subset A$ be a neighborhood of $a$ in $A$. Then $U=U_1\cap A$, where $U_1\subset X$ is open. There is an open set $a\in V_1\subset U_1$ such that $\dim\del V_1\leq n-1$. The set $V:=V_1\cap A$ is open in $A$, $a\in V\subset U$, and
\bea
\del_A V=\del_A(V_1\cap A)\subset(\del V_1)\cap A\subset\del V_1~~\sr{\txt{induction hypothesis}}{\Longrightarrow}~~\dim\del_A V\leq n-1.\nn\qedhere
\eea
\end{proof}

\begin{lmm}\label{IFSlmm7}
A regular space $X$ has dimension $\leq n$ $\iff$ every point $x\in X$ and any closed set $C\subset X$ not containing $x$ can be separated by a closed set of dimension $\leq n-1$.
\end{lmm}
\begin{proof}
($\Ra$): Assume $\dim X\leq n$. Since $X-C$ is a neighborhood of $x$, by regularity, there is a neighborhood $V\ni x$ such that $\ol{V}\subset X-C$. Take a neighborhood $x\in W\subset V$ such that $\del W\leq n-1$. Note that $\del W$ separates $x$ and $C$, since $X-\del W=W\cup(X-\ol{W})$ and $x\in W$, $C\subset X-\ol{W}$.

($\La$): Assume any point and any closed set not containing it can be separated by a closed set of dimension $\leq n-1$. Let $U\ni x$ be open. Then the closed set $X-U\not\ni x$ can be separated from $x$ by a closed set $B$ of dimension $\leq n-1$. That is, we have open sets $U',V'\subset X-B$ such that
\bea
X-B=U'\cup V',~~~~p\in U',~~~~X-U\subset V',~~~~U'\cap V'=\emptyset.\nn
\eea
Hence, $U'\subset U$ and $\del U'\subset B$ implies $\dim\del U'\leq n-1$.
\end{proof}

\begin{lmm}\label{IFSlmm8}
Let $X$ be a completely normal space. A subspace $A\subset X$ has dimension $\leq n$ $\iff$ for any $a\in A$ and any neighborhood $a\in U\subset X$, there exists a neighborhood $a\in V\subset U$ such that $\dim\big(\del V\cap A\big)\leq n-1$.
\end{lmm}
\begin{proof}
($\Ra$): Assume $\dim A\leq n$. Let $a\in A$ and $a\in U\subset X$ be open. Then $U'=U\cap A$ is open in $A$, and so there exists an open set $V'\subset A$ such that $a\in V'\subset U'$ and $\dim\del_A V'\leq n-1$.
Since $X$ is completely normal and $\ol{V'}$ and $A-\ol{V'}$ are separated, it follows there exists an open set $V\subset X$ such that $V'\subset V$ and $\ol{V}\cap(A-\ol{V'})=\emptyset$. Since $V'\subset U$, by replacing $V$ with $V\cap U$, we may assume wlog that $V\subset U$. By construction, we have $\del V\cap V'=\emptyset=\del V\cap(A-\ol{V'})$, and so
\begin{align}
\del V\cap A=\del V\cap(V'\cup\del_AV'\cup[A-\ol{V'}])=\del V\cap\del_AV'\subset\del_AV',~~\Ra~~\dim(\del V\cap A)\leq n-1.\nn
\end{align}

($\La$): Assume for any $a\in A$ and any neighborhood $a\in U\subset X$, there exists a neighborhood $a\in V\subset U$ such that $\dim\big(\del V\cap A\big)\leq n-1$. Let $x\in A$ and $U'\subset A$ a neighborhood of $x$ in $A$. Then $U'=U\cap A$ for an open set $U\subset X$. Thus, there exists an open set $x\in V\subset U$ such that $\dim\big(\del V\cap A\big)\leq n-1$. Let $V'=V\cap A$. Then $V'$ is open in $A$, $x\in V'\subset U'$, and
\bea
\del_AV'=\del_A(V\cap A)\subset\del V\cap A,~~\Ra~~\dim\del_AV'\leq n-1.\nn\qedhere
\eea
\end{proof}

\begin{lmm}\label{IFSlmm9}
If $X$ is a completely normal space, then any two subspaces $A,B\subset X$ satisfy
\bea
\dim(A\cup B)\leq1+\dim A+\dim B.\nn
\eea
\end{lmm}
\begin{proof}
We proceed by induction on $\dim A+\dim B$. The statement holds for $\dim A=\dim B=-1$. Let $m:=\dim A$, $n:=\dim B$. Assume the statement holds for the following cases. (i) $\dim A_1\leq m$, $\dim B_1\leq n-1$. (ii) $\dim A_1\leq m-1$, $\dim B_1\leq n$. (I.e., each of (i),(ii) implies $\dim(A_1\cup B_1)\leq 1+\dim A_1+\dim B_1$.)

We need to show the statement also holds for $\dim A_1\leq m$, $\dim B_1\leq n$. Let $x\in A\cup B$, and assume wlog that $x\in A$. Let $U$ be a neighborhood of $x$ in $X$. Since $\dim A\leq m$, it follows from Lemma \ref{IFSlmm8} that there exists an open set $V$ such that $x\in V\subset U$ and $\dim(\del V\cap A)\leq m-1$. Moreover, since $\del V\cap B\subset B$ and $\dim B\leq n$, we also have $\dim(\del V\cap B)\leq n$. Thus, by the induction hypotheses (i),(ii),
\bea
\dim(\del V\cap(A\cup B))\leq1+\dim(\del V\cap A)+\dim(\del V\cap B)\leq 1+(m-1)+n=m+n.\nn
\eea
It follows again by Lemma \ref{IFSlmm8} that $\dim(A\cup B)\leq m+n+1=1+\dim A+\dim B$.
\end{proof}

\begin{crl}\label{IFScrl4}
Let $X_0,X_1,\cdots,X_n\subset X$ be subspaces such that $\dim X_i\leq 0$ for all $i=0,...,n$. Then $\dim(X_0\cup X_1\cup\cdots\cup X_n)\leq n$.
\end{crl}

\begin{lmm}\label{IFSlmm10}
Fix an integer $n\geq -1$. Suppose every space that is a \ul{countable} union of \ul{closed} (or $F_\sigma$) sets of dimension $\leq n-1$ has dimension $\leq n-1$. (Assume all sets, including their unions, are separable normal Hausdorff spaces). Then any space of dimension $\leq n$ is the union of a subspace of dimension
$\leq n - 1$ and a subspace of dimension $\leq 0$.
\end{lmm}
\begin{proof}
Let $X$ be a space of dimension $\leq n$. Then we know $X$ has a base of sets whose boundaries have dimension $\leq n-1$. Since $X$ is separable, let $\{U_1,U_2,\cdots\}$ be a countable base of sets whose boundaries $\{B_1,B_2,\cdots\}$ have dimension $\leq n-1$ (where $B_i=\del U_i$). By hypothesis, the set {\small $B:=\bigcup_{i=1}^\infty B_i$} has dimension $\leq n$. It remains to show that {\small $\dim(X-B)\leq 0$}.

By Lemma \ref{IFSlmm8}, since $\dim(\del U_i\cap(X-B))=\dim\emptyset\leq -1$ for all $i$ (and $\{U_i\}$ is a base), it follows that $\dim(X-B)\leq -1+1=0$.
\end{proof}

\begin{thm}[\textcolor{OliveGreen}{Union theorem for dimension $n$: \cite[Theorem III.2]{HureWall41}}]\label{IFSthmIII2}
A space that is a \ul{countable} union of \ul{closed} (or $F_\sigma$) sets of dimension $\leq n$ has dimension $\leq n$. (Assume all sets, including their unions, are separable metric spaces.)
\end{thm}
\begin{proof}
We proceed by induction on $n$. Denote the statement of the theorem by $S_n$. $S_{-1}$ holds trivially. Also $S_0$ holds by Theorem \ref{IFSthmII2}. Assume $S_{n-1}$ holds. Then by Lemma \ref{IFSlmm10}, $X$ is the union of a subspace of dimension $\leq n-1$ and a subspace of dimension $\leq 0$.

Let $X=C_1\cup C_2\cup\cdots$, where each $i$ is closed and $\dim C_i\leq n$. We need to show $\dim X\leq n$. Let
\bea
\textstyle K:=C_1,~~~~K_i:=C_i-\bigcup_{j=1}^{i-1}C_j=C_i\cap\left(X-\bigcup_{j=1}^{i-1}C_j\right),~~~~i=2,3,\cdots.\nn
\eea
Then it is clear that $X=\bigcup K_i$ and $K_i\cap K_j=\emptyset$ if $i\neq j$. Also, $K_i$ is an $F_\sigma$ (being the intersection of a closed set and an open set in a metric space, which is thus $F_\sigma$), and $\dim K_i\leq n$ (since $K_i\subset C_i$). By Lemma \ref{IFSlmm10}, each $K_i$ can be expressed as follows:
\bea
K_i=M_i\cup N_i,~~~~\dim M_i\leq n-1,~~~~\dim N_i\leq 0.\nn
\eea
Let $M:=\bigcup M_i$, $N=\bigcup N$. Then $X=M\cup N$. Note that each $M_i$ is an $F_\sigma$ \ul{in $M$}, since $M_i\subset K_i$ and $K_i\cap K_j=\emptyset$ for $i\neq j$, and so
\bea
M_i=M_i\cap K_i=M\cap K_i~~~~\txt{(where we know $K_i$ is an $F_\sigma$ in $X$)}.\nn
\eea
That is $M$ is a union of $F_\sigma$ sets of dimension $\leq n-1$, and so {\small $\dim M\leq n-1$} (by induction hypotheses $S_{n-1}$). Similarly, each $N_i$ is an $F_\sigma$ of dimension $\leq 0$ in $N$, and so {\small $\dim N\leq 0$} (by $S_0$).

Hence, by Lemma \ref{IFSlmm9}, {\small $\dim(X)=\dim(M\cup N)\leq 1+\dim M+\dim N\leq 1+(n-1)+0=n$}.
\end{proof}

\begin{crl}\label{IFScrl5}
The union of two subspaces each of dimension $\leq n$, at least one of which is closed, has dimension $\leq n$. (Assume the spaces are separable metric spaces.)
\end{crl}
\begin{proof}
Let $A,B\subset X$ each have dimension $\leq n$, and assume wlog that $B$ is closed. Then $A=(A\cup B)-B$ is open in $A\cup B$. Thus, $A$ is $F_\sigma$ (as an open set in a metric space). It follows that $A\cup B$ is a countable union of $F_\sigma$ sets of dimension $\leq n$.
\end{proof}

\begin{crl}\label{IFScrl6}
A space of dimension $\leq n$ remains a space of dimension $\leq n$ after the adjunction of a single point. (Assume the spaces are separable metric spaces.)
\end{crl}

\begin{crl}[\textcolor{OliveGreen}{Consequence of Lemma \ref{IFSlmm8}}]\label{IFScrl7}
Let $X$ be a completely normal space. If a subspace $A\subset X$ has dimension $\leq n$, then for any $x\in X$ and any neighborhood $x\in U\subset X$, there exists a neighborhood $x\in V\subset U$ such that $\dim\big(\del V\cap A\big)\leq n-1$.
\end{crl}
\begin{proof}
For each point $x\in X$, $A\cup\{x\}$ has dimension $\leq n$. Hence the conclusion follows by Lemma \ref{IFSlmm8}.
\end{proof}

\begin{crl}\label{IFScrl8}
A space of dimension $\leq n$ is the union of a subspace of dimension $\leq n - 1$ and a subspace of dimension $\leq 0$. (Proof: This follows from Lemma \ref{IFSlmm10}.)
\end{crl}

\begin{thm}[\textcolor{OliveGreen}{Decomposition theorem for dimension $n$: \cite[Theorem III.3]{HureWall41}}]\label{IFSthmIII3}
A space has dimension $\leq n$ $\iff$ it is the union of $n + 1$ subspaces each of dimension $\leq 0$.
\end{thm}
\begin{proof}
This follows by repeated application of Corollary \ref{IFScrl8}, and then application of Corollary \ref{IFScrl4}.
\end{proof}

\begin{crl}\label{IFScrl9}
Let $X$ be a space with $\dim X=n$. Then for any numbers $p,q\geq -1$ satisfying $n=p+q+1$, we have $X=A\cup B$ for subspaces $A,B\subset X$ of dimensions $\dim A=p$, $\dim B=q$.
\end{crl}

\begin{thm}[\textcolor{OliveGreen}{Dimension of product space: \cite[Theorem III.4]{HureWall41}}]\label{IFSthmIII4}
If $X,Y$ are spaces, at least one of which is nonempty, then $\dim(X\times Y)\leq \dim X+\dim Y$.
\end{thm}
\begin{proof}
We proceed by induction on $\dim X+\dim Y$. The result is clear if either $\dim X=-1$ or $\dim Y=-1$. Let $m:=\dim X$, $n:=\dim Y$. Assume the result holds for the cases (i) $\dim X_1\leq m$, $\dim Y_1\leq n-1$ and (ii) $\dim X_1\leq m-1$, $\dim Y_1\leq n$. (I.e., these conditions imply $\dim(X_1\times Y_1)\leq\dim X_1+\dim Y_1$.)

Let $(x,y)\in X\times Y$. Then any open neighborhood of $(x,y)$ contains an open set of the form $U\times V\subset X\times Y$, where $U\subset X$ is open in $X$, $V\subset Y$ is open in $Y$, and by the definition of dimension we can assume $\dim\del U\leq m-1$, $\dim\del V\leq n-1$. Note that $\del(U\times V)=(\del U\times\ol{V})\cup(\ol{U}\times\del V)$.

By the hypotheses (i),(ii), each term on the right has dimension $\leq m+n-1$, and so $\dim\del(U\times V)\leq m+n-1$ by Theorem \ref{IFSthmIII2} (since the terms on the right are closed subspaces). Hence, $\dim(X\times Y)\leq m+n$.
\end{proof}

\section{Separation of sets in $n$-dimensional spaces}
\begin{lmm}\label{IFSlmm10a}
Let $X$ be a space, $C_1,C_2\subset X$ disjoint closed sets, and $A\subset X$ a set of dimension $\leq n$. Then there exists a closed set $C$ separating $C_1,C_2$, such that $\dim(A\cap C)\leq n-1$.
\end{lmm}
\begin{proof}
If $n=-1$, the statement holds trivially. Also, the statement holds for $n=0$ by Lemma \ref{IFSlmm5b}. So, assume $n\geq 1$. By Corollary \ref{IFScrl8}, $A=D\cup E$, where $\dim D\leq n-1$, $\dim E\leq 0$. By the $n=0$  case (i.e., Lemma \ref{IFSlmm5b}) we can separate $C_1,C_2$ by a closed set $C\subset X$ such that $E\cap C=\emptyset$ (which implies $A\cap C\subset D$).

Hence, $\dim(A\cap C)\leq n-1$, since $\dim D\leq n-1$.
\end{proof}

\begin{crl}\label{IFScrl9a}
Let $X$ be a $T_1$ space of dimension $\leq n$. Then any two disjoint closed sets $C_1,C_2\subset X$ can be separated by a closed set $C\subset X$ of dimension $\leq n-1$.
\end{crl}
Note that the converse of this result also holds: Indeed if any two disjoint closed sets in $X$ can be separated by a closed set of dimension $\leq n-1$, then for any $x\in X$ and any open set $U\ni x$, the disjoint closed sets $x,U^c$ can be separated by a closed set $C$ of dimension $\leq n-1$, i.e., $X-C=V\cup W$, where $x\in V,U^c\subset W$ and $V,W$ are disjoint and open in $X-C$ (hence also open in $X$). It follows that $V\subset W^c\subset U$ and $\del V\subset C$ implies $\dim\del V\leq n$.

\begin{lmm}\label{IFSlmm10b}
Let $X$ be a space of dimension $\leq n-1$. Let $C_i,C_i'\subset X$, $i=1,...,n$, be $n$ pairs of closed sets such that $C_i\cap C_i'=\emptyset$ for each $i=1,...,n$. Then there exist $n$ closed sets $B_1,...,B_n$ such that $B_i$ separates $C_i,C_i'$, and $B_1\cap \cdots\cap B_n=\emptyset$.
\end{lmm}
\begin{proof}
By Corollary \ref{IFScrl9a}, we have a closed set $B_1$ separating $C_1,C_1'$ such that
\bea
\dim B_1\leq (n-1)-1=n-2.\nn
\eea
By Lemma \ref{IFSlmm10a}, we then get a closed set $B_2$ separating $C_2,C_2'$ such that
\bea
\dim B_1\cap B_2\leq (n-2)-1=n-3.\nn
\eea
By repeated application of Lemma \ref{IFSlmm10a} as above, we get $n$ closed sets $B_1,...,B_n$ such that $B_i$ separates $C_1,C_i'$ and $\dim B_1\cap\cdots\cap B_k\leq n-k-1$ for $k=1,...,n$.

Hence, for $k=n$, we see that $B_1\cap\cdots\cap B_n=\emptyset$.
\end{proof}

\section{Dimension of Euclidean spaces}
\begin{lmm}\label{IFSlmm10c}
Let $I:=[-1,1]\subset\Real$. Let $C_i$ be the $i$th face of $I^n\subset\Real^n$ given by $C_i:=\{x\in I^n:x_i=1\}$, and $C_i'$ the face opposite to $C_i$ given by $C_i':=\{x\in I^n:x_i=-1\}$. If $K_i$ is a \ul{closed} set separating $C_i,C_i'$, then ~$K_1\cap\cdots\cap K_n\neq\emptyset$.
\end{lmm}
\begin{proof}
Since $K_i$ separates $C_i,C_i'$ in $I^n$, we have {\footnotesize $I^n-K_i=U_i\cup U_i'$, $C_i\subset U_i$, $C_i'\subset U_i'$, $U_i\cap U_i'=\emptyset$}, where $U_i,U_i'$ are open in $I^n-K_n$, and hence open in $I^n$. Define the cont. map
\bea
v:I^n\ra I^n,~~x\mapsto v(x),~~
v(x)_i:=\left\{
          \begin{array}{ll}
            -\dist(x,K_i), & x\in U_i \\
            0, & x\in K_i \\
            \dist(x,K_i), & x\in U_i'
          \end{array}
        \right\}.\nn
\eea
Define another continuous map $f:I^n\ra I^n$ by $f(x):=x+v(x)$, where
\bea
f(x)_i=\left\{
         \begin{array}{ll}
           x_i-\dist(x,K_i), & x\in U_i \\
           x_i, & x\in K_i \\
           x_i+\dist(x,K_i), & x\in U_i'
         \end{array}
       \right\}\in I=[-1,1].\nn
\eea
Since $I^n$ is homeomorphic to the closed $n$-ball, it follows from the Brouwer fixed point theorem that there exists $x^0\in I^n$ such that $f(x^0)=x^0$, i.e., $v(x^0)=0$. This implies $\dist(x^0,K_i)=0$, i.e., $x^0\in K_i$, for each $i$.
\end{proof}

\begin{lmm}\label{IFSlmm10d}
If $I^n=[-1,1]^n$, then $\dim I^n\geq n$.
\end{lmm}
\begin{proof}
Suppose that $\dim I^n\leq n-1$. Then there exist $n$ closed sets $B_i\subset I^n$, with $B_i$ separating a pair of disjoint faces $C_i,C_i'$, and $B_1\cap\cdots\cap B_n=\emptyset$ (a contradiction of Lemma \ref{IFSlmm10c}).
\end{proof}

\begin{thm}[\textcolor{OliveGreen}{\cite[Theorem IV1, p.41]{HureWall41}}]\label{IFSthmIV1}
The Euclidean space $\Real^n$ has dimension $=n$.
\end{thm}
\begin{proof}
\ul{$\dim\Real^n\leq n$}: From the definition of dimension, $\dim\Real=1$. Thus, it follows from Theorem \ref{IFSthmIII4} and induction that $\dim\Real^n\leq n$. \ul{$\dim\Real^n\geq n$}: This follows from Lemma \ref{IFSlmm10d} and the fact that $I^n\subset\Real^n$.
\end{proof}

\begin{crl}
If $I:=[-1,1]$, then $\dim I^n=n$. (This follows from Lemma \ref{IFSlmm10d} and the fact $I^n\subset\Real^n$.)
\end{crl}

\section{Imbedding theorems}
In this section, where necessary, assume the spaces are second countable completely normal Hausdorff (ScCnH-) spaces.
\begin{dfn}[\textcolor{blue}{\index{Diameter of a cover}{Diameter (or mesh) of a cover}, \index{Order of a cover}{Order of a cover}}]
Let $X$ be a space and $\U$ an open cover of $X$. (i) If $X$ is a metric space, the \ul{diameter} of $\U$ is ~$\diam\U:=\sup\{\diam U:U\in \U\}$. (ii) The \ul{order} of $\U$ (denoted $\Ord\U$) $:=$ (the largest number of elements of $\U$ with a common point)$-1$, i.e., $\Ord\U$ $:=$ the largest number $n$ such that there exist $n+1$ elements of $\U$ with nonempty intersection. Equivalently, if $\Ord_x\U$ (the order of $\U$ at $x$) is the number of elements of $\U$ containing $x\in X$, then ~$\Ord\U:=\sup_{x\in X}\Ord_x\U$.
\end{dfn}

\begin{lmm}\label{IFSlmm11}
Let $X$ be a ScCnH-space and $A\subset X$ a subspace of dimension $\leq 0$. If $U_1,U_2\subset X$ are open sets that cover $A$ (i.e., $A\subset U_1\cup U_2$), then there exist \ul{open sets} $V_1,V_2\subset X$ that \ul{cover} $A$ and also \ul{satisfy}
\bea
V_1\cap V_2=\emptyset,~~~~V_1\subset U_1,~~~~V_2\subset U_2.\nn
\eea
\end{lmm}
\begin{proof}
The case of $\dim A=-1$ (i.e., $A=\emptyset$) is clear. So, assume $\dim A=0$.
By replacing $X$ with $U_1\cup U_2$, we can assume $X=U_1\cup U_2$. Then, $C_1:=X-U_2$, $C_2:=X-U_1$ are disjoint closed sets. Since $\dim A=0$, by Lemma \ref{IFSlmm5b}, there exists a closed set $B$ such that $A\cap B=\emptyset$ and $B$ separates $C_1,C_2$, i.e., there exist open sets $V_1,V_2$ that (are the desired sets because they) satisfy
\bea
V_1\cap V_2=\emptyset,~~X-B=V_1\cup V_2,~~C_1\subset V_1,~~C_2\subset V_2.\nn\qedhere
\eea
\end{proof}

\begin{lmm}\label{IFSlmm12}
Let $X$ be a space and $A\subset X$ a subspace of dimension $\leq 0$. If $U_1,\cdots,U_r\subset X$ are open sets that cover $A$ (i.e., $A\subset \bigcup U_i$), then there exist open sets $V_1,\cdots,V_r\subset X$ that cover $A$ and also satisfy
\bea
V_i\cap V_j=\emptyset~~\txt{if}~~i\neq j,~~~~\txt{and}~~~~V_i\subset U_i~~\txt{for all}~~i=1,2,\cdots,r.\nn
\eea
\end{lmm}
\begin{proof}
We proceed by induction on $r$. The result holds for $r=1$ trivially and for $r=2$ by Lemma \ref{IFSlmm11}. Assume the statement holds for any $r-1$ open sets (induction hypothesis). Let $U_{r-1}':=U_{r-1}\cup U_r$.
Then we have the open covering of $A$ consisting of the open sets
\bea
U_1,...,U_{r-2},U_{r-1}'.\nn
\eea
By the induction hypothesis, there exist pairwise disjoint open sets $V_1,...,V_{r-2},V_{r-1}'$ that cover $A$, and
\bea
V_1\subset U_1,~~\cdots,~~V_{r-2}\subset U_{r-2},~~V_{r-1}'\subset U_{r-1}'.\nn
\eea
Now, $V_{r-1}'\cap A$ has dimension $\leq 0$ and is covered by $U_{r-1}$ and $U_r$, and so by Lemma \ref{IFSlmm11}, there exist pairwise disjoint open sets $V_{r-1},V_r$ that cover $V_{r-1}'\cap A$, and $V_{r-1}\subset U_{r-1}$, $V_r\subset U_r$.

Hence, the desired collection of sets is $V_1,\cdots,V_r$.
\end{proof}

\begin{thm}[\textcolor{OliveGreen}{\cite[Theorem V1, p.54]{HureWall41}}]\label{IFSthmV1}
Let $X$ be a space of dimension $\leq n$. Then every finite open covering $\U$ of $X$ has an open refinement $\V$ of order $\leq n$.
\end{thm}
\begin{proof}
By the decomposition theorem (Theorem \ref{IFSthmIII3}), $A=A_1\cup\cdots\cup A_{n+1}$, where $\dim A_i\leq 0$. Let $A_i$ be covered by $r(i)$ members $\U^i=\{U_1^i,\cdots,U_{r(i)}^i\}$ of $\U$. Then it follows from Lemma \ref{IFSlmm12} that there exists (for each $i$) a finite open cover $\V^i$ of $A_i$, with
\bea
\V^i=\{V_1^i,....,V_{r(i)}^i\},~~~~i=1,....,n+1,~~~~V_j^i\subset U_j^i,\nn
\eea
where the elements of each $\V^i$ are pairwise disjoint, i.e., $V_j^i\cap V_k^i=\emptyset$ if $j\neq k$. Let $\V$ be the cover of $X$ consisting of all these sets, i.e., $\V:=\V^1\cup\V^2\cup\cdots\cup\V^{n+1}=\{V^i_j:j=1,...,r(i),~i=1,...,n+1\}$.

Now, any elements $B_1,...,B_t\in\V$ with a common intersection point must come from different $\V^i$'s, and so $t\leq n+1$. That is, $\V$ has order $\leq n$.
\end{proof}

\begin{crl}\label{IFScrl10}
Let $X$ be a compact metric space of dimension $\leq n$. Then for any $\vep>0$, there exists an open cover of $X$ of diameter (or mesh) $<\vep$ and order $\leq n$.
\end{crl}
\begin{proof}
Since $X$ is compact, a finite number of balls $\B=\{B_1,...,B_r\}$ of radius $\vep/2$ covers $X$, and it follows from Theorem \ref{IFSthmV1} that $\B$ has an open refinement of order $\leq n$.
\end{proof}

\begin{dfn}[\textcolor{blue}{$\vep$-map, Set of $\vep$-maps}]
Let $X$ be a compact metric space, $Y$ a metric space, and $\vep>0$. Then a continuous map $g:X\ra Y$ is an $\vep$-map of $X$ to $Y$ if the primage $g^{-1}(y)$ of each point $y\in g(X)\subset Y$ has diameter $<\vep$. We will denote the set of all $\vep$-maps from $X$ to $Y$ by $G_\vep(X,Y)$.
\end{dfn}

\begin{lmm}\label{IFSlmm13}
Let $X$ be a compact metric space and $Y$ a metric space. Then a continuous map $g:X\ra Y$ is a $1/i$-map for each positive integer $i$ $\iff$ $g$ is an imbedding (i.e., a homeomorphism onto its image).
\end{lmm}
\begin{proof}
Assume $g$ is a $1/i$-map for each $i$. To show $g:X\ra g(X)$ is a homeomorphism, it is enough to show $g$ is injective since a bijective continuous map from a compact space to a Hausdorff space is a homeomorphism. Let $x\in X$. Then $\diam g^{-1}\big(g(x)\big)<1/i$ for all $i$, and so $\diam g^{-1}\big(g(x)\big)=0$, i.e., $g$ is injective.

Conversely, if {\footnotesize $X\sr{g}{\ra} g(X)$} is a homeomorphism, it is clear that $g$ is a $1/i$ map for each $i$.
\end{proof}

\begin{lmm}\label{IFSlmm14}
Let $X$ be a compact metric space, $Y$ a metric space, and $\vep>0$. Then the set $G_\vep(X,Y)$ is open in $C(X,Y)$.
\end{lmm}
\begin{proof}
Let $g\in G_\vep(X,Y)$. We need to find $\eta>0$ such that for all $f\in C(X,Y)$, $d_u(f,g)<\eta$ implies $f\in G_\vep$. Let $\delta:=\inf\left\{d\big(g(x),g(x')\big):d(x,x')\geq\vep\right\}$. Since $X$ (hence $X\times X$) is compact and $(x,x')\mapsto d\big(g(x),g(x')\big)$ is continuous, it follows that $\delta=d\big(g(x),g(x')\big)$ for some $x,x'$ such that $d(x,x')\geq\vep$. This implies $\delta>0$, otherwise, if $\delta=0$ or $g(x)=g(x')$, then $\{x,x'\}\subset g^{-1}(g(x))=g^{-1}(g(x'))$ implies $d(x,x')<\vep$ since $g$ is an $\vep$-map (which contradicts $d(x,x')\geq\vep$).

Let $f\in C(X,Y)$ satisfy $d_u(f,g)<\delta/2$. Let $x,x'\in X$ be such that $f(x)=f(x')$, i.e., $x,x'\in f^{-1}(y)$ for some $y\in f(X)$. Then $d(x,x')<\vep$ (so that $f$ is an $\vep$-map): This is because $d(x,x')\geq\vep$ implies $d\big(g(x),g(x')\big)\geq\vep$ (by the definition of $\delta$), meanwhile
{\small\begin{align}
d\big(g(x),g(x')\big)\leq d\big(g(x),f(x)\big)+d\big(f(x),f(x')\big)+d\big(f(x'),g(x')\big)\leq d_u(f,g)+0+d_u(f,g)<\delta.\nn
\end{align}}
\end{proof}

\begin{prp}\label{IFSlmm15}
Let $X$ be a compact metric space of dimension $\leq n$. For any $\vep>0$, the set $G_\vep(X,I^{2n+1})$ is dense in $C(X,I^{2n+1})$.
\end{prp}
\begin{proof}
Let $f\in C(X,I^{2n+1})$ and $\eta>0$. We need to find $g\in G_\vep(X,I^{2n+1})$ such that $d_u(f,g)<\eta$. Since $f$ is uniformly continuous (as a continuous map from a compact metric space to a metric space), there exists $0<\delta<\vep$ such that $d(x,x')<\delta$ implies $d\big(f(x),f(x')\big)<\eta/2$.

By Corollary \ref{IFScrl10}, $X$ has an open cover $\U=\{U_1,...,U_r\}$ of diameter $\max_i\diam U_i<\delta$ and $\Ord\U\leq n$. It follows that
\bea
\diam U_i<\delta~~~~\Ra~~~~\diam f(U_i)<\eta/2,~~~~\txt{for each}~~i=1,...,r.\nn
\eea

Select points $p_i\in I^{2n+1}$ that are (i) in general position (i.e., $\{p_1-p_i:i=2,...,r\}$ are linearly independent) in $\Real^{2n+1}$ and (ii) satisfy
\bea
\dist(p_i,f(U_i))<\eta/2,~~~~i=1,...,r.\nn
\eea

For each $i=1,\cdots,r$ define $w_i:X\ra\Real$ by $w_i(x):=\dist(x,X-U_i)$. It is clear that $w_i(x)=0$ for $x\not\in U_i$ and $w_i(x)>0$ for $x\in U_i$. Since the $U_i$'s cover $X$, it follows for each $x\in X$, at least one $w_i(x)>0$. Define $g:X\ra I^{2n+1}$ by $g(x):=\sum_iw_i(x)p_i/\sum_iw_i(x)$, which is continuous as a ratio of continuous functions since each $w_i$ is continuous and $h(x):=\sum_iw_i(x)p_i$ satisfies $\|h(x)-h(x')\|\leq\sum_i|w_i(x)-w_i(x')|$. We will now show $d_u(f,g)<\eta$ and $g\in G_\vep(X,I^{2n+1})$.

To show $d_u(f,g)<\eta$, fix $x\in X$. Renumber the $U_i$ so that $\{U_1,...,U_s\}$ is the set of all $U\in\U$ that contain $x$. Then $w_i(x)=\dist(x,X-U_i)>0$ for $i\leq s$, and $w_i(x)=\dist(x,X-U_i)=0$ for $i>s$. Thus, $g(x)=\sum_{i=1}^sw_i(x)p_i/\sum_{i=1}^sw_i(x)$, i.e., only $p_1,...,p_s$ terms contribute to $g(x)$. So, for each $i\leq s$,
{\small\begin{align}
d(p_i,f(x))\leq\inf_{x'\in U_i}\left[d(p_i,f(x'))+d(f(x'),f(x))\right]\leq\dist(p_i,f(U_i))+\diam f(U_i)<\eta/2+\eta/2=\eta.\nn
\end{align}}
It follows from the triangle inequality that $d(f(x),g(x))<\eta$, and so $d_u(f,g)<\eta$.

To show $g\in G_\vep(X,I^{2n+1})$, fix $x\in X$. Let $U_{i_1},...,U_{i_s}$ be all $U\in\U$ that contain $x$. Since $\Ord\U\leq n$, we have $s\leq n+1$. Consider the $(s-1)$-dimensional ``linear'' subspace of $I^{2n+1}$ given by
{\small\begin{align}
&\textstyle L(x):=I^{2n+1}\cap\Span_\Real\Conv\{p_{i_1},...,p_{i_s}\}=I^{2n+1}\cap\Span_\Real\left\{\sum_{k=1}^s\al_kp_{i_k}:\al_k\geq 0,~\sum\al_k=1\right\}\nn\\
&\textstyle~~~~=I^{2n+1}\cap\Span_\Real\left\{p_{i_1}+\sum_{k=2}^s\al_k(p_{i_k}-p_{i_1}):\al_k\geq 0,~\sum_{k=2}^s\al_k\leq 1\right\}.~~~~~~\txt{(note $g(x)\in L(x)$)}\nn
\end{align}}
Let $x'\in X$ be another point. Then with $U_{i'_1},...,U_{i'_{s'}}$ being all $U\in\U$ containing $x'$, we have
{\small\begin{align}
&\textstyle L(x'):=I^{2n+1}\cap\Span_\Real\Conv\{p_{i'_1},...,p_{i'_{s'}}\}=I^{2n+1}\cap\Span_\Real\left\{\sum_{k=1}^{s'}\beta_kp_{i'_k}:\beta_k\geq 0,~\sum\beta_k=1\right\}\nn\\
&\textstyle~~~~=I^{2n+1}\cap\Span_\Real\left\{p_{i'_1}+\sum_{k=2}^{s'}\beta_k(p_{i'_k}-p_{i'_1}):\beta_k\geq 0,~\sum_{k=2}^s\beta_k\leq 1\right\},\nn
\end{align}}
where $\Ord\U\leq n$ again implies $s'\leq n+1$. $\bullet$ Now, if $L(x)\cap L(x')\neq\emptyset$, then some $p_i\in L(x)\cap L(x')$, [hence $x,x'\in U_i$]. \ul{Indeed} if $0\neq y\in L(x)\cap L(x')$, then $y=\sum\al_kp_{i_k}=\sum\beta_{k'}p_{i'_{k'}}$ is a nontrivial relation among the $p$'s involved (``eliminating'' at least one $p$), and so $L(x,x'):=I^{2n+1}\cap\Span_\Real\Conv\{p_{i_1},...,p_{i_s}\}\cup\{p_{i'_1},...,p_{i'_{s'}}\}$ has dimension (compared to it dimension without intersection $\leq s+s'-1$)
\bea
\dim L(x,x')\leq (s+s'-1)-1=s+s'-2~~~~\big[~\leq (n+1)+(n+1)-2=2n~\big].\nn
\eea
But since all of the $p$'s are in general position, if none of the $p_{i_{k}}$ is equal to some $p_{i'_{k'}}$, then $\dim L(x,x')=s+s'-1>s+s'-2$ (a contradiction).

Finally, for any $x,x'\in X$, if $g(x)=g(x')$, i.e., $x,x'\in g^{-1}(t)$ for some $t\in g(X)$, then $L(x)\cap L(x')\neq\emptyset$, since it contains $g(x)=g(x')$. Thus $x,x'\in U_i$ for some $i$, and so ~$d(x,x')\leq\diam U_i<\delta<\vep$.
\end{proof}

\begin{thm}[\textcolor{OliveGreen}{\cite[Theorem V2, p.56]{HureWall41}}]\label{IFSthmV2}
Let $X$ be a compact metric space of dimension $\leq n$. Then $X$ imbeds into $I^{2n+1}$ (i.e., $X$ is homeomorphic to a subspace of $I^{2n+1}$).
\end{thm}
\begin{proof}
As before, let $G_{1/i}(X,I^{2n+1})$ be all $1/i$-maps from $X$ to $I^{2n+1}$. By Lemma \ref{IFSlmm13} the set
\bea
\textstyle H:=\bigcap_{i=1}^\infty G_{1/i}(X,I^{2n+1})~\subset~C(X,I^{2n+1})\nn
\eea
consists of imbeddings of $X$ into $I^{2n+1}$. Also, by Lemma \ref{IFSlmm14} each $G_{1/i}(X,I^{2n+1})$ is open (and thus a $G_\delta$ set), and by Proposition \ref{IFSlmm15} each $G_{1/i}(X,I^{2n+1})$ is dense in the complete metric space $C(X,I^{2n+1})$. It thus follows from Baire's category theorem that $H$ is also dense in $C(X,I^{2n+1})$.
\end{proof}

{\flushleft \dotfill}
\vspace{0.3cm}

We will now discuss the general case of the above theorem with compactness removed.

\begin{dfn}[\textcolor{blue}{Cover-map, Set of cover-maps}]
Let $X$ be a space and $\U=\{U_\al\}_{\al\in A}$ an open cover of $X$. A map of spaces $f:X\ra Y$ is a $\U$-map if every point $y\in Y$ has a neighborhood $V\ni y$ whose preimage $f^{-1}(V)\subset U_\al\in\U$ for some $\al$. We will denote the set of continuous $\U$-maps from $X$ to $Y$ by ~$G_\U(X,Y):=\{\txt{$\U$-maps}~g\in C(X,Y)\}$.
\end{dfn}

\begin{dfn}[\textcolor{blue}{Cover-neighborhood of a point}]
Let $X$ be a space and $\U$ an open cover of $X$. The $\U$-neighborhood (or $\U$-star) $S_\U(x)$ of a point $x\in X$ is the union of all elements of $\U$ containing $x$, i.e., ~$S_\U(x)~:=~\bigcup\{U\in\U:x\in U\}$.
\end{dfn}
\begin{dfn}[\textcolor{blue}{Basic sequence of covers}]
Let $X$ be a space. A countable collection of open covers $\U_1$, $\U_2$, $\cdots$ of $X$ is a basic sequence if every neighborhood $U$ of any point $x\in X$ contains at least one of the cover-neighborhoods (or stars) ~$S_{\U_1}(x)$, $S_{\U_2}(x)$, $\cdots$.
\end{dfn}

\begin{lmm}\label{BaSeqExist}
Every regular space $X$ with a countable base $U_1,U_2,\cdots$ has a basic sequence of covers.
\end{lmm}
\begin{proof}
Since $X$ is regular, for any $U_m\neq\emptyset$, with any point $x\in U_m$, we have an open set $V\ni x$ such that $\ol{V}\subset U_m$, and hence a $U_n\ni x$ such that $\ol{U}_n\subset U_m$.
For $n,m\geq 1$, consider all pairs of nonempty base elements $U_n,U_m$ satisfying $\ol{U}_n\subset U_m$. Then we have the countable collection of open covers $\U_{n,m}:=\{X-\ol{U}_n,U_m\}$. Fix $x\in X$. Then it is clear that $x\in U_n\subset\ol{U}_n\subset U_m\subset U$ for some $n,m$. Since $x\in U_n$, we have $S_{\U_{m,n}}(x)=U_m\subset U$. Hence $\{\U_{n,m}\}$ is a basic sequence of covers.
\end{proof}

\begin{lmm}\label{BaSeqHomeo}
Let $X$ be a regular Hausdorff space and $\U_1$, $\U_2$, $\cdots$ a basic sequence of open covers of $X$. If a continuous map of spaces $g:X\ra Y$ is a $\U_i$-map for all $i$, then $g$ is an imbedding (i.e., a homeomorphism onto its image).
\end{lmm}
\begin{proof}
Let $x\in X$ and $U\ni x$ an open set. Then some $S_{\U_i}(x)\subset U$ (since $\{\U_i\}$ is a basic sequence of covers). Moreover, for some open set $V\ni g(x)$, we have $g^{-1}(V)\subset U_i$ for some $U_i\in\U_i$. Since $x\in g^{-1}(V)\subset U_i$ implies $x\in U_i$, it follows that $g^{-1}(V)\subset U_i\subset S_{\U_i}(x)\subset U$. That is, for any nbd $U\ni x$, there exists a nbd $V\ni g(x)$ such that $g^{-1}(V)\subset U$.

It follows that $g$ is injective: Indeed if $x\neq x'$, let $N(x),N(x')$ be disjoint neighborhoods. If $g(x)=g(x')$, then there is an open set $V\ni g(x)=g(x')$ such that $g^{-1}(V)\subset N(x)\cap N(x')$, which is a contradiction.

Also it is clear that $g$ is an open map, which means $g^{-1}$ is continuous.
\end{proof}

\begin{lmm}[\textcolor{OliveGreen}{Lebesgue number lemma}]\label{LebNoLmm}
Let $(X,d)$ be a compact metric space and $\U=\{U_\al\}$ an open cover of $X$. Then there exists a number $\delta=\delta_\U>0$ such that for all $A\subset X$,
\bea
\txt{diam}~A~<\delta~~\Ra~~A\subset U_\al~~\txt{for some}~~\al.\nn
\eea
\end{lmm}
\begin{proof}
Since $\{U_\al\}$ is a cover of $X$, for every $x\in X$ some $U_{\al_x}\ni x$, and since $U_{\al_x}$ is open, some $B(x,\vep_x)\subset U_{\al_x}$, where $\vep_x>0$. The sets $\{B(x,\vep_x/2)\}_{x\in X}$ form an open cover of $X$, and so by compactness,
\bea
X=B(x_1,\vep_{x_1/2})\cup\cdots\cup B(x_n,\vep_{x_n/2})~~~~\txt{for some}~~~~x_1,...,x_n\in X.\nn
\eea
Let $\delta:=\min\{\vep_{x_1}/2,...,\vep_{x_n}/2\}$. If $A\subset X$ such that $\dim A<\delta$, pick any $a\in A$ and let $a\in B(x_{i_a},\vep_{x_{i_a}}/2)$. Then
\bea
A\subset B\left(x_{i_a},\delta+\vep_{x_{i_a}}/2\right)\subset B\left(x_{i_a},\vep_{x_{i_a}}/2+\vep_{x_{i_a}}/2\right)=B\left(x_{i_a},\vep_{x_{i_a}}\right)\subset U_{\al_{x_{i_a}}}.\nn\qedhere
\eea
\end{proof}

\begin{lmm}\label{IFSlmm16}
Let $X$ be a metric space and $Y$ a compact metric space. For any cover $\U$ of $X$, the set $G_\U(X,Y)$ is open in $C(X,Y)$.
\end{lmm}
\begin{proof}
Let $g\in G_\U(X,Y)$. That is, $g$ is continuous and every $y\in Y$ has a neighborhood $V_y\ni y$ such that $g^{-1}(V_y)\subset U$ for some $U\in\U$. Since $Y$ is compact, let $Y\subset \bigcup_{i=1}^kV_{y^i}$. By Lemma \ref{LebNoLmm}, let $\delta>0$ be a number such that every set $A\subset Y$ of diameter $\diam<\delta$ is contained in some $V_{y^i}$, and hence $g^{-1}(A)\subset g^{-1}(V_{y^i})\subset U$ for some $U\in\U$.

Let $f\in C(X,Y)$ be any continuous map such that $d_u(f,g)<\delta/3$, i.e., $f\in B_{\delta/3}(g)\subset C(X,Y)$. For $y\in Y$, let $A:=g^{-1}\big(B_{\delta/3}(y)\big)$. Then $g(A)\subset B_{\delta/3}(y)$ implies $\diam(A)<\delta$, and so $A\subset U$ for some $U\in\U$, which implies $f\in G_\U(X,Y)$ (i.e., $f$ is a continuous $\U$-map) since
\begin{align}
&f^{-1}\big(B_{\delta/3}(y)\big)=\{x:d(f(x),y)<\delta/3\}\subset \{x:d(g(x),y)\leq d(g(x),f(x))+d(f(x),y)<\delta\}\nn\\
&~~~~=g^{-1}\big(B_\delta(y)\big).\nn\qedhere
\end{align}
\end{proof}

\begin{prp}\label{IFSlmm17}
Let $X$ be a first countable regular Hausdorff space such that $\diam X\leq n$. For any open cover $\W$ of $X$, the set $G_\W(X,I^{2n+1})$ is dense in $C(X,I^{2n+1})$.
\end{prp}
\begin{proof}
Fix $f\in C(X,I^{2n+1})$ and $\eta>0$. We need to fine $g\in G_\W(X,I^{2n+1})$ such that $d_u(f,g)<\eta$. By its compactness, $I^{2n+1}$ has a finite open cover $\V=\{V_1,\cdots,V_r\}$ of diameter $<\vep/2$. By Theorem \ref{IFSthmV1}, the open covers $\W$ and $\{f^{-1}(V_1),\cdots,f^{-1}(V_r)\}$ of $X$ have a common open refinement $\U=\{U_1,\cdots,U_r\}$ of order $\leq n$, and since $U_i\subset f^{-1}(V_j)$ for some $j$, we have $f(U_i)\subset V_j$, and so
\bea
\diam f(U_i)\leq\diam V_j<\eta/2,~~~~i=1,...,r.\nn
\eea

Select points $p_i\in I^{2n+1}$ that are (i) in general position (i.e., $\{p_1-p_i:i=2,...,r\}$ are linearly independent) in $\Real^{2n+1}$ and (ii) satisfy
\bea
\dist(p_i,f(U_i))<\eta/2,~~~~i=1,...,r.\nn
\eea

For each $i=1,\cdots,r$ define $w_i:X\ra\Real$ by $w_i(x):=\dist(x,X-U_i)$. It is clear that $w_i(x)=0$ for $x\not\in U_i$ and $w_i(x)>0$ for $x\in U_i$. Since the $U_i$'s cover $X$, it follows for each $x\in X$, at least one $w_i(x)>0$. Define $g:X\ra I^{2n+1}$ by $g(x):=\sum_iw_i(x)p_i/\sum_iw_i(x)$, which is continuous as a ratio of continuous functions since each $w_i$ is continuous and $h(x):=\sum_iw_i(x)p_i$ satisfies $\|h(x)-h(x')\|\leq\sum_i|w_i(x)-w_i(x')|$. We will now show $d_u(f,g)<\eta$ and $g\in G_\U(X,I^{2n+1})\subset G_\W(X,I^{2n+1})$.

(The proof that $d_u(f,g)<\eta$ proceeds exactly as in the proof of Proposition \ref{IFSlmm15}.)

To show $d_u(f,g)<\eta$, fix $x\in X$. Renumber the $U_i$ so that $\{U_1,...,U_s\}$ is the set of all $U\in\U$ that contain $x$. Then $w_i(x)=\dist(x,X-U_i)>0$ for $i\leq s$, and $w_i(x)=\dist(x,X-U_i)=0$ for $i>s$. Thus, $g(x)=\sum_{i=1}^sw_i(x)p_i/\sum_{i=1}^sw_i(x)$, i.e., only $p_1,...,p_s$ terms contribute to $g(x)$. So, for each $i\leq s$,
{\small\begin{align}
d(p_i,f(x))\leq\inf_{x'\in U_i}\left[d(p_i,f(x'))+d(f(x'),f(x))\right]\leq\dist(p_i,f(U_i))+\diam f(U_i)<\eta/2+\eta/2=\eta.\nn
\end{align}}
It follows from the triangle inequality that $d(f(x),g(x))<\eta$, and so $d_u(f,g)<\eta$.

(The proof that $d_u(f,g)<\eta$ proceeds as in the proof of Proposition \ref{IFSlmm15}, except for the last paragraph.)

To show $g\in G_\U(X,I^{2n+1})$, fix $x\in X$. Let $U_{i_1},...,U_{i_s}$ be all $U\in\U$ that contain $x$. Since $\Ord\U\leq n$, we have $s\leq n+1$. Consider the $(s-1)$-dimensional ``linear'' subspace of $I^{2n+1}$ given by
{\small\begin{align}
&\textstyle L(x):=I^{2n+1}\cap\Span_\Real\Conv\{p_{i_1},...,p_{i_s}\}=I^{2n+1}\cap\Span_\Real\left\{\sum_{k=1}^s\al_kp_{i_k}:\al_k\geq 0,~\sum\al_k=1\right\}\nn\\
&\textstyle~~~~=I^{2n+1}\cap\Span_\Real\left\{p_{i_1}+\sum_{k=2}^s\al_k(p_{i_k}-p_{i_1}):\al_k\geq 0,~\sum_{k=2}^s\al_k\leq 1\right\}.~~~~~~\txt{(note $g(x)\in L(x)$)}\nn
\end{align}}
Let $x'\in X$ be another point. Then with $U_{i'_1},...,U_{i'_{s'}}$ being all $U\in\U$ containing $x'$, we have
{\small\begin{align}
&\textstyle L(x'):=I^{2n+1}\cap\Span_\Real\Conv\{p_{i'_1},...,p_{i'_{s'}}\}=I^{2n+1}\cap\Span_\Real\left\{\sum_{k=1}^{s'}\beta_kp_{i'_k}:\beta_k\geq 0,~\sum\beta_k=1\right\}\nn\\
&\textstyle~~~~=I^{2n+1}\cap\Span_\Real\left\{p_{i'_1}+\sum_{k=2}^{s'}\beta_k(p_{i'_k}-p_{i'_1}):\beta_k\geq 0,~\sum_{k=2}^s\beta_k\leq 1\right\},\nn
\end{align}}
where $\Ord\U\leq n$ again implies $s'\leq n+1$. $\bullet$ Now, if $L(x)\cap L(x')\neq\emptyset$, then some $p_i\in L(x)\cap L(x')$, [hence $x,x'\in U_i$]. \ul{Indeed} if $0\neq y\in L(x)\cap L(x')$, then $y=\sum\al_kp_{i_k}=\sum\beta_{k'}p_{i'_{k'}}$ is a nontrivial relation among the $p$'s involved (``eliminating'' at least one $p$), and so $L(x,x'):=I^{2n+1}\cap\Span_\Real\Conv\{p_{i_1},...,p_{i_s}\}\cup\{p_{i'_1},...,p_{i'_{s'}}\}$ has dimension (compared to it dimension without intersection $\leq s+s'-1$)
\bea
\dim L(x,x')\leq (s+s'-1)-1=s+s'-2~~~~\big[~\leq (n+1)+(n+1)-2=2n~\big].\nn
\eea
But since all of the $p$'s are in general position, if none of the $p_{i_{k}}$ is equal to some $p_{i'_{k'}}$, then $\dim L(x,x')=s+s'-1>s+s'-2$ (a contradiction).

Since there are only a finite number number of the linear subspaces $L(x)$, as the $p$'s are finite, there exists a number (which we can choose to be) $\eta>0$ such that any two of the subspaces $L(x),L(x')$ either meet or else have distance $\dist\big(L(x),L(x')\big)\geq\eta$ from each other. If $d\big(g(x),g(x')\big)<\eta$, then it is clear that $\dist\big(L(x),L(x')\big)\leq d\big(g(x),g(x')\big)<\eta$ since $g(x)\in L(x)$, $g(x')\in L(x')$, and so $L(x),L(x')$ meet. This implies $x,x'\in U_i$ for some $i$, which shows $g$ is a $\U$-map.
\end{proof}

\begin{thm}[\textcolor{OliveGreen}{\cite[Theorem V3]{HureWall41}}]\label{ImbedThm2}
Let $X$ be a metric space. If $\dim X\leq n$, then $X$ is homeomorphic to a subspace of $I^{2n+1}$. (The converse is trivially true as well). Moreover, the set of homeomorphisms $X\ra I^{2n+1}$ is dense in $C(X,I^{2n+1})$.
\end{thm}
\begin{proof}
Let $\U_1,\U_2,\cdots$ be a basic-sequence of covers of $X$, and let $H:=\bigcap_{i=1}^\infty G_{\U_i}(X,I^{2n+1})$. Each element of $H$ is an imbedding of $X$ into $I^{2n+1}$ by Lemma \ref{BaSeqHomeo}. Each $G_{\U_i}(X,I^{2n+1})$ is open (hence a $G_\delta$ set) in $C(X,I^{2n+1})$ by Lemma \ref{IFSlmm16}. By Proposition \ref{IFSlmm17} therefore, each $G_{\U_i}(X,I^{2n+1})$ is a dense $G_\delta$ set in $C(X,I^{2n+1})$, and so by Baire's category theorem, $H$ is dense in $C(X,I^{2n+1})$.
\end{proof}


}

\backmatter
\bibliographystyle{amsplain}

\begin{bibdiv}
\begin{biblist}
\bibitem{aguilar.etal2002} M. Aguilar, S. Gitler, C. Prieto, \emph{Algebraic Topology from a Homotopical Viewpoint}, Universitext, Springer (2002).
\bibitem{ProbList} AimPL: \emph{Mapping theory in metric spaces}, available at \url{http://aimpl.org/mappingmetric} (See Problem 1.4).
\bibitem{akofor2019} Akofor, E., \emph{On Lipschitz retraction of finite subsets of normed spaces}, Isr. J. Math. (2019). \url{https://doi.org/10.1007/s11856-019-1935-x}

\bib{ambro-prodi1993}{book}{
   author={Ambrosetti, Antonio},
   author={Prodi, Giovanni},
   title={A primer of nonlinear analysis},
   series={Cambridge Studies in Advanced Mathematics},
   volume={34},
   note={Corrected reprint of the 1993 original},
   publisher={Cambridge University Press, Cambridge},
   date={1995},
   pages={viii+171},
   isbn={0-521-48573-8},
}

\bib{artin-book}{book}{
   author={Artin, Michael},
   title={Algebra},
   publisher={Prentice Hall, Inc., Englewood Cliffs, NJ},
   date={1991},
   pages={xviii+618},
   isbn={0-13-004763-5},
}

\bib{baake-schlaegel2011}{article}{
   author={Baake, Michael},
   author={Schl\"{a}gel, Ulrike},
   title={The Peano-Baker series},
   journal={Tr. Mat. Inst. Steklova},
   volume={275},
   date={2011},
   number={Klassicheskaya i Sovremennaya Matematika v Pole Deyatel\cprime nosti
   Borisa Nikolaevicha Delone},
   pages={167--171},
   issn={0371-9685},
   translation={
      journal={Proc. Steklov Inst. Math.},
      volume={275},
      date={2011},
      number={1},
      pages={155--159},
      issn={0081-5438},
   },
}

\bib{bacac-kovalev2016}{article}{
   author={Ba\v{c}\'{a}k, Miroslav},
   author={Kovalev, Leonid V.},
   title={Lipschitz retractions in Hadamard spaces via gradient flow
   semigroups},
   note={[Paging previously given as 1--9]},
   journal={Canad. Math. Bull.},
   volume={59},
   date={2016},
   number={4},
   pages={673--681},
   issn={0008-4395},
}

\bib{ben-lind2000}{book}{
   author={Benyamini, Yoav},
   author={Lindenstrauss, Joram},
   title={Geometric nonlinear functional analysis. Vol. 1},
   series={American Mathematical Society Colloquium Publications},
   volume={48},
   publisher={American Mathematical Society, Providence, RI},
   date={2000},
   pages={xii+488},
   isbn={0-8218-0835-4},
}

\bib{borsuk-ulam1}{article}{
   author={Borsuk, Karol},
   author={Ulam, Stanislaw},
   title={On symmetric products of topological spaces},
   journal={Bull. Amer. Math. Soc.},
   volume={37},
   date={1931},
   number={12},
   pages={875--882},
   issn={0002-9904},
}

\bib{borsuk1949}{article}{
   author={Borsuk, Karol},
   title={On the third symmetric potency of the circumference},
   journal={Fund. Math.},
   volume={36},
   date={1949},
   pages={236--244},
   issn={0016-2736},
}

\bib{BorovIbrag2009}{article}{
   author={Borovikova, Marina},
   author={Ibragimov, Zair},
   title={The third symmetric product of $\Bbb R$},
   journal={Comput. Methods Funct. Theory},
   volume={9},
   date={2009},
   number={1},
   pages={255--268},
}

\bib{BorovEtal2010}{article}{
   author={Borovikova, Marina},
   author={Ibragimov, Zair},
   author={Yousefi, Hassan},
   title={Symmetric products of the real line},
   journal={J. Anal.},
   volume={18},
   date={2010},
   pages={53--67},
   issn={0971-3611},
}

\bib{bott1952}{article}{
   author={Bott, R.},
   title={On the third symmetric potency of $S_1$},
   journal={Fund. Math.},
   volume={39},
   date={1952},
   pages={264--268 (1953)},
   issn={0016-2736},
}

\bib{brudnyis2012v1}{book}{
   author={Brudnyi, Alexander},
   author={Brudnyi, Yuri},
   title={Methods of geometric analysis in extension and trace problems.
   Volume 1},
   series={Monographs in Mathematics},
   volume={102},
   publisher={Birkh\"{a}user/Springer Basel AG, Basel},
   date={2012},
   pages={xxiv+560},
   isbn={978-3-0348-0208-6},
}

\bib{brudnyis2012v2}{book}{
   author={Brudnyi, Alexander},
   author={Brudnyi, Yuri},
   title={Methods of geometric analysis in extension and trace problems.
   Volume 2},
   series={Monographs in Mathematics},
   volume={103},
   publisher={Birkh\"{a}user/Springer Basel AG, Basel},
   date={2012},
   pages={xx+414},
   isbn={978-3-0348-0211-6},
}

\bib{BBI}{book}{
   author={Burago, Dmitri},
   author={Burago, Yuri},
   author={Ivanov, Sergei},
   title={A course in metric geometry},
   series={Graduate Studies in Mathematics},
   volume={33},
   publisher={American Mathematical Society, Providence, RI},
   date={2001},
   pages={xiv+415},
   isbn={0-8218-2129-6},
}

\bib{chernoff-1992}{article}{
   author={Chernoff, Paul R.},
   title={A simple proof of Tychonoff's theorem via nets},
   journal={Amer. Math. Monthly},
   volume={99},
   date={1992},
   number={10},
   pages={932--934},
   issn={0002-9890},
}

\bib{chinen2015}{article}{
   author={Chinen, Naotsugu},
   title={Symmetric products of the Euclidean spaces and the spheres},
   journal={Comment. Math. Univ. Carolin.},
   volume={56},
   date={2015},
   number={2},
   pages={209--221},
   issn={0010-2628},
}

\bib{chinen-koyama2010}{article}{
   author={Chinen, Naotsugu},
   author={Koyama, Akira},
   title={On the symmetric hyperspace of the circle},
   journal={Topology Appl.},
   volume={157},
   date={2010},
   number={17},
   pages={2613--2621},
   issn={0166-8641},
}

\bib{chinen2018}{article}{
   author={Chinen, Naotsugu},
   title={On isometries of symmetric products of metric spaces},
   journal={Topology Appl.},
   volume={248},
   date={2018},
   pages={24--39},
   issn={0166-8641},
}

\bib{conway}{book}{
   author={Conway, John B.},
   title={A course in functional analysis},
   series={Graduate Texts in Mathematics},
   volume={96},
   edition={2},
   publisher={Springer-Verlag, New York},
   date={1990},
   pages={xvi+399},
   isbn={0-387-97245-5},
}

\bib{CoronaEtal2017}{article}{
   author={Corona-V\'{a}zquez, F.},
   author={Qui\~{n}ones-Estrella, R. A.},
   author={Villanueva, H.},
   title={Embedding cones over trees into their symmetric products},
   journal={Topology Appl.},
   volume={231},
   date={2017},
   pages={77--91},
   issn={0166-8641},
}

\bib{deimling}{book}{
   author={Deimling, Klaus},
   title={Nonlinear functional analysis},
   publisher={Springer-Verlag, Berlin},
   date={1985},
   pages={xiv+450},
   isbn={3-540-13928-1},
}

\bib{dugundji1951}{article}{
   author={Dugundji, J.},
   title={An extension of Tietze's theorem},
   journal={Pacific J. Math.},
   volume={1},
   date={1951},
   pages={353--367},
   issn={0030-8730},
}

\bib{dugundji1966}{book}{
   author={Dugundji, James},
   title={Topology},
   publisher={Allyn and Bacon, Inc., Boston, Mass.},
   date={1966},
   pages={xvi+447},
}

\bib{dunkl1964}{article}{
   author={Dunkl, C. F.},
   author={Williams, K. S.},
   title={Mathematical notes: a simple norm inequality},
   journal={Amer. Math. Monthly},
   volume={71},
   date={1964},
   number={1},
   pages={53--54},
   issn={0002-9890},
}



\bib{enderton1977}{book}{
   author={Enderton, Herbert B.},
   title={Elements of set theory},
   publisher={Academic Press [Harcourt Brace Jovanovich, Publishers], New
   York-London},
   date={1977},
   pages={xiv+279},
}

\bib{espinola-khamsi}{article}{
   author={Esp\'{\i}nola, R.},
   author={Khamsi, M. A.},
   title={Introduction to hyperconvex spaces},
   conference={
      title={Handbook of metric fixed point theory},
   },
   book={
      publisher={Kluwer Acad. Publ., Dordrecht},
   },
   date={2001},
   pages={391--435},
}

\bib{espinola.etal}{article}{
   author={Esp\'{\i}nola, Rafa},
   author={Madiedo, \'{O}scar},
   author={Nicolae, Adriana},
   title={Borsuk-Dugundji type extension theorems with Busemann convex
   target spaces},
   journal={Ann. Acad. Sci. Fenn. Math.},
   volume={43},
   date={2018},
   number={1},
   pages={225--238},
   issn={1239-629X},
}

\bib{ferry-weinbg2013}{article}{
   author={Ferry, Steve},
   author={Weinberger, Shmuel},
   title={Quantitative algebraic topology and Lipschitz homotopy},
   journal={Proc. Natl. Acad. Sci. USA},
   volume={110},
   date={2013},
   number={48},
   pages={19246--19250},
   issn={0027-8424},
}

\bib{folland-book}{book}{
   author={Folland, Gerald B.},
   title={Real analysis},
   series={Pure and Applied Mathematics (New York)},
   edition={2},
   note={Modern techniques and their applications;
   A Wiley-Interscience Publication},
   publisher={John Wiley \& Sons, Inc., New York},
   date={1999},
   pages={xvi+386},
   isbn={0-471-31716-0},
}

\bibitem{goldrei1996} Derek Goldrei, \emph{Classic Set Theory: For Guided Independent Study}, Chapman $\&$ Hall/CRC, 1996.

\bib{hakobyan-herron2008}{article}{
   author={Hakobyan, Hrant},
   author={Herron, David A.},
   title={Euclidean quasiconvexity},
   journal={Ann. Acad. Sci. Fenn. Math.},
   volume={33},
   date={2008},
   number={1},
   pages={205--230},
   issn={1239-629X},
}

\bib{hatcher2001}{book}{
   author={Hatcher, Allen},
   title={Algebraic topology},
   publisher={Cambridge University Press, Cambridge},
   date={2001},
   pages={xii+544},
   isbn={0-521-79160-X},
   isbn={0-521-79540-0},
}

\bib{heinonen}{book}{
   author={Heinonen, Juha},
   title={Lectures on Lipschitz analysis},
   series={Report. University of Jyv\"{a}skyl\"{a} Department of Mathematics and
   Statistics},
   volume={100},
   publisher={University of Jyv\"{a}skyl\"{a}, Jyv\"{a}skyl\"{a}},
   date={2005},
   pages={ii+77},
   isbn={951-39-2318-5},
}

\bib{heinonen2}{book}{
   author={Heinonen, Juha},
   title={Lectures on analysis on metric spaces},
   series={Universitext},
   publisher={Springer-Verlag, New York},
   date={2001},
   pages={x+140},
   isbn={0-387-95104-0},
}

\bib{howes1995}{book}{
   author={Howes, Norman R.},
   title={Modern analysis and topology},
   series={Universitext},
   publisher={Springer-Verlag, New York},
   date={1995},
   pages={xxx+403},
   isbn={0-387-97986-7},
}

\bib{HureWall41}{book}{
   author={Hurewicz, Witold},
   author={Wallman, Henry},
   title={Dimension Theory},
   series={Princeton Mathematical Series, v. 4},
   publisher={Princeton University Press, Princeton, N. J.},
   date={1941},
   pages={vii+165},
}

\bib{illanes-nadler}{book}{
   author={Illanes, Alejandro},
   author={Nadler, Sam B., Jr.},
   title={Hyperspaces},
   series={Monographs and Textbooks in Pure and Applied Mathematics},
   volume={216},
   note={Fundamentals and recent advances},
   publisher={Marcel Dekker, Inc., New York},
   date={1999},
   pages={xx+512},
   isbn={0-8247-1982-4},
}

\bib{IllanMart2017}{article}{
   author={Illanes, Alejandro},
   author={Mart\'{\i}nez-de-la-Vega, Ver\'{o}nica},
   title={Symmetric products as cones},
   journal={Topology Appl.},
   volume={228},
   date={2017},
   pages={36--46},
   issn={0166-8641},
}

\bib{kelley1975}{book}{
   author={Kelley, John L.},
   title={General topology},
   note={Reprint of the 1955 edition [Van Nostrand, Toronto, Ont.];
   Graduate Texts in Mathematics, No. 27},
   publisher={Springer-Verlag, New York-Berlin},
   date={1975},
   pages={xiv+298},
}

\bib{koskela94}{article}{
   author={Koskela, Pekka},
   title={The degree of regularity of a quasiconformal mapping},
   journal={Proc. Amer. Math. Soc.},
   volume={122},
   date={1994},
   number={3},
   pages={769--772},
   issn={0002-9939},
}

\bib{kovalev2016}{article}{
   author={Kovalev, Leonid V.},
   title={Lipschitz retraction of finite subsets of Hilbert spaces},
   journal={Bull. Aust. Math. Soc.},
   volume={93},
   date={2016},
   number={1},
   pages={146--151},
   issn={0004-9727},
}

\bib{kovalev2015}{article}{
   author={Kovalev, Leonid V.},
   title={Symmetric products of the line: embeddings and retractions},
   journal={Proc. Amer. Math. Soc.},
   volume={143},
   date={2015},
   number={2},
   pages={801--809},
   issn={0002-9939},
}

\bibitem{lewin1991} Lewin, J. (1991). \emph{A simple proof of zorn's lemma.} The American Mathematical Monthly, 98(4), pp. 353-354.

\bib{Luu-Sak98}{article}{
   author={Luukkainen, Jouni},
   author={Saksman, Eero},
   title={Every complete doubling metric space carries a doubling measure},
   journal={Proc. Amer. Math. Soc.},
   volume={126},
   date={1998},
   number={2},
   pages={531--534},
   issn={0002-9939},
}

\bib{Luu98}{article}{
   author={Luukkainen, Jouni},
   title={Assouad dimension: antifractal metrization, porous sets, and
   homogeneous measures},
   journal={J. Korean Math. Soc.},
   volume={35},
   date={1998},
   number={1},
   pages={23--76},
}

\bib{malig2006}{article}{
   author={Maligranda, Lech},
   title={Simple norm inequalities},
   journal={Amer. Math. Monthly},
   volume={113},
   date={2006},
   number={3},
   pages={256--260},
   issn={0002-9890},
}

\bib{mckemie1987}{article}{
   author={McKemie, M. Jean},
   title={Quasiconformal groups with small dilatation},
   journal={Ann. Acad. Sci. Fenn. Ser. A I Math.},
   volume={12},
   date={1987},
   number={1},
   pages={95--118},
   issn={0066-1953},
}

\bib{michael1953J}{article}{
   author={Michael, Ernest},
   title={Some extension theorems for continuous functions},
   journal={Pacific J. Math.},
   volume={3},
   date={1953},
   pages={789--806},
   issn={0030-8730},
}

\bib{michael1953}{article}{
   author={Michael, Ernest},
   title={A note on paracompact spaces},
   journal={Proc. Amer. Math. Soc.},
   volume={4},
   date={1953},
   pages={831--838},
   issn={0002-9939},
}

\bib{mostovoy2004}{article}{
   author={Mostovoy, Jacob},
   title={Lattices in $\Bbb C$ and finite subsets of a circle},
   journal={Amer. Math. Monthly},
   volume={111},
   date={2004},
   number={4},
   pages={357--360},
   issn={0002-9890},
}

\bib{munkres}{book}{
   author={Munkres, James R.},
   title={Topology},
   note={Second edition},
   publisher={Prentice Hall, Inc., Upper Saddle River, NJ},
   date={2000},
   pages={xvi+537},
   isbn={0-13-181629-2},
}

\bib{papado2014}{book}{
   author={Papadopoulos, Athanase},
   title={Metric spaces, convexity and nonpositive curvature},
   series={IRMA Lectures in Mathematics and Theoretical Physics},
   volume={6},
   publisher={European Mathematical Society (EMS), Z\"{u}rich},
   date={2005},
   pages={xii+287},
   isbn={3-03719-010-8},
}

\bib{PWZ1961}{article}{
   author={Penrose, R.},
   author={Whitehead, J. H. C.},
   author={Zeeman, E. C.},
   title={Imbedding of manifolds in euclidean space},
   journal={Ann. of Math. (2)},
   volume={73},
   date={1961},
   pages={613--623},
   issn={0003-486X},
}

\bib{petryshyn1970}{article}{
   author={Petryshyn, W. V.},
   title={A characterization of strict convexity of Banach spaces and other
   uses of duality mappings},
   journal={J. Functional Analysis},
   volume={6},
   date={1970},
   pages={282--291},
}

\bib{phelps57}{article}{
   author={Phelps, R. R.},
   title={Convex sets and nearest points},
   journal={Proc. Amer. Math. Soc.},
   volume={8},
   date={1957},
   pages={790--797},
   issn={0002-9939},
}

\bib{plichko-yost2000}{article}{
   author={Plichko, Anatolij M.},
   author={Yost, David},
   title={Complemented and uncomplemented subspaces of Banach spaces},
   note={III Congress on Banach Spaces (Jarandilla de la Vera, 1998)},
   journal={Extracta Math.},
   volume={15},
   date={2000},
   number={2},
   pages={335--371},
   issn={0213-8743},
}

\bib{salazar2004}{article}{
   author={Salazar, Jos\'{e} M.},
   title={Fixed point index in symmetric products},
   journal={Trans. Amer. Math. Soc.},
   volume={357},
   date={2005},
   number={9},
   pages={3493--3508},
   issn={0002-9947},
}

\bib{schori68}{article}{
   author={Schori, R. M.},
   title={Hyperspaces and symmetric products of topological spaces},
   journal={Fund. Math.},
   volume={63},
   date={1968},
   pages={77--88},
   issn={0016-2736},
}

\bibitem{sims} B. Sims, ``Projection methods in geodesic metric spaces'' (\url{https://carma.newcastle.edu.au/DRmethods/sims13-ACFFTO.pdf}), 7th Asian Conference on Fixed Point Theory and Optimization, Kasetsart University, Thailand.

\bib{shvartsman2004}{article}{
   author={Shvartsman, P.},
   title={Barycentric selectors and a Steiner-type point of a convex body in
   a Banach space},
   journal={J. Funct. Anal.},
   volume={210},
   date={2004},
   number={1},
   pages={1--42},
   issn={0022-1236},
}

\bib{stone1948}{article}{
   author={Stone, A. H.},
   title={Paracompactness and product spaces},
   journal={Bull. Amer. Math. Soc.},
   volume={54},
   date={1948},
   pages={977--982},
   issn={0002-9904},
}

\bibitem{tao2011} T. Tao, ``Brouwer's fixed point and invariance of domain theorems, and Hilbert's fifth problem'' (\url{https://terrytao.wordpress.com/2011/06/13}).

\bib{thele1974}{article}{
   author={Thele, R. L.},
   title={Some results on the radial projection in Banach spaces},
   journal={Proc. Amer. Math. Soc.},
   volume={42},
   date={1974},
   pages={483--486},
   issn={0002-9939},
}

\bib{tuffley2002}{article}{
   author={Tuffley, Christopher},
   title={Finite subset spaces of $S^1$},
   journal={Algebr. Geom. Topol.},
   volume={2},
   date={2002},
   pages={1119--1145},
   issn={1472-2747},
}

\bib{tukia81}{article}{
   author={Tukia, Pekka},
   title={A quasiconformal group not isomorphic to a M\"{o}bius group},
   journal={Ann. Acad. Sci. Fenn. Ser. A I Math.},
   volume={6},
   date={1981},
   number={1},
   pages={149--160},
   issn={0066-1953},
}

\bib{tyson-wu2005}{article}{
   author={Tyson, Jeremy T.},
   author={Wu, Jang-Mei},
   title={Characterizations of snowflake metric spaces},
   journal={Ann. Acad. Sci. Fenn. Math.},
   volume={30},
   date={2005},
   number={2},
   pages={313--336},
   issn={1239-629X},
}

\bib{wheeden-zygmund}{book}{
   author={Wheeden, Richard L.},
   author={Zygmund, Antoni},
   title={Measure and integral},
   series={Pure and Applied Mathematics (Boca Raton)},
   edition={2},
   publisher={CRC Press, Boca Raton, FL},
   date={2015},
   note={An introduction to real analysis},
   pages={xvii+514},
   isbn={978-1-4987-0289-8},
}

\bib{whitehead1949}{article}{
   author={Whitehead, J. H. C.},
   title={Combinatorial homotopy. I},
   journal={Bull. Amer. Math. Soc.},
   volume={55},
   date={1949},
   pages={213--245},
   issn={0002-9904},
}

\bib{wu-wen47}{article}{
   author={Wu, Wen-Ts\"{u}n},
   title={Note sur les produits essentiels sym\'{e}triques des espaces
   topologiques},
   language={French},
   journal={C. R. Acad. Sci. Paris},
   volume={224},
   date={1947},
   pages={1139--1141},
   issn={0001-4036},
}

\end{biblist}
\end{bibdiv}
{\footnotesize
\printindex
}

\hrulefill
\vspace{1cm}
\begin{center}
\bit
\item[] \hspace{0cm}{\Large About the author}
\vspace{0.5cm}
\item[] \textbf{Name of author:}~~ Earnest Akofor
\vspace{0.3cm}

%

\item[] \textbf{Degrees awarded:}~
\vspace{0.2cm}
\bit
\item[] Ph.D. Electrical and Computer Engineering, Syracuse University, 2016
\vspace{0.2cm}

\item[] Ph.D. Physics, Syracuse University, 2010
\vspace{0.2cm}

\item[] Diploma. Mathematical Sciences, AIMS, 2004
\vspace{0.2cm}

\item[] Diploma. High Energy Physics, ICTP, 2003
\vspace{0.2cm}

\item[] B.Sc. Physics, University of Buea, 2001
\eit
\vspace{0.5cm}

\item[] \textbf{Professional experience:}~
\vspace{0.2cm}
\bit
\item[] Instructor \& Teaching Assistant, Syracuse University, 2016-2020
\vspace{0.2cm}

\item[] Peer Review, IJDSN (2015) and IEEE Transactions (2015-2016)
\vspace{0.2cm}

\item[] Research Assistant, Syracuse University, 2011-2016
\vspace{0.2cm}

\item[] Adjunct Instructor, ITT Technical Institute, 2011
\vspace{0.2cm}

\item[] Teaching \& Research Assistant, Syracuse University, 2004-2009
\eit
\eit
\end{center}

\end{document}